\documentclass[a4paper,10pt]{book}
\usepackage[utf8]{inputenc}
\DeclareUnicodeCharacter{00A0}{}
\usepackage{mathrsfs}
\usepackage{lmodern}
 \usepackage[francais]{babel}
 \input{xypic}
\DeclareUnicodeCharacter{00A0}{}
\usepackage{amssymb}
\newcommand{\RR}{{\mathbb R}}
\usepackage{hyperref} 
\usepackage[all]{xy}
\newcommand{\deux}[1]{\refstepcounter{subsection}\label{#1}\medskip\noindent {\bf (\thesubsection)}\hspace{.1cm}}
\newcommand{\trois}[1]{\refstepcounter{subsubsection}\label{#1}\medskip\noindent {\bf
    (\thesubsubsection)}\hspace{.1cm}}

\newcommand{\got}[1]{{\mathfrak #1}}
\newcommand{\red}[1]{\widetilde{#1}}
\renewcommand{\geq}{\geqslant}
\newcommand{\pos}{_+}
\newcommand{\posh}{_+^{\rm hom}}
\newcommand{\homog}{^{\rm hom}}
\newcommand{\Aff}{{\mathbb A}}
\newcommand{\proj}{{\rm Proj}\;}
\newcommand{\ml}[2]{\mathsf L^{#1}_{#2}}
\newcommand{\isol}[2]{|\mathsf L^{#1}_{#2}|}
\newcommand{\Hom}{{\rm Hom}}

\newcommand{\ZZ}{{\Bbb Z}}
\newcommand{\NN}{{\Bbb N}}
\newcommand{\QQ}{{\Bbb Q}}

\newcommand{\VV}{{\Bbb V}}
\renewcommand{\epsilon}{\varepsilon}
\renewcommand{\H}{{\rm H}}
\newcommand{\FF}{{\Bbb F}}
\newcommand{\DD}{{\Bbb D}}
\newcommand{\CC}{{\Bbb C}}
\newcommand{\PP}{{\Bbb P}}
\renewcommand{\leq}{\leqslant}
\renewcommand{\phi}{\varphi}
\newcommand{\sch}[1]{\mathscr #1}
\renewcommand{\hom}{{\rm Hom}}
\newcommand{\ens}{\mathsf{Ens}}
\newcommand{\gp}{\mathsf{Gp}}
\newcommand{\ab}{\mathsf{Ab}}
\newcommand{\ann}{\mathsf{Ann}}
\newcommand{\aalg}{A\text{-}\mathsf{Alg}}
\newcommand{\amod}{A\text{-}\mathsf{Mod}}
\renewcommand{\bmod}{B\text{-}\mathsf{Mod}}
\renewcommand{\top}{\mathsf{Top}}
\newcommand{\point}{\mathsf{TopPt}}
\newcommand{\hompres}{\mathsf{Top}/\mathsf{h}}
\newcommand{\pointhom}{\mathsf{TopPt}/\mathsf {h}}
\newcommand{\ti}{^{\times}}
\newcommand{\limind}{\lim\limits_{\longrightarrow}}
\newcommand{\limproj}{\lim\limits_{\longleftarrow}}
\newcommand{\spec}{{\rm Spec}\;}
\newcommand{\pref}[1]{\mathsf{Pref}_{#1}}
\newcommand{\fasc}[1]{\mathsf{Faisc}_{#1}}
\newcommand{\oxmod}{\sch O_X\text{-}\mathsf{Mod}}
\newcommand{\homfsc}{\underline{\hom}}
\newcommand{\fscend}{\underline{\rm End}\;}
\makeatletter
\let\old@subsection\subsection
\newcommand{\subsection@star}[1]{\old@subsection*{#1}\addcontentsline{toc}{subsection}{\hspace{.2cm}#1}}
\renewcommand{\subsection}{\@ifstar{\subsection@star}{\old@section}}
\makeatother
\title{Introduction à la théorie des schémas
\\ \small Polycopié issu de deux cours du master {\em Mathématiques fondamentales}
de l'UPMC}
\author{Antoine Ducros}
\date{Premier semestre 2013-2014}
\begin{document}

\maketitle
\pagestyle{myheadings}
\tableofcontents

\chapter*{Introduction}
\markboth{Introduction}{Introduction}
\addcontentsline{toc}{chapter}{Introduction}

La théorie des schémas fut développée par Grothendieck et son école dans les années cinquante et soixante. Elle est exposée dans
un gigantesque corpus de textes, répartis en deux familles : 

\medskip
$\bullet$ les {\em Éléments de géométrie algébrique}
(EGA), qui ont été publiés sous forme de volumes entiers des {\em Publications mathématiques de l'IHES} -- mentionnons que EGA I a donné lieu ultérieurement à un livre ; 

$\bullet$ les notes du {\em Séminaire de géométrie algébrique}
du Bois-Marie (SGA), publiées dans la collection {\em Lecture Notes in Mathematics}, et dont la SMF
procède aujourd'hui à la réédition (saisie des manuscrits en LaTeX, corrections, commentaires, etc.). 

\medskip
La première motivation de ce travail d'ampleur exceptionnelle était la mise au point d'outils permettant 
de démontrer la
{\em conjecture de Weil}, ce qui advint effectivement, la pierre finale à l'édifice ayant été apportée
par Pierre Deligne en 1973. À titre purement
culturel\footnote{Ce cours ne permettra pas malheureusement pas d'aborder ni 
même d'effleurer ces questions},
indiquons en quelques mot (une partie de)
ce que dit cette conjecture, ou plutôt ce théorème. 

\subsection*{La conjecture de Weil}
Donnons-nous un système fini~$X$
d'équations polynomiales {\em homogènes} en~$n$ variables, 
à coefficients dans~$\ZZ$. Pour tout nombre
premier~$p$, il définit par réduction modulo~$p$ un système d'équations polynomiales 
homogènes à coefficients
dans~$\FF_p$ ; si~$k$ est une extension finie de~$\FF_p$, on notera~$X(k)$ l'ensemble
des solutions de ce système
dans~$k^n\setminus\{(0,\ldots, 0)\}$ {\em modulo la multiplication par un scalaire non nul}
(comme les équations sont homogènes, si un~$n$-uplet est solution, il en va de même de 
tous ses multiples par un même scalaire). 

De façon analogue, 
on note~$X(\CC)$ l'ensemble des éléments de~$\CC^n\setminus\{(0,\ldots, 0)\}$
solution de~$X$, modulo la multiplication par un scalaire non nul ; il hérite d'une topologie naturelle, 
déduite de celle de~$\CC$ et pour laquelle il est compact. On peut associer
pour tout~$i$ à
l'espace topologique~$X(\CC)$ un~$\QQ$-espace vectoriel de cohomologie
\footnote{Il y a plusieurs définitions possibles, toutes équivalentes ({\em via} les cochaînes singulières, 
{\em via}
les complexes de \v Cech, ou encore {\em via}
la théorie des foncteurs dérivés), que nous n'expliciterons pas ici ; 
nous renvoyons le lecteur intéressé à un cours ou un ouvrage de
topologie algébrique.}
$\H^i(X(\CC),\QQ)$ qui contient des informations sur la
«forme»
de~$X(\CC)$, et dont on démontre
qu'il est de dimension finie. 

On peut définir la dimension algébrique~$d$ de~$X$ ; la dimension topologique de~$X(\CC)$ est alors
égale à~$2d$ (car~$\CC$ est de dimension réelle égale à~$2$ : la droite affine complexe
est un plan réel, une courbe algébrique complexe donne
lieu à une surface de Riemann, une surface algébrique complexe à un espace de dimension réelle 4, etc.).

\medskip
On fait enfin une hypothèse technique sur~$X$, qui en pratique s'avère raisonnable : on suppose qu'il est 
{\em lisse}. Nous ne donnerons pas la définition précise ici ; indiquons
simplement que s'il consiste en une équation~$f$, cela signifie que les dérivées partielles de~$f$ ne s'annulent
pas simultanément sur le lieu des zéros de~$f$ dans~$\CC^n\setminus\{(0,\ldots, 0)\}$ ; en général, cela implique
que~$X(\CC)$ a une structure naturelle de variété différentielle, et même de variété analytique complexe. 

\medskip
Fixons un nombre premier~$p$. Pour tout~$n\geq 1$, il existe à isomorphisme
(non canonique)
près une unique extension~$\FF_{p^n}$ de~$\FF_p$ de degré~$n$ ; notons~$x_{p,n}$ le cardinal
de l'ensemble fini~$X(\FF_{p^n})$, et 
posons
$$Z_p=\exp\left(\sum_{n\geq 1} \frac{x_{p,n} T^n}n\right)\in \QQ((T)).$$

On démontre alors (Weil, Dwork, Grothendieck, Deligne) les assertions suivantes : 

\medskip
1) Pour tout~$p$, la série $Z_p$ est une fraction rationnelle. 

2) Pour tout~$p$ suffisamment grand, on peut plus précisément écrire
$$Z_p=\frac{\prod\limits_{0\leq i\leq 2d, \;i\;{\rm impair}} R_{i,p}}{\prod\limits_{0\leq i\leq 2d,\; i\;{\rm pair}} R_{i,p}}$$
où~$R_{i,p}$ est pour tout~$i$ un polynôme unitaire à coefficients dans~$\ZZ$ de degré~$\dim_{\QQ} \H^i(X(\CC),\QQ)$ 
dont toutes les racines complexes ont pour module~$p^{i/2}$. 

\medskip
On voit en particulier qu'il existe un lien profond,
lorsque~$p$ est assez grand, entre le nombre 
de solutions de~$X$ dans les~$\FF_{p^n}$ (pour~$n$ variable)
et la topologie de~$X(\CC)$.  

\subsection*{L'intérêt des schémas} 

Au-delà de ce succès majeur qu'a représenté la preuve 
de la conjecture de Weil, la théorie des schémas s'est imposée
comme un outil 
à peu près indispensable pour qui souhaite faire de la géométrie
algébrique sur un corps,
et plus encore sur un anneau, quelconques. 

\medskip
Elle a l'inconvénient, comme nous le verrons, d'être d'un accès ardu : la mise en place
des définitions et la démonstration des propriétés de base sont longues et
parfois délicates. 

\medskip
Mais une fois
franchis ces premiers obstacles un peu âpres, elle s'avère d'une extrême souplesse. 
Et elle a un 
immense avantage : elle apporte de l'intuition géométrique dans des situations qui
pouvaient
{\em a priori}
sembler purement algébriques, essentiellement parce qu'elle permet de 
penser à
{\em n'importe quel anneau}
comme à un anneau de «fonctions»
sur un objet géométrique.

\medskip
Par exemple, 
reprenons le système d'équations~$X$
à coefficients dans~$\ZZ$ considéré au paragraphe précédent. 
La théorie des schémas lui associe
une sorte de {\em fibration}
dont les différentes fibres sont, {\em grosso modo},
les variétés algébriques obtenues à partir de~$X$ d'une part
en le réduisant modulo~$p$ pour
chacun des nombres premiers~$p$, d'autre part
en le voyant comme un système d'équations à coefficients dans~$\QQ$. 
Cette fibration fournit ainsi
un certain {\em liant}
entre les différentes caractéristiques, qui s'avère
très utile pour
comprendre
dans quelle mesure
ce qui se passe modulo~$p$ peut
avoir un rapport avec ce qui se passe
en caractéristique nulle -- et rend moins mystérieuse
la relation entre cardinal de~$X(\FF_{p^n})$ et 
topologie de~$X(\CC)$.

\subsection*{Les outils indispensables}

La théorie des schémas repose de manière cruciale sur un certain nombre
d'outils et notions, auxquels nous consacrons une première partie, 
elle-même divisée en trois chapitres distincts.  

\medskip
$\bullet$ La premier porte sur
les {\em catégories}. Comme vous le verrez, on ne vous y présente
pas véritablement une théorie\footnote{On n'y établit pour ainsi dire qu'un seul énoncé, 
le {\em lemme de Yoneda},
dont la preuve est essentiellement triviale, même si elle
peut être très déroutante à la première lecture.}, mais plutôt un {\em langage}
très commode. Il
permet, en dégageant un certain nombre de propriétés
formelles qui leur sont communes, de donner une description unifiée de
situations rencontrées dans des domaines extrêmement divers. On peut {\em en
principe} l'utiliser dans à peu près n'importe quelle 
branche des mathématiques ; 
{\em en pratique}, les géomètres algébristes à la Grothendieck
en sont particulièrement friands.

\medskip
$\bullet$ Le seconde est le
plus difficile sur le plan technique. Il est consacré à {\em l'algèbre commutative},
c'est-à-dire à
l'étude des anneaux commutatifs, et des idéaux {\em de}
et 
modules {\em sur}
ces derniers. L'algèbre commutative joue en géométrie algébrique un rôle
absolument crucial, 
analogue à celui de l'analyse réelle en géométrie différentielle : elle constitue en quelque sorte
la partie {\em locale}
de la théorie. 

Nous commençons  par présenter des notions
et résultats très généraux : localisation, anneaux locaux et lemme de Nakayama,
produit tensoriel, modules projectifs, 
algèbres finies et entières, dimension de Krull, 
lemme de {\em going-up}. 
Puis nous en venons à des théorèmes plus spécifiques
et nettement plus délicats, qui concernent les algèbres de type fini sur un corps :
normalisation
de Noether, {\em Nullstellensatz}, et
calcul de la dimension de Krull
d'une telle algèbre. 

\medskip
$\bullet$ Le dernier
présente les définitions et propriétés de base des
{\em faisceaux} sur un espace topologique. Ceux-ci ont été initialement 
introduits par Leray en topologie algébrique
et c'est Serre qui, dans son 
article fondateur {\em Faisceaux algébriques cohérents}, a le premier
mis en évidence les services qu'ils pouvaient rendre en géométrie algébrique ; Grothendieck
les a ensuite placés au cœur de toute sa théorie. 

\medskip
Celle-ci repose ainsi de façon essentielle sur la notion d'{\em espace localement annelé}\footnote
{Les schémas sont ainsi définis comme des espaces localement annelés satisfaisant une condition
supplémentaire ; signalons par ailleurs que les objets géométriques
plus classiques (variétés différentielles, variétés analytiques complexes ou réelles...)
sont aussi de manière naturelle des espaces localement annelés.}
(c'est un espace topologique muni d'un faisceau d'un certain type),
qui est au cœur de notre chapitre faisceautique\footnote{Le lecteur
trouvera sans doute
avec raison
que l'adjectif «faisceautique»
est très laid ; d'un point de vue strictement
linguistique, le bon terme aurait probablement été «fasciste», 
mais il n'est évidemment
plus utilisable.}, lequel se conclut 
par l'étude de certaines propriétés 
de {\em faisceaux de modules}
particuliers
sur un espace localement annelé, et notamment de ceux
qui sont {\em localement libres}
de rang~$1$ et jouent un rôle absolument central en géométrie algébrique. 

\setcounter{chapter}{-1}
\chapter{Prérequis et rappels}

\section{Anneaux}
\deux{conv-ann-comm}
{\bf Convention}.
Dans tout ce qui suit,
 et sauf mention  expresse
 du contraire, «anneau»
 signifiera «anneau commutatif unitaire», 
 «algèbre»
 signifiera «algèbre commutative unitaire», et 
 un morphisme d'anneaux ou d'algèbres sera toujours supposé
 envoyer l'unité de la source sur celle du but.

\deux{pre-requis} Le lecteur sera supposé familier avec les définitions
 d'anneau, d'idéal
 et d'anneau quotient.... ainsi qu'avec les
 propriétés élémentaires de ces objets, 
 que nous ne rappellerons pas ici pour la plupart. Nous allons toutefois
 insister sur quelques points sans doute connus, mais qui sont importants
 et au sujet desquels on peut commettre facilement quelques erreurs.  
 
 \trois{ann-nul} Dans la définition d'un anneau~$A$, on n'impose pas à~$1$
 d'être différent de~$0$. En fait, l'égalité~$1=0$ se produit dans un et seul cas,
 celui où~$A$ est {\em l'anneau nul}
 $\{0\}$. 
 
 \trois{def-ati}
 Si~$A$ est un anneau, on notera~$A\ti$ l'ensemble des éléments
 inversibles de~$A$ ; il est stable par multiplication
 et~$(A\ti, \times)$ est un groupe 
 
 \deux{ann-integre} Un anneau~$A$ est dit
 {\em intègre} s'il est non nul
 et si l'on a pour tout couple~$(a,b)$
 d'éléments de~$A$ l'implication
 $$(ab=0)\Rightarrow (a=0\;{\rm ou}\;b=0).$$
 
 \trois{attention-non-nul}
 On prendra garde de ne jamais oublier 
 de vérifier la  première de ces deux conditions : {\em un anneau intègre est 
 par définition non nul} (l'expérience a montré qu'on avait tout intérêt à imposer cette restriction
 pour éviter une profusion de cas particuliers à distinguer dans les définitions, énoncés  et
 démonstrations ultérieurs). 
 
 \trois{ann-integre-bourb} Le lecteur amateur de facéties bourbakistes appréciera certainement
 la définition alternative suivante : un anneau~$A$ est intègre si et seulement si 
 {\em tout produit fini d'éléments non nuls de~$A$ est non nul.}
 Elle contient en effet la non-nullité
 de~$A$, puisqu'elle implique
 que l'unité~$1$, qui n'est autre que le produit {\em vide}
 d'éléments de~$A$, est non nulle. 
 
 \deux{def-corps} On dit qu'un anneau~$A$ est un corps s'il est non nul 
 et si tout élément non nul de~$A$ est inversible ; il revient au même de
 demander que~$A$ ait exactement deux idéaux, à savoir~$\{0\}$ et~$A$. 
 Si~$A$ est un corps, il est intègre et~$A\ti =A\setminus\{0\}$. 
 
 \deux{mor-corps}
 Soit~$f$ un morphisme d'un corps~$K$ vers un anneau
 {\em non nul}
 $A$. Comme~$A$ est non nul, $1\notin {\rm Ker}\;f$ ;
 puisque les seuls idéaux de~$K$ sont~$\{0\}$ et~$K$,
 il vient~${\rm Ker}\;f=\{0\}$ et~$f$ est injectif. 
 
 \medskip
 En particulier, tout morphisme de corps est injectif. 
 
 \deux{alg-rap}
 Soit~$A$ un anneau. Une
 {\em $A$-algèbre}
 est un anneau~$B$ muni d'un morphisme~$f : A\to B$. Bien
 que~$f$ fasse partie des données, il sera très souvent omis
 (on dira simplement «soit~$B$
 une~$A$-algèbre»). Il arrivera même que l'on écrive
 abusivement~$a$ au lieu de~$f(a)$ pour~$a\in A$ ; 
mais cette entorse à la rigueur peut
 être dangereuse, surtout lorsque~$f$ n'est pas injective : si on 
 la commet, il faut en avoir conscience et y mettre fin lorsque la situation
 l'exige.

 \deux{rem-nilrad}
Soir~$A$ un anneau. Un élément~$a$ de~$A$
est dit
{\em nilpotent}
s'il existe~$n\geq 0$ tel que~$a^n=0$. L'ensemble
des éléments nilpotents 
de~$A$ est un idéal de~$A$ (nous vous laissons la preuve
en exercice, appelé le {\em nilradical}
de~$A$. On dit que~$A$ est
{\em réduit}
si son nilradical est nul, c'est-à-dire encore si~$A$
n'a pas d'élément idempotent non trivial. 

 \section{Modules} 
 \markboth{Algèbre commutative}{Modules}

 \deux{pas-rap-modules}
 Soit~$A$ un anneau. Nous ne rappellerons pas ici les définitions des objets de base de la théorie des 
 $A$-modules, à savoir les~$A$-modules eux-mêmes, les sous-modules, les applications~$A$-linéaires, 
 les familles libres et génératrices, les bases, les supplémentaires.... Ce sont
 {\em mutatis mutandis}
 les mêmes qu'en algèbre linéaire. 
 
 \deux{conv-somme}
 Il arrivera souvent dans la suite qu'on manipule des expressions de
 la forme~$\sum_{i\in I} m_i$, où les~$m_i$ sont des éléments
 d'un~$A$-module~$M$ fixé. Il sera toujours implicitement supposé, 
 dans une telle écriture, 
 que {\em presque tous les~$m_i$ sont nuls.} 
 Elle n'aurait sinon aucun sens : {\em en algèbre, on ne sait faire
 que des sommes finies} ; pour donner un sens à des sommes infinies, 
 il est nécessaire d'introduire des structures de nature topologique.

 \deux{attention-ev}
 {\bf Attention !} On prendra garde que certains énoncés usuels 
 portant sur les espaces vectoriels deviennent faux en général pour les modules sur un anneau quelconque
 (ce qui empêche leurs preuves de s'étendre à ce nouveau contexte est
 le plus souvent qu'elles font appel à un moment
 ou un autre à l'inversion d'un scalaire non nul). 
 Donnons quelques exemples. 
 
 \trois{pas-base} Il est faux en général qu'un module possède une base. Par exemple, 
 le~$\ZZ$-module~$\ZZ/2\ZZ$ n'en possède pas. En effet, s'il en admettait une elle serait non 
 vide (puisqu'il est non nul), et comprendrait donc au moins un élément qui serait
 annulé par~$2$, 
 contredisant ainsi
 sa liberté.

 \trois{rang-module}
 Soit~$A$ un anneau. On dit qu'un~$A$-module~$M$
 est {\em libre}
 s'il possède une base. On vient de voir que ce n'est pas automatique ; 
 mais lorsque c'est le cas {\em et lorsque~$A$
 est non nul} on démontre que, comme en algèbre linéaire, 
 toutes les bases de~$M$ ont même cardinal, appelé
 {\em rang}
 de~$M$.

 \medskip
 Il faut faire attention au cas de l'anneau nul~$\{0\}$. Le
 seul module sur celui-ci est le module trivial~$\{0\}$, 
 et {\em toute}
 famille~$(e_i)_{i\in I}$ d'éléments de ce module
 (qui vérifie nécessairement~$e_i=0$ pour tout~$i$) en est une base, 
 indépendamment du cardinal de~$I$ (il peut être vide, fini, infini dénombrable ou non, etc.). 
 Nous laissons au lecteur qu'amusent les manipulations logiques dans les cas un peu extrêmes le 
 soin de prouver cette assertion. 
 
 \medskip
 Notez par contre que sur un anneau non nul, un élément~$e$ d'une famille libre
 n'est jamais nul (sinon, il satisferait la relation non triviale $1\cdot e=0$). 
 
 \medskip
 Pour éviter de fastidieuses distinctions de cas, on se permet d'appeler
 {\em module libre de rang~$n$} (où~$n$ est un entier, ou même un cardinal)
 tout module libre ayant une base de cardinal~$n$ : cela permet d'inclure
 le module nul sur l'anneau nul, qui est ainsi libre de tout rang. 
 
 \trois{inj-pas-sur}
 Soit~$A$ un anneau et soit~$M$ un~$A$-module libre de rang fini.
 Il est faux en général qu'un endomorphisme injectif de~$M$ 
 lui-même soit bijectif\footnote{Nous
 verrons par
 contre un peu plus loin qu'un
 endomorphisme {\em surjectif}
 d'un tel~$M$
 est toujours bijectif.}. Par exemple, 
$\ZZ$
est un~$\ZZ$-module libre de rang~$1$ sur lui-même, 
et la multiplication par~$2$
en est un endomorphisme injectif non surjectif. 

\trois{pas-supplem}
Soit~$A$ un anneau. Il est faux en général que tout sous-module
d'un~$A$-module~$M$ admette un supplémentaire
dans~$M$, même si~$M$
est libre. Par exemple, le lecteur pourra démontrer à titre d'exercice que
le sous-module~$2\ZZ$ de~$\ZZ$ n'a pas de supplémentaire dans~$\ZZ$.

\deux{type-pres-fini}
Soit~$A$ un anneau. Un~$A$-module~$M$ est dit {\em de type fini}
s'il possède une famille génératrice finie, c'est-à-dire encore s'il existe un entier~$n$
et une surjection linéaire~$A^n\to M$. 

Il est dit {\em de présentation finie} s'il existe une telle surjection {\em possédant un noyau 
de type fini}. 

\trois{type-fin-pres-fin}
Tout~$A$-module de type fini est de présentation finie. La réciproque est vraie
si~$A$ est noethérien, car on démontre que dans ce cas tout sous-module d'un~$A$-module
de type fini est de type fini ; elle est fausse en général (pour un contre-exemple, {\em cf.}
\ref{rem-pas-pf} {\em infra}).

\trois{libre-tf}
Si~$M$ est un~$A$-module
libre et de type fini, il est de rang fini. Pour le voir, on choisit une base~$(e_i)_{i\in I}$
de~$M$. Comme~$M$ est de type fini, il est engendré par une famille finie de vecteurs, et chacun 
d'eux est combinaison linéaire d'un nombre fini de~$e_i$. Il existe donc un ensemble fini 
d'indices~$J\subset I$ tel que~$(e_i)_{i\in J}$ soit génératrice, et donc soit une base de~$M$
(si~$A\neq \{0\}$ on a alors nécessairement~$I=J$, car un élément~$e_i$ avec~$i\notin J$ est forcément
non nul et ne peut donc être combinaison linéaire des~$e_i$ avec~$i\in J$ ; mais si~$A=\{0\}$ notez que~$I$
peut contenir strictement~$J$, {\em cf.}~\ref{rang-module}). 

\deux{supplem}
{\bf Sommes directes externes et internes.}
Soit~$A$ un anneau.

\trois{somme-directe-interne}
{\em La somme directe interne.}
Soit~$M$ un~$A$-module et~soit~$(M_i)$
une famille de sous-modules de~$M$. La {\em somme des~$M_i$}
est le sous-module de~$M$ constitué des éléments de la forme
$\sum m_i$ où~$m_i\in M_i$ pour tout~$i$ ; on le note~$\sum M_i$. 

\medskip
On dit que la somme des~$M_i$ est {\em directe}, et l'on écrit~$\sum M_i=\bigoplus M_i$, si 
pour tout élément~$m$ de~$\sum M_i$ l'écriture~$m=\sum m_i$ est unique. Il suffit de le vérifier
pour~$m=0$, c'est-à-dire de s'assurer que~$(\sum m_i=0)\Rightarrow (\forall i\;m_i=0)$. 

 \trois{somme-directe-externe}
 {\em La somme
 directe externe.}
 Soit~$(M_i)_{i\in I}$ une famille de~$A$-modules, qui ne sont pas
 {\em a priori}
 plongés dans un même~$A$-module. 
 Soit~$N$
 le sous-module de~$\prod M_i$ formé des familles~$(m_i)$
 telles que~$m_i=0$ pour presque tout~$i$. Pour tout~$i\in I$, on dispose d'une injection naturelle
 $h_i: M_i\hookrightarrow N$
 qui envoie un élément~$m$ sur la famille~$(m_j)$ avec~$m_j=0$ si~$j\neq i$ et~$m_i=m$. Il est immédiat
 que~$N=\bigoplus h_i(M_i)$. On a ainsi construit un module qui contient une copie de chacun des~$M_i$, et est égal
 à la somme
 directe desdites copies. On dit que~$N$ est la somme directe {\em externe}
 des~$M_i$, et on le note encore~$\bigoplus M_i$.

\trois{liens-interne-externe}
{\em Liens entre les sommes directes interne et externe.}
Soit~$M$ un~$A$-module et soit~$(M_i)$ une famille de sous-modules
de~$M$. On dispose pour tout~$i$ de l'inclusion~$u_i : M_i\hookrightarrow M$. 
La famille des~$u_i$ définit une application linéaire
$$u : (m_i)\mapsto \sum_i u_i(m_i)$$
de la somme directe {\em externe}~$\bigoplus M_i$ vers~$M$. On vérifie aussitôt
que les~$M_i$ sont en somme directe dans~$M$ au sens de~\ref{somme-directe-interne}
si et seulement si~$u$ est injective, et que~$M=\bigoplus M_i$ au sens de~\ref{somme-directe-interne}
si et seulement si~$u$
est bijectif. 

\part{Les outils de la géométrie algébrique} 
\chapter{Le langage des catégories}

\section{Définitions et premiers exemples}
\markboth{Le langage des catégories}{Définitions, exemples}

\deux{defcat} {\bf Définition.}
Une {\em catégorie}~$\mathsf C$ consiste en les données suivantes. 

\medskip
$\bullet$ Une classe d'objets~${\rm Ob}\;\mathsf C$.

$\bullet$ Pour tout couple~$(x,y)$ d'objets de~$\mathsf C$, un ensemble
$\Hom_{\mathsf C}(x,y)$ dont les éléments sont appelés les {\em 
morphismes de source~$x$ et de but~$y$}, ou encore les {\em morphismes de
$x$ vers~$y$}. 

$\bullet$ Pour tout~$x\in {\rm Ob}\;\mathsf C$, un élément~${\rm Id}_x$ de~$\Hom_{\mathsf C}(x,x)$. 

$\bullet$ Pour tout triplet~$(x,y,z)$ d'objets de~$\mathsf C$, une application

$$\Hom_{\mathsf C}(x,y)\times \Hom_{\mathsf C}(y,z)\to \Hom_{\mathsf C}(x,z), \;\;(f,g)\mapsto g\circ f\;.$$

Ces données sont sujettes aux deux axiomes suivants. 

\medskip
$\diamond$ {\em Associativité} : on a~$(f\circ g)\circ h=f\circ (g\circ h)$ pour tout triplet~$(f,g,h)$ de morphismes
tels que ces compositions aient un sens. 

$\diamond$ {\em Neutralité des identités} : pour tout couple~$(x,y)$ d'objets de~$\mathsf C$ et tout morphisme~$f$ de~$x$ vers~$y$, 
on a~$f\circ {\rm Id}_x={\rm Id}_y\circ f=f$.

\deux{cat-fondements}
{\bf Commentaires.} La définition ci-dessus est 
un peu vague : nous n'avons pas précisé ce que signifie «classe»~-- ce n'est
pas une notion usuelle de théorie des ensembles. Cette imprécision est volontaire : dans le cadre 
de ce cours, il ne nous a pas paru souhaitable d'entrer dans le détail des problèmes de fondements 
de la théorie des catégories. Mais il est possible de développer celle-ci rigoureusement, de plusieurs façons différentes : 

- on peut travailler avec une variante de la théorie des ensembles dans laquelle la notion de classe a un sens (une classe
pouvant très bien ne pas être un ensemble) ; 

- on peut imposer à~${\rm Ob}\;\mathsf C$ d'être un ensemble.

\medskip
Comme on le verra, c'est plutôt la première approche que nous suivrons implicitement : dans les exemples que nous
donnons ci-dessous, ${\rm Ob}\; C$ n'est en général pas un ensemble. 

Pour suivre la deuxième, il faudrait modifier
la définition de toutes nos catégories en ne gardant que les objets qui appartiennent à un certain ensemble fixé au préalable, 
et absolument gigantesque : il
doit être assez gros pour qu'on puisse y réaliser toutes les constructions du cours. C'est ce que Grothendieck
et son école ont fait dans
SGA IV (parce que
certaines questions abordées dans cet ouvrage
requièrent de manière {\em impérative}
de rester dans le cadre
ensembliste traditionnel) ; ils qualifient ce
type d'ensembles gigantesques 
d'{\em univers}.

\deux{ex-cat}
{\bf Exemples de catégories.}

\trois{cat-class} Commençons par des exemples classiques, qui mettent en jeu de «vrais»
objets et de «vrais» morphismes. 

\medskip
$\bullet$ La catégorie $\ens$ : ses objets sont les ensembles, et ses morphismes les applications. 

$\bullet$ La catégorie~$\gp$ : ses objets sont les groupes, et ses morphismes les morphismes de groupes. 

$\bullet$ La catégorie~$\ab$ : ses objets sont les groupes {\em abéliens}, et ses morphismes les morphismes
de groupes. 

$\bullet$ La catégorie~$\ann$ : ses objets sont les anneaux, et ses morphismes les morphismes d'anneaux. 

$\bullet$ La catégorie~$\amod$, où~$A$ est un anneau : ses objets sont les~$A$-modules, et ses morphismes
les applications~$A$-linéaires. 

$\bullet$ La catégorie~$\top$ : ses objets sont les espaces topologiques, et ses morphismes sont
les applications continues. 

\trois{cat-semiclass}
Partant d'une catégorie, on peut en définir d'autres par un certain nombre de procédés standard. 

\medskip
\begin{itemize}
\item[$\bullet$] Commençons par
un exemple abstrait. Soit~$\mathsf C$ une catégorie et soit~$S$ un objet de~$\mathsf C$. On note~$\mathsf C/S$ la catégorie définie comme suit. 

\begin{itemize}
\item[$\diamond$] Ses objets sont les couples~$(X,f)$ 
où~$X$ est un objet de~$\mathsf C$ est où~$f: X\to S$ est un morphisme.

\item[$\diamond$] Un morphisme de~$(X,f)$ vers~$(Y,g)$ est un morphisme~$\phi : X\to Y$ te que le diagramme
$$\diagram X \rrto^\phi\drto_f&&Y\dlto_g
\\&S&\enddiagram$$
commute. 
\end{itemize}

\item[$\bullet$] La construction duale existe : on peut définir~$S\backslash \mathsf C$ 
comme la catégorie dont les objets sont les couples~$(X,f)$ 
où~$X$ est un objet de~$\mathsf C$ est où~$f: S\to X$ est un morphisme, et dont les flèches sont définies
comme le lecteur imagine (ou devrait imaginer).

\item[$\bullet$] Donnons 
deux exemples
explicites 
de catégories de la forme~$S\backslash \mathsf C$. 

\medskip

\begin{itemize}

\item[$\diamond$] Si~$\mathsf C=\ann$ et si~$A\in {\rm Ob}\;\mathsf C$ alors~$A\backslash \mathsf C$ n'est autre que la catégorie
$\aalg$ des~{\em $A$-algèbres}.

\item[$\diamond$] Si~$\mathsf C=\top$ et si~$S=\{*\}$ la catégorie~$S\backslash \mathsf C$
est appelée catégorie des {\em espaces topologiques pointés}
et sera notée~$\point$. Comme se donner
une application continue
de~$\{*\}$ dans un espace topologique~$X$ revient à choisir un point de~$X$, 
la catégorie~$\point$ peut se décrire comme suit : 

\medskip
\begin{itemize}

\item[-] ses objets 
sont les couples~$(X,x)$ 
où~$X$ est un espace topologique et où~$x\in X$ (d'où son nom) ; 

\item[-] un morphisme d'espaces topologiques pointés
de~$(X,x)$ vers~$(Y,y)$ est une application 
continue~$\phi : X\to Y$ telle que~$\phi(x)=y$.
\end{itemize}
\end{itemize}

\end{itemize}

\trois{cat-quotient}
On peut aussi, partant d'une catégorie, conserver ses objets mais ne plus considérer ses morphismes
que modulo une certaine relation d'équivalence. Nous n'allons pas détailler le formalisme général, mais 
simplement illustrer ce procédé par un exemple. Si~$X$ et~$Y$ sont deux espaces topologiques
et si~$f$
et~$g$ sont deux applications continues de~$X$ vers~$Y$, on dit qu'elles sont {\em homotopes}
s'il existe une application~$h$ continue de~$[0;1]\times X$ vers~$Y$ telle que~$h(0,.)=f$
et~$h(1,.)=g$ ; on définit ainsi
ainsi une relation d'équivalence sur l'ensemble des applications continues de~$X$ vers~$Y$.

\medskip
On construit alors la catégorie~$\hompres$ des
{\em espaces topologiques à homotopie près}
comme suit : ses objets sont les espaces topologiques ; et si~$X$ et~$Y$ 
sont deux espaces topologiques, $\hom_{\hompres}(X,Y)$ est le quotient 
de l'ensemble des applications continues de~$X$ vers~$Y$ par la relation d'homotopie
(les morphismes de~$X$ vers~$Y$ dans~$\hompres$ ne sont donc plus 
tout à fait de «vrais»~morphismes). 

\medskip
On peut combiner cette construction avec celle des espaces topologiques
pointés, et obtenir ainsi la catégorie $\pointhom$
des {\em espaces topologiques pointés
à homotopie près} : ses objets sont les espaces topologiques pointés ; 
et si~$(X,x)$ et~$(Y,y)$ 
sont deux espaces topologiques pointés, $\hom_\pointhom(X,Y)$ est le quotient 
de l'ensemble~$\hom_{\point}((X,x),(Y,y))$ par la relation d'homotopie, qui est définie dans ce contexte exactement comme
ci-dessus avec la condition supplémentaire~$h(t,x)=y$ pour tout~$t$. 

\trois{cat-decret} On peut également construire des catégories par décret, sans chercher 
à donner une interprétation tangible des objets et morphismes ; donnons quelques exemples. 

\medskip
\begin{itemize}
\item[$\bullet$]  Soit~$G$ un groupe. On lui associe traditionnellement deux catégories : 

\begin{itemize}
\item[$\diamond$] la catégorie~$BG$ ayant un seul objet~$*$ avec~$\hom(*,*)=G$
(la composition est définie comme étant égale à la loi interne de~$G$) 

\item[$\diamond$]  la catégorie~$\widetilde{BG}$
telle que~${\rm Ob}\;\widetilde{BG}=G$ et telle qu'il y ait un et un seul morphisme entre deux objets donnés 
de~$\widetilde{BG}$. 
\end{itemize}

\item[$\bullet$] Soit~$k$ un corps. 
On définit comme suit la catégorie~$\VV_k$ .

\begin{itemize}
\item[$\diamond$] ${\rm Ob}\;\VV_k={\mathbb N}$.

\item[$\diamond$]  $\hom_{\VV_k}(m,n)=M_{n,m}(k)$ pour tout~$(m,n)$, 
la composition de deux morphismes
étant définie comme égale au produit des deux matrices correspondantes.

\end{itemize}
\end{itemize} 

\trois{cat-opp}
{\bf Catégorie opposée.} Si~$\mathsf C$ est une catégorie, on 
définit sa catégorie {\em opposée}
~$\mathsf C^{\rm op}$ comme suit : ${\rm Ob}\;\mathsf C^{\rm op}={\rm Ob}\; C$, 
et~$\hom_{C^{\rm op}}(X,Y)=\hom_{\mathsf C}(Y,X)$ pour tout couple~$(X,Y)$
d'objets de~$\mathsf C$. 

\deux{def-iso} Soit~$\mathsf C$
une catégorie et soient~$X$
et~$Y$ deux objets de~$\mathsf C$. 

\trois{def-endo} Un {\em endomorphisme}
de~$X$ est un élément de~$\hom_{\mathsf C}(X,X)$.

\trois{def-iso} On dit qu'un 
morphisme~$f:X\to Y$ est un {\em isomorphisme}
s'il existe un morphisme~$g: Y\to X$ tel que~$f\circ g={\rm Id}_Y$
et~$g\circ f={\rm Id}_X$. 
On vérifie immédiatement que si un tel~$g$ existe, il est unique ; et on le
note alors en général~$f^{-1}$. 

\medskip
Un {\em automorphisme}
de~$X$ est un isomorphisme de~$X$ vers~$X$. 

\trois{def-fleche} Il arrivera souvent que l'on emploie
le terme
{\em flèche}
au lieu de morphisme. 

\section{Foncteurs}
\markboth{Le langage des catégories}{Foncteurs}

\deux{def-fonct} {\bf Définition.}
Soient~$\mathsf C$ et~$\mathsf D$ deux catégories. Un {\em foncteur}
$F$ de~$\mathsf C$ vers~$\mathsf D$ est la donnée : 

\medskip
$\bullet$ pour tout~$X\in {\rm Ob}\;\mathsf C$, d'un objet~$F(X)$ de~$\mathsf D$ ; 

$\bullet$ pour tout morphisme~$f: X\to Y$ de~$\mathsf C$,
d'un morphisme~$F(f) : F(X)\to F(Y)$ de~$\mathsf D$. 

\medskip
On impose de plus les propriétés suivantes : 

\medskip
$\bullet$ $F({\rm Id}_X)={\rm Id}_{F(X)}$ pour tout~$X\in {\rm Ob}\;\mathsf C$ ; 

$\bullet$ $F(f\circ g)=F(f)\circ F(g)$ pour tout couple~$(f,g)$ de flèches composables
de~$\mathsf C$. 

\deux{comment} {\bf Commentaires.}

\trois{fonct-iso}
Il est immédiat qu'un foncteur transforme un isomorphisme en un isomorphisme. 

\trois{rem-contra} La notion définie au~\ref{def-fonct}
est en fait la notion de foncteur {\em covariant}. Il existe une notion  de
foncteur {\em contravariant} ; la définition est la même, 
à ceci près que si~$f\in \hom_{\mathsf C}(X,Y)$ 
alors~$F(f)\in \hom_{\mathsf D}(F(Y), F(X))$ (en termes imagés, $F$ renverse le sens des flèches), 
et que~$F(f\circ g)=F(g)\circ F(f)$ pour tout couple~$(f,g)$ de flèches composables. 

\trois{contra-op} La plupart des résultats et notions que nous présenterons ci-dessous
seront relatifs aux foncteurs {\em covariants}, mais ils se transposent
{\em mutatis mutandis} au
cadre des foncteurs contravariants (nous laisserons ce soin
au lecteur, et utiliserons librement à l'occasion
ces transpositions) : il suffit de renverser certaines flèches et/ou l'ordre 
de certaines compositions. 

Cette assertion peut paraître imprécise, mais on peut lui donner un sens rigoureux
-- et la justifier --
en remarquant qu'un foncteur contravariant de~$\mathsf C$ vers~$\mathsf D$ n'est autre, par définition, qu'un foncteur
covariant de~$\mathsf C^{\rm op}$ vers~$\mathsf D$ (ou de~$\mathsf C$ vers~$\mathsf D^{\rm op}$).

\deux{ex-foncteurs}
{\bf Exemples.}

\trois{fonct-idiot} Sur toute catégorie~$\mathsf C$ on dispose d'un foncteur identité~${\rm Id}_{\mathsf C} :\mathsf C\to C$
défini par les égalités~${\rm Id}_{\mathsf C}(X)=X$ et~${\rm Id}_{\mathsf C}(f)=f$ pour tout objet~$X$ et toute flèche~$f$ 
de~$\mathsf C$. 

Si~$\mathsf C,D$ et~$E$ sont trois catégories, si~$F$ est un foncteur de~$\mathsf C$ vers~$\mathsf D$ et si~$G$ est un foncteur de~$\mathsf D$
vers~$E$ on définit de façon évidente le foncteur composé~$G\circ F : \mathsf C\to E$. La composition des foncteurs est associative, 
et les foncteurs identité sont neutres pour celle-ci. 

Le composé de deux foncteurs de même variance est covariant, celui de deux foncteurs de variances opposées est contravariant. 

\trois{fonct-oubli} {\bf Les foncteurs d'oublis.}
Ce sont des foncteurs covariants dont l'application
revient, comment leur nom
l'indique, à oublier certaines des structures en jeu. 

Par exemple, on dispose d'un foncteur d'oubli de~$\gp$ vers~$\ens$ : il associe à un groupe l'ensemble sous-jacent, 
et à un morphisme de groupes l'application ensembliste sous-jacente. 
On dispose de même d'un foncteur d'oubli de~$\ann$ vers~$\ens$, de~$\top$ vers~$\ens$.... ou encore, un anneau commutatif unitaire
$A$ étant donné, de~$\amod$
vers~$\ab$ (on n'oublie alors qu'une partie de la structure). 

\trois{fonct-rep} {\bf Les deux foncteurs associés à un objet.}
Soit~$\mathsf C$ une catégorie et soit~$X$ un objet de~$\mathsf C$. On lui associe naturellement
deux foncteurs de~$\mathsf C$ vers~$\ens$ : 

\medskip
$\bullet$ Le foncteur covariant~$h_X$, qui envoie un objet~$Y$ sur~$\hom_{\mathsf C}(X,Y)$ 
et une flèche~$f: Y\to Y'$ sur l'application~$\hom_{\mathsf C}(X,Y)\to \hom_{\mathsf C}(X,Y'), g\mapsto f\circ g$. 

$\bullet$ Le foncteur contravariant~$h^X$, qui envoie un objet~$Y$ sur~$\hom_{\mathsf C}(Y,X)$ 
et une flèche~$f: Y\to Y'$ sur l'application~$\hom_{\mathsf C}(Y',X)\to \hom_{\mathsf C}(Y,X), g\mapsto g\circ f$.

\trois{dualite} Soit~$A$ un anneau commutatif unitaire. En appliquant la construction
du~\ref{fonct-rep}
ci-dessus à l'objet~$A$ de~$\amod$, on obtient un foncteur contravariant
$h^A : \amod \to \ens$ qui envoie~$M$ sur~$\hom_{\amod}(M,A)$.

En fait, $\hom_{\amod}(M,A)$ n'est pas un simple ensemble : il possède une structure
naturelle de~$A$-module, et la formule~$M\mapsto \hom_{\amod}(M,A)$ définit
un foncteur contravariant de~$\amod$ dans elle-même, souvent noté~$M\mapsto M^\vee$
(en termes pédants, le foncteur~$h^A$ ci-dessus est la composée
de~$M\mapsto M^\vee$ avec le foncteur oubli de~$\amod$
vers~$\ens$).

\trois{groupe-fond} {\bf Le groupe fondamental}. 
On définit en topologie algébrique un foncteur $(X,x)\mapsto \pi_1(X,x)$, qui va
de la catégorie~$\point$ vers celle des groupes (on appelle~$\pi_1(X,x)$
le {\em groupe fondamental}
de~$(X,x)$). Par construction, deux applications
continues homotopes entre espaces topologiques
pointés induisent le même morphisme entre
les groupes fondamentaux ; autrement dit, on peut voir
$(X,x)\mapsto \pi_1(X,x)$ comme un foncteur de~$\pointhom$
vers~$\ens$. 

\trois{cov-bg} Nous allons décrire deux foncteurs mettant en jeu les
catégories un peu artificielles du~\ref{cat-decret}. 

\medskip

$\bullet$ Soit~$G$ un groupe. On définit comme suit un foncteur
covariant de~$\widetilde {BG}$
vers~$BG$ : il envoie n'importe quel élément~$g$ de~$G$
sur~$*$, et l'unique flèche entre deux éléments~$g$ et~$h$ de~$g$ sur~$gh^{-1}$.

$\bullet$ On définit comme suit un foncteur
covariant
de~$\VV_k$ dans la catégorie des espaces vectoriels
de dimension finie sur~$k$ : il envoie~$n$
sur~$k^n$, et il envoie une matrice~$M$ sur 
l'application linéaire de matrice~$M$
dans les bases canoniques. 

\medskip
\deux{plein-fidele}
{\bf Définition.} Soient~$\mathsf C$ et~$\mathsf D$ deux catégories, et soit~$F : \mathsf C\to \mathsf D$ un foncteur
covariant. On dit
que~$F$ est {\em fidèle} (resp. {\em plein}, resp. {\em pleinement fidèle})
si pour tout couple~$(X,Y)$ d'objets de~$\mathsf C$, 
l'application~$f\mapsto F(f)$ de~$\hom_{\mathsf C}(X,Y)$ vers
$\hom_{\mathsf D}(F(X), F(Y))$ est injective (resp. surjective, resp. bijective). 

\deux{sous-cat}
{\bf Définition}. Soit~$\mathsf C$ une catégorie. Une {\em sous-catégorie}
de~$\mathsf C$ est une catégorie~$\mathsf D$ telle que~${\rm Ob }\;D\subset {\rm Ob}\;\mathsf C$ 
et telle que~$\hom_{\mathsf D}(X,Y)\subset \hom_{\mathsf C}(X,Y)$ pour tout couple~$(X,Y)$
d'objets de~$\mathsf D$. 

Elle est dit {\em pleine}
si~$\hom_{\mathsf D}(X,Y)= \hom_{\mathsf C}(X,Y)$ pour tout couple~$(X,Y)$
d'objets de~$\mathsf D$. 

\deux{inclusion-sous-cat}
Si~$\mathsf D$ est une sous-catégorie d'une catégorie~$\mathsf C$, on dispose d'un foncteur 
covariant naturel d'inclusion de~$\mathsf C$ dans~$\mathsf D$. Il est fidèle, et pleinement
fidèle si~$\mathsf D$ est une sous-catégorie pleine de~$\mathsf C$. 

\section{Morphismes de foncteurs et équivalences de catégories}
\markboth{Le langage des catégories}{Morphismes de foncteurs}

\deux{def-morph-fonct} {\bf Définition.}
Soient~$\mathsf C$ et~$\mathsf D$
deux catégories et soient~$F$ et~$G : \mathsf C\to \mathsf D$ deux foncteurs covariants. Un
{\em morphisme}
(ou {\em transformation naturelle})
$\phi$ de~$F$ vers~$G$ est la donnée, pour tout
objet~$X$
de~$\mathsf C$, 
d'un morphisme
$$\phi(X) : F(X)\to G(X)$$ de la catégorie~$\mathsf D$, de sorte que pour
toute flèche~$f: X\to Y$ de~$\mathsf C$, le diagramme
$$\diagram F(X)\rto ^{\phi(X)}\dto_{F(f)}&G(X)\dto^{G(f)}\\
F(Y)\rto^{\phi(Y)}&G(Y)\enddiagram$$
commute. 

\trois{comment-vocab}
On résume parfois ces conditions en disant simplement qu'un morphisme
de~$F$ vers~$G$ est la donnée pour tout~$X$ d'un morphisme
de~$F(X)$ dans~$G(X)$ qui est {\em fonctoriel en~$X$.} 


\trois{id-compo} L'identité~${\rm Id}_F$ du foncteur~$F$ est le morphisme de~$F$ dans lui-même
induit par la collections des ${\rm Id}_{F(X)}$ où~$X$ parcourt~${\rm Ob}\;\mathsf C$. 
Les morphismes de foncteurs se composent de façon évidente. Attention toutefois : on 
ne peut pas dire que les foncteurs de~$\mathsf C$ vers~$\mathsf D$
constituent eux-mêmes une
catégorie, car rien n'indique {\em a priori}
que les morphismes entre deux tels foncteurs
(qui mettent en jeu une collection de données paramétrée par~${\rm Ob}\;\mathsf C$)
forment un ensemble.

\deux{ex-mor-fonct}
{\bf Exemples}.

\trois{mvvv}
Soit~$A$ un anneau commutatif unitaire. Pour tout~$A$-module~$M$, 
on dispose d'une application naturelle~$\iota(M) : M\to M^{\vee\vee}$
définie par la formule
$$m\mapsto (\phi\mapsto \phi(m)).$$ La collection des~$\iota(M)$ pour~$M$
variable définit un morphisme~$\iota$
du foncteur~${\rm Id}_{\amod}$ vers le foncteur~$M\mapsto M^{\vee\vee}$
de~$\amod$ dans elle-même (notons que ce dernier foncteur est bien covariant, en tant que
composée de deux foncteurs contravariants). 

\trois{hxhy}
Soit~$\mathsf C$ une catégorie, soient~$X$ et~$Y$ deux objets de~$\mathsf C$
et soit~$f$ un morphisme de~$X$ vers~$Y$. Pour tout~$Z\in {\rm Ob}\;\mathsf C$, 
on dispose de deux applications naturelles fonctorielles en~$Z$

$$\hom(Y,Z)\to \hom(X,Z), g\mapsto g\circ f\;\;{\rm et}\;
\hom(Z,X)\to \hom(Z,Y), g\mapsto f\circ g,$$
qui constituent deux morphismes de foncteurs~$f^*:h_Y\to h_X$ et~$f_*:h^X\to h^Y$. 

\medskip
On vérifie aisément qu'on a les égalités~$(f\circ g)_*=f_*\circ g_*$
et~$(f\circ g)^*=g^*\circ f^*$ à chaque fois qu'elles ont un sens. On
a aussi~$({\rm Id}_X)^*={\rm Id}_{h_X}$ et~$({\rm Id}_X)_*={\rm Id}_{h^X}$. 

\deux{iso-foncteur} {\bf Définition.}
Soient~$\mathsf C$ et~$\mathsf D$
deux catégories et soient~$F$ et~$G $
deux foncteurs
de~$\mathsf C$ vers~$\mathsf D$. On dit qu'un morphisme~$\phi : F\to G$
est un {\em isomorphisme} s'il existe~$\psi : G\to F$ 
tel que~$\psi \circ \phi={\rm Id}_F$ et~$\phi \circ \psi ={\rm Id}\;G$. 
Un tel~$\psi$ est dans ce cas unique, et est noté~$\phi^{-1}$. 

Le morphisme~$\phi$ est un isomorphisme si et seulement si~$\phi(X)$
est un isomorphisme pour tout~$X\in {\rm Ob}\;\mathsf C$, et l'on a alors 
$\phi^{-1}(X)=\phi(X)^{-1}$ pour tout
tel~$X$.

\deux{ex-bidual}
{\bf Exemple : la bidualité.} 
Soit~$A$ un anneau commutatif unitaire, et soit~$\mathsf C$ la sous-catégorie
pleine de~$\amod$ constituée des modules libres de rang fini 
({\em i.e.} qui possèdent une base finie, ou encore qui sont isomorphes à~$A^n$ pour
un certain~$n$). Si~$M\in {\rm Ob}\;\mathsf C$ alors~$M^{\vee \vee}\in {\rm Ob}\; C$, 
et le morphisme~$\iota$ du~\ref{mvvv}
ci-dessus induit un {\em isomorphisme}
entre le foncteur~${\rm Id}_{\mathsf C}$ et le foncteur~$M\mapsto M^{\vee \vee}$
de~$\mathsf C$ dans
$\mathsf C$. 

Cette assertion est simplement l'énoncé rigoureux traduisant le fait
q'un module libre de rang fini est {\em canoniquement}
isomorphe à son bidual. 

\deux{def-equiv}
{\bf Définitions.}
Soient~$\mathsf C$ et~$\mathsf D$ deux catégories, et soit~$F$ un foncteur de~$\mathsf C$ vers~$\mathsf D$. 
Un {\em quasi-inverse}
de~$F$ est un foncteur~$G : \mathsf D\to \mathsf C$ tel que~$F\circ G\simeq {\rm Id}_{\mathsf D}$ et~$G\circ F\simeq {\rm Id}_{\mathsf C}$. 

\medskip
On dit que~$F$ est une {\em équivalence de catégories}
s'il admet un quasi-inverse. 

\trois{tauto-equiv-cat}
Il résulte immédiatement des définitions que si~$G$ est un quasi-inverse
de~$F$ alors~$F$ est un quasi-inverse de~$G$,
et que la composée de deux équivalences de catégories est une équivalence de catégories. 

\medskip
Si~$G$ et~$H$ sont deux quasi-inverses de~$F$ ils sont isomorphes {\em via}
la composition d'isomorphismes 
$$G=G\circ {\rm Id}_{\mathsf D}\simeq G\circ (F\circ H)=(G\circ F)\circ H\simeq {\rm Id}_{\mathsf C}\circ H=H.$$

\trois{equiv-pf-es}
{\em Exercice.}
Montrez qu'un foncteur~$F: \mathsf C\to \mathsf D$ est une équivalence de catégories
si et seulement si~$F$
est pleinement fidèle et {\em essentiellement surjectif}, 
ce qui veut dire que pour tout objet~$Y$ de~$\mathsf D$ il existe un objet~$X$
de~$\mathsf C$ tel que~$F(X)$ soit {\em isomorphe}
à~$Y$. 

\deux{ex-equiv-cat}
{\bf Exemples d'équivalences de catégories.}

\trois{ex-iso-cat}
{\bf Les isomorphismes de catégories.} Un foncteur~$F:\mathsf C\to \mathsf D$
est appelé un {\em isomorphisme}
de catégories s'il possède un inverse, c'est-à-dire un foncteur~$G:\mathsf D\to \mathsf C$ tel
que~$F\circ G$ et~$G\circ F$ soient respectivement {\em égaux} (et pas seulement isomorphes)
à~${\rm Id}_{\mathsf D}$ et~${\rm Id}_{\mathsf C}$.Un tel~$G$ est unique s'il existe. 

Les isomorphismes de catégories sont des cas particuliers d'équivalences de catégories
que l'on rencontre très rarement en pratique. Citons deux exemples ; le premier est trivial,
le second
est nettement plus intéressant.  

\medskip
$\bullet$ Si~$\mathsf C$ est une catégorie alors~${\rm Id}_{\mathsf C}$ est un isomorphisme de catégories.

$\bullet$ Soit~$A$ un anneau
commutatif unitaire. Soit~$\mathsf C$ la catégorie 
des~$A[T]$-modules, et soit~$ D$
la catégorie dont les objets sont les couples~$(M,u)$ où~$M$
est un~$A$-module et~$u$ un endomorphisme de~$M$, et où
un morphisme de~$(M,u)$ vers~$(N,v)$ est une 
application~$A$-linéaire~$f: M\to N$ telle que~$v\circ f=f\circ u$. 
Soit~$F$ le foncteur qui envoie un~$A[T]$-module~$M$
sur le~$A$-module sous-jacent à~$M$ muni de~$m\mapsto Tm$ ; soit~$G$ le foncteur
qui envoie un couple~$(M,u)$ sur le~$A[T]$-module de groupe
additif sous-jacent~$(M,+)$ et de loi externe~$(P,m)\mapsto P(u)(m)$. 
Le foncteur~$F$ induit alors
un isomorphisme entre les catégories~$\mathsf C$
et~$\mathsf D$, et~$G$ est son inverse.

\trois{vn-vect-k}
Soit~$k$ un corps On a défini au~\ref{cov-bg}
un foncteur covariant
de la catégorie~$\VV_k$ (\ref{cat-decret}) 
vers celle des espaces vectoriels de dimension finie sur~$k$. 
{\em Ce foncteur est une équivalence de catégories}
-- cette assertion est essentiellement 
une reformulation conceptuelle (ou pédante) de l'existence de bases. 

En effet, 
on peut construire son quasi-inverse comme suit. On commence par choisir
\footnote{Cette opération requiert une forme redoutablement puissante d'axiome du choix,
puisque  les~$k$-espaces vectoriels de dimension finie
ne constituent pas un ensemble.}
pour tout~$k$-espace vectoriel~$E$ de dimension finie, une base~$b(E)$
de~$E$. On
envoie alors un~$k$-espace vectoriel~$E$
de dimension finie sur~$\dim E$, et une application linéaire~$f: E\to F$ entre deux
tels espaces sur~${\rm Mat}_{b(E), b(F)}f$. 

\trois{equipi1}{\bf La théorie des revêtements}
Nous allons maintenant indiquer (sans la moindre démonstration)
comment la théorie qui relie lacets tracés sur un espace topologique
et revêtements de ce dernier peut être reformulée en termes d'équivalence 
de catégories. 

\medskip
On se donne donc un espace topologique pointé~$(X,x)$. La catégorie
des~$\pi_1(X,x)$-ensembles est celle dont les objets sont les ensembles munis
d'une action à droite
de~$\pi_1(X,x)$, et les flèches les applications~$\pi_1(X,x)$-équivariantes. 

\medskip
On suppose que~$X$
est connexe, et qu'il est 
par ailleurs semi-localement simplement connexe\footnote
{Bourbaki propose de remplacer cette expression
peu engageante par {\em délaçable}.}. 

\medskip
Un {\em revêtement}
de~$X$ est un espace topologique~$Y$ 
muni d'une application continue~$f: Y\to X$ tel que tout point 
de~$X$ possède un voisinage ouvert~$U$ pour lequel 
il existe un ensemble discret~$E_U$ et un diagramme
commutatif
$$\diagram f^{-1}(U)\rto^\simeq\dto&E_U\times U\dlto\\
U&\enddiagram.$$ Un morphisme
entre deux revêtements~$(Y\to X)$ et~$(Z\to X)$
est une application continue~$Y\to Z$ telle que
$$\diagram Y\drto\rrto&&Z\dlto\\
&X&\enddiagram$$ commute. 

\medskip
Si~$Y$ est un revêtement de~$X$, sa fibre~$Y_x$ en~$x$ hérite
d'une action naturelle de~$\pi_1(X,x)$ à droite (un lacet~$\ell$
et un point~$y$ de~$Y_x$ étant donnés, on relève~$\ell$
en partant de~$y$ ; à l'arrivée, on ne retombe pas nécessairement sur ses pieds : 
on atteint un antécédent~$z$ de~$x$ qui en général diffère de~$y$, et l'on
pose~$y.\ell=z$). 

\medskip
Le lien entre lacets et revêtements peut alors
s'exprimer ainsi : le foncteur~$Y\mapsto Y_x$ établit une
équivalence entre la catégorie des revêtements de~$X$ et celle
des~$\pi_1(X,x)$-ensembles. 

{\em Description d'un quasi-inverse}. On choisit un revêtement 
universel~$\tilde X$ de~$X$ (c'est-à-dire un 
revêtement de~$X$ connexe, non vide et simplement connexe). 
Le foncteur qui envoie un~$\pi_1(X,x)$-ensemble~$E$
sur le produit contracté
$$\tilde X\times^{\pi_1(X,x)}E:=\tilde X\times E/((y.\ell, e)\sim (y, e.\ell))$$
est alors un quasi-inverse de~$Y\mapsto Y_x$. 

\trois{trans-gelfand}
{\bf La transformation de Gelfand}. Soit~$\mathsf{Comp}$
la sous-catégorie pleine de~$\top$ formée
des espaces topologiques compacts. 
Soit~$X$ un espace topologique compact. 
L'algèbre~$\sch C^0(X,\CC)$  est une~$\mathbb C$-algèbre 
commutative de Banach (pour la norme~$f\mapsto \sup_{x\in X} |f(x)|$) ; 
elle possède une involution~$f\mapsto \overline f$
qui prolonge
la conjugaison complexe,
et l'on a par construction~$||f||=\sqrt{|| f\overline f||}$ pour tout~$f$. 

\medskip
Soit~$\mathsf C$ la
catégorie définie comme suit. 

\medskip
$\bullet$ Un objet de~$\mathsf C$ est une~$\mathbb C$-algèbre
de Banach commutative~$A$
 munie d'une involution~$a\mapsto \overline a$
prolongeant la conjugaison complexe et telle que~$||a||=\sqrt{||a\overline a||}$
pour tout~$a\in A$. 

$\bullet$ Si~$A$ et~$B$ sont deux objets de~$\mathsf C$ un morphisme 
de~$A$ dans~$B$ est un morphisme d'algèbres~$f: A \to B$ tel que
~$||f(a)||\leq ||a||$ et~$f(\overline a)=\overline{f(a)}$ pour tout~$a\in A$. 

\medskip
Il est immédiat que~$\sch C : X\mapsto \sch C^0(X,\CC)$ définit un foncteur 
contravariant de~$\mathsf{Comp}$ vers~$\mathsf C$.

\medskip
On démontre en analyse fonctionnelle (théorie de Gelfand)
que le foncteur~$\sch C$ est anti-équivalence de catégories («anti» signifiant simplement qu'il
est contravariant). 

\medskip
Décrivons partiellement
un quasi-inverse~${\rm Sp}$ (le «spectre»)
de~$\sch C$. Si~$A\in C$, l'{\em ensemble} 
sous-jacent à l'espace topologique
compact~${\rm Sp}\;A$ est l'ensemble des idéaux maximaux de~$A$ (nous ne préciserons
pas ici la construction de la topologie de~${\rm Sp}\;A$). 

\medskip
{\em Remarque.} Un point d'un espace 
topologique ~$X$ définit bien un idéal maximal 
de l'algèbre~$\sch C^0(X,\CC)$, à savoir l'ensemble des fonctions
qui s'y annulent ; vous pouvez à titre d'exercice vérifier que 
si~$X$ est compact, on obtient ainsi une {\em bijection} 
entre~$X$ et l'ensemble des idéaux maximaux de~$\sch C^0(X,\CC)$. 

%
%
%
%
%

\section{Foncteurs représentables et lemme de Yoneda}
 
\markboth{Le langage des catégories}{Foncteurs représentables}

\deux{def-fonct-rep}
{\bf Définition.}
Soit~$\mathsf C$ une catégorie. Un foncteur covariant
de~$F : \mathsf C\to \ens$ est dit
{\em représentable}
s'il existe~$X\in {\rm Ob}\;\mathsf C$
et un isomorphisme~$h_X\simeq F$.

\deux{ex-fonct-rep}
{\bf Exemples.}
Dans ce qui suit, 
les anneaux et algèbres 
sont commutatifs   et unitaires,
et les morphismes d'anneaux ou algèbres 
sont unitaires. 

\trois{asurp} Soit~$A$ un anneau
et soit~$P$ une partie de~$A$. Le foncteur~$F$
qui envoie un anneau~$B$
sur l'ensemble des morphismes de~$A$
dans~$B$ s'annulant sur les éléments de~$P$ est représentable : en effet, 
si~$\pi$ désigne l'application quotient~$A\to A/\langle P\rangle$, 
alors~$f\mapsto f\circ \pi$ établit une bijection, fonctorielle en~$B$, 
entre~${\rm  Hom}_{\ann}(A/\langle P \rangle, B)$ et~$F(B)$. 

\trois{asurp-alg}
Nous allons donner une présentation un peu différente 
de l'exemple~\ref{asurp}
ci-dessus, dont nous conservons les notation~$A$ et~$P$. Soit~$\Phi$ le foncteur qui envoie une~$A$-algèbre~$(B, f\colon A \to B)$
sur un singleton~$\{*\}$ si~$P$ est contenu dans~${\rm Ker}\; f$ et sur~$\varnothing$ sinon ; il est représentable. 
En effet, soit~$(B, f\colon A\to B)$ une~$A$-algèbre. Le morphisme~$f$
admet une
factorisation par~$A\mapsto A/\langle P\rangle$ si et seulement si~$P\subset {\rm Ker}\; f$, et si c'est le cas cette factorisation est unique ; 
nous la noterons~$\bar f$. On peut reformuler cette assertion en disant que
l'ensemble des morphismes de~$A$-algèbres de~$A/\langle P\rangle$ dans~$A$
est le singleton~$\{\bar f\}$ si~$P\subset {\rm Ker}\;f$, et est vide sinon. 
On dispose dès lors d'une bijection, fonctorielle en~$B$, entre~$\hom_{\aalg}(A/\langle P\rangle,B)$ et~$\Phi(B)$ : c'est 
la bijection entre~$\{\bar f\}$ et~$\{*\}$ si~$P\subset {\rm Ker}\; F$, et~${\rm Id}_{\emptyset}$ dans le cas contraire.

\medskip

\trois{apuissancen}
Soit~$A$ un anneau et
soit~$n$
un entier. Le foncteur qui envoie 
une~$A$-algèbre~$B$ sur~$B^n$ est représentable. En effet,
l'application~$f\mapsto (f(T_1),\ldots, f(T_n))$ établit
une bijection fonctorielle en~$B$ 
entre~${\rm Hom}_{\aalg}(A[T_1,\ldots, T_n], B)$ et~$B^n$. 

\trois{apuissancenmodp}
Soit~$A$ un anneau,  soit~$n$
un entier et
soit~$(P_i)$ une famille
de polynômes appartenant à~$A[T_1,\ldots, T_n]$. 
Le foncteur~$G$~qui envoie 
une~$A$-algèbre~$B$ sur

$$\{(b_1,\ldots, b_n)\in B^n,P_i(b_1,\ldots, b_n)=0 \;\forall i\}$$
est représentable. En effet,
l'application~$f\mapsto (f(\overline{T_1}),\ldots, f(\overline{T_n}))$ établit
une bijection fonctorielle en~$B$ 
entre~${\rm Hom}_{\aalg}(A[T_1,\ldots, T_n]/(P_i)_i, B)$ et~$G(B)$.

\trois{groupe-fond-rep}
Le foncteur~$(X,x)\mapsto \pi_1(X,x)$
de~$\pointhom$ vers~$\ens$
est représentable : il s'identifie naturellement, par sa définition
même, à~${\rm Hom}_{\pointhom}((S_1,o),.)$, où~$S_1$ est
le cercle et~$o$ un point quelconque choisi sur~$S_1$. 

\deux{pre-intro-yoneda} Montrer
qu'un 
foncteur~$F$ est représentable revient à exhiber
un objet~$X$ {\em et un isomorphisme}~$h_X\simeq F$. 
Ce dernier point peut sembler délicat, étant donnée la définition
d'un morphisme de foncteurs comme une gigantesque collection de morphismes
(sujette à des conditions de compatibilité). On va voir un
peu plus bas (\ref{repres-yoneda})
qu'il n'en est rien, et qu'un
tel isomorphisme~$h_X\simeq F$
admet toujours une description concrète 
simple et maniable. 

\deux{intro-yoneda}
Soit~$\mathsf C$ une catégorie et soit~$X\in {\rm Ob}\;\mathsf C$. Soit~$F$
un foncteur covariant de~$\mathsf C$ dans~$\ens$, et soit~$\xi\in F(X)$. On vérifie
immédiatement que la collection d'applications
$$h_X(Y)=\hom(X,Y)\to F(Y), \; f\mapsto F(f)(\xi)$$ définit, lorsque~$Y$
parcourt~${\rm Ob}\;\mathsf C$, un morphisme de~$h_X$ dans~$F$ ; nous le noterons~$\phi_\xi$. 

\deux{lemme-yoneda}
{\bf Lemme de Yoneda.}
{\em Soit~$\psi$ un morphisme de~$h_X$ dans~$F$. Il existe un
unique~$\xi \in F(X)$ tel que~$\psi=\phi_\xi$, à savoir~$\psi(X)({\rm Id}_X)$.}

\medskip
{\em Démonstration.} Commençons par remarquer que l'énoncé a bien un sens : 
on a~${\rm Id}_X\in \hom(X,X)=h_X(X)$, d'où il résulte que~$\psi(X)({\rm Id}_X)$ est un
élément bien défini de~$F(X)$. 

\medskip
{\em Preuve de l'unicité.}
Soit~$\xi$ tel que~$\psi=\phi_\xi$. On a alors
$$\psi(X)({\rm Id}_X)=\phi_\xi(X)({\rm Id}_X)=F({\rm Id}_X)(\xi)={\rm Id}_{F(X)}(\xi)=\xi,$$
d'où l'assertion requise. 

\medskip
{\em Preuve de l'existence.}
On pose~$\xi=\psi(X)({\rm Id}_X)$, et l'on va démontrer
que~$\psi=\phi_\xi$. Soit~$Y\in {\rm Ob}\;\mathsf C$ et
soit~$f: X\to Y$ un élément de~$h_X(Y)$. 

On a par définition d'un morphisme de foncteurs un diagramme
commutatif 
$$\diagram
h_X(X)\dto_{\psi(X)}\rto^{h_X(f)}&h_X(Y)\dto^{\psi(Y)}\\
F(X)\rto^{F(f)}&F(Y)\enddiagram$$

Comme~$h_X(f)$ 
est la composition avec~$f$,
on
a
$h_X(f)({\rm Id}_X)=f$. L'image de~${\rm Id}_X$ par~$\psi(X)$
est par ailleurs égale à~$\xi$ par définition de ce dernier. Par commutativité du diagramme
il vient
alors~$\psi(Y)(f)=F(f)(\xi)=\phi_\xi(Y)(f)$, ce qu'il fallait démontrer.~$\Box$

\trois{mor-hx-f-ens}
{\bf Remarque.}
En particulier, les morphismes de~$h_X$ vers~$F$ constituent un ensemble
(en bijection naturelle avec~$F(X)$), ce qui n'était pas évident
{\em a priori}. 

\trois{mor-hx-hy}
Plaçons-nous
dans le cas particulier où~$F=h_Y$ pour un certain~$Y\in {\rm Ob}\;\mathsf C$. L'élément~$\xi$
du lemme appartient alors à~$h_Y(X)=\hom(Y,X)$, et le morphisme
$\psi=\phi_\xi$ est donné par la formule
$f\mapsto h_Y(f)(\xi)= f\circ \xi$. En conséquence, $\psi$ est le morphisme
$\xi^*: h_X\to h_Y$ induit par~$\xi$ ({\em cf.}~\ref{hxhy}).

\medskip
Il s'ensuit que {\em $\xi\mapsto \xi^*$ établit une bijection entre~$\hom(Y,X)$
et l'ensemble des morphismes de foncteurs de~$h_X$
vers~$h_Y$.}

\trois{mor-composition} Soit~$\xi \in \hom(Y,X)$. C'est un isomorphisme
si et seulement si~$\xi^*$ est un isomorphisme. En effet, si~$\xi$ est un isomorphisme
d'inverse~$\zeta$, on a

$$\zeta^*\circ \xi^*=(\xi \circ \zeta)^*={\rm Id}_Y^*={\rm Id}_{h_Y},$$ et de même
$\xi^*\circ \zeta^*={\rm Id}_{h_X}$. 

\medskip
Réciproquement, supposons que~$\xi^*$ soit un isomorphisme. Sa réciproque est alors
d'après~\ref{mor-hx-hy}
de la forme~$\zeta^*$ pour un certain~$\zeta\in \hom(X,Y)$. On a
$$(\xi \circ \zeta)^*=\zeta^*\circ \xi^*={\rm Id}_{h_X}={\rm Id}_X^*$$
et donc~$\xi \circ \zeta={\rm Id}_X$ d'après~{\em loc. cit.} De même, $\zeta\circ \xi={\rm Id}_Y$.

\deux{repres-yoneda}
Soit~$\mathsf C$ une catégorie et soit~$F$ un foncteur covariant de~$\mathsf C$ dans
$\ens$. Le foncteur~$F$ est représentable si et seulement si il existe
un objet~$X$ de~$\mathsf C$ et un isomorphisme~$h_X\simeq F$. En vertu du lemme de Yoneda, 
cela revient à demander qu'il existe un objet~$X$
de~$\mathsf C$ et un élément~$\xi$
de~$F(X)$ tel que~$\phi_\xi$ soit
un isomorphisme, c'est-à-dire tel que
$f\mapsto F(f)(\xi)$ établisse pour tout~$Y\in{\rm Ob}\;\mathsf C$ une bijection
$\hom(X,Y)\simeq F(Y)$. Nous dirons qu'un tel couple~$(X,\xi)$ est un
{\em représentant} de~$F$.

\deux{repres-canon}
{\bf Canonicité du représentant.}
Soient~$F$ et~$G$ deux foncteurs covariants d'une catégorie~$\mathsf C$ dans~$\ens$. 
On les suppose représentables ; soit~$(X,\xi)$
un représentant de~$F$ et
soit~$(Y,\eta)$ un représentant de~$G$.

\trois{repres-morph} Soit~$\rho$ un morphisme de foncteurs
de~$F$ dans~$G$. Il existe un unique morphisme~$\psi : h_X\to h_Y$
faisant commuter le diagramme
$$\diagram h_X\dto_{\phi_\xi}^\simeq\rto^\psi&
h_Y\dto_{\phi_\eta}^\simeq \\
F\rto^\rho&G\enddiagram\;\;\;,$$ à savoir
$\phi_\eta^{-1}\circ \rho\circ \phi_\xi$ ; comme~$\phi_\eta^{-1}$ et~$\phi_\xi$
sont des isomorphismes, $\psi$ est un isomorphisme si et seulement si~$\rho$ est un isomorphisme.

D'après~\ref{mor-hx-hy}, on peut reformuler ceci 
en disant qu'il existe un unique~$\lambda$
appartenant à~$\hom(Y,X)$
tel
que
le diagramme~$$\diagram h_X\dto_{\phi_\xi}^\simeq\rto^{\lambda^*}&
h_Y\dto_{\phi_\eta}^\simeq \\
F\rto^\rho&G\enddiagram\;\;\;$$
commute. En vertu de~\ref{mor-composition}, 
$\lambda$ est un isomorphisme si et seulement si~$\lambda^*$
est un isomorphisme, c'est-à-dire par
ce qui précède
si et seulement si~$\rho$
est un isomorphisme. 

\medskip
En appliquant
la commutativité du diagramme
à l'élément~${\rm Id}_X$
de~$ \hom(X,X)=h_X(X)$, et en
remarquant que~$\lambda^*(X)({\rm Id}_X)=\lambda$, 
on obtient l'égalité
$$\phi_\eta(X)(\lambda)=\rho(X)\left(\phi_\xi(X) ({\rm Id}_X)\right),$$
qui {\em caractérise}~$\lambda$ 
puisque~$\phi_\eta(X)$ est une bijection ; par définition
de~$\phi_\eta$ et~$\phi_\xi$, elle se récrit
$$G(\lambda)(\eta)=\rho(X)(\xi)$$
(notons que c'est une égalité entre éléments de~$G(X)$). 

%

\trois{deduction-canon}
En appliquant
ce qui précède
lorsque~$G=F$
et~$\rho={\rm Id}_F$, on obtient l'assertion suivante : 
{\em il existe
un unique
morphisme~$\lambda$ de~$Y$ vers~$X$ tel que~$F(\lambda)(\eta)=\xi$ ; 
ce morphisme est un isomorphisme, et peut également être caractérisé par 
la commutativité du diagramme
$$\diagram h_X\drto_{\phi_\xi}^\simeq\rrto^{\lambda^*}&&h_Y
\dlto^{\phi_\eta}_\simeq\\
&F&\enddiagram\;\;.$$}

\medskip
Il existe en particulier un unique isomorphisme~$\lambda : Y\simeq X$
tel que~$F(\lambda)(\eta)=\xi$ ; les couples~$(X,\xi)$ et~$(Y,\eta)$ sont donc
en un sens {\em canoniquement isomorphes}. Pour cette raison, nous nous permettrons
souvent de parler par abus
{\em du}
représentant~$(X,\xi)$ d'un foncteur représentable~$F$ ; et nous dirons que~$\xi$
est {\em l'objet universel relatif à~$F$}. 

\medskip
Il arrivera qu'on omette de mentionner~$\xi$ et qu'on 
se contente de dire que~$X$ est le représentant de~$F$ ; mais attention, 
c'est un peu imprudent, car pour un même~$X$, plusieurs~$\xi$ peuvent convenir.

\trois{reformulation-yoneda}
On peut alors reformuler
le lemme de Yoneda, ou plus précisément
la déclinaison qui en est faite en~\ref{repres-morph}
en disant que {\em se donner un morphisme entre foncteurs représentables,
c'est se donner un morphisme entre leurs représentants.} 

\deux{fonct-rep-revisit}
Nous allons revisiter
certains des exemples
vus en~\ref{ex-fonct-rep}
en donnant à chaque fois le représentant du foncteur considéré,
et en insistant sur son objet universel.

\trois{asurp-rev} Soit~$A$ un anneau
et soit~$P$ une partie de~$A$. Le foncteur
qui envoie un anneau~$B$
sur l'ensemble des morphismes de~$A$
dans~$B$ s'annulant sur les éléments de~$P$ est représentable. 
Si~$\pi$ désigne l'application quotient~$A\to A/\langle P\rangle$, 
alors~$(A/(P), \pi)$ est son représentant. 

\medskip
Le morphisme~$\pi$ est le «morphisme universel s'annulant sur~$P$». 

\trois{apuissancen-rev}
Soit~$A$ un anneau et
soit~$n$
un entier. Le foncteur qui envoie 
une~$A$-algèbre~$B$ sur~$B^n$ est représentable. 
Le couple~$(A[T_1,\ldots, T_n], (T_1,\ldots, T_n))$ est son représentant. 

\medskip
Le~$n$-uplet~$(T_1,\ldots, T_n)$ est le «$n$-uplet universel d'éléments d'une~$A$-algèbre». 

\trois{apuissancenmodp-rev}
Soit~$A$ un anneau,  soit~$n$
un entier et
soit~$(P_i)$ une famille
de polynômes appartenant à~$A[T_1,\ldots, T_n]$. 
Le foncteur
qui envoie 
une~$A$-algèbre~$B$ sur
$$\{(b_1,\ldots, b_n)\in B^n,P_i(b_1,\ldots, b_n)=0 \;\forall i\}$$
est représentable. 
Le couple~$(A[T_1,\ldots, T_n]/(P_i)_i,(\overline{T_1},\ldots, \overline{T_n}))$
est son représentant. 

\medskip
Le~$n$-uplet~$(\overline{T_1},\ldots, \overline{T_n})$ est le «$n$-uplet universel d'éléments d'une~$A$-algèbre
en lequel les~$P_i$ s'annulent» . 

\trois{groupe-fond-rep-rev} Le foncteur~$(X,x)\mapsto \pi_1(X,x)$
de~$\pointhom$ vers~$\ens$
est représentable ;  le couple
$((S_1,o),{\rm Id}_{(S_1, o)})$ est son représentant. 

\medskip
L'application identité de~$(S_1,o)$
est le «lacet universel à homotopie près».

\deux{mor-bn-bm} {\bf Exemple d'application du lemme de Yoneda}. 
Soient~$m$ et~$n$ deux entiers. 
Soit~$A$ un anneau commutatif unitaire, et soient~$F$ et~$G$ les foncteurs~$B\mapsto B^n$ et~$B\mapsto B^m$
sur la catégorie des~$A$-algèbres. Ils sont
représentables, de représentants respectifs~$(A[T_1,\ldots,T_n],(T_1,\ldots, T_n))$
et~$(A[S_1,\ldots,S_m],(S_1,\ldots, S_m))$. Soit~$\phi$
un morphisme de~$F$ vers~$G$. Par le lemme de Yoneda, il provient d'un unique
morphisme de~$A[S_1,\ldots,S_m]$ vers~$A[T_1,\ldots, T_n]$. Un tel morphisme
est lui-même donné par~$m$ polynômes~$P_1,\ldots, P_m$ appartenant
à~$A[T_1,\ldots, T_n]$ (les images des~$S_i$). On vérifie immédiatement 
que~$\phi$ est alors décrit par les formules
$$(b_1,\ldots, b_n)\mapsto (P_1(b_1,\ldots,b_n),\ldots, P_m(b_1,\ldots, b_n)).$$

Ce résultat est un sens assez intuitif : la seule façon d'associer de façon naturelle
à tout~$n$-uplet d'éléments d'une~$A$-algèbre~$B$ un~$m$-uplet d'éléments de~$B$
consiste à utiliser une formule polynomiale à coefficients dans~$A$ ; on pouvait s'y attendre,
puisque les seules opérations que l'on sache effectuer dans une~$A$-algèbre générale 
sont l'addition, la multiplication interne, et la multiplication par les éléments de~$A$. 

\deux{ob-iniit-fin}
{\bf Objet initial, objet final}. Soit~$\mathsf C$ une catégorie. Comme
il existe une et une seule application d'un singleton dans un autre, la
flèche~$\mathsf \mathsf C\to \mathsf {Ens}$ qui envoie un objet~$X$ sur un singleton~$\{*\}$
peut être vue d'une unique manière comme un foncteur covariant~$F$, et d'une unique manière comme 
un foncteur contravariant~$G$. 

Si~$F$ est représentable, on appelle {\em objet initial} de~$\mathsf C$
tout représentant de~$F$. Un objet~$X$ de~$\mathsf C$ est initial si et seulement si~${\rm Hom}(X,Y)$ est un 
singleton pour tout objet~$Y$ de~$\mathsf C$. 

Si~$G$ est représentable, on appelle {\em objet final} de~$\mathsf C$
tout représentant de~$F$. Un objet~$X$ de~$\mathsf C$ est final si et seulement si~${\rm Hom}(Y,X)$ est un 
singleton pour tout objet~$Y$ de~$\mathsf C$. 

\deux{ex-init-fin}
{\bf Exemples.} 

\trois{ens-in-fin}
Dans la catégorie des ensembles, $\varnothing$ est initial, et~$\{*\}$ est final. 

\trois{gp-in-fin}
Dans la catégorie des groupes, $\{e\}$ est initial et final. 

\trois{amod-in-fin}
Dans la catégorie des modules sur un anneau
commutatif unitaire~$A$ donné,~$\{0\}$ est initial et final. 

\trois{ann-in-fin}
Dans la catégorie des anneaux,~$\mathbb Z$ est initial, et~$\{0\}$ est final. 

\trois{corps-in-fin}
La catégorie des corps n'a ni objet initial, ni objet final. Celle des corps de
caractéristique nulle a~$\mathbb Q$ comme objet initial, et n'a pas d'objet final ; celle
des corps de caractéristique~$p>0$ a~${\mathbb F}_p$ comme objet initial, et n'a pas
d'objet final. 

\section{Produits fibrés et sommes amalgamées}
\markboth{Le langage des catégories}{Produits cartésiens et sommes amalgamées}

\subsection*{Produits cartésiens et produits fibrés}

\deux{prod-cartesien}
{\bf Le produit cartésien}. Le produit cartésien
 ensembliste~$X\times Y$ de deux
ensembles~$X$ et~$Y$ est par définition l'ensemble des couples~$(x,y)$
où~$x\in X$ et~$y\in Y$. Soit~$\mathsf C$ une catégorie et soient~$X$ et~$Y$ 
deux objets de~$\mathsf C$. Si
le foncteur contravariant~
$$h^X\times h^Y:=T\mapsto h^X(T)\times h^Y(T)=\hom(T,X)\times \hom(T,Y)$$ est représentable,
son représentant est en général noté~$X\times Y$
et est appelé
le {\em produit cartésien}
de~$X$ et~$Y$. 

\medskip
On prendra garde que ce représentant est 
en réalité 
 constitué de~$X\times Y$ {\em et d'un couple~$(p,q)$
 de morphismes où~$p$ 
 va de~$X\times Y$ vers~$X$ 
et~$q$ de~$X\times Y$ vers~$Y$.}
Ces morphismes n'ont pas  de notation standard ; on les appelle
 les première et seconde projections. 
 
 \medskip
 La définition du
 couple~$(X\times Y), (p,q))$ 
 comme représentant de~$h^X\times h^Y$ signifie que pour tout objet~$T$ de~$\mathsf C$, 
 tout morphisme~$\lambda: T\to X$
 et tout morphisme~$\mu: T\to Y$, il existe un et un seul morphisme~$\pi : T\to X\times Y$
 tel que~$\lambda =p\circ \pi$ et~$\mu=q\circ \pi$. Ou encore que pour tout objet~$T$ de~$\mathsf C$,
 l'application~$\pi \mapsto (p\circ \pi, q\circ \pi)$ établit une bijection entre
 $\hom (T,X\times Y)$ et~$\hom(T,X)\times \hom(T,Y)$. En termes plus imagés : 
 se donner un morphisme d'un objet~$T$ vers~$X\times Y$, c'est se donner un morphisme de~$T$ 
 vers~$X$ et un morphisme de~$T$ vers~$Y$. 
 
 \deux{ex-prod-cart}
 {\bf Quelques exemples de catégories dans lesquelles les produits cartésiens
 existent}. 
 
 \trois{prod-cart-end}
 Dans la catégorie des ensembles, le produit cartésien de deux ensembles
 est leur produit cartésien ensembliste. 
 
\trois {prod-cart-top} Dans la catégorie des espaces topologiques, le produit cartésien de deux espaces
 topologiques est leur produit cartésien ensembliste {\em muni de la topologie produit.} 
 
 \trois{prod-carta-alg}
 Dans la catégorie des groupes (resp. des anneaux, resp. des modules
 sur un anneau), le produit 
 cartésien de deux objets est leur produit cartésien 
 ensembliste, la ou les opérations étant définies coordonnée
 par coordonnée. 
 
 \deux{prod-fib}
 {\bf Le produit fibré}. 
 
 \trois{prod-fib-en}
 {\bf Le cas des ensembles.}
 Soient~$X,Y$ et~$S$
 trois ensembles, et soit $f: X\to S$ et~$g: Y\to S$ deux applications. 
 Le {\em produit fibré ensembliste}
 $X\times_{f,g}Y$ (ou
 plus simplement~$X\times_S Y$ s'il n'y a pas d'ambiguïté sur
 la définition des flèches~$f$ et~$g$) 
 est par définition l'ensemble des couples~$(x,y)$
où~$x\in X$, où~$y\in Y$ et où~$f(x)=g(y)$. Il est muni d'une application 
naturelle~$h$ de but~$S$, qui envoie~$(x,y)$
sur~$f(x)=g(y)$, et l'on a pour tout~$s$
appartenant à~$S$ l'égalité~$h^{-1}(s)=f^{-1}(s)\times g^{-1}(s)$ : 
ainsi,~$X\times_S Y$ est un «produit cartésien fibre à fibre», d'où l'expression
produit fibré.

\trois{prod-fib-cat}
{\bf Le cas général.}
Soit~$\mathsf C$ est une catégorie, soient~$X,Y$ 
et~$S$ 
trois objets de~$\mathsf C$, et soient~$f: X\to S$ et~$g: Y\to S$
deux morphismes. Si
le foncteur contravariant
$$h^X\times_{h^S}h^Y:=T\mapsto {\rm Hom}(T,X)\times_{{\rm Hom}(T,S)}  {\rm Hom}(T,Y)$$ est représentable, 
on note en général~$X\times_S Y$ ou~$X\times_{f,g}Y$ son représentant,
que l'on appelle {\em produit fibré de~$X$ et~$Y$
au-dessus de~$S$ (ou le long de~$(f,g)$).}
On prendra garde que ce représentant est en réalité 
 constitué de~$X\times_SY$ {\em et d'un couple~$(p,q)$de morphismes où~$p$ 
 va de~$X\times_S Y$ vers~$X$ 
et~$q$ de~$X\times_SY$ vers~$Y$, et où~$f\circ p=g\circ q$}. 
Ces morphismes n'ont pas  de notation standard ; on les appelle
 les première et seconde projections. 
 
 La définition du couple~$((X\times_SY), (p,q))$
 comme
 représentant du
 foncteur~$h^X\times_{h^S}h^Y$
 signifie 
 que pour tout objet~$T$ de~$\mathsf C$, tout morphisme
 $\lambda: T\to X$
 et tout morphisme~$\mu : T\to Y$ tel que~
$ f\circ \lambda=g\circ \mu$, 
 il existe un et un seul morphisme~$\pi : T\to X\times_S Y$
 tel que~$\lambda=p\circ \pi$ et~$\mu=q\circ \pi$. Ou encore que pour tout objet~$T$ de~$\mathsf C$,
 l'application~$\pi \mapsto (p\circ \pi, q\circ \pi)$ établit une bijection entre
 $\hom (T,X\times_SY)$ et~$\hom(T,X)\times_{\hom (T,S)} \hom(T,Y)$. En termes plus imagés : 
 se donner un morphisme d'un objet~$T$ vers~$X\times_S Y$, c'est se donner un morphisme de~$T$ 
 vers~$X$ et un morphisme de~$T$ vers~$Y$ dont les composés respectifs avec~$f$ et~$g$
 coïncident.

 \trois{diag-prod-fib}
 On peut schématiser ce qui précède par le diagramme commutatif suivant, 
 dans lequel les flèches en dur sont données, et où la flèche en pointillés
 est celle fournie
 par la propriété universelle.
  
 $$\xymatrix{T\ar@{.>}|-{\pi}[rd]\ar@/^/^\mu[rrd]\ar@/_/_\lambda[rdd]&&
 \\& X\times_S Y \ar[d]_p\ar[r]^q&Y\ar[d]^g\\
 &X\ar[r]_f&S
 }$$

\trois{prod-fib-verss}
Notons que~
le produit fibré~$X\times_SY$ est muni d'un morphisme naturel vers~$S$, 
 à savoir~$f\circ p=g\circ q$. 
  
 \deux{ex-prod-fib}
 {\bf Quelques exemples de catégories dans lesquelles les produits fibrés
 existent}. 
 
 \trois{prod-fib-ens}
 Dans la catégorie des ensembles, le produit fibré de deux ensembles
au-dessus d'un troisième  est leur produit fibré ensembliste. 
 
 \trois{prod-fib-top}
 Dans la catégorie des espaces topologiques, le produit fibré de deux espaces
 topologiques au-dessus d'un troisième
 est leur produit fibré ensembliste {\em muni de la topologie induite
 par la topologie produit.} 
 
 \trois{prod-fib-alg}
 Dans la catégorie des groupes (resp. des anneaux, resp. des
 modules sur un anneau~$A$), le produit 
 fibré de deux objets au-dessus d'un troisième
 est leur produit
 fibré ensembliste, la ou les opérations étant définies coordonnée
 par coordonnée.

 \subsection*{Quelques tautologies} 
 
 \medskip
 Les démonstrations
 des assertions qui suivent sont laissées en exercice, à l'exception d'une
 d'entre elles qui est rédigée à titre d'exemple.

 \deux{prod-fib-fonct}Si~$\mathsf C$ est une catégorie dans laquelle les produits fibrés existent
 et si~$Y$ est un objet de~$\mathsf C$, alors~$X\mapsto X\times Y$ est de manière naturelle
 un foncteur covariant en~$X$
 de~$\mathsf C$ dans elle-même.
 
 \deux{cart-ob-fin} Si~$\mathsf C$ est une catégorie admettant un objet final~$S$,
 et si~$X$ et~$Y$ sont deux objets de~$\mathsf C$ tels que~$X\times Y$ existe,
 alors~$X\times_S Y$ existe\footnote{En toute rigueur, il
 faudrait préciser quels 
 sont les morphismes~$X\to S$
 et~$Y\to S$ que l'on considère ; mais comme~$S$
 est final, il n'y a pas le choix.}
 et s'identifie à~$X\times Y$ : les produits cartésiens
 sont donc des produits fibrés, pour peu qu'il existe un objet final.

\deux{prod-c-sur-s}
Soit~$\mathsf C$ une catégorie et 
 soit~$S$ un objet de~$\mathsf C$. Rappelons que
 la catégorie~$\mathsf C/S$ a été définie au~\ref{cat-semiclass}.

  \trois{fib-cart-ca}
  Soient~$X\to S$ et~$Y\to S$ deux objets de~$\mathsf C/S$
  (\ref{cat-semiclass}). 
 Alors
  le produit cartésien de~$(X\to S)$
  et~$(Y\to S)$ existe
  dans~$\mathsf C/S$ si et seulement si~$X\times_SY$ existe
  dans~$\mathsf C$, et si c'est le
  cas alors
  $$ (X\to S)\times(Y\to S)=X\times_SY$$ (ce
  dernier est vu comme objet de~$ C/S$ {\em via} son morphisme
  naturel vers~$S$, {\em cf.} la fin de~\ref{prod-fib}).
  
  \medskip
  On en déduit, à l'aide de~\ref{prod-fib-fonct}
  que si les produits fibrés existent dans~$\mathsf C$ alors pour tout
  objet~$(Y\to S)$ de~$\mathsf C/S$ la formule
  $(X\to S)\mapsto X\times_SY$ définit un foncteur covariant en~$(X\to S)$
  de~$\mathsf C/S$ dans elle-même.

  \trois{rest-cs} Soit~$X$ un objet de~$\mathsf C$. Supposons que~$X\times S$ existe,
  et considérons-le comme un objet de~$\mathsf C/S$ {\em via} la seconde projection. 
  La restriction
  à~$\mathsf C/S$ du foncteur~$h^X$ est alors égale à~${\rm Hom}_{\mathsf C/S}(.,X\times S)$.
  
\deux{xyzt} Soit~$\mathsf C$ une catégorie 
  dans laquelle les produits fibrés existent, et soient~$X,Y$, $Z$ et~$T$
  des objets de~$\mathsf C$. Supposons donnés un morphisme~$f:X\to Z$, un 
  morphisme~$g:Y\to Z$, et un morphisme~$h:T\to Y$. On a alors un isomorphisme
  naturel~$$(X\times_ZY)\times_Y T\simeq X\times_Z T,$$ où~$T$ est 
  vu comme muni du morphisme~$\lambda : T\to Z$ composé de~$T\to Y$ et~$Y\to Z$. 

\medskip
Nous allons justifier cette assertion. L'idée de la preuve consiste à se ramener
grâce au lemme de Yoneda à l'assertion ensembliste
correspondante\footnote{On trouve quelque part dans
les SGA l'expression «procédé ensembliste breveté»
à propos de ce type de méthode.}.

\trois{red-cas-ens}
{\em Réduction au cas ensembliste.} 
Grâce au lemme de Yoneda, il suffit, pour exhiber un isomorphisme entre
deux objets, d'en exhiber un entre les foncteurs qu'ils représentent. On va donc
chercher à montrer que
$$(h^X\times_{h^Z}h^Y)\times_{h^Y}h^T\simeq h^X\times_{h^Z}h^T,$$
c'est-à-dire encore qu'il existe, 
pour tout objet~$S$ de~$\mathsf C$,
une bijection fonctorielle en~$S$ 

$$(\hom (S,X)\times_{\hom (S,Z)}\hom(S,Y))\times_{\hom(S,Y)}\hom(S,T)$$
$$\simeq
\hom(S,X)\times_{\hom(S,Z)}\hom (S,T).$$ Il suffit donc de montrer
l'existence d'un isomorphisme comme dans l'énoncé lorsque~$\mathsf C$ est la catégorie
des ensembles, en s'assurant en plus qu'il est bien fonctoriel en les données. 

\trois{le-cas-ens}
{\em Le cas ensembliste.}

On a
$$ (X\times_ZY)\times_Y T=\{((x,y),t), f(x)=g(y)\;{\rm et}\;h(t)=y\}$$
et~$$X\times_ZT=\{(x,t), f(x)=\lambda(t)=g(h(t))\}.$$ 

La bijection
cherchée est alors $((x,y),t)\mapsto (x,t)$, 
sa réciproque
étant égale à~$(x,t)\mapsto ((x,h(t)),t)$ ; la fonctorialité
en les données est évidente. 

\subsection*{Sommes disjointes et sommes amalgamées}

Il s'agit des notions duales de celles
  de produit cartésien et de produit fibré.

\deux{somme-dis}
{\bf La somme disjointe}.  La somme
 disjointe ensembliste~$X\coprod Y$ de deux
ensembles~$X$ et~$Y$ est par définition la réunion
d'une copie de~$X$ et d'une copie de~$Y$, rendues disjointes
(pour
le faire proprement, on peut par exemple la définir comme le sous-ensemble
de~$X\cup Y)\times \{0,1\}$ formé
des couples de la forme~$(x,0)$ avec~$x\in X$ ou~$(y,1)$
avec~$y\in Y$). 

\medskip
Si~$\mathsf C$ est une catégorie et si~$X$ et~$Y$ sont
deux objets de~$\mathsf C$, et si
le foncteur covariant~$$h_X\times h_Y $$ est représentable,
on note en général~$X\coprod Y$ son représentant,
que l'on appelle la {\em somme disjointe}
de~$X$ et~$Y$. Attention : ce représentant est en réalité 
 constitué de~$X\coprod Y$ et d'un couple~$(i,j)$
 de morphismes où~$i$ 
 va de~$X$ vers~$X\coprod Y$
et~$j$ de~$Y$ vers~$X\coprod Y$.
Ces morphismes n'ont pas  de notation standard.  

La définition du
couple~$(X\coprod Y, (i,j))$
comme représentant de~$h_X\times h_Y$
signifie 
que pour tout objet~$T$ de~$\mathsf C$, tout morphisme~$\lambda: X\to T$
 et tout morphisme~$\mu: Y\to T$, il existe un et un seul morphisme~$\sigma : X\coprod Y \to T$
 tel que~$\lambda =\sigma\circ i$ et~$\mu=\sigma \circ j$. Ou encore que pour tout objet~$T$ de~$\mathsf C$,
 l'application~$\pi \mapsto (\pi\circ i, \pi\circ j)$ établit une bijection entre
 $\hom (X\coprod Y,T)$ et~$\hom(X,T)\times \hom(Y,T)$. En termes plus imagés : 
 se donner un morphisme de~$X\coprod Y$ vers un objet~$T$,
c'est se donner un morphisme de~$X$ 
 vers~$T$ et un morphisme de~$Y$ vers~$T$.

 \deux{ex-som-dis}
 {\bf Quelques exemples de catégories dans lesquelles les sommes
 disjointes existent}. 
 
 \trois{som-dis-ens}
 Dans la catégorie des ensembles, la somme disjointe
 de deux ensembles
 est leur somme disjointe ensembliste. 
 
\trois{som-dis-top}
Dans la catégorie des espaces topologiques, la somme disjointes
de deux espaces topologiques~$X$ et~$Y$ est leur 
somme disjointe ensembliste~$X\coprod Y$ munie de la topologie définie comme suit : 
sa restriction à~$X$ (resp.~$Y$) identifié à une partie de~$X\coprod Y$ est la topologie de~$X$
(resp. de~$Y)$ ; les parties~$X$ et~$Y$ de~$X\coprod Y$ sont toutes deux ouvertes. 
 
\trois{som-dis-mod}
Dans la catégorie des modules sur un anneau~$A$,
 la somme disjointe de deux modules est leur somme directe externe. 
 {\em L'ensemble sous-jacent
 n'est donc pas leur somme disjointe ensembliste} (on a identifié les zéros des deux modules, 
 et rajouté les sommes d'éléments). 
  
 \trois{som-dis-gp}
 Dans la catégorie des groupes, le même type de phénomène se reproduit. Si~$G$ 
 et~$H$ sont deux groupes, leur somme disjointe dans la catégorie des groupes
 est le {\em produit libre~$G*H$ de~$G$ et~$H$}, défini comme suit. Un élément
 de~$G*H$ est une suite finie~$x_1\ldots x_n$
 où chaque~$x_i$ est un élément non trivial
 de~$G\coprod H$, et telle qu'éléments de~$G$ et éléments de~$H$ alternent
 (une telle suite est souvent appelée un {\em mot}, dont les~$x_i$
 sont les {\em lettres}.) 
 On fait le
 produit de deux mots
 en les concaténant, puis en contractant
 le résultat
 obtenu selon l'algorithme suivant. 
  
\begin{itemize}
\item[$\bullet$] Si le mot contient une séquence
 de la forme~$xy$ avec~$x$ et~$y$
 appartenant tous deux à~$G$ (resp.~$H$), et si l'on note~$z$
 le produit~$xy$ calculé dans~$G$ (resp.~$H$) alors : 
 
 \begin{itemize}
 \item[$\diamond$] si~$z\neq e$ on remplace la séquence~$xy$
 par la lettre~$z$ ; 
 
 \item[$\diamond$] si~$z=e$ on supprime la séquence~$xy$. 
 
 \end{itemize} 
 \item[$\bullet$] On recommence l'opération jusqu'à ce que le
 mot ne contienne plus
 de séquence
 de la forme~$xy$ avec~$x$ et~$y$
 appartenant tous deux à~$G$
 ou tous deux à~$H$.  
 
 \end{itemize}
 
 \medskip
 On prendra garde que les~$x_i$
 sont des éléments de la réunion
 {\em disjointe}
 de~$G$ et~$H$ : même si~$G=H$ on considère
 dans cette construction~$G$
 et~$H$ comme deux groupes {\em distincts}, ce qui veut dire qu'on ne simplifie {\em jamais} un produit~$xy$
 avec~$x\in G$ et~$y\in H$ ou le contraire. Notons que l'élément
 neutre de~$G*H$ est le mot
 vide. 
 
 \medskip
 Par exemple, si~$G=H=\mathbb Z$, et si l'on écrit~$G=a^{\mathbb Z}$
 et~$H=b^{\mathbb Z}$ pour distinguer les deux groupes, 
 on voit que~$G*H$ est constitué des mots en les lettres~$a, a^{-1}, b$
 et~$b^{-1}$ : c'est le groupe
 libre sur deux générateurs. 
 
\deux{som-amalg}
 {\bf La somme amalgamée}.

 \trois{som-amalg-ens}
 {\bf Le cas ensembliste.} 
 Soient~$X,Y$ et~$S$
 trois ensembles, et soit~$f: S\to X$ et~$g: S\to Y$ deux applications. 
 La {\em somme amalgamée ensembliste}~$X\coprod_{f,g}Y$ (ou
 plus simplement~$X\coprod_SY$ s'il n'y a pas d'ambiguïté sur
 la définition des flèches~$f$ et~$g$) 
 est par définition le quotient
de~$X\coprod Y$
par la relation d'équivalence engendrée
par les relations~$f(s)\sim g(s)$ pour~$s$ parcourant~$S$.
On dispose de deux applications
naturelles~$X\to X\coprod_S Y$ et~$Y\to X\coprod_S Y$.  

\trois{som-amalg-cat}
{\bf Le cas général.} 
Si~$\mathsf C$ est une catégorie, si~$X,Y$ 
et~$S$ sont
trois objets de~$\mathsf C$, si~$f: S\to X$ et~$g: S\to Y$
sont deux morphismes,
et si
le foncteur covariant~$h_X\times_{h_S}h_Y$ est représentable,
on note en général~$X\coprod_S Y$ ou~$X\coprod_{f,g}Y$
son représentant, que l'on appelle
la {\em somme amalgamée de~$X$ et~$Y$ le long de~$S$
(ou le long de~$(f,g)$)}. 
Attention : ce représentant est en réalité 
 constitué de~$X\coprod_S Y$ et d'un couple~$(i,j)$ de morphismes
 où~$i$ 
 va de~$X$ vers~$X\coprod_S Y$
et~$j$ de~$Y$ vers~$X\coprod_S Y$, et où~$i\circ f=j\circ g$. 
Ces morphismes n'ont pas  de notation standard. 
 
 La définition du
 couple~$(X\coprod_S Y, (i,j))$
 comme représentant du foncteur
 $h_X\times_{h_S}h_Y$ signifie 
que
pour tout objet~$T$ de~$\mathsf C$, tout morphisme~$\lambda: X\to T$
 et tout morphisme~$\mu : Y\to T$ tel que~
$ \lambda\circ i=\mu\circ j$, 
 il existe un et un seul morphisme~$\sigma : X\coprod_S Y\to T$
 tel que~$\lambda=\sigma\circ i$ et~$\mu=\sigma\circ j$. Ou encore que pour tout objet~$T$ de~$\mathsf C$,
 l'application~$\pi \mapsto (\pi\circ i, \pi \circ j)$ établit une bijection entre
 $\hom (X\coprod_SY,T)$ et~$\hom(X,T)\times_{\hom (S,T)} \hom(Y,T)$. En termes plus imagés : 
 se donner un morphisme de~$X\coprod_S Y$
 vers un objet~$T$,
 c'est se donner un morphisme de~$X$ 
 vers~$T$ et un morphisme de~$Y$ vers~$T$ dont les composés respectifs avec~$i$ et~$j$
 coïncident. 

\trois{diag-som-amalg}
 On peut schématiser ce qui précède par le diagramme commutatif suivant, 
 dans lequel les flèches en dur sont données, et où la flèche en pointillés
 est celle fournie
 par la propriété universelle. 
  
 $$\xymatrix{T&&
 \\& X\coprod_S Y\ar@{.>}|-{\sigma}[ul]&Y\ar_j[l]\ar@/_/_\mu[ull]\\
 &X\ar^i[u]\ar@/^/^\lambda[uul]&S\ar^f[l]\ar_g[u]
 }$$

\trois{s-vers-somamalg}
Notons que~$S$ est muni d'un morphisme
 naturel vers~$X\coprod_SY$,
  à savoir~$i\circ f=j\circ g$.

 \deux{ex-som-amalg}
 {\bf Quelques exemples de catégories dans lesquelles les sommes amalgamées
 existent}. 
 
 \trois{som-amalg-ens}
 Dans la catégorie des ensembles, la somme amalgamée de deux ensembles
le long d'un troisième  est leur somme amalgamée ensembliste. 
 
 \trois{som-amalg-top}
 Dans la catégorie des espaces topologiques, le somme amalgamée
 de deux espaces
 topologiques~$X$ et~$Y$
 le long d'un troisième espace~$S$ 
 est leur somme amalgamée ensembliste, munie de la topologie quotient
 (la somme disjointe~$X\coprod Y$ étant munie quant à elle de sa topologie décrite plus
 haut). 
 
 \trois{som-amalg-gp}
 Dans la catégorie des groupes, la somme amalgamée
 de deux groupes~$G$ et~$H$ le long d'un groupe~$K$
 muni de deux morphismes~$g: K\to G$ et~$h:K\to H$
 existe, et est notée~$G*_KH$. 
 On peut par exemple la construire
 comme
 le quotient du produit libre~$G*H$
 par son plus petit sous-groupe distingué~$D$ contenant
 les éléments de la forme~$g(k)h(k)^{-1}$ pour~$k\in K$. Cette
 présentation a l'avantage de rendre triviale la vérification de la propriété universelle, 
 mais l'inconvénient de ne pas être très maniable en pratique : elle ne dit pas
 comment décider si deux mots de~$G*H$ ont même classe modulo~$D$. 
 
\medskip

Toutefois, {\em si~$g$
et~$h$
sont injectives}, 
on peut décrire~$G*_KH$
d'une façon plus tangible
(cela revient
peu ou prou
à exhiber un système 
de représentants explicite de représentants de
la congruence modulo~$D$), que nous allons 
maintenant
détailler. Pour alléger les notations, nous utilisons
les injections~$g$ et~$h$ pour identifier~$K$ à un sous-groupe
de~$G$ aussi bien que de~$H$.

\medskip
On fixe un système de représentants~$S$ des classes
 non triviales de~$G$ modulo~$K$, et un système de représentants~$T$
 des classes non triviales de~$H$ modulo~$K$. 
 Un élément
 de~$G*_KH$ est
 alors
 un mot~$x_1\ldots x_n$
 où chaque~$x_i$ est un élément
 de~$(K\setminus\{e\}) \coprod S\coprod T$
 et où : 
 
 \medskip
 - l'ensemble des indices~$i$ tels que~$x_i \in K\setminus \{e\}$ est vide ou égal à~$\{1\}$ ; 
 
- pour tout indice~$i\leq n-1$, on a~$x_{i+1}\in T$ si~$x_i\in S$, et~$x_{i+1}\in S$ si~$x_i\in T$. 

\medskip
On fait le
 produit de deux mots
 en les concaténant, puis en
 transformant le mot obtenu
 par l'algorithme suivant. 
  
\medskip
\begin{itemize}

\item[$\bullet$] Si le mot contient une séquence de la forme~$kk'$
où~$k$ et~$k'$ appartiennent à~$K$, et si~$z$ désigne le produit~$kk'$ 
calculé dans~$K$, alors : 

\begin{itemize}
 \item[$\diamond$] si~$z=e$ on supprime la séquence~$kk'$ ; 
 
 \item[$\diamond$] si~$z\neq e$ on remplace la séquence~$kk'$
 par la lettre~$z$. 
\end{itemize}

\item[$\bullet$] Si le mot contient une séquence
 de la forme~$sy$ avec~$s\in S$ et~$y\in K\setminus\{e\}$
 ou~$y\in S$ 
on pose~$z=sy\in G$ ; puis on 
écrit~$z=ks'$ avec~$k\in K$ et~$s'\in S\cup\{e\}$, et l'on procède comme suit : 
 
 \begin{itemize}
 \item[$\diamond$] si~$k=e$ et~$s'=e$ on supprime la séquence~$sy$ ; 
 
 \item[$\diamond$] si~$k=e$ et~$s'\neq e$ on remplace la séquence~$sy$
 par la lettre~$s'$ ; 
 
 \item[$\diamond$] si~$k\neq e$ et~$s'=e$ on remplace la séquence~$sy$
 par la lettre~$k$ ; 
 
\item[$\diamond$] si~$k\neq e$ et~$s'\neq e$ on remplace la séquence~$sy$
 par la séquence~$ks'$.

 \end{itemize}
 
 \item[$\bullet$] Si le mot contient une séquence
 de la forme~$ty$ avec~$t\in T$ et~$y\in K\setminus\{e\}$
 ou~$y\in T$ 
on pose~$z=ty\in H$ ; puis on 
écrit~$z=kt'$ avec~$k\in K$ et~$t'\in T\cup\{e\}$, et l'on procède comme suit :

 \begin{itemize}
 \item[$\diamond$] si~$k=e$ et~$t'=e$ on supprime la séquence~$ty$ ; 
 
 \item[$\diamond$] si~$k=e$ et~$t'\neq e$ on remplace la séquence~$ty$
 par la lettre~$t'$ ; 
 
 \item[$\diamond$] si~$k\neq e$ et~$t'=e$ on remplace la séquence~$ty$
 par la lettre~$k$ ; 
 
\item[$\diamond$] si~$k\neq e$ et~$t'\neq e$ on remplace la séquence~$ty$
 par la séquence~$kt'$.

 \end{itemize}
 
 \item[$\bullet$] On recommence jusqu'à ce qu'il n'y ait plus de séquence de l'une des trois
 formes évoquées ; le mot est alors sous la forme requise. 
 \end{itemize}
 
\medskip
L'élément
neutre de~$G*_KH$ est le mot vide. 
   
 \trois{som-amalg-ann}
 Dans la catégorie des anneaux
 commutatifs unitaires, la somme amalgamée de deux
 anneaux~$B$ et~$C$ le long d'un troisième anneau~$A$ existe, 
 c'est le produit tensoriel~$B\otimes_A C$, que nous verrons un peu
 plus loin. 
 
\medskip
 {\bf J'invite le lecteur à écrire lui-même 
 à propos des sommes disjointes et amalgamées, les «quelques tautologies»
 duales de celles vues plus haut sur  les produits cartésiens
 et produits fibrés. }
 
 \section{Adjonction} 
 \markboth{Le langage des catégories}{Adjonction}
 
 \deux{def-adj}
 {\bf Définition.}
 Soient~$\mathsf C$
 et~$\mathsf D$ deux catégories, et soient~$F : \mathsf C\to \mathsf D$ et~$G: \mathsf D\to \mathsf C$
 deux foncteurs covariants. On dit que~$F$ est un
 {\em adjoint à gauche de~$G$},
 ou que~$G$ est un {\em adjoint à droite de~$F$},
 ou que~$(F,G)$ est un {\em couple de foncteurs adjoints}, 
 s'il existe 
 une collection d'isomorphismes 
 $$\iota_{(A,B)} : \hom_{\mathsf D}(F(A), B)\simeq \hom_{\mathsf C}(A, G(B)),$$
 paramétrée par~$(A,B) \in {\rm Ob}\;\mathsf C\times {\rm Ob}\;\mathsf D$
 et fonctorielle en~$(A,B)$. 
 
 \trois{comm-adj}
 {\bf Commentaires.} La fonctorialité en~$(A,B)$ signifie que pour toute
 flèche~$u : A\to A'$ de~$\mathsf C$ et toute flèche~$v : B\to B'$ de~$\mathsf D$, 
 le diagramme
 
 $$\diagram\hom_{\mathsf D}(F(A'), B)\rrto_{\iota_{(A',B)}}^\simeq\dto_{\circ F(u)}&&\hom_{\mathsf C}(A',G(B))\dto^{\circ u}\\
 \hom_{\mathsf D}(F(A), B)\rrto_{\iota_{(A,B)}}^\simeq\dto_{v\circ}&&\hom_{\mathsf C}(A, G(B))\dto^{G(v)\circ}\\
 \hom_{\mathsf D}(F(A), B')\rrto_{\iota_{(A,B')}}^\simeq&&\hom_{\mathsf C}(A,G(B'))\enddiagram$$
 commute. 
 
 \trois{exo-sys-adj}
 {\bf Exercice.}
 Décrire la donnée des~$\iota_{(A,B)}$ 
 comme un morphisme de foncteurs. 
 
 \deux{lien-adj-yoneda}
 Soient~$\mathsf C$ et~$\mathsf D$ deux catégories, et soit~$G: \mathsf D\to \mathsf C$ 
 un foncteur covariant. 
 
 \trois{ad-impl-rep}
 Supposons que~$G$ admet un adjoint à gauche~$F$, et fixons un objet~$A$
 de~$\mathsf C$. Le foncteur covariant
 $$B\mapsto \hom_{\mathsf C}(A, G(B))$$
 de~$\mathsf D$ dans~$\ens$ est alors isomorphe
 à~$$B\mapsto \hom_{\mathsf D}(F(A), B).$$ Il est donc représentable, et~$F(A)$ est son 
 représentant.

\trois{rep-impl-ad} Supposons réciproquement que pour tout
objet~$A$
de~$\mathsf C$,
le foncteur~$B\mapsto \hom_{\mathsf C}(A,G(B))$
soit représentable, et notons~$F(A)$ son représentant. 

\medskip
Soit maintenant~$u : A\to A'$ une flèche de~$\mathsf C$. 
Elle induit un morphisme de foncteurs 
$$[B\mapsto \hom_{\mathsf C}(A',G(B))]\longrightarrow [B\mapsto \hom_{\mathsf C}(A,G(B))],$$
d'où par le lemme de Yoneda un morphisme~$F(u) : F(A)\to F(A')$. Ainsi,
$A\mapsto F(A)$ apparaît comme un foncteur covariant de~$\mathsf C$
vers~$\mathsf D$, qui est par construction un
adjoint 
à gauche de~$G$.

\deux{comment-adj}
{\bf Commentaires.}

\trois{comment-abus-adj}
Le lecteur attentif aura peut-être remarqué que nous sommes allés un peu vite : par construction, 
on dispose pour tout objet~$A$ d'un isomorphisme 
$$\hom_{\mathsf D}(F(A),B)\simeq \hom_{\mathsf C}(A, G(B))$$
qui est fonctoriel en~$B$, mais nous n'avons pas vérifié explicitement la fonctorialité en~$A$. Celle-ci
peut se justifier {\em grosso
modo}
ainsi -- nous vous invitons à écrire vous-mêmes les détails.

\medskip
Lorsque
nous
disons
que~$F(A)$ est le représentant de
$B\mapsto \hom_{\mathsf C}(A,G(B))$, nous commettons un abus. En effet, 
le représentant
consiste en réalité en la donnée de~$F(A)$ {\em et d'un objet universel},
à savoir ici un morphisme de~$A$ vers~$G(F(A))$ ; et le lemme de Yoneda
que nous avons appliqué ci-dessus fournit un morphisme~$F(u) : F(A)\to F(A')$
{\em compatible aux objets universels livrés avec~$A$ et~$A'$}. C'est cette compatibilité
qui garantit la fonctorialité attendue. 

\trois{unique-adj}
Supposons que~$G$ admette deux adjoints à gauche~$F_1$ et~$F_2$ ; donnons-nous
deux systèmes d'isomorphismes~$(\iota^1_{(A,B)})$
et~$(\iota^2_{(A,B)})$ décrivant l'adjonction respective des couples
$(F_1,G)$ et~$(F_2,G)$. On déduit alors du lemme de Yoneda (appliqué, 
comme au~\ref{comment-abus-adj}
ci-dessus, pour chaque~$A\in {\rm Ob}\;\mathsf C$, avec prise en compte des objets universels)
qu'il existe un unique isomorphisme~$F_1\simeq F_2$ compatible
aux systèmes $(\iota^1_{(A,B)})$
et~$(\iota^2_{(A,B)})$, dans un sens que nous laissons au lecteur le soin
de préciser. Nous nous permettrons pour cette raison de
parler de {\em l'}
adjoint à gauche de~$G$.

\trois{ad-g-d}
Nous venons d'évoquer diverses questions
liées à l'existence éventuelle d'un adjoint
{\em à gauche}
d'un foncteur. Nous vous invitons à
énoncer les résultats analogues concernant l'adjonction
{\em à droite}.

\deux{ex-fonct-adj}
{\bf Exemples de foncteurs adjoints.}

\trois{ex-adj-oubli-amod} Soit~$A$ un anneau commutatif
unitaire et soit~$O_1 : \amod\to \ens$ le foncteur d'oubli. Pour tout ensemble~$X$, 
notons~$L_1(X)$ le~$A$-module libre engendré par~$X$, c'est-à-dire
le~$A$-module~$\bigoplus_{x\in X}A\cdot e_x.$

Pour tout~$A$-module~$M$ et tout ensemble~$X$ on a un isomorphisme 
naturel $$\hom_{\ens}(X,O_1(M))\simeq \hom_{\amod}(L_1(X),M), $$
qui associe à une application~$\phi : X\to M$ son «prolongement
par linéarité»,
{\em i.e.}  l'unique application~$A$-linéaire
de~$L_1(X)$ dans~$M$ envoyant~$e_x$ sur~$\phi(x)$ pour tout~$x$. Il est visiblement
fonctoriel en~$X$
et~$M$ ; en conséquence, $L_1$ est l'adjoint à gauche
de~$O_1$. 

\trois{ex-adj-oubli-gp}
Soit~$O_2: \gp \to \ens$ le foncteur d'oubli. Pour tout ensemble~$X$, 
on note~$L_2(X)$ le
{\em groupe libre}
engendré par~$X$, défini comme suit. On introduit pour chaque
élément~$x$ de~$X$ un symbole~$x^{-1}$, et l'on note~$X^{-1}$
l'ensemble des~$x^{-1}$ pour~$x$ parcourant~$X$. Un élément de~$L_2(X)$ est alors un
«mot {\em réduit}
sur l'alphabet~$X\cup X^{-1}$», 
c'est-à-dire une suite finie d'éléments de~$X\cup X^{-1}$ ne comprenant aucune
séquence de deux termes consécutifs de la forme~$xx^{-1}$ ou~$x^{-1}x$
pour~$x\in X$. On multiplie deux mots en les concaténant
puis en réduisant le mot obtenu, c'est-à-dire en éliminant 
les séquences~$xx^{-1}$ ou~$x^{-1}x$
tant qu'on en rencontre ; l'élément neutre est le mot vide.

\medskip
Pour tout
groupe~$G$
et tout ensemble~$X$ on a un isomorphisme 
naturel $$\hom_{\ens}(X,O_2(G))\simeq \hom_\gp(L_2(X),G), $$
qui associe à une application~$\phi : X\to G$ son «prolongement
homomorphique»,
{\em i.e.}  l'unique morphisme de groupes
de~$L_2(X)$ dans~$G$ envoyant~$x$ sur~$\phi(x)$ pour tout~$x\in X$
(défini par ma formule que le lecteur imagine). Il est visiblement
fonctoriel en~$X$
et~$G$ ; en conséquence, $L_2$ est
l'adjoint à gauche de~$O_2$. 

\trois{ex-adj-oubli-alg}
Soit~$A$
un anneau commutatif unitaire, et soit~$O_3$ le foncteur
d'oubli~$B\mapsto (B,\times)$ de la catégorie des~$A$-algèbres vers celle
des monoïdes commutatifs avec unité. Soit~$L_3$ le foncteur qui 
associe à un tel monoïde~$\Gamma$
l'algèbre~$\bigoplus_{\gamma
\in \Gamma} A_{e_\gamma}$
où~$e_{\gamma}\cdot e_{\gamma'}=e_{\gamma\gamma'}$
pour tout~$(\gamma, \gamma')$. 

\medskip
Pour toute~$A$-algèbre~$B$ et tout monoïde commutatif avec unité~$\Gamma$, 
on a un isomorphisme 
naturel $$\hom_{\mathsf{mon}}(\Gamma,O_3(B))\simeq \hom_{A-\mathsf{Alg}}(L_3(\Gamma),B), $$
qui associe à un morphisme~$\phi : \Gamma \to (B,\times)$ 
l'unique morphisme de~$A$-algèbres de~$L_3(\Gamma)$ dans~$B$ envoyant~$e_\gamma$
sur~$\phi(\gamma)$ pour tout~$\gamma$.  Il est visiblement
fonctoriel en~$\Gamma$
et~$B$ ; en conséquence, $L_3$ est
l'adjoint à gauche de~$O_3$.

\trois{ex-adj-oubli-ab}
Soit~$I$ le foncteur d'inclusion de~$\ab$ dans~$\gp$, 
et soit~$A$ le foncteur de~$\gp$ dans~$\ab$ qui envoie un groupe~$G$
sur son abélianisé~$G/<ghg^{-1}h^{-1}>_{(g,h)\in G^2}$.  Soit~$H$ un groupe
abélien
et soit~$G$ un groupe. On a un isomorphisme naturel

$$\hom_\gp(G, I(H))\simeq \hom_\ab(A(G), H),$$
qui envoie un morphisme~$\phi : G\to H$ vers l'unique morphisme
$A(G)\to H$ par lequel il se factorise (son
existence est assurée par
la propriété
universelle du quotient). Il est visiblement
fonctoriel en~$G$
et~$H$ ; en conséquence, $A$ est l'adjoint à gauche
de~$I$. 

\trois{ex-adj-top}
Soit~$O_4$ le foncteur d'oubli 
de~$\top$
dans~$\ens$, et soit~${\rm Dis}$ 
(resp.~${\rm Gro}$)
le foncteur de~$\ens$ vers~$\top$
consistant à munir un ensemble
de la topologie discrète (resp. grossière). 

\medskip
Toute application ensembliste entre espaces
topologiques est automatiquement continue dès lors que 
la topologie de sa source
(resp. son but) est discrète (resp. grossière). 

\medskip
On dispose donc pour toute ensemble~$X$ et tout espace 
topologique~$T$ de deux bijections naturelles

$$\hom_\ens(X,O_4(T))\simeq \hom_\top ({\rm Dis}(X), T)$$
$${\rm et}\;\;\hom_\ens(O_4(T),X)\simeq \hom_\top(T, {\rm Gro}(X)).$$ Elles
sont visiblement fonctorielles en~$X$ et~$T$ ; en conséquence, 
${\rm Dis}$ est l'adjoint à gauche de~$O_4$, et~${\rm Gro}$ 
est son adjoint à droite. 

\trois{ex-adj-hom-prod}
Soient~$X,Y$ et~$Z$ trois ensembles. On dispose
d'une bijection 
$$\hom_\ens(X\times Y,Z)\simeq \hom_\ens(X, \hom_\ens(Y,Z)) : $$
elle envoie~$f$ sur~$x\mapsto [y\mapsto f(x,y)]$, et sa réciproque 
envoie~$g$ sur~$(x,y)\mapsto g(x)(y)$. Elle est visiblement fonctorielle
en~$X,Y,Z$.
En conséquence, $Y$ étant fixé, le foncteur
$X\mapsto X\times Y$ est l'adjoint à gauche de~$Z\mapsto \hom_\ens(Y,Z)$.

\trois{ex-adj-qi}
{\bf Exercice.}
Montrez que si un foncteur covariant~$F$ admet un quasi-inverse~$G$ alors~$G$ est 
à la fois adjoint à~$F$ à gauche et à droite. 

\section{Limites inductives et projectives}
\markboth{Le langage des catégories}{Limites inductives et projectives}

\deux{intro-lim-indproj}
Cette section va être consacrée à deux techniques
très générales, et en un sens duales l'une de l'autre, 
de construction d'objets dans une catégorie : la formation
(lorsque c'est possible)
de la {\em limite projective} et de la {\em limite projective}
d'un {\em diagramme}. 

\medskip
Les sommes disjointes et amalgamées ainsi que l'objet initial apparaîtront 
{\em a posteriori}
comme des exemples de limites inductives ; les produits cartésiens et fibrés ainsi 
que l'objet final apparaîtront {\em a posteriori}
comme des exemples de limites projectives. Si nous avons
choisi d'en faire une présentation directe précédemment,
c'est en raison de leur importance
particulière, et également 
pour vous permettre de vous faire la main
sur des cas un peu plus concrets et explicites que ceux que nous allons
maintenant aborder.  

\deux{def-diagramma}
Soit~$\mathsf C$ une catégorie.
Un
{\em diagramme}
dans~$\mathsf C$ est la donnée d'une famille~$(X_i)_{i\in I}$ 
d'objets de~$\mathsf C$ et pour tout~$(i,j)$ d'un ensemble~$E_{ij}$ de flèches
de~$X_i$ vers~$X_j$. Mentionnons incidemment 
qu'un diagramme~$((X_i), (E_{ij}))$ est
dit~{\em commutatif}
si pour tout triplet~$(i,j,k)$, toute~$f\in E_{ij}$, toute~$g\in E_{jk}$ et toute~$h\in E_{ik}$ 
on a~$h=g\circ f$. 

\deux{mor-diag}
Soit~$\sch D=((X_i), (E_{ij}))$ un diagramme dans~$\mathsf C$ et soit~$Y$ un objet de~$\mathsf C$. 

\trois{def-mor-d-y}
Un {\em morphisme de~$\sch D$ vers~$Y$} est une famille~$(a_i : X_i\to Y)_i$ de flèches 
telle que l'on ait pour tout couple~$(i,j)$ et pour tout~$f\in E_{ij}$ l'égalité~$a_j\circ f=a_i$. On note
$\hom_{\mathsf C} (\sch D, Y)$ l'ensemble des morphismes de~$\sch D$ vers~$Y$.

\trois{def-mor-y-d}
Un {\em morphisme de~$Y$ vers~$\sch D$}
est une famille~$(b_i : Y\to X_i)_i$ de flèches 
telle que l'on ait pour tout couple~$(i,j)$ et pour tout~$f\in E_{ij}$ l'égalité~$f\circ b_i =b_j$. On note
$\hom_{\mathsf C} (Y, \sch D)$ l'ensemble des morphismes de~$Y$ vers~$\sch D$. 

\deux{def-lims}
{\bf Définition des limites inductive et projective}. 
Soit~$\sch D$ un diagramme de~$\mathsf C$. 

\trois{def-lim-ind}
Si le foncteur covariant de~$\mathsf C$ dans~$\ens$ qui envoie un objet~$Y$ sur~$\hom_{\mathsf C}(\sch D,Y)$
est représentable, son représentant est appelé
la
{\em limite inductive}
du diagramme~$\sch D$
et est noté~$\limind \sch D$. Notons qu'il est fourni avec une famille
de morphismes~$(\lambda_i : X_i\to \limind \sch D)$. 

\trois{def-lim-proj}
Si le foncteur contravariant de~$\mathsf C$ dans~$\ens$ qui envoie un objet~$Y$ sur~$\hom_{\mathsf C}(Y,\sch D)$
est représentable, son représentant est appelé
la
{\em limite projective}
du diagramme~$\sch D$
et est noté~$\limproj \sch D$. Notons qu'il est fourni avec une famille
de morphismes~$(\mu_i : \limproj \sch D \to X_i)$.

\deux{ex-lim-ind}
{\bf Premiers exemples de limite inductive.}
Les faits suivants découlent tautologiquement des définitions. 

\trois{lim-diag-vide}
La limite inductive du {\em diagramme vide}
existe si et seulement si~$\mathsf C$ a un objet initial, et si c'est le cas
les deux coïncident.

\trois{lim-diag-somme-disj}
Soient~$X$ et~$Y$ deux objets de~$\mathsf C$. La limite
inductive du diagramme
$$X\;\;\;Y$$ (deux objets, pas de morphismes) existe si et seulement si~$X\coprod Y$
existe, et si c'est le cas les deux coïncident. 

\trois{lim-diag-somme-disj-gen}
Les deux exemples précédents se généralisent comme suit : soit
$\sch D=(X_i)_{i\in I}$
un diagramme constitué d'une famille d'objets, sans morphismes. Dire que~$\limind \sch D$
existe revient à dire que le foncteur~$Y\mapsto \prod \hom_{\mathsf C} (X_i,Y)$ est représentable ; si c'est le 
cas la limite est notée~$\coprod X_i$ et est appelée la {\em somme disjointe}
des~$X_i$. Se donner un morphisme {\em depuis}
$\coprod X_i$, c'est se donner un morphisme depuis chacun des~$X_i$. 

\trois{lim-diag-somme-amalg}
Soient~$S\to X$ et~$S\to Y$ deux morphismes de~$\mathsf C$. La limite
inductive du diagramme
$$
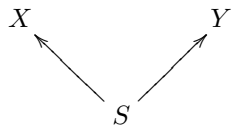
\xymatrix{
X&&Y\\
&S\ar[ul]\ar[ur]&}$$ existe si et seulement si~$X\coprod_SY$
existe, et si c'est le cas les deux coïncident. 

\deux{ex-lim-proj}
{\bf Premiers exemples de limite projective.}
Les faits suivants découlent tautologiquement des définitions. 

\trois{limpro-diag-vide}
La limite projective du {\em diagramme vide}
existe si et seulement si~$\mathsf C$ a un objet final, et si c'est le cas
les deux coïncident.

\trois{lim-diag-prod-fib}
Soient~$X$ et~$Y$ deux objets de~$\mathsf C$. La limite
projective du diagramme
$$X\;\;\;Y$$ existe si et seulement si~$X\times Y$
existe, et si c'est le cas les deux coïncident. 

\trois{lim-diag-prod-fib-gen}
Les deux exemples précédents se généralisent comme suit : soit
$\sch D=(X_i)_{i\in I}$
un diagramme constitué d'une famille d'objets, sans morphismes. Dire que~$\limproj \sch D$
existe revient à dire que le foncteur~$Y\mapsto \prod \hom_{\mathsf C} (Y,X_i)$ est représentable ; si c'est le 
cas la limite est notée~$\prod X_i$ et est appelée le {\em produit}
des~$X_i$. Se donner un morphisme {\em vers}
$\prod X_i$, c'est se donner un morphisme vers chacun des~$X_i$. 

\trois{lim-diag-somme-amalg}
Soient~$X\to S$ et~$Y\to S$ deux morphismes de~$\mathsf C$. La limite
projective du diagramme
$$\xymatrix{
X\ar[dr]&&Y\ar[dl]\\
&S&}$$ existe si et seulement si~$X\times_SY$
existe, et si c'est le cas les deux coïncident. 

\deux{tauto-sous-singleton}
Soit~$\mathsf C$ une catégorie, et soit~$\sch D=((X_i), (E_{ij})$ un diagramme
dans~$\mathsf C$. Supposons que pour tout objet~$Y$ de~$\mathsf C$ et tout~$i$,
l'ensemble
$\hom(X_i,Y)$ soit ou bien vide, ou bien un singleton ; c'est par exemple le cas si~$\mathsf C$ est 
la catégorie des algèbres sur un anneau~$A$ et si~$X_i$ est pour tout~$i$ un quotient de~$A$. 

\medskip
Dans ce cas, pour tout objet~$Y$ de~$\mathsf C$, l'ensemble $\hom_{\mathsf C}(\sch D,Y)$ est vide s'il existe~$i$
tel que~$\hom_{\mathsf C}(X_i,Y)$ soit vide, et est un singleton sinon. En effet, il est clair que s'il existe~$i$ tel que~$\hom_{‘\mathsf C}(X_i,Y)=\varnothing$
alors~$\hom_{\mathsf C}(\sch D,Y)=\varnothing$. Dans le cas contraire, $\hom_{\mathsf C}(X_i,Y)$ est pour tout~$i$ un singleton~$\{f_i\}$, et le seul
élément éventuel de~$\hom_{\mathsf C}(\sch D,Y)$ est donc la famille~$(f_i)$ ; mais celle-ci appartient effectivement à~$\hom_{\mathsf C}(\sch D,Y)$ : comme
$f_i$ est le seul élément de~$\hom_{\mathsf C}(X_i,Y)$ pour tout~$i$, les conditions du type~$f_j\circ f=f_i$ sont {\em automatiquement}
vérifiées. 

\medskip
On voit donc que~$\hom_{\mathsf C}(\sch D,Y)$ ne dépend que de~$(X_i)$, et pas de la famille de morphismes~$(E_{ij})$ ; il en va dès lors de même
de l'existence de~$\limind \sch D$, et de sa valeur le cas échéant. 

\medskip
Nous laissons le lecteur énoncer les assertions duales à propos des limites projectives. 

\subsection*{Exemples de catégories admettant
des limites inductives}

\deux{exemples-limind}
Soit~$\mathsf C$ une catégorie
et soit~$\sch D=((X_i), (E_{ij}))$ un diagramme dans~$\mathsf C$. 
Nous allons
mentionner un certain nombre de cas dans lesquels~$\limind \sch D$ existe
et décrire celle-ci, en laissant les vérifications (éventuellement fastidieuses, mais
triviales) au lecteur. 

\trois{limind-ens}
{\em Supposons que
$\mathsf C=\ens$.} Dans ce cas, $\limind \sch D$ est 
le quotient de
l'union disjointe ensembliste~$\coprod X_i$ par la relation d'équivalence 
engendrée par les relations~$x\sim f(x)$ où~$x\in X_i$ et~$f\in E_{ij}$ pour un certain 
couple~$(i,j)$, et la flèche structurale~$X_j\to \limind \sch D$  est pour tout~$j$ l'application
composée $\lambda_j \colon X_j\to \coprod X_i\to \limind \sch D$.
On remarque que~$\limind \sch D=\bigcup \lambda_i(X_i)$.

\trois{limind-top}
{\em Supposons que~$\mathsf C=\top$.}
La limite
inductive de~$\sch D$ coïncide alors en tant qu'ensemble avec celle
construite au~\ref{limind-ens}
ci-dessus ; la topologie de~$\limind\sch D$
est obtenue en munissant~$\coprod X_i$ de la topologie
d'union disjointe (pour laquelle une partie est ouverte si et seulement si sa trace sur chacun des~$X_i$ est ouverte),
puis en passant à la topologie quotient. Pour tout~$j$, le morphisme
structural~$X_j\to \limind \sch D$ est l'application~$\lambda_j$ de~\ref{limind-ens}
(qui est continue). 

\trois{limind-ann}
{\em Supposons que~$\mathsf C=\ann$ ou que~$\mathsf C=\aalg$ pour un certain anneau~$A$.}
La limite inductive de~$\sch D$ sera décrite plus loin, lorsque nous
aurons vu le produit tensoriel. 

\trois{limin-damod}
{\em Supposons que~$\mathsf C=\amod$ pour un certain anneau~$A$.}
Pour tout~$i$, notons~$h_i$ le morphisme canonique de~$X_i$ vers
la somme directe (externe)
$\bigoplus X_i$. La limite inductive de~$\sch D$ est alors
le
quotient de~$\bigoplus X_i$ par son sous-module engendré par les
éléments $h_j(f(x))-h_i(x)$
où~$(i,j)\in I^2$, où~$x\in X_i$ et où~$f\in E_{ij}$ ; pour tout~$j$, la flèche structurale
$\mu_j \colon X_j\to \limind\sch D$ est la composée de~$h_j$ et du morphisme quotient ; 
on remarque que~$\limind \sch D$ est la somme (interne, et non directe en général) des~$\mu_i(M_i)$. 

\trois{limind-gp}
{\em Supposons que~$\mathsf C=\gp$.} Nous allons
nous contenter d'indiquer
ici succinctement la construction de~$\limind \sch D$ (nous ne nous en servirons pas). 
On commence
par construire
le {\rm produit libre} $\star_i X_i$ (qui sera la somme disjointe des~$X_i$ dans la catégorie des groupes) : 
c'est l'ensemble des mots dont les lettres sont des éléments non triviaux des~$X_i$, deux lettre consécutives
n'appartenant jamais au même~$X_i$. Pour tout~$j$, soit~$\iota_j$
le morphisme canonique~$X_j\hookrightarrow \star_i X_i$. 
La limite inductive de~$\sch D$ est alors le quotient de~$\star_i X_i$ par son plus petit
sous-groupe distingué contenant~$\iota_i(x)\iota_j(f(x))^{-1}$ pour tout couple~$(i,j)$, 
tout~$x$ appartenant à~$X_i$ et toute~$f\in E_{ij}$. Pour tout~$j$, la flèche structurale
$\nu_j \colon X_j\to \limind \sch D$ est la composée de~$\iota_j$ et du morphisme quotient ; 
on remarque que~$\limind \sch D$ est engendré par les~$\nu_i(X_i)$. 

\deux{syst-filtr-ind}
{\bf Diagrammes commutatifs filtrants.}
Soit~$\mathsf C$ une catégorie et soit~$I$ un ensemble
{\em préordonné}, c'est-à-dire muni d'une relation~$\leq$
que l'on suppose réflexive et transitive, mais pas nécessairement
antisymétrique ; l'exemple qu'on peut
avoir en tête est celui d'un anneau muni de la relation de divisibilité. 

\trois{desc-idiagramme}
On peut voir~$I$ comme une catégorie, dans laquelle 
l'ensemble~$\hom(i,j)$ est un singleton si~$i\leq j$, et est vide sinon. 
Donnons-nous un foncteur de~$I$ vers~$\mathsf C$, c'est-à-dire : 

\medskip
$\bullet$ pour tout~$i\in I$, un objet~$X_i$ de~$\mathsf C$ ; 

$\bullet$ pour tout~$(i,j)\in I^2$ avec~$i\leq j$, une flèche~$f_{ij}$ de~$X_i$ vers~$X_j$, 
de sorte que~$f_{ii}={\rm Id}_{X_i}$ pour tout~$i$,
et que~$f_{jk}\circ f_{ij}=f_{ik}$
dès que~$i\leq j\leq k$. 

\medskip
Si l'on pose~$E_{ij}=\{f_{ij}\}$ si~$i\leq j$ et~$E_{ij}=\varnothing$
sinon, la famille~$\sch D=((X_i), (E_{ij}))$ est alors un diagramme commutatif dans
la catégorie~$\mathsf C$, que l'on dira {\em défini}
par le foncteur~$I\to \mathsf C$ dont on est parti. 

\trois{def-diag-filtr}
On suppose
de surcroît
que l'ensemble d'indices~$I$ est {\em filtrant}, 
c'est-à-dire qu'il est non vide et que pour tout couple~$(i,j)$ d'éléments de~$I$ il existe~$k\in I$
avec~$i\leq k$ et~$j\leq k$ (c'est par exemple le cas d'un anneau muni de la divisibilité : considérer
le produit de deux éléments) ; il revient au même de demander que toute partie finie de~$I$
admette un majorant, grâce au petit jeu usuel : 
l'existence d'un majorant de la partie vide garantira que~$I$ est non vide. 
Nous
résumerons la situation en disant
que~$\sch D$ est un diagramme {\em commutatif filtrant}.
On suppose que sa limite inductive existe et pour tout~$i$
l'on notera
encore~$\lambda_i : X_i\to \limind \sch D$
la flèche canonique.

\trois{syst-filtr-ind-alg}
{\em On suppose que~$\mathsf C=\ens$, $\gp$, $\ann$, $\aalg$ ou~$\amod$ (pour un certain anneau~$A$). }
Nous allons tout d'abord donner une description
de la limite inductive {\em ensembliste}
de~$\sch D$ qui est spéficique au contexte filtrant et diffère (en apparence
seulement, évidemment) de celle
de~\ref{limind-ens}. Soit~$\sch R$
la relation sur la somme disjointe ensembliste~$\coprod X_i$ 
définie comme suit : si~$a\in X_i$ et si~$b\in X_j$ alors~$a\sch R b$ si et seulement
si il existe~$k$ majorant~$i$ et~$j$ tel que~$f_{ik}(a)=f_{jk}(b)$. 
C'est une relation d'équivalence en vertu du caractère filtrant de~$I$ (qui sert pour la transitivité,
la réflexivité et la symétrie étant évidentes) ; on vérifie que le quotient~$L:=\coprod X_i/\sch R$, muni 
pour tout~$j$ de la composée~$\ell_j \colon X_j \to \coprod X_i \to L$, s'identifie à la limite inductive
ensembliste de~$\sch D$ (c'est-à-dire qu'il satisfait la propriété universelle correspondante). 

\medskip
On peut résumer cette construction en disant que la limite inductive ensembliste~$(L, (\ell_i))$
de~$\sch D$ est caractérisée
par l'égalité~$L=\bigcup \ell_i(X_i)$ et par le fait que pour tout~$(i,j)$, tout~$a\in X_i$ et tout~$b\in X_j$, on a 
$\ell_i(a)=\ell_j(b)$ si et seulement si il existe~$k$ majorant~$i$ et~$j$ tel que~$f_{ik}(a)=f_{jk}(b)$.

\medskip
Nous allons maintenant expliquer brièvement comment munir~$L$ d'une structure d'objet de~$\mathsf C$ 
pour laquelle les~$\ell_i$ seront des morphismes, et qui fera de~$(L, (\ell_i))$ la limite inductive de~$\sch D$ 
dans la catégorie~$\mathsf C$.

\medskip
Il n'y a rien à faire si~$\mathsf C=\ens$. Dans les autres cas, il faut définir une ou plusieurs lois
sur~$L$. Indiquons par exemple comment on procède lorsque~$\mathsf C$
est la catégorie des anneaux, les autres situations
se traitant de façon analogue. Soient~$x$ et~$y$ deux éléments de~$L$ ; on peut écrire~$x=\ell_i(x_i)$ et~$y=\ell_j(y_j)$
pour~$i$ et~$j$ convenables, avec~$x_i\in X_i$ et~$y_j\in X_j$. On choisit alors un indice~$k$
majorant~$i$ et~$j$ ; on a~$x=\ell_k(x_k) $ et~$y=\ell_k(y_k)$, où l'on a posé
$x_k=f_{ik}(x_i)$ et~$y_k=f_{jk}(y_j)$. On vérifie sans peine
que 
l'élément~$\ell_k(x_k+y_k)$ de~$L$ ne dépend que de~$x$ et~$y$, et on le note~$x+y$ ; de même, 
l'éléments~$\ell_k(x_ky_k)$ de~$L$ ne dépend
que de~$x$ et~$y$, et on le note~$xy$. Il est facile de voir que~$(L,+,\times)$ 
est un anneau, dont les éléments~$0$ et~$1$
sont respectivement égaux à~$\ell_i(0)$ et~$\ell_i(1)$ pour n'importe quel~$i\in I$ (il est donc vital pour
les définir qu'il existe un tel~$i$ : ici
intervient le fait que~$I\neq \varnothing$). Les~$\ell_i$ sont des morphismes d'anneaux par construction, 
et il est immédiat que~$(L, (\ell_i))$ satisfait la propriété universelle de la limite inductive dans~$\ann$. 

\trois{filtrant-commute}
{\em Remarque.} Dans les catégories algébriques (groupes, anneaux, modules, algèbres), 
l'ensemble sous-jacent à une limite inductive {\em filtrante}
est donc la limite inductive des ensembles sous-jacents. 
C'est faux en général pour les limites inductives quelconques. Pensez par exemple au produit libre de deux groupes, ou même 
tout simplement aux objets initiaux : 
l'objet initial de~$\ann$ est~$\ZZ$, celui de~$\gp$ 
est~$\{e\}$, et
celui de~$\ens$ est~$\varnothing$.

\subsection*{Exemples de catégories admettant
des limites projectives}

\deux{exemples-limproj}Soit~$\mathsf C$ une catégorie
et soit~$\sch D=((X_i), (E_{ij}))$ un diagramme dans~$\mathsf C$. 
Nous allons
mentionner un certain nombre de cas dans lesquels~$\limproj \sch D$ existe
et décrire celle-ci, en laissant une fois encore
les vérifications au lecteur.

\trois{limproj-ens}
{\em Supposons que
$\mathsf C=\ens$.} Dans ce cas, $\limproj \sch D$ est 
le sous-ensemble de~$\prod X_i$ 
formée des familles~$(x_i)$ telles que
$f(x_i)=x_j$ pour tout couple~$(i,j)$ et toute~$f\in E_{ij}$.
Pour tout~$j$, la flèche structurale~$\limproj \sch D\to X_j$ est simplement
la restriction de la~$j$-ième projection.

\trois{limproj-top}
{\em Supposons que~$\mathsf C=\top$.}
La limite
projective de~$\sch D$ coïncide alors en tant qu'ensemble avec celle
construite au~\ref{limproj-ens}.
La topologie de~$\limproj\sch D$
est induite par la topologie produit sur~$\prod X_i$, et la flèche
structurale $\limproj \sch D\to X_j$ est encore la
restriction de la~$j$-ième projection, qui est continue. 

\trois{limproj-ann}
{\em Supposons que~$\mathsf C=\ann, \gp, \amod$ ou~$\aalg$ pour un certain anneau~$A$.}
Une fois encore, l'ensemble sous-jacent à~$\limproj \sch D$
coïncide avec la limite projective ensembliste de~$\sch D$ décrite au~\ref{limproj-ens}. 
Sa structure d'objet de~$\mathsf C$ s'obtient 
en munissant le produit ensembliste~$\prod X_i$ de la structure
d'objet de~$\mathsf C$
définie composante par composante, puis en remarquant que
~$\limproj \sch D$ en est un sous-objet
au sens de~$\mathsf C$ ; 
et pour tout~$j$, la flèche structurale~$\limproj \sch D\to X_j$ est encore
la restriction de la~$j$-ième projection, qui est un morphisme de~$\mathsf C$. 

\deux{lim-sur-s}
Soit~$\mathsf C$ une catégorie et soit~$\sch D=((X_i), (E_{ij}))$
un diagramme de~$\mathsf C$. Soit~$S$ un objet de~$\mathsf C$. 

\trois{sur-s-ind}
Supposons donné pour tout~$i$ un morphisme~$p_i : X_i\to S$
de sorte que~$p_j\circ f=p_j$ pour tout~$(i,j)$ et toute~$f\in E_{ij}$. La
famille~$(p_i)$ induit alors un morphisme~$p \colon \limind \sch D\to S$. 

\medskip
Grâce aux~$p_i$, on peut voir~$\sch D$ comme un diagramme
dans la catégorie~$\mathsf C/S$ (\ref{cat-semiclass}). Il est alors
tautologique (faites l'exercice) que la limite inductive de~$\sch D$ dans la
catégorie~$\mathsf C/S$ existe et s'identifie à~$p \colon \limind \sch D \to S$. 

\trois{sur-s-proj}
Supposons donné pour tout~$i$ un morphisme~$q_i : S\to X_i$
de sorte que~$f\circ q_i=q_j$ pour tout~$(i,j)$ et toute~$f\in E_{ij}$. La
famille~$(p_i)$ induit alors un morphisme~$q\colon S\to \limproj \sch D$. 

\medskip
Grâce aux~$q_i$, on peut voir~$\sch D$ comme un diagramme
dans la catégorie~$S\backslash \mathsf C$ (\ref{cat-semiclass}). Il est alors
tautologique (faites l'exercice) que la limite projective de~$\sch D$ dans la
catégorie~$S\backslash \mathsf C$ existe et s'identifie à~$q \colon S \to \limproj \sch D.$ 

\subsection*{Adjonction et passage à la limite}

\deux{lim-et-adj}
Soient~$\mathsf C$ et~$\mathsf D$ deux catégories, et soit~$\sch D=((X_i), (E_{ij}))$
un diagramme dans~$\mathsf C$. Soit~$F$ un foncteur de~$\mathsf C$ vers~$\mathsf D$. Supposons
que~$F$ admet un adjoint à droite~$G$, et que~$\limind \sch D$ existe
dans~$\mathsf C$. Soit~$Y$ un objet de~$D$, et soit~$F(\sch D)$ 
le diagramme de~$\mathsf D$ déduit de~$\sch D$ en appliquant~$F$
à chacun de ses constituants (objets et morphismes). 
Notons~$E$ l'ensemble des triplets~$(i,j,f)$ avec~$f\in E_{ij}$. 
L'ensemble~$\hom_{\mathsf D}(F(\sch D), Y)$ est égal à
$$\{(f_i : F(X_i)\to Y), f_i\circ F(f)=f_j\;\;\;\forall(i,j,f)\in E\}$$
$$=\{(g_i : X_i\to G(Y)), g_i\circ f=g_j\;\;\;\forall(i,j,f)\in E\}\;\;\;\text{(par adjonction)}$$ 
$$=\hom_{\mathsf C}(\sch D, G(Y))$$
$$=\hom_{\mathsf C}(\limind \sch D,G(Y))\;\;\;\text{(par définition de la limite inductive)}$$
$$=\hom_{\mathsf D}(F(\limind \sch D), Y)\;\;\;\text{(par adjonction)}.$$
En conséquence, $\limind F(\sch D)$ existe et est égal à~$F(\limind \sch D)$ : le foncteur~$F$
préserve les limites inductives. 

\medskip
Nous laissons le lecteur démontrer par une méthode
analogue (ou en utilisant les catégories opposées)
que si~$F$ admet un adjoint à gauche, il préserve les limites projectives. 

\deux{comment-preserv}
{\bf Adjonction et préservations des limites : exemples.}

\trois{adj-reserv-ind}
Nous avons vu au~\ref{limind-top}
que dans la catégorie des espaces
topologiques, l'ensemble sous-jacent à une limite inductive
est la limite inductive des ensembles sous-jacents ; autrement
dit, le foncteur oubli de~$\top$ vers~$\ens$ préserve les limites inductives. 
Il a une excellente raison pour ce faire : il admet en effet un adjoint à droite
(\ref{ex-adj-top}). 

\medskip
Nous avons par contre signalé
en~\ref{filtrant-commute}
qu'en général, l'ensemble
sous-jacent à une limite inductive d'anneaux, groupes, modules ou algèbres n'est 
pas la limite inductive des ensembles sous-jacents ; il en résulte que les foncteurs
oubli correspondant n'ont pas d'adjoint à droite. 

\trois{adj-reserv-proj}
Nous avons vu au~\ref{limproj-top}
et au~\ref{limproj-ann}
que dans les catégories des espaces
topologiques, des groupes, des anneaux, 
des modules et des algèbres l'ensemble sous-jacent à une limite projective
est la limite projective
des ensembles sous-jacents ; autrement
dit, les foncteur oubli de~$\top$, $\ann$ et~$\gp$ 
vers~$\ens$ préservent les limites projectives.  
Ils ont une excellente raison pour ce faire : ils admettent
en effet tous un adjoint à gauche. Celui-ci est décrit
en~\ref{ex-adj-top} pour les espaces topologiques, en~\ref{ex-adj-oubli-amod}
pour les $A$-modules et en~\ref{ex-adj-oubli-gp}
pour les groupes ; en 
ce qui concerne les anneaux et les algèbres, nous
laissons le lecteur construire lui-même les anneau et algèbre libres sur un ensemble~$X$ 
dans l'esprit de ce qui a été fait en~\ref{ex-adj-oubli-amod}, \ref{ex-adj-oubli-gp}
et~\ref{ex-adj-oubli-alg}.

%
%
%
%
%

\chapter{Algèbre commutative}

 \section{Localisation}
 
 \deux{pre-intro-loc}
 Soit~$A$ un anneau. Lorsqu'on se donne un sous-ensemble~$P$
  de~$A$, on sait construire un anneau «défini à partir de~$A$, en {\em décrétant}
  que les éléments de~$P$ sont nuls, et en n'imposant aucune autre contrainte
  que celle-ci, et ses conséquences découlant de la théorie générale des anneaux» : 
  c'est  le quotient~$A/(P)$. Celui-ci est caractérisé par sa propriété universelle,
  c'est-à-dire encore par le foncteur (ici, covariant) qu'il représente.   
  
  \medskip
  C'est en fait une illustration d'un phénomène assez général : à chaque fois
  lorsqu'on veut intuitivement {\em imposer} une contrainte, et {\em seulement}
  cette contrainte, la construction rigoureuse qui répond à ce caprice s'exprime 
  en termes de propriété universelle, ou encore de foncteur à représenter. 
  
 \deux{intro-loc}
 Nous allons en voir un nouvel exemple avec ce qu'on appelle
  la {\em localisation}. Soit~$A$ un anneau, et soit~$S$ un sous-ensemble
  de~$A$. Le but intuitif est de construire un objet à partir de~$A$ en imposant aux
  éléments de~$S$ d'être inversibles -- et rien d'autre. Techniquement, 
  on s'intéresse au foncteur covariant
  $$F: B\mapsto \{f\in {\rm Hom}(A,B), f(s)\in B\ti\;\;\forall \;s\in S\},$$
  et nous allons montrer de deux façons différentes qu'il est représentable, c'est-à-dire qu'il
  existe un anneau~$C$ et un morphisme~$g : A\to C$ tels que :
  
  \medskip
  $\bullet$ $g(S)\subset C\ti$
   ; 
  
  $\bullet$ pour tout~$f: A\to B$ tel que~$f(S)\subset B\ti$, il existe un unique~$h : C\to B$
  faisant commuter le diagramme
  $$\diagram A\rto^f\dto_g&B\\
  C\urto_h&\enddiagram\;\;.$$
  
 \deux{methode-bourin}
 {\em Première preuve.}   Posons~$$C=A[T_s]_{s\in S}/(sT_s-1)_{s\in S}$$
 et~$g=a\mapsto \overline a$. 
 
 \trois{preuve-bourin}
 Soit~$f : A\to B$ un morphisme
 tel que~$f(S)\subset B\ti$ ; nous allons montrer l'existence d'un unique~$h : C\to B$
 tel que~$f=h\circ g$. Voyant~$B$ et~$C$ comme des~$A$-algèbres {\em}
 via~$f$ et~$g$ respectivement, cela revient à montrer l'existence d'un unique
 morphisme~$h$ de~$A$-algèbres de~$C$ vers~$B$. 
 
 \medskip
 {\em Unicité.}
 Soit~$h$ un tel morphisme. On a alors 
 pour tout~$s\in S$ les égalités
$$1= h(\overline{T_ss})=h(\overline{T_s})h(\overline s)=h(\overline{T_s})h(g(s))
=h(\overline{T_s})f(s),$$
et donc~$h(\overline{T_s})=f(s)^{-1}$. Comme les~$\overline{T_s}$
engendrent la~$A$-algèbre~$C$, il y a au plus un tel morphisme~$h$. 

\medskip
{\em Existence.}
Soit~$\phi$ l'unique morphisme de~$A$-algèbres de~$A[T_s]_{s\in S}$
vers~$B$ qui envoie~$T_s$ sur~$f(s)^{-1}$ pour tout~$s\in S$. 
On a~$\phi(T_ss-1)=0$ 
pour tout~$s\in S$ et~$\phi$ induit donc par passage au quotient un 
morphisme d'algèbres~$h : C\to B$.

\trois{comment-bourin} Cette preuve est la plus économique,
  et en un sens la plus naturelle : on a {\em forcé}
  les éléments
  de~$S$ à être inversibles, en adjoignant 
  formellement à~$A$ un 
  symboles~$T_s$ pour chaque élément~$s$ de~$S$, et
  en imposant par décret l'égalité~$sT_s=1$. 
    
  \medskip
 Mais elle
 présente un défaut : il est {\em en pratique}
 extrêmement difficile d'arriver à dire quoi que ce soit
 sur 
 la~$A$-algèbre
 $A[T_s]_{s\in S}/(sT_s-1)_{s\in S}$ ; c'est un cas où
 la connaissance du foncteur représenté par un objet qui, 
 {\em en théorie}, caractérise l'objet en question à isomorphisme près,
 n'est pas suffisante. 
 
 Par exemple,  il semble {\em a priori} impossible
 de donner un critère simple permettant 
 de savoir si~$A[T_s]_{s\in S}/(sT_s-1)_{s\in S}$
 est nul ou non. Nous allons donc donner une autre construction du représentant
 de~$F$. 
 
 \trois{intro-produit-en-croix}
{\em Réduction au cas d'une
partie multiplicative.}
Convenons
 de dire qu'une partie~$T$ de~$A$ est {\em multiplicative} si elle
 contient~$1$ et si~$ab\in T$ dès que~$a\in T$ et~$b\in T$ (on
 peut une fois encore condenser la définition en style bourbakiste,
 en demandant simplement
 que~$T$ soit stable par produit fini, ce qui la force à contenir~$1$ puisque ce dernier
 est le produit vide). L'ensemble~$\widehat S$ des produits finis d'éléments de~$S$ (en incluant~$1$
 qui est le produit vide) est visiblement la plus petite partie multiplicative
 de~$A$ contenant~$S$ ; nous dirons que c'est la partie 
 multiplicative
 {\em engendrée par~$S$.}
 Si~$B$ est un anneau et si~$f: A\to B$ est un morphisme,
 il est immédiat que~$f(s)$ est inversible pour tout~$s\in S$ si et seulement
 si~$f(s)$ est inversible pour tout~$s\in \widehat S$. On peut donc, pour étudier 
 le foncteur~$F$, 
 remplacer~$S$ par~$\widehat S$ ; autrement dit,
 on s'est ramené au cas où~$S$ est multiplicative.  
 
 \trois{def-produit-en-croix}
 On définit alors sur~$A\times S$ la relation~$\mathscr R$ suivante : $(a,s){\mathscr R}(b,t)$ si et seulement
 si il existe~$r\in S$ tel que~$r(at-bs)=0$. 
 On vérifie que c'est une relation d'équivalence, et l'on note~$S^{-1}A$ le quotient correspondant. 
 
 Les formules
 $$((a,s); (b,t))\mapsto (at+bs, st)\;{\rm et}\;((a,s);(b,t))\mapsto (ab,st)$$ 
 passent au quotient, et définissent deux lois~$+$ et~$\times$ sur~$S^{-1}A$ qui en font un anneau commutatif. 
 
 Si~$(a,s)\in A\times S$, on écrira~$\frac a s$ au lieu de~$\overline{(a,s)}$. Cette notation permet de disposer
 des formules naturelles
 $$\frac a s + \frac b t=\frac{at+bs}{st}\;{\rm et} \;\frac a s \cdot \frac b t=\frac{as}{bt},$$
 et l'on a~$$\frac a s =\frac b t \;\iff\;\exists r\in S, r(at-bs)=0.$$
 
 L'application~$a\mapsto \frac a 1$ est un morphisme d'anneaux de~$A$ dans~$S^{-1}A$, et 
 si~$s\in S$ alors~$\frac s 1$ est inversible, d'inverse~$\frac 1 s$. 
 
 \medskip
 {\em Le couple~$(S^{-1}A,a\mapsto \frac a 1)$ représente le foncteur~$F$.}
 En effet, soit~$f: A\to B$ un morphisme tel que~$f(S)\subset B\ti$ ; nous allons montrer
 qu'il existe un unique morphisme~$h$
 de~$S^{-1}A$ dans~$B$
 tel que~$h(\frac a 1)=f(a)$ pour tout~$a\in A$. 
 
 \medskip
 {\em Unicité.} Soit~$h$ un tel morphisme. On a alors pour tout~$(a,s)\in A\times S$
 les égalités
 $$h\left(\frac a s\right)=h\left(\frac a 1 \cdot \frac 1 s\right)=h\left(\frac a 1\right)h\left(\left(\frac s 1\right)^{-1}\right)=f(a)f(s)^{-1},$$
 d'où l'unicité. 
 
 {\em Existence.}
 On vérifie immédiatement que l'application~
 $$A\times S\to B, (a,s)\mapsto f(a)f(s)^{-1}$$
 passe au quotient par~$\sch R$. Elle induit donc une application~$h: S^{-1}A\to B$, 
 qui envoie toute fraction~$\frac a s$
 sur~$f(a)f(s)^{-1}$. Par un calcul explicite, on s'assure que~$h$ est un morphisme d'anneaux,
 et l'on a bien $h(\frac a 1)=f(a)$ pour tout~$a\in A$.

\medskip
On dit
 que~$S^{-1}A$ est le {\em localisé} de~$A$ par rapport à la partie
 multiplicative~$S$.  
 
 \deux{comment-localis}
 {\bf Commentaires.}
 
 \trois{croix-complique}
 La condition d'égalité entre fractions de~$S^{-1}A$
 est plus compliquée que le bon vieux produit en croix traditionnel ; c'est le prix
 à payer pour travailler avec des anneaux quelconques, {\em i.e.} non nécessairement intègres ni réduits. 
 Notons toutefois que si~$S$ ne contient pas de diviseurs de zéro
--
c'est par exemple le cas si~$A$ est intègre et si~$0\notin S$ -- la condition
 «il existe~$r\in S$ tel que~$r(at-bs)=0$» équivaut à la relation usuelle~«$at-bs=0$». 
 
 \trois{asur1-pas-in}
 La flèche~$A\to S^{-1}A$ n'est pas injective en général, c'est la raison pour 
 laquelle on préfère souvent écrire~$\frac a 1$ et non~$a$. Son noyau est facile à décrire : 
 c'est l'ensemble des éléments~$a$ de~$A$ tels qu'il existe~$r\in S$ vérifiant l'égalité~$ra=0$. 
 Une fois encore, les choses se simplifient si~$S$ ne contient pas de diviseurs 
 de zéros (et donc en particulier si~$A$ est intègre et si~$0\notin S$) :
 on voit immédiatement que sous ces hypothèses,~$a\mapsto \frac a 1$ est injective. 
 
\trois{smoins1a-nul}
L'anneau~$S^{-1}A$ est nul si et seulement si~$1=0$ dans~$S^{-1}A$, c'est-à-dire
 encore si et seulement si~$\frac 1 1=0$, donc si et seulement si il existe~$r\in S$ tel que~$r\cdot1=0$. 
 {\em Autrement dit,~$S^{-1}A$ est nul si et seulement si~$0\in S$}. 
 
 \trois{singleton-as}
 On a défini
 {\em l'anneau} $S^{-1}A$
 comme représentant du foncteur~$F$ défini en~\ref{intro-loc}.
 Par un raisonnement
 en tout point analogue à celui tenu en~\ref{asurp-alg},
 on en déduit que 
 pour toute~$A$-algèbre~$(B, f\colon A \to B)$, 
 l'ensemble des morphismes de~$A$-algèbres de~$S^{-1}A$ dans~$B$ est un singleton 
 si~$f(S)\subset B\ti$, et est vide sinon ; puis
 que la~{\em $A$-algèbre}~$S^{-1}A$
 représente le foncteur de~$\aalg$
 dans~$\ens$ qui envoie~$(B, f \colon  \to B)$ sur~$\{*\}$
 si~$f(S)\subset B\ti$ et sur~$\varnothing$ sinon.

\deux{ex-smoins1A}
{\bf Exemples.}
 
 \trois{corps-fractions}
 Soit~$A$ un anneau  intègre. Le sous-ensemble~$S:=A\setminus\{0\}$ en 
 est une partie multiplicative. Le localisé~$S^{-1}A$ est non nul
 puisque~$0\notin S$, et si~$\frac a s$ est un élément non nul
 de~$S^{-1}A$ alors~$a\neq 0$ ; en conséquence, $a\in S$ et~$\frac a s$
 est inversible d'inverse~$\frac s a$. Il s'ensuit
 que~$S^{-1}A$ est un corps, 
 appelé {\em corps des fractions de~$A$}
 et souvent noté~${\rm Frac}\;A$.

 \medskip
 Puisque~$0\notin S$, l'anneau intègre~$A$ s'injecte dans~${\rm Frac}\;A$. 
 Ce dernier est précisément le plus petit corps contenant~$A$, dans le sens suivant : 
 pour tout corps~$K$ et tout morphisme
 {\em injectif}
 $A\hookrightarrow K$, il existe un unique plongement~${\rm Frac}\;A\hookrightarrow K$
 tel que le diagramme
 $$\diagram A\rto\dto&K\\
 {\rm Frac}\;A\urto&\enddiagram$$
 commute. En effet, la flèche injective~$A\to K$ envoie tout élément
 de~$S=A\setminus\{0\}$ sur un élément non nul, et partant inversible, de~$K$ ; l'assertion
 requise est alors un cas particulier 
 de la propriété universelle de~$S^{-1}A$. 
 
 \trois{a-f}
 Soit~$A$ un anneau et soit~$f\in A$. La partie multiplicative~$S$ engendrée
 par~$f$ est~$\{f^n\}_{n\in \NN}$, et le localisé correspondant est le plus souvent
 noté~$A_f$. On déduit de la construction par quotient
 décrite en~\ref{methode-bourin}
 que~$A_f\simeq A[T]/(fT-1)$. 
 
 \medskip
 En vertu de~\ref{smoins1a-nul}, 
 l'anneau~$A_f=0$ est nul si et seulement si~$S$
 contient~$0$, 
 c'est-à-dire si et seulement si~$f$ est nilpotent. 
 
 Donnons une preuve alternative
 de ce fait. L'anneau~$A_f$ est nul si et seulement
 si~$A[T]/(Tf-1)=0$, c'est-à-dire si et seulement si~$1-Tf$ est inversible
 dans~$A[T]$. Or~$1-Tf$ est inversible dans~$A[[T]]$, 
 d'inverse~$g=\sum f^iT^i$. Par unicité de l'inverse (lorsqu'il existe), 
 on voit que~$1-Tf$ est inversible dans~$A[T]$ si et seulement
 si~$g\in A[T]$, c'est-à-dire si et seulement si~$f$ est nilpotent.

\deux{smoinsun-fonct}
{\bf Fonctorialité.}
Soit~$f:A\to B$ un morphisme d'anneaux. Soit~$S$ une partie multiplicative
 de~$A$, et soit~$T$ une partie multiplicative de~$B$ telle que~$f(S)\subset T$ (par
 exemple,~$T=f(S)$). La flèche composée~$A\to B\to T^{-1}B$ envoyant chaque
 élément de~$S$ sur un inversible, elle induit une flèche
 de~$S^{-1}A$ vers~$T^{-1}B$, donnée par les
 formules~$$\frac a s \mapsto \frac {f(a)}{f(s)}.$$ 
 
 \deux{mise-en-garde-mult}
 {\bf À propos des éléments inversibles de~$S^{-1}A$.}
 Soit~$A$ un anneau et soit~$S$ une partie multiplicative de~$A$. 
 Soit~$S^\flat$ l'ensemble
 des éléments de~$A$ qui deviennent inversibles dans~$S^{-1}A$
 (ce n'est pas une notation standard, nous ne nous en servirons que dans la brève discussion
 qui suit). 
 C'est une partie multiplicative de~$A$ qui contient~$S$ par définition,
 mais on prendra garde que l'inclusion~$S\subset S^\flat$
 est en général {\em stricte} ; par exemple, $S^\flat$ 
 contient automatiquement~$A\ti$, qui n'a aucune raison
 {\em a priori}
 d'être contenu dans~$S$. Si~$T$ est une partie multiplicative
 contenant~$S$ on a~$S^\flat \subset T^\flat$, puisque~$A\to T^{-1}A$
 se factorise par~$S^{-1}A$. 
 
 \trois{desc-s-flat}
 Il est en fait possible de décrire précisément~$S^\flat$ : 
 c'est l'ensemble des éléments~$a\in A$ tels qu'il existe~$b$ 
 vérifiant~$ab\in S$. En effet, s'il existe un tel~$b$ alors $\frac{ab}1=\frac a 1\cdot \frac b 1$ 
 est inversible dans~$S^{-1}A$, ce qui force~$\frac a 1$ (et~$\frac b 1$) à l'être aussi, et~$a\in S^\flat$.

\medskip
Réciproquement, si~$a\in S^\flat$, il existe~$\alpha\in A$ et~$s\in S$ tels que $\frac a 1\cdot \frac\alpha s = \frac 1 1$,
ce qui veut dire qu'il existe~$t\in S$ tel que~$t(a\alpha-s)=0$ ; on a donc~$t\alpha a=st\in S$, ce qui achève
de prouver l'assertion requise. 
 
\trois{desc-s-flat-funct}
 On peut également donner une description fonctorielle de~$S^\flat$ : c'est l'ensemble~$E$
 des éléments~$a$
 de~$A$ tel que~$\phi(a)\in B\ti$ pour tout morphisme d'anneaux~$\phi \colon A\to B$ 
 vérifiant~$\phi(S)\subset B\ti$. En effet, si~$a\in E$ l'image de~$a$ dans~$S^{-1}A$ est inversible, puisque
 $A\to S^{-1}A$ 
 envoie chaque élément de~$S$ sur un inversible ; ainsi, $E\subset S^\flat$. 
 
 Réciproquement, si~$a\in S^\flat$ et si~$\phi : A\to B$ est un morphisme d'anneaux tel que
 $\phi(S)\subset B\ti$ alors comme~$\phi$ se factorise (d'une unique manière)
 par
 la flèche~$A\to S^{-1}A$ et comme~$a$ s'envoie par définition de~$S^\flat$
 dans~$(S^{-1}A)\ti$, on
 a~$\phi(a)\in B\ti$ et~$a\in E$. 
 
 \trois{sflat-versuss}
 Comme~$S\subset S^\flat$ on a un morphisme naturel de~$S^{-1}A$ dans~$(S^\flat)^{-1}A$,
 morphisme qui induit en vertu de~\ref{desc-s-flat-funct}
 un isomorphisme entre les foncteurs représentés par ces deux~$A$-algèbres ; en conséquence,
 ce morphisme est un isomorphisme. On peut bien entendu s'en assurer de manière plus terre-à-terre
 à partir de la description explicite de~$S^\flat$ donnée au~\ref{desc-s-flat} ; l'exercice est laissé au lecteur. 
 
 \medskip
 On déduit de ce fait, ou directement du~\ref{desc-s-flat-funct}
 ci-dessus, que~$(S^\flat)^\flat=S^\flat$. 
 
 \deux{limind-smoinsun-a}
 {\bf Localisation et limite inductive : le cas filtrant.} 
 Soit~$A$ un anneau, soit~$I$ un ensemble
 préordonné filtrant (\ref{syst-filtr-ind} et~\ref{syst-filtr-ind-alg}) et soit~$\Sigma$ une partie multiplicative
 de~$A$. Supposons donnée pour tout~$i\in I$ une partie multiplicative
 de~$A$ contenue dans~$\Sigma$, et faisons les hypothèses suivantes : 
 
 \medskip
 $\bullet$ $\Sigma=\bigcup S_i$ ; 
 
 $\bullet$ pour tout~$(i,j)\in I^2$ avec~$i\leq j$, les éléments de~$S_i$ deviennent inversibles dans~$S_j^{-1}A$ (il suffit pour cela que~$S_i\subset S_j$, 
 mais ce n'est pas nécessaire, {\em cf.}
~\ref{mise-en-garde-mult}
 {\em et sq.}). 
 
 \medskip
 Posons~
 $$\sch D=((S_i^{-1}A)_{i\in I}, (S_i^{-1}A\to S_j^{-1}A)_{i\leq j}).$$ C'est un
 diagramme commutatif filtrant dans la catégorie des~$A$-algèbres. La famille
des flèches canoniques~$S_i^{-1}A\to  \Sigma^{-1}A$ définit un morphisme
de~$\sch D$ dans~$\Sigma^{-1}A$, dont nous allons montrer
qu'il induit un isomorphisme

$$\limind \sch D\simeq \Sigma^{-1}A.$$

En vertu de la description explicite des limites inductives filtrantes donnée en~\ref{syst-filtr-ind-alg}, cet énoncé
équivaut à la
validité des deux assertions qui suivent. 

\trois{asi-asigma-sur}
{\em Soit~$\alpha\in \Sigma^{-1}A$. Il existe~$i$ tel que~$\alpha$ provienne de~$S_i^{-1}A$.}
Mais c'est évident : par définition, $\alpha$ s'écrit~$\frac a s$ avec~$a\in A$ et~$s\in \Sigma$. Comme~$\Sigma$
est la réunion des~$S_i$, il existe~$i$ tel que~$s\in S_i$, et~$\alpha$ est dès lors égal à l'image de l'élément~$\frac a s$ de~$S_i^{-1}A$. 

\trois{asi-sigma-inj}
{\em Soit~$(i,j)\in I^2$, soit~$\alpha\in S_i^{-1}A$ et soit~$\beta\in S_j^{-1}A$. Supposons que~$\alpha$ et~$\beta$ ont même image dans~$\Sigma^{-1}A$ ; 
il existe alors un majorant~$k$ de~$\{i,j\}$ tel que~$\alpha$ et~$\beta$ aient déjà même image dans~$S_k^{-1}A$.}
Pour le voir, on écrit~$\alpha=\frac a s$ et~$\beta =\frac b t$ avec~$(a,b)\in A^2, s\in S_i$ et~$t\in S_j$. Comme~$\alpha$ et~$\beta$
ont même image dans~$\Sigma^{-1}A$, il existe~$\sigma\in \Sigma$ tel que~$\sigma(at-bs)=0$. Puisque~$\Sigma$ est la réunion des~$S_i$, il existe~$\ell$
tel que~$\sigma\in S_{\ell}$. Choisissons un majorant~$k$ de~$\{i,j,\ell\}$. Comme~$k\leq \ell$, l'élément~$\sigma$ est inversible dans~$S_k^{-1}A$, et l'égalité
$\sigma(at-bs)=0$ implique donc que~$at-bs=0$ dans~$S_k^{-1}A$, et partant que~$\frac a s=\frac b t$ dans~$S_k^{-1}A$, ce qui termine la preuve. 

\trois{limind-loc-multicat}
{\em Remarque.}
Comme la limite inductive de~$\sch D$ dans la catégorie des~$A$-algèbres
«est» aussi sa limite inductive dans la catégorie des anneaux ainsi que dans celle des~$A$-modules
(cela découle de~\ref{syst-filtr-ind-alg}), il
résulte de ce qui précède que~$\Sigma^{-1}A$ s'identifie
également à la limite inductive de~$\sch D$
dans la catégorie des anneaux et dans celle des~$A$-modules. 

 \deux{limind-smoinsun-fonct}
 {\bf Localisation et limite inductive : le cas général}. Nous allons donner 
 dans ce qui suit une preuve plus conceptuelle du fait que le morphisme de~$A$-algèbres
 $\limind \sch D\to \Sigma^{-1}A$ de~\ref{limind-smoinsun-a}
 est un isomorphisme ; cette nouvelle démonstration a l'avantage de marcher sous des hypothèses
 nettement plus faibles que nous allons maintenant énoncer. 
 
 \medskip
 On conserve les notations~$A$ et~$\Sigma$ de~\ref{limind-smoinsun-a}. On désigne
 par contre maintenant par~$I$ un ensemble quelconque, et l'on se donne pour tout~$i\in I$
 une partie multiplicative~$S_i$ de~$A$ contenue dans~$\Sigma$ ; nous supposons simplement
 que les~$S_i$ engendrent multiplicativement~$\Sigma$. 
 
 \medskip
 On se donne un diagramme~$\sch D$ dans la catégorie des~$A$-algèbres dont la famille
 d'objets est $(S_i^{-1}A)_i$ ; {\em on n'impose rien à la famille des flèches de~$\sch D$}.

 \trois{preuve-limind-loc-fonct}
 Soit~$(B, f \colon A \to B)$ une~$A$-algèbre ; pour toute partie multiplicative~$S$
 de~$A$, l'ensemble~$\hom_{\aalg}(S^{-1}A,B)$ est un singleton si~$f(S)\subset B\ti$, et est vide sinon. 
 Il s'ensuit en vertu de~\ref{tauto-sous-singleton} que~$\hom_{\aalg}(\sch D,B)$ est un singleton si~$f(S_i)\subset B\ti$ 
 pour tout~$i$, et est vide sinon, et ce {\em quelles que soient les flèches de~$\sch D$}. 
 
 \medskip
 Comme les~$S_i$ engendrent multiplicativement~$\Sigma$, on a~$f(S_i)\subset B\ti$ pour tout~$i$ si et seulement si~$f(\Sigma)\subset B\ti$. 
 En conséquence, on dispose d'un isomorphisme fonctoriel en~$B$ entre~$\hom_{\aalg}(\sch D,B)$ et~$\hom_{\aalg} (\Sigma^{-1}A,B)$, 
 qui montre que~$\Sigma^{-1}A$ s'identifie à la limite inductive du diagramme~$\sch D$.

\trois{comment-limind-loc-fonct}
{\em Commentaires.}
Insistons à nouveau 
sur le fait que ce qui précède vaut pour
{\em tout}
diagramme~$\sch D$ dont la famille d'objets est~$(S_i^{-1}A)_i$. C'est par exemple
le cas du diagramme sans flèches : la~$A$-algèbre~$\Sigma^{-1}A$ est ainsi
en particulier la somme
disjointe des~$S_i^{-1}A$. 

\medskip
On prendra garde qu'ici, contrairement à ce qui valait plus haut, 
({\em cf.}~\ref{limind-loc-multicat}),  l'identification entre~$\Sigma^{-1}A$ et~$\limind \sch D$
n'est avérée {\em a priori}
que dans la catégorie des~$A$-algèbres, mais pas dans celle des anneaux ou des~$A$-modules ; en effet,
la «coïncidence» des limites inductives dans les différentes catégories est une spécificité du cas filtrant, 
prise en défaut en général. 

\medskip
Ainsi, on déduit de ce qui précède que~$\ZZ[1/6]$ est la somme disjointe
de~$\ZZ[1/2]$ et~$\ZZ[1/3]$ dans la catégorie des~$\ZZ$-algèbres, c'est-à-dire des anneaux ; 
mais leur somme disjointe comme~$\ZZ$-modules est égale à~$\ZZ[1/2]\bigoplus \ZZ[1/3]$, qui est un
sous-module strict de~$\ZZ[1/6]$.

\section{Idéaux premiers et maximaux}\label{IDEAUX}
\markboth{Algèbre commutative}{Idéaux premiers et maximaux}

Cette section ne
contient à proprement parler aucun 
résultat nouveau. Elle vise simplement 
à présenter une approche des idéaux premiers et maximaux
qui est 
sans doute un peu différente de celle dont vous avez l'habitude,
et imprègne
(le plus souvent implicitement)
la géométrie algébrique à la Grothendieck.

\deux{rappels-defprime}
{\bf Rappels des définitions.}
Soit~$A$ un anneau. Un idéal~$\got p$
de~$A$ est dit {\em premier}
s'il est strict et si~$ab\in \got p\Rightarrow a\in \got p$
ou~$b\in \got p$. Il revient au même de demander
que~$A/\got p$ soit un anneau intègre. 

\medskip
Un idéal~$\got m$
de~$A$ est dit
{\em maximal}
s'il est strict et s'il est maximal {\em en tant
qu'idéal strict}. Cela revient à demander que~$A/\got m$
ait exactement deux idéaux, à savoir~$\{0\}$ et~$A/\got m$ ; 
autrement dit, $\got m$ est maximal si et seulement si~$A/\got m$ est un corps. 
On déduit de cette dernière caractérisation que tout idéal maximal est premier. 

\deux{exist-premier}
Soit~$A$ un anneau et soit~$I$ un idéal strict de~$A$. On déduit
immédiatement du lemme de Zorn que~$I$ est contenu dans un idéal
maximal. Si~$A$ est non nul il possède donc
un idéal maximal : appliquer
ce qui précède avec~$I=\{0\}$, qui est alors strict. 

\medskip
On voit en particulier 
que tout anneau non nul possède un idéal {\em premier}. Notons que cette propriété
est {\em a priori}
plus faible que l'existence d'un idéal maximal, mais 
elle ne peut pas à ma connaissance être établie
directement. 
 
\deux{premiers-morphismes}
{\bf Idéaux premiers
et morphismes vers les corps.}
Soit~$A$ un anneau et soit~$f$ un morphisme de~$A$ vers un corps~$K$. 
 Le
 noyau de~$f$ est visiblement
 un idéal premier. Réciproquement, soit~$\mathfrak p$ un idéal premier
 de~$A$ ; la flèche composée~$A\to A/{\mathfrak  p}\hookrightarrow {\rm Frac}(A/{\mathfrak p})$ a pour
 noyau~$\mathfrak p$. Ainsi, les idéaux premiers sont exactement
 {\em les noyaux de morphismes dont le but est un corps.} 
 
 \trois{desc-memenoyau}
 On peut donc décrire un idéal premier de~$A$ comme 
 une classe d'équivalence de 
 morphismes~$(A\to K)$  où~$K$ est un corps, pour la relation d'équivalence
 «avoir même noyau». 
 
 \trois{desc-sous-extension}
 Cette relation admet une description alternative : 
 si~$K$ et~$L$ sont deux corps, deux morphismes~$A\to K$ et~$A\to L$ ont même noyau
 si et seulement si il existe un corps~$F$ et un diagramme
 commutatif~
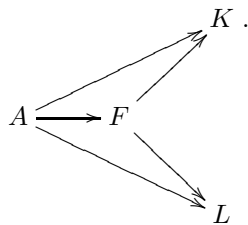
$$\xymatrix{ &&K\\
 A\ar[r]\ar[rrd]\ar[rru]&F\ar[ur]\ar[dr]&\\
 && L}.$$ 
 
 En effet, si un tel diagramme existe 
 alors~$${\rm Ker}(A\to K)={\rm Ker}(A\to L)={\rm Ker}(A\to F)$$
 puisque~$F\to K$ et~$F\to L$ sont injectifs en tant que morphismes de corps. 
 
 Réciproquement, supposons que~$A\to K$ et~$A\to L$ aient même noyau~$\mathfrak p$. 
 En vertu des propriétés universelles du quotient et du corps des fractions, la 
 flèche~$A\to K$ admet une unique factorisation sous la
 forme~$$A\to {\rm Frac}(A/{\mathfrak p})\hookrightarrow K,$$
 et il en va de même de~$A\to L$. Il existe donc un  
 diagramme comme ci-dessus avec~$F={\rm Frac}(A/{\mathfrak p})$. 
 
 \trois{rem-asurp-pluspetit}
 On a en fait montré au~\ref{desc-sous-extension}
 ci-dessus
 qu'un morphisme~$A\to K$ appartient à la classe
 qui correspond à~$\mathfrak p$
 si et seulement si il se factorise par
 la flèche~$A\to {\rm Frac}\;(A/\got p)$, 
 et que
 si c'est le cas cette factorisation est unique. En d'autres termes, 
 le morphisme canonique~$A\to {\rm Frac}\;(A/\got p)$
 est le {\em plus petit}
 élément de la classe de morphismes~$A\to K$ associée à~$\got p$.

 \deux{idemax-morphisme}
 {\bf Idéaux maximaux et surjection vers un corps.}
 Soit~$\mathfrak m$ un idéal maximal
 de~$A$. Le quotient~$A/\got m$ est un corps, et la
 flèche~$A\to A/{\mathfrak m}$ est surjective. 
 
 \medskip
 Réciproquement, si~$K$ est un corps et si~$f: A\to K$ est surjective, 
 alors comme~$K$ s'identifie à~$A/{\rm Ker} f$, le noyau de~$f$ est un
 idéal maximal de~$A$. 
 
 \trois{conclu-idmax}
 Ainsi, un idéal maximal de~$A$ peut être vu comme une classe d'équivalence
 de {\em surjections}~$A\to K$, où~$K$ est un corps, pour la relation d'équivalence
 «avoir même noyau». Et si~$A\to K$ et~$A\to L$ sont deux
 surjections ayant même noyau~$\mathfrak m$, les corps~$K$ et~$L$ s'identifient tous deux
 à~$A/{\mathfrak m}$ comme~$A$-algèbres. 
 Il y a donc en fait à isomorphisme
 canonique près {\em une seule} surjection dans la classe d'équivalence qui correspond
 à un idéal maximal donné~$\got m$ : c'est la surjection quotient
 de~$A$ vers~$A/\got m$.

\trois{idmax-idpremier}
{\em Idéaux maximaux au sein des idéaux premiers.} Donnons-nous un idéal
premier~$\got p$ de~$A$. Il correspond à une classe d'équivalence
de morphismes~$A\to K$, où~$K$ est un corps, 
classe
dont~$A\to {\rm Frac}(A/\got p)$ est le plus petit élément. Par ce qui précède,
l'idéal~$\got p$ est maximal si et seulement si il existe, dans la classe d'équivalence
qui lui correspond, un morphisme surjectif. Mais cela revient à demander que le plus petit
morphisme de la classe, à savoir~$A\to {\rm Frac}(A/\got p)$, soit surjectif, c'est-à-dire encore
que~$ {\rm Frac}(A/\got p)=A/\got p$, et donc que~$A/\got p$ soit un corps ; 
on retrouve bien
(heureusement !)
la définition
usuelle.  

\deux{exemples-idprimes}{\bf Exemple : le cas de~$\ZZ$.}
Nous donnons ci-dessous la liste des idéaux premiers et maximaux de~$\ZZ$, 
en donnant leur description 
du point
de vue des morphismes vers les corps. 

\medskip
$\bullet$ L'idéal~$(0)$ ; il correspond à la classe des morphismes injectifs~$\ZZ\to K$, c'est-à-dire
des morphismes~$\ZZ\to K$ où~$K$ est un corps de caractéristique nulle. Le plus petit morphisme
de cette classe est l'inclusion~$\ZZ\hookrightarrow \QQ$, laquelle n'est pas surjective :~$(0)$ n'est pas
maximal. 

$\bullet$ Pour tout nombre premier~$p$, l'idéal~$(p)$ ;  il correspond à la classe des morphismes~$\ZZ\to K$
de noyau~$p\ZZ$, c'est-à-dire
des morphismes~$\ZZ\to K$ où~$K$ est un corps de caractéristique~$p$. Le plus petit morphisme
de cette classe est la flèche naturelle~$\ZZ\to \FF_p$, qui est surjective :~$(p)$ est maximal. 

\deux{spec-contrav}
{\bf Fonctorialité contravariante du spectre.} Si~$A$ est un anneau, on note~${\rm Spec}\; A$ le
{\em spectre} de~$A$, c'est-à-dire l'ensemble
des idéaux premiers de~$A$ (nous verrons plus tard, lors du cours sur les schémas, que~${\rm Spec}\; A$ 
peut être muni d'une topologie, et même d'une structure supplémentaire).

La flèche~$A\mapsto {\rm Spec}\;A$ est de manière naturelle un foncteur {\em contravariant} ; 
nous allons  
donner deux descriptions de la flèche
de~${\rm Spec}\;B$ vers~${\rm Spec}\;A$
induite par un morphisme d'anneaux~$f: A\to B$.

  \medskip
  1) {\em Description dans le langage classique}. À un idéal premier~$\got q$ de~$B$, on associe l'idéal~$f^{-1}(\got q)$ de~$A$,
  dont on vérifie qu'il est premier. 
  
  2) {\em Description du point de vue des morphismes vers les corps}. Si~$K$ est un corps et~$B\to K$ un morphisme,
  le noyau de la flèche composée~$A\to B\to K$ ne dépend que de celui de~$B\to K$ (c'est son
  image réciproque dans~$A$). On peut ainsi sans ambiguïté associer à la classe
  d'équivalence de~$B\to K$ la classe d'équivalence de la composée~$A\to B\to K$.

 \subsection*{Anneaux locaux}

%
%
%
 \deux{prop-def-annloc}
 {\bf Proposition-définition.}
 {\em
 Soit~$A$ un anneau. Les assertions suivantes sont équivalentes.  
 
  \medskip
 i) $A$ possède un et un seul idéal maximal. 
 
 ii) L'ensemble des éléments non inversibles de~$A$ est un idéal de~$A$. 
 
 \medskip
Si elles sont satisfaites
on dit que~$A$ est un anneau 
{\em local}. Son unique idéal maximal est alors précisément l'ensemble de ses éléments non inversibles}.

\medskip
{\em Démonstration.}
Supposons que~i) est vraie, et soit~$\got m$ l'unique idéal 
maximal de~$A$. Si un élément de~$a$ appartient à~$\got m$,
il n'est pas inversible puisque~$\got m$ est strict par définition. 
Réciproquement, si~$a$ n'est pas inversible, l'idéal~$(a)$ est strict, et est donc
contenu dans un idéal maximal qui ne peut être que~$\got m$ ; ainsi~$a\in \got m$,
et l'ensemble des éléments non inversibles de~$A$ est exactement~$\got m$. 

\medskip
Supposons maintenant que~ii) est vraie, et soit~$\got m$ l'ensemble des éléments
non inversibles de~$A$. Comme~$1$ est inversible, il n'appartient pas à~$\got m$, 
qui est donc un idéal strict. 
Par ailleurs, si~$I$ est un idéal strict de~$A$, il ne contient aucun élément inversible et est donc contenu 
dans~$\got m$. Il s'ensuit aussitôt que ce dernier est l'unique idéal maximal de~$A$.~$\Box$

\deux{corps-annloc}
{\bf Exemple trivial.}
Tout corps est un anneau local, dont~$(0)$ est l'unique
idéal maximal. 

\deux{ann-local-geodiff}
{\bf Exemple géométrique.} Nous  donner un exemple
qui illustre
la pertinence de l'épithète «local». 
Soit~$X$ un espace
 topologique et soit~$x\in X$. On considère l'ensemble des
 couples~$(U,f)$ où~$U$ est un voisinage ouvert de~$x$ et~$f\in {\mathscr C}^0(U,\RR)$,
 sur lequel on met la relation d'équivalence suivante : $(U,f)\sim(V,g)$ si et seulement
 si il existe un voisinage ouvert~$W$ de~$x$ dans~$U\cap V$ tel que~$f_{|W}=g_{|W}$. 
 L'ensemble quotient~$A$ hérite alors d'une structure d'anneau naturelle. En bref,~$A$ est 
 l'ensemble des fonctions continues (à valeurs réelles) définies au voisinage de~$x$, deux fonctions
 appartenant à~$A$
 étant considérées comme égales si elles coïncident au voisinage de~$x$ ; on dit aussi que~$A$ est l'anneau
 des {\em germes de fonctions continues en~$x$.} L'évaluation en~$x$ induit un morphisme~$f\mapsto f(x)$
 de~$A$ dans~$\RR$. 
 
 Ce morphisme est surjectif, grâce aux fonctions constantes. Son noyau~$\got m$ est donc
 un idéal maximal de~$A$. Nous allons montrer que c'est le seul ; il suffit, par le critère donné
 ci-dessus, de vérifier que~$\got m$ est exactement l'ensemble des éléments
 non inversibles de~$A$. Soit~$f\in A\setminus \got m$. Choisissons un voisinage ouvert~$U$
 de~$x$ sur lequel~$f$ est définie. Comme~$f\notin \got m$, on a~$f(x)\neq 0$. Comme~$f$
 est continue, il existe un voisinage ouvert~$V$ de~$x$ dans~$U$ sur lequel~$f$ ne s'annule pas. 
 L'inverse~$g$ de~$f$ est alors une fonction continue sur~$V$, et l'on a~$fg=1$ dans l'anneau~$A$.
 Ainsi~$f$ est inversible, ce qui achève la preuve. 
 
 \deux{autre-exemple-geom}
 {\bf Remarque.}  On aurait pu tout aussi bien
 remplacer~$X$ par une variété différentiable (resp. analytique complexe)
 et~$A$ par l'anneau des germes
 de fonctions~${\mathscr C}^\infty$ (resp. holomorphe). 
 
 \deux{corps-resid}
 Si~$A$ est un anneau local d'idéal maximal~$\got m$, 
 le corps~$A/\got m$ sera appelé le
 {\em corps résiduel}
 de~$A$. 
 
 \subsection*{Localisation et idéaux premiers}
 
 \deux{rem-noyau}Soit~$A$ un anneau et soit~$S$ une partie multiplicative
 de~$A$. Soit~$f: A\to B$ un morphisme tel que~$f(S)\subset B\ti$ ; il induit un 
 morphisme~$g: S^{-1}A\to B$, donné par la formule
 $\frac a s\mapsto f(a)f(s)^{-1}$. 
 
 \medskip
 Un calcul explicite montre que le noyau de~$g$ ne dépend que de celui de~$f$, 
 et réciproquement. Plus précisément : 
 
 \medskip
 $\bullet$ ${\rm Ker}\;g=\{\frac a s\}_{a\in {\rm Ker}\;f}$ ; 
 
 $\bullet$ ${\rm Ker}\;f=\{a\;{\rm t.q.}\;\frac a 1 \in {\rm Ker}\;g\}$. 
 
\deux{ideaux-smoinsuna}
{\bf Idéaux premiers de~$S^{-1}A$}. 
Se donner un morphisme de~$S^{-1}A$ vers un corps~$K$ revient à se
 donner un morphisme de~$A$ vers~$K$ qui envoie chaque élément de~$S$ sur un 
 élément inversible
 de~$K$, c'est-à-dire sur un élément non nul de~$K$ ; cela revient donc à se donner un morphisme
 de~$A$ vers~$K$ dont le noyau ne rencontre pas~$S$. 
 
 Compte-tenu de la description des idéaux premiers en termes de morphismes vers un corps, 
 et de la 
 description explicite des noyaux donnée
 au~\ref{rem-noyau}
 ci-dessus, on en déduit que
 ~$$\got p\mapsto\got pS^{-1}A= \left\{\frac a s, a\in \got p, s\in S\right\}\;{\rm et}\;\got q\mapsto \left\{a\;{\rm t.q.}\;\frac a 1 \in \got q\right\}$$ établissent
 une bijection (visiblement croissante) entre l'ensemble des idéaux premiers de~$A$ ne rencontrant pas~$S$
 et l'ensemble des idéaux premiers de~$S^{-1}A$.

 \medskip
 On peut également formuler cette 
 dernière assertion comme suit : l'application
 ${\rm Spec}\; S^{-1}A\to {\rm Spec}\;A$ induite par~$A\to S^{-1}A$ 
 ({\em cf.}~\ref{spec-contrav})
 est injective, et a pour image
 l'ensemble des idéaux premiers de~$A$ qui ne rencontrent pas~$S$.

\deux{appl-nilpo}
{\bf Lemme.}
{\em Soit~$A$ un anneau et soit~$f\in A$. Les assertions suivantes
sont équivalentes : 

i) $f$ est nilpotent ; 

ii) pour tout corps~$K$ et tout morphisme~$\phi : A\to K$ 
on a~$\phi(f)=0$ ; 

iii) $f$ appartient à tous les idéaux premiers de~$A$.

\medskip
En d'autres termes, 
le nilradical de~$A$
est l'intersection
de tous les idéaux premiers de~$A$.}

\medskip
{\em Démonstration.} 
L'équivalence de~ii)
et~iii) résulte de la caractérisation des idéaux premiers
comme noyaux de morphismes vers un corps. L'implication~i)$\Rightarrow$ii)
est évidente. Supposons maintenant que~iii) est vraie,
et montrons~i). 

\medskip
L'ensemble des idéaux premiers de~$A_f$ est d'après~\ref{ideaux-smoinsuna}
en bijection avec l'ensemble des idéaux premiers de~$A$ qui ne rencontrent pas~$\{f^n\}_{n\in \NN}$, 
c'est-à-dire avec l'ensemble des idéaux premiers de~$A$ qui ne contiennent pas~$f$. Puisqu'on 
est sous l'hypothèse~iii), cet ensemble est vide. 

En conséquence,
$A_f$ n'a aucun idéal premier, ce qui signifie qu'il est nul. Il s'ensuit 
en vertu~\ref{a-f}
que~$f$ est nilpotent.~$\Box$

\deux{intro-ap}
{\bf Localisé d'un anneau en un idéal premier.}
Soit~$A$ un anneau, et soit~$\got p$ un idéal premier
de~$A$. Le sous-ensemble~$S=A\setminus \got p$
de~$A$ en est une partie multiplicative, 
et le localisé~$S^{-1}A$ est le plus souvent noté~$A_{\got p}$. On l'appelle
le
{\em localisé de~$A$
en l'idéal~$\got p$.}

\trois{ideaux-ap}
En vertu de~\ref{ideaux-smoinsuna}, l'ensemble des
idéaux premiers de~$A_{\got p}$ est en bijection croissante
avec l'ensemble des idéaux premiers de~$A$ ne rencontrant
pas~$S$, c'est-à-dire contenus dans~$\got p$. Or cet ensemble admet
un plus grand élément, à savoir~$\got p$. On en déduit que~$A_{\got p}$ 
possède un et un seul idéal maximal : celui qui correspond à~$\got p$. D'après
la description explicite de la bijection évoquée (voir~{\em loc. cit.}), cet idéal
est~$\got pA_{\got p}=\{\frac a s, a\in \got p, s\notin \got p\}\subset A_{\got p}$. 

\trois{corps-resid-ap}
Le morphisme
composé~$A\to A/\got p\hookrightarrow {\rm Frac}\;(A/\got p)$ envoie tout élément
de~$S$ sur un élément non nul, et partant inversible, de~${\rm Frac}\;(A/\got p)$. Il se
factorise donc de manière unique par~$A_{\got p}$. Le morphisme correspondant
de~$A_{\got p}$ vers~${\rm Frac}\;(A/\got p)$ est par construction donné par la formule
$\frac a s\mapsto \frac{\bar a}{\bar s}$ ; on voit immédiatement que son noyau est~$\got pA_{\got p}$. 
Il est par ailleurs surjectif, puisque tout élément de ${\rm Frac}\;(A/\got p)$
est de la forme~$\frac{\bar a}{\bar s}$ 
avec~$a\in A$ et~$s\notin \got p$. 

En conséquence,
le corps résiduel $A_{\got p}/\got p A_{\got p}$ s'identifie naturellement 
à~${\rm Frac}\;(A/\got p)$. 

\trois{ap-limind-af}
{\em Expression de~$A_{\got p}$ comme limite inductive
filtrante.}
La relation de divisibilité fait de~$A\setminus \got p$
un ensemble pré-ordonné filtrant (si~$f$ et~$g$ sont deux éléments
de~$A\setminus \got p$, leur produit est un multiple commun 
à~$f$ et~$g$ dans~$A\setminus \got p$).
Si~$f$ et~$g$ sont deux éléments de~$A$ tels que~$f|g$ alors~$f$
est inversible dans~$A_g$, et il existe donc un morphisme
de~$A$-algèbres~$A_f\to A_g$. 

\medskip
Il en résulte l'existence d'un diagramme
commutatif filtrant 

$$\sch D:=((A_f)_{f\in A\setminus \got p}, (A_f\to A_g)_{f|g})$$
dans la catégorie des~$A$-algèbres. Comme~$A\setminus \got p$
est évidemment égal à la la réunion de ses sous-parties multiplicatives de la forme~$\{f^n\}_{n\in \NN}$ pour~$f$
parcourant~$A\setminus \got p$, il résulte de~\ref{limind-smoinsun-a}
et de la remarque~\ref{limind-loc-multicat}
que~$A_{\got p}$ s'identifie à la limite inductive de~$\sch D$ dans la catégorie des~$A$-algèbres (et des
anneaux, et des $A$-modules).  

\deux{ex-ap}
{\bf Exemples.}

 \trois{azero}
 Supposons~$A$ intègre.

\medskip
Si~$\got p=\{0\}$, l'anneau~$A_{\got p}$ n'est autre par définition que le corps
 des fractions de~$A$. 
 
 \medskip
 En général, comme~$0\notin S$, la relation des produits en croix
 qui définit l'égalité dans~$A_{\got p}$ est la même que celle
 qui définit l'égalité dans~${\rm Frac}\;A$ ; ainsi, $A_{\got p}$ apparaît
 comme le {\em sous-anneau}
 de~${\rm Frac}\;A$ constitué des fractions qui admettent une écriture avec un dénominateur
 n'appartenant pas à~$\got p$. 
 
 \trois{zp}
 Soit~$p$ un nombre premier. Le localisé~$\ZZ_{(p)}$ est 
 d'après ce qui précède le sous-anneau de~$\QQ$
égal à
$$\left\{\frac a b, a\in \ZZ, b\in \ZZ\setminus p\ZZ\right\}.$$

\deux{rem-ap-germes}
{\bf Remarque.}
Le langage des schémas
 permet, pour tout anneau~$A$ et tout
 idéal premier~$\got p$ de~$A$, d'interpréter~$A_{\got p}$ comme un 
 anneau de germes de fonctions, analogue à ceux vus plus haut
 (exemple~\ref{ann-local-geodiff}
 et remarque~\ref{autre-exemple-geom}), et donc d'y penser en termes
 géométriques.

 \section{Endomorphismes d'un module et lemme de Nakayama}
\markboth{Algèbre commutative}{Lemme de Nakayama}

\deux{poly-annul-alin}
{\bf Proposition.}
{\em Soit~$A$ un anneau, 
soit~$n\in \NN$
et soit~$M$ un~$A$-module possédant une famille génératrice de cardinal~$n$. 
Soit~$I$ un idéal
de~$A$, et soit~$u$ un endomorphisme de~$M$ tel que~$u(M)\subset IM:=\{\sum a_i m_i, a_i\in I, m_i\in M\}$. 
Il existe alors une famille~$(a_1,\ldots, a_n)$ telle que~$a_j$
appartienne à~$ I^j$ pour tout~$j$, et telle que
$$u^n+a_1u^{n-1}+\ldots+a_1u+a_n{\rm Id}=0.$$}

\trois{rem-ideal-a}
{\bf Remarque.}
Lorsque~$I=A$, la condition~$u(M)\subset IM$ est automatiquement satisfaite. La proposition assure
donc entre autres que tout endomorphisme de~$M$ annule un polynôme unitaire de degré~$n$ à coefficients
dans~$A$. 

\trois{demo-poly)annul}
{\em Démonstration de la proposition~\ref{poly-annul-alin}}.
Choisissons une famille génératrice~$(e_1,\ldots, e_n)$
de~$M$. Comme~$u(M)=IM$, on a 
pour tout~$m\in M$ une égalité de la forme~$u(m)=\sum a_\ell m_\ell$ avec~$a_\ell \in I$
pour tout~$\ell$. En écrivant chacun des~$m_\ell$ comme combinaison 
linéaire des~$e_i$, on voit qu'on peut écrire~$u(m)$ comme combinaison linéaire des~$e_i$
{\em à coefficients dans~$I$.}

\medskip
En particulier, il existe une famille~$(a_{ij})$ d'éléments de~$I$ tels que
l'on ait~$u(e_i)=\sum_i a_{ij}e_i$ pour tout~$j$. Soit~$X$ la matrice~$(a_{ij})\in M_n(A)$ ; c'est en quelque 
sorte {\em une}
matrice de~$u$ dans la {\em famille génératrice}
$(e_1,\ldots, e_n)$. 

\medskip
Un calcul immédiat (le même que celui effectué en algèbre linéaire)
montre qu'on a pour tout~$(\lambda_1,\ldots, \lambda_n)\in A^n$ l'égalité
$u(\sum \lambda_ie_i)=\sum \mu_i e_i$ avec
$$\left(\begin{array}{c}\mu_1\\ .
\\
.
\\
. \\
\mu_n\end{array}\right)=
X\cdot\left(\begin{array}{c}\lambda_1\\ .
\\
.
\\
.
\\
\lambda_n\end{array}\right).$$

Par récurrence, on en déduit que l'on a pour tout entier~$r$ et pour tout
$(\lambda_1,\ldots, \lambda_n)\in A^n$ l'égalité
$u^r(\sum \lambda_ie_i)=\sum \nu_i e_i$ avec
$$\left(\begin{array}{c}\nu_1\\ .
\\
.
\\
. \\
\nu_n\end{array}\right)=
X^r\cdot\left(\begin{array}{c}\lambda_1\\ .
\\
.
\\
.
\\
\lambda_n\end{array}\right).$$

En vertu
du théorème de Cayley-Hamilton
\footnote{Vous ne l'avez peut-être rencontré que sur un corps, mais
sa validité dans ce cadre entraîne sa validité pour tout anneau. En effet, s'il est
vrai 
sur tout corps, il est vrai en particulier pour la matrice~$(X_{ij})\in M_n(\QQ(X_{ij}))$ ; il s'énonce
dans ce cas précis
comme une identité polynomiale à coefficients dans~$\ZZ$ en les~$(X_{ij})$.
Cette identité débouche
par spécialisation pour toute matrice~$(\alpha_{ij})$ à coefficients
dans un anneau quelconque sur la «même»
identité pour les~$\alpha_{ij}$... 
laquelle est précisément le théorème de Cayley-Hamilton pour~$(\alpha_{ij})$.},
on a~$\chi_X(X)=0$ et donc par ce qui précède~$\chi_X(u)=0$. Mais comme les~$a_{ij}$ appartiennent
à~$I$, le polynôme~$\chi_X$ est de la 
forme~$T^n+a_1T^{n-1}+\ldots+a_n$ avec~$a_j \in I^j$ pour tout~$j$, ce
qui achève la
démonstration.~$\Box$ 

\deux{nakayama}
{\bf Lemme de Nakayama.}
{\em Soit~$A$
un anneau, soit~$I$ un idéal de~$A$ et soit~$M$
un~$A$-module
de type fini. Les assertions suivantes sont équivalentes :

i) il existe un élément~$a$ de~$A$ congru à~$1$ modulo~$I$ et tel que~$aM=\{0\}$ ; 

ii) $M=IM$.}

\medskip
{\em Démonstration.}
Supposons que~i) soit vraie, écrivons~$a=1+b$ avec~$b\in I$. On a pour 
tout~$m\in M$ l'égalité~$(1+b)m=0$, et donc~$m=-bm$. Ainsi, $M=IM$.

\medskip
Supposons que~ii) soit vraie, et appliquons 
la proposition~\ref{poly-annul-alin}
avec~$u={\rm Id}_M$ (c'est possible puisque~$M$ est de type fini). 
Elle assure l'existence d'une famille~$(a_j)$ avec~$a_j\in I^j$ pour tout~$j$ telle que
$${\rm Id}_M^n+a_1{\rm Id}_M^{n-1}+\ldots+a_0{\rm Id}_M=0.$$ En l'appliquant
à un élément~$m$
de~$M$, on obtient
$(1+a_1+\ldots+a_n)m=0$ ; ainsi, i) est vraie avec~$a=a_1+\ldots+a_n$.~$\Box$

\medskip
Ce lemme est surtout utile en pratique {\em via}
son corollaire suivant -- qui n'est autre que la version originelle 
du lemme de Nakayama.

\deux{coro-nakayama}
{\bf Corollaire.}
{\em Soit~$A$ un anneau local d'unique idéal maximal~$\got m$,
et soit~$M$ un~$A$-module de type fini tel que~$M=\got m M$. 
On a alors~$M=\{0\}$.}

\medskip
{\em Démonstration.}
Le lemme de Nakayama assure qu'il existe un élément~$a$ congru à~$1$ modulo~$\got m$
tel que~$aM=\{0\}$. Étant non nul modulo~$\got m$, l'élément~$a$ appartient à~$A\ti$ ; 
il s'ensuit que~$M$ est trivial.~$\Box$ 

\deux{surj}
Donnons une conséquence très utile de ce corollaire ; on désigne toujours par~$A$
un anneau local d'idéal maximal~$\got m$.

\trois{surj-test-residuel}
Soit~$M$ un~$A$-module
de type fini, soit~$N$
un~$A$-module et soit~$f: N\to M$ une application~$A$-linéaire. Pour que~$f$ soit surjective,
il faut et il suffit que l'application induite~$N/\got m N\to M/\got m N$ le soit. 

\medskip
C'est en effet clairement nécessaire. Supposons maintenant que
la flèche $N/\got m N\to M/\got m N$ est surjective, et soit~$m\in M$. Par hypothèse, 
on peut écrire
$$m=f(n)+\sum a_i m_i$$
où~$n\in N$, où les~$m_i$
appartiennent à~$M$ et où les~$a_i$
appartiennent à~$\got m$. Il s'ensuit que~$m$ est égal
à~$\sum a_i m_i$ modulo~$f(N)$. En conséquence, le module
quotient~$M/f(N)$ vérifie l'égalité $$M/f(N)=\got m M/f(N).$$ Comme il est de type
fini puisque c'est déjà le cas de~$M$, il est nul d'après le corollaire~\ref{coro-nakayama}
ci-dessus. Ainsi, $f(N)=M$ et~$f$ est surjective. 

\trois{gen-test-residuel}
Soit~$(e_i)_{i\in i}$ une famille d'éléments de~$M$. Elle engendre le module~$M$
si et seulement si les~$\overline{e_i}$ engendrent le~$A/\got m$-espace vectoriel~$M/\got m M$. 

C'est en effet une simple application du~\ref{surj-test-residuel}
ci-dessus, au cas où~$N$ est le module~$A^{(I)}$ formé des familles~$(a_i)_{i\in I}$
de~$A^I$ dont presque tous les éléments sont nuls, et où~$f$ est l'application
$(a_i)\mapsto \sum a_i e_i$. 

\medskip
Nous allons
maintenant donner
une application astucieuse et très frappante
du lemme de Nakayama. 

\deux{prop-surj-inj}
{\bf Corollaire.}
{\em Soit~$A$ un anneau, soit~$M$ un~$A$-module de type fini et soit~$u$ un endomorphisme
surjectif de~$M$. L'endomorphisme~$u$ est alors bijectif et~$u^{-1}$ est un polynôme en~$u$.}

\medskip
{\em Démonstration.} La loi
externe
$$(P,m)\mapsto P(u)(m)$$
définit sur le groupe abélien~$(M,+)$
une structure de~$A[X]$-module
qui prolonge
celle de~$A$-module, et
la multiplication par~$X$ est égale à l'endomorphisme~$u$. 
Par hypothèse, ~$M$ est de type fini comme~$A$-module ; 
il l'est {\em a fortiori} comme~$A[X]$-module. 

La surjectivité de~$u$ signifie
que pour tout~$m\in M$, 
il existe~$n\in M$ tel que~$Xn=m$. En conséquence,~$M=(X)M$. 

Le lemme de Nakayama assure alors qu'il existe un polynôme~$P$
congru à~$1$ modulo~$X$ tel que~$PM=0$. Écrivons~$P=1+XQ$, 
avec~$Q\in A[X]$. Soit~$m\in M$. 
On a $Pm=0$, soit~$P(u)(m)=0$,
soit encore~$({\rm Id}+uQ(u))(m)=0$. Ceci valant
pour tout~$m$, il vient~${\rm Id}=u(-Q(u))$. Comme deux polynômes
en~$u$ commutent, on aussi~${\rm Id}=(-Q(u))u$. Ainsi,~$u$ est bijectif
et~$u^{-1}=-Q(u)$.~$\Box$

%
%

\section{Le produit tensoriel : cas de deux modules}
\markboth{Algèbre commutative}{Produit tensoriel de deux modules}

{\em On fixe pour toute cette section un anneau~$A$.}

\subsection*{Définition, exemples et premières propriétés}

\deux{intro-tens}
Soient~$M$ et~$N$ deux~$A$-modules. Avant de
définir rigoureusement le produit tensoriel
de~$M$
et~$N$, expliquons intuitivement le
but de sa construction. On cherche à fabriquer la loi bilinéaire {\em la plus générale
possible} de source~$M\times N$, c'est-à-dire à donner un sens au produit d'un élément
de~$M$ par un élément de~$N$, en ne lui imposant rien d'autre que la bilinéarité. 

Comme à chaque fois que l'on cherche à construire
un objet obéissant à une 
liste limitative de contraintes, 
la définition rigoureuse de l'objet en question s'exprime au moyen d'une propriété
universelle ou, si l'on préfère, du foncteur qu'il représente. 

\medskip
Pour tout~$A$-module~$P$, 
on note~${\rm Bil}_A(M\times N,P)$ l'ensemble
des applications bilinéaires de~$M\times N$
vers~$P$.

\deux{def-tens}
{\bf Définition -- proposition.}
{\em Le foncteur covariant~$P\mapsto {\rm Bil}_A(M\times N,P)$
de~$\amod$
dans~$\ens$  est représentable, 
et son représentant est noté
$$(M\otimes_A N, (m,n)\mapsto m\otimes n).$$
On dit que~$M\otimes_A N$ est le {\em produit tensoriel} de~$M$ et~$N$
au-dessus de~$A$. }

\medskip
{\em Démonstration.}
On part
d'un~$A$-module libre~$L$ 
ayant une base~$(e_{m,n})$ indexée par les éléments de~$M\times N$, 
et l'on note~$L_0$ le sous-module de
$L$ engendré par les 

$$e_{m,n+\lambda n'}-e_{m,n}-\lambda e_{m,n'}\;\;{\rm et}\;
e_{m+\lambda m',n}-e_{m,n}-\lambda e_{m',n}$$
pour~$(m,m',n,n',\lambda)$ parcourant~$M^2\times N^2\times A$. 
On pose alors

$$M\otimes_AN=L/L_0, \;{\rm et}\;m\otimes n=\overline{e_{m,n}}\;
{\rm pour}\;{\rm tout}\;(m,n)\in M\times N.$$

Notons que par construction,
les~$m\otimes n$
engendrent le~$A$-module
$M\otimes_AN$. 

\medskip
Montrons que
$(M\otimes_A N, (m,n)\mapsto m\otimes n)$
représente~$F$. Soit~$P$ un~$A$-module 
et soit~$ b \in {\rm Bil}_A(M\times N,P)$. 
une application bilinéaire. Il s'agit de 
prouver
qu'il existe une et une seule application
linéaire~$\lambda : M\otimes_AN\to P$ 
telle que~$\lambda(m\otimes n)=b(m,n)$ 
pour tout~$(m,n)$. 

\medskip
{\em Unicité.}
Elle provient simplement du fait que la famille~$(m\otimes n)$ 
est génératrice. 

\medskip
{\em Existence.}
Soit~$\phi$ l'unique application~$A$-linéaire
de~$L$ dans~$P$ envoyant~$e_{m,n}$
sur~$b(m,n)$ pour tout~$(m,n)$. Comme~$b$ est bilinéaire, 
l'application~$\phi$ s'annule sur les éléments de~$L_0$ ; elle induit
donc par passage au quotient une application linéaire
$\lambda : M\otimes_A N\to P$, et l'on a pour tout~$(m,n)\in M\times N$
les égalités 
$$\lambda(m\otimes n)=\lambda(\overline{e_{m,n}})=\phi(e_{m,n})=b(m,n),$$
ce qui achève la démonstration.~$\Box$

\deux{comment-prod-tens}
{\bf Commentaires et premières propriétés.}

\trois{const-prod-tens}
La construction du produit tensoriel
n'est guère subtile ; elle consiste à imposer
par décret les propriétés requises. 
Elle n'est en pratique {\em jamais}
utilisée, et il faut à tout prix
éviter de penser au produit tensoriel comme au
quotient d'un module libre monstrueux. 

Il y a toutefois une chose à en retenir : le fait
que~$M\otimes_AN$ est engendré
 (comme~$A$-module,
 ou même comme groupe abélien puisque l'on a
  pour tout~$(a,m,n)$
  les égalités~$a(m\otimes n)=(am)\otimes n$)
 par les éléments de la forme~$m\otimes n$, qu'on appelle les {\em tenseurs purs}.

 \trois{prop-petit-tenseur}
 L'application bilinéaire universelle~$(m,n)\mapsto m\otimes n$ a tendance à coder
 les propriétés (de nature linéaire) vérifiées par {\em toutes}
 les applications
 bilinéaires de source~$M\times N$. Illustrons cette pétition
 de principe par un exemple : nous allons montrer que
 ~$\sum m_i\otimes n_i=0$
 si et seulement si pour
 {\em tout}
 $A$-module~$P$ et {\em toute} application 
 bilinéaire~$b:M\times N\to P$, on a~$\sum b(m_i,n_i)=0$. 
 
 \medskip
 {\em Supposons que pour tout $A$-module~$P$ et
 toute
 application 
 bilinéaire~$b:M\times N\to P$, on ait~$\sum b(m_i,n_i)=0$}. C'est en particulier
 le cas pour l'application~$\otimes$, et il vient~$\sum m_i\otimes n_i=0$. 
 
 \medskip
 {\em Supposons que~$\sum m_i\otimes n_i=0$.}
 Soit~$P$ un~$A$-module
 et soit~$b \in {\rm Bil}_A(M\times N,P)$. L'application~$b$
 induit une application~$A$-linéaire ~$\lambda : M\otimes_A N\to P$
 telle que~$\lambda(m\otimes n)=b(m,n)$
 pour tout~$n$. On a donc
 
 $$\sum b(m_i,n_i)=\sum \lambda(m_i\otimes n_i)=\lambda(\sum m_i \otimes n_i)
 =0,$$
 ce qu'il fallait démontrer.
 
 \trois{nullite-tens}
 {\bf Exercice.}
 Dans le même esprit, montrez que~$M\otimes_AN$ est nul si et seulement
 si toute application bilinéaire de source~$M\times N$ est nulle. 
 
 \deux{ex-prod-tens}
 {\bf Premiers exemples.}
 
 \trois{m-tenseur-0}
 {\bf Exemple trivial}.
 Si~$M$ est un~$A$-module quelconque alors
 $$\{0\}\otimes M=M\otimes\{0\}=\{0\}\;:$$ 
 cela vient du fait que le produit tensoriel est engendré par un tenseur pur, et qu'un tenseur
 pur dont l'un des deux facteurs est nul est lui-même nul. 
 
 \trois{prod-tens-sym}
 {\bf Symétrie du produit tensoriel.}
 Soient~$M$ et~$N$ 
 deux~$A$-modules. L'application
 de~$M\times N$ dans~$N\otimes_AM$ qui envoie~$(m,n)$
 sur~$n\otimes m$ est bilinéaire, et induit donc une application
 $A$-linéaire~$M\otimes_AN\to N\otimes_AM$ qui
 envoie~$m\otimes n$ sur~$n\otimes m$ pour tout~$(n,m)$. 
 
 On a de même une application~$A$-linéaire~$v : N\otimes_AM\to M\otimes_AN$
 qui envoie
 $n\otimes m$ sur~$m\otimes n$. Il est immédiat
 que~$u\circ v={\rm Id}_{N\otimes_AM}$
 et~$v\circ u={\rm Id}_{M\otimes_AN}$ 
 (le vérifier sur les tenseurs purs). Ainsi, $u$
 et~$v$ sont deux isomorphismes réciproques l'un de l'autre. 
 
 \trois{ex-prod-tens-sym}
 {\bf Construction
 fonctorielle de~$u$.} Soit~$P$ un~$A$-module. L'application
 $$b\mapsto [(n,m)\mapsto b(m,n)]$$ induit un isomorphisme
 fonctoriel en~$P$ entre l'ensemble des applications
 bilinéaires de~$M\times N$ vers~$P$ et celui des applications
 bilinéaires de~$N\times M$ vers~$P$. Par le lemme de Yoneda, 
 cet isomorphisme est induit par une bijection~$A$-linéaire
 de~$M\otimes_AN$ vers~$N\otimes_AM $ ; on vérifie
 immédiatement que cette bijection n'est autre que~$u$. 
  
 \trois{otimes-pas-integre}
On prendra garde qu'en général, le
produit tensoriel de deux modules non nuls peut
très bien être nul. Nous allons montrer par
exemple que
$$(\ZZ/2\ZZ)\otimes_{\ZZ}\ZZ/3\ZZ=0.$$ Pour
 cela, il suffit de montrer que~$a\otimes b=0$ pour tout~$a\in \ZZ/2\ZZ$
 et tout~$b\in \ZZ/3\ZZ$. Donnons-nous donc un tel couple~$(a,b)$.
 On a~$$a\otimes b=(3-2)a\otimes b=3a\otimes b-2a\otimes b=a\otimes 3b-2a\otimes b=0,$$
 puisque~$2a=0$ et~$3b=0$. 
 
 \medskip
 Plus généralement,$(\ZZ/p\ZZ)\otimes_{\ZZ}\ZZ/q\ZZ=0$ dès que~$p$ et~$q$
 sont premiers entre eux : on raisonne comme ci-dessus, 
 en remplaçant l'égalité~$3-2=1$ par une relation de Bezout entre~$p$
 et~$q$. 

\deux{prod-tens-libre1}
{\bf Produit tensoriel par un module libre
de rang~$1$.}
Soit~$N$ un~$A$-module, et soit~$M$
un~$A$-module libre de rang~$1$. Soit~$e$
une base de~$M$. 

\trois{m-tenseur-ae}
Soit~$\phi$
l'application
linéaire
de~$N$ vers~$M\otimes_A N$ 
donnée par la formule~$n\mapsto e\otimes n$. 
Nous allons montrer que c'est un isomorphisme en exhibant 
 réciproque. 

Soit~$b$ l'application de~$M\times N$ dans~$M$ 
qui envoie un couple~$(ae, n)$ sur~$an$ (elle est bien définie
car~$e$ est une base de~$M$) ; elle est bilinéaire, donc induit
une application~$A$-linéaire~$\psi$ de~$M\otimes_AN$ vers~$N$
qui envoie~$ae\otimes n$ sur~$an$ pour tout~$(a,n)$. 

On vérifie immédiatement par leur effet sur les tenseurs purs que~$\phi$
et~$\psi$ sont réciproques l'une de l'autre. 

\trois{m-tenseur-a}
{\bf Un cas particulier important.}
On déduit de ce qui précède que
pour tout~$A$-module~$N$, l'application
linéaire
$n\mapsto 1\otimes n$ induit
un isomorphisme~$N\simeq A\otimes_A N$. 

\trois{m-tenseur-ae-fonct}
{\bf Construction fonctorielle de~$\phi$.} 
Soit~$P$ un~$A$-module. L'application
$$b \mapsto [n\mapsto b(e,n)]$$ définit une bijection
fonctorielle en~$P$ entre~${\rm Bil}_A(M\times N,P)$
et~$\hom_A(N,P)$,
de réciproque
$$\lambda\mapsto [(ae,n)\mapsto a\lambda(n)].$$
Par le lemme de Yoneda, cette
collection de bijections 
est induite par une application~$A$-linéaire de~$N$ vers~$N\otimes_AM$, dont
on vérifie qu'elle n'est autre que~$\phi$. 

\deux{ae-tenseur-af}
{\bf Produit tensoriel de deux modules libres de rang~$1$}. 
Soient maintenant~$M$ et~$N$ deux~$A$-modules libres de rang~$1$. 
Donnons-nous une base~$e$ de~$M$ et une base~$f$ de~$N$.

\trois{base-ae-tenseur-af}
Il résulte de~\ref{m-tenseur-ae}
que
la formule~$n\mapsto n\otimes e$ définit
un isomorphisme~$N\simeq M\otimes_A N$.
Comme~$a\mapsto af$ définit un isomorphisme~$A\simeq M$, 
on voit que~$a\mapsto a e\otimes f$
définit un isomorphisme~$\iota : A\simeq M\otimes_A N$. En d'autres
termes~$M\otimes_AN$ est libre de rang~$1$ de base~$e\otimes f$. 

\trois{ae-tenseur-af-fonct}
{\bf Construction fonctorielle de~$\iota$.}
Soit~$P$ un~$A$-module. L'application
$$b\mapsto b(e,f)$$
définit une bijection fonctorielle en~$P$
entre~${\rm Bil}_A(M\times N,P)$ et~$P$,
de réciproque
$$p\mapsto [(ae, bf)\mapsto abp].$$ 
Comme par ailleurs~$p\mapsto [a\mapsto ap]$
définit une bijection fonctorielle en~$P$ entre~$P$
et~$\hom_A(A,P)$
(de réciproque~$\lambda\mapsto \lambda(1)$),
on obtient par composition 
une bijection fonctorielle en~$P$ entre
${\rm Bil}_A(M\times N,P)$ et~$\hom_A(A,P)$.

Par le lemme de Yoneda, 
cette collection de bijections est induite par une bijection~$A$-linéaire
de~$A$ vers~$M\otimes_AN$, dont on vérifie qu'elle n'est autre que~$\iota$. 

\deux{fonct-prod-tens}
{\bf Fonctorialité du produit tensoriel en ses deux arguments.}
Soient~$M, M', N$ 
et~$N'$ quatre~$A$-modules, 
et soient~$f : M\to M'$ et~$g : N\to N'$
deux applications~$A$-linéaires. 

\trois{f-tenseur-g}
L'application
de~$M\times N$ vers~$M'\otimes_AN'$ donnée
par la formule~$(m,n)\mapsto (f(m)\otimes g(n))$
est bilinéaire. Elle induit donc une application~$A$-linéaire
$f\otimes g : M\otimes_AN\to M'\otimes_AN'$,
telle que~$(f\otimes g)(m\otimes n)=f(m)\otimes g(n)$
pour tout~$(m,n)$. On vérifie que~$(f,g)\mapsto f\otimes g$
est elle-même une application bilinéaire
de~$\hom_A(M,M')\times \hom_A(N,N')$
vers~$\hom_A(M\otimes_AN,M'\otimes_AN')$. 

\trois{f-tenseur-g-fonct}
{\bf Description fonctorielle de~$f\otimes g$.}
Soit~$P$ un~$A$-module. La formule
$$b\mapsto [(m,n)\mapsto b(f(m), g(n))]$$
définit une application 
de~${\rm Bil}_A(M'\times N',P)$
vers~${\rm Bil}_A(M\times N,P)$
qui est fonctorielle en~$P$. 

Par le lemme de Yoneda, cette collection
d'applications est induite par une application~$A$-linéaire
de~$M\otimes_AN$
vers~$M'\otimes_AN'$, dont on vérifie qu'elle n'est autre que~$f\otimes g$.

\deux{adj-prod-tens}
{\bf Propriétés d'adjonction.}
Soient~$L, M$ et~$N$ trois~$A$-modules, et soit~$f \colon N\to \hom(M,L)$ 
une application~$A$-linéaire. L'application de~$M\times N$ dans~$L$
qui envoie~$(m,n)$ sur~$f(n)(m)$ est bilinéaire, et induit donc
une application linéaire~$M\otimes_A N\to L$, qui dépend manifestement
linéairement de~$f$ ; on a donc construit une application~$A$-linéaire
$$p \colon \hom(N, \hom(M,L))\to \hom (M\otimes_A N,L).$$
On définit par ailleurs une application~$A$-linéaire 
$$q\colon \hom (M\otimes_A N,L)\to \hom(N, \hom(M,L))$$
par la formule~$q(g)=n\mapsto [m\mapsto g(m\otimes n)]$,
et l'on vérifie immédiatement que~$p$ et~$q$ sont des bijections
et réciproques l'une de l'autre, fonctorielles en~$M, L$ et~$N$. 
En particulier, si l'on considère~$M$ comme fixé, 
on dispose d'un isomorphisme
$$\hom(N, \hom(M,L))\simeq \hom (M\otimes_A N,L)$$
qui est fonctoriel en~$N$ et~$L$ ; en conséquence, 
$L\mapsto \hom(M,L)$ est adjoint à droite à~$N\mapsto M\otimes_A N$. 

\medskip
Comme le foncteur~$N\mapsto M\otimes_A N$ admet un adjoint à droite, {\em il
commute aux limites inductives} (\ref{lim-et-adj}).

\deux{pro-tens-somme-dir}
{\bf Produit tensoriel et somme directe.}
Soit~$M$ un~$A$-module et soit~$(N_i)$
une famille de~$A$-modules. Pour tout~$i$, 
on note~$u_i$ l'injection naturelle de~$N_i$ 
dans~$\bigoplus N_i$. 

\trois{tens-somme-dir}
La famille des~${\rm Id}_M\otimes u_i : M\otimes_AN_i\to M\otimes_A(\bigoplus N_i)$
induit un morphisme~$\chi : (\bigoplus M\otimes_A N_i)\to M\otimes_A (\bigoplus N_i)$. Nous allons
montrer que c'est un isomorphisme. Comme la somme directe est un cas particulier de limite inductive, 
on peut pour ce faire se contenter d'invoquer la commutation de~$M\otimes$ aux limites inductives
({\em cf.}~\ref{adj-prod-tens}
ci-dessus). Mais nous allons également donner deux preuves directes, la première consistant
à exhiber la bijection réciproque par une formule, et la seconde à utiliser le lemme de Yoneda. 

\trois{tens-somme-dir-down}
{\em Construction de la bijection réciproque de~$\chi$ par une formule.}
L'application
de~$M\times(\bigoplus N_i)$ dans~$\bigoplus M\otimes_A N_i$
donnée par la formule
$(m, (n_i))_i\mapsto (m\otimes n_i)_i$ est bilinéaire. 
Elle induit
dès lors
une application linéaire~$\rho$
de~$M\otimes_A(\bigoplus N_i)$
vers~$\bigoplus M\otimes_A  N_i$. On vérifie immédiatement
que~$\chi$
et~$\rho$ sont inverses l'une de l'autre. 

\trois{tens-somme-dir-fonct}
{\em Preuve de la bijectivité de~$\chi$
{\em via}
le lemme de Yoneda.}
La somme directe~$\bigoplus M\otimes_A N_i$
représente le foncteur~$\prod_i h_{M\otimes_A N_i}$, 
c'est-à-dire encore le foncteur qui envoie un~$A$-module~$P$
sur~$\prod {\rm Bil}_A(M\times N_i, P)$. 

\medskip
Soit~$P$ un~$A$-module. La formule
$$b\mapsto (b_{|M\times N_i})_i$$ définit une bijection
fonctorielle en~$P$ de~${\rm Bil}_A(M\times(\bigoplus N_i),P)$
vers
le produit~$\prod {\rm Bil}_A(M\times N_i,P)$, de réciproque
$$(b_i)\mapsto [(m,(n_i)_i)\mapsto \sum b_i(m,n_i)].$$

Par le lemme de Yoneda, cette collection
de bijections est induite par une bijection~$A$-linéaire
de~$\bigoplus (M\otimes_AN_i)$
vers~$M\otimes_A(\bigoplus N_i)$, 
dont on vérifie qu'elle n'est autre que~$\chi$. 

\deux{prod-tens-mod-libres}
{\bf Produit tensoriel de modules libres.}
Soient~$M$ et
$N$
deux~$A$-modules libres, 
de bases respectives
$(e_i)$
et~$(f_j)$ ; on a
les égalités $M=\bigoplus A\cdot e_i$
et~$N=\bigoplus A \cdot f_j$. 

On déduit alors
de~\ref{tens-somme-dir}
que~$M\otimes_A N\simeq \bigoplus_j M\otimes_A (A\cdot f_j)$. 
En réappliquant~\ref{tens-somme-dir}
à chacun des sommandes (et en utilisant la symétrie du produit tensoriel, 
{\em cf.}
\ref{prod-tens-sym}), 
il vient~$M\otimes_A N\simeq \bigoplus_{i,j}(A\cdot e_i)\otimes_A(A\cdot f_j)$. 

\medskip
Mais en vertu de~\ref{base-ae-tenseur-af}, le module~$(A\cdot e_i)\otimes_A(A\cdot f_j)$
est pour tout~$(i,j)$ libre de base~$e_i\otimes f_j$. Il s'ensuit que~$M\otimes_A N$
est libre de base~$(e_i\otimes f_j)_{i,j}$. 

\medskip
\deux{conclu-tens-lib}
Ainsi, le produit tensoriel de deux~$A$-modules libres~$M$
et~$N$ est libre, et si~$A$ est non nul le rang de~$M\otimes_AN$
est égal au produit du rang de~$M$
et du rang de~$N$. 

Il s'ensuit
que
{\em si~$A$ est non nul, le produit tensoriel
de deux~$A$-modules libres non nuls est toujours non nul}. Notez 
un cas particulier important : le produit tensoriel de deux
espaces vectoriels non nuls sur un corps~$k$ est non nul. Nous
aurons plusieurs fois l'occasion de l'utiliser. 

\deux{prod-tens-ensfini}
{\bf Produit tensoriel d'une famille de modules.}
Ce qu'on a vu pour deux modules se généralise sans peine à une famille quelconque~$(M_i)_{i\in I}$
de~$A$-modules : le foncteur qui envoie un~$A$-module~$P$ sur l'ensemble des applications
multilinéaires de~$\prod_{i\in I}M_i$ vers~$P$ est représentable par un~$A$-module
noté~$\bigotimes_{i\in I}M_i$, fourni avec une application multilinéaire universelle
$(x_i)\mapsto \otimes_{i\in I}x_i$ de~$\prod M_i$ vers~$\bigotimes_{i\in I}M_i$. 
La construction est analogue à celle donnée à la proposition~\ref{def-tens} : on part d'un 
module libre de base indexée par~$\prod M_i$ que l'on quotiente par les relations
exigées. 

\medskip
Si~$I$ est réunion disjointe de deux ensembles~$I'$ et~$I''$ on a
$$\bigotimes_{i\in I}M_i\simeq \left(\bigotimes_{i\in I'}M_i\right)\otimes\left(\bigotimes_{i\in I''}M_i\right).$$
Notez que le produit tensoriel {\em vide}
de modules a un sens, et est égal au~$A$-module~$A$ : le lecteur est invité
à vérifier que c'est une conséquence de sa définition comme représentant d'un foncteur, 
et que c'est par ailleurs bien ce que donne la construction esquissée ci-dessus. 

\subsection*{Propriétés d'exactitude}

\deux{rapp-suite-exacte}
{\bf Brefs rappels sur les suites exactes de~$A$-modules.}
Soient~$n^-$ et~$n^+$ deux éléments de~$\ZZ\cup\{-\infty,+\infty\}$
et soit
$$S=\ldots \to M_i\to M_{i+1}\to M_{i+2}\to \ldots$$
une suite de morphismes de~$A$-modules, où~$i$ 
parcourt l'ensemble~$I$ des entiers relatifs compris entre~$n^-$ et~$n^+$. 

\medskip
Soit~$i$ un élément de~$I$ tel que~$i-1$ et~$i+1$ appartiennent à~$I$. On dit
que la suite~$S$
est {\em exacte en~$M_i$}
si le noyau de~$M_i\to M_{i+1}$
est égal à l'image de~$M_{i-1}\to M_i$. On dit que~$S$ est {\em exacte}
si elle est exacte en~$M_i$ pour tout~$i$ tel que~$i-1$
et~$i+1$ appartiennent à~$I$ (les indices extrêmes, s'ils existent,
ne comptent donc pas). 

Il résulte de la définition que dans une suite
exacte, la composée de deux flèches successives est toujours nulle. 

\medskip
Donnons quelques exemples. 

\trois{ex-exa-d}
La suite
$$\diagram M'\rto^f&M\rto^g&M''\rto&0\enddiagram $$
est exacte si et seulement si~$g$ est surjective et~${\rm Ker}\;g={\rm Im}\;f$. 

\medskip
Le lecteur est invité à vérifier que cela peut se reformuler en termes plus
catégoriques de la façon suivante : {\em le triplet
$$(g\colon M\to M'', \;\;0\colon \{0\}\to M'', \;\;g\circ f \colon M'\to M'')$$
définit un morphisme du diagramme
$$\sch D=\xymatrix{
\{0\}\ar[rd]^0&
\\
&M\\
{M'}\ar[ur]_f&}$$
vers~$M''$, et identifie~$M''$ à la limite inductive de~$\sch D$.}

\trois{ex-exa-g}
La suite
$$\diagram 0\rto &M'\rto^f&M\rto^g&M''\enddiagram$$
est exacte si et seulement si~$f$ est injective
et~${\rm Ker}\;g={\rm Im}\;f$. 

\medskip
Le lecteur est invité à vérifier que cela peut se reformuler en termes plus
catégoriques de la façon suivante : {\em le triplet
$$(f\colon M'\to M, \;\;0\colon M'\to \{0\}, \;\;g\circ f \colon M'\to M'')$$
définit un morphisme de~$M'$ vers le diagramme
$$\sch D=
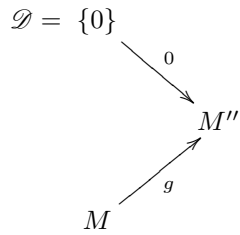
\xymatrix{
\{0\}\ar[rd]^0&
\\
&{M''}\\
M\ar[ur]_g&}$$
et identifie~$M'$ à la limite projective de~$\sch D$.}

\trois{ex-exa-dg}
La suite
$$\diagram 0\rto &M'\rto^f&M\rto^g&M''\rto&0\enddiagram$$
est exacte si et seulement si~$f$ est injective, $g$ est
surjective
et~${\rm Ker}\;g={\rm Im}\;f$. 

\deux{fonct-exact}
Soit~$B$ un anneau, et soit~$F$ un foncteur covariant
de~$\amod$ vers~$\bmod$. 

\trois{def-fonct-exact}
On dit que~$F$
est {\em exact à gauche} (resp. {\em exact 
à droite}, resp. {\em exact})
si et seulement si il satisfait les conditions suivantes : 

\medskip
$\bullet$ $F$ est {\em additif},
c'est-à-dire que pour tout couple~$(M,N)$
de~$A$-modules l'application
naturelle
$\hom_A(M,N)\to \hom_B(F(M),F(N))$
est un morphisme de groupes (à titre d'exercice, 
vous pouvez vérifier
que cela entraîne
la commutation de~$F$ aux sommes
directes finies) ; 

$\bullet$ pour toute suite exacte~$S$ 
de la forme
$$0\to M'\to M\to M''\;({\rm resp.} \;M'\to M\to M''\to 0, \;{\rm resp.}\;0\to M'\to M\to M''\to 0),$$
la suite~$F(S)$ est exacte.

\trois{equiv-font-exact}
{\bf Lemme.}
{\em Les propriétés suivantes sont équivalentes : 

\medskip
i) $F$ est exact ; 

ii) $F$ transforme toute suite exacte en une suite exacte ;

iii) $F$ est exact à gauche et transforme les surjections en surjections ; 

iv) $F$ est exact à droite et transforme les injections en injections.}

\medskip
{\em Démonstration.}
Il est clair que~iii)$\Rightarrow$i), que~iv)$\Rightarrow$i), et que~ii) entraîne 
iii) et~iv). Il reste à montrer
que~i) entraîne~ii). Il suffit par définition de s'assurer que si~$F$ est exact,
il transforme toute suite exacte~$M'\to M\to M''$ en une suite exacte. Mais cela résulte du fait 
que la suite exacte~$M'\to M\to M''$ se dévisse en suite exactes 
$$0\to K\to M'\to P\to 0, 0\to P\to M\to Q\to 0\;\;\text{et}\;0\to Q\to M''\to R\to 0$$
(prendre pour~$K$ le noyau de~$M'\to M$, pour~$P$ son image, pour~$Q$ le quotient de~$M$ par~$P$,
et pour~$R$ le conoyau de~$M\to M''$). Chacune d'elle reste par définition exacte quand on applique~$F$, 
et en recollant les suites obtenues on voit que~$F(M')\to F(M)\to F(M'')$ est exacte.~$\Box$

\deux{prod-tens-exd}
{\bf Proposition.}
{\em Soit~$M$ un~$A$-module. Le foncteur~$N\mapsto M\otimes_A N$ est exact à droite.}

\medskip
{\em Démonstration.} 
Que~$N\mapsto M\otimes_AN$
soit un foncteur additif résulte de~\ref{f-tenseur-g}. Soit maintenant
$$\diagram N\rto^f&L\rto^g&P\rto&0\enddiagram$$
une suite exacte. Nous allons montrer que
$$\diagram M\otimes_AN\rrto^{{\rm Id}_M\otimes f}&&
M\otimes_AL\rrto^{{\rm Id}_M\otimes g}&&M\otimes_AP\rto&0\enddiagram$$
est exacte. Nous allons commencer par une preuve conceptuelle extrêmement concise
mais peu explicite, puis donner une démonstration directe «à la main». 

\trois{prod-tens-exd-direct}
{\em La preuve conceptuelle.}
Le produit tensoriel commutant aux limites inductives (\ref{adj-prod-tens}), 
l'assertion requise découle aussitôt de la caractérisation catégorique
de l'exactitude à droite d'une suite (\ref{ex-exa-d}).  

\trois{prod-tens-surj}
{\em Preuve «à la main»
de la surjectivité de~${\rm Id}_M\otimes g$.}
Soit~$(m,p)$
appartenant à~$M\times P$. Comme~$g$
est surjective, l'élément~$p$ de~$P$
a un antécédent~$\ell$
dans~$L$ par~$g$. 

On a alors
$({\rm Id}_M\otimes g)(m\otimes \ell)=m\otimes g(\ell)=m\otimes p$. 
Ainsi, l'image de~${\rm Id}_M\otimes g$
contient tous les tenseurs purs ; en conséquence, elle est égale
à~$M\otimes_AP$ tout entier.

\trois{prod-tens-milieu}
{\em Preuve «à la main»
de l'égalité~${\rm Ker}({\rm Id}_M\otimes g)={\rm Im}({\rm Id}_M\otimes f)$.}
On a~$g\circ f=0$ ; il s'ensuit que~$({\rm Id_M}\otimes g)\circ ({\rm Id}_M\otimes f)=0$ ; 
autrement dit, ${\rm Im}({\rm Id}_M\otimes f)\subset {\rm Ker}({\rm Id}_M\otimes g)$. 
L'application~${\rm Id}_M\otimes g$ induit donc une surjection

$$M\otimes_AL/({\rm Im}({\rm Id}_M\otimes f))\to M\otimes_A P.$$
Pour montrer
que~${\rm Im}({\rm Id}_M\otimes f)$
est égale à~${\rm Ker}({\rm Id}_M\otimes g)$, il suffit de montrer que
cette surjection est un isomorphisme ; nous allons pour ce faire
exhiber sa réciproque. 

\medskip
 Soit~$m\in M$, soit~$p\in P$ et soit~$\ell$
 un antécédent de~$p$
 par~$g$. La classe de~$m\otimes \ell$
 modulo~${\rm Im}({\rm Id}_M\otimes f)$
 ne dépend alors pas du choix de~$\ell$. En effet, si~$\ell'$
 est un (autre)
 antécédent de~$p$ alors~$\ell-\ell'\in {\rm Ker}\;g={\rm Im}\;f$. En conséquence,
 il existe~$n\in N$ tel que~$\ell-\ell'=f(n)$,
 et l'on a donc
 $$m\otimes \ell-m\otimes \ell'=m\otimes (\ell-\ell')=m\otimes f(n)=({\rm Id}_M\otimes f)(m\otimes n),$$
 d'où l'assertion. 
 
 \medskip
 L'application de~$M\times P$ 
 vers~$M\otimes_AL/({\rm Im}({\rm Id}_M\otimes f))$
 qui envoie~$(m,p)$ sur la classe de~$m\otimes \ell$
 pour n'importe quel antécédent~$\ell$
 de~$p$ est donc bien définie. On voit immédiatement 
 qu'elle est bilinéaire ; elle induit donc une application
 $A$-linéaire~$\sigma : M\otimes_AP\to M\otimes_AL/({\rm Im}({\rm Id}_M\otimes f))$. 
 On vérifie sur les tenseurs purs que~$\sigma$ est bien un inverse à gauche et à droite
 de la surjection
 $$M\otimes_AL/({\rm Im}({\rm Id}_M\otimes f))\to M\otimes_A P$$
 induite par~${\rm Id}_M\otimes g$, 
 ce qui achève la démonstration.~$\Box$ 
 
 \deux{prod-tens-pas-exact}
 {\bf Remarque.}
  Le foncteur~$N\mapsto M\otimes_AN$ n'est pas exact à gauche en général
  (c'est-à-dire qu'en général, il ne préserve pas l'injectivité des flèches),
  comme le
 montre le contre-exemple suivant.
 
 \medskip
 On se place dans le cas où~$A=\ZZ$.
 Soit~$f: \ZZ\to \ZZ$
 la multiplication par~$2$ ; c'est une application~$\ZZ$-linéaire 
 injective. Pour tout~$\ZZ$-module~$M$,
 l'endomorphisme~${\rm Id}_M\otimes f$
 de~$M\otimes_{\ZZ}\ZZ\simeq M$
 est la multiplication par~$2$. 
 
 Lorsque~$M=\ZZ/2\ZZ$, celle-ci coïncide
 avec l'application nulle, et n'est en particulier pas injective. 
 
 \deux{def-module-plat}
 On dit qu'un~$A$-module~$M$ est
 {\em plat}
 si le foncteur~$N\mapsto M\otimes_AN$ est exact, c'est-à-dire
 s'il transforme les injections en injections. 
 
 \trois{comment-plat}
La platitude n'apparaîtra guère
{\em dans ce cours}, et 
c'est essentiellement à titre culturel que nous la mentionnons. 
Mais il s'agit d'une notion absolument cruciale en théorie des schémas, 
qui en dépit de sa définition purement algébrique un peu sèche a un
sens géométrique profond, et joue de surcroît un rôle technique majeur.

\trois{module-libre-plat}
Soit~$M$ un~$A$-module libre ; il est alors plat. En effet, choisissons
une base~$(e_i)$
de~$M$, et donnons-nous une injection~$A$-linéaire
$N\hookrightarrow  N'$.   

\medskip
On a~$M=\bigoplus A\cdot e_i$. On a
En conséquence,
on dispose d'après~\ref{pro-tens-somme-dir}
d'isomorphismes canoniques
$M\otimes_A N\simeq \bigoplus (A\cdot e_i\otimes_A N)$, 
et~$M\otimes_A N'\simeq  \bigoplus (A\cdot e_i\otimes_A N')$.
Il résulte par ailleurs de~\ref{prod-tens-libre1}
que l'on a pour tout
indice~$i$ des isomorphismes naturels
$A\cdot e_i\otimes_A N\simeq N$, et~$A\cdot e_i\otimes_A N'\simeq N'$.
Il s'ensuit que~$A\cdot e_i\otimes_A N\hookrightarrow A\cdot e_i\otimes_A N'$
pour tout~$i$, puis que~$M\otimes_A N\hookrightarrow M\otimes_A N'$.

\trois{esp-vect-corps}
Notons un cas particulier important de ce qui précède : {\em tout espace
vectoriel sur un corps est plat.}

\subsection*{Quelques objets classiques revisités}
\deux{m-tenseur-mdual}
Soient~$M$ et~$N$ deux~$A$-modules. L'application
de~$M^\vee \times N$ dans~$\hom_A(M,N)$ définie
par la formule
$$(\phi,n)\mapsto [m\mapsto \phi(m)n]$$ est bilinéaire, 
elle induit donc une application~$A$-linéaire~$\phi$ 
de~$M^\vee\otimes_A N$ vers~$\hom_A(M,N)$. 

\trois{mmdual-libre}
Supposons que~$M$ et~$N$ soient tous deux
{\em libres de rang fini.}
On choisit une base~$(e_i)_{1\leq i\leq m}$ de~$M$,
et une base~$(f_j)_{1\leq j\leq n}$
de~$N$. On déduit de la description matricielle des
applications entre modules libres que~$\hom_A(M,N)$ est libre
de rang~$nm$, une base étant donnée par la famille~$(u_{ij})$
où~$u_{ij}$ est caractérisé par les égalités
$u_{ij}(e_\ell)=\delta_{\ell i}f_j$ pour tout~$\ell$. 

\medskip
En appliquant cette remarque lorsque~$N=A$, on
voit que~$M^\vee$ est libre de base~$(e_i^*)$, où~$e_i^*$
désigne pour tout~$i$ la~$i$-ème forme coordonnée dans la base~$(e_1,\ldots, e_m)$. 

\medskip
On en déduit grâce à~\ref{prod-tens-mod-libres}
que le
module~$M^\vee\otimes_A N$ est libre de base~$(e_i^*\otimes f_j)_{ij}$. 

Fixons~$i$ et~$j$ et soit~$\ell\in \{1,\ldots,m\}$. On a par définition de~$\phi$ l'égalité
$$\phi(e_i^*\otimes f_j)(e_\ell)=e_i^*(e_\ell)f_j=\delta_{\ell i}f_j.$$
Autrement dit, $\phi(e_i^*\otimes f_j)=u_{ij}$. Ainsi, $\phi$ transforme une base
de~$M^\vee\otimes_A N$ en une base de~$\hom_A(M,N)$. Par conséquent,
$\phi$ est bijective. 

\trois{comment-mdual-m-hom}
{\em Commentaires.}
Ce qui précède est une illustration d'une démarche très fréquente en algèbre commutative,
(on l'a d'ailleurs déjà implicitement rencontrée à propos de la bidualité, {\em cf.}
l'exemple~\ref{ex-bidual}) : 

\medskip
$\bullet$ on commence par construire une application linéaire de manière complètement
naturelle (sans aucun choix à effectuer) ; 

$\bullet$ on montre ensuite, sous l'hypothèse qu'un ou plusieurs des modules en jeu sont
libres de rang fini, que cette application est bijective ; et pour ce faire, on
{\em choisit} une base 
dans laquelle on effectue les calculs. 

\medskip
Signalons par ailleurs que la bijection réciproque de~$\phi$
(lorsque~$M$ et~$N$ sont libres de rang fini) n'admet pas de description 
naturelle au moyen d'une formule explicite. 

\trois{trace-tenseur}
On se place maintenant dans le cas où~$N=M$, et on suppose toujours que~$M$
est libre de rang fini, de base~$(e_i)_{i=1,\ldots, m}$. 

\medskip
L'application
de~$M^\vee\times M$ vers~$A$ qui envoie~$(\phi,m)$ sur~$\phi(m)$
étant bilinéaire, elle induit une application linéaire~$\lambda$
de~$M^\vee\otimes_A M$
vers~$A$, et il existe une unique application linéaire~$\tau$ telle que le diagramme

$$\diagram
M^\vee\otimes_A M\dto_\lambda \rto^\phi_\simeq&{\rm End}_A(M)\dlto^\tau
\\A&\enddiagram$$ 
commute. Soient~$i$ et~$j$ deux entiers compris entre~$1$ et~$m$. 
On a
$$\tau(u_{ij})=\tau (\phi(e_i^*\otimes e_j))=\lambda (e_i^*\otimes e_j)=e_i^*(e_j)=\delta_{ij}={\rm Tr}\;(u_{ij}).$$
Ceci valant pour tout~$(i,j)$, la forme linéaire~$\tau$ coïncide avec la trace, dont on a ainsi donné une définition
intrinsèque (ne faisant pas intervenir un choix de base). 

\section{Produit tensoriel d'un module et d'une algèbre}
\markboth{Algèbre commutative}{Produit tensoriel d'un module et d'une algèbre}

\subsection*{Définitions, exemples et premières propriétés}

\deux{prod-mod-alg-intro}
On désigne toujours par~$A$ un anneau, et l'on se donne une~$A$-algèbre~$B$. 
Si~$M$ est un~$B$-module, il possède une structure naturelle
({\em i.e.}, fonctorielle en~$M$) de~$A$-module, 
obtenue par «restriction des scalaires» à~$A$. S'il n'y a pas de
risque de confusion, on notera encore~$M$ ce~$A$-module ; dans le
cas contraire, on écrira~$_AM$. 

\deux{def-struct-bmodule}
Soit~$M$ un~$A$-module. Nous allons montrer que le~$A$-module
~$B\otimes_A M$ possède une unique structure de~$B$-module, étendant sa structure
de~$A$-module et telle que~$\beta\cdot (b\otimes m)=(\beta b)\otimes m$ pour
tout~$(\beta,b,m)\in B^2\times M$. 

\trois{unique-mtenseurab}
L'unicité est claire : elle provient du fait que les tenseurs purs engendrent~$B\otimes_AM$
comme groupe abélien. 

\trois{exist-mtenseurab}
Montrons maintenant l'existence. Soit~$\beta \in B$. 
L'application de~$B\times M$ dans~$B\otimes_AM$ qui envoie
$(b,m)$ sur~$\beta b\otimes m$ est bi-$A$-linéaire. Elle induit donc 
une application~$A$-linéaire~$\mu_\beta$
de~$B\otimes_AM$ dans lui-même. On vérifie aussitôt 
(en testant comme d'habitude
les propriétés requise sur les tenseurs purs) que
l'application~$(\beta,v)\mapsto \mu_\beta(v)$ 
de~$B\times_A (B\otimes_AM)$ vers~$B\otimes_AM$ 
définit une structure de~$B$-module sur~$B\otimes_AM$
répondant à nos exigences. 

\trois{fonctorialite-btenseura-m}
Si~$f: M\to N$ est une application~$A$-linéaire, il 
est immédiat que~${\rm Id}_B\otimes f : B\otimes_AM\to B\otimes_AN$
est~$B$-linéaire. On peut donc voir~$M\mapsto B\otimes_AM$
comme un foncteur de~$\amod$
vers~$\bmod$. 

\trois{intuition-b-tenseur-m}
On dit que
le~$B$-module~$B\otimes_A M$ est déduit de~$M$ par {\em extension des scalaires
de~$A$ à~$B$.} Intuitivement,~$B\otimes_AM$ est le~$B$-module le
plus général fabriqué à partir de~$M$, en autorisant  la multiplication
externe par les éléments de~$B$, et non plus simplement de~$A$. 

Comme toujours, ce type de description informelle
se traduit rigoureusement en terme de propriété universelle,
ou encore de représentation d'un foncteur ; c'est l'objet du lemme ci-dessous. 

\deux{prop-univ-btensm}
{\bf Lemme.}
{\em Soit~$M$ un~$A$-module. 
Le couple$$(B\otimes_A M, m\mapsto 1\otimes m)$$
représente le foncteur covariant de~$\bmod$
vers~$\ens$
qui envoie~$P$ sur~${\rm Hom}_A(M,P)$.}

\medskip
{\em Démonstration.}
Notons pour commencer 
que~$m\mapsto 1\otimes m$ est bien~$A$-linéaire, et donc que l'énoncé a un sens. 

\medskip
Soit~$P$ un~$B$-module et soit~$f$ une application~$A$-linéaire de~$M$ dans~$P$. 
Il s'agit de montrer qu'il existe une unique application~$B$-linéaire~$g: B\otimes_AM\to P$
telle que~$g(1\otimes m)=f(m)$ pour tout~$m\in M$. 

\trois{uniquelin-btenseur}
{\em Unicité.}
Supposons qu'une telle~$g$ existe. On a alors pour tout~$(b,m)\in B\times N$ 
les égalités
$$g(b\otimes m)=g(b\cdot(1\otimes m))=bg(1\otimes m)=bf(m),$$
et
comme les tenseurs purs engendrent~$B\otimes_AM$ l'application~$g$
est bien uniquement déterminée. 

\trois{existelin-btenseur}
{\em Existence.}
On s'inspire de la formule exhibée dans la preuve de l'unicité. L'application
de~$B\times M$ vers~$M$ qui envoie~$b\otimes m$
sur~$bf(m)$ est bi-$A$-linéaire, et induit donc une application~$A$-linéaire~$g: B\otimes_AM\to P$, 
qui envoie~$b\otimes m$ sur~$bf(m)$ pour tout~$(b,m)$. On vérifie immédiatement que~$g$ est~$B$-linéaire, 
et l'on a bien par construction~$g(1\otimes m)=f(m)$ pour tout~$m$. ~$\Box$

\deux{reformule-btenseura-m}
On peut reformuler le lemme ci-dessus de différentes façons. 

\trois{adjonction-restriction-sclaaire}
{\bf Reformulation catégorique.} Le foncteur
$M\mapsto B\otimes_AM$ de~$\amod$
vers~$\bmod$ est adjoint à
gauche au foncteur~$N\mapsto _AN$ de~$\bmod$
vers~$\amod$. 

\trois{informel-extension-scalaire}
{\bf Reformulation informelle.}
Si~$M$ est un~$A$-module, 
se donner une application~$B$-linéaire de~$B\otimes_A M$ dans
un~$B$-module~$P$
revient à se donner une application~$A$-linéaire de~$M$ dans~$P$.

\deux{exemple-ext-scalaires}
Nous allons maintenant décrire explicitement~$B\otimes_A M$ 
dans un certain nombre de cas particuliers. 

\trois{ext-scalaire-libre}
Soit~$M$ un~$A$-module libre, et soit~$(e_i)$
une base de~$M$. On a la décomposition~$M=\bigoplus A\cdot e_i$. Par commutation
du produit tensoriel aux sommes directes, il vient
$B\otimes_A M\simeq \bigoplus B\otimes_A A\cdot e_i$.

Par ailleurs,
le~$A$-module~$A\cdot e_i$ est 
pour tout~$i$
libre de base~$e_i$ ; on en déduit 
grâce à~\ref{prod-tens-libre1}
que
$b\mapsto b\otimes e_i$ établit 
une bijection~$A$-linéaire entre~$B$
et~$B\otimes_A A\cdot e_i$. Comme~$b\otimes e_i=b\cdot(1\otimes e_i)$
pour tout~$(b,i)$, 
on voit finalement que~$B\otimes_A M$ est libre de base~$(1\otimes e_i)$. 

\trois{a-parenth-i-tenseurb}
Soit~$I$ un ensemble. Le~$A$-module~$A^{(I)}$ est libre ; soit~$(\theta_i)$
sa base canonique ($\theta_i$ est pour tout~$i$
 la famille~$(\delta_{ij})_j$ de~$A^{(I)}$). 

\medskip
Par ce qui précède,~$B\otimes_A A^{(I)}$ est libre de base~$1\otimes \theta_i$. Cela signifie
que
$$(b_i)\mapsto \sum b_i\cdot (1\otimes \theta_i)=\sum b_i \otimes \theta_i$$ établit un isomorphisme
entre~$B^{(I)}$ et~$B\otimes_A A^{(I)}$.

Modulo cet isomorphisme,
l'application naturelle
$$(a_i)\mapsto 1\otimes (a_i)=1\otimes \sum a_i\theta_i=\sum a_i \otimes \theta_i$$
s'identifie à la flèche~$A^{(I)}\to B^{(I)}$ déduite du morphisme structural de~$A$
vers~$B$. 

\trois{extenson-scal-gen-rel}
Soit maintenant~$M$ un~$A$-module quelconque
et soit~$(e_i)$ une famille génératrice de~$M$. 
L'unique
application linéaire
de~$A^{(I)}$ dans~$M$ qui envoie~$\theta_i$
sur~$e_i$
pour tout~$i$ (c'est celle 
qui envoie toute famille~$(a_i)$
sur~$\sum a_i e_i$)
est alors surjective ; soit~$(f_\ell)_{\ell \in \Lambda}$
une famille génératrice de son noyau. 
On dispose d'une suite exacte 
$$\diagram A^{(\Lambda)}\rrto^{(a_\ell)\mapsto \sum a_\ell f_\ell}&&A^{(I)}\rrto^{(a_i)\mapsto \sum a_i e_i}&&M\rto&0\enddiagram,$$
c'est-à-dire encore d'un isomorphisme
$[A^{(I)}/(f_\ell)_\ell]\simeq M$ envoyant~$\overline{\theta_i}$
sur~$e_i$ pour tout~$i$.

\medskip
Par exactitude à droite
du produit tensoriel
et en vertu du~\ref{a-parenth-i-tenseurb}
ci-dessus, cette suite induit
{\em via}
la tensorisation avec~$B$ une suite exacte

$$
\diagram
B^{(\Lambda)}\rrto^{(b_\ell)\mapsto \sum b_\ell f_\ell}&&B^{(I)}\rrto^{(b_i)\mapsto \sum b_i (1\otimes e_i)}&&B\otimes_AM\rto&0\enddiagram,$$
c'est-à-dire un isomorphisme
$[B^{(I)}/(f_\ell)_\ell]\simeq B\otimes_AM$
envoyant~$\overline{\theta_i}$
sur~$1\otimes e_i$ pour tout~$i$ 
(par abus, on désigne encore
par~$f_\ell$ et~$\theta_i$
les images respectives de~$f_\ell$ et~$\theta_i$ dans~$B^{(I)}$
par la flèche~$A^{(I)}\to B^{(I)}$) ; notons en particulier que~$1\otimes e_i$
est une famille génératrice de~$B\otimes_AM$. 

\medskip
De manière un peu informelle, 
on voit que le~$A$-module~$M$
et le~$B$-module~$B\otimes_AM$ 
admettent la
«même» description par générateurs (les~$\theta_i$) et relations
(les~$f_\ell$). 

\medskip
On peut résumer cela
par le
slogan suivant, vague mais assez
intuitif : {\em $B\otimes_AM$ est à~$B$ ce que~$M$ est à~$A$.}

\trois{formules-al-extension}
Ce principe s'applique aussi aux applications linéaires. Plus précisément,
donnons-nous deux~$A$-modules~$M$ et~$N$, 
et choisissons une famille génératrice~$(e_i)$ de~$M$
et une famille génératrice~$(f_j)$ de~$N$. Soit~$u$ une application~$A$-linéaire
de~$M$ vers~$N$. Pour tout~$i$, il existe une famille~$(a_{ij})$ d'éléments
de~$A$ telle que~$u(e_i)=\sum a_{ij}f_j$, et elle détermine entièrement~$u$ :
on a~$u(\sum \lambda_i e_i)=\sum_j(\sum _i \lambda_i a_{ij})f_j$ pour toute
famille~$(\lambda_i)$ de scalaires. 

\medskip
Il résulte de~\ref{extenson-scal-gen-rel}
que~$(1\otimes e_i)$ est une famille génératrice de~$B\otimes_A M$, 
et
que~$(1\otimes f_j)$ est une famille génératrice de~$B\otimes_A N$. 
On a de plus pour tout~$i$ les égalités
$$({\rm Id}_B\otimes u)(1\otimes e_i)=1\otimes u(e_i)=1\otimes(\sum a_{ij}f_j)=\sum a_{ij}(1\otimes f_j).$$
On voit ainsi
que~${\rm Id}_B\otimes u : B\otimes_A M\to B\otimes_A N$ est décrite,
dans les familles génératrices~$(1\otimes e_i)$ et~$(1\otimes f_j)$, par les «mêmes»
formules que~$u$ dans les familles génératrices~$(e_i)$
et~$(f_j)$.

\deux{m-tensb-tensc}
 {\bf Extensions des scalaires successives.}
 Soit~$M$
 un~$A$-module, soit~$B$ une~$A$-algèbre, 
 et soit~$C$ une~$B$-algèbre. On dispose alors d'un isomorphisme 
 naturel de~$C$-modules 
 
 $$C\otimes_B(B\otimes_AM)\simeq C\otimes_AM.$$
 
 On peut le voir de deux façons différentes.
 
 \trois{m-tensb-tensc-mains}
 {\em Première méthode.}
 Soit~$c\in C$. L'application
de~$B\times M$ dans~$ C\otimes_AM$ qui envoie~$(b,m)$
sur~$cb\otimes m$
est bilinéaire, et induit donc une application~$A$-linéaire~$m_c$ 
de~$B\otimes_A M$ vers~$C\otimes_A M$. 

L'application de~$C\times (B\otimes_AM)$
vers~$(C\otimes_AM)$ qui envoie~$(c,v)$ sur~$m_c(v)$ est~$B$-bilinéaire, 
donc induit une application~$B$-linéaire
$$\phi : C\otimes_B(B\otimes_AM)\to C\otimes_A M.$$

\medskip
L'application
de~$C\times M$ vers~$C\otimes_B(B\otimes_AM)$ qui envoie~$(c,m)$
sur~$c\otimes (1\otimes m)$
est bilinéaire, et induit donc une application~$A$-linéaire
$$\psi : C\otimes_A M
\to 
C\otimes_B(B\otimes_AM).$$
On vérifie
aisément
que~$\phi$ et~$\psi$ sont~$C$-linéaires et réciproques
l'une de l'autre. 

\trois{m-tensb-tensc-fonct}
{\em Preuve fonctorielle.}
Soit~$P$ un~$C$-module. On dispose d'isomorphismes naturels

$$\hom_C(C\otimes_B(B\otimes_A M), P)\simeq \hom_B(B\otimes_A M,P)\simeq \hom_A(M,P)$$
$$\simeq \hom_C(C\otimes_A M, P)$$
qui sont fonctoriels en~$P$ et~$M$. 
Par composition on obtient un isomorphisme 
naturel
 $\hom_C(C\otimes_B(B\otimes_A M), P)\simeq \hom_C(C\otimes_AM,P)$
 fonctoriel en~$P$ et~$M$.
 Le lemme de Yoneda assure qu'il
 provient d'une bijection~$C$-linéaire de~$C\otimes_AM$ vers~$C\otimes_B(B\otimes_AM)$,
dont on vérifie qu'elle coïncide avec~$\psi$. 

\subsection*{Comportement vis-à-vis des localisations et quotients}

\deux{ext-scal-loc-quot}
Soit~$M$ un~$A$-module. Nous allons décrire l'extension des scalaires de~$M$
à deux types de~$A$-algèbres particulières, à savoir les quotients et les localisations. 

\trois{ext-scal-quot}
{\em Le cas des quotients.}
Soit~$I$ un idéal de~$A$. La
structure de~$A$-module du quotient~$M/IM$
est induite par une structure de~$A/I$-module
sur ce dernier (la multiplication externe par un scalaire~$a$ ne dépend
dans ce module que de la classe de~$a$ modulo~$I$). 

\medskip
La surjection~$M\to M/IM$ étant~$A$-linéaire, elle induit un
morphisme de~$A/I$-modules~$p: (A/I)\otimes_AM\to M/IM$. 

Par ailleurs, l'application~$A$-linéaire~$m\mapsto 1\otimes m$
de~$M$ dans~$(A/I)\otimes_A M$ s'annule sur~$IM$ : en effet, 
si les~$a_i$ sont des éléments de~$I$ et les~$m_i$ des éléments de~$M$, on a
$$1\otimes \sum a_im_i=\sum a_i \otimes m_i=0,$$ puisque les facteurs de gauche
vivent dans~$A/I$. 

Elle induit donc
une application~$A$-linéaire~$s : M/IM\to (A/I)\otimes_AM$, qui comme 
toute application~$A$-linéaire entre~$A/I$-modules est automatiquement~$A/I$-linéaire
(par surjectivité de~$A$ vers~$A/I$). On vérifie immédiatement que~$p$ et~$s$
sont inverses l'une de l'autre. 

\medskip
On a ainsi construit un isomorphisme canonique de~$A/I$-modules

$$(A/I)\otimes_A M\simeq M/IM.$$

\trois{ext-scal-localis}
{\em Le cas des localisations.}
Soit~$S$ une partie multiplicative de~$A$. 
On définit sur~$M\times S$ la relation~$\mathscr R$ suivante : $(m,s){\mathscr R}(n,t)$ si et seulement
 si il existe~$r\in S$ tel que~$r(tm-sb)=0$. 
 On vérifie que c'est une relation d'équivalence, et l'on note~$S^{-1}M$ le quotient correspondant. 
 Les formules
 $$((m,s),(n,t))\mapsto (tm+sn, st)\;{\rm et}\;((a,s),(m,t))\mapsto (am,st)$$ 
 passent au quotient, et définissent une loi
 interne~$+$ sur~$S^{-1}M$ 
 et une  loi
 externe~$\times : S^{-1}A\times S^{-1}M\to S^{-1}M$ qui font
 de~$S^{-1}M$ un~$S^{-1}A$-module (la preuve de ce fait,
 aussi triviale que fastidieuse, est laissée au lecteur). 
  
 Si~$(m,s)\in M\times S$, on écrira~$\frac m s$ au lieu de~$\overline{(m,s)}$. Cette notation permet de disposer
 des formules naturelles
  $$\frac m s + \frac n t=\frac{sn+tm}{st}\;{\rm et} \;\frac a s \cdot \frac m t=\frac{am}{bt},$$
 et l'on a~$$\frac m s =\frac n t \;\iff\;\exists r\in S, r(tm-sn)=0.$$
 
 Si~$f$ est une application~$A$-linaire 
 de~$M$ vers un~$A$-module~$N$, on vérifie que la
 formule~$(m,s)\mapsto \frac{f(m)}s$ passe au quotient par~$\sch R$
 et induit une
 application~$S^{-1}A$-linéaire de~$S^{-}M$ vers~$S^{-1}N$, qui envoie
 par construction une fraction~$\frac m s$
 sur la fraction~$\frac {f(m)}s$. Ainsi, $M\mapsto S^{-1}M$
 apparaît comme un foncteur de~$\amod$ vers~$S^{-1}A\text{-}\mathsf{Mod}$. 
 
 \medskip
L'application~$m\mapsto \frac m 1$ de~$M$ dans~$S^{-1}M$ est~$A$-linéaire ; elle induit
donc une application~$S^{-1}A$-linéaire~$\phi$ de~$S^{-1}A\otimes_A M$ dans~$S^{-1}M$. 

Par ailleurs, si~$m\in M$ et~$s\in S$, on vérifie immédiatement que l'élément~$\frac 1 s\otimes m$
de~$S^{-1}A\otimes_A M$ ne dépend que de la classe de~$(m,s)$ modulo~$\sch R$, que l'application
$\psi : S^{-1}M\to S^{-1}A\otimes_A M$ construite par ce biais est~$S^{-1}A$-linéaire, et que~$\phi$
et~$\psi$ sont inverses l'une de l'autre. 

\medskip
On a ainsi construit un isomorphisme de~$S^{-1}A$-modules
$$S^{-1}A\otimes_A M\simeq S^{-1}M$$ qui est visiblement fonctoriel en~$M$. 

 \trois{smoinsun-m-plat}
 {\em Une application importante.} Soit~$f$ une injection
 $A$-linéaire de~$M$ dans un~$A$-module~$N$,
 soit~$m\in M$ et soit~$s\in S$. Supposons que~$\frac{f(m)}s=0$ ; 
 cela signifie qu'il existe~$r\in S$ tel que~$rf(m)=0$, ou encore tel que~$f(rm)=0$. 
 Mais comme~$f$ est injective, il vient~$rm=0$ puis~$\frac m s=0$. Ainsi,
 l'application
 linéaire~$S^{-1}M\to S^{-1}N$ induite par~$f$ est injective. 
 
 \medskip
 {\em En conséquence, le~$A$-module~$S^{-1}A$
 est plat}. 
 
 \trois{mp-limind-mf}
 {\em Localisation des modules et limites inductives.}
 Soit~$\Sigma$ une partie multiplicative de~$A$ et soit~$(I\leq)$ un ensemble
pré-ordonné filtrant. 
 Soit~$(S_i)_{i\in I}$ une famille de parties multiplicatives de~$A$ contenues dans~$\Sigma$, 
 telles que tout élément de~$S_i$ soit inversible dans~$S_j$ dès que~$i\leq j$ ; on suppose
 de plus que les~$S_i$ engendrent
 $\Sigma$ multiplicativement. Soit~$\sch D$ le diagramme de~$A$-algèbres
 dont la famille d'objets
 est~$(S_i^{-1}A)_{i\in I}$ et dont les flèches sont les morphismes canoniques 
 $S_i^{-1}A\to S_j^{-1}A$ pour~$i\leq j$, et soit~$M\otimes_A \sch D$ l'image de~$\sch D$ par
 le foncteur~$M\otimes_A\bullet $ ; les objets de~$M\otimes_A \sch D$ sont les~$S_i^{-1}M$, 
 et ses morphismes
 sont les flèches canoniques $S_i^{-1}M\to S_j^{-1}M$ pour~$i\leq j$. 
 Le morphisme canonique de~$\sch D$ dans~$\Sigma^{-1}A$ identifie ce dernier
 à la limite inductive de~$\sch D$ dans la catégorie des~$A$-modules 
 (\ref{limind-smoinsun-a}
 et remarque~\ref{limind-loc-multicat}) ; puisque le produit tensoriel
 commute aux limites inductives, la flèche naturelle $M\otimes_A\sch D\to \Sigma^{-1}M$
 identifie ce dernier à la limite inductive de~$M\otimes_A\sch D$ dans la catégorie des~$A$-modules.

 \medskip
 Le lecteur que rebuterait l'invocation de la commutation aux limites inductives pourra donner une
 démonstration directe de ce fait. Il suffit en effet de s'assurer 
 que les assertions
\ref{asi-asigma-sur}
et~\ref{asi-sigma-inj}
restent vraies lorsqu'on remplace partout les localisés de~$A$ par ceux de~$M$, ce qui se fait sans
la moindre difficulté, en reprenant leurs preuves
{\em mutatis mutandis}. 

\medskip
Mentionnons un cas particulier important. Soit~$\got p$ un idéal premier de~$A$, et soit~$\Delta$
le diagramme
$$((M_f)_{f\in A\setminus \got p}, (M_f\to M_g)_{f|g})$$ 
(qui se déduit par tensorisation avec~$M$ de celui considéré au~\ref{ap-limind-af}). Sa limite inductive dans la catégorie des~$A$-modules
s'identifie alors à~$M_{\got p}$.

\section{Modules projectifs}
\markboth{Algèbre commutative}{Modules projectifs}

{\em On fixe pour toute cette section un anneau~$A$.}

\subsection*{Propriétés se testant sur une famille couvrante de localisés}

\deux{fam-compl-loc}
Soit~$(S_i)_{i\in I}$ une famille de parties multiplicatives de~$A$. Nous dirons que la famille~$(S_i)$ est {\em couvrante}
si pour tout idéal premier~$\got p$ de~$A$, il existe~$i$ tel que~$\got p$ 
 ne rencontre pas~$S_i$. 
 
 \trois{comment-couvrante}
 Il revient au même de demander que pour tout idéal premier~$\got p$ de~$A$, il existe~$i$ tel que~$S_i$ s'envoie dans les éléments
 inversibles de~$A_{\got p}$, c'est-à-dire encore tel qu'il existe un morphisme de~$A$-algèbres de~$S_i^{-1}A$ dans~$A_{\got p}$. 
 
 \trois{exemple-couvrante}
 Les deux exemples fondamentaux de famille
 couvrante à avoir en tête sont les suivants. 
 
 \medskip

a) La famille~$(A\setminus \got p)_{\got p\in \spec A}$ est couvrante par définition. 
 
 b) Soit~$(f_i)_{i\in I}$ une famille d'éléments de~$A$ telle que l'idéal engendré par les~$f_i$
 soit égal à~$A$. La famille~$(\{f_i^n\}_{n\in \NN})_i$ est alors couvrante. 
 En effet, si~$\got p$ est un idéal 
 premier de~$A$, il ne peut contenir toutes les~$f_i$ puisqu'elles engendrent~$A$, et 
 si l'on choisit~$i$ tel que~$f_i\notin \got p$
 alors~$\got p$ ne rencontre pas~$\{f_i^n\}_{n\in \NN}$. 
 
 \deux{fixe-couvrante}
 On fixe une famille couvrante~$(S_i)$ de parties multiplicatives de~$A$. Le but de ce qui suit est de montrer
 que certaines propriétés (d'un module, d'un morphisme...) sont vraies si et seulement si elles sont vraies après 
 localisation par chacune des~$S_i$.

 \trois{lemme-null-msi}
 {\bf Lemme.}
 {\em Soit~$M$ un~$A$ module. La flèche naturelle~$M\to \prod_i S_i^{-1}M$ 
 est injective.}
 
 \medskip
 {\em Démonstration.}
 Soit~$m$ un élément tel que~$\frac m 1=0$ dans~$S_i^{-1}M$ pour tout~$i$ et soit~$J$ l'idéal 
 annulateur de~$M$. Fixons~$i$ ; par hypothèse, il existe~$s_i \in S_i$ tel que~$s_im=0$ ; en conséquence, 
 $J$ rencontre tous les~$S_i$. Il n'est dès lors contenu dans aucun idéal premier de~$A$, ce qui veut dire 
 qu'il est égal à~$A$ ; il vient~$m=1\cdot m=0$.~$\Box$ 
 
 \trois{coroll-nul}
 {\bf Corollaire.}
 {\em Le module~$M$ est nul si et seulement si tous les~$S_i^{-1}M$ sont nuls.}~$\Box$
 
 \trois{coroll-reduit}
 {\bf Corollaire.}
 {\em L'anneau~$A$ est réduit si et seulement si~$S_i^{-1}A$ est réduit pour tout~$i$.}
 
 \medskip
 {\em Démonstration.}
 Supposons~$A$ réduit, et soit~$i\in I$. Soient~$a\in A$ et~$s\in S_i$ tels que
 l'élément~$\frac a s$ de~$S_i^{-1}A$ soit nilpotent. Il existe~$n\geq 1$ tel que~$\frac {a^n}{s^n}=0$ dans~$S_i^{-1}A$, 
 ce qui veut dire qu'il existe~$t\in S_i$ tel que~$ta^n=0$. On a {\em a fortiori}
 $(ta)^n=0$, et comme~$A$ est réduit $ta=0$, ce qui entraîne que~$\frac a s=0$ dans~$S_i^{-1}A$. Ainsi, ce dernier est réduit.
 
 \medskip
 Réciproquement, supposons~$S_i^{-1}A$ réduit pour tout~$i$, et soit~$a$ un élément nilpotent de~$A$. 
 Son image dans chacun des~$S_i^{-1}A$ est nilpotente, donc nulle et le lemme~\ref{lemme-null-msi}
 assure alors que~$a=0$.~$\Box$ 
 
 \trois{lemme-ex-msi}
 {\bf Lemme.}
 {\em Soit 
 $$\sch D=\xymatrix{{M'}\ar[r]^u&M\ar[r]^v&{M''}}$$ un diagramme de~$\amod$. C'est une suite exacte si et seulement
 si $$S_i^{-1}\sch D:=\xymatrix{{S_i^{-1}M'}\ar[rr]^{S_i^{-1}u}&&{S_i^{-1}M}\ar[rr]^{S_i^{-1}v}&&{S_i^{-1}M''}}$$ est exacte pour tout~$i$.}
 
 \medskip
 {\em Démonstration.}
 Le sens direct provient de la platitude du~$A$-module
 $S_i^{-1}A$ pour tout~$i$ (\ref{smoinsun-m-plat}). Supposons maintenant que
 $S_i^{-1}\sch D$ soit une suite exacte pour tout~$i$ et soit~$m\in M$. Par hypothèse, $v\circ u(m)$ s'annule dans~$S_i^{-1}M''$ pour tout~$i$
 (puisque~$S_i^{-1}v\circ S_i^{-1}u=0$), 
 et est donc
 nul d'après le lemme~\ref{lemme-null-msi} ; ainsi, $v\circ u=0$. 
 
 \medskip
 Soit~$P$ le conoyau de~$u$ ; la flèche~$v$ induit d'après
 ce qui précède une flèche~$P\to M''$, et il s'agit de montrer qu'elle
 est injective. Fixons~$i$. L'exactitude à droite du produit tensoriel assure que
 $$S_i^{-1}M'\to S_i^{-1}M\to S_i^{-1}P\to 0$$ est exacte ; cela signifie que~$S_i^{-1}P$
 s'identifie au conoyau de~$S_i^{-1}u$,  et l'exactitude de~$S_i^{-1}\sch D$ entraîne alors
 l'injectivité de~$S_i^{-1}P\to S_i^{-1}M''$. 
 
 \medskip
 Soit maintenant~$p\in P$ un élément dont l'image dans~$M''$ est nulle. Puisque~$S_i^{-1}P\to S_i^{-1}M''$
 est injective pour tout~$i$, l'élément~$p$ appartient au noyau de~$P\to \prod S_i^{-1}P$, et est dès lors
 nul d'après le lemme~\ref{lemme-null-msi}.~$\Box$ 
 
 \deux{intro-type-fini}
 Nous allons maintenant énoncer un résultat de la même veine
 que les précédents, mais qui requiert que la famille
 des~$S_i$ soit finie (il ne pourra donc pas s'utiliser en général avec
 la famille~$(A\setminus \got p)_{\got p\in \spec A}$). 
 
 \trois{type-fini-si}
 {\bf Lemme.}
 {\em Supposons que l'ensemble d'indice~$I$ est fini, et soit~$M$ un~$A$-module. 
 Il est de type fini si et seulement si~$S_i^{-1}M$ est un~$S_i^{-1}A$-module de type fini pour tout~$i$.}
 
 \medskip
 {\em Démonstration.}
 L'implication directe provient du fait que l'extension des scalaires préserve 
 en vertu de l'exactitude à droite du produit tensoriel
 la propriété d'être de type fini ({\em cf.} \ref{extenson-scal-gen-rel})
pour le raisonnement précis).
On suppose maintenant que~$S_i^{-1}M$ est de type fini
 pour tout~$i$. Fixons~$i$. Il existe
 par hypothèse des éléments~$m_{i1},\ldots, m_{in_i}$
 de~$M$ et des éléments~$s_{i1},\ldots, s_{in_i}$ de~$S_i$ tels que
 la famille~$(\frac{m_{i1}}{s_{i1}},\ldots, \frac{m_{in_i}}{s_{in_i}})$
 engendre le~$S_i^{-1}A$-module~$S_i^{-1}M$. 
 Comme les~$s_{ij}$ sont inversibles dans~$S_i^{-1}A$, la famille des~$m_{ij}$ engendre encore~$S_i^{-1}M$
 comme~$S_i^{-1}A$-module. 
 
 \medskip
 Soit~$L$ un~$A$-module libre de base~$(e_{ij})_{i\in I, 1\leq j\leq n_i}$ et soit~$f\colon L\to M$ l'application~$A$-linéaire qui envoie
 $e_{ij}$ sur~$m_{ij}$ pour tout~$(i,j)$. Par construction, $S_i^{-1} f\colon S_i^{-1} L\to S_i^{-1} M$ est surjective pour tout~$i$, 
 et~$f$ est donc surjective en vertu du lemme~\ref{lemme-ex-msi}. Comme~$I$ est fini, la famille~$(e_{ij})$ est finie et~$M$
 est de type fini.~$\Box$ 
 
 \trois{rem-tf-loc}
 {\em Remarque.}
 Le corollaire ci-dessus est faux en général sans hypothèse de finitude sur~$I$ ; donnons un contre-exemple. 
 On se place sur l'anneau~$\ZZ$, et l'on considère la famille de parties multiplicatives~$(\ZZ\setminus(p))_{p\;\text{premier}}$,
 laquelle est couvrante : 
 un idéal premier de la forme~$p\ZZ$ avec~$p$ premier ne rencontre pas~$\ZZ\setminus(p)$, et~$(0)$ ne rencontre
 quant à lui aucune
 des~$\ZZ\setminus(p)$. 
 
 \medskip
 Soit~$M$ le~$\ZZ$-module~$\bigoplus_p \ZZ/p\ZZ$, la somme étant prise sur tous les nombres premiers. 
 Un calcul immédiat (faites-le) montre que le~$\ZZ_{(p)}$-module~$M_{(p)}$ est isomorphe
 pour tout nombre premier~$p$ à~$\ZZ/p\ZZ=\ZZ_{(p)}/p\ZZ_{(p)}$ ; il est donc de type fini
 (et même de présentation finie). Pourtant, $M$ n'est pas de type fini : il ne comprend que des éléments
 d'ordre fini, donc serait fini s'il était de type fini. 
 
 \subsection*{Suites exactes scindées, modules projectifs}
 
 \deux{intro-split}
 Soit
 $$\xymatrix{0\ar[r]&{M'}\ar[r]^i&M\ar[r]^p &{M''}\ar[r]&0}$$ une suite exacte de~$A$-modules. 
 
 \trois{equiv-split}
 Les assertions
 suivantes sont équivalentes : 
 
 \medskip
 i) la surjection~$p$ admet une {\em section}, c'est-à-dire une application~$A$-linéaire
$s: M''\to M$ telle que~$p\circ s={\rm Id}_{M''}$ ; 
 
 ii) il existe un isomorphisme~$\theta\colon M'\oplus M''\simeq M$ tel que
 le diagramme 
 
 $$\xymatrix{
 &&&M\ar[rrd]^p&&&\\
 0\ar[r]&{M'}\ar[rru]^i\ar[rrd]&&&&{M''}\ar[r]&0\\
 &&&{M'\oplus M''}\ar[rru] \ar[uu]^{\simeq}_\theta&&&&}$$
 commute. 
 
 \medskip
 En effet, supposons que~i) soit vraie et montrons que
 $$\theta \colon (m',m'')\mapsto i(m')+s(m'')$$
 convient. Soit~$m\in M$. On a~$p(m-s(p(m))=p(m)-p(m)=0$,
 et~$m-s(p(m))$ appartient donc à~${\rm Ker}\;p={\rm Im}\;i$. Si l'on pose~$m''=p(m)$ et si l'on note~$m'$ l'unique élément
 de~$M'$ tel que~$m-s(p(m))=i(m')$ on a donc~$m=i(m')+s(m'')$ et~$\theta$ est surjectif. 
 
 Soit maintenant~$(m',m'')\in M'\times M''$ tel que~$\theta(m',m'')=i(m')+s(m'')=0$.
 En appliquant~$p$ il vient~$0=p(i(m'))+p(s(m''))=m''$. 
 On a alors~$i(m')=0$, et partant~$m'=0$ par injectivité de~$i$. Ainsi, $(m',m'')=(0,0)$
 et~$\theta$ est injectif. 
 
 L'application~$A$-linéaire~$\theta$ est en conséquence
 un isomorphisme ; qu'elle fasse
 commuter le diagramme résulte du fait que~$p\circ s={\rm Id}_{M''}$. 
 
 \medskip
 Réciproquement, supposons que~ii) soit vraie. On vérifie aussitôt que
 l'application linéaire $m''\mapsto \theta(m'',0)$ est une section
 de~$p$.

 \trois{def-scind}
 Lorsque ces conditions équivalents sont satisfaites, 
 on dit que la suite exacte 
 $$\xymatrix{0\ar[r]&{M'}\ar[r]^i&M\ar[r]^p &{M''}\ar[r]&0}$$
 est {\em scindée}. 
 
 \trois{exo-scind}
 À titre  d'exercice, montrez que $$\xymatrix{0\ar[r]&{M'}\ar[r]^i&M\ar[r]^p &{M''}\ar[r]&0}$$
 est scindée si et seulement si~$i$ admet une {\em rétraction}, c'est-à-dire
 une application linéaire~$r\colon M\to M'$ telle que~$r\circ i={\rm Id}_{M'}$. 
 
 \trois{contrex-scind}
 La propriété d'être scindée n'a rien d'automatique : ainsi,
 la suite exacte de~$\ZZ$-modules
 $$\xymatrix{0 \ar[r]&\ZZ\ar[rr]^{\times 2}&&\ZZ\ar[r]&{\ZZ/2\ZZ}\ar[r]&0}$$ n'est pas 
 scindée, par exemple parce qu'il n'existe aucun isomorphisme entre~$\ZZ$ et~$\ZZ\oplus \ZZ/2\ZZ$,
 puisque~$\ZZ$ n'admet pas d'élément de~$2$-torsion non trivial. 
 
 \deux{exact-hom}
 Soit~$P$ un~$A$-module. Nous laissons le lecteur vérifier que le foncteur covariant
 $M\mapsto \hom(P, M)$ est exact à gauche. Nous dirons que~$P$ est
 {\em projectif}
 si ce foncteur est exact, c'est-à-dire encore si pour toute surjection~$p\colon M\to N$, la flèche 
 $u \mapsto p\circ u$ de~$\hom (P, M)$ vers~$\hom (P, N)$ est surjective ; en termes plus imagés,
 cela signifie que toute application linéaire de~$P$ vers~$N$ se relève à~$M$. 
 
 \trois{ex-base}
 {\em Si le module~$P$ est libre, il est projectif.}
 En effet, supposons que~$P$ admette une base~$(e_i)_{i\in I}$, et soit~$p \colon M\to N$ une surjection
 linéaire. Donnons-nous
 une application linéaire~$u\colon P\to N$. 
 Choisissons\footnote{Si~$I$ est infini, cela requiert l'axiome du choix.}
 pour tout~$i$ un antécédent~$m_i$ de~$u(e_i)$ dans~$N$. Soit~$v$ l'unique
 application~$A$-linéaire de~$P$ dans~$M$ envoyant~$e_i$ sur~$m_i$ pour tout~$i$. 
 On a alors pour tout~$i$
 les égalités~$p(v(e_i))=p(m_i)=u(e_i)$ ; comme une application linéaire de source~$P$ est connue
 dès qu'on connaît son effet sur une base, il vient~$p\circ v=u$, et~$P$ est donc bien projectif. 
 
 \trois{libre-proj-pasequiv}
 {\em Remarque.}
 La réciproque de~\ref{ex-base}
 est fausse : il y a des exemples de modules projectifs qui ne sont pas libres. Par
 exemple, si~$A$ est un anneau de Dedekind ({\em e.g.}
 $A$ est un anneau d'entiers de corps de nombres) et~$I$ un idéal de~$A$ alors~$I$ est un module
 projectif, qui est libre si et seulement si il est principal (nous esquisserons plus loin
 la preuve de ce fait) ; et il y a des exemples d'anneaux de Dedekind non principaux. 

 \trois{suite-scind}
 Supposons que~$P$ est projectif, et donnons-nous une suite exacte
 $$\xymatrix{0\ar[r]&L\ar[r]^i&M\ar[r]^p &P\ar[r]&0}.$$ Elle est alors
 automatiquement scindée. En effet, l'application linéaire~$\hom(P, M)\to \hom(P,P)$
 induite par~$p$ est 
 surjective par projectivité de~$P$. En particulier, ${\rm Id}_P$ a un antécédent~$s\in \hom(P,M)$ ; 
 par définition, cela signifie que~$p\circ s={\rm Id}_P$, et~$s$ est donc une section de~$p$, ce qui achève la
 preuve. 

\trois{proj-loc}
{\bf Lemme.}
{\em Soit~$A$ un anneau local et soit~$P$ un module projectif de type fini sur~$A$. 
Le module~$P$ est libre (de rang fini, {\em cf.}~\ref{libre-tf}).}

\medskip
{\em Démonstration.}
Soit~$\got m$ l'idéal  maximal de~$A$. Le quotient~$P/\got m P$ est un~$(A/\got m)$-espace
vectoriel de dimension finie. Choisissons une famille~$(e_1,\ldots, e_n)$ 
de~$P$ dont les classes modulo~$\got mP$ constituent
une base de~$P/\got m P$, et soit~$p : A^n\to P$ le morphisme~$(a_i)\mapsto \sum a_i e_i$. 

\medskip
Par construction, la flèche~$(A/\got m)^n\to P/\got mP$ induite par~$p$ est bijective, 
et en particulier surjective. On déduit du lemme de Nakayama, ou plus précisément
de l'un de ses avatars (\ref{surj-test-residuel})
que~$p$ est surjective. Soit~$K$ son noyau. Comme~$P$ est projectif, 
il résulte de~\ref{suite-scind}
qu'il existe un isomorphisme~$A^n\simeq K\oplus P$ modulo lequel~$p$
est la seconde projection. On en déduit  en quotientant  modulo~$\got m$
un isomorphisme~$(A/\got m)^n\simeq K/ \got mK \oplus P/\got mP$ modulo
lequel~$(A/\got mA)^n\to P/\got mP$ est la seconde projection. Mais
on a signalé ci-dessus que~$(A/\got mA)^n\to P/\got mP$ est un isomorphisme ; 
en conséquence,~$K/\got mK=0$. L'isomorphisme~$A^n\simeq K\oplus P$ assure
par ailleurs que~$K$ s'identifie à un quotient de~$A^n$, et est en particulier de type fini ; 
on déduit alors du lemme de Nakayama, et cette fois-ci plus
précisément du corollaire~\ref{coro-nakayama}, que~$K$ est nul. Ainsi, $A^n\simeq P$.~$\Box$

 \deux{equi-proj-not}
 Soit~$P$ un~$A$-module. Choisissons une famille génératrice~$(p_i)_{i\in I}$ de~$P$, 
 soit~$L$ un module libre quelconque de base~$(e_i)_{i\in I}$ paramétrée
 par~$I$, soit~$\pi : L\to P$ la surjection~$\sum a_i e_i \mapsto \sum a_i p_i$, et soit~$K$ son noyau. 
 
 \deux{equi-proj-theo}
 {\bf Théorème.}
 {\em Les assertions suivantes sont équivalentes : 
 
 \medskip
 i) $P$ est projectif ; 
 
 ii) toute suite exacte de la forme
 $$0\to N\to M\to P\to 0$$ est sicndée ; 
 
 iii)  la suite exacte
 $$ \xymatrix{0\ar[r]&K\ar@^{{(}->}[r]&L\ar[r]^\pi &P\ar[r]&0}$$
 est scindée ; 
 
 iv) Il existe un isomorphisme~$L\simeq P\oplus K$ ; 
 
 v) $P$ est {\em facteur direct} d'un module libre, c'est-à-dire qu'il existe un module~$\Lambda$
 tel que~$\Lambda \oplus P$ soit libre.}
 
 \medskip
 {\em Démonstration.}
 L'implication i)$\Rightarrow$ii) a été vue au~\ref{suite-scind}. Les implication~ii)$\Rightarrow$iii), iii)$\Rightarrow$iv)
 et~iv)$\Rightarrow$v) sont évidentes. 
 
 \medskip
 Supposons maintenant que~v) soit vraie. 
 Soit~$\pi \colon M\to N$ une surjection linéaire, et soit~$u \colon P\to N$
 une application linéaire. Comme~$\Lambda\oplus P$ est libre, il est projectif (\ref{ex-base}), et
 l'application~$0\oplus u$ de~$\Lambda\oplus P$ vers~$N$ se relève donc en une 
 application linéaire~$v \colon \Lambda \to M$ ; par construction, $v|_P$ relève~$u$ et~$P$ est projectif.~$\Box$
 
 \trois{coro-proj-extscal}
 {\bf Corollaire.}
 {\em Soit~$M$ un~$A$-module projectif et soit~$B$ une~$A$-algèbre. 
 Le~$B$-module~$B\otimes_A M$ est projectif.}
 
 \medskip
 {\em Démonstration.}
 Le théorème ci-dessus assure qu'il existe un~$A$-module~$\Lambda$
 tel que~$M\oplus\Lambda$ soit libre. Le~$B$-moodule
 $B\otimes_A(M\oplus \Lambda)=B\otimes_A M\oplus B\otimes_A \Lambda$ est alors
 libre, et~$B\otimes_A M$ est donc projectif, là encore par me théorème ci-dessus.~$\Box$ 
 
 \trois{coro-proj-prin}
 {\bf Corollaire.}
 {\em Soit~$A$ un anneau principal et soit~$M$ un~$A$-module projectif
 de type fini. Le module~$M$ est libre (de rang fini, {\em cf.}~\ref{libre-tf}).}
 
 \medskip
 {\em Démonstration.}
 Le théorème ci-dessus assure qu'il existe un~$A$-module~$\Lambda$
 tel que~$M\oplus\Lambda$ soit libre. Cela assure en particulier que~$M$ est sans torsion ; 
 la théorie générale des modules sur les anneaux principaux garantit alors que~$M$ est
 libre.~$\Box$

 \subsection*{Modules de présentation finie}
 
 \deux{pro-pres-fin}
 Si~$P$ est un module {\em projectif}
 de type fini, il est automatiquement de présentation finie. Pour le voir, on choisit une famille génératrice
 finie~$(p_1,\ldots, p_n)$  de~$P$, et l'on note~$\pi : A^n\to P$
 la surjection $(a_i)\mapsto \sum a_i p_i$. En vertu du théorème~\ref{equi-proj-theo}, 
 on a un isomorphisme~$A^n\simeq P\oplus {\rm Ker}\;\pi$. Modulo cet isomorphisme 
 on a~${\rm Ker}\;\pi=A^n/P$, et~${\rm Ker}\;\pi$ est en particulier de type fini. 
 
 \deux{intro-pres-fini}
 Par définition, si un module~$M$ est de présentation finie,
 {\em il existe}
 une surjection d'un module libre type fini vers~$M$ dont le noyau est lui-même de type fini. 
 Mais ce sera en fait le cas de {\em toute}
 telle surjection, comme le montre la proposition suivante. 
 
 \trois{pro-presfin}
 {\bf Proposition.}
 {\em Soit~$M$ un~$A$-module de présentation finie, 
soit~$L$ un~$A$-module de type fini et soit~$p \colon L\to M$ une
surjection linéaire. Le noyau~$K$ de~$p$ est alors de type fini.}

\medskip
{\em Démonstration.}
Comme~$M$ est de présentation finie, il existe~$n\geq 0$ et une surjection
$q\colon A^n\to M$ dont le noyau est de type fini. Comme~$L$ est de type fini, il existe
$m\geq 0$ et une surjection linéaire~$\pi \colon A^m\to L$. 

\medskip
Comme~$A^n$ est libre, il est projectif ; il 
s'ensuit qu'il existe~$s\colon A^n\to L$ tel que~$p\circ s=q$. 
De même, $A^m$ est projectif et il existe donc~$t\colon A^m\to A^n$ tel que~$q\circ t=p\circ \pi$. 
Le diagramme
$$\xymatrix{
{A^m\oplus A^n}\ar[rr]^{t\oplus {\rm Id}}\ar[d]_{\pi \oplus s}&&{A^n}\ar[d]^q\\
L\ar[rr]^p&&M}$$ est alors commutatif, et les flèches~$\pi\oplus s$ et~$t\oplus {\rm Id}$ 
sont surjectives. Soit~$K'$ l'image réciproque de~$K$
par la surjection~$\pi\oplus s$ ; le module~$K'$
se surjecte sur~$K$, et il suffit dès lors de prouver que~$K'$
est de type fini. 
On déduit du diagramme commutatif

$$\xymatrix{K'\ar[r]\ar[d]&
{A^m\oplus A^n}\ar[rr]^{t\oplus {\rm Id}}\ar[d]_{\pi \oplus s}&&{A^n}\ar[d]^q\\
K\ar[r]&L\ar[rr]^p&&M}$$
et des définitions de~$K$ et~$K'$ 
que~$K'=(t\oplus{\rm Id})^{-1}({\rm Ker}\;q)$. Par hypothèse, 
${\rm Ker}\;q$ est de type fini ; puisque~$t\oplus {\rm Id}$
est surjective, il existe une famille finie~$(f_1,\ldots, f_r)$ d'éléments
de~$A^m\oplus A^n$ dont les images par~$t\oplus {\rm Id}$ engendrent~${\rm Ker}\;q$, 
et il est alors immédiat que
$$K'=(t\oplus{\rm Id})^{-1}({\rm Ker}\;q)=\langle f_1,\ldots, f_r\rangle +{\rm Ker}\;(t\oplus {\rm Id}).$$

\medskip
Le caractère surjectif de~$t\oplus {\rm Id}$ et le caractère projectif de~$A^n$ assurent
que~$A^m\oplus A^n$ est isomorphe à~$A^n\oplus {\rm Ker}\;(t\oplus {\rm Id})$ ; en particulier, 
${\rm Ker}\;(t\oplus {\rm Id})$ est un quotient de~$A^n\oplus A^m$ et est en conséquence
de type fini ; il en résulte que~$K'$ est de type fini.~$\Box$

\trois{rem-pas-pf}
On peut maintenant donner un exemple de module
de type fini qui n'est pas de présentation finie. 
Soit~$A$ un anneau non noethérien soit~$I$
un idéal de~$A$ 
qui n'est pas de type fini
(par exemple, 
on peut prendre~$A=\CC[X_n]_{n\in \NN}$ et~$I=(X_n)_{n\in \NN}$). Le~$A$-module quotient~$A/I$ est de type 
fini (il est engendré par~$\bar 1$) mais n'est pas de présentation finie : sinon, 
le noyau~$I$ de la flèche quotient~$A\to A/I$ serait de type fini d'après la proposition~\ref{pro-presfin}
ci-dessus. 

\trois{coro-descente-pf}
{\bf Corollaire.}
{\em Soit~$(S_i)$ une famille finie et couvrante (\ref{fam-compl-loc}) de parties multiplicatives de~$A$
et soit~$M$ un~$A$-module. Le~$A$-module~$M$ est de présentation finie si et seulement si
le~$S_i^{-1}A$-module~$S_i^{-1}M$ est de présentation finie pour tout~$i$.}

\medskip
{\em Démonstration.}
L'implication directe provient du fait que le caractère «de présentation finie»
est stable par extension des scalaires en vertu de
l'exactitude à droite du produit tensoriel ({\em cf.}~\ref{extenson-scal-gen-rel}
pour le raisonnement précis). 
Supposons maintenant que~$S_i^{-1}M$ est de présentation finie
pour tout~$i$, et montrons qu'il en va de même de~$M$. 

\medskip
 On déduit du lemme~\ref{type-fini-si}
que~$M$ est de type fini. Choisissons une partie génératrice finie~$(m_j)_{1\leq j\leq n}$
de~$M$ ; 
soit~$p \colon A^n\to M$ la surjection~$(a_j)\mapsto \sum a_j m_j$
et soit~$K$ le noyau de~$p$. 

\medskip
Pour tout~$i$, la suite
$$0\to S_i^{-1}K\to (S_i^{-1}A)^n\to S_i^{-1}M\to 0$$ 
induite par~$p$ est exacte par platitude
de~$S_i^{-1}A$. Comme~$S_i^{-1}M$ est de présentation finie par hypothèse, 
$S_i^{-1}K$ est de type fini d'après
la	proposition~\ref{pro-presfin}. Ceci vaut pour tout~$i$ ; en appliquant
une fois encore le lemme~\ref{type-fini-si},  on voit que~$K$ est de type fini,
et donc que~$M$ est de présentation finie.~$\Box$ 

\trois{comment-desc-presf}
{\em Remarque.}
La finitude de la famille~$(S_i)$ est indispensable à la validité du corollaire~\ref{coro-descente-pf},
{\em cf.} la remarque~\ref{rem-tf-loc}. 

\deux{present-fini-mor}
Soient~$M$ et~$N$ deux~$A$-modules et soit~$B$ une~$A$-algèbre. 
Toute application~$A$-linéaire de~$M$ vers~$N$ induit une application~$B$-linéaire
de~$B\otimes_A M$ vers~$B\otimes_A N$ ; on définit
par ce biais une application
$A$-linéaire de~$\hom_A(M,N)$ vers~$\hom_B(B\otimes_A M, B\otimes_A N)$, 
puis une application~$B$-linéaire
$$B\otimes_A\hom_A(M,N)\to \hom_B(B\otimes_A M, B\otimes_A N)$$
par la propriété universelle de l'extension des scalaires. 

\trois{pres-fini-mor-prop}
{\bf Proposition.}
{\em Supposons
que~$M$ est de présentation finie et que~$B$ est plat sur~$A$. La flèche
$$B\otimes_A\hom_A(M,N)\to \hom_B(B\otimes_A M, B\otimes_A N)$$
est alors bijective.}

\medskip
{\em Démonstration.}
Fixons un isomorphisme~$M\simeq L/(\sum_i a_{ij}e_i)_j$, où~$L$ est un module
libre de base finie~$(e_i)_{i\in I}$
et où les~$a_{ij}$ sont des scalaires, les indices~$j$ parcourant un ensemble fini~$J$. 

\medskip
Soit~$C$ une~$A$-algèbre quelconque. 
On a par exactitude à droite du produit tensoriel
un isomorphisme
$$C\otimes_A M\simeq C\otimes_A L/\left(\sum_j a_{ij}\cdot 1\otimes e_i\right)_i,$$
et les~$1\otimes e_i$ forment une base de~$C\otimes_A L$. On en 
déduit un isomorphisme naturel, et fonctoriel en~$C$, entre 
$\hom_C(C\otimes_AM, C\otimes_A N)$ et le noyau de
$$(C\otimes_A N)^I \to (C\otimes_A N)^J,\; \;(n_i)\mapsto \left(\sum_i a_{ij}n_j\right)_j$$
(prendre pour~$n_i$ l'image de~$\overline{1\otimes e_i}$). 

\medskip
En particulier, 
$\hom_A(M,N)$ est le noyau de~$N^I\to N^J, (n_i)\mapsto (\sum_i a_{ij}n_j)_j$,
et~$\hom_B(B\otimes_AM, B\otimes_A N)$ est le noyau de
$$(B\otimes_A N)^I \to (B\otimes_A N)^J,\; \;(n_i)\mapsto \left(\sum_i a_{ij}n_j\right)_j.$$

Le foncteur~$B\otimes \bullet$ commute aux produits finis (qui sont des sommes directes), 
et à la formation du noyau car~$B$ est plat. Il s'ensuit que
$$\hom_B(B\otimes_A M, B\otimes_A N)\simeq B\otimes_A \hom_A(M,N),$$
ce qu'il fallait démontrer.~$\Box$ 

\trois{prop-iso-mp}
{\bf Corollaire.}
{\em Soient~$M$ et~$N$ deux~$A$-modules de présentation finie
et soit~$\got p$ un idéal premier de~$A$ tel que
les~$A_{\got p}$-modules
$M_{\got p}$ et~$N_{\got p}$
soient isomorphes. Il
existe alors~$f\notin \got p$ tel que
les~$A_f$-modules~$M_f$
et~$N_f$
soient isomorphes.}

\medskip
{\em Démonstration.}
Les localisés de~$A$ étant plats sur~$A$, on déduit
de la proposition~\ref{pres-fini-mor-prop}
les égalités
$$\hom_{A_{\got p}}(M_{\got p}, N_{\got p})=\hom_A(M,N)_{\got p}\;\;\text{et}\;
\hom_{A_{\got p}}(N_{\got p}, M_{\got p})=\hom_A(N,M)_{\got p}$$
et
$$\hom_{A_f}(M_f, N_f)=\hom_A(M,N)_f\;\;\text{et}\;
\hom_{A_f}(N_f, M_f)=\hom_A(N,M)_f$$ pour tout~$f\in A$. 
Rappelons par ailleurs que si~$\Lambda$ et un~$A$-module
quelconque,
$\Lambda_{\got p}$ s'identifie à la limite inductive filtrante des~$\Lambda_f$
pour~$f\notin \got p$ (\ref{mp-limind-mf}). 

\medskip
Choisissons un isomorphisme~$u \colon M_{\got p}\to N_{\got p}$ et
soit~$v$ sa réciproque. Par ce qui précède, on peut voir~$u$ et~$v$ comme des éléments
respectifs de~$\hom_A(M,N)_{\got p}$ et~$\hom_A(N,M)_{\got p}$, 
et il existe~$g\in A\setminus \got p$ tel que~$u$ et~$v$ proviennent
respectivement de deux éléments (encore notés~$u$ et~$v$)
de~$\hom_A(M,N)_g=\hom_{A_g}(M_g,N_g)$ et~$\hom_A(N,M)_g=\hom_{A_g}(N_g,M_g)$.

Les égalités~$u\circ v={\rm Id}_{N_{\got p}}$ et~$v\circ u={\rm Id}_{M_{\got p}}$ peuvent être vues
comme des égalités entre éléments de~$\hom_A(M,M)_{\got p}$ et~$\hom_A(N,N)_{\got p}$, 
et il existe donc un élément~$f\in A\setminus \got p$, multiple de~$g$, telle que les
égalités en question vaillent déjà dans~$\hom_A(M,M)_f=\hom_{A_f}(M_f,M_f)$ et~$\hom_A(N,N)_f=\hom_{A_f}(N_f,N_f)$. 
En conséquence, les~$A_f$-modules~$M_f$ et~$N_f$ sont isomorphes.~$\Box$ 

\subsection*{Retour aux modules projectifs}
\deux{lemme-fond-moduleslibres}
{\bf Lemme.}
{\em Soit~$M$ un~$A$-module et soit~$(S_i)_{i\in I}$ une famille couvrante et finie
de parties multiplicatives de~$A$. Les assertions suivantes sont équivalentes : 

\medskip
i) le~$A$-module~$M$ est projectif de type fini ; 

ii) pour tout~$i$,  le~$S_i^{-1}$-module~$S_i^{-1}M$ est projectif
de type fini. }

\medskip
{\em Démonstration.}
L'implication~i$\Rightarrow$ii) provient du fait que le caractère projectif
et le fait d'être de type fini 
sont préservés par extension des scalaires. 

\medskip
Supposons maintenant que~ii) soit vraie ; 
on déduit de~\ref{pro-pres-fin}
que~$S_i^{-1}M$ est pour tout~$i$ un~$S_i^{-1}A$-module
de {\em présentation finie}, 
et le lemme~\ref{type-fini-si} assure alors que~$M$ est de présentation finie.
Choisissons une famille
génératrice finie~$(m_i)_{i\in I}$
de~$M$ et soit~$p$
la surjection~$A^I\to M, (a_i)\mapsto \sum a_i m_i$ ; elle induit une
application~$\hom_A(M, A^I)\to \hom_A(M,M)$, d'où pour tout~$i$
une application
$$S_i^{-1}\hom_A(M,A^I)\to S_i^{-1}\hom_A(M,M)$$
qui s'identifie
d'après la proposition~\ref{pres-fini-mor-prop}
(et le fait que~$M$ est de présentation finie)
à la flèche
$$\hom_{S_i^{-1}A}(S_i^{-1}M, (S_i^{-1}A)^I)\to \hom_{S_i^{-1}A}(S_i^{-1}M,S_i^{-1}M),$$
laquelle est surjective puisque~$S_i^{-1}M$ est projectif. 
On déduit alors du lemme~\ref{lemme-ex-msi}
que~$p$ est surjective. En particulier, ${\rm Id}_M$ a un antécédent~$s$
par~$p$, qui fournit une section de~$p$. Ainsi, $M$ est projectif 
({\em cf.}
la condition~iii) du th.~\ref{equi-proj-theo}).~$\Box$

\deux{theo-proj-fond}
{\bf Théorème.}
{\em Soit~$M$ un~$A$-module.
Les propositions suivantes sont équivalentes : 

\medskip
i) Le module~$M$ est projectif et de type fini. 

ii) Le module~$M$ est de présentation finie, et pour tout idéal premier~$\got p$ de~$A$, 
il existe un entier~$n_{\got p}$ tel que~$M_{\got p}\simeq A_{\got p}^{n_{\got p}}$. 

iii) Il existe une famille finie~$(f_i)$ d'éléments de~$A$ engendrant~$A$ comme idéal, et pour tout~$i$ un 
entier~$n_i$ tel que~$M_{f_i}\simeq A_{f_i}^{n_i}$.}

\medskip
{\em Démonstration.}
Supposons~i) vraie. On sait que~$M$ est de présentation finie d'après~\ref{pro-pres-fin}. Soit~$\got p$
un idéal premier de~$A$. Le~$A_{\got p}$-module~$M_{\got p}$ est projectif et de type fini, 
et est dès lors libre de rang fini d'après
le lemme~\ref{proj-loc}. Ainsi, ii) est vraie. 

\medskip
Supposons~ii) vraie et soit~$\got p$ un idéal premier de~$A$. La 
proposition~\ref{prop-iso-mp}
assure qu'il existe~$f_{\got p}\in A\setminus \got p$
tel que~$M_{\got p}\simeq A_{f_{\got p}}^{n_{\got p}}$. 

Par construction, l'idéal engendré par les~$f_{\got p}$ n'est contenu 
dans aucun idéal premier de~$A$, et il coïncide dès lors avec~$A$. Cela signifie qu'il existe 
une famille~$(a_{\got p})$ de scalaires {\em presque tous nuls}
tels que~$\sum a_{\got p}f_{\got p}=1$. Si~$I$ désigne l'ensemble
des idéaux premiers~$\got p$ tels que~$a_{\got p}\neq 0$, la famille~$(f_i)_{i\in I}$
vérifie les conditions requises par~iii).

\medskip
Supposons~iii) vraie. Pour tout~$i$, le~$A_{f_i}$-module~$M_{f_i}$ est libre de rang fini,
et est en particulier projectif et de présentation finie. Comme la famille~$(\{f_i^n\}_n)_i$ est couvrante
({\em cf.} l'exemple b) de~\ref{exemple-couvrante}), 
le lemme~\ref{lemme-fond-moduleslibres}
assure que~$M$ est projectif de type fini.~$\Box$ 

\deux{ex-dedekind}
{\bf L'exemple des anneaux de Dedekind.}
Il y a différentes caractérisation des anneaux de Dedekind, et celle que nous
utiliserons est la suivante. Un anneau~$A$ est de Dedekind si et seulement si 
il est intègre, noethérien, et possède la propriété suivante : si~$\got p$ 
est un idéal premier non nul de~$A$ alors~$A_{\got p}$ est un anneau de valuation
discrète, c'est-à-dire un anneau principal ayant (à équivalence près) un unique élément irréductible~$\pi$, qui est alors
nécessairement le générateur de son unique idéal maximal. 

\deux{prop-dedek}
{\bf Proposition.}
{\em Soit~$A$ un anneau de Dedekind et soit~$K$ son corps des fractions. Soit~$M$ un~$A$-module. Les assertions suivantes sont équivalentes : 

\medskip
i) $M$ est projectif de type fini et~$M_{\got p}$ est libre de rang~$1$ sur~$A_{\got p}$ pour tout idéal
premier~$\got p$ de~$A$. 

ii) $M$ est projectif de type fini et~$K\otimes_A M$ est de dimension~$1$ sur~$K$. 

iii) $M$ est isomorphe comme~$A$-module à un idéal non nul de~$A$. 

\medskip
Par ailleurs un idéal de~$A$ est libre comme~$A$-module si et seulement si il est principal.}

\medskip
{\em Démonstration.}
Commençons par une remarque que nous allons utiliser constamment : si~$S$ est une partie
multiplicative de~$A$ ne contenant pas~$0$ et si~$N$ est un~$A$-module sans torsion,
la flèche naturelle~$N\to S^{-1}N=S^{-1}_A \otimes_A N$ est injective (c'est une conséquence triviale
de la condition de nullité d'une fraction). 

\medskip
Si~i) est vraie, ii) est vraie : prendre~$\got p=\{0\}$. Supposons maintenant que~ii)
est vraie. Comme~$M$ est projectif
il est facteur direct d'un module libre, et en particulier sans torsion. En conséquence, $M$ s'injecte dans le~$K$-espace vectoriel~$K\otimes_AM$ 
qui est de dimension~$1$, ce qui fournit (une fois choisie une base de~$K\otimes_AM$) un plongement~$A$-linéaire de~$M$ dans~$K$, donc une identification de~$M$
à un sous-$A$-module de~$K$. Comme~$M$ est de type fini, il est engendré par un nombre fini de fractions ; soit~$\alpha$ un multiple commun non nul de leurs dénominateurs. 
Le sous~$A$-module~$\alpha M$ de~$K$ est isomorphe à~$M$, et contenu dans~$A$ ; c'est donc un idéal de~$A$, 
nécessairement non nul puisque~$K\otimes_A M$ est non nul, d'où~iii). 

\medskip
Montrons maintenant~iii)$\Rightarrow$i). On peut supposer que~$M$ est un idéal non nul de~$A$ ; 
il est donc de présentation finie par noethérianité de~$A$. Soit~$\got p$ un idéal premier de~$A$. 
L'injection~$M\hookrightarrow A$ induit par platitude de~$A_{\got p}$ 
une injection~$M_{\got p}\hookrightarrow A_{\got p}$, et~$M_{\got p}$ est par ailleurs non nul
puisque~$M$ est non nul et s'injecte dans~$M_{\got p}$ (étant contenu dans
l'anneau intègre~$A$, le module~$M$ est sans torsion). 

En conséquence~$M_{\got p}$ s'identifie pour tout~$\got p$ à un idéal non nul
de l'anneau~$A_{\got p}$, lequel est principal (c'est évident si~$\got p=\{0\}$, et c'est dû au fait
que~$A$ est un anneau de Dedekind sinon). Il s'ensuit que le~$A_{\got p}$-module~$M_{\got p}$ est libre de rang~$1$. 
On en déduit alors~i) grâce au théorème~\ref{theo-proj-fond}
dont la condition~ii) est ici vérifiée.

\medskip
Il reste à justifier la dernière assertion. Si~$M$ est un idéal principal de~$A$ il est 
engendré par un élément~$a$ et est donc libre, de rang~$0$ si~$a=0$ et de rang~$1$ sinon. 
Réciproquement, supposons que~$M$ soit un idéal de~$A$, et qu'il soit libre comme~$A$-module. 
Si~$a$ et~$b$ sont deux éléments non nuls de~$M$, on a l'égalité~$b\cdot a-a\cdot b=0$, ce qui exclut qu'ils puissent
figurer tous deux dans une même famille libre. En conséquence, le rang de
$M$
est égal
à zéro
(auquel cas~$M$ est nul)
ou à un (auquel cas~$M$ est engendré par un élément non nul de~$A$). Dans tous les cas, $M$
est principal.~$\Box$

\section{Produit tensoriel de deux algèbres}
 \markboth{Algèbre commutative}{Produit tensoriel de deux algèbres}

\subsection*{Définition, exemples, premières propriétés}

\deux{intro-b-tensa-c}
On désigne toujours par~$A$ un anneau. 
Soient~$B$ et~$C$ deux~$A$-algèbres. Nous allons démontrer qu'il existe une unique
loi interne~$\cdot$ 
sur $B\otimes_AC$  telle que
~$(B\otimes_AC,+,\cdot)$ soit
un anneau et telle que

 $$(b\otimes c)\cdot (\beta \otimes \gamma)=(b\beta)\otimes (c\gamma)$$
 pour tout~$(b,\beta,c,\gamma)\in B^2\times C^2$. 
 
 \medskip
 \trois{unique-bac}
 {\em Unicité.}
 Elle résulte du fait que les tenseurs purs engendrent
 $B\otimes_A C$ comme groupe abélien. 
 
 \trois{exist-bac}
 {\em Existence}. Soit~$(b,c)\in B\times C$. L'application
 de~$B\times C$ vers~$B\otimes_AC$ qui envoie~$(\beta,\gamma)$ sur~$b\beta\otimes c\gamma$
 est bilinéaire. Elle induit donc une application~$A$-linéaire
 $m_{b,c}$ de~$B\otimes_A C$ dans lui-même. On vérifie aussitôt
 que~$(b,c)\mapsto m_{b,c}$ est
 elle-même~$A$-bilinéaire. Elle induit donc une application~$A$-linéaire~$\chi$ 
 de~$B\otimes_AC$ dans~${\rm End}_A(B\otimes_AC)$. 
 
 \medskip
 On définit une loi interne sur~$B\otimes_AC$ par la formule
 $(v,w)\mapsto v\cdot w:=\chi(v)(w)$. Il découle de sa construction 
 qu'elle est bi-$A$-linéaire et
 satisfait les  égalités
 $$(b\otimes c)\cdot (\beta\otimes\gamma)=b\beta\otimes c\gamma$$ 
 pour tout~$(b,\beta,c,\gamma)\in B^2\times C^2$. On déduit
 sans difficulté de ces formules
 que~$(B\otimes_AC,+,\cdot)$ est un anneau.

 \deux{prop-elem-btensc}
 Par définition de la loi~$\cdot$
 sur~$B\otimes_A C$, les applications~$b\mapsto b\otimes 1$ et~$c\mapsto 1\otimes c$ sont des
 morphismes d'anneaux, et les composées~$A\to B\to B\otimes_AC$
 et~$A\to C\to B\otimes_AC$ coïncident : par~$A$-bilinéarité de~$\otimes$,
 on a en effet~$a\otimes 1=1\otimes a$ pour tout~$a\in A$ (comme d'habitude, on 
 note encore~$a$ les images de~$a$ dans $B$ et~$C$). 
 
 \medskip
 L'anneau~$B\otimes_AC$ hérite ainsi d'une structure de~$B$-algèbre et d'une structure
 de~$C$-algèbre, qui induisent la même structure de~$A$-algèbre. On vérifie que ses structures
 de~$A$-module,~$B$-module et~$C$-module sont précisément les structures sous-jacentes 
 à ces structures d'algèbre. 
 
 \deux{intuition-btensc}
 Intuitivement,~$B\otimes_AC$ est la~$A$-algèbre la plus
 générale fabriquée à partir de~$B$ et~$C$, en définissant «artificiellement» la
 multiplication d'un élément~$b$ de~$B$ par un élément~$c$ de~$C$ (c'est~$b\otimes c$).
 Comme toujours, 
 ce type de description se traduit rigoureusement en termes de propriété universelle, 
 ou de foncteur représenté.
 
 \deux{prop-univ-btensc}
 {\bf Proposition.}
 {\em Le couple $(B\otimes_AC, (b\mapsto b\otimes 1, c\mapsto 1\otimes c))$
 fait de~$B\otimes_AC$ la somme disjointe de~$B$ et~$C$ dans la catégorie des~$A$-algèbres.}
 
 \trois{tauto-a-tens-b}
 {\bf Remarque.}
 Il revient {\em tautologiquement}
 au même d'affirmer que~$(B\otimes_AC, (b\mapsto b\otimes 1, c\mapsto 1\otimes c))$
 fait de~$B\otimes_AC$ la somme amalgamée de~$B$ et~$C$ le long de~$A$ dans la
 catégorie des anneaux : c'est la version duale de~\ref{fib-cart-ca}. 
 
 \trois{demo-propuniv-bc}
 {\em Démonstration de la proposition~\ref{prop-univ-btensc}}.
 Soit~$D$ une~$A$-algèbre
 et soient $f: B\to D$ et~$g: C\to D$ deux morphismes de~$A$-algèbres. Il s'agit
 de montrer qu'il existe un unique morphisme de~$A$-algèbres~$h$
 de~$B\otimes_AC$
 vers~$D$ tel que le diagramme
 
 $$
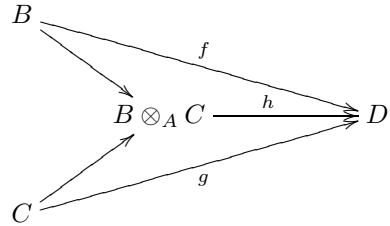
\xymatrix{B\ar[rd]\ar[rrrd]^f&&&\\&B\otimes_AC\ar[rr]^h&&D
\\C\ar[ru]\ar[rrru]_g&&&}
$$ 
commute (ou encore un unique morphisme
d'anneaux de~$B\otimes_AC$ dans~$D$ qui soit à la
fois un morphisme de~$B$-algèbres et de~$C$-algèbres).

\medskip
{\em Unicité.}
Si~$h$ existe, on a nécessairement pour tout~$(b,c)\in B\times C$
les égalités
$$h(b\otimes c)=h((b\otimes 1) \cdot (1\otimes c))=h(b\otimes 1) h(1\otimes c)=f(b)g(c),$$
et l'unicité de~$h$ découle du fait que les tenseurs purs engendrent le groupe abélien~$B\otimes_A C$. 

\medskip
{\em Existence.}
On s'inspire de la seule formule possible obtenue en démontrant l'unicité. L'application de~$B\times C$
dans~$D$ qui envoie~$(b,c)$ sur~$f(b)g(c)$ est bi-$A$-linéaire, et elle induit donc une application
$A$-linéaire~$h : B\otimes_AC\to D$, qui satisfait par construction les égalités
$h(b\otimes c)=f(b)g(c)$ pour tout~$(b,c)\in B\times C$. On déduit de celles-ci que~$h$ est un morphisme
de~$A$-algèbres répondant aux conditions posées.~$\Box$

\deux{philo-btenseurc}
Selon le contexte, il y a plusieurs façons d'envisager~$B\otimes_AC$. On peut
y penser comme à un objet {\em symétrique en~$B$ et~$C$} : c'est 
par exemple le cas lorsqu'on le décrit 
informellement comme en~\ref{intuition-btensc}
ou, plus rigoureusement, lorsqu'on le caractérise
par le foncteur qu'il représente (prop.~\ref{prop-univ-btensc}). 

Mais on peut faire psychologiquement
jouer un rôle différent à~$B$ et~$C$, en considérant qu'on part d'une~$A$-algèbre~$C$
et qu'on la transforme en une~$B$-algèbre~$B\otimes_AC$ (ou l'inverse, évidemment). Le slogan
à retenir lorsqu'on aborde les choses de ce point de vue est,
à analogue
à celui vu plus haut pour les modules : 
{\em la $B$-algèbre~$B\otimes_AC$ est à~$B$
ce que~$C$ est à~$A$.}  Nous allons l'illustrer par différents
exemples. 

\deux{b-tenseur-at-intro}
On désigne toujours par~$B$ une~$A$-algèbre ; soit~$I$ un ensemble
d'indices. 

\trois{b-tenseur-at}
La donnée des deux morphismes naturels
de~$A$-algèbres
$$B\to B[T_i]_{i\in I}\;{\rm et}\;A[T_i]_{i\in I}\to B[T_i]_{i\in I}$$
induit un morphisme~$\phi$
de~$B\otimes_A A[T_i]_{i\in I}$ vers
$B[T_i]_{i\in I}$, qui est à la fois un morphisme
de~$B$-algèbres et un morphisme de~$A[T_i]_{i\in I}$-algèbres. 

\medskip
La propriété universelle de la~$B$-algèbre~$B[T_i]_{i\in I}$ assure
par ailleurs l'existence d'un unique morphisme
de~$B$-algèbres $\psi \colon B[T_i]_{i\in I}
\to B\otimes_A A[T_i]_{i\in I}$ qui envoie~$T_i$
sur~$1\otimes T_i$ pour tout~$i$. On vérifie aussitôt que~$\phi$ et~$\psi$ sont inverses
l'un de l'autre. 

\medskip
On a donc construit un isomorphisme 
$$B\otimes_A A[T_i]_{i\in I}\simeq B[T_i]_{i\in I},$$ compatible
aux structures de~$B$-algèbres et de~$A[T_i]_{i\in I}$-algèbres sur
ses source et but. 

\trois{b-tenseur-at-alter}
Donnons une autre construction de ces isomorphismes. Le~$A$-module~$A[T_i]_{i\in I}$
est libre de base~$(\prod_{i\in I}T_i^{e(i)})_{e}$, où~$e$ parcourt l'ensemble des applications
de~$I$ dans~$\NN$ s'annulant presque partout. 

Le~$B$-module~$B\otimes_A A[T_i]_{i\in I}$ est donc
libre
de base~$(1\otimes \prod_{i\in I}T_i^{e(i)})_{e}$. Compte-tenu de la définition 
de la loi d'anneau sur~$B\otimes_A A[T_i]_{i\in I}$, on a par ailleurs
$1\otimes \prod_{i\in I}T_i^{e(i)}=\prod_{i\in I} (1\otimes T_i)^{e(i)}$ pour tout~$e$. Ainsi,
$B\otimes_A A[T_i]_{i\in I}$ est une algèbre de polynômes
en les~$1\otimes T_i$ : on retrouve donc
l'isomorphisme
du~\ref{b-tenseur-at}. 

\trois{comment-sym-dissym}
En particulier~$A[S]\otimes_A A[T]\simeq (A[S])[T]=A[S,T]=(A[T])[S]$. 
On voit bien sur
cet exemple les différentes façons dont on peut penser au produit tensoriel : la
première écriture est à~$A[S]$ ce que~$A[T]$ est à~$A$, la seconde est symétrique
en les facteurs, la troisième est à~$A[T]$ ce que~$A[S]$ est à~$A$. 

\deux{modules-algebres}
Certains isomorphismes de modules
exhibés lors de l'étude du produit tensoriel
d'un module par une algèbre se trouvent en fait, lorsque le
module en jeu est lui-même une algèbre, être des isomorphismes
d'algèbres -- on le vérifie immédiatement à l'aide des formules
explicites qui les décrivent. Donnons deux exemples. 

\trois{asuri-tens-b}
Soit~$I$ un idéal de~$A$ et soit~$B$ une~$A$-algèbre. L'isomorphisme

$$A/I\otimes_A B\simeq B/IB$$ ({\em cf.} \ref{ext-scal-quot})
est alors un isomorphisme de~$B$-algèbres
et de~$A/I$-algèbres. 

\trois{tens-alg-localis}
Soit~$S$ une partie multiplicative de~$A$ et soit~$B$ une~$A$-algèbre ; soit~$\phi \colon B\to A$
le morphisme structural. 
La notation~$S^{-1}B$ est {\em a priori}
ambiguë : elle pourrait désigner le localisé 
du~$A$-module~$B$ par la partie multiplicative~$S$ de~$A$
ou bien, modulo notre abus usuel consistant à oublier~$\phi$, l'anneau
localisé de~$B$ par rapport à sa partie multiplicative~$\phi(S)$. 
Mais {\em a posteriori}, il n'y a aucun problème : nous laissons le 
lecteur vérifier qu'il existe un isomorphisme de~$S^{-1}A$-modules
entre le premier et le second de ces objets, donné par la formule
$$\frac b s \mapsto \frac b {\phi(s)}\;\;;$$ 
et que modulo celui-ci,
l'isomorphisme
$$S^{-1}A\otimes_A B\simeq S^{-1}B$$
({\em cf.} \ref{ext-scal-localis}) est un isomorphisme
de~$B$-algèbres et de~$S^{-1}A$-algèbres. 

\trois{b-tensc-tensd}
Soient~$B$ et~$C$ deux~$A$-algèbres, et soit~$D$ une~$B$-algèbre. L'isomorphisme 

$$D\otimes_B(B\otimes_A C)\simeq D\otimes_AC$$
({\em cf.}~\ref{m-tensb-tensc})
est alors un isomorphisme de~$D$-algèbres et de~$C$-algèbres. 

\trois{b-tens-csuri}
En vertu de~\ref{asuri-tens-b}
et~\ref{b-tensc-tensd}, il existe
pour toute~$A$-algèbre~$B$, toute~$A$-algèbre~$C$, 
et tout idéal~$I$ de~$C$ des isomorphismes naturels

$$B\otimes_AC/I\simeq (B\otimes_A C)\otimes_C (C/I) \simeq (B\otimes_AC)/I(B\otimes_AC).$$

\trois{b-tenseur-atsurp}
En vertu de~\ref{b-tens-csuri}
et~\ref{b-tenseur-at}, il existe pour toute~$A$-algèbre~$B$, 
tout ensemble d'indices~$I$ et toute famille~$(P_j)$
de polynômes appartenant à~$A[T_i]_{i\in I}$ un isomorphisme naturel
(de~$B$-algèbres aussi bien que de~$A[T_i]_{i\in I}$-algèbres)

$$B\otimes_A (A[T_i]_{i\in I}/(P_j))\simeq B[T_i]_{i\in I}/(P_j),$$
(où l'on note encore~$P_j$ l'image de~$P_j$ dans~$B[T_i]_{i\in I}$). 

\medskip
De manière un peu informelle, 
on voit que la~$A$-algèbre~$A[T_i]_{i\in I}/(P_j)$
et la~$B$-algèbre~$B\otimes_A(A[T_i]_{i\in I}/(P_j))$ 
admettent la
«même» description par générateurs (les~$T_i$) et relations
(les~$P_j$).

\trois{b-tenseur-atsurp}
{\em Exercice.}
Construire directement l'isomorphisme ci-dessus par une méthode analogue
à celle suivie au~\ref{b-tenseur-at}.

\deux{ex-tenseur-alg}
{\bf Exemples.}

\trois{zx-modp}
Soit~$A$ la~$\ZZ$-algèbre~$\ZZ[X]/(6X^2+18X -3)$. Pour toute~$\ZZ$-algèbre~$B$,
on a~$B\otimes_\ZZ A\simeq B[X]/(6X^2+18X-3)$. L'allure de cette dernière~$B$-algèbre
dépend beaucoup de~$B$. 
Ainsi : 

\medskip
$\bullet$ si~$B=\QQ$, elle est égale à~$\QQ[X]/(6X^2+18X-3)$ qui est un corps 
de degré 2 sur~$Q$, car~$6X^2+18X-3$ est irréductible sur~$\QQ[X]$ (son discriminant
est~$324+72=396=4\times 9\times 11$ qui n'est pas un carré dans~$\QQ$) ; 

$\bullet$ si~$B=\FF_2$ elle est égale~à
$\FF_2[X]/(3)=\{0\}$ car~$3$ est inversible modulo~$2$ ; 

$\bullet$ si~$B=\FF_3$ elle est égale à~$\FF_3[X]/(0)=\FF_3[X]$ ; 

$\bullet$ si~$B=\FF_5$ elle est égale à
$$\FF_5[X]/(X^2-2X+2)=\FF_5[X]/(X+1)(X+2)$$
$$\simeq \FF_5[X]/(X+1)\times \FF_5[X]/(X+2)\simeq \FF_5\times \FF_5; $$

$\bullet$ si~$B=\FF_{11}$, elle est égale à
$$\FF_{11}[X]/(6X^2-4X-3)=\FF_{11}[X]/(2(6X^2-4X-3))=\FF_{11}[X]/(X^2+3X-6)$$
$$=\FF_{11}[X](X-4)^2\simeq \FF_{11}[Y]/Y^2.$$

\medskip
On voit qu'en tensorisant
 la même~$\ZZ$-algèbre par différents corps on a
obtenu un corps, l'anneau nul, un anneau de polynômes, 
un produit de deux corps, et un anneau non réduit.

\trois{c-tenseur-c}
Nous  allons maintenant 
décrire l'anneau~$\CC\otimes_{\RR}\CC$. 

\medskip
Comme~$\CC\simeq \RR[X]/(X^2+1)$, 
cet anneau s'identifie à 
$$\CC[X]/(X^2+1)\simeq \CC[X]/(X-i)(X+i)\simeq \CC[X]/(X-i)\times \CC[X]/(X+i)\simeq \CC\times \CC.$$

À titre d'exercice, vérifiez que l'isomorphisme~$\CC\otimes_{\RR}\CC\simeq \CC\times\CC$ ainsi construit
envoie~$b\otimes \beta$ sur~$(b\beta,b\overline \beta)$. 

\medskip
On voit à travers cet exemple qu'un produit tensoriel de deux
corps au-dessus d'un troisième n'est pas nécessairement un corps, ni même un anneau intègre. 
Nous allons voir qu'il peut même arriver qu'un tel produit tensoriel ne soit pas réduit.

\trois{kp-tenseur-kp}
Soit~$k$ un corps de caractéristique~$p>0$ non parfait, 
et soit~$a$ un élément de~$k$ qui n'est pas une puissance~$p$-ième. On démontre
(nous laissons la vérification 
au lecteur à titre d'exercice) que~$X^p-a$ est alors un polynôme irréductible. Soit~$L$ le corps~$k[X]/(X^p-a)$
et soit~$\alpha$ la classe de~$X$ dans~$L$. 

On a~$$L\otimes_kL=L\otimes_kk[X]/(X^p-a)\simeq L[X]/(X^p-a)= L[X]/(X^p-\alpha^p)=L[X]/(X-\alpha)^p$$
car l'élévation à la puissance~$p$ est un morphisme d'anneaux en caractéristique~$p$. La classe de~$X-\alpha$
modulo~$(X-\alpha)^p$ fournit alors un élément nilpotent non nul de~$L\otimes_kL$. 

\subsection*{Limites inductives dans la catégorie des anneaux}

\deux{prod-tens-alg-fini}
Soit~$(A_i)_{i\in I}$ une famille de~$A$-algèbres. 
On vérifie facilement, par un raisonnement analogue à celui tenu 
aux~\ref{intro-b-tensa-c} 
{\em et sq.}, qu'il existe une unique structure de~$A$-algèbre sur 
le $A$-module~$\bigotimes_{i\in I} A_i$ (\ref{prod-tens-ensfini})
telle que
$$(\otimes_ix_i)\cdot (\otimes_iy_i)=\otimes_ix_iy_i$$ pour tout
couple~$((x_i), (y_i))$ d'éléments de~$\prod A_i$. 

\trois{limin-alg-ifini}
Il est facile
de voir que cette construction 
fournit la somme disjointe des~$A$-algèbres $A_i$ {\em lorsque~$I$
est fini} (moralement, la raison de cette restriction est qu'on ne sait pas quel sens
donner à un produit infini d'éléments dans un anneau, sauf s'ils sont presque tous égaux à 1). 

\trois{limind-alg-inonfini}
Pour construire la somme disjointe 
des~$A_i$ lorsque l'ensemble~$I$ est quelconque, on
procède comme suit. Si~$J$ et~$J'$ sont deux sous-ensembles
finis de~$I$, on note~$f_{JJ'}$ le morphisme
de~$A$-algèbres de~$\bigotimes_{i\in J}A_i$ vers~$\bigotimes_{i\in J'}A_i$ induit par les flèches structurales
$A_i\to \bigotimes_{i\in J'}A_i$ lorsque~$i$ parcourt~$J$. Si~$\otimes_{i\in J} x_i$
est un tenseur pur de~$\bigotimes_{i\in J}A_i$ alors~$f_{JJ'}(\otimes_{i\in J} x_i)=\otimes_{i\in J'} y_i$
où~$y_i=x_i$ si~$i\in I$ et~$y_i=1$ sinon. 

\medskip
L'ensemble~$E$ des parties finies de~$I$ ordonné par l'inclusion est filtrant, 
et
$$\sch D=((A_J)_{J\in E}, (f_{JJ'})_{J\subset J'})$$ est un diagramme commutatif
filtrant dans la catégorie des~$A$-algèbres. On montre sans peine quels~$\limind \sch D$
est la somme disjointe des~$A_i$ dans la catégories des~$A$-algèbres. 

\trois{exercice-prodtens-infini}
{\em Exercice.}
Vérifiez que le sous-module de~$\bigotimes_{i\in I}A_i$
engendré par les tenseurs purs dont presque toutes les composantes sont égales à~$1$
en est une sous-algèbre, et qu'elle s'identifie canoniquement à~$\limind \sch D$. 

\trois{limind-ann-infini}
{\em Changement de notation.}
On choisit désormais de noter~$\bigotimes_{i\in I}A_i$ la somme disjointe des~$A_i$ dans la catégorie 
des~$A$-algèbres -- cette notation est compatible avec la précédente si~$I$ est fini
ou si presque tous les~$A_i$ sont nuls, mais ne l'est pas sinon.  

\deux{lim-ind-ai}
Soit maintenant~$\Delta=((A_i), (E_{ij}))$ un diagramme 
dans la catégorie des~$A$-algèbres ; pour tout~$i$, soit~$\lambda_i \colon A_i\to \bigotimes_{i\in I} A_i$
le morphisme canonique. 

\trois{desc-limind-ai}Nous laissons le lecteur vérifier que 
$\limind \Delta$ existe et s'identifie au quotient de~$ \bigotimes_{i\in I} A_i$
par son idéal engendré par les~$\lambda_i(x)-\lambda_j(f(x))$ où~$(i,j)\in I^2$, où~$x\in A_i$
et où~$f\in E_{ij}$. 

\trois{desc-limind-anneaux}
{\em Remarque.}
On sait en donc
particulier construire les limites inductives quelconques dans la catégorie des~$\ZZ$-algèbres,
c'est-à-dire des anneaux. 

\trois{revisit-btensc}
Soient~$B$ et~$C$ deux~$A$-algèbres. Le produit tensoriel~$B\otimes_A C$ est la
somme amalgamée de~$B$ et~$C$ le long de~$A$ dans la catégorie des anneaux. Il
résulte de~\ref{lim-ind-ai} (et de la remarque~\ref{desc-limind-anneaux})
que~$B\otimes_AC$ s'identfie
au quotient
$$(B\otimes_{\ZZ}A\otimes_{\ZZ}C)/(a\otimes 1\otimes 1-1\otimes a\otimes 1, 1\otimes a\otimes 1-1\otimes 1\otimes a)_{a\in A}.$$
Vérifiez-le directement à titre d'exercice. 

\subsection*{Applications 
à la théorie des corps}

\deux{prop-ext-compose}
{\bf Proposition}.
{\em Soient~$k\hookrightarrow F$ 
et~$k\hookrightarrow L$
deux extensions de corps. Il existe alors un corps~$K$ et deux plongements~$F\hookrightarrow K$
et~$L\hookrightarrow K$ tels que le diagramme
$$
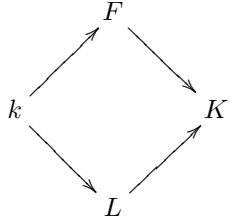
\xymatrix{&F\ar[dr]&\\k\ar[ru]\ar[rd]&&K\\
&L\ar[ur]&}$$ commute.}

\medskip
{\em Démonstration.}
Les~$k$-espaces vectoriels~$F$ et~$L$ sont non nuls, puisque ce sont des corps. 
Comme ils sont libres sur~$k$
(c'est le cas de tout espace vectoriel),
leur produit tensoriel est {\em non nul.} 
La~$k$-algèbre~$F\otimes_kL$ étant non nulle, 
elle possède un idéal maximal. Si l'on note~$K$ le corps quotient correspondant, les flèches
composées~$F\to F\otimes_kL\to K$ 
et~$L\to F\otimes_kL\to K$ satisfont les conditions requises.~$\Box$

\deux{corps-bourbaki}Indiquons maintenant
deux applications de ce fait à la théorie des extensions de corps.

\trois{corps-dec}
{\bf Unicité du corps de décomposition.} Soit~$P$ un polynôme non nul à coefficients
dans~$k$, et soient~$F$ et~$L$ deux corps de décompositions de~$P$ sur~$k$ 
(un corps de décomposition de~$P$ sur~$k$ est une extension de~$k$ dans
laquelle~$P$ est scindé, et qui est engendrée par les racines de~$P$).

{\em Les corps~$F$ et~$L$ sont alors~$k$-isomorphes}. En effet, par ce qui précède, il existe une
extension~$K$ de~$k$ et deux~$k$-plongements~$i : F\hookrightarrow K$ et~$j: L\hookrightarrow K$. 
L'image~$i(F)$ est isomorphe à~$F$, et est donc un corps de décomposition de~$P$ sur~$k$. Il s'ensuit
que~$P$ est scindé dans~$K$ et que~$i(F)$ est le sous-corps de~$K$ engendré par~$k$ et les racines de~$P$. 

De même,~$j(L)$ est le sous-corps de~$K$ engendré par~$k$ et les racines de~$P$. Par conséquent,~$j(L)=i(K)$. 
On a donc deux~$k$-isomorphismes~$F\simeq i(F)$ et~$i(F)=j(L)\simeq L$, d'où un~$k$-isomorphisme~$F\simeq L$.

\trois{clot-alg}
{\bf Unicité de la clôture algébrique.}
Soient~$F$ et~$L$ deux clôtures
algébriques de~$k$
(une clôture algébrique de~$k$ est une extension algébrique
de~$k$ qui est algébriquement close).

{\em Les corps~$F$ et~$L$ sont alors~$k$-isomorphes}. En effet, par ce qui précède, il existe une
extension~$K$ de~$k$ et deux~$k$-plongements~$i : F\hookrightarrow K$ et~$j: L\hookrightarrow K$. 
L'image~$i(F)$ est isomorphe à~$F$, et est donc une clôture algébrique de~$k$. Il s'ensuit
que~$i(F)$ est nécessairement le sous-corps de~$K$ formé des éléments algébriques sur~$k$.  

De même,~$j(L)$ est le sous-corps de~$K$ formé des éléments algébriques sur~$k$.  
Par conséquent,~$j(L)=i(K)$. 
On a donc deux~$k$-isomorphismes~$F\simeq i(F)$ et~$i(F)=j(L)\simeq L$, d'où un~$k$-isomorphisme~$F\simeq L$. 

\deux{rem-corps-bourbaki}
{\bf Remarque.}
Vous connaissez peut-être des preuves de l'unicité du corps de décomposition
d'un polynôme~$P$ consistant à construire l'isomorphisme entre deux tels corps~$F$
et~$L$ 
en choisissant successivement des racines de diviseurs convenables de~$P$ dans~$F$ et~$L$. On peut se demander
où est passé ce choix de racines dans la preuve proposée ci-dessus, qui peut
donner l'impression que «rien ne se passe». En fait, tout 
est caché dans le choix de l'idéal maximal de~$F\otimes_kL$, effectué lorsqu'on veut exhiber un corps~$K$. 

\medskip
La même chose se produit pour les clôtures algébriques : se donner un isomorphisme entre deux clôtures
algébriques~$F$ et~$L$ de~$k$ revient peu ou prou à faire des choix compatibles de racines, dans~$F$
et~$L$, de {\em tous} les polynômes irréductibles de~$k[X]$. Là encore, ce choix est pudiquement dissimulé
derrière celui de l'idéal maximal de~$F\otimes_kL$.

\section{Algèbres finies et algèbres entières}
 \markboth{Algèbre commutative}{Algèbres finies et entières}

On fixe un anneau~$A$. 

\subsection*{Définitions, exemples, premières propriétés}

\deux{def-alg-ent}
{\bf Définition.}
Soit~$B$ une~$A$-algèbre. On dit que~$B$ est {\em finie} si~$B$ est de type fini {\em comme~$A$-module}.

\deux{alg-ent-sorites}
{\bf Exemples et premières propriétés.}

\trois{asuri-fini} Si~$I$ est un idéal de~$A$ alors~$A/I$ est engendré par~$\bar 1$
comme~$A$-module. C'est donc une~$A$-algèbre finie. 

\trois{fini-impl-tf}
Une~$A$-algèbre finie est de type fini comme~$A$-module ; elle l'est {\em a fortiori}
comme~$A$-algèbre. 

\trois{fin-trans}
Soit~$B$ une~$A$-algèbre finie et soit~$C$ une~$B$-algèbre finie. La~$A$-algèbre~$C$ est alors finie ; nous laissons
au lecteur le soin de rédiger la preuve, qui repose essentiellement sur un «principe de la {\em famille génératrice}
télescopique».

\deux{equiv-entier}
{\bf Proposition-définition.} {\em Soit~$B$ une~$A$-algèbre et soit~$x\in B$. Les assertions suivantes sont équivalentes.

\medskip
i) L'élément~ $x$ annule un polynôme {\em unitaire} appartenant à
~$A[X]$. 

ii) La~$A$-algèbre~$A[x]$ est finie.

iii) Il existe une sous-$A$-algèbre finie~$C$ de~$B$ contenant~$x$. 

\medskip
Lorsqu'elles sont satisfaites, on dit que~$x$ est {\em entier} sur~$A$. Si tout élément
de~$B$ est entier sur~$A$, on dit que~$B$ est {\em entière} sur~$A$.} 

\medskip
{\em Démonstration.}
Supposons que~i) soit vraie ; il existe alors~$n\geq 0$ et~$a_0,\ldots, a_{n-1}\in A$
tels que
$$x^n+a_{n-1}x^{n-1}+\ldots+a_0=0.$$
Si~$M$ désigne le sous-$A$-module de~$B$ engendré par les~$x^i$
pour~$0\leq i\leq n-1$, on en déduit que~$x^n\in M$ puis, par récurrence, 
que~$x^m\in M$ pour tout~$m\geq n$. Ainsi, l'algèbre~$A[x]$ coïncide 
avec le~$A$-module de type fini~$M$, et elle est en conséquence finie, 
ce qui prouve~ii). 

\medskip
Si~ii) est vraie, il est clair que~iii) est vraie (prendre~$C=A[x]$).

\medskip
Supposons que~iii) soit vraie, et soit~$u : C\to C$ la multiplication
par~$x$. Comme~$C$ est un~$A$-module de type fini, 
le théorème~\ref{poly-annul-alin} -- appliqué ici avec~$I=A$--
assure 
l'existence d'un polynôme unitaire~$P$ à coefficients dans~$A$
(dont le degré peut être choisi égal au cardinal d'une famille
génératrice
finie
fixée du~$A$-module~$C$) tel que~$P(u)=0$. On a en particulier~$P(x)=P(u)(1)=0$, et~i)
est vraie.~$\Box$

\deux{rem-equiv-ent}
{\bf Quelques remarques.} 

\trois{fini-impl-ent}
On déduit de la caractérisation des éléments
entiers par la propriété~iii) ci-dessus que toute algèbre finie est entière. Nous allons voir
ci-dessous que la réciproque est vraie pour les algèbres de type fini, mais fausse en général. 

\trois{deg-poly-ent} Au cours de la preuve de~iii)$\Rightarrow$i), on
a vu que si~$C$ est engendré comme~$A$-module
par~$n$ éléments,
alors on peut trouver un polynôme unitaire de~$A[X]$ annulant~$x$ {\em et de degré~$n$}. 

\trois{sorites-ent-mor}
Soit~$ f: B\to C$ un morphisme de~$A$-algèbres, et soit~$x$ un élément
de~$B$
entier sur~$A$. Il est immédiat que~$f(x)$ est entier sur~$A$ aussi 
(tout polynôme de~$A[X]$ annulant~$x$ annule~$f(x)$). 

\trois{ent-et-quot}
Soit~$I$ un idéal de~$A$
et soit~$B$ une~$A/I$-algèbre ; on peut aussi voir~$B$
comme une~$A$-algèbre. Un élément de~$B$ est entier sur~$A$
si et seulement si il est entier sur~$A/I$, et~$B$ est finie sur~$A$ si et seulement si elle est
finie sur~$A/I$ (cela provient du fait que si~$\overline a$ est la classe
modulo~$I$ d'un élément~$a$ de~$A$
alors~$\bar ab=ab$ pour tout~$b\in B$). 

\trois{ent-et-doublequot}
Soit~$B$ une~$A$-algèbre, soit~$J$ un idéal de~$B$ et soit~$I$ un idéal de~$A$
dont l'image dans~$B$ est contenue dans~$J$. Si~$x$ est un élément de~$B$ entier
sur~$A$, alors son image~$\overline x$ dans~$B/J$ est entière sur~$A/I$ : on peut ou bien le
voir directement (si~$P$ annule~$x$ alors la réduction de~$P$ modulo~$I$ annule $\overline x$)
ou bien utiliser~\ref{sorites-ent-mor}
pour conclure que~$\bar x$ est entier sur~$A$, puis~\ref{ent-et-quot}
pour en déduire qu'il est entier sur~$A/I$.

\deux{ent-fin-equiv}
{\bf Proposition.}
{\em Soit~$B$ une~$A$-algèbre. 
Les assertions suivantes sont équivalentes : 

\medskip
1) $B$ est finie. 

2) $B$ est entière et de type fini. 

3) $B$ est engendrée comme~$A$-algèbre par un nombre
fini d'éléments entiers sur~$A$.}

\medskip
{\em Démonstration.}
On a déjà vu qu'une~$A$-algèbre finie est de type fini, et entière (\ref{fini-impl-ent}). L'implication
2)$\Rightarrow$3) est évidente ; il reste à montrer que~3)$\Rightarrow$1). 

Supposons que~$B$ soit engendrée par une famille finie~$b_1,\ldots b_r$ d'éléments
entiers sur~$A$. Nous allons montrer par récurrence sur~$r$ que~$B$ est finie. 

\medskip
Si~$r=0$ alors la flèche structurale~$A\to B$ est surjective, ce qui veut dire que~$B$ est de la forme
$A/I$ pour un certain idéal~$I$ ; elle est dès lors finie (\ref{asuri-fini}). 

\medskip
Supposons~$r>0$ et la propriété vraie au rang~$r-1$. 
Notons~$C$
la~$A$-algèbre~$A[b_1,\ldots, b_{r-1}]$. D'après l'hypothèse 
de récurrence, $ C$ est finie sur~$A$. L'élément~$b_r$
de~$B$ est entier sur~$A$
par hypothèse ; il l'est {\em a fortiori}
sur~$C$. La~$C$-algèbre~$C[b_r]=B$ est donc finie en vertu de la 
proposition~\ref{equiv-entier} ; par transitivité, $B$ est finie sur~$A$. ~$\Box$

\deux{coro-fermeture-integrale}
{\bf Corollaire-définition.}
{\em Soit~$B$ une~$A$-algèbre. Le sous-ensemble~$C$
de~$B$ formé des éléments entiers sur~$A$ est une sous-$A$-algèbre de~$B$, 
que l'on appelle {\em fermeture intégrale}
de~$A$ dans~$B$.
En particulier, si~$B$ est engendrée par des éléments entiers sur~$A$,
elle est entière sur~$A$.} 

\medskip
{\em Démonstration.}
Soient~$x$ et~$y$ appartenant à~$C$. Il résulte de
la proposition~\ref{ent-fin-equiv}
ci-dessus, et plus précisément de l'implication
3)$\Rightarrow$2) de son énoncé, que~$A[x,y]$ est entière sur~$A$,
ce qui achève la démonstration.~$\Box$ 

\deux{coroll-transitiv}
{\bf Corollaire.}
{\em Soit~$B$ une~$A$-algèbre entière, soit~$C$ une~$B$-algèbre et soit~$x$ un 
élément de~$C$. Si~$x$ est entier sur~$B$ il est entier sur~$A$ ; en particulier si~$C$
est entière sur~$B$ elle est entière sur~$A$.}

\medskip
{\em Démonstration.}
Comme~$x$ est entier sur~$B$, il existe~$n\geq 0$ et~$b_0,\ldots, b_{n-1}\in B$
tels que~$x^n+b_{n-1}x^{n-1}+\dots+b_0=0$. Soit~$B'$ la sous-$A$-algèbre
de~$B$ engendrée par les~$b_i$. Comme ceux-ci sont entiers sur~$A$
(car~$B$ est entière sur~$A$), la~$A$-algèbre~$B'$ est finie 
d'après la proposition~\ref{ent-fin-equiv}. 

Par construction, $x$ est entier sur~$B'$, et~
$B'[x]$ est donc finie sur~$B'$. Par transitivité, $B'[x]$ est finie sur~$A$,
et~$x$ est 
en conséquence entier sur~$A$.~$\Box$ 

\deux{zbarre-pasfini}
{\bf Exemple d'algèbre entière non finie.}
Notons~$\overline{\ZZ}$ la fermeture intégrale de~$\ZZ$
dans~$\CC$. C'est par définition une~$\ZZ$-algèbre entière. Nous allons
montrer par l'absurde qu'elle n'est pas finie. 

\medskip
Si elle l'était, il résulterait
de~\ref{deg-poly-ent}
qu'il existe un entier~$N$ tel que tout élément
de~$\overline {\ZZ}$
soit annulé par un polynôme unitaire
de degré~$N$ à coefficients 
dans~$\ZZ$. 

\medskip
Soit maintenant~$n\geq 1$. L'élément~$\sqrt[n]2$ de~$\CC$
appartient à~$\overline{\ZZ}$ car il est annulé par~$X^n-2$. Ce dernier est irréductible
sur~$\QQ$
en vertu du critère d'Eisenstein, et est donc le polynôme minimal
de~$\sqrt[n]2$ sur~$\QQ$. En conséquence, $\sqrt [n]2$ n'est racine d'aucun
polynôme non nul à coefficients rationnels de degré~$<n$, et on aboutit ainsi à une contradiction
en prenant~$n>N$.

\deux{entier-changebase}
{\bf Lemme}.
{\em Soit~$C$ une~$A$-algèbre et soit~$B$ une~$A$-algèbre. Si
la~$A$-algèbre~$C$ est finie (resp. entière) alors la~$B$-algèbre
$B\otimes_A C$ est finie (resp. entière).}

\medskip
{\em Démonstration.}
Supposons~$C$ finie sur~$A$, et soit~$(e_i)$ une famille
génératrice finie de~$C$ comme~$A$-module. Comme~$(1\otimes e_i)$ engendre~$B\otimes_AC$
comme~$B$-module, la~$B$-algèbre $B\otimes_A C$ est finie. 

\medskip
Supposons~$C$ entière sur~$A$. La~$B$-algèbre~$B\otimes_A C$ est engendrée
comme~$B$-algèbre (et même comme~$B$-module) par les~$1\otimes c$ pour~$c$ 
parcourant~$C$. Or si~$c\in C$, l'élément~$1\otimes c$ de~$C$ est entier sur~$A$
(\ref{sorites-ent-mor}), 
et {\em a fortiori}
sur~$B$. Il s'ensuit, en vertu du corollaire~\ref{coro-fermeture-integrale},
que~$B\otimes_A C$ est entière sur~$B$.~$\Box$ 

\deux{def-int-clos}
Soit~$A$ un anneau intègre. On appelle
{\em clôture intégrale de~$A$}
la fermeture intégrale de~$A$ dans son corps des fractions. 
On dit que~$A$ est {\em intégralement clos},
ou {\em normal},
s'il est égal à sa clôture intégrale. 

\trois{fact-int-clos}
{\bf Exercice.}
Démontrez qu'un anneau factoriel est intégralement clos. Démontrez que l'anneau des fonctions
holomorphes de~$\CC$
dans~$\CC$ est intègre, et intégralement clos. 

\trois{contrex-int-clos}
L'anneau~$\ZZ[\sqrt 5]$ n'est pas intégralement clos. Son corps
des fractions est en effet~$\QQ(\sqrt 5)$, lequel contient
le nombre d'or
$$\frac {1+\sqrt 5}2,$$ qui n'appartient pas à~$\ZZ[\sqrt 5]$ mais est entier
sur celui-ci, puisqu'il est racine de~$X^2-X-1$. 

\medskip
L'anneau~$\CC[T^2,T^3]\subset \CC[T]$ n'est pas intégralement clos. Son corps
des fractions est en effet~$\CC(T)$ (car~$T=T^3/T^2$). 
Il contient
$T$ qui n'appartient pas à~$\CC[T^2,T^3]$ mais est entier
sur celui-ci, puisqu'il est racine de~$X^2-T^2$.

\deux{intro-alg-ent-corps}
Nous allons maintenant énoncer un lemme très simple
dont la preuve est élémentaire, et sur lequel repose
{\em in fine}
le lemme crucial dit
de «going-up»
que nous verrons un peu plus loin. 

\deux{alg-ent-corps}
{\bf Lemme.}
{\em Soit~$B$ un anneau intègre et soit~$A$ un sous-anneau de~$B$. On 
suppose que~$B$ est entier sur~$A$. Les assertions
suivantes sont alors équivalentes. 

\medskip
i) $A$ est un corps. 

ii) $B$ est un corps. }

\medskip
{\em Démonstration.}
Supposons que~$A$ est un corps, et soit~$x$ un élément non nul de~$B$. 
Comme~$x$ est entier sur~$A$, la~$A$-algèbre $A[x]$ est un~$A$-espace
vectoriel de dimension finie. Comme~$B$ est intègre et~$x$ non nul, 
l'endomorphisme~$y\mapsto xy$
de~$A[x]$ est injectif, et partant surjectif. Il existe en particulier~$y\in A[x]\subset B$ tel
que~$xy=1$, et~$x$ est inversible
dans~$B$. Ainsi, $B$ est un corps. 

\medskip
Supposons maintenant que~$B$ est un corps, et soit~$x$ un élément
non nul de~$A$. Comme~$B$ est un corps, $x$ possède un inverse
$1/x$ dans~$B$. Comme~$B$ est entière sur~$A$, il existe~$n\in \NN$ 
et~$a_0,\ldots, a_{n-1}$ dans~$A$ tels que
$$\frac 1 {x^n} + a_{n-1}\frac 1 {x^{n-1}}+\ldots +a_0=0.$$ En multipliant par~$x^n$, il vient
$$1=x(-a_{n-1}-a_{n-2}x-\ldots -a_0x^{n-1}),$$
et donc
$$\frac 1 x = -a_{n-1}-a_{n-2}x-\ldots -a_0x^{n-1},$$ 
qui appartient à~$A$. Ainsi, $x$ est inversible dans~$A$ et~$A$ est un corps.~$\Box$ 

\subsection*{Degré de transcendance}

\deux{corps-alg-pasentier}
Soit~$k\hookrightarrow L$ une extension de corps. 

\trois{def-alg-corps}
Au lieu
de dire qu'un élément donné de~$L$ est entier sur~$k$, on dit plutôt qu'il est {\em algébrique}
sur~$k$. Si tout élément de~$L$ est algébrique sur~$k$, on dit que~$L$ elle-même 
est algébrique sur~$k$. 

\trois{fermeture-algebr}
L'ensemble des éléments de~$L$ algébriques sur~$K$
est une~$k$-algèbre d'après
le corollaire~\ref{coro-fermeture-integrale}, qui est intègre et entière sur~$k$ ; c'est donc
un corps en vertu du lemme~\ref{alg-ent-corps} (vous connaissiez certainement ces faits ; le but 
de ces remarques est simplement de montrer qu'on peut les retrouver
comme des cas particuliers 
de ce qu'on a établi plus haut au sujet des algèbres entières).

\deux{def-alg-indep}
Soit~$k$
un corps, soit~$A$ une~$k$-algèbre,
et soit~$(x_i)$ une famille
d'éléments
de~$A$. On dit que
les~$x_i$ sont {\em algébriquement indépendants} sur~$k$ si le 
morphisme~$k[X_i]_i\to k[x_i]_i$ qui envoie~$X_i$ sur~$x_i$ pour tout~$i$
est bijectif. Il est toujours surjectif ; par conséquent, les~$x_i$ sont algébriquement
indépendants si et seulement si il est injectif, c'est-à-dire si et seulement si
les~$x_i$ n'annulent aucun polynôme non trivial à coefficients dans~$k$. 

\medskip

Si~$A$ est non nulle il existe toujours au moins une famille d'éléments
de~$A$ algébriquement indépendants sur~$k$ : la famille
{\em vide}. Notons par contre que si~$A=\{0\}$ elle
ne possède
aucune famille d'éléments algébriquement indépendants puisque~$k$ ne s'injecte pas
dans~$\{0\}$. 

\deux{corps-al-indep}
Soit~$L$ une extension de~$k$. 

\trois{singl-alg-indep}
Si~$x\in L$, la famille singleton~$\{x\}$ est algébriquement indépendante
sur~$k$ 
si et seulement si~$x$ {\em n'est pas} algébrique
sur~$k$. On dit alors que~$x$ est {\em transcendant}
sur~$k$. 

\trois{alg-indep-max}
Soit~$(x_i)$ une famille d'éléments de~$L$ algébriquement indépendants sur~$k$. 
Le corps~$k(x_i)_{i\in I}$ engendré par~$k$ et les~$x_i$ 
s'identifie au corps de fractions rationnelles~$k(X_i)_{i\in I}$ en les indéterminées~$X_i$. 

Il résulte immédiatement des définitions que la famille~$(x_i)_{i\in I}$ 
est maximale, en tant que famille d'éléments de~$L$ algébriquement indépendants sur~$k$, si et seulement
si~$L$ est algébrique sur~$k(x_i)_{i\in I}$.

\trois{def-base-transc}
On appelle {\em base de transcendance} de~$L$
sur~$k$ une famille d'éléments de~$L$ algébriquement indépendants
sur~$k$ et maximale pour cette propriété. La théorie de l'indépendance algébrique
et des bases de transcendance est tout à fait analogue
à celle de l'indépendance linéaire et des bases. Il existe en fait
une théorie 
générale qui couvre les deux, à savoir celle des {\em relations de dépendance abstraites} 
(voir
par exemple à ce sujet {\em Basic Algebra II}, de Jacobson).  
On démontre ainsi les faits suivants.

\medskip

$\bullet$  Toute 
famille d'éléments de~$L$ algébriquement indépendants de~$k$ est contenue dans une base de transcendance
de~$L$ sur~$k$ ; 
en particulier, il existe une base de transcendance 
de~$L$ sur~$k$ (appliquer ce qui précède
à la famille vide). 

$\bullet$ Si~$(x_i)$ est une famille d'éléments de~$L$ telle que~$L$
soit algébrique sur~$k(x_i)_i$, elle contient une base de transcendance
de~$L$ sur~$k$. 

$\bullet$ Toutes les bases de transcendance
de~$L$ sur~$k$ ont même cardinal, appelé {\em degré de transcendance de~$L$ sur~$k$} ; 

$\bullet$ Si le degré de transcendance de~$L$ sur~$K$ est fini, toute famille
d'éléments de~$L$ algébriquement indépendants sur~$K$ qui est de cardinal~$\deg {\rm tr.} (L/K)$
est une base de transcendance 
de~$L$ sur~$k$ ; et toute famille~$(x_i)_{i\in I}$ de cardinal $\deg {\rm tr.} (L/K)$
telle que~$L$ soit algébrique sur~$k(x_i)_i$ est une base de transcendance
de~$L$ sur~$k$.

\trois{ex-base-tr}
{\bf Exemples.} Le
degré de transcendance de~$k(X_1,\ldots, X_n)$ sur~$k$ est égal à~$n$,
et~$(X_1,\ldots, X_n)$ en est une base de transcendance. 

\medskip
Soit~$f$ une fraction rationnelle non constante dans~$k(X)$. L'élément~$X$ est alors
algébrique sur~$k(f)$ (exercice facile). En conséquence~$f$ est transcendant
(sinon,~$X$ serait algébrique sur~$k$), et~$\{f\}$ est une base de transcendance
de~$k(X)$ sur~$k$.

\subsection*{Lemme de
{\em going-up}
et dimension de Krull}

\deux{intro-going-up}
Nous allons maintenant démontrer le
lemme dit
de {\em going-up}. Il est extrêmement utile en
en algèbre commutative et géométrie algébrique (le
langage des schémas permettra d'en donner une interprétation géométrique), mais 
est également intéressant pour sa preuve. Celle-ci consiste en effet essentiellement
à se ramener au lemme~\ref{alg-ent-corps} ci-dessus , lui-même élémentaire, 
par de judicieux passages au quotient et localisations ; elle peut donc 
aider à mieux comprendre comment utiliser en pratique ces opérations.

\deux{going-up}
{\bf Lemme de {\em going-up}.} {\em Soit~$f: A\to B$ un morphisme d'anneaux faisant de~$B$ une~$A$-algèbre entière. 

\medskip
1) Si~$\got q$ est un idéal premier
de~$B$, l'idéal premier~$f^{-1}(\got q)$ de~$A$ est maximal
si et seulement si~$\got q$ est maximal. 

2) Si~$\got q$ et~$\got q'$ sont
deux idéaux premiers~{\em distincts} de~$B$ tels que
l'on ait~$f^{-1}(\got q)=f^{-1}(\got q')$,
alors~$\got q$ et~$\got q'$ sont non comparables pour l'inclusion. 

3) Si~$J$ est un idéal de~$B$ et si l'on pose~$I=f^{-1}(J)$, 
il existe
pour tout idéal premier~$\got p$ de~$A$ contenant~$I$
un idéal premier~$\got q$
de~$B$ contenant~$J$ tel que~$f^{-1}(\got q)=\got p$. }

\medskip
{\em Remarque.} Dans le cas où~$f$ est injectif, on a~$f^{-1}(\{0\})=\{0\}$ et l'assertion~3) ci-dessus affirme
alors que pour tout idéal premier~$\got p$ de~$A$, il existe un idéal premier~$\got q$
de~$B$ tel que~$f^{-1}(\got q)=\got p$ (en fait, 
comme on le verra ci-dessous, on {\em ramène}
la preuve de 3) à celle
de ce cas particulier). Autrement dit, si~$f: A\to B$ est une injection entière
alors~${\rm Spec}\; B\to {\rm Spec}\; A$ est surjective. 

\deux{preuve-going-up}
{\em Démonstration du lemme
de going-up.}
On commence par établir l'assertion~1), qui sera elle-même utilisée dans la preuve de~2) et~3). 

\trois{going-up1}
{\em Preuve de~1)}. Soit~$\got q$ un idéal premier de~$B$. Posons~$\got p=f^{-1}(\got q)$. 
La flèche~$f : A\to B$ induit une injection~$A/\got p\hookrightarrow B/\got q$, qui fait
de~$B/\got q$ une~$A/\got p$-algèbre entière (\ref{ent-et-doublequot}). Comme~$\got q$ est premier, 
$B/\got q$ est intègre. Il résulte alors du lemme~\ref{alg-ent-corps}
que~$B/\got q$ est un corps
si et seulement si~$A/\got p$ est un corps ; autrement dit,
$\got q$ est maximal si et seulement si~$\got p$ est maximal, ce qui achève
de prouver~1). 

\trois{going-up2}
{\em Preuve de~2)}. Soient~$\got p$ un idéal premier de~$A$. Il s'agit de
montrer que les éléments de~$\{\got q\in \spec B,f^{-1}(\got q)=\got p\}$ sont deux à deux non comparables
pour l'inclusion, c'est-à-dire encore que les antécédents de~$\got p$ pour l'application
$\spec B\to \spec A$ induite par~$f$ sont deux à deux non comparables pour l'inclusion. 

\medskip
Posons~$S=A\setminus \got p$. Le
localisé~$S^{-1}A$ est l'anneau local~$A_{\got p}$
d'idéal maximal~$\got pA_{\got p}$, et l'on a par ailleurs
$$f(S)^{-1}B\simeq B\otimes_A A_{\got p}; $$
en particulier, $f(S)^{-1}B$ est entier sur~$A_{\got p}$. On a un diagramme
commutatif 
$$\diagram
\spec B \dto &\spec f(S)^{-1}B\dto\lto\\
\spec A &\spec A_{\got p}\lto\enddiagram.$$

La flèche~$\spec A_{\got p} \to \spec A$ induit une
bijection d'ensembles ordonnés entre~$\spec A_{\got p}$
et  l'ensemble des idéaux premiers de~$A$ contenus
dans~$\got p$ (elle fait correspondre~$\got p$
à~$\got pA_{\got p}$) ; la flèche
$\spec f(S)^{-1}B\to \spec B$ induit une bijection d'ensembles
ordonnés entre~$\spec f(S)^{-1}B$
et l'ensemble des idéaux premiers de~$B$ ne rencontrant pas~$f(S)$. 

Par ailleurs, soit~$\got q$ un idéal premier de~$B$ situé au-dessus de~$\got p$, 
c'est-à-dire tel que~$f^{-1}(\got q)=\got p$. Cette dernière égalité assure que
$\got q$ ne rencontre pas~$f(S)$, et donc que~$\got q$ appartient à l'image de~$\spec f(S)^{-1}B$. 
Il s'ensuit qu'il existe une bijection d'ensembles ordonnés entre l'ensemble des idéaux premiers
de~$B$ situés au-dessus de~$\got p$, et celui des idéaux premiers de~$f(S)^{-1}B$ situés au-dessus
de~$\got p A_{\got p}$. 

\medskip
Mais comme~$\got pA_{\got p}$ est l'idéal maximal de~$A_{\got p}$, il résulte de l'assertion~1),
appliquée à la~$A_{\got p}$-algèbre entière~$f(S)^{-1}B$, que l'ensemble des idéaux premiers de 
~$f(S)^{-1}B$ situés au-dessus
de~$\got p A_{\got p}$ est exactement l'ensemble des idéaux
{\em maximaux}
de~$f(S)^{-1}B$ ; or ceux-ci sont deux à deux non comparables pour l'inclusion, 
et il en va donc de même des idéaux premiers
de~$B$ situés au-dessus de~$\got p$, d'où~2). 

\trois{going-up3}
{\em Preuve de 3)}.
La flèche~$A\to B$ induit une injection~$A/I\to B/J$. 
On a par ailleurs un diagramme commutatif
$$\diagram \spec B/J\rto \dto&\spec B\dto\\
\spec A/I\rto &\spec A \enddiagram.$$
La flèche~$\spec A/I\to \spec A$ induit une bijection 
entre~$\spec A/I$ et l'ensemble des idéaux premiers de~$A$
contenant~$I$, et la flèche~$\spec B/J\to \spec B$ induit une bijection 
entre~$\spec B/J$ et l'ensemble des idéaux premiers de~$A$
contenant~$J$ ; on peut donc, quitte à remplacer~$A$ par~$A/I$ et~$B$
par~$B/J$, se ramener au cas où~$I$ et~$J$
sont nuls et où~$f$ est injective 

\medskip
Il s'agit maintenant, un idéal premier~$\got p\in \spec A$ étant donné,
de montrer l'existence d'un idéal premier~$\got q\in \spec B$ au-dessus de~$\got p$. 

Posons~$S=A\setminus \got p$.
En se fondant sur le diagramme commutatif
considéré au~\ref{going-up2}
ci-dessus lors de la preuve de~2), on voit qu'il suffit de montrer
l'existence d'un idéal premier de~$f(S)^{-1}B$ situé au-dessus de~$\got pA_{\got p}$. 

\medskip
Par définition de~$S$, cette partie ne contient pas~$0$. Comme~$f $ est injective, 
$f(S)$ ne contient pas non plus~$0$ ; par conséquent, $f(S)^{-1}B$ est non nul. Il possède
dès lors un idéal maximal~$\got m$. En vertu de l'assertion~1),
appliquée à la~$A_{\got p}$-algèbre entière~$f(S)^{-1}B$, l'idéal~$\got m$ est situé
au-dessus d'un idéal maximal de~$A_{\got p}$, et donc de~$\got pA_{\got p}$. Ceci
achève la démonstration. ~$\Box$

\deux{def-krull}
{\bf Définition.}
Soit~$A$ un anneau. On appelle
{\em dimension de Krull}
de~$A$ la borne supérieure de l'ensemble~$\sch E$
des entiers~$n$ tels qu'il existe 
une chaîne

$$\got p_0\subsetneq \got p_1\subsetneq\ldots \subsetneq \got p_n$$
où les~$\got p_i$ sont des idéaux premiers de~$A$ (attention : notez bien que la numérotation
commence à~$0$). 

\trois{valeurs-krull}
{\bf Remarque.}
L'ensemble~$\sch E$ 
peut être vide. C'est le cas si et seulement
si~$A$ n'a pas d'idéaux premiers, c'est-à-dire si et seulement si~$A=\{0\}$ ; il y a alors
une ambiguïté dans la définition de la dimension de Krull\footnote{En effet, 
nous invitons le lecteur à vérifier que la borne supérieure
de la partie vide est
par définition {\em le plus petit
élément de l'ensemble ordonné
dans lequel on travaille} (s'il existe). Il faut donc préciser ici quel est l'ensemble
en question ; les conventions adoptées reviennent à décider qu'on travaille
dans~$\RR\cup\{-\infty,+\infty\}$.},
qu'on lève en posant par convention
${\rm dim}_{\rm Krull}\{0\}=-\infty.$ 

\medskip
Lorsque~$A$ est non nul, l'ensemble
non vide~$\sch E$
peut être fini,
auquel cas la dimension de Krull appartient à~$\NN$, ou infini -- cela signifie qu'il existe
des chaînes strictement croissantes arbitrairement longues d'idéaux premiers de~$A$,
et l'on a alors~${\rm dim}_{\rm Krull}\;A=+\infty$. 

\trois{comment-krull}
{\bf Commentaire sur la terminologie.}
Nous verrons lors du cours sur les schémas que le terme {\em dimension}
est bien choisi : la dimension de Krull d'un anneau peut
en effet s'interpréter comme la dimension du schéma qui lui est associé. 

\trois{ex-krull-corps}
{\bf Anneaux de dimension nulle.}
Un anneau~$A$ est de dimension de Krull nulle si et seulement si il possède un et un seul
idéal premier. C'est notamment le cas lorsque~$A$ est un corps. 

\trois{ex-krull-dim1}
{\bf Anneaux intègres de dimension 1}. Un anneau intègre~$A$ est de
dimension~$1$ si et seulement si les deux conditions suivantes sont satisfaites : 

$\bullet$ $A$ possède un idéal premier non nul ;

$\bullet$ tout idéal 
premier non nul de~$A$ est maximal.

\medskip
C'est notamment le cas lorsque~$A$ est un anneau de Dedekind qui n'est pas un corps. 
En particulier, tout anneau principal qui n'est pas un corps est de dimension de Krull égale
à~$1$.

\trois{ex-krul-infini}
{\bf Un anneau de dimension de Krull infinie.}
Soit~$A$
l'anneau~$k[X_i]_{1\leq i}$. Pour tout~$j\in \NN$, l'idéal~$\got p_j:=(X_i)_{1\leq i\leq j}$ 
de~$A$ est premier : en effet, le quotient~$A/\got p_j$ s'identifie naturellement
à l'anneau intègre~$k[X_i]_{i>j}$. Pour tout
entier~$n$, on dispose d'une chaîne
strictement croissante de~$n+1$
idéaux premiers de~$A$,
à savoir 
$$(0)\subsetneq \got p_1\subsetneq \got p_2\subsetneq\ldots \subsetneq \got p_n.$$ 
En conséquence, ${\rm dim}_{\rm Krull}A=+\infty$. 

\deux{dim-krull-culture}
On démontre que si~$A$ est un anneau
{\em local noethérien}
dont l'idéal
maximal
possède un système de générateurs de cardinal~$n$, 
alors~${\rm dim}_{\rm Krull} A\leq n$
(en particulier, ${\rm dim}_{\rm Krull}A$ est finie). 

\trois{regulier-culture}
Un anneau local~$A$
est dit
{\em régulier}
s'il est 
noethérien
et si son idéal maximal possède un système de générateurs
de cardinal exactement égal à~${\rm dim}_{\rm Krull}A$. 
Les anneaux locaux réguliers jouent un rôle majeur en géométrie 
algébrique. 

\trois{krull-nagata}
Nagata a construit un exemple d'anneau noethérien (non local)
dont la dimension de Krull est infinie.

\deux{krull-transfert}
{\bf Proposition.}
{\em Soit~$A$ un anneau et soit~$B$ une~$A$-algèbre. On suppose que la flèche~$f:A\to B$
est injective et que~$B$ est entière sur~$A$. On a alors~${\rm dim}_{\rm Krull} B={\rm dim}_{\rm Krull} A$.}

\medskip
{\em Démonstration.} On va montrer l'égalité des dimensions
par double majoration. 

\trois{b-majore-a}
{\em Prouvons que~${\rm dim}_{\rm Krull}B\geq {\rm dim}_{\rm Krull} A$.}
Soit~$\got p_0\subsetneq \got p_1\subsetneq\ldots \subsetneq \got p_n$
une chaîne strictement croissante d'idéaux premiers de~$A$. 
Il existe alors une chaîne
$\got q_0\subset \got q_1\subset \ldots \subset \got q_n$ d'idéaux premiers de~$B$
tels que~$f^{-1}(\got q_i)=\got p_i$ pour tout~$i$. 
Pour le voir, on raisonne par récurrence sur~$n$. 

\medskip
Si~$n=0$, c'est vrai par la remarque
suivant l'énoncé du lemme de going-up
(c'est ici qu'intervient l'injectivité de~$f$). 

Supposons maintenant~$n>0$ et les~$\got q_i$ construits
pour~$i\leq n-1$. On applique alors l'assertion~3)
du lemme de going-up avec~$J=\got q_{n-1}, I=f^{-1}(\got q_{n-1})
=\got p_{n-1}$, et~$\got p=\got p_n$. Elle assure l'existence d'un idéal premier~$\got q_n$
de~$B$ contenant~$\got q_{n-1}$ et tel que~$f^{-1}(\got q_n)=\got p_n$, qui est ce qu'on souhaitait. 

\medskip
Comme~$f^{-1}(\got q_i)=\got p_i$ pour tout~$i$, et comme les~$\got p_i$ 
sont deux à deux distincts, les~$\got q_i$ sont deux à deux distincts. La chaîne
des~$\got q_i$ est ainsi {\em strictement}
croissante, et a même longueur que celle des~$\got p_i$. On en déduit
l'inégalité requise $${\rm dim}_{\rm Krull}B\geq {\rm dim}_{\rm Krull} A.$$

\trois{b-minore-a}
{\em Prouvons que~${\rm dim}_{\rm Krull}A\geq {\rm dim}_{\rm Krull} B$.}
Soit~$\got q_0\subsetneq \got p_1\subsetneq\ldots \subsetneq \got q_n$
une chaîne strictement croissante d'idéaux premiers de~$B$. Pour tout
indice~$i$,
posons~$\got p_i=f^{-1}(\got q_i)$. Les~$\got p_i$ sont des idéaux premiers de~$A$. 

\medskip
Soit~$i\leq n-1$ ; nous allons montrer 
que~$\got p_i\subsetneq \got p_{i+1}$. L'inclusion
large  $\got p_i\subset \got p_{i+1}$ provient immédiatement des définitions (et de l'inclusion large
de~$\got q_i$ dans~$\got q_{i+1}$). Il reste à s'assurer que~$\got p_i\neq \got p_{i+1}$ ; 
mais s'ils étaient égaux, on aurait
alors~$f^{-1}(\got q_i)=f^{-1}(\got q_{i+1})$ avec~$\got q_i\subsetneq\got q_{i+1}$,
contredisant l'assertion~2) du lemme de going-up. 

\medskip
La chaîne
des~$\got p_i$ est ainsi {\em strictement}
croissante, et a même longueur que celle des~$\got q_i$. On en déduit
l'inégalité requise $${\rm dim}_{\rm Krull}A\geq {\rm dim}_{\rm Krull} B,$$
ce qui achève la démonstration.~$\Box$ 

\section{Algèbres de type fini sur un corps : normalisation 
de Noether, {\em Nullstellensatz}}
 \markboth{Algèbre commutative}{Normalisation de Noether et {\em Nullstellensatz}}

{\em On fixe pour toute cette section un corps~$k$.}

\subsection*{Le lemme de normalisation de Noether}

\deux{objectif-noether}
Notre premier objectif est de démontrer le
«lemme de normalisation de Noether», qui permet
dans certaines situations de ramener l'étude d'une~$k$-algèbre
de type fini~$A$ 
à celle d'une algèbre de polynômes (le fait de savoir écrire~$A$
comme un {\em quotient}
d'une telle algèbre est
le plus souvent insuffisant, tant la théorie des idéaux 
de~$k[X_1,\ldots, X_n]$ est complexe dès que~$n>1$).

\medskip
Pour ce faire, nous avons choisi d'isoler le cœur technique
de la démonstration en lui donnant un statut de proposition 
indépendante, par laquelle nous allons commencer. 

\deux{prop-pre-noether}
{\bf Proposition.}
{\em Soit~$n\geq 0$ et soit~$I$ un idéal {\em strict}
de~$k[X_1,\ldots, X_n]$. Il existe~$n$ éléments~$Y_1,\ldots, Y_n$ 
de~$k[X_1,\ldots, X_n]$ tels que les propriétés suivantes
soient satisfaites : 

\medskip
1) les~$Y_i$ sont algébriquement indépendants sur~$k$ ; 

2) $k[X_1,\ldots, X_n]$ est entière sur~$k[Y_1,\ldots, Y_n]$ ; 

3) l'idéal $I\cap k[Y_1,\ldots, Y_n]$ de~$k[Y_1,\ldots, Y_n]$ 
est engendré par~$(Y_1,Y_2,\ldots, Y_m)$ pour un certain
entier~$m$ compris entre~$0$ et~$n$. }

\medskip
{\em Démonstration.}
On procède en plusieurs étapes. 

\trois{comment-pre-noether}
{\bf Remarque préalable.} Supposons avoir construit
une famille~$(Y_i)_{1\leq i\leq n}$ d'éléments de~$k[X_1,\ldots, X_n]$ 
telle que~$k[X_1, \ldots, X_n]$ soit entière sur~$k[Y_1,\ldots, Y_n]$. Le corps
$k(X_1,\ldots, X_n)$ est alors engendré par des éléments algébriques sur~$k(Y_1,\ldots, Y_n)$, et
est donc lui-même algébrique sur~$k(Y_1,\ldots, Y_n)$. En
conséquence, $(Y_1,\ldots, Y_n)$ contient une base de transcendance de~$k(X_1,\ldots, X_n)$ sur~$k$ ; 
mais comme le degré de transcendance de $k(X_1,\ldots, X_n)$ sur~$k$ 
vaut~$n$, cela signifie que~$(Y_1,\ldots, Y_n)$ est une base de transcendance de~$k(X_1,\ldots, X_n)$ sur~$k$. 
En particulier, les~$Y_i$ sont algébriquement indépendants. Ainsi, si l'on prouve~2), l'assertion~1)
sera {\em automatiquement}
vérifiée. 

\trois{debut-preuve-prenoeth}
Pour démontrer la proposition, on procède par récurrence sur~$n$. Si~$n=0$ alors~$k[X_1,\ldots, X_n]=k$. 
Comme~$I$ est strict, il est nécessairement nul. On vérifie alors aussitôt qu'on peut
prendre pour~$(Y_i)$ la {\em famille vide}. 

\medskip
Supposons maintenant~$n>0$, et la proposition vraie au rang~$n-1$. Si~$I=\{0\}$
on peut prendre~$Y_i=X_i$ pour tout~$i$. On se place désormais
dans le cas où~$I\neq 0$. On choisit
alors un élément~$P\neq 0$ dans~$I$ ; comme~$I$ est strict, $P$ est non constant. 

\trois{construction-zi}
On pose~$Z_1=P$. La
prochaine étape de la preuve
consiste à 
construire~$Z_2,\ldots,Z_n$
de sorte que~$k[X_1,\ldots, X_n]$ soit entier sur~$k[Z_1,\ldots, Z_n]$. 

\medskip
Nous allons en fait chercher les~$Z_i$ sous une forme très
particulière. Plus précisément, donnons-nous
une famille~$r_2,\ldots, r_n$ d'entiers~$>0$. Nous allons chercher
des conditions suffisantes sur les~$r_i$ pour que
$k[X_1,\ldots, X_n]$ soit entier sur~$k[Z_1,\ldots, Z_n]$
avec~$Z_i=X_i+X_1^{r_i}$ pour~$i\geq 2$. 

\medskip
Notons déjà que comme~$X_i=Z_i-X_1^{r_i}$ pour tout~$i\geq 2$, et 
comme~$Z_i$ est évidemment entier sur~$k[Z_1,\ldots, Z_n]$ pour tout~$i$
il suffit, pour obtenir ce qu'on souhaite, de faire en sorte que~$X_1$ soit entier sur~$k[Z_1,\ldots, Z_n]$. 

\medskip
Lorsqu'on remplace dans l'écriture de~$P$ chacun des~$X_i$ pour~$i\geq 2$
par~$Z_i-X_1^{r_i}$, un monôme donné de~$P$ de multi-degré~$(i_1,\ldots, i_n)$ 
donne lieu à un polynôme en~$X_1$ à coefficients dans~$k[Z_1,\ldots, Z_n]$ dont le terme dominant
est de la forme~$\alpha X_1^{i_1+r_2i_2+\ldots+r_n i_n}$, avec~$\alpha\in k\ti$. 

\trois{supposition-iiri}
Supposons que lorsque~$(i_1,\ldots, i_n)$ parcourt l'ensemble des multi-degrés
des monômes de~$P$, les entiers~$i_1+r_2i_2+\ldots+r_ni_n$ soient deux à deux distincts. 
Dans ce cas, le maximum~$N$ de ces entiers est atteint pour un et un seul multi-degré
$(i_1,\ldots, i_n)$, et~$N$ est nécessairement non nul puisque~$P$ n'est pas constant. Il résulte 
alors de ce qui précède que l'égalité~$P-Z_1=0$
peut se récrire
$$\alpha X_1^N+Q(X_1)-Z_1=0,$$ où~$\alpha \in k\ti$
et où~$Q$ est un polynôme de degré~$<N$ en~$X_1$, à coefficients dans~$k[Z_1,\ldots, Z_n]$. 
Comme~$N>0$, le polynôme~$Q(X_1)-Z_1$ est encore un polynôme de degré~$<N$ en~$X_1$, à coefficients dans~$k[Z_1,\ldots, Z_n]$. 
En multipliant l'égalité ci-dessus par~$\alpha^{-1}$, on obtient
$$X_1^N+\alpha^{-1}(Q(X_1)-Z_1)=0,$$ et~$X_1$ est entier sur~$k[Z_1,\ldots, Z_n]$. 

\trois{conclusion-zi}
Il suffit donc, pour construire des polynômes $Z_2,\ldots, Z_n$
satisfaisant la propriété souhaitée, 
d'exhiber une famille~$(r_2,\ldots, r_n)$ d'entiers
strictement positifs telle que
les entiers~$i_1+r_2i_2+\ldots+r_ni_n$ soient deux à deux distincts
lorsque~$(i_1,\ldots, i_n)$ parcourt l'ensemble des multi-degrés
des monômes de~$P$. 

\medskip
Soit~$p$ un entier strictement supérieur au maximum des degrés intervenant
dans les monômes de~$P$. On peut alors prendre $r_2=p, r_3=p^2,\ldots, r_n=p^{n-1}$ : 
en effet, ils satisfont la condition exigée en vertu de l'unicité du développement
d'un entier en base~$p$.

\trois{passage-zi-yi}
{\em Construction des~$Y_i$ et conclusion.}
On déduit
de~\ref{comment-pre-noether}
que les~$Z_i$ sont algébriquement indépendants sur~$k$. En particulier,
$k[Z_2,\ldots, Z_n]$ est une algèbre de polynômes en~$n-1$
variables. Comme~$I$ est strict, $1\notin I\cap k[Z_2,\ldots, Z_n]$
et ce dernier est donc un idéal strict de $k[Z_2,\ldots, Z_n]$. 
Il s'ensuit, par l'hypothèse
de récurrence, qu'il existe
une famille~$(Y_2,\ldots, Y_n)$
d'éléments de~$k[Z_2,\ldots, Z_n]$, algébriquement indépendants
sur~$k$ et tels que :

\medskip
$\bullet$ $k[Z_2,\ldots, Z_n]$ est entier
sur~$k[Y_2,\ldots, Y_n]$ ;

$\bullet$ l'idéal $I\cap k[Y_2,\ldots, Y_n]$ 
de~$k[Y_2,\ldots, Y_n]$ est
engendré par~$(Y_2,\ldots, Y_m)$
pour un certain~$m$. 

\medskip
Posons~$Y_1=Z_1$. L'anneau~$k[Z_1,\ldots, Z_n]$ est alors
engendré par des éléments entiers sur~$k[Y_1,\ldots, Y_n]$. Il est donc entier
sur~$k[Y_1,\ldots, Y_n]$ ; par transitivité, $k[X_1,\ldots, X_n]$ est entier sur~$k[Y_1,\ldots, Y_n]$. 
Cela force les~$Y_i$ à être algébriquement indépendants sur~$k$ (\ref{comment-pre-noether}). 

\medskip
Pour conclure, il suffit de montrer que
l'idéal~$I\cap k[Y_1,\ldots, Y_n]$ de~$k[Y_1,\ldots, Y_n]$ 
est engendré par~$(Y_1,\ldots, Y_m)$. On a~$Y_1=Z_1=P\in I$, et~$Y_j\in I$
pour~$2\leq j\leq m$ par choix des~$Y_j$. 
Par conséquent, l'idéal de~$k[Y_1,\ldots, Y_n]$ engendré
par les~$Y_j$ pour~$1\leq j\leq m$
est contenu~$I\cap k[Y_1,\ldots, Y_n]$. 

\medskip
Il reste à vérifier l'inclusion réciproque. Soit~$R\in I\cap k[Y_1,\ldots, Y_n]$. 
On peut écrire~$R=Y_1S+T$, où~$T\in k[Y_2,\ldots, Y_n]$ (on met à part les monômes
contenant~$Y_1$). Comme~$R\in I$ et~$Y_1\in I$, on a~$T\in I\cap k[Y_2,\ldots, Y_n]$. 
Par conséquent, $T$ s'écrit~$\sum_{2\leq i\leq m} A_i Y_i$ avec~$A_i\in k[Y_2,\ldots, Y_m]$ pour tout~$i$,
et
on a
alors~$R=SY_1+\sum_{2\leq i\leq m} A_i Y_i$. Ainsi, $R$ appartient à l'idéal
de~$k[Y_1,\ldots, Y_n]$ engendré par~$(Y_1,\ldots, Y_m)$, ce qui achève la démonstration.~$\Box$ 

\deux{normal-noeth}
{\bf Lemme de normalisation de Noether.}
{\em Soit~$A$ une~$k$-algèbre de type fini non nulle. Il
existe une famille~$(y_1,\ldots, y_d)$ d'éléments de~$A$
algébriquement indépendants sur~$k$
et tels que~$A$ soit entière sur~$k[y_1,\ldots, y_d]$. }

\medskip
{\em Démonstration.}
Comme~$A$ est de type fini, elle s'identifie à~$k[X_1,\ldots, X_n]/I$ pour un 
certain~$n$ et un certain idéal~$I$ de~$k[X_1,\ldots, X_n]$, qui est strict
puisque~$A$ est non nulle. 

\medskip
On peut donc appliquer la proposition~\ref{prop-pre-noether}
ci-dessus. Elle assure l'existence d'une famille~$(Y_1,\ldots, Y_n)$
d'éléments de~$k[X_1,\ldots, X_n]$ algébriquement indépendants sur~$k$
et d'un entier~$m\leq n$
tels que~$k[X_1,\ldots, X_n]$ soit entière sur~$k[Y_1,\ldots, Y_n]$,
et tels que l'idéal~$I\cap k[Y_1,\ldots, Y_n]$ de~$k[Y_1,\ldots, Y_n]$ 
soit engendré par~$Y_1,Y_2,\ldots, Y_m$. Notons~$y_1,\ldots, y_d$ les
classes respectives de~$Y_{m+1},\ldots, Y_n$ modulo~$I$ (on a donc~$d=n-m$). 
Comme~$k[X_1,\ldots, X_n]$ est entière sur~$k[Y_1,\ldots, Y_n]$,
l'anneau quotient~$A=k[X_1,\ldots, X_n]/I$ est entier sur
sa sous-algèbre~$k[Y_1,\ldots, Y_n]/(I\cap k[Y_1,\ldots, Y_n])$, 
qui est engendrée par les~$y_i$ 
car chacun des~$Y_j$ pour~$j\leq m$ appartient à~$I$ ; par conséquent, $A$
est entière
sur~$k[y_1,\ldots, y_d]$. 

\medskip
Il reste à s'assurer que les~$y_i$ sont algébriquement indépendants sur~$k$. 
Soit donc~$Q\in k[T_1,\ldots, T_d]$ un polynôme s'annulant en les~$y_i$. 
Cela signifie que~$Q(Y_{m+1}\ldots, Y_n)\in I\cap k[Y_1,\ldots, Y_n]$. 
Cet idéal de~$k[Y_1,\ldots, Y_n]$
étant engendré par les~$Y_j$ pour~$j\leq m$, le polynôme~$Q(Y_{m+1},\ldots, Y_n)$
peut s'écrire~$\sum_{1\leq j\leq m}A_j Y_j$ où les~$A_j$ appartiennent à~$k[Y_1,\ldots, Y_n]$. 
Mais comme~$Q(Y_{m+1},\ldots, Y_n)$ ne
fait intervenir que les variables~$Y_j$ pour~$j>m$, il ne peut admettre
une telle écriture que s'il est nul. Il en résulte que~$Q=0$, ce qu'il fallait démontrer.~$\Box$

\deux{comment-noether}
{\bf Commentaires.}

\trois{entier-fini-noether}
Comme~$A$ est de type fini comme~$k$-algèbre, elle est
de type fini comme~$B$-algèbre pour toute sous-algèbre~$B$
de~$A$. Il résulte alors de la proposition~\ref{ent-fin-equiv}
que l'on obtiendrait un énoncé équivalent à celui du lemme de normalisation
de Noether en remplaçant «entière sur~$k[y_1,\ldots, y_d]$»
par «finie sur~$k[y_1,\ldots, y_d]$». Nous avons opté pour l'épithète
«entière»
car ce que notre preuve du lemme de normalisation de Noether montre {\em effectivement},
c'est bien le caractère entier de~$A$ sur~$k[y_1,\ldots, y_d]$. 

\trois{integre-dimtr}
Supposons~$A$ intègre. Dans ce cas, son corps des fractions est engendré
par des éléments algébriques sur~$k(y_1,\ldots, y_d)$ et est donc lui-même
algébrique sur~$k(y_1,\ldots, y_d)$. En conséquence, $(y_1,\ldots, y_d)$ est une base
de transcendance de~${\rm Frac}\;A$ sur~$k$, et l'entier~$d$ est 
dès lors nécessairement
égal au degré de transcendance de~${\rm Frac}\;A$ sur~$k$. 

\subsection*{Le {\em Nullstellensatz}}
\deux{intro-nullst}
Nous allons maintenant donner trois variantes 
du {\em Nullstellensatz}\footnote{En allemand : théorème du lieu des zéros.},
également appelé «théorème des zéros de Hilbert». Les deux premières s'énoncent dans un langage
purement algébrique, la troisième évoque la géométrie algébrique et justifiera l'appellation du théorème. 

\deux{nullst-1}
{\bf Première variante du~{\em Nullstellensatz}.}
{\em Soit~$A$ une~$k$-algèbre de type fini. Si~$A$ est un corps, c'est une
extension finie de~$k$.}

\medskip
{\em Démonstration.}
Le lemme de normalisation de Noether assure qu'il existe
des éléments~$y_1,\ldots, y_d$ de~$A$, algébriquement indépendants
sur~$k$ et tels que~$A$ soit entière (ou finie, {\em cf.}~\ref{entier-fini-noether})
sur~$k[y_1,\ldots, y_d]$. 

\medskip
Si~$A$ est un corps,
l'anneau de polynômes $k[y_1,\ldots,y_d]$ est lui aussi un corps
d'après le lemme~\ref{alg-ent-corps}, ce qui force~$d$ à être nul, et
$k[y_1,\ldots, y_d]$ à être égale à~$k$ ; le corps~$A$
est alors une~$k$-algèbre finie, ce qu'il fallait démontrer.~$\Box$ 

\deux{nullst-2}
{\bf Deuxième
variante du~{\em Nullstellensatz}.}
{\em Soit~$A$ un~$k$-algèbre de type fini et soit~$\got p$
un idéal de~$A$. les assertions suivantes sont équivalentes : 

\medskip
i) $A/\got p$ est une extension finie de~$k$ ; 

ii) ${\rm Frac}(A/\got p)$ est une extension finie de~$k$ ; 

iii) $\got p$ est maximal.} 

\medskip
{\em Démonstration.} Il est clair que~i)$\Rightarrow$ii). Supposons que~ii) soit vraie. 
L'algèbre~$A/\got p$ est alors, en tant que sous-algèbre de ${\rm Frac}(A/\got p)$,
une algèbre intègre et entière sur~$k$. D'après le lemme~\ref{alg-ent-corps}, c'est un corps
et~$\got p$ est maximal. 

\medskip
Enfin, si~$\got p$ est maximal alors~$A/\got p$ est une~$k$-algèbre de type fini
qui est un corps, et c'est donc une extension finie de~$k$ en vertu de la variante
précédente du {\em Nullstellensatz}
(\ref{nullst-1}).~$\Box$

\deux{nullst-3}
{\bf Troisième variante du {\em Nullstellensatz}}. 
{\em Soit~$n\in \NN$ et soit~$(P_i)_{i\in I}$ une famille
de polynômes appartenant à~$k[T_1,\ldots, T_n]$. Les assertions suivantes sont
équivalentes : 

\medskip
1) Il existe une extension {\em finie}~$L$ de~$k$ et un~$n$-uplet
$(x_1,\ldots, x_n)\in L^n$
tel que~$P_i(x_1,\ldots, x_n)=0$ pour tout~$i$. 

2) Il existe une extension~$L$ de~$k$ et un~$n$-uplet
$(x_1,\ldots, x_n)\in L^n$
tel que~$P_i(x_1,\ldots, x_n)=0$ pour tout~$i$. 

3) Pour toute famille~$(Q_i)_{i\in I}$
de polynômes presque tous nuls appartenant à~$k[X_1,\ldots, X_n]$, 
on a~$\sum Q_i P_i\neq 1$.}

\medskip
{\em Démonstration.}
L'implication~1)$\Rightarrow$2) est évidente. Supposons
que~2) soit vraie, et donnons-nous une famille~$(Q_i)$ de polynômes
presque tous nuls appartenant à~$k[X_1,\ldots, X_n]$. 
On a
$$(\sum Q_i P_i)(x_1,\ldots, x_n)=\sum Q_i(x_1,\ldots, x_n)P_i(x_1,\ldots, x_n)=0,$$
ce qui entraîne que~$\sum Q_i P_i\neq 1$. 

\medskip
Supposons enfin que~3)
soit vraie. Cela signifie que l'idéal~$J$ de~$k[X_1,\ldots, X_n]$
engendré par les~$P_i$ ne contient pas~1 ; autrement dit,
il est strict. 
L'algèbre quotient~$A:=k[X_1,\ldots, X_n]/J$ est alors non nulle, 
et elle admet de ce fait un idéal maximal~$\got m$. 
Comme~$A$ est de type fini sur~$k$ par construction, 
le corps quotient~$L:=A/\got m$ est
une extension finie de~$k$ en vertu de la deuxième
variante du {\em Nullstellensatz}
(\ref{nullst-2}). Si l'on note~$x_i$ l'image de~$\overline {X_i}$
par la surjection canonique~$A\to L$, le~$n$-uplet~$(x_1,\ldots, x_n)$ 
annule par construction chacun des~$P_i$, d'où~1).~$\Box$ 

\deux{comment-nullst}
{\bf Commentaires.}

 \trois{nullst-algclos}
 Si le corps~$k$ est algébriquement clos, l'assertion~1) ci-dessus peut se reformuler en 
 disant qu'il existe un~$n$-uplet~$(x_1,\ldots,x_n)$ d'éléments
 {\em du corps~$k$ lui-même}
 tels que~$P_i(x_1,\ldots, x_n)=0$ pour tout~$i$. 
 
 \trois{nullst-algclos-gros}
 Supposons que les trois conditions équivalentes~1), 2) et~3) ci-dessus soient
 satisfaites, et soit~$F$ une extension {\em algébriquement close}
 de~$k$. Il existe alors un~$n$-uplet
$(x_1,\ldots, x_n)\in F^n$
tel que~$P_i(x_1,\ldots, x_n)=0$ pour tout~$i$. 

\medskip
En effet, commençons par choisir~$L$
et~$(x_1,\ldots, x_n)$ comme dans~1). Comme~$L$ est une extension
finie de~$k$, elle admet un~$k$-plongement dans~$F$ (faites l'exercice si
vous ne connaissez pas ce résultat ; on peut le démontrer «à la main» 
par récurrence sur~$[L:k]$, ou 
par des méthodes analogues à celles décrites au~\ref{corps-bourbaki}). 
Ce plongement permet de voir~$L$ comme un sous-corps de~$F$,
et~$(x_1,\ldots, x_n)$ comme un~$n$-uplets d'éléments de~$F$, d'où l'assertion.

\subsection*{Un calcul de dimension de Krull}

\deux{intro-krull-k}
Le but de cette section est de calculer explicitement 
la dimension de
Krull d'une~$k$-algèbre intègre de type fini, en se fondant
de manière cruciale sur le lemme de
normalisation de Noether. On désigne toujours par~$k$ un corps. 

\deux{prop-dimkrull-transc}
{\bf Théorème.}
{\em Soit~$A$ une~$k$-algèbre intègre de type fini. La dimension
de Krull de~$A$ est finie, et coïncide avec le degré de transcendance de~${\rm Frac}\;A$
sur~$k$. }

\medskip
{\em Démonstration.} 
Le lemme de normalisation de Noether assure l'existence
d'éléments~$y_1,\ldots, y_d$ de~$A$, algébriquement indépendants sur~$k$
et tels que~$A$ soit entière sur~$k[y_1,\ldots, y_d]$. De 
plus, l'entier~$d$ est nécessairement égal au degré de transcendance de~${\rm Frac}\;A$
sur~$k$ (\ref{integre-dimtr}). On raisonne désormais par récurrence sur~$d$. 

\trois{krull-tr-d0}
Si~$d=0$ alors~$k[y_1,\ldots,y_d]=k$, et~$A$ est donc une algèbre intègre entière sur~$k$. 
C'est dès lors un corps en vertu du lemme~\ref{alg-ent-corps} ; par conséquent, sa dimension de Krull
est nulle. 

\trois{krull-tr-d}
Supposons~$d>0$ et la proposition vraie pour les entiers~$<d$. Comme~$k[y_1,\ldots, y_d]$
s'injecte dans~$A$ et comme~$A$ est entière sur~$k[y_1,\ldots, y_d]$, il résulte 
de la proposition~\ref{krull-transfert}
que la dimension de Krull de~$A$ est égale à celle de l'algèbre de polynômes
$k[y_1,\ldots, y_d]$. Il suffit donc de montrer que la dimension de Krull
de~$k[y_1,\ldots, y_d]$
est égale à~$d$. 

\medskip
Pour tout~$i\leq d$, notons~$\got p_i$ l'idéal~$(y_1	,\ldots, y_i)$ de~$k[y_1,\ldots, y_d]$. Il est
premier (le quotient~$k[y_1,\ldots, y_d]/\got p_i$ est l'anneau intègre~$k[y_j]_{i<j\leq d}$), et l'on
a ainsi construit une chaîne strictement croissante d'idéaux premiers

$$(0)=\got p_0\subsetneq \got p_1\subsetneq \ldots\subsetneq \got p_d$$ 
(que les inclusions soient strictes résulte du fait évident que~$y_i \notin \got p_{i-1}$ pour tout~$i\geq 1$). 
Par conséquent, la dimension de Krull de~$k[y_1,\ldots, y_d]$ est au moins égale à~$d$. 

\medskip
Nous allons montrer maintenant qu'elle vaut au plus~$d$. Soit
$$\got q_0\subsetneq \got q_1\subsetneq \ldots \subsetneq \got q_m$$
une chaîne strictement croissante d'idéaux premiers de~$A$ ; nous allons prouver que~$m\leq d$, 
ce qui permettra de conclure. 

\medskip
Quitte à rajouter~$(0)$ en début de chaîne, on peut toujours supposer que
l'on a~$\got q_0=(0)$. Si~$m=0$
alors~$m\leq d$ ; on suppose maintenant que~$m>0$. L'idéal~$\got q_1$ est alors non nul, il possède
donc un élément~$f\in k[y_1,\ldots, y_d]$ qui est non constant (puisque~$\got q_1$ est strict). 
Comme~$\got q_1$ est premier, il existe un diviseur irréductible~$g$
de~$f$ qui appartient à~$\got q_1$. La chaîne
$$\got q_1/(g)\subsetneq \got q_2/(g)\subsetneq \ldots \subsetneq \got q_m/(g)$$
est une chaîne strictement croissantes d'idéaux premiers de~$k[y_1,\ldots, y_d]/(g)$. Il suffit
maintenant pour obtenir la majoration souhaitée de prouver que la dimension de Krull
de~$k[y_1,\ldots, y_d]/(g)$
est majorée par~$d-1$. 

\medskip
L'anneau~$B:=k[y_1,\ldots, y_d]/(g)$ est intègre (puisque~$g$ est irréductible, et puisque~$k[y_1,\ldots, y_d]$ est factoriel). 
Il est engendré par les~$\overline {y_i}$ comme~$k$-algèbre. Le corps~${\rm Frac}\;B$ est donc engendré par les~$\overline{y_i}$
comme extension de~$k$ ; en particulier, la famille~$(\overline{y_i})$ contient une base de transcendance de~${\rm Frac}\;B$
sur~$k$. 

\medskip
Par ailleurs, les~$\overline{y_i}$ {\em ne sont pas}
algébriquement indépendants sur~$k$, 
puisque~$g(\overline{y_1},\ldots, \overline{y_d})=0$ par construction. En conséquence, 
$(\overline{y_i})$ n'est pas elle-même une base de transcendance de~${\rm Frac}\;B$ sur~$k$. 
Le degré de transcendance~$\delta$ de~${\rm Frac}\;B$ sur~$k$ est donc
strictement inférieur à~$d$. 
En vertu de l'hypothèse de récurrence, la dimension de Krull de~$B$
est égale à~$\delta\leq d-1$, ce qui achève la démonstration.~$\Box$ 

\deux{comment-krull}
{\bf Commentaires.}

\trois{krull-polyn}
Le théorème ci-dessus assure en particulier que
la dimension de Krull de
la~$k$-algèbre $k[X_1,\ldots, X_d]$ 
(qui
n'est autre du point de vue de la géométrie algébrique
que~«l'anneau des fonctions
sur l'espace affine de dimension~$d$»)
est égale à~$d$ ; on a d'ailleurs démontré explicitement
cette égalité au cours de la preuve. 

\trois{exo-kx-surf}
{\bf Exercice.}
En reprenant les notations du~\ref{krull-tr-d}
ci-dessus, montrez que le degré de transcendance
de~${\rm Frac}\;B$ est {\em exactement}
égal
à~$d-1$.

\chapter{Théorie des faisceaux}

\section{Préfaisceaux et faisceaux}
\markboth{Théorie des faisceaux}{Préfaisceaux et faisceaux}

\deux{intro-cat-c}
Nous allons définir dans ce qui suit les notions de
{\em préfaisceau}, puis plus tard de {\em faisceau}, sur un espace topologique.
Il en existe plusieurs variantes : on peut manipuler des faisceaux ou préfaisceaux d'ensembles, 
mais aussi de groupes, d'anneaux...Aussi, afin 
d'éviter des répétitions fastidieuses ou une profusion pénible d'abréviations
«resp.», nous fixons une catégorie~$\mathsf C$, qui peut être 
$\ens$, $\gp$, $\ann$, $\amod$ ou~$\aalg$ pour un certain anneau~$A$.

\medskip
Nous utiliserons relativement souvent des constructions
faisant appel à la notion de limite inductive {\em filtrante}
({\em cf.}
\ref{syst-filtr-top}) dans la catégorie~$\mathsf C$. Nous attirons votre attention sur le fait qu'une 
telle limite est «indépendante» de la catégorie $\mathsf C$, dans le sens précis suivant : 
son ensemble sous-jacent coïncide avec la limite inductive
ensembliste du diagramme considéré.  

\subsection*{Préfaisceaux} 
\deux{conv-espace-topo}
Soit~$X$ un espace topologique. Soulignons
qu'on ne fait
{\em aucune}
hypothèse sur~$X$ : on ne suppose pas qu'il est séparé, ni même que tous ses points sont fermés...
Cette généralité ne complique ni les définitions, ni les preuves, et elle est indispensable en géométrie algébrique : 
en effet, le plus souvent,
un schéma possède des points non fermés, et n'est {\em a fortiori}
pas topologiquement séparé.

\deux{def-prefaisc}
{\bf Définition.}
Un
{\em préfaisceau}
$\sch F$ sur~$X$
à valeurs dans~$\mathsf C$ 
consiste en les données suivantes. 

\medskip
$\bullet$ Pour tout ouvert~$U$ de~$X$, un objet~$\sch F(U)$
de~$\mathsf C$  dont les éléments
sont parfois appelés les {\em sections}
de~$\sch F$ sur~$U$ ; 

$\bullet$ Pour tout couple~$(U,V)$ d'ouverts de~$X$ avec~$V\subset U$, 
un morphisme
$$r_{U\to V} : \sch F(U)\to \sch F(V)$$ dit
de {\em restriction}, et parfois notée~$s\mapsto s|_V$.

\medskip
Ces données sont sujettes aux deux axiomes suivants : 

\medskip
$\diamond$ $r_{U\to U}={\rm Id}_{\sch F(U)}$ pour tout
ouvert~$U$ de~$X$ ; 

$\diamond$ $r_{V\to W}\circ r_{U\to V}=r_{U\to W}$ pour tout triplet~$(U,V,W)$
d'ouverts de~$X$
avec~$W\subset V\subset U$.

\trois{convention-cat-c}
{\em Convention terminologique.}
Lorsqu'il sera nécessaire de préciser la catégorie dans laquelle nos préfaisceaux
prennent leurs valeurs, nous parlerons de préfaisceaux de groupes, d'ensembles, etc. 
Lorsque nous dirons «préfaisceaux»
sans référence à une catégorie particulière, cela signifiera «préfaisceaux à valeurs dans~$\mathsf C$». 

\trois{variante-prefaisc}
On peut définir un préfaisceau de façon plus concise, qui vous paraîtra
peut-être un peu pédante (mais a l'avantage de se généraliser à bien d'autres cadres 
que celui de la topologie). Soit~$\mathsf{Ouv}_X$ la catégorie dont les objets sont les ouverts 
de~$X$ et les flèches les {\em inclusions} ;
un préfaisceau sur~$X$ est alors simplement un foncteur contravariant
de~$\mathsf{Ouv}_X$ vers~$\mathsf C$.

\medskip
\deux{ex-prefaisc}
{\bf Exemples.}
Nous allons donner quelques exemples de préfaisceaux, que nous allons
décrire en nous contentant de donner leurs valeurs sur les ouverts : la définition 
des flèches de restriction est à chaque fois évidente.

\trois{fonct-cont}
Si~$X$ est un espace topologique,~$U\mapsto \mathscr C^0(U,\RR)$ est un préfaisceau
de~$\RR$-algèbres sur~$X$. 

\trois{fonct-smooth}
Si~$X$ est une variété différentielle,~$U\mapsto \mathscr C^\infty(U,\RR)$ est un préfaisceau 
de~$\RR$-algèbres sur~$X$. 

\trois{fonct-hol}
Si~$X$ est une variété analytique complexe,~$U\mapsto {\mathscr H}(U,\CC)$ est un préfaisceau 
de~$\CC$-algèbres sur~$X$, 
où~${\mathscr H}(U,\CC)$ désigne l'anneau des fonctions holomorphes sur~$U$. 

\trois{pref-const}
Si~$X$ est un espace topologique et~$E$ un ensemble,~$U\mapsto E$ est un préfaisceau d'ensembles
sur~$X$, appelé le {\em préfaisceau constant associé à~$E$.} 
Si~$E$ est un groupe (resp. ...), le préfaisceau constant associé hérite d'une structure naturelle
de préfaisceau de groupes (resp. ...).

\trois{pref-idiot}
Terminons par un exemple un peu plus artificiel que les précédents : si~$X$ est un espace topologique,
les formules
$$U\mapsto \left\{\begin{array}{l} \ZZ\;{\rm si}\;U=X\\
\{0\}\;{\rm sinon}\end{array}\right.$$ définissent un préfaisceau de groupes abéliens sur~$X$. 

\deux{fibres-pref}
Soit~$X$ un espace topologique, soit~$\sch F$ un préfaisceau sur~$X$,
soit~$x$
un point de~$X$
et soit~$\sch V$ 
l'ensemble des voisinages ouverts de~$x$,
qui est filtrant (pour l'ordre {\em opposé}
à l'inclusion). Soit~$\sch D$ le diagramme commutatif filtrant 
$$\sch D=((\sch F(U))_{U\in \sch V}, (r_{U\to V})_{V\subset U}))$$
dans la catégorie~$\mathsf C$. 
La limite inductive de~$\sch D$ est notée
~$\sch F_x$ et est appelée
la
{\em fibre}
de~$\sch F$ en~$x$. 

\medskip
Concrètement, $\sch F_x$ est le quotient
de l'ensemble des couples~$(U,s)$, où~$U\in \sch V$ et~où~$s$ est
une section
de~$\sch F$ sur~$U$, 
par la relation d'équivalence suivante : {\em $(U,s)\sim(V,t)$ si et seulement 
si il existe un voisinage ouvert~$W$ de~$x$ dans~$U\cap V$ tel que~$s|_W=t|_W$.}

De manière un peu plus informelle,~$\sch F_x$ est 
 l'ensemble des sections de~$\sch F$ définies au voisinage de~$x$, deux sections
 appartenant à~$\sch F_x$
 étant considérées comme égales si elles coïncident au voisinage de~$x$. 

\trois{germe-pref}
Soit~$U$ un voisinage ouvert de~$x$ et soit~$s\in \sch F(U)$. L'image
de~$s$ dans~$\sch F_x$ est appelée le
{\em germe de~$s$ en~$x$}
et est notée~$s_x$. 

\trois{ex-germes}
{\em Exemples de fibres.} Si~$\sch F$ est le préfaisceau 
constant associé à un ensemble~$E$
alors~$\sch F_x\simeq E$. 

Supposons que~$X=\CC$, que~$x$ est l'origine et que~$\sch F$
est le préfaisceau des fonctions holomorphes. La fibre~$\sch F_x$ 
s'identifie alors à la~$\CC$-algèbre~$\CC\{z\}$
des séries formelles en~$z$ de rayon~$>0$. 

\deux{mor-pref}
Soit~$X$ un espace topologique. Les préfaisceaux sur~$X$
forment 
une catégorie $\pref X$, les flèches se décrivant
comme suit. Si~$\sch F$ et~$\sch G$ sont
deux tels préfaisceaux, un morphisme de~$\sch F$ vers~$\sch G$
consiste en la donnée, pour tout ouvert~$U$
de~$X$, d'un morphisme de~$\sch F(U)$ vers~$\sch G(U)$, de sorte
que le diagramme
$$\diagram \sch F(U)\dto_{r_{U\to V}}\rto&\sch G(U)\dto^{r_{U\to V}}\\
\sch F(V)\rto&\sch G(V)\enddiagram$$
commute pour tout couple~$(U,V)$ d'ouverts de~$X$ tels que~$V\subset U$. 

\medskip
On dit que~$\phi$ est {\em injectif}
(resp. {\em surjectif})
si~$\sch F(U)\to \sch G(U)$ est injectif (resp. surjectif)
pour tout~$U$.

\trois{mor-pref-pedant}
Si l'on considère les préfaisceaux~$\sch F$ 
et~$\sch G$ 
comme des foncteurs contravariants de~$\mathsf{Ouv}_X$
vers~$\mathsf C$, un morphisme de~$\sch F$ vers~$\sch G$ n'est autre qu'un morphisme
de foncteurs de~$\sch F$ vers~$\sch G$. 

En particulier, un morphisme~$\sch F\to \sch G$ est un isomorphisme
si et seulement si~$\sch F(U)\to \sch G(U)$ est bijectif pour tout~$U$.

\trois{mor-fibre}
Si~$\phi : \sch F\to \sch G$ est un morphisme de préfaisceaux sur~$X$, 
il induit pour tout~$x$ un morphisme~$\phi_x : \sch F_x\to\sch G_x$.

\trois{mor-im}
Soit~$\phi : \sch F\to \sch G$ un morphisme
de préfaisceaux sur~$X$. 
On définit le préfaisceau~${\rm Im}\;\phi$ 
par la formule
$$U\mapsto {\rm Im}\;(\sch F(U)\to \sch G(U))$$
pour tout ouvert~$U$ de~$X$. On dit que~${\rm Im}\;\phi$
est l'{\em image}
de~$\phi$ ; c'est un sous-préfaisceau de~$\sch G$, en un sens évident ; 
il est égal à~$\sch G$ si et seulement si~$\phi$
est surjectif. 

\trois{mor-ker}
{\em On suppose que~$\mathsf C=\gp$
ou~$\amod$.}
Soit~$\phi : \sch F\to \sch G$ un morphisme
de préfaisceaux de groupes
sur~$X$. 
On définit le préfaisceau~${\rm Ker}\;\phi$ 
par la formule
$${U\mapsto {\rm Ker}\;(\sch F(U)\to \sch G(U)) }$$
pour tout ouvert~$U$ de~$X$. On dit que~${\rm Ker}\;\phi$
est le {\em noyau}
de~$\phi$ ; c'est un sous-préfaisceau de~$\sch F$,
et~$\phi$ est injectif si et seulement si
son noyau est le préfaisceau trivial. 

\trois{lim-projind-pref}
Soit~$\sch D$ un diagramme dans
$\pref X$ ; pour tout ouvert~$U$ de~$X$, on note~$\sch D(U)$
le diagramme dans~$\mathsf C$ 
déduit de~$\sch D$ par évaluation en~$U$ de ses constituants.  Il est immédiat
que
$$U\mapsto \limind \sch D(U)\;\;\text{et}\;\;U\mapsto \limproj \sch D(U)$$
définissent deux préfaisceaux sur~$X$,  
et que ceux-ci sont respectivement les limites inductive et projective de~$\sch D$.

\deux{appli-cont-pref}
Soit~$f : Y\to X$ une application continue entre espaces topologiques. 

\trois{f-lower-star}
Soit~$\sch G$ un préfaisceau sur~$Y$. 
La formule
$$U\mapsto \sch G(f^{-1}(U))$$ définit un préfaisceau sur~$X$, 
que l'on note~$f_*\sch G$. 

\trois{f-moins-un}
Soit~$\sch F$ un préfaisceau sur~$X$. 
Si~$V$ est un ouvert de~$Y$, on note~$\got E_V$ l'ensemble
des voisinages ouverts de~$f(V)$, et~$\sch D_V$ le diagramme
commutatif filtrant
$$((\sch F(U))_{U\in \got E_V}, (r_{U\to U'})_{U'\subset U}).$$

\medskip
On note
alors~$f^{-1}\sch F$ le préfaisceau~$V\mapsto \limind \sch D_V$.

\trois{pf-pb-fonct}
On vérifie aussitôt que~$f_*$ est de façon naturelle
un foncteur covariant de
$\pref Y$
vers~$\pref X$, 
et que $f^{-1}$ est de façon naturelle
un foncteur covariant de~$\pref X$
vers~$\pref Y$. On notera que
ces foncteurs sont «indépendants de~$\mathsf C$» :
pour~$f_*$, c'est évident par sa définition même,
et pour~$f^{-1}$, {\em cf.}~\ref{intro-cat-c}.

\trois{compose-pf-pb}
Si~$Z$ est un espace topologique et~$g\colon Z\to Y$ une application continue,
il existe des isomorphismes canoniques de foncteurs
(construisez-les à titre d'exercice !)
$$(f\circ g)_*\simeq f_*\circ g_*\;\;\;\text{et}\;\;\;(f\circ g)^{-1}\simeq g^{-1}\circ f^{-1}.$$

\deux{exemples-pullbacks}
{\bf Deux cas particuliers intéressants.}

\trois{ex-restriction}
Soit~$X$ un espace topologique, soit~$U$ un ouvert de~$X$
et soit~$j$ l'inclusion~$U\hookrightarrow X$. Il résulte immédiatement
de la définition
que pour tout préfaisceau~$\sch F$ sur~$X$, le préfaisceau~$j^{-1}\sch F$
n'est autre que la {\em restriction}~$\sch F|_U$ de~$\sch F$ à~$U$,
c'est-à-dire le préfaisceau~$V\mapsto \sch F(V)$ où~$V$ se contente de parcourir
l'ensemble des ouverts {\em de~$U$}. 

\trois{pull-back-singleton}
Soit~$X$ un espace topologique, soit~$x\in X$ et soit~$i$
la flèche canonique~$\{x\}\hookrightarrow X$. Il résulte immédiatement
des définitions que l'on a pour tout préfaisceau~$\sch F$ sur~$X$
l'égalité~$i^{-1}\sch F(\{x\})=\sch F_x$.

\deux{fibres-bon-pullback}
Soit~$f\colon Y\to X$ une application 
continue entre espaces topologiques.
La définition du foncteur~$f^{-1}$ est sensiblement plus compliquée que celle
de~$f_*$, mais
$f^{-1}$ se comporte
plus simplement que~$f_*$
en ce qui concerne les fibres. En effet, si~$y\in Y$ et si~$x$ désigne son image 
sur~$X$, on dispose pour tout préfaisceau~$\sch F$ sur~$X$
d'un
isomorphisme naturel~$f^{-1}\sch F_y \simeq \sch F_x$ ; 
le lecteur le déduira à peu près formellement de~\ref{compose-pf-pb}
et~\ref{pull-back-singleton}. 

\medskip
Par contre, si~$\sch G$ est un préfaisceau sur~$Y$, la fibre~$f_*\sch G_x$ n'admet pas 
à notre connaissance de description maniable. 

\deux{adj-fmoin-fstar-pre}
Soit~$f: Y\to X$ une application continue entre espaces topologiques. Le
couple~$(f^{-1},f_*)$ est alors un couple de foncteurs adjoints. La preuve
détaillée de ce fait est laissée au lecteur ; nous allons nous
contenter de décrire brièvement les isomorphismes
d'adjonction. On se donne 
un préfaisceau~$\sch F$ sur~$X$ et un préfaisceau~$\sch G$ sur~$Y$. 

\trois{adj-pre1}
{\em Description d'un isomorphisme~$\hom(f^{-1}\sch F,\sch G)\simeq \hom (\sch F, f_*\sch G)$}.
Soit~$\phi$ un morphisme de~$f^{-1}\sch F$ vers~$\sch G$. On lui associe
un morphisme~$\psi$ de~$\sch F$ vers~$f_*\sch G$ comme suit. 

\medskip
Soit~$U$ un ouvert de~$X$, et soit~$s\in \sch F(U)$. Par définition de~$f^{-1}\sch F$, le couple
$(U,s)$ donne lieu à une section~$s'$ de~$f^{-1}\sch F$ sur~$f^{-1}(U)$. Son image~$\phi(s')$
est une section de~$\sch G$ sur~$f^{-1}(U)$, et donc par définition une section de~$f_*\sch G$ sur~$U$. 
On pose alors~$\psi(s)=\phi(s')\in f_*\sch G(U)$. 

\trois{adj-pre2}
{\em Description de l'isomorphisme $\hom (\sch F, f_*\sch G)\simeq \hom (f^{-1}\sch F,\sch G)$
réciproque du précédent.}
Soit~$\psi$ un morphisme de~$\sch F$ vers~$f_*\sch G$. On lui associe un morphisme
$\phi$ de $f^{-1}\sch F$ vers~$\sch G$
comme suit. 

\medskip
Soit~$V$ un ouvert de~$Y$, et soit~$\tau\in f^{-1}\sch F(V)$. La section~$\tau$ 
provient d'une section~$t$ de~$\sch F$ définie sur un ouvert~$U$ de~$X$ contenant~$f(V)$. 
Son image~$\psi(t)$ est une section de~$f_*\sch G$ sur~$U$, c'est-à-dire par définition
une section de~$\sch G$ sur~$f^{-1}(U)\supset V$. On vérifie que
l'élément~$\psi(t)|_V$
de~$\sch G(V)$
ne dépend que de~$\tau$ et pas du choix de~$(U,t)$, et il est dès lors
licite de poser~$\phi(\tau)=\psi(t)|_V$.

\subsection*{Faisceaux}

\deux{def-fasc}
{\bf Définition.}
Soit~$X$ un espace topologique. 
On dit qu'un préfaisceau~$\sch F$ sur~$X$
(à valeurs, comme d'habitude, dans la catégorie~$\mathsf C$
de~\ref{intro-cat-c}), 
est un {\em faisceau}
s'il satisfait la condition suivante : 
{\em pour tout ouvert~$U$ de~$X$, pour tout
recouvrement ouvert~$(U_i)_{i\in I}$ de~$U$ 
et toute famille~$(s_i)\in \prod_{i\in I} \sch F(U_i)$ telle
que~$s_i|_{U_i\cap U_j}=s_j|_{U_i\cap U_j}$ 
pour tout~$(i,j)$, il existe une et une seule section~$s$
de~$\sch F$ sur~$U$ 
telle que~$s|_{U_i}=s_i$ pour tout~$i$.}

\trois{def-recoll} Si~$s$ et les~$s_i$
sont comme ci-dessus, on dit que~$s$ est obtenue
{\em par recollement des~$s_i$.} 

\trois{fasc-vide}
Soit~$\sch F$ un faisceau sur~$X$.
L'ouvert vide de~$X$ est alors recouvert par
la {\em famille vide}
d'ouverts ; en appliquant à celle-ci 
et à la {\em famille vide} de sections la définition d'un 
faisceau, on voit qu'il existe une et une seule section
de~$\sch F$ sur~$\varnothing$. Ainsi~$\sch F(\varnothing)$
est un singleton. Lorsque~$\mathsf C=\gp$ (resp.~$\amod$, $\ann$
ou~$\aalg$) cela signifie que~$\sch F(\varnothing)$ est le groupe trivial
(resp. le module, l'anneau ou l'algèbre nuls). 

\medskip
Le lecteur que rebuterait (bien à tort) ce style de gymnastique mentale
peut prendre l'égalité~$\sch F(\varnothing)=\{*\}$ comme un
axiome supplémentaire à rajouter à la définition d'un faisceau.

\deux{mor-fasc}
On voit les faisceaux sur~$X$
comme une sous-catégorie {\em pleine}
$\fasc X$
de celle des préfaisceaux sur~$X$ : si~$\sch F$ et~$\sch G$ sont deux faisceaux sur~$X$,
un morphisme de faisceaux de~$\sch F$ vers~$\sch G$ est simplement un morphisme
de préfaisceaux de~$\sch F$ vers~$\sch G$.

\deux{fasc-germes}
Soit~$X$ un espace topologique. Il résulte
immédiatement de la définition qu'une section~$s$
d'un faisceau sur~$X$ est entièrement déterminée
par la famille~$(s_x)_{x\in X}$
de ses germes. On en déduit qu'un 
morphisme
d'un préfaisceau~$\sch F$ sur~$X$ vers un faisceau~$\sch G$
sur~$X$ 
est entièrement déterminé par
la famille~$(\sch F_x\to \sch G_x)_{x\in X}$ de
morphismes induits au niveau des fibres. 

\deux{ex-fasc}
{\bf Exemples.}

\trois{fonct-cont-fasc}
Si~$X$ est un espace topologique,~$U\mapsto \mathscr C^0(U,\RR)$ est un faisceau
de~$\RR$-algèbres sur~$X$. 

\trois{fonct-smooth-fasc}
Si~$X$ est une variété différentielle,~$U\mapsto \mathscr C^\infty(U,\RR)$ est un faisceau 
de~$\RR$-algèbres sur~$X$. 

\trois{fonct-hol-fasc}
Si~$X$ est une variété analytique complexe,~$U\mapsto {\mathscr H}(U,\CC)$ est un faisceau 
de~$\CC$-algèbres sur~$X$, 
où~${\mathscr H}(U,\CC)$ désigne l'anneau des fonctions holomorphes sur~$U$. 

\trois{pref-const}
Si~$X$ est un espace topologique et~$\{*\}$ un singleton,
le préfaisceau constant~$U\mapsto \{*\}$ est un faisceau
d'ensembles (et aussi d'ailleurs de groupes, anneaux, $A$-modules ou~$A$-algèbres si l'on y tient).

\trois{noyau-fasc}
Si~$\phi : \sch F \to \sch G$ est un morphisme
de faisceau de groupes sur un espace topologique~$X$, 
le préfaisceau noyau~${\rm Ker}\;\phi$ (défini au~\ref{mor-im}) 
est un faisceau. 

\trois{rest-fasc}
Soit~$X$ un espace topologique. Pour tout faisceau~$\sch F$ sur~$X$,
et tout ouvert~$U$ de~$X$, la restriction de~$\sch F$ à~$U$ est un faisceau
sur~$U$. 

\deux{contrex-fasc}
{\bf Contre-exemples.}

\trois{contre-idiot}
Soit~$X$ un espace topologique. Si~$X$ admet
un recouvrement~$(U_i)$
par des ouverts {\em stricts}
\footnote{Il existe des espaces
topologiques naturels du point de vue de la géométrie 
algébrique qui sont non vides et pour lesquels cette condition n'est pas vérifiée : 
par exemple l'espace~$\{a,b\}$ dont les ouverts
sont~$\varnothing, \{a\}$ et~$\{a,b\}$. Notons
toutefois que si~$X$ est 
de cardinal au moins 2 et si tous ses points
sont fermés,
$X$ est bien recouvert par des ouverts
stricts.}, le préfaisceau de groupes
$$U\mapsto \left\{\begin{array}{l} \ZZ\;{\rm si}\;U=X\\
\{0\}\;{\rm sinon}\end{array}\right.$$ n'est pas un faisceau : les
sections globales~$1$ et~$0$ sont distinctes, mais ont toutes deux mêmes
restrictions à chacun des~$U_i$. 

\trois{const-pas-fasc}
Soit~$X$ un espace topologique et soit~$E$ un ensemble
{\em non singleton}. Le préfaisceau d'ensembles
constant associé à~$E$
envoie en particulier~$\varnothing$
sur~$E$ et n'est donc pas un faisceau en vertu de~\ref{fasc-vide}. 

\deux{pre-contrex-log}
{\bf Autour des images préfaisceautiques.}
Soit~$X$ un espace topologique et soit~$\sch F$ un faisceau sur~$X$. 
On vérifie immédiatement qu'un sous-préfaisceau~$\sch F'$ de~$\sch F$ est un faisceau
si et seulement si «l'appartenance à~$\sch F'$ est une propriété locale», {\em i.e.}
si pour tout ouvert~$U$ de~$X$ et tout recouvrement ouvert~$(U_i)$
de~$U$, une section~$s$ de~$\sch F$ sur~$U$
appartient à~$\sch F'(U)$ dès que~$s|_{U_i}$ appartient 
à~$\sch F'(U_i)$ pour tout~$i$. 

\trois{im-inj-faisc}
Soit~$\phi : \sch F\to \sch G$ un
morphisme
de faisceaux, et soit~$s$ une section de~$\sch G$ sur un ouvert~$U$
de~$X$ ; soit~$(U_i)$ un recouvrement ouvert de~$U$ tel que
$s|_{U_i}$ appartiennent à~${\rm Im}\;\phi$ pour tout~$i$. 
Choisissons pour tout~$i$ un antécédent~$t_i$
de~$s|_{U_i}$ dans~$\sch F(U_i)$. 

\medskip
{\em Supposons que~$\phi$ est injective
(il n'y a alors qu'un choix
possible pour les~$t_i$).}
Dans ce cas,~$t_i|_{U_i\cap U_j}$ et~$t_j|_{U_i\cap U_j}$ 
coïncident pour tout~$(i,j)$, puisque ce sont deux antécédents
de~$s|_{U_i\cap U_j}$ par une flèche injective. Il s'ensuit que les~$t_i$ se recollent
en une section~$t$ de~$\sch F$ sur~$U$, qui satisfait par construction 
l'égalité~$\phi(t)=s$. Ainsi, ${\rm Im}\;\phi$ est un sous-faisceau de~$\sch G$. 

\medskip
{\em On ne suppose plus~$\phi$ injective.}
Dans ce cas, rien ne garantit {\em a priori}
que le système~$(t_i)$
d'antécédents puisse être choisi de sorte que les~$t_i$ se recollent, et le contre-exemple
ci-dessous montre que~${\rm Im\;\phi}$
n'est pas un faisceau en général.

\trois{contrex-log}
Soit~$\sch H$ le faisceau
des fonctions holomorphes sur~$\CC$, et soit~$d$
la dérivation~$f\mapsto f'$, qui est un endomorphisme
du faisceau de~$\CC$-espaces vectoriels~$\sch H$. 
La fonction~$z\mapsto 1/z$
appartient à~$\sch H(\CC\ti)$. Comme n'importe quelle fonction
holomorphe, elle admet localement des primitives et appartient donc localement
à~${\rm Im}\;d$. Par contre, elle n'admet pas de primitive sur~$\CC\ti$ (il n'existe pas
de logarithme complexe continu sur~$\CC\ti$) ; elle n'appartient donc pas
à~${\rm Im}\;d(\CC\ti)$. En conséquence, 
{\em ${\rm Im}\;d$ n'est pas un faisceau}. 

\deux{intro-faisceautisation}
On observe là un cas particulier d'un phénomène général : 
lorsqu'on applique aux faisceaux des constructions «naïves»
(c'est-à-dire définies ouvert par ouvert), on obtient des préfaisceaux
qui n'ont en général aucune raison d'être des faisceaux. Qu'à cela ne tienne : 
on y remédie grâce au procédé de {\em faisceautisation}, que nous allons maintenant
décrire. Il consiste {\em grosso modo}
à modifier un préfaisceau donné pour en faire un
faisceau, sans
l'altérer davantage que ne le requiert cet objectif ; cela
va se traduire rigoureusement en terme de foncteur à représenter. 

\subsection*{La faisceautisation}

\deux{def-faisceautis}
{\bf Proposition-définition.}
{\em Soit~$X$ un espace topologique et soit~$\sch F$ un préfaisceau sur~$X$. Le foncteur
covariant de~$\fasc X$ vers~$\ens$
qui envoie un faisceau~$\sch G$ sur~$\hom_{\mathsf{Pref}}(\sch F,\sch G)$ est
représentable par un couple~$(\widehat {\sch F}, \pi : \sch F\to \widehat{\sch F})$.
Le faisceau~$\widehat{\sch F}$ est appelé le {\em faisceautisé}
de~$\sch F$, ou encore le {\em faisceau associé}
à~$\sch F$.}

\deux{faisceauti-propuniv}
{\bf Formulation équivalente.}
Il revient au même de dire qu'il existe un faisceau~$\widehat{\sch F}$
sur~$X$
et un
morphisme de préfaisceaux~$\pi : \sch F\to \widehat{\sch F}$ tel que pour tout faisceau~$\sch G$
sur~$X$ et tout morphisme de préfaisceaux~$\phi : \sch F\to \sch G$, il existe un unique
morphisme~$\psi : \widehat{\sch F}\to \sch G$ faisant
commuter le diagramme $$\diagram \sch F\dto_\pi\rto^\phi&\sch G\\
\widehat{\sch F}\urto_\psi&\enddiagram\;\;\;.$$

\deux{demo-faisceautis}
{\em Démonstration de la proposition~\ref{def-faisceautis}}.
L'idée qui préside à la construction de~$\widehat{\sch F}$ est simple : 
puisqu'une section d'un faisceau est caractérisée par ses germes, 
on va {\em définir}
une section de~$\widehat{\sch F}$ comme
une collection «raisonnable»
de germes de sections de~$\sch F$. 
Plus précisément, soit~$\got G$ l'ensemble
$\coprod_{x\in X}\sch F_x$. Pour tout
ouvert~$U$ de~$X$, on note~$\widehat{\sch F}(U)$ l'ensemble
des applications~$f : U\to \got G$ satisfaisant
les conditions suivantes : 

\medskip
a) $f(x)\in \sch F_x$ pour tout~$x\in U$ ; 

b) pour tout~$x\in U$ il existe un voisinage ouvert~$V$
de~$x$ dans~$U$ et une section~$s$ de~$\sch F$ sur~$V$
telle que~$f(y)=s_y$ pour tout~$y\in U$. 

\medskip
Il est immédiat que~$\widehat{\sch F}$ est un faisceau
(de fonctions à valeurs dans~$\got G$). On note~$\pi$
le morphisme de préfaisceaux de~$\sch F$ vers~$\widehat{\sch F}$
qui pour tout ouvert~$U$
de~$X$ associe à un élément~$s$
de~$\sch F(U)$
l'élément~$(x\mapsto s_x)$
de~$\widehat{\sch F}(U)$. 

\medskip
Si~$U$ est un ouvert de~$X$, si~$x\in U$
et si~$f\in \widehat{\sch F}(U)$, il existe un voisinage
ouvert~$V$
de~$x$ dans~$U$ et une section~$t\in \sch F(U)$ telle que~$f|_V=\pi(t)$ : 
c'est une reformulation de~b).

\trois{faisceautis-fibr}
Nous allons
maintentant
établir un résultat qui nous sera utile pour 
établir la propriété universelle de~$(\widehat{\sch F},\pi)$ et qui présente 
par ailleurs un intérêt intrinsèque : {\em pour tout~$x\in X$,
la flèche~$\pi_x : \sch F_x\to \widehat{\sch F}_x$
est {\em bijective}.}
Soit donc~$x\in X$.

\medskip
{\em Surjectivité de~$\sch F_x\to \widehat{\sch F}_x$.}
Soit~$\sigma\in \widehat{\sch F}_x$ ; par définition, $\sigma$ est le germe
en~$x$
d'une section~$\tau$ de~$\widehat {\sch F}$
définie
sur un voisinage ouvert~$U$
de~$x$, et l'on peut restreindre~$U$
de sorte qu'il existe~$t\in \sch F(U)$ vérifiant l'égalité~$\pi(t)=\tau$.
On a alors~$\pi_x(t_x)=\tau_x=\sigma$ ; ainsi,
$\pi_x$
est surjective. 

\medskip
{\em Injectivité de~$\pi_x$.}
Soient~$s$ et~$s'$
deux éléments de~$\sch F_x$ dont les images 
dans~$\widehat{\sch F}_x$ coïncident. Il existe un voisinage
ouvert~$U$ de~$x$ et deux sections~$t$ et~$t'$
de~$\sch F$ sur~$U$ telles que~$t_x=s$ et~$t'_x=s'$ ; comme
$\pi(t)_x=\pi(t')_x$ par hypothèse,  on peut restreindre~$U$
de sorte que~$\pi(t)=\pi(t')$. Par définition de~$\pi$, cela signifie que~$t_y=t'_y$ pour tout~$y\in U$. 
En appliquant ceci avec~$y=x$, il vient~$s=s'$, et
la flèche~$\sch F_x\to \widehat{\sch F}_x$
est injective. 

\medskip
{\em Calcul des germes d'une section de~$\widehat{\sch F}$.}
Soit~$U$ un ouvert de~$X$, soit~$f\in \widehat{\sch F}(U)$ et soit~$x\in U$. 
Il existe un voisinage ouvert~$V$ de~$x$ dans~$U$ et une section~$t$
de~$\sch F$ sur~$V$
telle
que~$f|_V=\pi(t)$. On a alors~$f_x=\pi_x(t_x)$ ; et par ailleurs
la définition même de~$\pi$ assure que~$t_x=f(x)$.

En conséquence, $f_x=\pi_x(f(x))$ ; on peut reformuler
cette égalité en disant que
{\em $f_x=f(x)$ modulo la bijection naturelle
$\pi_x : \sch F_x\simeq \widehat{\sch F}_x$}. 

\trois{preuve-propuniv-faisceautis}
Prouvons maintenant que~$(\widehat{\sch F}, \pi)$ satisfait la propriété universelle requise. Soit~$\sch G$ un faisceau sur~$X$,
et soit~$\phi : \sch F\to \sch G$ un morphisme de préfaisceaux. Il s'agit de montrer l'existence et l'unicité de~$\psi: \widehat{\sch F}\to \sch G$
tel que le diagramme
$$\diagram \sch F\dto_\pi\rto^\phi&\sch G\\
\widehat{\sch F}\urto_\psi&\enddiagram$$
commute. 

\medskip
{\em Unicité de~$\psi$.}
Si~$\psi$ est un morphisme faisant commuter le diagramme, 
on a pour tout~$x\in X$ l'égalité
$\psi_x\circ \pi_x=\phi_x$, et donc~$\psi_x=\phi_x\circ \pi_x^{-1}$ (rappelons que~$\pi_x$ est bijective, 
{\em cf.}~\ref{faisceautis-fibr}
{\em supra}). 
L'unicité de~$\psi$ découle alors du fait qu'un morphisme d'un préfaisceau vers un faisceau est entièrement
déterminé par son effet sur les fibres. 

\medskip
{\em Existence de~$\psi$.}
Soit~$U$ un ouvert de~$X$ et soit~$f\in \widehat{\sch F}(U)$. 
Nous allons tout d'abord montrer qu'il
existe une section~$t$ de~$\sch G$
sur~$U$, nécessairement unique, telle que~$t_x=\phi_x(f(x))$
pour tout~$x\in U$. 

Il existe un recouvrement ouvert~$(U_i)$
de~$U$ et pour tout~$i$ une section~$s_i$ de~$\sch F$
sur~$U_i$ telle que~$f(x)=s_{i,x}$ pour tout~$x\in U_i$. Posons~$t_i=\phi(s_i)$. 
On a alors par construction pour tout~$x\in U_i$ l'égalité~$t_{i,x}=\phi_x(s_{i,x})=\phi_x(f(x))$. 
Cette dernière écriture ne dépend plus de~$i$. Il s'ensuit que si~$x\in U_i\cap U_j$ alors~$t_{i,x}=t_{j,x}$. 
Comme~$\sch G$ est un faisceau, ceci entraîne que~$t_i|_{U_i\cap U_j}=t_j|_{U_i\cap U_j}$ pour tout~$(i,j)$, 
puis que les~$t_i$ se recollent en une section~$t$ de~$\sch G$ ; on a bien
par construction~$t_x=\phi_x(f(x))$
pour tout~$x\in U$. 

On pose alors~$\psi(f)=t$. On vérifie 
immédiatement que
$\psi$ est bien un morphisme de préfaisceaux. On a pour tout ouvert~$U$ de~$X$, tout point~$x$ de~$U$
et toute section~$f\in \sch F(U)$ les égalités
$$\underbrace{\psi_x(f_x)=\phi_x(f(x))}_{{\rm par~construction~de}\;\psi}=\phi_x(\pi_x^{-1}(f_x)).$$ 
On a donc  $\psi_x=\phi_x\circ \pi_x^{-1}$ pour tout~$x\in X$, 
soit encore~$\psi_x\circ \pi_x=\phi_x$. Comme un morphisme d'un 
préfaisceau vers un faisceau est entièrement déterminé par son effet sur les fibres,
il vient~$\psi \circ \pi=\phi$, ce qui achève la démonstration.~$\Box$ 

\deux{adj-faisceautis}
Soit~$\sch F$ un préfaisceau sur~$X$, et soit~$\sch G$ un faisceau sur~$X$. On a
une bijection naturelle
$$\hom_{\pref X}(\sch F,\sch G)\simeq \hom_{\fasc X}(\widehat{\sch F},\sch G),$$
fonctorielle en~$\sch F$ et~$\sch G$. En conséquence,
$\sch F\mapsto \widehat{\sch F}$ est l'adjoint à gauche du foncteur d'inclusion
de~$\pref X$ dans~$\fasc X$.

\deux{exem-faisceautis}
{\bf Exemples.}

\trois{triv-faisceautis} Si~$\sch F$ est un faisceau
sur un espace topologique~$X$  alors
$\widehat {\sch F}=\sch F$ ; on le déduit ou bien de la construction
de~$\widehat{\sch F}$, ou bien du fait que~$(\sch F,{\rm Id}_{\sch F})$
satisfait visiblement la propriété universelle requise. 

\trois{faisceautis-const}
Soit~$E$ un ensemble et soit~$X$ un espace topologique. 
Pour tout~$x$
appartenant à~$X$, la fibre en~$x$ du préfaisceau
d'ensembles
constant~$U\mapsto E$ est égale à~$E$. Son faisceautisé 
s'identifie donc, d'après notre construction, 
au faisceau des applications localement
constantes sur~$X$ à valeurs dans~$E$. On l'appelle le {\em 
faisceau constant
associé à~$E$}, et on le note~$\underline E$. 

\medskip
Notez que si~$X$ est localement connexe, le faisceau~$\underline E$ envoie
un ouvert~$U$ de~$X$ sur~$E^{\pi_0(U)}$, où~$\pi_0(U)$ est l'ensemble
des composantes connexes de~$U$. Mais en général la description de~$\underline E$
est un peu plus compliquée : se donner une section de~$\underline E$ sur un ouvert~$U$ revient à 
se donner une partition de~$U$ en ouverts fermés, et à assigner à chacun d'eux un élément de~$E$.  

\trois{faisceautis-sous}
Soit~$\sch F$ un faisceau, et soit~$\sch G$ un sous-{\em préfaisceau}
de~$\sch F$. L'inclusion~$\sch G \hookrightarrow \sch F$
induit un morphisme~$\widehat{\sch G}\hookrightarrow \sch F$
dont on vérifie (l'exercice est laissé au lecteur)
qu'il est injectif et identifie~$\widehat{\sch G}$ au
sous-faisceau de~$\sch F$
formé des sections qui appartiennent localement à~$\sch G$. En termes un
peu plus précis,
$\widehat{\sch G}(U)$ est pour tout ouvert~$U$
l'ensemble des sections~$s\in \sch F(U)$ satisfaisant
la condition suivante : pour tout~$x\in U$,
il existe un voisinage ouvert~$V$ de~$x$ dans~$U$ tel que~$s_{|V}\in \sch G(V)$. 

\subsection*{Images faisceautiques et défaut d'exactitude}

\deux{im-faisceaut}
{\bf Image faisceautique}.
Soit~$X$ un espace topologique. Soit
$\phi : \sch F\to \sch G$
un morphisme de faisceaux sur~$X$. Ce qu'on notera 
désormais ${\rm Im}\;\phi$, ce sera l'image
{\em faisceautique}
de~$\phi$, c'est-à-dire le {\em faisceau
associé} à son image préfaisceautique
$U\mapsto \{\phi(s), s\in \sch F(U)\}\subset \sch G(U)$. 
D'après
le~\ref{faisceautis-sous}
ci-dessus, une section de~${\rm Im}\;\phi$
sur un ouvert~$U$ de~$X$ est une section
de~$\sch G$ sur~$U$ qui admet {\em localement}
un antécédent par~$\phi$. 

\deux{ex-im-faisc}
{\bf Exemples.}

\trois{im-faisceautique-inj}
Si~$\phi$ est un morphisme
{\em injectif}
entre deux faisceaux sur un espace topologique, ${\rm Im}\;\phi$
coïncide avec l'image préfaisceautique
de~$\phi$, puisque celle-ci
est déjà un faisceau 
(\ref{im-inj-faisc}). 

\trois{im-faisceautique-holo}
Ce
n'est pas le cas en général :
on a vu
au~\ref{contrex-log}
que si~$\sch H$
désigne le faisceau des fonctions
holomorphes sur~$\CC$
et~$d: \sch H\to \sch H$ la dérivation,
l'image préfaisceautique de~$d$ n'est pas
un faisceau. 

Par ailleurs, comme toute fonction 
holomorphe admet {\em localement}
une primitive, 
on a~${\rm Im}\;d=\sch H$.

\deux{surj-faisc}
Soit~$X$ un espace topologique et
soit~$\phi : \sch F\to \sch G$ un morphisme
de faisceaux sur~$X$. De même
qu'on a modifié la définition de l'image, on modifie celle
de la surjectivité : on dira désormais que~$\phi$ est
surjectif si son image (faisceautique !) est égale à~$\sch G$. 
Ainsi, {\em la dérivation est un endomorphisme surjectif du faisceau 
des fonctions holomorphes sur~$\CC$}. 

\trois{inj-faisc}
Insistons sur le fait qu'on n'a par contre pas modifié la définition de l'injectivité, 
ni celle du noyau pour un morphisme de faisceaux de groupes :
elles restent définies ouvert par ouvert 
comme pour les morphismes de préfaisceaux.

\trois{rappel-iso-faisc}
Rappelons que~$\phi$ est un isomorphisme si
et seulement si
la flèche~$\sch F(U)\to \sch G(U)$ 
est bijective pour tout
ouvert~$U$ de~$X$ (\ref{mor-pref-pedant}). 
Il s'ensuit, en vertu de~\ref{im-faisceautique-inj},
que~$\phi$ est un isomorphisme si et seulement si il est à la fois
injectif et surjectif.

\trois{inj-surj-fibres}
Nous laissons au lecteur le soin de démontrer
l'assertion suivante : $\phi$ est injectif (resp. surjectif, resp. bijectif)
si et seulement si~$\phi_x : \sch F_x\to \sch G_x$
est injectif (resp. surjectif, resp. bijectif)
pour tout~$x\in X$. Elle a l'avantage de remettre injectivité et surjectivité
sur le même plan, alors qu'on pouvait avoir l'impression d'une certaine
dissymétrie entre elles -- l'injectivité étant définie de manière naïve quand
la surjectivité ne se teste qu'après faisceautisation de l'image.

\deux{exacte-faisc}
{\bf Suites exactes de faisceaux de groupes.}
Soit~$X$ un espace topologique, 
soient~$A$ et~$B$ deux éléments de~$\ZZ\cup\{-\infty,+\infty\}$
et soit

$$S=\ldots \to \sch F_i\to \sch F_{i+1}\to \sch F_{i+2}\to \ldots$$
une suite de morphismes de
faisceaux de groupes sur~$X$, où~$i$ 
parcourt l'ensemble~$I$ des entiers relatifs compris entre~$A$ et~$B$.

Soit~$i$ un élément de~$I$ tel que~$i-1$ et~$i+1$ appartiennent à~$I$. On dit
que la suite~$S$
est {\em exacte en~$\sch F_i$}
si le noyau de~$\sch F_i\to \sch F_{i+1}$
est égal à l'image de~$\sch F_{i-1}\to \sch F_i$. On dit que~$S$ est {\em exacte}
si elle est exacte en~$\sch F_i$ pour tout~$i$ tel que~$i-1$
et~$i+1$ appartiennent à~$I$ (les indices extrêmes, s'ils existent,
ne comptent donc pas). 

Il résulte de la définition que dans une suite
exacte, la composée de deux flèches successives est toujours nulle. 

\trois{exa-fa-stalk}
Démontrez que l'exactitude d'une suite de faisceaux de groupes
se teste sur les fibres. 

\trois{ex-exafa-d}
La suite
$$\diagram \sch F'\rto^f&\sch F\rto^g&\sch F''\rto&0\enddiagram $$
est exacte si et seulement si~$g$ est surjective et~${\rm Ker}\;g={\rm Im}\;f$. 

\trois{ex-exafa-g}
La suite
$$\diagram 0\rto &\sch F'\rto^f&\sch F\rto^g&\sch F''\enddiagram$$
est exacte si et seulement si~$f$ est injective
et~${\rm Ker}\;g={\rm Im}\;f$. 

\trois{ex-exafa-dg}
La suite
$$\diagram 0\rto &\sch F'\rto^f&\sch F\rto^g&\sch F''\rto&0\enddiagram$$
est exacte si et seulement si~$f$ est injective, $g$ est
surjective
et~${\rm Ker}\;g={\rm Im}\;f$.

\deux{ex-suitex-faisc}
{\bf Exemples.} On note~$\sch H$
le faisceau des fonctions holomorphes 
sur~$\CC$, et~$\sch H\ti$ celui des fonctions holomorphes
inversibles. 

\trois{exa-deri-holo}
On a vu au~\ref{im-faisceautique-holo}
que la dérivation~$d$
est un endomorphisme surjectif de~$\sch H$
(vu comme faisceau de~$\CC$-espaces vectoriels).
Par ailleurs, une fonction
holomorphe sur un ouvert~$U$ de~$\CC$ a une dérivée nulle
si et seulement si elle est constante
{\em sur chaque composante connexe de~$U$.}
On a donc une suite exacte 
naturelle de faisceaux de~$\CC$-espaces vectoriels
sur l'espace topologique $\CC$ :

$$\diagram
0\rto&\underline{\CC}\rto &\sch H\rto^d&\sch H\rto&0\enddiagram.$$

Décrivons la suite exacte qui lui correspond au niveau des fibres. Soit~$x$
un point de~$\CC$.  
Le développement en série entière en la variable~$u=z-x$ fournit un isomorphisme
de~$\CC$-algèbres entre~$\sch H_x$ et l'anneau $\CC\{u\}$ des séries
entières de rayon~$>0$. La fibre en~$x$ de la suite exacte précédente
est la suite exacte de~$\CC$-espaces vectoriels
$$\diagram 0\rto&  {\CC}\rto &\CC\{u\} \rto^{\partial/\partial u}&\CC\{u\}\rto&0\enddiagram.$$

\trois{exa-exp-holo}
Toute fonction holomorphe inversible est localement le
logarithme d'une fonction holomorphe. Par ailleurs,
l'exponentielle d'une fonction holomorphe~$f$
sur un ouvert~$U$
de~$\CC$ est égale à~$1$ si et seulement si~$f$
est constante de valeur appartenant à~$2i\pi\ZZ$
sur chaque composante connexe de~$U$. 
On a donc une suite exacte 
naturelle de faisceaux de
groupes abéliens
sur l'espace topologique $\CC$,
appelée {\em suite exponentielle} : 

$$\diagram
0\rto&\underline{2i \pi \ZZ}\rto &\sch H\rto^{\exp}&\sch H\ti \rto&1\enddiagram.$$

Décrivons la suite exacte qui lui correspond au niveau des fibres. Soit~$x$
un point de~$\CC$.  
Le développement en série entière en la variable~$u=z-x$ fournit un isomorphisme
de~$\CC$-algèbres entre~$\sch H_x$ et l'anneau $\CC\{u\}$ des séries
entières de rayon~$>0$. La fibre en~$x$ de la suite exacte précédente
est la suite exacte de groupes abéliens 
$$\diagram 0\rto&  2i\pi\ZZ \rto &\CC\{u\} \rto^{\exp}&\CC\{u\}\ti\rto&1\enddiagram.$$

\deux{exactitude-gamma}
Soit~$X$ un espace topologique et soit
$$\diagram 1\rto&\sch F'\rto^u& \sch F\rto^v &\sch F''\enddiagram$$ une suite exacte de faisceaux de groupes
sur~$X$. Pour tout ouvert~$U$ 
de~$X$, la suite~
$$1\to \sch F'(U)\to \sch F(U)\to \sch F''(U)$$ est exacte. En effet :

\medskip
$\bullet$ $\sch F'(U)\to \sch F(U)$
est injective par définition de l'injectivité d'un morphisme de faisceaux ; 

$\bullet$ l'exactitude en~$\sch F$ signifie que~${\rm Ker}\;v={\rm Im}\;u$. Mais 
comme~$u$ est injective, ${\rm Im}\;u$ coïncide avec l'image préfaisceautique
de~$u$, et l'exactitude en~$\sch F(U)$ en découle aussitôt. 

\trois{gamma-pas-exact}
Ainsi, le foncteur $\sch F\mapsto \sch F(U)$ est exact à gauche. Il n'est pas
exact en général, car il ne transforme pas nécessairement les surjections en 
surjections, comme en atteste notre sempiternel contre-exemple : 
si~$\sch H$ désigne le faisceau des fonctions holomorphes sur~$\CC$
la dérivation~$d : \sch H\to \sch H$ est surjective, mais l'application induite
$\sch H(\CC\ti)\to \sch H(\CC\ti)$ ne l'est pas, puisque son 
image ne contient pas~$z\mapsto 1/z$. 

Donnons-en un autre, qui traduit le même phénomène
(l'absence de logarithme continu sur~$\CC\ti$) : 
la flèche~$\exp : \sch H\to \sch H\ti$ est surjective, mais
l'application induite~$\sch H(\CC\ti)\to \sch H\ti(\CC\ti)$ ne l'est pas, car son image
ne contient pas l'identité. 

\trois{philo-pas-exact}
Ce défaut d'exactitude -- dont la mesure précise constitue
l'objet de ce qu'on appelle la {\em cohomologie} -- est, en un sens, le 
principal intérêt de la théorie des faisceaux : il traduit en effet les difficultés de
recollement d'antécédents, elles-mêmes liées à la «forme»
de l'espace topologique
considéré (présence ou non de «trous», etc.) ; il permet
donc d'une certaine manière de décrire
cette forme algébriquement. 

Ainsi, les deux contre-exemples du~\ref{gamma-pas-exact}
ci-dessus sont intimement liés au fait que~$\CC\ti$ n'est pas simplement
connexe. 

\deux{quot-faisc}
{\bf Quotients.}
Soit~$X$ un espace topologique, 
soit~$\sch G$ un faisceau de groupes sur~$X$, et soit~$\sch H$
un sous-faisceau de groupes de~$\sch G$. Le {\em faisceau quotient}
$\sch G/\sch H$ est le faisceau {\em d'ensembles}
associé
au préfaisceau $U\mapsto \sch G(U)/\sch H(U)$, et l'on dispose
d'une surjection naturelle~$\sch G\to \sch G/\sch H$. 

\medskip
{\em Supposons que~$\sch H$ est distingué
dans~$\sch G$ ({\em i.e.}
~$\sch H(U)$ est distingué dans
$\sch G(U)$ pour tout ouvert~$U$).}
Le faisceau quotient~$\sch G/\sch H$ hérite alors d'une
structure naturelle de faisceau de groupes, et~$\sch G\to \sch G/\sch H$ 
est un morphisme de faisceaux de groupes ; 
nous laissons au lecteur le soin d'énoncer et prouver la propriété
universelle du couple~$(\sch G/\sch H, \sch G\to \sch G/\sch H)$ sous ces hypothèses. Indiquons
simplement que la suite

$$1\to \sch H\to \sch G\to \sch G/\sch H\to 1$$ est exacte ; et qu'inversement,
si~$$1\to \sch H\to \sch G\to \sch G'\to 1$$ est une suite
exacte de faisceaux de groupes,
$\sch H$ est distingué dans~$\sch G$ et
le groupe~$\sch G'$ s'identifie
naturellement à~$\sch G/\sch H$.

\deux{lim-pref}
{\bf Limites dans la catégorie des faisceaux}.
Soit~$X$ un espace topologique et soit~$\sch D$ un diagramme
dans~$\fasc X$ ; pour tout ouvert~$U$
de~$X$, on note~$\sch D(U)$ le diagramme de~$\mathsf C$
obtenu en évaluant en~$U$
les constituants de~$\sch D$. 

\medskip
Les limites projectives et inductives
de~$\sch D$ existent ; nous allons brièvement indiquer leur construction,
en laissant le détail des vérifications au lecteur --elles sont
élémentaires. 

\trois{lim-proj-fasci}
{\em Construction de la limite projective.}
Elle est très simple : le préfaisceau~$U\mapsto \limproj \sch D(U)$ se trouve
être un faisceau, qui s'identifie à~$\limproj \sch D$. 

\trois{lim-ind-faisc}
{\em Construction de la limite inductive.}
Elle est un tout petit peu moins simple : le préfaisceau
$U\mapsto \limind \sch D(U)$ n'est pas un faisceau en général,
mais son faisceautisé s'identifie à~$\limind \sch D$ (combiner la propriété
universelle de la limite inductive préfaisceautique~$U\mapsto \limind \sch D(U)$ 
et celle du faisceautisé). 

\trois{lim-faisc-exemples}
{\em Exercice.}
Traduire en termes de limite
inductive (resp. projective) l'exactitude
de la suite de~\ref{ex-exafa-d} (resp. \ref{ex-exafa-g}), 
dans l'esprit de~\ref{ex-exa-d}
et~\ref{ex-exa-g} ; montrez que le faisceau quotient~$\sch G/\sch H$ de~\ref{quot-faisc}
peut s'interpréter comme la limite inductive
d'un diagramme convenable. 

\trois{limind-parfois-facile}
Il peut arriver dans certains cas qu'il ne
soit pas nécessaire de faisceautiser pour 
obtenir une limite inductive : par exemple, 
le lecteur vérifiera que la somme directe
{\em préfaisceautique}
d'une famille {\em finie}
de faisceaux de~$A$-modules est un faisceau. 
Nous l'invitons à montrer par un contre-exemple
que cette assertion est fausse en général sans hypothèse de
finitude. 

\trois{pourquoi-limproj-facile}
Il y a une explication conceptuelle au fait que
la limite projective préfaisceautique soit toujours un faisceau. Cette
assertion peut en effet se reformuler en disant que le foncteur d'inclusion~$I_X\colon \fasc X \hookrightarrow \pref X$
 commute aux limites projectives, et il a une excellente
raison pour ce faire : il admet un adjoint à gauche, à savoir la faisceautisation. 

\medskip
Le foncteur~$I_X$
n'admet pas d'adjoint à droite {\em en général}, faute de commuter aux limites inductives
({\em cf. supra}). 
Il peut toutefois en admettre dans certains cas très simples : par exemple, si~$X=\varnothing$
et si~$\mathsf C=\amod$, la catégorie~$\pref X$
s'identifie à celle des~$A$-modules, et la catégorie~$\fasc X$
à la catégorie~$\mathsf{Nul}$
dont le seul objet est le module nul ; et le foncteur~$M\mapsto \{0\}$ est à
la fois adjoint
{\em à gauche et à droite}
à l'inclusion~$\mathsf{Nul}\hookrightarrow \amod$. 

\medskip
Notons par contre que si~$\mathsf C=\ens$, le foncteur~$I_X$ n'admet
{\em jamais}
d'adjoint à droite : en effet, l'objet initial de~$\fasc X$ 
envoie~$U$ sur~$\varnothing$ si~$U\neq \varnothing$ et sur~$\{*\}$ sinon, 
alors que l'objet initial de~$\pref X$ envoie tout ouvert
{\em y compris~$\varnothing$}
sur~$\varnothing$. Le foncteur~$I_X$ ne commute donc pas aux limites inductives,
et n'admet dès lors pas d'adjoint à droite. 

\deux{fonct-faisc}
{\bf Fonctorialité.} Soit~$f: Y\to X$ une application continue entre espaces topologiques. 

\trois{f-lowerstar-ok}
Si~$\sch G$ est un faisceau sur~$Y$, on vérifie
immédiatement que 
le préfaisceau~$f_*\sch G$ est un faisceau. 

\trois{f-moinsu-pasok}
Par contre, si~$\sch F$ est un faisceau sur~$X$, le préfaisceau~$f^{-1}\sch F$ n'est pas
un faisceau en général. {\em C'est désormais son faisceautisé que l'on désignera par~$f^{-1}\sch F$}. 
Comme la faisceautisation ne modifie pas les fibres, on a encore pour tout~$y\in Y$ d'image~$x$
sur~$X$ un isomorphisme~$f^{-1}\sch F_y\simeq \sch F_x$. 

\trois{fonct-cont-faisc}
On a ainsi défini deux foncteurs : le foncteur~$f^{-1}$ 
de~$\fasc X$ vers~$\fasc Y$, et le foncteur~$f_*$
de~$\fasc Y$ vers~$\fasc X$.

\medskip
Si~$Z$ est un espace topologique et si~$g\colon Z\to Y$
est une application continue, on a des isomorphismes naturels de foncteurs

$$ (f\circ g)_*\simeq f_*\circ g_*\;\;\;\text{et}\;\;\;(f\circ g)^{-1}\simeq g^{-1}\circ f^{-1}.$$
C'est en effet évident pour les images directes ; et en ce qui concerne les images réciproques, 
les définitions fournissent un morphisme naturel~
de~$(f\circ g)^{-1}$ vers~$g^{-1}\circ f^{-1}$, dont on vérifie
sur les fibres que c'est un isomorphisme. 

\medskip
En combinant l'assertion préfaisceautique correspondante
et la propriété universelle du faisceautisé, on démontre que~$(f^{-1},f_*)$
est un couple de foncteurs adjoints. 

\trois{imm-ouverte-faisc}
Si~$U$ est un ouvert de~$X$ et
si~$j: U\hookrightarrow X$ est la flèche d'inclusion,
on vérifie que pour tout faisceau~$\sch F$ sur~$X$, le faisceau~$j^{-1}\sch F$
n'est autre que la restriction de~$\sch F$ à~$U$. 

\trois{im-sinlgeton-faisc}
Si~$x\in X$ et si~$i\colon \{x\}\hookrightarrow X$ est l'inclusion, 
on vérifie que pour tout faisceau~$\sch F$ sur~$X$, le faisceau~$i^{-1}\sch F$
envoie~$\{x\}$ sur~$\sch F_x$ (et~$\varnothing$ sur~$\{*\}$, nécessairement).

\section{Espaces annelés}
\markboth{Théorie des faisceaux}{Espaces annelés}

\subsection*{Définition, exemples, premières propriétés}

\deux{def-esp-ann}
{\bf Définition.}
Un {\em espace annelé}
est un couple~$(X,\sch O_X)$ où~$X$
est un espace topologique et~$\sch O_X$ un faisceau d'anneaux sur~$X$,
que l'on appelle parfois le faisceau
{\em structural}. 

\deux{exem-esp-ann}
{\bf Exemples.} 

\trois{esp-top-ann}
Soit~$X$ un espace topologique, et soit~$\sch O_X$ le faisceau 
des fonctions continues à valeurs réelles sur~$X$. Le couple~$(X,\sch O_X)$ est un
espace annelé (en~$\RR$-algèbres).

\trois{var-diff-ann}
Soit~$X$ une variété différentielle,
et soit~$\sch O_X$ le faisceau des fonctions~$\mathscr C^\infty$
sur~$X$. Le couple~$(X,\sch O_X)$ est un
espace annelé (en~$\RR$-algèbres).

\trois{var-hol-ann}
Soit~$X$ une variété analytique complexe,
et soit~$\sch O_X$ le faisceau des fonctions holomorphes sur~$X$. 
 Le couple~$(X,\sch O_X)$ est un
espace annelé (en~$\CC$-algèbres).

\trois{pref-const}
Soit~$X$ un espace topologique et soit~$\sch O_X$
le faisceau constant~$\underline{\ZZ}$ sur~$X$. Le couple~$(X,\sch O_X)$
est un espace annelé. 

\trois{rest-ouv-ann}
Soit~$(X,\sch O_X)$ un espace annelé, et soit~$U$
un ouvert de~$X$. Le couple~$(U,\sch O_X|_U)$ est un espace annelé ; 
sauf mention expresse du contraire, on considèrera toujours~$U$ comme
étant muni de cette structure d'espace annelé, et on écrira
à l'occasion~$\sch O_U$ au lieu de~$\sch O_X|_U$.

\deux{def-mor-espann}
{\bf Définition.} Soient~$(Y, \sch O_Y)$ et~$(X,\sch O_X)$ deux
espaces annelés. Un
{\em morphisme d'espaces annelés}
de~$(Y,\sch O_Y)$ vers~$(X,\sch O_X)$ est constitué d'une application continue~$\phi : Y\to X$ et d'une donnée
supplémentaire que l'on peut présenter de trois façons différentes, dont l'équivalence 
résulte des définitions et de l'adjonction entre~$\phi^{-1}$ et~$\phi_*$ : 

\medskip
a) un morphisme de faisceaux d'anneaux de~$\sch O_X$ vers~$\phi_*\sch O_Y$ ; 

b) un morphisme de faisceaux d'anneaux de~$\phi^{-1}\sch O_X$ vers~$\sch O_Y$ ; 

c) pour tout couple~$(U,V)$ formé d'un ouvert~$U$ de~$X$ et d'un ouvert~$V$ de~$Y$
tel que~$\phi(V)\subset U$, un morphisme d'anneaux~$\phi^*:\sch O_X(U)\to \sch O_Y(V)$, en exigeant que 
si~$U$ et~$U'$ sont deux ouverts de~$X$ avec~$U'\subset U$, 
et~$V$ et~$V'$ deux ouverts de~$Y$ avec~$V'\subset V$,~$\phi(V)\subset U$
et~$\phi(V')\subset U'$ alors le diagramme
$$\diagram \sch O_X(U)\rto^{\phi^*}\dto&\sch O_Y(V)\dto
\\ \sch O_X(U')\rto^{\phi^*} &\sch O_Y(V')\enddiagram$$
commute. 

\trois{comment-mor-espann}
On prendra
garde
que si~$(U,V)$ est comme au~c)
l'application~$\phi$ va de~$V$ vers~$U$, mais l'application~$\phi^*$
entre 
anneaux de sections va «dans l'autre sens», 
à savoir de~$\sch O_X(U)$ vers~$\sc O_Y(V)$. 

\trois{espann-categ}
En s'appuyant sur la définition~\ref{def-mor-espann}, on fait des espaces
annelés une catégorie -- la définition des identités et de la composition des morphismes
est laissée au lecteur. 

\deux{ex-mor-ann}
{\bf Exemples.}

\trois{esp-top-morann}
Soient~$Y$ et~$X$ deux espaces topologique, respectivement munis 
de leurs faisceaux de fonctions continues à valeurs réelles, et soit~$\phi$ 
une application continue de~$Y$ vers~$X$. Elle induit naturellement un 
morphisme d'espaces annelés entre~$Y$ et~$X$ : pour tout ouvert~$U$
de~$X$, tout ouvert~$V$ de~$Y$
tel que~$\phi(U)\subset V$ et toute fonction continue~$f: U\to \RR$, 
on pose~$\phi^*f=f\circ \phi$. 

\trois{var-diff-morann}
Soient~$Y$ et~$X$ deux
variétés différentielles, respectivement munies 
de leurs faisceaux de fonctions continues à valeurs réelles, et soit~$\phi$ 
une application~$\mathscr C^\infty$
de~$Y$ vers~$X$. Elle induit naturellement un 
morphisme d'espaces annelés entre~$Y$ et~$X$ : pour tout ouvert~$U$
de~$X$, tout ouvert~$V$ de~$Y$
tel que~$\phi(U)\subset V$ et toute fonction~$\mathscr C^\infty$~$f: U\to \RR$, 
on pose~$\phi^*f=f\circ \phi$. 

\trois{var-hol-morann}
Soient~$Y$ et~$X$ deux
variétés analytiques complexes, respectivement munies 
de leurs faisceaux de fonctions holmomorphes, et soit~$\phi$ 
une application
holomorphe de~$Y$ vers~$X$. Elle induit naturellement un 
morphisme d'espaces annelés entre~$Y$ et~$X$ : pour tout ouvert~$U$
de~$X$, tout ouvert~$V$ de~$Y$
tel que~$\phi(U)\subset V$ et toute fonction
holomorphe~$f: U\to \CC$, 
on pose~$\phi^*f=f\circ \phi$.

\trois{rest-ouv-morann}
Soit~$(X,\sch O_X)$ un espace annelé, et soit~$U$
un ouvert de~$X$. L'immersion~$j : U\hookrightarrow X$ est
sous-jacente à un morphisme naturel d'espace annelés : si~$U'$ est un ouvert
de~$U$ et~$X'$ un ouvert de~$X$ contenant~$U'$, le morphisme
$j^*: \sch O_X(X')\to \sch O_U(U')=\sch O_X(U')$ est simplement la restriction. 

\medskip
Soit~$\phi : (Y,\sch O_Y)\to (X,\sch O_X)$ un morphisme d'espaces annelés
tel que~$\phi(Y)$
soit contenu dans~$U$. L'application continue~$Y\to U$ induite
par~$\phi$ est sous-jacente à un morphisme d'espaces annelés
de~$(Y,\sch O_Y)$ vers~$(U,\sch O_U)$ : si~$W$ est un ouvert de~$Y$
et si~$V$ est un ouvert de~$U$ contenant~$\phi(Y)$, 
le morphisme
d'anneaux correspondant~$\sch O_U(W)=\sch O_X(W)\to \sch O_Y(V)$
n'est autre que~$\phi^*$. 

En d'autres termes, toute factorisation
{\em ensembliste}
par~$U$ est automatiquement {\em morphique}. 

\medskip
On vérifie aisément que ce morphisme
d'espace annelés~$(Y,\sch O_Y)\to (U,\sch O_U)$ est le seul
dont la composée avec~$j$ soit égale à~$\phi$. Cela signifie que
$((U,\sch O_U),j)$ représente le foncteur qui envoie un espace annelé~$(Y,\sch O_Y)$ 
sur
$$\{\phi\in \hom_{\mathsf{Esp-ann}}((Y,\sch O_Y), (X,\sch O_X)),\;\;\;\phi(Y)\subset U\}.$$

\trois{imm-ouverte-espann}
Nous dirons qu'un morphisme d'espaces annelés~$\phi \colon Y\to X$
est une {\em immersion ouverte}
s'il induit un isomorphisme entre~$Y$ et un ouvert de~$X$. 

\trois{morann-fibres}
Soit~$\phi : (Y,\sch O_Y)\to (X,\sch O_X)$
un morphisme d'espaces annelés ; soit~$y$
un point de~$Y$
et soit~$x$ son image sur~$X$. Le
morphisme~$\phi$ induit alors de manière naturelle un morphisme d'anneaux
de~$\sch O_{X,x}$ vers~$\sch O_{Y,y}$, souvent encore noté~$\phi^*$. 

\deux{lim-ind-espann}
{\bf Limites inductives d'espaces annelés.}
Soit~$\sch D=((X_i), (E_{ij}))$ un diagramme dans la catégorie des espaces
annelés. Il admet une limite inductive~$\limind \sch D$
que nous allons brièvement décrire -- la justification
est laissée en exercice.

\medskip
En tant qu'espace topologique,
$\limind \sch D$ coïncide avec la limite
inductive de~$\sch D$ dans~$\top$, construite 
au~\ref{limind-top}. Il reste à décrire le faisceau structural
de~$\limind \sch D$. Soit~$U$ un ouvert de~$\limind \sch D$ ; 
son image réciproque de~$U$ sur~$\coprod X_i$ est de la forme~$\coprod U_i$,
où~$U_i$ est pour tout~$i$ un ouvert de~$X_i$, et où~$f(U_i)\subset U_j$ pour tout~$(i,j)$ et toute~$f\in E_{ij}$. 
On définit alors~$\sch O_{\limind \sch D}(U)$
comme l'ensemble des familles~$(s_i)\in \prod \sch O_{X_i}(U_i)$
telles que~$f^*s_j=s_i $ pour tout~$(i,j)$ et toute~$f\in E_{i,j}.$
On note~$\lambda_i \colon X_i\to \limind \sch D$ la flèche structurale.

\subsection*{Les~$\sch O_X$-modules}

\deux{ox-mod}
Soit~$(X,\sch O_X)$ un espace annelé. Un {\em $\sch O_X$-module}
$\sch M$  est un faisceau en groupes abéliens~$\sch M$ sur~$X$ muni,
pour tout ouvert~$U$ de~$X$, d'une loi externe~$\sch O_X(U)\times \sch M(U)\to
\sch M(U)$ qui fait du groupe abélien~$\sch M(U)$ un~$\sch O_X(U)$-module, ces données
étant sujettes à la condition suivante : pour tout ouvert~$U$ de~$X$, tout ouvert~$V$
de~$U$, toute section~$s\in \sch M(U)$ et toute~$f\in \sch O_X(U)$,
on a
$$(fs)|_V=(f|_V)(s|_V).$$

Un morphisme de~$\sch O_X$-modules est un morphisme
de faisceaux en groupes qui est~$\sch O_X$-linéaire en un sens évident. 

\deux{ex-ox-mod}
Donnons 
un premier exemple
venu de la géométrie différentielle : 
si~$(X,\sch O_X)$ est une variété différentielle munie de son faisceau
des fonctions~$\sch C^\infty$, alors le faisceau des champs de vecteurs,
qui associe à un ouvert~$U$ de~$X$ les sections du fibré tangent de~$X$
au-dessus de~$U$ (ou, de façon
équivalente, les dérivations de la~$\RR$-algèbre~$\sch O_X(U)$)
est de manière naturelle
un~$\sch O_X$-module.

\deux{fibres-mod}
Soit~$(X,\sch O_X)$ un espace annelé et soit~$\sch M$
un~$\sch O_X$-module. Pour tout~$x$
appartenant à~$X$, la fibre
$\sch M_x$ hérite d'une structure naturelle de~$\sch O_{X,x}$-module. 

\deux{op-usuelles-alg}
{\bf Faisceautisation d'opérations usuelles sur les modules}. 
Les notions usuelles en théorie des modules se faisceautisent souvent,
donnant ainsi lieu à des notions analogues en théorie des~$\sch O_X$-modules. 
Donnons quelques exemples ; on fixe un espace annelé~$(X,\sch O_X)$. 

\trois{somm-ox}
{\em Les limites inductives et projectives}.
Soit~$\sch D$ un diagramme dans~$\oxmod$ ; pour
tout~$U$, notons~$\sch D(U)$ le diagramme des~$\sch O_X(U)\text{-}\mathsf{Mod}$
obtenu par évaluation en~$U$ des constituants de~$\sch D$. 

\medskip
Il admet
une limite projective et une limite inductive, construites comme
à la section précédente : le préfaisceau~$U\mapsto \limproj \sch D(U)$
est naturellement un faisceau, admettant une structure évidente de~$\sch O_X$-module, 
et il s'identifie à~$\limproj \sch D$ ; le préfaisceau~$U\mapsto \limind\sch D(U)$
n'est quant à lui en général pas un faisceau, mais quand on le faisceautise on obtient
un~$\sch O_X$-module qui s'identifie à~$\limind \sch D$. 

\medskip
Il peut arriver que pour certaines limites inductives l'opération 
de faisceautisation ne soit pas nécessaire : on vérifie par exemple
que la somme directe préfaisceautique d'une famille finie de~$\sch O_X$-modules
est déjà
un faisceau. 

\trois{tens-ox-mod}
Soient~$\sch M$ et~$\sch N$
deux~$\sch O_X$-modules. Le préfaisceau
$$U\mapsto \sch M(U)\otimes_{\sch O_X(U)}\sch N(U)$$ n'est pas
un faisceau en général ; son faisceautisé est noté~$\sch M\otimes_{\sch O_X}\sch N$,
et hérite d'une structure naturelle de~$\sch O_X$-module. 

\medskip
Le lecteur établira sans la moindre difficulté les faits suivants (en raisonnant ouvert 
par ouvert, et en appliquant la propriété universelle du faisceautisé). On dispose d'un
morphisme bi-$\sch O_X$-linéaire canonique
de~$\sch M\times \sch N$ dans~$\sch M\otimes_{\sch O_X}\sch N$, notée~$(f,g)\mapsto f\otimes g$ ; pour tout~$\sch O_X$-module~$\sch P$
et tout morphisme bi-$\sch O_X$-linéaire $b$~$\sch M\times \sch N$ vers~$\sch P$, il existe une unique
application~$\sch O_X$-linéaire $\ell\colon \sch M\otimes_{\sch O_X}\sch N$ telle que~$\ell(f\otimes g)=b(f,g)$ pour tout $(f,g)$.

\trois{prod-caslibre}
Soit~$n$ un entier.
Le préfaisceau~$U\mapsto \sch M(U)\otimes_{\sch O_X(U)}\sch O_X(U)^n$
est un faisceau, qui s'identifie à~$\sch M^n$ ; 
le produit tensoriel~$\sch M\otimes_{\sch O_X}\sch O_X^n$ 
s'obtient donc directement, sans faisceautisation, et est naturellement
isomorphe à~$\sch M^n$. 

\trois{fibre-prodtens-oxmod}
On a pour tout~$x\in X$ un isomorphisme canonique
$$(\sch M\otimes_{\sch O_X}\sch N)_x\simeq 
\sch M_x\otimes_{\sch O_{X,x}}\sch N_x.$$ 

Pour s'en convaincre, on peut ou bien procéder «à la main»
à l'aide
des constructions explicites de la limite inductive filtrante et du produit tensoriel, 
ou bien remarquer que les deux~$\sch O_{X,x}$-modules en jeu représentent
le même foncteur covariant, à savoir celui qui envoie un~$\sch O_{X,x}$-module~$P$
sur
$$\limproj {\rm Bil}_{\sch O_X(U)}(\sch M(U)\times \sch N(U),P)$$ où~$U$
parcourt l'ensemble des voisinages ouverts de~$x$. 

\medskip
Le~$\sch O_X$-module $ \sch M \otimes_{\sch O_X}\sch N$ peut
être caractérisé par une propriété universelle ; nous
vous suggérons à titre d'exercice de l'énoncer et de la prouver. 

\trois{ox-alg}
Une {\em $\sch O_X$-algèbre}
$\sch A$ est un faisceau d'anneaux sur~$X$ muni d'un morphisme
de faisceaux d'anneaux~$\sch O_X\to \sch A$. Une~$\sch O_X$-algèbre
hérite d'une structure naturelle de~$\sch O_X$-module. 

\medskip
Si~$\sch A$ et~$\sch B$ sont deux 
$\sch O_X$-algèbres, un morphisme de~$\sch O_X$-algèbres
de~$\sch A$ vers~$\sch B$ est un morphisme de faisceau d'anneaux
de~$\sch A\to \sch B$
faisant commuter le diagramme
$$\diagram \sch A\rto&\sch B\\
\sch O_X\uto\urto&\enddiagram.$$

\medskip
Si~$\sch A$ est une~$\sch O_X$-algèbre et si~$\sch N$ est un~$\sch A$-module,
$\sch N$ hérite {\em via}
la flèche~$\sch O_X\to \sch A$ d'une structure naturelle de~$\sch O_X$-module.
On dit que cette structure
est obtenue par {\em restriction
des scalaires}
à de~$\sch A$ à~$\sch O_X$. 
 
\medskip
Si~$\sch A$ est une~$\sch O_X$-algèbre et si~$\sch M$ est un~$\sch O_X$-module,
le produit tensoriel~$\sch A\otimes_{\sch O_X}\sch M$ hérite d'une structure naturelle
de~$\sch A$-module ; on dit que le~$\sch A$-module~$\sch A\otimes_{\sch O_X}\sch M$
est déduit de~$\sch M$
par {\em 
extension des scalaires}
de~$\sch A$
à~$\sch O_X$.  Nous vous laissons
en exercice la démonstration du fait suivant : l'extension des scalaires
de~$\sch O_X$ à~$\sch A$ est l'adjoint à gauche de la restriction
des scalaires de~$\sch A$ à~$\sch O_X$. 

\medskip
Si~$\sch A$ et~$\sch B$ sont deux~$\sch O_X$-algèbres, 
leur produit tensoriel~$\sch A\otimes_{\sch O_X}\sch B$ est de manière 
naturelle une~$\sch O_X$-algèbre ; nous vous invitons à prouver
que~$\sch A\otimes_{\sch O_X}\sch B$
est la somme disjointe de~$\sch A$
et~$\sch B$ dans la catégorie des~$\sch O_X$-algèbres. 

\deux{fonct-oxmod}
{\bf Fonctorialité}.
Soit~$\phi : (Y,\sch O_Y)\to (X,\sch O_X)$ un
morphisme d'espaces annelés. Soit~$\sch N$ un~$\sch O_Y$-module
et soit~$\sch M$ un~$\sch O_X$-module. 

\trois{f-lowerstar-oxmod}
La structure de~$\sch O_Y$-module sur~$\sch N$ induit de manière naturelle
une structure de~$\phi_*\sch O_Y$-module sur~$\phi_*\sch N$. Le morphisme
$\phi$ est par définition fourni avec un morphisme de faisceaux d'anneaux
$\sch O_X\to \phi_*\sch O_Y$, par le biais duquel
le~$\phi_*\sch O_Y$-module~$\phi_*\sch N$
peut être vu comme un~$\sch O_X$-module. 

\trois{f-upperstar-oxmod}
Les choses se passent un peu mois bien concernant le foncteur~$\phi^{-1}$ : 
le faisceau~$\phi^{-1}\sch M$ hérite d'une structure naturelle
de~$\phi^{-1}(\sch O_X)$-module, mais la flèche~$\phi^{-1}(\sch O_X)\to \sch O_Y$
va dans le mauvais sens et
ne permet pas de faire de~$\phi^{-1}\sch M$
un~$\sch O_Y$-module. Elle permet par contre
de le
{\em transformer}
de manière universelle en un~$\sch O_Y$-module, par tensorisation ; 
le~$\sch O_Y$-module~$\sch O_Y\otimes_{\phi^{-1}\sch O_X}\sch M$
ainsi obtenu est noté~$\phi^*\sch M$. 

\medskip
Si~$y\in Y$
et si~$x$ désigne son image
sur~$X$, il existe un isomorphisme
naturel~$(\phi^*\sch M)_y\simeq \sch O_{Y,y}\otimes_{\sch O_{X,x}}\sch M_x$ :
on le voit en combinant le bon comportement vis-à-vis des fibres
de~$\phi^{-1}$ (\ref{f-moinsu-pasok})
et du produit tensoriel (\ref{fibre-prodtens-oxmod}).

\trois{fonct-oxmod}
On vérifie que~$\phi_*$ est de façon naturelle un foncteur
de la catégorie des~$\sch O_Y\text{-}\mathsf{Mod}$
vers~$\oxmod$, et que~$\phi^*$ est de façon naturelle un foncteur
de la catégorie de~$\oxmod$ vers~$\sch O_Y\text{-}\mathsf{Mod}$. 

\trois{compose-oxmod}
Si~$Z$ est un espace annelé et~$\psi \colon Z\to Y$ un morphisme, 
on dispose d'isomorphismes naturels de foncteurs

$$(\phi\circ \psi)_*\simeq \phi_*\circ \psi_*\;\;\text{et}\;\;(\phi\circ \psi)^*\simeq \psi^*\circ \phi^* \;.$$
C'est en effet évident pour les images directes ; et en ce qui concerne les images réciproques, 
les définitions fournissent un morphisme 
de~$(f\circ g)^*$ vers~$g^*\circ f^*$, dont on vérifie
sur les fibres que c'est un isomorphisme.

\trois{pullback-prod-tens}
Soit~$\sch M$ et~$\sch N$ deux~$\sch O_X$-modules ; on a un isomorphisme
naturel
$$\phi^*\sch M\otimes_{\sch O_Y}\phi^*\sch N\simeq \phi^*(\sch M\otimes_{\sch O_X}\sch N).$$
En effet, les définitions fournissent un morphisme naturel du terme de gauche vers celui de droite,
et l'on vérifie sur les fibres que c'est un isomorphisme. 

\trois{adj-ocmod}
Le couple~$(\phi^*,\phi_*)$ est un couple
de foncteurs adjoints : c'est une conséquence
formelle
des propriétés d'adjonction
du couple~$(\phi^{-1},\phi_*)$, et du couple
formé de l'extension des scalaires et de la restriction
des scalaires
(faisceautiques). 

\trois{ex-phietoile}
Il résulte immédiatement de la définition que~$\phi^*\sch O_X=\sch O_Y$. 
Par ailleurs, comme~$\phi$ a un adjoint à droite, il commute aux limites inductives
et en particulier aux sommes directes quelconques (le lecteur pourra le vérifier directement). 
En particulier, $\phi^*(\sch O_X^n)=\sch O_Y^n$ pour tout entier~$n$.

\deux{intro-loclibre}
{\bf Le faisceau des homomorphismes}. 
Soit~$(X,\sch O_X)$ un espace annelé. Si~$\sch F$ et~$\sch G$
sont deux~$\sch O_X$-modules, on
vérifie immédiatement
que
$$U\mapsto \hom(\sch F|_U,\sch G|_U)$$ est un faisceau, qui possède
lui-même une structure
naturelle de~$\sch O_X$-module ; on le note~$\homfsc (\sch F,\sch G)$. 
Le faisceau~$\homfsc(\sch F,\sch F)$ sera également noté~$\fscend \sch F$ ; il a
une structure naturelle ({\em via}
la composition des endomorphismes ouvert par ouvert)
de~$\sch O_X$-algèbre {\em non commutative en général}. 

\trois{homint-fonct}
La flèche~$(\sch F,\sch G)\mapsto \homfsc (\sch F,\sch G)$ est de manière 
naturelle un foncteur covariant en~$\sch G$ et
un foncteur contravariant en~$\sch F$.

\trois{exemple-homint}
{\em Exemple}. 
Soit~$\sch F$ un~$\sch O_X$-module et soit~$n\in \NN$. Soit~$U$ un ouvert de~$X$
 et soient~$e_1,\ldots, e_n$ des sections
 de~$\sch F$ sur~$U$. La formule
 $$(a_1,\ldots, a_n)\mapsto \sum a_i e_i$$ définit un morphisme
 de~$\sch O_U$-modules de~$\sch O_U^n$ vers~$\sch F|_U$, que l'on dira {\em induit
 par les~$e_i$}. 
 Réciproquement, tout morphisme~$\phi$
 de~$\sch O_U$-modules
 de~$\sch O_U^n$ vers~$\sch F|_U$ est de cette forme : prendre
 $e_i=\phi(\underbrace{0,\ldots, 1,\ldots, 0}_{\text{le}\;1\;\text{est \`a la place}\;i})$.

\medskip
En faisant varier~$U$
dans cette construction, on obtient
un isomorphisme canonique~$\homfsc(\sch O_X^n,\sch F)\simeq \sch F^n$. 
En particulier, $\homfsc(\sch O_X^n,\sch O_X^m)\simeq \sch O_X^{nm}$ pour tout~$m$. 

\trois{homoth-faisc}
Soit~$\sch F$ un~$\sch O_X$-module. 
À toute section~$f$ de~$\sch O_X$ sur un ouvert~$U$ de~$X$
est associée de manière naturelle
un endomorphisme de~$\sch F|_U$, à savoir {\em l'homothétie 
de rapport~$f$}, donnée par la formule~$s\mapsto fs$. 

\medskip
En faisant varier~$U$ dans cette construction, on obtient
un morphisme canonique de~$\sch O_X$-algèbres
$\sch F\to \fscend \sch F$. Lorsque~$\sch F=\sch O_X$,
ce morphisme est un isomorphisme en vertu de~\ref{exemple-homint}. 

\deux{dual-interne}
Soit~$\sch F$ un~$\sch O_X$-module. On note~$\sch F^\vee$ 
le~$\sch O_X$-module $\homfsc (\sch F,\sch O_X)$, que l'on appelle
aussi le {\em dual}
de~$\sch F$. 

\trois{bidual-fasci}
Soit~$\sch F$ un~$\sch O_X$-module. Soit~$U$ un ouvert de~$X$
et soit~$s\in \sch F(U)$. La section~$s$ définit de manière naturelle
un morphisme de~$\sch F^\vee|_U$ vers~$\sch O_U$, 
donné par la formule~$\phi\mapsto \phi(s)$. 

\medskip
En faisant varier~$U$
dans cette construction,  on obtient un morphisme canonique
de~$\sch F$ dans son bidual~$\sch F^{\vee \vee}$. 

\medskip
On vérifie immédiatement à l'aide de~\ref{exemple-homint}
que~$\sch O_X^n\to (\sch O_X^n)^{\vee \vee}$
est un isomorphisme pour tout~$n$. 

\trois{endos-dual-fasc}
Soient~$\sch F$ et~$\sch G$ deux~$\sch O_X$-modules. Soit~$U$ un ouvert de~$X$, soit~$\phi$
appartenant à~$\sch F^\vee(U)$ 
et soit~$s\in \sch G(U)$. Le couple~$(\phi,s)$ définit un morphisme de~$\sch F|_U$ vers~$\sch G|_U$, donné
par la formule~$t\mapsto \phi(t)s$.  

\medskip
En faisant varier~$U$
dans cette construction, on obtient un morphisme bi-$\sch O_X$-linéaire de~$\sch F^\vee\times \sch G$
vers~$\homfsc (\sch F,\sch G)$, qui induit un morphisme de~$\sch O_X$-modules
de~$\sch F^\vee \otimes_{\sch O_X}\sch G$ vers~$\homfsc(\sch F,\sch G)$. 

\trois{dual-endos-libre}
Soient~$n$ et~$m$ deux entiers. 
En vertu de~\ref{exemple-homint}, la flèche naturelle
$$(\sch O_X^n)^\vee(U) \otimes_{\sch O_X(U)}\sch O_X^m(U)
\to \hom(\sch O_U^n,\sch O_U^m)$$ est un isomorphisme pour tout ouvert~$U$
de~$X$.  En conséquence, 
$$(\sch O_X^n)^\vee\otimes_{\sch O_X}\sch O_X^m\to \homfsc(\sch O_X^n,\sch O_X^m)$$
est un isomorphisme. 

\deux{exo-homint}
{\bf Exercice}. Soit~$\sch F$ un~$\sch O_X$-module. Montrez que
$\sch H\mapsto \homfsc(\sch F,\sch H)$ est adjoint à droite à
$\sch G\mapsto \sch G\otimes_{\sch O_X}\sch F$.

%
%
%
%
%
%

\section{Espaces localement annelés}
\markboth{Théorie des faisceaux}{Espaces localement annelés}

\subsection*{Définition, exemples, premières propriétés}
\deux{def-locann}
{\bf Définition.}
On appelle
{\em espace localement annelé}
un espace annelé~$(X,\sch O_X)$
tel que~$\sch O_{X,x}$ soit pour tout~$x\in X$
un anneau local. 

\deux{ex-locann}
{\bf Exemples.}

\trois{esp-top-locann}
{\bf Le cas des faisceaux de fonctions.} 
Soit~$k$ un corps, soit~$X$ un espace topologique, et soit~$\sch O_X$
un sous-faisceau de~$k$-algèbres du faisceau de toutes les fonctions de~$X$ vers~$k$
(en particulier, $\sch O_X$ contient les fonctions constantes).
Supposons que~$\sch O_X$ possède la propriété suivante : {\em pour tout ouvert~$U$ de~$X$
et tout~$x\in U$, une fonction~$f\in \sch O_X(U)$ telle que~$f(x)\neq 0$ admet un inverse dans~$\sch O_X(V)$
pour un certain voisinage ouvert~$V$ de~$x$ dans~$U$.}

Sous ces hypothèses, $(X,\sch O_X)$ est localement annelé ; pour tout~$x\in X$, l'idéal
maximal de~$\sch O_{X,x}$ est le  noyau de la surjection~$x\mapsto f(x)$
de~$\sch O_{X,x}$ sur~$k$.

\medskip
La preuve est {\em
mutatis mutandis} celle donnée
au~\ref{ann-local-geodiff}
lorsque~$k=\RR$ et lorsque~$\sch O_X$ est le faisceau
des fonctions continues de~$X$ dans~$\RR$, mais l'assertion plus
générale que nous présentons ici s'applique dans bien d'autres cas intéressants : 

\medskip
$\bullet$ le corps $k$ est égal à~$\RR$, l'espace~$X$
est une variété différentielle et $\sch O_X$ est
le faisceau des fonctions~$\mathscr C^\infty$
sur~$X$ ; 

$\bullet$ le corps $k$ est égal à~$\CC$, l'espace~$X$
est une variété analytique complexe et $\sch O_X$ est
le faisceau des fonctions
holomorphes
sur~$X$ ;

$\bullet$ le corps~$k$ est algébriquement clos, $X$
est une variété algébrique sur~$k$ au sens des articles FAC et GAGA de
Serre (qui est aussi celui adopté par Perrin dans son cours de géométrie algébrique),
et~$\sch O_X$ est le faisceau des fonctions régulières sur~$X$. 

\trois{rest-ouv-locann}
{\bf Stabilité par restriction à un ouvert.}
Soit~$(X,\sch O_X)$ un espace localement
annelé, et soit~$U$
un ouvert de~$X$. L'espace annelé~$(U,\sch O_U=\sch O_X|_U)$ est 
localement annelé : cela provient du fait que  l'on a
pour tout~$x\in U$ l'égalité~$\sch O_{U,x}=\sch O_{X,x}$. 

\deux{comment-locann}
Soit
~$(X,\sch O_X)$ un espace localement annelé.

\trois{not-eval}
Soit~$x\in X$. Notons~$\kappa(x)$ le corps résiduel de l'anneau local~$\sch O_{X,x}$. 
On dit aussi que c'est le corps résiduel {\em du point~$x$}, 
et l'on dispose d'une surjection canonique~$\sch O_{X,x}\to \kappa(x)$. 
{\em Pour des raisons psychologiques}, on choisit de noter cette
surjection~$f\mapsto f(x)$, et d'appeler ce morphisme «évaluation en~$x$». 
On a ainsi l'équivalence

$$f(x)\neq 0\iff \;f\;\text{est inversible dans}\;\sch O_{X,x}.$$

\trois{eval-u}
Soit~$U$ un ouvert de~$X$ et soit~$x\in U$. La
flèche composée
$$\sch O_X(U)\to \sch O_{X,x}\to \kappa(x)$$ sera encore notée~$f\mapsto f(x)$.
Remarquons que si~$f$ est un élément inversible de~$\sch O_X(U)$, son image~$f(x)$ par le morphisme
d'évaluation est un élément inversible du corps~$\kappa(x)$, et est en particulier non nulle. 

\trois{locann-non-vide}
Si~$X=\varnothing$ on a~$\sch O_X(X)=\{0\}$ puisque~$\sch O_X$ est un faisceau. Supposons 
maintenant que~$X$ soit non vide, et soit~$x\in X$. L'évaluation
en~$x$ fournit un morphisme de~$\sch O_X(X)$ vers le corps~$\kappa(x)$, ce qui force~$\sch O_X(X)$
à être {\em non nul} : si~$1$ était nul dans~$\sch O_X(X)$ on aurait $1(x)=0$, c'est-à-dire
$1=0$ dans~$\kappa(x)$, ce qui est absurde.

\trois{lemme-inv-locann}
{\bf Lemme.}
{\em Soit~$U$ un ouvert de~$X$ et soit~$f\in \sch O_X(U)$. L'ensemble~$D(f)$ des points~$x\in U$ en lesquels~$f$ ne s'annule pas est un ouvert de~$U$, 
et~$f$ est inversible dans~$\sch O_X(U)$ si et seulement si~$D(f)=U$.}

\medskip
{\em Démonstration.}
Soit~$x\in D(f)$. Comme~$f(x)\neq 0$, on déduit de
l'équivalence mentionnée en~\ref{not-eval}
que~$f$ est inversible dans~$\sch O_{X,x}$, c'est-à-dire sur un voisinage ouvert~$V$de~$x$.
On a alors~$f(y)\neq 0$ pour tout~$y\in V$
d'après~\ref{eval-u} ; en conséquence,~$V\subset D(f)$
et~$D(f)$ est ouvert. 

\medskip
Si~$f$ est inversible dans~$\sch O_X(U)$ alors~$f(x)\neq 0$ pour tout~$x\in U$
d'après~\ref{eval-u}, et~$U=D(f)$. 

\medskip
Réciproquement, supposons que~$U=D(f)$ et soit~$x\in U$. 
Comme~$f(x)\neq 0$, on déduit de
l'équivalence mentionnée en~\ref{not-eval}
que~$f$ est inversible dans~$\sch O_{X,x}$, c'est-à-dire au voisinage de~$x$. Ceci valant
pour tout~$x\in U$, il existe un recouvrement ouvert~$(U_i)$ de~$U$ et pour tout~$i$
un inverse~$g_i$ de~$f|_{U_i}$ dans~$\sch O_X(U_i)$. 

Pour tout couple~$(i,j)$, chacune des
deux restrictions de~$g_i$ et~$g_j$ à~$U_i\cap U_j$ 
est un inverse de~$f|_{U_i\cap U_j}$ ; par unicité de l'inverse, elles coïncident. 
Les section~$g_i$ du faisceau~$\sch O_X$ se recollent donc en une section~$g$ 
de~$\sch O_X$ sur~$U$, qui satisfait les égalités~$gf|_{U_i}=1$ pour tout~$i$ ; en conséquence,
$gf=1$ et~$f$ est inversible
dans~$\sch O_X(U).$$\Box$ 

\deux{comment-espann}
{\bf Commentaires.}
On voit que
le faisceau structural
d'un espace 
localement
annelé {\em quelconque}
ressemble par certains aspects aux faisceaux
de fonctions à valeurs dans un corps
tels que décrits en~\ref{esp-top-locann} : ses sections peuvent être évaluées en tout point
(le résultat vivant dans un corps), le lieu des points en lesquels une section ne s'annule pas
est un ouvert, et une section est inversible si et seulement si elle ne
s'annule pas. 

Pour cette raison, on pense assez souvent aux sections 
du faisceau structural comme à des fonctions,
et il arrive fréquemment d'ailleurs qu'on les qualifie
(un peu abusivement) ainsi. 

\medskip
Nous attirons toute fois l'attention sur deux points importants
qui montrent les limites de l'intuition «fonctionnelle»
appliquée aux espaces localement annelés généraux. 

\trois{corps-varie}
{\em Premier point.}
Dans un espace localement annelé~$(X,\sch O_X)$,
le corps~$\kappa(x)$ {\em dépend {\em a priori}
de~$x$.} Dans 
la situation 
considérée au~\ref{esp-top-locann}
il était constant
mais en général, il peut effectivement varier. 

\medskip
Nous n'avons pas pour le moment d'exemple naturel (c'est-à-dire non construit
exprès) d'espace localement annelé sur lequel cela se produit. Indiquons simplement que
cela arrive fréquemment sur un schéma, et donnons en attendant un exemple
«artificiel»
très simple : on prend pour~$X$ un ensemble à deux éléments~$x$ et~$y$, muni
de la topologie discrète et du faisceau

$$\varnothing \mapsto \{0\},\;\;\{x\}\mapsto \CC,\;\;\{y\}\mapsto \ZZ/2\ZZ,\;\,\;{\rm et}\; \{x,y\}\mapsto \CC\times \ZZ/2\ZZ.$$
On vérifie immédiatement que~$X$ est un espace localement annelé, que~$\kappa(x)=\CC$
et que~$\kappa(y)=\ZZ/2\ZZ$.

\trois{nilp-locann}
{\em Second point.}
Si~$(X,\sch O_X)$ 
est un espace localement annelé, si~$U$ est un ouvert de~$X$ et si~$f\in \sch O_X(U)$, 
{\em il se peut que $f(x)=0$ pour tout~$x\in U$ sans que~$f$ soit nulle.}

\medskip
Par exemple, supposons
que~$f$ soit nilpotente et non nulle. Dans ce
cas,~$f(x)$ est pour tout~$x\in X$ un élément nilpotent
d'un corps, et est donc trivial. 

\medskip
Nous allons donner un exemple très simple où une telle~$f$
existe. Fixons un corps~$k$, et considérons 
un espace topologique singleton~$\{x\}$. Se donne une structure d'espace
localement annelé sur~$\{x\}$ revient à choisir un anneau local. Soit~$\sch O$ 
le quotient~$k[T]/(T^2)$ et soit~$f$ la classe de~$T$ ; elle est non nulle.

L'ensemble
des idéaux premiers de~$\sch O$ est en bijection avec l'ensemble 
des idéaux premiers de~$k[T]$ contenant~$(T^2)$ ; il n'y en a qu'un,
à savoir~$(T)$. L'anneau~$\sch O$ est donc local, 
son unique idéal maximal est~$(f)$,
et son corps résiduel est~$k[T]/(T)=k$. 

On a ainsi bien 
défini une structure d'espace localement annelé sur~$x$. Comme~$f$
est nilpotente, on a~$f(x)=0$. On pouvait d'ailleurs
le voir directement ici, puisque l'évaluation
en~$x$ est la réduction modulo l'idéal maximal de~$\sch O$, 
c'est-à-dire justement modulo~$(f)$. D'une manière générale,
si~$g$ est un élément de~$\sch O$, il s'écrit~$a+bf$, où~$a$
et~$b$ sont deux éléments uniquement déterminés de~$k$, et
on a alors~$g(x)=a$. 

\trois{pouquoi-nilp}
{\bf Commentaires.} On peut se demander pourquoi autoriser
ce genre d'horreurs, alors qu'on a fait en sorte, pour ce qui 
concerne la non-annulation, que les propriétés usuelles soient satisfaites. 
La raison est que la présence de «fonctions» 
nilpotentes non nulles
peut avoir un sens géométrique profond,
et c'est notamment le cas dans l'exemple
que l'on vient de traiter.  

\medskip
En effet, considérons, dans le plan affine sur~$k$
en coordonnées~$S$ et~$T$, la parabole~$P$ d'équation~$T^2=S$
et la droite~$D$ d'équation~$S=0$. Leur intersection naïve
est le point~$x$
de coordonnées~$(0,0)$. En théorie des schémas, cette intersection est un peu 
plus riche que~$\{x\}$ : on garde en mémoire le corps de base et les équations, 
et l'intersection sera donc l'espace topologique~$\{x\}$ muni
du faisceau (ou de l'anneau, si l'on préfère)~$k[S,T]/(S,T^2-S)\simeq k[T]/T^2$ :
on retrouve l'espace localement annelé évoqué plus haut. 

La présence de nilpotents
non triviaux parmi les fonctions sur~$P\cap D$ s'interprète intuitivement
comme suit : l'intersection~$P\cap D$
est égale au point~$x$  {\em infinitésimalement
épaissi} parce que~$P$ et~$D$ sont tangentes en~$x$ ;  
le point d'intersection~$x$ 
est en quelque sorte 
double, et c'est cette multiplicité qui
est codée algébriquement par l'existence de nilpotents non triviaux. 

\medskip
Cet exemple est significatif : c'est 
pour prendre en compte les multiplicités dans la théorie
que Grothendieck a décidé d'admettre les «fonctions» nilpotentes non nulles . 
Cela
se révèle un outil extraordinairement souple, mais il y a un prix à payer : il faut
autoriser une
«fonction»
à s'annuler en tout point
sans être globalement nulle. D'où le choix du formalisme abstrait des espaces localement annelés, 
qui mime 
en partie le point de vue fonctionnel classique, 
mais permet ce genre de fantaisies finalement très utiles.

\deux{digres-alg}
{\bf Digression algébrique}. Soient~$A$ et~$B$ deux anneaux locaux
d'idéaux maximaux respectifs~$\got m$ et~$\got n$. Soit~$f$ un morphisme
de~$A$ vers~$B$. Si~$a\in A$ et si~$f(a)\in\got n$ alors~$f(a)$ n'est pas
inversible, et~$a$ n'est donc pas inversible non plus ; autrement dit,~$a\in \got m$. 
On dit que~$f$ est {\em local} 
si la réciproque est vraie, c'est-à-dire si~$f(\got m)\subset \got n$, ou encore
si~$f^{-1}(\got n)=\got m$.

\deux{def-mor-locann}
{\bf Définition.}
Soient~$(Y,\sch O_Y)$ et~$(X,\sch O_X)$ deux espaces localement annelés. Un {\em morphisme
d'espaces localement annelés} de~$(Y,\sch O_Y)$ vers~$(X,\sch O_X)$ est un morphisme~$\phi$
d'espaces  annelés tel que~$\sch O_{X,\phi(y)}\to \sch O_{Y,y}$ soit pour tout~$y\in Y$
un morphisme local.

\trois{recriture-morlocann}
Nous allons récrire cette condition de façon plus suggestive. Soit~$y$
un point de~$Y$. 
Dire que~$\sch O_{X,\phi(y)}\to \sch O_{Y,y}$ est local signifie que si~$f$ est un
élément de~$\sch O_{X,\phi(y)}$, 
alors~$f$ appartient à l'idéal maximal
de~$\sch O_{X,\phi(y)}$ si et seulement si~$f^*\phi$ appartient à l'idéal maximal de~$\sch O_{Y,y}$. En termes plus imagés,
cela se traduit par l'équivalence
~$$f(\phi(y))=0\iff (\phi^*f)(y)=0.$$

\trois{ex-morlocann}
{\bf Exemples.} Dans chacun des exemples
classiques~\ref{esp-top-morann},
\ref{var-diff-morann}
et~\ref{var-hol-morann}, l'application~$\phi^*$ est simplement 
$f\mapsto f\circ \phi$ : on a donc tautologiquement~$(\phi^*f)(y)=0\iff f(\phi(y))=0$,
et~$\phi$ est dès lors à chaque fois un morphisme d'espaces
localement annelés. 

\trois{philo-morlocann}
Soit~$\phi: (Y,\sch O_Y)\to (X,\sch O_X)$
 un morphisme d'espaces localement annelés. 
On ne peut pas espérer que la formule
~$\phi^*f=f\circ \phi$ soit valable sans hypothèse supplémentaire : celle-ci n'a
en effet simplement {\em aucun sens}
en général, puisque~$\sch O_X$ n'est pas
nécessairement un faisceau de fonctions à valeurs
dans un corps fixé, pour les deux raisons évoquées ci-dessus
(\ref{corps-varie}
et~\ref{nilp-locann}). 
Mais on va voir qu'elle est tout de même, d'une certaine manière,
aussi valable qu'il est possible.

 \medskip
 Soit~$y$
 un point de~$Y$. 
 Comme~$\phi$
 est un morphisme d'espaces localement annelés, on
 a~$f(\phi(y))=0\iff( \phi^*f)(y)=0$ pour tout~$f \in \sch O_{X,\phi(y)}$. En conséquence,
 $\phi^* : \sch O_{X,\phi(y)}\to \sch O_{Y,y}$ induit par passage au quotient un
 plongement~$\kappa(\phi(y))\hookrightarrow \kappa(y)$, de sorte que le diagramme
 $$\diagram \sch O_{Y,y}\rto&\kappa(y)
 \\ \sch O_{X,\phi(y)}\uto\rto&\kappa(\phi(y))\uto\enddiagram$$ 
 commute. Autrement dit, 
 modulo le plongement de~$\kappa(\phi(y))$
 dans~$\kappa(y)$, on a pour toute~$f\in \mathscr O_{X,\phi(y)}$
 l'égalité 
 $$(\phi^*f)(y)=f(\phi(y)).$$
 Elle évoque irrésistiblement,  comme annoncé, l'égalité~$f^*\phi=f\circ \phi$ ; mais répétons
 qu'il serait illicite de la traduire ainsi puisque~$f$ ne peut pas en général s'interpréter
 comme une vraie fonction naïvement composable avec~$\phi$. 
 
\trois{morlocann-fonct}
Par contre, si~$k$ est un corps et si~$\sch O_Y$
 et~$\sch O_X$ sont des faisceaux de fonctions à valeurs
 dans~$k$ comme dans~\ref{esp-top-locann}, alors 
 l'égalité ci-dessus signifie précisément que~$\phi^*f=f\circ \phi$. 
 Dans ce contexte, un morphisme d'espaces localement 
 annelés est donc simplement une application continue~$\phi : Y\to X$
 telle que
 la fonction~$f\circ \phi$ appartienne à~$\sch O_Y$ pour fonction~$f$
 appartenant à~$\sch O_X$, et~$\phi^*$ est obligatoirement
 donné par la formule~$\phi^*f=f\circ \phi$. 
 
 \deux{imm-ouverte-esplocann}
 Soit~$X$ un espace localement annelé et soit~$U$ un ouvert de~$X$. 
 L'immersion canonique d'espaces annelés $j:U\hookrightarrow X$ 
 ({\em cf.}
 \ref{rest-ouv-morann})
 est alors un morphisme d'espaces localement annelés (les
 morphismes induits au niveau des fibres sont des isomorphismes, car
 si~$x\in U$
 l'anneau local~$\sch O_{U,x}$ s'identifie canoniquement à~$\sch O_{X,x}$). 
 
 Si~$\phi : Y\to X$ est un morphisme
 d'espaces localement annelés tel que~$\phi(Y)$
 soit contenu dans~$U$,
 l'unique morphisme d'espaces annelés~$\psi : Y\to U$ tel que
 $j\circ \psi=\phi$ est en fait un morphisme d'espaces
 localement annelés (là encore parce que~$\sch O_{U,x}=\sch O_{X,x}$
 pour tout~$x\in U$). 
 
 Autrement dit, dans la catégorie
 des espaces localement annelés, on observe
 le phénomène
 déjà constaté dans la catégorie des espaces annelés : toute factorisation
 {\em ensembliste}
 par un ouvert est automatiquement
{\em morphique}. 

Le morphisme
$\psi$ est le seul
morphisme d'espaces localement
annelés dont la composée avec~$j$ soit égale à~$\phi$
(puisque c'est déjà le cas dans la catégorie des espaces annelés, {\em cf.}~\ref{rest-ouv-morann}).
Cela signifie que
$((U,\sch O_U),j)$ représente le foncteur qui envoie un espace
localement annelé~$(Y,\sch O_Y)$ 
sur
$$\{\phi\in \hom_{\mathsf{Esp-loc-ann}}((Y,\sch O_Y), (X,\sch O_X)),\;\;\;\phi(Y)\subset U\}.$$

\trois{imm-ouverte-esplocann}
Un morphisme d'espaces
localement annelés~$\phi \colon Y\to X$
est appelé une {\em immersion ouverte}
s'il induit un isomorphisme entre~$Y$ et un ouvert de~$X$. 

\deux{lim-ind-locann}
Soit~$\sch D=((X_i), (E_{ij})$ un diagramme
dans la catégorie des espaces localement annelés. 
Le lecteur vérifiera que sa limite inductive dans la catégorie des espaces annelés, 
construite au~\ref{lim-ind-espann}, est en fait localement annelée,
et s'identifie à la limite inductive de~$\sch D$ dans la catégorie des espaces
localement annelés. Soit~$\lambda_i \colon X_i \to \limind \sch D$ la flèche structurale.

\deux{limind-locann-comment}
{\bf Remarque.}
Ainsi, l'inclusion de la catégorie des espaces localement annelés
dans celle des espaces annelés commute aux limites inductives. Elle a une bonne
raison pour ce faire : elle admet un adjoint à droite, que nous décrirons plus loin
pour la curiosité du lecteur, mais dont nous ne nous servirons pas. 

\deux{lim-proj-locann}
Exactement comme dans le cas des espaces annelés, on déduit
de~\ref{imm-ouverte-esplocann}
que si~$X$ est un espace localement annelé et si~$U$ et~$V$
sont deux ouverts de~$X$ alors~$U\cap V$
s'identifie au produit fibré~$U\times_X V$. 
  
 \subsection*{Une conséquence géométrique du lemme de Nakayama}

 \deux{intro-oxn-f}
 Soit~$(X,\sch O_X)$ 
 un espace
 localement
 annelé. Si~$\sch F$ est un~$\sch O_X$-module, si~$x\in X$
 et si~$f$ est une section de~$\sch F$ sur un voisinage ouvert~$U$ de~$x$, 
 on se permettra, lorsque le contexte est clair, 
 de noter encore~$f$ l'image~$f_x$ de~$\sch F$ dans~$\sch F_x$. 
 
 \deux{f-kappax}
 Soit~$\sch F$ un~$\sch O_X$-module, et soit~$x\in X$. On
 désigne par~$\got m_x$ l'idéal maximal de
 $\sch O_{X,x}$. 
 Le~$\kappa(x)$-espace
 vectoriel~$\kappa(x)\otimes_{\sch O_{X,x}}\sch F_x=\sch F_x/\got m_x\sch F_x$
 sera plus simplement noté~$\kappa(x)\otimes_x \sch F$. Si~$f$ est une section
 de~$\sch F$ définie au voisinage de~$x$, son image dans
 $\kappa(x)\otimes_x \sch F$ sera notée~$f(x)$ (cette notation est compatible
 avec celle déjà utilisée lorsque~$\sch F=\sch O_X$).  L'application
 $f\mapsto f(x)$ induit par sa définition même
 une surjection de~$\sch F_x$
 vers~$\kappa(x)\otimes_x \sch F$.

 \deux{surj-oxn-f}
 Soit~$(e_1,\ldots, e_n)$ une famille de sections de~$\sch F$ sur~$X$, 
 et soit~$\phi : \sch O_X^n\to \sch F$ le morphisme induit. 
 On dit que~$(e_1,\ldots, e_n)$ 
 {\em engendre~$\sch F$}
 si~$\phi$ est surjectif. Ce signifie 
 que
 le morphisme
 induit~$\sch O_{X,x}^n\to \sch F_x$ est surjectif pour tout~$x$, c'est-à-dire encore
 que les~$e_i$ engendrent le~$\sch O_{X,x}$-module $\sch F_x$
 pour tout~$x$. 
 Si c'est le cas, les~$e_i(x)$ engendrent
 pour tout~$x$ l'espace vectoriel~$\kappa(x)\otimes_x \sch F=\sch F_x/\got m_x\sch F_x$, 
 qui est donc de dimension au plus~$n$. 
 
 \deux{def-loctf-oxmod}
 {\bf Définition.}
 Un~$\sch O_X$-module~$\sch F$ est dit
 {\em localement de type fini}
 si pour tout~$x\in X$ il existe un voisinage 
 ouvert~$U$ de~$x$ et une famille finie de sections
 de~$\sch F$ sur~$U$ qui engendrent~$\sch F|_U$. 
 
 Si c'est le cas, il résulte de~\ref{surj-oxn-f}
 que~$\kappa(x)\otimes_x \sch F$ est de dimension 
 finie pour tout~$x$. 
 
 \deux{nakayama-geom}
 {\bf Proposition.}
 {\em Soit~$\sch F$ un~$\sch O_X$-module
 localement de type fini, et soit~$x$ un point de~$X$. 
 Soit~$(e_1,\ldots, e_n)$ 
 une famille
 de
 sections de~$\sch F$
 sur un voisinage
 ouvert~$U$ de~$x$. 
 Les assertions suivantes sont équivalentes :
 
 \medskip
 i) les~$e_i(x)$ engendrent le~$\kappa(x)$-espace vectoriel~$\kappa(x)\otimes
 \sch F_x$ ; 
  
  ii) il existe un voisinage ouvert~$V$ de~$x$
  dans~$U$ tel que les ~$e_i$
 engendrent~$\sch F|_V$. } 
 
 \medskip
 {\em Démonstration.}
 L'implication~ii)$\Rightarrow$i) a été vue au~\ref{surj-oxn-f}. 
 Supposons maintenant que~i) soit vraie. 
 Les~$e_i(x)$ engendrant~$\kappa(x)\otimes_x \sch F=\sch F_x/\got m_x\sch F$,
 le lemme de Nakayama assure que les~$e_i$ engendrent
 le~$\sch O_{X,x}$-module
 $\sch F_x$. 
 
 \medskip
 Par ailleurs, comme~$\sch F$
 est localement de type fini, il existe un voisinage
 ouvert~$V$ de~$x$ dans~$U$, et une famille~$(f_1,\ldots, f_m)$ de
 sections de~$\sch F$ sur~$V$
 qui engendrent~$\sch F|_V$. 
 
 \medskip
 Comme les~$e_i$ engendrent~$\sch F_x$, il existe une famille~$(a_{ij})$ 
 d'éléments de~$\sch O_{X,x}$ tels que~$f_j=\sum_i a_{ij}e_i$ pour tout~$j$. 
 Quitte à restreindre~$V$, on peut supposer que les~$a_{ij}$ sont définies sur~$V$,
 et que l'égalité~$f_j=\sum a_{ij}e_i$ vaut dans~$\sch F(V)$. 
 
 \medskip
 Soit~$y\in V$. L'égalité~$f_j=\sum a_{ij}e_i$
 vaut dans~$\sch F_y$ ; ce dernier étant engendré par les~$f_j$
 (puisqu'elle engendrent~$\sch F|_V$), il
 est dès lors également engendré par les~$e_i$. Ceci étant vrai pour tout~$y\in V$, 
 les~$e_i$ engendrent~$\sch F|_V$, ce qui achève la démonstration.~$\Box$  
 
 \trois{comment-naka-geom}
 {\bf Commentaires.}
 L'étape cruciale de la preuve ci-dessus, celle durant laquelle
 «il se passe vraiment quelque chose», est 
 l'utilisation du lemme de Nakayama pour garantir que les~$e_i$ engendrent~$\sch F_x$ ; 
 le reste n'est qu'une application directe de la définition des germes en~$x$, couplée à un tout petit
 peu d'algèbre linéaire. 
 
\medskip
On peut donc
considérer la proposition~\ref{nakayama-geom}
comme une traduction géométrique 
du lemme de Nakayama, traduction qui se présente essentiellement
sous la forme d'un
{\em passage du ponctuel au local}
(pour le caractère générateur d'une famille finie de sections). 

\trois{null-ponct-loc}
Mentionnons un cas particulier important
de la proposition~\ref{nakayama-geom}, qui met
particulièrement bien en lumière
ce passage du ponctuel au local : l'espace
vectoriel~$\kappa(x)\otimes_x \sch F$ est nul si et seulement
si il existe un voisinage ouvert~$V$ de~$x$ dans~$X$ tel
que~$\sch F|_V$
soit nul (appliquer la proposition
à la {\em famille vide}
de sections).  

\deux{coro-naka-geom}
{\bf Corollaire.}
{\em Soit~$(X,\sch O_X)$ un espace localement
annelé et soit~$\sch F$ un~$\sch O_X$-module 
localement de type fini. La fonction
$$r: x\mapsto \dim_{\kappa(x)} \kappa(x)\otimes_x \sch F$$
est semi-continue supérieurement, c'est-à-dire
que pour tout entier~$d$, l'ensemble des points~$x$
de~$X$ tels que~$r(x)\leq d$ est ouvert.}

\medskip
{\em Démonstration.}
Soit~$x$ un point de~$X$ en lequel~$r(x)\leq d$. 
Choisissons une famille~$e_1,\ldots, e_n$ de sections de~$\sch F$,
définies sur un voisinage ouvert~$U$ de~$x$, et telles
que les~$e_i(x)$ forment une base de~$\kappa(x)\otimes_x \sch F$. 
Comme~$r(x)\leq d$, on a~$n\leq d$. En vert de la proposition~\ref{nakayama-geom},
il existe un voisinage ouvert~$V$ de~$x$
dans~$U$ tel que les~$e_i$ engendrent~$\sch F|_V$ ;  en conséquence, 
on a~$r(y)\leq n\leq d$ pour tout~$y\in V$ (\ref{surj-oxn-f}).~$\Box$

\trois{comment-zero-ouvert}
Notons un cas particulier fondamental, dont l'énoncé
peut apparaître contre-intuitif au premier abord : le sous-ensemble~$U$
de~$X$ formé des point~$x$ tels que~$r(x)=0$, c'est-à-dire encore
tel que~$\kappa(x)\otimes_x \sch F=\{0\}$, est un {\em ouvert}. Remarquons
de surcroît que
le faisceau~$\sch F|_U$ a toutes ses fibres nulles d'après
le cas particulier de la proposition~\ref{nakayama-geom}
signalé au~\ref{null-ponct-loc}, et est donc lui-même nul. 

\trois{exemple-r-scs}
{\bf Un exemple.}
Soit~$X$ une variété différentielle 
munie du faisceau~$\sch O_X$ des fonctions~$\sch C^\infty$. 
Soit~$x\in X$ et soit~$U$ un ouvert de~$X$. On note~$\sch I(U)$ l'idéal de~$\sch O_X(U)$
défini comme suit : 

\medskip
$\bullet$ si~$x\in U$ alors~$\sch I(U)$ est l'ensemble des fonctions 
appartenant à~$\sch O_X(U)$ et
s'annulant en~$x$ ; 

$\bullet$ si~$x\notin U$ alors~$\sch I(U)=\sch O_X(U)$. 

\medskip
Il est immédiat que~$\sch I$ est un sous-faisceau de~$\sch O_X$ ; soit~$\sch F$ le quotient
$\sch O_X/\sch I$. Il est (localement)
de type fini par construction ; nous allons déterminer la fonction~$r: y\mapsto \dim_{\RR}\RR\otimes_y \sch F$
(notez que~$\kappa(y)=\RR$ pour tout~$y\in X$). 

\medskip
Soit~$i$ l'inclusion de~$\{x\}$ dans~$X$. Le faisceau~$i_*\underline{\RR}$ envoie un ouvert~$U$
de~$X$
sur~$\RR$ si~$U$ contient~$x$, et sur~$\{0\}$ sinon
(c'est un «faisceau gratte-ciel supporté en~$x$»). Il hérite d'une structure naturelle de~$\sch O_X$-module, définie
comme suit : sur un ouvert~$U$ ne contenant pas~$x$, il n'y a rien à faire ; sur un ouvert~$U$ contenant~$x$, on fait agir
$\sch O_X(U)$ sur $i_*\underline{\RR}(U)=\RR$ par la formule~$(f,\lambda)\mapsto f(x)\lambda$. 

\medskip
On dispose d'une surjection
$\sch O_X$-linéaire
naturelle de~$\sch O_X$
sur~$i_*\underline{\RR}$ : 
là encore, sur un
ouvert~$U$ ne contenant pas~$x$, il n'y a rien à faire ; et sur un ouvert~$U$ contenant~$x$, on 
envoie une fonction~$f\in \sch O_X(U)$ sur~$f(x)\in \RR=i_*\underline{\RR}(U)$. Par construction, le noyau de cette surjection
est~$\sch I$. En conséquence, $\sch F\simeq i_*\underline{\RR}$. Il s'ensuit que~$\RR\otimes_y\sch F=\{0\}$
si~$y\neq x$, et que~$\RR\otimes_x \sch F=\RR$ ; il vient
$$r(y)=0\;\text{si}\; y\neq x\;\;\text{et}\;\;r(x)=1.$$

\section{Faisceaux localement libres de rang~$1$}
\markboth{Espaces localement annelés}{Faisceaux localement libres de rang~$1$}

\subsection*{Définition, exemples, premières propriétés}

\deux{def-faisc-loc}
Soit~$X$ un espace localement annelé
et soit~$n$ un entier. 
Un~$\sch O_X$-module~$\sch F$ est dit
{\em localement libre de rang fini} (resp. {\em de rang~$n$})
si tout point de~$X$ possède un voisinage ouvert~$U$ tel que~$\sch F|_U$ soit
isomorphe à~$\sch O_U^m$ pour un certain entier~$m$ (resp. à~$\sch O_U^n$).

\trois{ex-faisc-loclibre}
Soit~$\sch F$ un~$\sch O_X$-module localement libre de rang~$n$. 
Il est dit {\em trivial}
s'il est isomorphe à~$\sch O_X^n$. 
On dira qu'une famille~$(U_i)$ d'ouverts de~$X$
{\em trivialise}
$\sch F$ si~$\sch F|_{U_i}$ est trivial pour tout~$i$.

\trois{rang-faisc-loclibre}
Soit~$\sch F$ un~$\sch O_X$-module
localement libre sur~$X$, et soit~$x\in X$. Par hypothèse, 
il existe un voisinage ouvert~$U$ de~$x$ dans~$X$, un entier~$m$ et un isomorphisme
$\sch F|_{U}\simeq \sch O_U^m$. On a alors~$\sch F_x\simeq \sch O_{X,x}^m$, et~$\kappa(x)\otimes \sch F
\simeq \kappa(x)^m$. L'entier~$m$ est ainsi uniquement déterminé : c'est le rang du module libre
$\sch F_x$ sur l'anneau non nul~$\sch O_{X,x}$, ou encore la dimension du~$\kappa(x)$-espace vectoriel
$\kappa(x)\otimes \sch F$ ; on dit que~$m$ est le {\em rang de~$\sch F$ en~$x$}.
Le rang de~$\sch F$ en tout point de~$U$
est encore égal à~$m$ par définition ; le rang de~$\sch F$ 
apparaît ainsi comme une {\em fonction localement constante} de~$X$
dans~$\NN$. Si~$X$ est connexe, cette fonction
est nécessairement constante, et a donc une valeur bien définie si
$X$ est
de surcroît
non vide, valeur
que l'on appelle encore le rang de~$\sch F$. 

\trois{fonctorialite-loclibre}
Soit~$\psi \colon Y \to X$ un morphisme
d'espaces localement annelés. Si~$\sch F$ est un~$\sch O_X$-module
localement libre de rang fini (resp. de rang~$n$) alors
$\psi^*\sch F$ est un~$\sch O_Y$-module localement libre de rang fini (resp.~$n$) : 
c'est une conséquence immédiate de~\ref{adj-ocmod}. 

\deux{exemple-loclibre-nontriv}
Nous allons donner un premier exemple de~$\sch O_X$-module
localement libre non trivial, dans le contexte de la géométrie analytique 
complexe -- nous en rencontrerons d'autres plus loin en théorie des schémas. 

\medskip
Soit~$S$ la sphère de Riemann $\CC\cup\{\infty\}$, munie du faisceau~$\sch O_S$ 
des fonctions holomorphes. Soit~$\sch I\subset \sch O_S$ le faisceau d'idéaux défini
comme suit : si~$U$ est un ouvert de~$S$ ne contenant pas~$0$ (resp. contenant~$0$)
alors~$\sch I(U)=\sch O_S(U)$ (resp.~$\sch I(U)$ est l'ensemble des fonctions~$f\in \sch O_S(U)$
s'annulant en~$0$). 

\trois{i-loclibre-1}
{\em Le~$\sch O_S$-module~$\sch I$ est localement libre de rang~$1$}.
En effet, on a d'une part par définition
$\sch I|_{S\setminus\{0\}}=\sch O_{S\setminus\{0\}}$ ; et d'autre part,~$f\mapsto zf$ 
définit un isomorphisme de~$\sch O_{\CC}$ sur~$\sch I|_{\CC}$
(une fonction holomorphe au voisinage de l'origine s'y annule si et seulement si
elle est multiple de~$z$). 

\trois{i-pas-trivial}
{\em Le~$\sch O_S$-module localement libre~$\sch I$ n'est pas trivial}. 
En effet, rappelons que les seules fonctions holomorphes définies sur~$S$ toute entière
sont les constantes (c'est le «principe du maximum»). 
Autrement dit, on a~$\sch O_S(S)=\CC$, et~$\sch I(S)=\{0\}$ (une fonction constante s'annule
à l'origine si et seulement si elle est nulle !). Par conséquent, $\sch I$ {\em n'est pas}
isomorphe à~$\sch O_S$.

\deux{pro-faisc-loclib}
Soient~$\sch F$ et~$\sch G$ deux~$\sch O_X$-modules localement libres
de rang fini. 

\trois{hom-loclibre}
On déduit de~\ref{exemple-homint}
que le~$\sch O_X$-module~$\homfsc (\sch F,\sch G)$ est
localement libre de rang fini
et que son rang (comme fonction sur~$X$, {\em cf. supra})
est égal au produit des rangs de~$\sch F$ et~$\sch G$. 

\medskip
En particulier, le dual~$\sch F^\vee$ est localement libre
de même rang que~$\sch F$. 

\trois{tens-loc-libre}
Il est immédiat que le produit tensoriel~$\sch F\otimes_{\sch O_X}\sch G$
est localement libre de rang égal au produit des rangs de~$\sch F$ et~$\sch G$ : c'est 
une simple conséquence du fait 
que~$(\sch O_X)^n\otimes_{\sch O_X}(\sch O_X^m)\simeq \sch O_X^{nm}$. 

\medskip
On déduit par ailleurs de~\ref{dual-endos-libre}
que~$\sch F^\vee \otimes_{\sch O_X}\sch G \to \homfsc(\sch F,\sch G)$ est un isomorphisme,
puisque cette propriété se teste localement. 

\trois{loc-libre-bidual}
On déduit de même de~\ref{bidual-fasci}
que~$\sch F\to \sch F^{\vee\vee}$ est un isomorphisme. 

\deux{loc-libre-rg1}
Soit~$\sch L$ un~$\sch O_X$-module localement libre de rang~$1$. Il résulte
de~\ref{homoth-faisc}
(toujours {\em via}
un raisonnement local) que la flèche naturelle~$\sch O_X\to \fscend \sch L$ est un isomorphisme
de~$\sch O_X$-algèbres. En composant
sa réciproque avec l'isomorphisme
entre~$\sch L\otimes_{\sch O_X}\sch L^\vee$
et~$\fscend \sch L$
fourni par~\ref{tens-loc-libre}, 
on obtient un isomorphisme
$$\sch L\otimes_{\sch O_X}\sch L^\vee\simeq
\sch O_X.$$ Il découle immédiatement 
des définitions des différentes flèches en jeu que cet
isomorphisme
est simplement
donné par la formule
$$s\otimes \phi\mapsto \phi(s).$$
On peut bien entendu vérifier directement que celle-ci
définit bien un isomorphisme : en raisonnant localement on se ramène au cas où~$\sch F=\sch O_X$,
pour lequel c'est évident. 

\deux{pic-esplocann}
Si~$\sch L$ et~$\sch L'$ sont deux~$\sch O_X$-modules localement
libres de rang~$1$, le~$\sch O_X$-module~$\sch L\otimes_{\sch O_X}\sch L'$ est 
lui aussi
localement libre de rang~$1$
d'après~\ref{tens-loc-libre}. 

\medskip
Le produit tensoriel induit de ce fait une loi de composition 
sur l'ensemble~${\rm Pic}\;X$
des classes d'isomorphie de~$\sch O_X$-modules localement libres
de rang~$1$. Elle est associative, commutative, et possède un élément neutre : la 
classe de~$\sch O_X$. Il résulte de~\ref{loc-libre-rg1}
que la classe~$[\sch L]$ d'un~$\sch O_X$-module localement
libre~$\sch L$ admet un symétrique, à savoir~$[\sch L^\vee]$. En conséquence, 
{\em ${\rm Pic}\;X$ est un groupe abélien}, appelé le {\em groupe de Picard}
de l'espace localement annelé~$X$. 

\subsection*{Sections inversibles et trivialisations}

\deux{sectinv-trivial} 
Soit~$\sch L$ un~$\sch O_X$-module localement libre de rang~$1$. Si~$s$
est une section globale de~$\sch L$, on note~$D(s)$ l'ensemble des points~$x$
tels que~$s(x)\neq 0$ ; cette notation est compatible avec celle déjà utilisée lorsque~$\sch L=\sch O_X$. 

\trois{pro-tens-section}
Soit~$\sch L'$ un second~$\sch O_X$-module localement libre de rang~$1$, et soient~$s$ et~$s'$
des sections globales respectives de~$\sch L$ et~$\sch L'$.
Le produit tensoriel~$s\otimes s'$ est une section de~$\sch L\otimes_{\sch O_X}\sch L'$, et l'on a~$D(s\otimes s')=D(s)\cap D(s')$. 
En effet, comme la propriété est locale sur~$X$, on peut supposer~$\sch L=\sch L'=\sch O_X$. Comme le morphisme
canonique
$$\sch O_X\otimes_{\sch O_X}\sch O_X\to \sch O_X, f\otimes g\mapsto fg$$ est un isomorphisme, il en résulte que~$D(s\otimes s')=D(ss')=D(s)\cap D(s')$, comme
annoncé. 

\trois{sect-inv-iso}
Soit~$s$ une section globale de~$\sch L$. Il existe un unique morphisme~$\ell$ 
de~$\sch O_X$ vers~$\sch L$ qui envoie~$1$ sur~$s$ : celui donné par la formule
$f\mapsto fs$. Nous allons montrer que~$\ell$ est un isomorphisme si et seulement
si~$D(s)=\varnothing$ ; si c'est le cas, nous dirons que~$s$ est {\em inversible}. 

\medskip
Si~$\ell$ est un isomorphisme, on a~$D(s)=D(1)=\varnothing$. Réciproquement, supposons
que~$D(s)=\varnothing$. Pour montrer que~$\ell$ est un isomorphisme, on peut raisonner
localement et donc supposer que~$\sch L=\sch O_X$. Dans ce cas, $s$ est une {\em fonction}
inversible, 
et~$\ell$ est donc bien un isomorphisme de réciproque~$f\mapsto f/s$. 

\trois{recap-sectinv}
La flèche~$\ell\mapsto \ell(1)$ établit donc une bijection
entre l'ensemble des isomorphismes
de~$\sch O_X$ vers~$\sch L$ et l'ensemble des sections inversibles de~$\sch L$. Si~$s$
est une telle section, l'isomorphisme qui lui correspond est~$f\mapsto fs$, et sa réciproque
sera notée~$t\mapsto t/s$. 

\medskip
Soit~$s$ une section inversible de~$\sch L$ et
soit~$s'$ une section inversible de~$\sch L'$. La
section~$s\otimes s'$ de~$\sch L\otimes_{\sch O_X}\sch L'$ est inversible d'après~\ref{pro-tens-section}. 
Si~$f$ et~$g$ sont deux sections de~$\sch O_X$ (sur un ouvert de~$X$) on a~$fs\otimes gs'=fg (s\otimes s')$ ; 
il en résulte que si~$t$ et~$t'$ sont des sections respectives de~$\sch L$ et~$\sch L'$ alors
$$(t\otimes t')/(s\otimes s')=(t/s)\cdot (t'/s').$$

\trois{mor-secinv}
Soit~$s$ une section inversible de~$\sch L$ et soit~$s'$ une section de~$\sch L'$. 
Il est immédiat qu'il existe un et un seul morphisme de~$\sch L$ vers~$\sch L'$ envoyant~$s$ sur~$s'$, donné
par la formule~$t\mapsto (t/s)s'$. Ce morphisme est un isomorphisme si et seulement si~$s'$ est inversible. 
Cette condition est en effet clairement nécessaire, et si elle est satisfaite on vérifie aussitôt que~$\tau \mapsto (\tau/s')s$
est un inverse à gauche et à droite du morphisme considéré. 

\deux{intuition-loclibre-rg1}
On peut penser à
un~$\sch O_X$-module~$\sch L$
localement libre de rang~$1$
comme
à une famille d'espaces vectoriels de dimension~$1$ (celle des~$\kappa(x)\otimes \sch L$
lorsque~$x$ parcourt~$X$), les sections inversibles correspondant aux familles d'éléments non nuls, 
c'est-à-dire de bases. 

\trois{loclibre-rg1-canon}
Un~$\sch O_X$-module
localement libre de rang~$1$
n'admet pas nécessairement
de section inversible (puisqu'il n'est pas forcément trivial, {\em cf.}
\ref{exemple-loclibre-nontriv}). Modulo l'interprétation
donnée au~\ref{intuition-loclibre-rg1}, cela correspond au fait que si un espace
vectoriel abstrait de dimension~$1$ sur un corps admet toujours une base, 
il n'en admet aucune qui soit canonique -- et il n'y a donc en général aucune raison de
pouvoir faire
un choix «cohérent» de bases dans une famille de tels espaces vectoriels. 

\trois{loclibre-tens-unites}
Plusieurs faits mentionnés ci-dessus 
(\ref{pro-tens-section}, \ref{recap-sectinv})
font apparaître une certaine parenté entre le produit
tensoriel des sections de~$\sch O_X$-modules localement libres 
de rang~$1$ et le produit classique des fonctions. C'est, là encore, la
déclinaison en famille d'un phénomène qui existe déjà au niveau des espaces vectoriels
de dimension 1, que nous allons maintenant expliquer. 

\medskip
Soient~$L$ et~$L'$ deux espaces vectoriels de dimension 1 sur un corps~$k$. 
{\em Si l'on choisit une base de~$L$ et une base de~$L'$}, on obtient deux isomorphismes
$L\simeq k$ et~$L'\simeq k$, modulo lesquels il devient possible de multiplier un élément de~$L$ 
par un élément de~$L'$ (on obtient un scalaire). Cette opération est évidemment hautement non canonique : 
elle dépend de façon cruciale du choix des bases. 

\medskip
Pour «multiplier» les éléments de~$L$ par ceux de~$L'$ de manière intrinsèque, 
on utilise le produit tensoriel
$$\otimes \colon L\times L'\to L\otimes_k L'.$$
Il y a bien entendu
un lien entre nos deux constructions d'un produit : si~$v$ est une base de~$L$ et~$v'$ une base de~$L'$ alors~$v\otimes v'$ est une
base de~$L\otimes_k L'$ ; chacune de ces bases identifie l'espace vectoriel correspondant à~$k$, et le diagramme
$$\xymatrix{
{L\times L'} \ar[rrr]^{(\lambda, \mu)\mapsto \lambda\otimes \mu}&&&{L\otimes_k L'}\\
{k\times k}\ar[rrr]_{(a,b)\mapsto ab}\ar[u]^{(a,b)\mapsto (av, bv')}_\simeq &&&k\ar[u]_{a\mapsto a v\otimes v'}^\simeq}$$
commute. 

\medskip
Vous avez déjà maintes fois rencontré, probablement de façon implicite, ce genre de considérations en... physique. 
La mesure d'une grandeur y est en effet le plus souvent non un scalaire bien déterminé, mais un élément d'un
$\RR$-espace vectoriel réel de dimension~$1$ (celui des temps, celui des longueurs...), dont le choix d'une base revient à celui 
d'une unité de référence. Une phrase courante comme «lorsqu'on multiplie deux longueurs, on obtient une aire»
évoque une opération qui, {\em conceptuellement}, ne consiste pas à 
multiplier deux nombres réels (même si {\em en pratique}, c'est évidemment ce que l'on
fait, {\em une fois choisi un système d'unités}), mais à appliquer le produit tensoriel~$L\times L\to L^{\otimes 2}$, où~$L$ est l'espace
vectoriel des longueurs (et~$L^{\otimes 2}$ celui des aires). Si l'on note~$m$ la base de~$L$ correspondant
au choix du mètre comme unité de longueur, la base~$m\otimes m$ de~$L^{\otimes 2}$ est celle qui 
correspond au mètre carré comme unité d'aire. 

La division d'une longueur par un temps (non nul) pour obtenir une vitesse est un tout
petit peu plus délicate à décrire en termes intrinsèques : si~$T$ désigne l'espace des temps, 
elle consiste à associer à un couple~$(\ell, t)$ de~$L\times (T\setminus \{0\})$ l'élément
$l\otimes t^{-1}$ de l'espace des vitesses~$L\otimes_{\RR} T^\vee$, où~$t^{-1}$ est la forme linéaire
$\tau \mapsto \tau/t$ sur l'espace vectoriel~$T$ (si l'on préfère, on peut décrire~$t^{-1}$ 
comme la base duale de~$t$). Si~$s$ désigne la base de~$T$ correspondant au choix de la
seconde comme unité de temps, la base~$m\otimes s^{-1}$ de~$L\otimes_{\RR}T^\vee$ est 
celle qui correspond au mètre par seconde comme unité de vitesse. 

\subsection*{Cocycles}

\deux{cocy-cobo}
Nous allons maintenant définir des objets (qui sont des cas
particuliers de constructions cohomologiques très générales) dont nous nous servirons
pour fabriquer des~$\sch O_X$-modules localement libres de rang~$1$ par recollement.

\trois{def-cocycle}
Soit~$(U_i)_{i\in I}$ un recouvrement ouvert de~$X$.
Un {\em cocycle
subordonné au recouvrement~$(U_i)$} est la donnée, pour tout couple~$(i,j)$ 
d'indices, d'un élément~$f_{ij}\in \sch O_X(U_i\cap U_j)^\times$, satisfaisant les conditions suivantes
(on pourrait déduire~ii) de~i) et~iii), mais nous avons 
préféré la faire figurer explicitement) : 

\medskip
i) $f_{ii}=1$ pour tout~$i$ ; 

ii) $f_{ij}=f_{ji}^{-1}$ pour tout~$(i,j)$ ; 

iii) $f_{ij}\cdot f_{jk}=f_{ik}$
pour tout~$(i,j,k)$ là où cette égalité a un sens, c'est-à-dire sur~$U_i\cap U_j\cap U_k$. 

\medskip
L'ensemble~$Z_{(U_i)}$ des cocycles
subordonnés à~$(U_i)$ hérite d'une structure naturelle de groupe, induite par la multiplication
des fonctions. 

\trois{def-cobord}
{\em Un premier exemple : les cobords}.
Donnons-nous pour tout~$i$ un élément~$a_i\in \sch O_X(U_i)^\times$. 
Pour tout~$(i,j)$, on note~$f_{ij}$
l'élément inversible~$(a_i/a_j)$ de~$\sch O_X(U_i\cap U_j)\ti$. 
La famille
$(f_{ij})$ est alors un cocycle subordonné à~$(U_i)$. Les cocycles de cette forme sont appelés
{\em cobords} ; l'ensemble~$B_{(U_i)}$ des cobords subordonnés à~$(U_i)$ est un sous-groupe
de~$Z_{(U_i)}$. 

\deux{fibres-et-cocycles}
{\bf Cocyles et~$\sch O_X$-modules localement libres de rang~$1$.}
On fixe un~$\sch O_X$-module $\sch L_0$, localement libre de rang~$1$.
Soit~$(U_i)$
un recouvrement ouvert de~$X$. Le but de ce qui suit est de construire une bijection entre
le groupe~$Z_{(U_i)}/B_{(U_i)}$ et l'ensemble~$G_{(U_i),\sch L_0}\subset {\rm Pic}\;X$
des classes d'isomorphie de~$\sch O_X$-modules
$\sch L$ localement libres de rang~$1$ et tels que~$\sch L|_{U_i}\simeq \sch L_0|_{U_i}$ pour tout~$i$ ; nous dirons
plus simplement qu'un tel~$\sch L$ est {\em $\sch L_0$-trivialisé par~$(U_i)$.}
Notons qu'un~$\sch O_X$-module localement libre de rang~$1$
est~$\sch O_X$-trivialisé par~$(U_i)$ si et seulement si il est trivialisé par~$(U_i)$ au sens
de~\ref{ex-faisc-loclibre}.

\trois{fibre-vers-cocycle}
Soit donc~$\sch L$ un~$\sch O_X$-module localement libre de rang~$1$
qui est~$\sch L_0$-trivialisé par~$(U_i)$ ; nous allons
expliquer comment lui associer un élément de~$Z_{(U_i)}/B_{(U_i)}$. 

\medskip
Choisissons pour tout~$i$ un isomorphisme~$\ell_i \colon \sch L_0|_{U_i}\to \sch L|_{U_i}$. Pour tout
couple~$(i,j)$, la composée~$\ell_i\circ \ell_j^{-1}$ est un 
automorphisme de~$\sch L|_{U_i\cap U_j}$, c'est-à-dire, en vertu de~\ref{loc-libre-rg1}, 
une homothétie de rapport~$f_{ij}$ pour une certaine fonction~$f_{ij}\in \sch O_X(U_i\cap U_j)^\times$, 
uniquement déterminée. On vérifie immédiatement que~$(f_{ij})$ est un cocycle subordonné à~$(U_i)$. 

Donnons-nous une seconde collection d'isomorphismes
$\lambda_i  \colon \sch L_0|_{U_i}\to \sch L|_{U_i}$, et soit~$(g_{ij})$ le cocycle qui lui 
est associé par le procédé ci-dessus. Pour tout~$i$, la composée~$\ell_i\circ \lambda_i^{-1}$
est un 
automorphisme de~$\sch L|_{U_i}$, c'est-à-dire, en vertu de~\ref{loc-libre-rg1}, 
une homothétie de rapport~$a_i$ pour une certaine fonction~$a_i\in \sch O_X(U_i)\ti$, 
uniquement déterminée. On vérifie alors aussitôt que~$f_{ij}=g_{ij}\cdot (a_i/a_j)$ pour tout~$(i,j)$. 

Ainsi,
la classe du cocycle~$(f_{ij})$ modulo~$B_{(U_i)}$ ne dépend pas du choix du système~$(\ell_i)$, mais seulement
de~$\sch L$ ; on la note~$\got h(\sch L)$.

\trois{rem-h-classe-iso}
{\em Remarque.}
Si~$\sch L'$ est second~$\sch O_X$-module
localement libre de rang~$1$ et si~$\theta \colon \sch L\to \sch L'$ est 
un isomorphisme, il résulte aussitôt des définitions que le cocycle associé
par la collection d'isomorphismes~$(\theta\circ \ell_i \colon \sch L_0|_{U_i}\simeq \sch L'|_{U_i})$ 
est égal à~$(f_{ij})$. Ainsi, $\got h(\sch L)$ ne dépend que de la classe d'isomorphie de~$\sch L$. 
On peut donc voir~$\got h$ comme une application de~$G_{(U_i),\sch L_0}$ vers~$Z_{(U_i)}/B_{(U_i)}$. 

\trois{cocycle-vers-fibre}
Nous allons maintenant construire une application~$\got l$
de~$Z_{(U_i)}/B_{(U_i)}$ vers~$G_{(U_i),\sch L_0}$,
puis montrer que~$\got l$ et~$\got h$ sont deux bijections réciproques l'une de l'autre.

\medskip
Soit~$(f_{ij})$ un cocycle
subordonné à~$(U_i)$.
Soit~$U$ un ouvert de~$X$. On définit~$\sch L(U)$ comme l'ensemble
des familles~$(s_i)$
appartenant à~$\prod_i \sch L_0(U\cap U_i)$ telles que
$s_i=f_{ij}s_j$ pour tout~$(i,j)$ là où cette égalité a un sens, c'est-à-dire sur~$U\cap U_i\cap U_j$
(en termes informels, on «tord» la condition de coïncidence usuelle par le cocycle~$(f_{ij})$). 

On voit aussitôt
que~$\sch L$ est un faisceau, qui possède une structure naturelle de~$\sch O_X$-module
(la multiplication externe se fait composante par composante). Fixons un indice~$k$. Pour tout ouvert~$U\subset U_k$ et toute
section~$s$
appartenant à~$\sch L_0(U)$, la famille~$(s_i)\in \prod \sch L_0(U\cap U_i)$ 
définie par la formule
$s_i=f_{ik}s$ est une section de~$\sch L$ sur~$U$
(nous laissons le lecteur le vérifier -- cela repose de manière essentielle sur le fait que~$(f_{ij})$ est un cocycle).

On définit par ce biais un morphisme~$\ell_k \colon \sch L_0|_{U_k}\simeq \sch L|_{U_k}$. Il est immédiat
que~$\ell_k$ est un isomorphisme de réciproque~$(s_i)\mapsto s_k$. {\em Autrement dit, sur l'ouvert~$U_k$
on peut identifier~$\sch L$ à~$\sch L_0$ en ne regardant que la composante d'indice~$k$}. Ainsi, $\sch L$ est~$\sch L_0$-trivialisé
par~$(U_i)$. On dit que~$\sch L$ est obtenu en
{\em tordant}
$\sch L_0$ par le cocycle~$(f_{ij})$. 

\trois{rem-bon-cocycle}
{\em Remarque}.
Il résulte des définitions que le cocycle associé au système d'isomorphismes~$(\ell_i)$
par le procédé décrit au~\ref{fibre-vers-cocycle} est précisément~$(f_{ij})$.

\trois{rem-cobord-sanseffet}
Donnons-nous pour tout~$i$ une fonction inversible~$a_i$ sur~$U_i$. Soit~$(g_{ij})$ le cocycle
défini par la formule~$g_{ij}=(a_i/a_j) f_{ij}$, et soit~$\sch M$ le~$\sch O_X$-module localement
libre de rang~$1$ associé à~$(g_{ij})$ par le procédé décrit ci-dessus. 
On montre sans difficulté que la
formule
$(s_i)\mapsto (a_is_i)$ définit
un isomorphisme de~$\sch L$ sur~$\sch M$, de réciproque~$(t_i)\mapsto (a_i^{-1}t_i)$. 

\medskip
La construction du~\ref{cocycle-vers-fibre}
permet ainsi d'associer à toute classe~$h$
appartenant à~$Z_{(U_i)}/B_{(U_i)}$ une 
classe~$\got l(h)\in G_{(U_i),\sch L_0}$. La remarque~\ref{rem-bon-cocycle} 
ci-dessus assure que~$\got h\circ \got l={\rm Id}_{Z_{(U_i)}/B_{(U_i)}}$. 

\trois{fin-cocycles-fibre}
Il reste à montrer
que~$\got l\circ \got h={\rm Id}_{G_{(U_i),\sch L_0}}$. 

\medskip
Soit~$\sch L$ un~$\sch O_X$-module localement libre de rang~$1$
qui est~$\sch L_0$-trivialisé
par~$(U_i)$. On se donne un système d'isomorphismes
$$(\ell_i \colon \sch L_0|_{U_i}\simeq \sch L|_{U_i})\;\;;$$ on lui associe un cocycle~$(f_{ij})$ comme
au~\ref{fibre-vers-cocycle}, puis on associe à~$(f_{ij})$ un~$\sch O_X$-module $\sch L'$ par
le procédé décrit au~\ref{cocycle-vers-fibre}. Il suffit pour conclure de démontrer que~$\sch L'\simeq \sch L$. 

\medskip
Soit~$U$ un ouvert de~$X$ et soit~$(s_i)$ une section de~$\sch L'$ sur~$U$. Par construction, on a pour tout~$(i,j)$
l'égalité~$\ell_i(s_i)=\ell_j(s_j)$ là où elle a un sens, c'est-à-dire sur~$U\cap U_i\cap U_j$. Les~$\ell_i(s_i)$ se recollent donc
en une section de~$\sch L$ sur~$U$. En faisant varier~$U$, on obtient un morphisme~$\sch L'\to \sch L$ ; nous laissons le lecteur
vérifier qu'il s'agit d'un isomorphisme, de réciproque~$s\mapsto (\ell_i^{-1}(s))_i$. 

\trois{conclu-iso-fibres}
Les applications~$\got l$ et~$\got h$ mettent donc comme annoncé
$Z_{(U_i)}/B_{(U_i)}$ et $G_{(U_i),\sch L_0}$ en bijection. 

\medskip
Supposons que~$\sch L_0=\sch O_X$. Il est immédiat que~$G_{(U_i),\sch O_X}$ est un {\em sous-groupe}
de~${\rm Pic}\;X$. Nous laissons au lecteur
le soin de prouver que~$\got h$ et~$\got l$ 
sont alors  
des isomorphismes
{\em de groupes} ; cela repose essentiellement sur le bon comportement du produit tensoriel vis-à-vis
de la multiplication des fonctions (\ref{pro-tens-section}-\ref{recap-sectinv}).

\part{La théorie des schémas}

\chapter{Le spectre comme espace topologique}

\section{Spectre d'un anneau}
\markboth{Le spectre comme espace topologique}{Le spectre d'un anneau}

\subsection*{Motivation et définition} 

\deux{intro-spec}
Soit~$A$ un anneau. La géométrie algébrique à la Grothendieck se propose
de lui associer un objet de nature géométrique, en partant du postulat suivant,
conforme à l'intuition provenant de théories
classiques (géométrie différentielle, 
géométrie analytique complexe, géométrie algébrique au sens des articles FAC et GAGA de Serre...) :  
{\em en géométrie, un «point» est quelque chose en lequel on peut évaluer des fonctions,
le résultat étant à valeurs dans un corps.} 

\deux{predef-spec}
On décide donc d'associer à tout morphisme~$A\to K$, où~$K$ est un corps, un point de notre objet géométrique
à construire --l'idée étant qu'on doit penser au morphisme en question comme à l'évaluation en le point correspondant. 

\trois{relat-equiv-spec}
Mais il y a beaucoup trop de tels morphismes~$A\to K$, lorsque~$K$ varie (au point que ceux-ci ne constituent même pas
un ensemble). Il faut donc en identifier certains pour obtenir un objet raisonnable. On décide ainsi que pour tout diagramme
commutatif 
$$
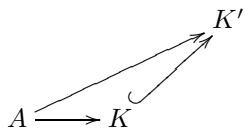
\xymatrix{
&&K'\\
A\ar[rru]\ar[r]&K\ar@{^{(}->}[ru]&}
$$
les morphismes~$A\to K$
et~$A\to K'$ définissent le même point. C'est naturel : il s'agit simplement de dire qu'on ne change
pas un point en agrandissant artificiellement le corps sur lequel il est défini. 

Par exemple, l'évaluation
$P\mapsto P(0)$ est un morphisme de~$\RR[T]$ dans~$\RR$ ; on peut
toujours s'amuser à le voir
comme
un morphisme
de~$\RR[T]$ dans~$\CC$, mais il s'agira encore de l'évaluation en l'origine : 
qu'on considère celle-ci comme un point réel ou un point complexe importe
peu, c'est le «même» point. 

\trois{prem-def-spec}
L'objet que l'on souhaite associer à~$A$ peut donc être défini comme le quotient
de~$\{A\to K\}_{K\;{\rm corps}}$ par la relation qu'engendrent les identifications
mentionnées ci-dessus. 

\medskip
Ce n'est certes pas une définition très tangible. Mais il résulte de~\ref{IDEAUX} 
que ce quotient est en bijection naturelle avec l'ensemble~$\spec A$
des idéaux premiers de~$A$, de la façon suivante :

\medskip
$\bullet$ à la classe d'un morphisme~$A\to K$ on associe le noyau de~$A\to K$ ; 

$\bullet$ à un idéal premier~$\got p$ on fait correspondre
la classe de
la flèche composée~$A\to A/\got p\hookrightarrow {\rm Frac}\;A/\got p$. 

\medskip
De plus, d'après
{\em loc. cit.}, $A\to {\rm Frac}\;A/\got p$
est le plus petit morphisme de sa classe : 
tout morphisme~$A\to K$ appartenant à celle-ci
admet une unique factorisation~$A\to {\rm Frac}\;A/\got p \hookrightarrow K$. 

\deux{spec-def}
L'objet de base associé à un anneau~$A$ par la théorie
des schémas
est donc
l'ensemble~$\spec A$ de ses idéaux premiers. Toutefois, 
pour favoriser l'intuition géométrique, il est préférable de penser aux éléments de~$\spec A$ comme
à des {\em points}, et 
de se rappeler qu'à tout point~$x$
de~$\spec A$ 
{\em correspond}, selon 
les besoins : 

\medskip
$\bullet$ un idéal premier~$\got p$ de~$A$ ; 

$\bullet$ une classe de morphismes~$A\to K$, où~$K$ est un corps, qui admet un plus
petit élément~$A\to \kappa(x)$ que l'on note suggestivement~$f\mapsto f(x)$. 

\medskip
Le lien entre les deux se déduit de~\ref{prem-def-spec} : on
a~$$\got p=\{f\in A, f(x)=0\},$$ le corps~$\kappa(x)$ est
égal à~${\rm Frac}\;A/\got p$ 
et~$A\to \kappa(x)$ est la flèche canonique
de~$A$ vers~${\rm Frac}\;A/\got p$, composée
de la flèche quotient~$A\to A/\got p$ et de l'injection
de~$A/\got p$ dans son corps des fractions. On dit que~$\kappa(x)$ est le
{\em corps résiduel}
du point~$x$. 

\deux{spectre-fonct}
Ainsi, $A$ apparaît comme une sorte d'anneaux de fonctions sur~$\spec A$, au moins
dans le sens où l'on dispose pour tout~$x\in \spec A$ d'un morphisme d'évaluation
$f\mapsto f(x)$, à valeurs dans le corps~$\kappa(x)$ {\em qui dépend
{\em a priori}
de~$x$.}

\trois{invers-speca}
{\em Inversibilité : tout se passe bien.}
Soit~$f\in A$. L'élément~$f$ appartient à~$A\ti$
 si et seulement si il n'appartient à aucun idéal 
 premier de~$A$ ; autrement dit, $f$ est inversible
 si et seulement si~$f(x)\neq 0$ pour tout~$x\in \spec A$ : en ce
 qui concerne l'inversibilité, $A$ se comporte effectivement
 comme un anneau de fonctions classiques sur~$\spec A$. 
 
 \trois{null-speca}
 {\em Annulation en tout point : les limites du point de vue fonctionnel.} 
 Soit~$f\in A$. On a~$f(x)=0$ pour tout~$x\in \spec A$ si et seulement 
 si~$f$ appartient à tous les idéaux premiers d~$A$, c'est-à-dire si et seulement si~$f$
 est nilpotent (lemme~\ref{appl-nilpo}). 
 
 Ainsi, lorsque~$A$ n'est pas réduit, $f$ peut s'annuler en tout point
 sans être elle-même nulle, et il est donc abusif de qualifier
 les éléments de~$A$ de fonctions ; on le fait tout de même
 parfois en pratique, soit parce qu'on travaille avec des anneaux réduits, soit
 pour le confort de l'analogie -- mais il faut garder en tête le problème des nilpotents ! 
 
 \subsection*{La topologie de Zariski}
 
 Nous allons maintenant définir une topologie sur~$\spec A$, et établir ses 
 propriétés de base, avant d'en 
 venir aux premiers exemples de spectres. 
 
 \deux{def-ve}
 Soit~$E$ une partie de~$A$. On note~$V(E)$ l'ensemble 
 des points~$x$ de~$\spec A$ tels que~$f(x)=0$ pour tout~$f\in E$.
 
 \deux{inter-union-zar}
 Si~$(E_i)_{i\in I}$ est une famille de sous-ensembles de~$A$, on voit immédiatement
 que~$\bigcap V(E_i)=V(\bigcup E_i)$. Supposons maintenant que~$I$ est fini, et soit~$F$ l'ensemble
 des éléments de~$A$ de la forme~$\prod_{i\in I}f_i$ où~$f_i\in E_i$ pour tout~$i$ ; on vérifie là encore
 sans problèmes que~$V(F)=\bigcup V(E_i)$.

\deux{sorites-ve}
Soit~$E\in A$. Les faits suivants découlent sans difficulté des définitions. 
 
 \trois{ve-vie}
 Si~$I$ désigne l'idéal engendré par~$E$ alors~$V(E)=V(I)$. 
 
 \trois{zar-premier}
 Soit~$x\in \spec A$ et soit~$\got p$ l'idéal premier correspondant. Le 
 point~$x$ appartient à~$V(E)$ si et seulement si~$E\subset \got p$. 
 
 \deux{top-zar} Il résulte de~\ref{inter-union-zar}
 que les parties de la forme~$V(E)$ pour~$E\subset A$ sont les fermés d'une topologie sur~$\spec A$,
 dite {\em de Zariski}. Par construction, les ouverts de~$\spec A$ sont les parties qui sont 
 réunion de sous-ensembles de la forme
 $$D(f):=\{x\in A, f(x)\neq 0\}$$ où~$f\in A$. Notons que~$D(fg)=D(f)\cap D(g)$ pour tout~$(f,g)\in A^2$,
 et que si~$x$ est un point de~$\spec A$ correspondant à un idéal premier~$\got p$ alors
 $x\in D(f)$ si et seulement si~$f\notin \got p$. 
 
 \deux{pas-separe}
 La topologie de Zariski ne ressemble guère aux topologies usuelles. Par exemple, en général
 les points de~$\spec A$ ne sont pas tous fermés (et~$\spec A$ n'est
 {\em a fortiori}
 pas séparé). 
 
 \medskip
 Plus précisément, soit~$x\in \spec A$
 et soit~$\got p$ l'idéal premier correspondant. 
 Il découle tautologiquement
 de la définition de la topologie de Zariski
 que~$\overline{\{x\}}$ est l'ensemble des points~$y$ tels que
 $f(y)=0$ pour toute~$f\in A$ s'annulant en~$x$. Autrement dit, 
 ~$\overline{\{x\}}=V(\got p)$. Cela signifie que
 si~$y$ est un point de~$\spec A$
 correspondant à un idéal premier~$\got q$, le point~$y$
 appartient à~$\overline{\{x\}}$ si et seulement si~$\got p\subset \got q$. 
 En particulier,~{\em le point~$x$ est fermé si et seulement si~$\got p$ est maximal}
 (et on a alors~$x=V(\got p)$). 
 
 \deux{def-qc}
{\bf Définition.}
On dit qu'un espace topologique~$X$
est {\em quasi-compact}
si de tout recouvrement ouvert de~$X$ on peut extraire un sous-recouvrement fini. 

\medskip
Attention : la différence entre «quasi-compact» et~«compact» est que l'on ne requiert pas
qu'un espace quasi-compact soit séparé. 

\deux{speca-qc}
{\bf Lemme.}
{\em L'espace topologique~$\spec A$ est quasi-compact.}

\medskip
{\em Démonstration.}
Soit~$(U_i)$ un recouvrement ouvert de~$\spec A$. Pour montrer qu'on 
peut en extraire un sous-recouvrement fini, on peut toujours le raffiner, 
ce qui autorise à supposer~$U_i$ de la forme~$D(f_i)$ pour tout~$i$. 

\medskip
Dire que~$\spec A=\bigcup D(f_i)$ signifie que pour tout~$x\in \spec A$,
il existe~$i$ tel que~$f_i(x)\neq 0$. En termes d'idéaux premiers, cela se traduit comme suit : 
pour tout idéal premier~$\got p$ de~$A$, il existe~$i$ tel que~$f_i\notin \got p$. Cela revient
à dire que l'idéal~$(f_i)$ n'est contenu dans aucun idéal premier de~$A$, c'est-à-dire encore
que~$(f_i)=A$, ou que~$1\in (f_i)$. Mais cette dernière condition équivaut à demander que~$1$ s'écrive comme
une combinaison {\em finie}
$\sum a_i f_i$ ; si~$J$ désigne l'ensemble des indices apparaissant effectivement dans cette écriture, 
on a~$1\in (f_i)_{i\in J}$ et partant~$\spec A=\bigcup_{i\in J}D(f_i)$
(on remonte la chaîne d'équivalences que l'on a mise en évidence).~$\Box$ 

\deux{comment-qc}
{\em Commentaires.}
La quasi-compacité est une propriété de finitude qui
{\em n'est pas}
un avatar schématique raisonnable de la compacité. Un tel avatar existe, 
c'est la
{\em propreté}
que nous rencontrerons plus bas,
et qui ne peut pas
être définie en termes purement 
topologiques. 

\subsection*{Premiers exemples}

\deux{spec-vide}
Si~$A$ est un anneau $\spec A=\varnothing$ si et seulement si~$A$
n'a pas d'idéal premier, c'est-à-dire si et seulement si~$A=\{0\}$. Dans le cas
contraire, $A$ admet un idéal maximal, et~$\spec A$ possède donc
au moins un point fermé (\ref{pas-separe}).

\deux{spec-corps}
Soit~$k$ un corps. Il possède un unique idéal premier, à savoir~$\{0\}$. 
Son spectre est donc un singleton~$\{x\}$, et~$\kappa(x)={\rm Frac}\;(k/\{0\})=k$ ; 
l'évaluation en~$x$ est bien entendu l'identité de~$k$. 

\deux{spec-z}
{\bf Description de~$\spec \ZZ$.}
Donnons la liste des idéaux  premiers de~$\ZZ$. 

\trois{id-max-z}
{\em Les idéaux maximaux.}
Ils sont de la forme~$(p)$ avec~$p$ premier.
À tout nombre premier~$p$ est donc
associé un point fermé~$x_p$ de~$\spec \ZZ$, qui s'identifie à~$V(p)$ : c'est le lieu 
d'annulation de~$p$ {\em vu comme fonction sur~$\spec \ZZ$.}

\medskip
On a~$\kappa(x_p)=\FF_p$, et~$f\mapsto f(x_p)$ est simplement
la réduction modulo~$p$. 

\trois{point-gen-specz}
{\em L'idéal~$(0)$.}
Soit~$\eta$ le point correspondant de~$\spec \ZZ$. On a
$$\overline{\{\eta\}}=V(0)=\spec \ZZ.$$ 
Le point~$\eta$ est donc {\em dense}
dans~$\spec \ZZ$ ; on dit aussi qu'il est {\em générique.}

\medskip
On a~$\kappa(\eta)=\QQ$, et~$f\mapsto f(\eta)$ est simplement
l'injection canonique~$\ZZ\hookrightarrow \QQ$. 

\deux{top-specz}
{\bf Les fermés de~$\spec \ZZ$.}
Soit~$F$ un fermé de~$\spec \ZZ$. Il est de la forme~$V(I)$ pour un certain 
idéal~$I$, que l'on peut écrire~$(a)$ avec~$a\in \NN$ (puisque~$\ZZ$ est principal) ; 
on a~$F=V(a)$. 

\trois{cas-a-nul}
Si~$a=0$ alors~$F=\spec \ZZ$. 

\trois{cas-a-pasnul}
Si~$a\geq 1$, on peut écrire~$a=\prod p_i ^{n_i}$ où
les~$p_i$ sont des nombres premiers deux à deux distincts et les~$n_i$
des entiers strictement positifs. On a alors~$F=\{x_{p_1}, \ldots, x_{p_n}\}$. 

\trois{conclu-fermez}
En conclusion les fermés de~$\spec \ZZ$ sont d'une part~$\spec \ZZ$ lui-même,
d'autre part les ensembles finis de points fermés.

\deux{spec-kt}
{\bf Description de~$\spec k[T]$.}
Soit~$k$ un corps. 
Donnons la liste des idéaux  premiers de~$k[T]$. 

\trois{id-max-kt}
{\em Les idéaux maximaux.}
Ils sont de la forme~$(P)$ avec~$P$ irréductible
unitaire. 
À tout polynôme
irréductible unitaire~$P\in k[T]$  est donc
associé un point fermé~$x_P$ de~$\spec k[T]$, qui s'identifie à~$V(P)$ : c'est le lieu 
d'annulation de~$P$ {\em vu comme fonction sur~$\spec k[T]$.}

\medskip
On a~$\kappa(x_P)=k[T]/(P)$, et~$f\mapsto f(x_P)$ est simplement
la réduction modulo~$P$. 

\trois{point-gen-speckt}
{\em L'idéal~$(0)$.}
Soit~$\eta$ le point correspondant de~$\spec k[T]$. On a
$$\overline{\{\eta\}}=V(0)=\spec k[T].$$ 
Le point~$\eta$ est donc {\em dense}
dans~$\spec k[T]$ ; on dit aussi qu'il est {\em générique.}

\medskip
On a~$\kappa(\eta)=k(T)$, et~$f\mapsto f(\eta)$ est 
l'injection canonique~$k[T]\hookrightarrow k(T)$.

\trois{points-class-it}
{\em Les points «naïfs».}
Le spectre de~$k[T]$ est la variante schématique de la droite affine. Un point de celle-ci 
n'est autre qu'un élément~$\lambda$ de~$k$. Or si~$\lambda \in k$, le point naïf correspondant
peut être vu comme un point fermé de~$k[T]$ : avec les notations de~\ref{id-max-kt},
c'est simplement le point~$x_{T-\lambda}$. En effet, celui-ci
est précisément le lieu d'annulation de~$T-\lambda$, 
et l'on dispose d'un isomorphisme $$\kappa(x_{T-\lambda})=k[T]/(T-\lambda)\simeq k$$
modulo lequel l'évaluation~$f\mapsto f(x_{T-\lambda})$ est
simplement l'évaluation classique
$f\mapsto f(\lambda)$. 

\medskip
Lorsque~$k$ est algébriquement clos, tout polynôme irréductible
de~$k[T]$ 
est de la forme~$T-\lambda$, et tous les points fermés de~$\spec k[T]$ sont donc 
des points naïfs : hormis le point générique, $\spec k[T]$ n'est constitué que de «vrais»
points. 

\deux{top-speckt}
{\bf Les fermés de~$\spec k[T]$.}
Soit~$F$ un fermé de~$\spec k[T]$. Il est de la forme~$V(I)$ pour un certain 
idéal~$I$, que l'on peut écrire~$(Q)$ où~$Q$ 
est un élément de~$k[T]$ nul ou unitaire (puisque~$k[T]$ est principal) ; 
on a~$F=V(Q)$. 

\trois{cas-q-nul}
Si~$Q=0$ alors~$F=\spec k[T]$. 

\trois{cas-Q-pasnul}
Si~$Q\neq 0$, on peut écrire~$Q=\lambda \prod_{i=1}^n P_i ^{n_i}$ où
les~$P_i$ sont des polynômes irréductibles unitaire 
deux à deux distincts et où les~$n_i$
sont 
des entiers strictement positifs. On a alors~$F=\{x_{P_1}, \ldots, x_{P_n}\}$. 

\trois{conclu-fermekt}
En conclusion les fermés de~$\spec k[T] $ sont d'une part~$\spec k[T]$ lui-même,
d'autre part les ensembles finis de points fermés.

\deux{t2+1}
{\bf Un exemple de point fermé non naïf.}
Le polynôme~$T^2+1$
de~$\RR[T]$ étant irréductible, il définit un point
fermé~$x_{T^2+1}$
de~$\spec \RR[T]$,
qui est le lieu d'annulation de~$T^2+1$. 
Son corps résiduel~$\RR[T]/(T^2+1)$ est isomorphe
à~$\CC$ (comme extension de~$\RR$) 
de deux manières différentes : on peut envoyer~$T$ sur~$i$ ou~$(-i)$ ; le
morphisme d'évaluation~$f\mapsto f(x_{T^2+1})$ s'identifie à l'évaluation 
classique~$f\mapsto f(i)$ dans le premier cas, et à~$f\mapsto f(-i)$ dans le second cas. 

\medskip
On voit que les points complexes naïfs~$i$ et~$(-i)$ de la droite affine induisent
le même point fermé
de~$\spec \RR[T]$ : cela traduit le fait que~$i$ et~$(-i)$ sont
en quelque sorte
{\em $\RR$-indiscernables.}

\deux{cas-k-alg}
{\bf Spectre d'une~$k$-algèbre de type fini.}
Soit~$k$ un corps, soit~$A$ une~$k$-algèbre de type fini et soit~$X$
son spectre. Si~$x\in X$, l'évaluation
$f\mapsto f(x)$ induit un plongement~$k\hookrightarrow \kappa(x)$
qui fait de~$\kappa(x)$ une extension de~$k$.

\trois{def-xl}
Pour toute extension~$L$ de~$k$ on note~$X(L)$ l'ensemble~$\hom_k(A, L)$. 
Cela peut paraître abusif, puisque~$X$ semble dès lors
désigner à la fois un schéma et un foncteur~$L\mapsto X(L)$, 
mais nous verrons plus loin qu'une telle notation est tout à fait justifiée ; nous nous permettons
donc de l'utiliser dès maintenant, car elle va être commode et n'induira aucune confusion. 

\medskip
Fixons une présentation~$A\simeq k[T_1,\ldots, T_n]/(P_1,\ldots,P_r)$ de~$A$. 
Elle induit une bijection 

$$X(L)\simeq \{(x_1,\ldots, x_n)\in L^n \;\;{\rm t.q.}\;\;\forall j\; \;P_j(x_1,\ldots, x_n)=0\}$$
qui est fonctorielle en~$L$
(à un~$n$-uplet~$(x_1,\ldots, x_n)$ correspond le morphisme
de~$A$ vers~$L$ induit par l'évaluation
des polynômes en~$(x_1,\ldots,x_n)$). 

Ainsi, on peut voir~$X(L)$
comme «l'ensemble des~$L$-points
de la variété algébrique d'équations~$P_1=0,\ldots, P_r=0$», dont~$X=\spec A$ est censé être la déclinaison
schématique. 

\medskip
Pour toute extension~$L$ de~$k$, on dispose d'une application naturelle de~$X(L)$ vers~$X$
(\ref{prem-def-spec}, \ref{spec-def}). 

\trois{pts-fermes-a}
Soit~$X_0$ l'ensemble des points fermés de~$X$ ; il est non vide dès que~$A$ est non nulle
(\ref{spec-vide}). On déduit du
{\em Nullstellensatz}, et plus précisément de sa variante donnée par l'énoncé~\ref{nullst-2}, que~$X_0$ est l'ensemble des points~$x\in X$ tels que~$\kappa(x)$ soit fini sur~$k$. 

\trois{kbarre-x0}
On fixe une clôture algébrique~$\bar k$ de~$k$, et on désigne
par~$G$ le groupe de Galois de~$\bar k/k$. 
Soit~$x\in X_0$. Comme~$\kappa(x)$ est une extension finie de~$k$, il admet 
un~$k$-plongement dans~$\bar k$. La composée~$A\to \kappa(x)\hookrightarrow \bar k$ 
est un élément de~$X(\bar k)$ dont l'image sur~$X$
est par construction égale à~$x$. 

\medskip
Réciproquement, donnons-nous un~$k$-morphisme~$A\to \bar k$. Comme~$A$ est de type fini, 
son image est engendrée par un nombre fini d'éléments, et est donc une extension finie~$L$ de~$k$. 
Le morphisme~$A\to \bar k$
se factorisant par la flèche surjective~$A\to L$, son image~$x$ sur~$X$ appartient à~$X_0$
et vérifie~$\kappa(x)=L$. 

\medskip
Ainsi, la flèche canonique~$X(\bar k)\to X$ a pour image~$X_0$. 

\trois{action-g}
Le but est maintenant de décrire le «noyau»
de cette flèche, ou plus précisément son défaut d'injectivité. 

\medskip
Pour commencer, remarquons qu'il y a une action naturelle 
de~$G$ sur~$X(\bar k)$ : si~$g\in G$ et si~$\phi : A\to \bar k$ est un élément
de~$X(\bar k)$, on pose~$g\cdot \phi= g\circ \phi$. Si l'on identifie~$X(\bar k)$
au sous-ensemble de~$\bar k^n$ formé des~$n$-uplets en lesquels tous les~$P_j$
s'annulent, le lecteur vérifiera aisément 
que~$g\cdot(x_1,\ldots, x_n)=(g(x_1),\ldots, g(x_n))$ pour tout~$g\in G$ et 
tout~$n$-uplet~$(x_1,\ldots, x_n)$ en lequel les~$P_j$ s'annulent
(qu'ils s'annulent aussi en le~$n$-uplet $(g(x_1),\ldots, g(x_n))$
résulte du fait qu'ils sont à coefficients dans~$k$). 

\medskip
Nous allons maintenant démontrer que~$X(\bar k)\to X_0$ identifie~$X_0$
au quotient~$X(\bar k)/G$. On retrouve le phénomène
entrevu au~\ref{t2+1} : deux points de~$X(\bar k)$ induisent le même 
point de~$X$ si et seulement si ils sont conjugués sous~$G$,
c'est-à-dire en un sens~$k$-indiscernables. 

\medskip
{\em Deux éléments de~$X(\bar k)$ conjugués sous l'action de~$G$
ont même image sur~$X_0$.}
En effet, donnons-nous~$\phi : A\to \bar k$ et~$g\in G$. L'existence du diagramme
commutatif
$$\xymatrix{
& {\bar k}\\
A\ar[ru]^{g\circ \phi}\ar[r]^\phi&{\bar k}\ar[u]_\simeq^g}$$
assure que les images de~$\phi$ et~$g\circ \phi$ sur~$X$
coïncident, ce qu'on souhaitait établir. 

\medskip
{\em Deux éléments de~$X(\bar k)$ 
qui ont même image sur~$X_0$ sont conjugués sous l'action de~$G$.}
En effet, donnons-nous deux morphismes~$\phi$
et~$\psi$ de~$A$ vers~$\bar k$ qui ont même image~$x$ sur~$X$. 
Cela signifie que~$\phi$ et~$\psi$ se factorisent tous deux par la surjection 
canonique~$A\to \kappa(x)$ ; autrement dit, $\phi$ est induit
par un~$k$-plongement~$\phi' : \kappa(x)\hookrightarrow \bar k$, 
et~$\psi$ par un~$k$-plongement~$\psi' : \kappa(x)\hookrightarrow \bar k$. 
Chacun de ces deux plongements fait de~$\bar k$ une clôture algébrique de~$\kappa(x)$. 
Comme deux telles clôtures algébriques sont isomorphes, 
il existe un automorphisme~$g$ de~$\bar k$ tel que~$g\circ \phi'=\psi'$ ; un tel~$g$ 
est automatiquement un~$k$-morphisme (car~$\phi'$ et~$\psi'$ sont des~$k$-morphismes), 
ce qui veut dire qu'il appartient à~$G$. On a par construction~$\psi=g\circ \phi$, ce qui achève la
démonstration. 

\trois{cas-part-kpoints}
Comme~$G$ fixe~$k$, il agit trivialement
sur le sous-ensemble~$X(k)$ de~$X(\bar k)$. On déduit alors
de~\ref{action-g}
que~$X(k)\to X_0$ est injectif, et bijectif si~$k$ est algébriquement clos. 

Ainsi, ce qu'on avait remarqué au~\ref{spec-kt}
dans un cas particulier vaut en général : l'ensemble~$X(k)$
des points «naïfs» de la~$k$-variété
algébrique d'équations~$P_1=0,\ldots, P_r=0$
se plonge dans~$X_0$ et s'identifie à
celui-ci lorsque~$k$ est algébriquement clos. 
Sous cette dernière hypothèse, on peut également
donner, dans une certaine mesure, une interprétation
classique des points «non naïfs»
de~$X$, {\em cf.} \ref{points-schem-interp}
{\em et sq. infra}.

\trois{xk-corpsk}
L'image de~$X(k)$ dans~$X_0$ est 
exactement
l'ensemble des points~$x$
tels que~$\kappa(x)=k$. En effet, si~$\phi : A\to k$ est un 
élément de~$X(k)$ il est nécessairement surjectif (considérer
les constantes), et le point~$x$ qu'il induit a donc pour corps résiduel~$k$. 

Réciproquement, si~$x\in X$ est tel que~$\kappa(x)=k$,
l'évaluation~$f\mapsto f(x)$ est un morphisme de~$A$
dans~$k$ qui induit~$x$, d'où notre assertion. 

\trois{desc-x0-naif}
Soit~$x\in X$ un point de corps résiduel~$k$. Par ce qui précède,
il provient d'un unique point de~$X(k)$, point que l'on voit 
comme un~$n$-uplet $(x_1,\ldots, x_n)$ en lequel les~$P_j$ s'annulent. 

L'évaluation~$A\to \kappa(x)=k$ est alors par construction
induite par l'évaluation classique des polynômes en~$(x_1,\ldots, x_n)$. 
L'idéal maximal associé à~$x$ est le noyau de~$A\to \kappa(x)$ ; nous invitons
le lecteur à vérifier qu'il est engendré par les~$\overline{T_i}-x_i$. 
En particulier, $\{x\}=V(\overline{T_1}-x_1,\ldots, \overline{T_n}-x_n)$.

\subsection*{Fonctorialité du spectre}

\deux{spec-fonct}
Soit~$\phi : A\to B$ un morphisme d'anneaux. On a mentionné
en~\ref{spec-contrav}
fait que~$\phi$ induit une application
$\psi$
de~$\spec B$
vers~$\spec A$, que l'on peut décrire de deux façons : 

\medskip
$\bullet$ au niveau des idéaux premiers, elle envoie~$\got q$ sur~$\phi^{-1}(\got q)$ ; 

$\bullet$ au niveau des morphismes dont le but est un corps, elle envoie la classe de~$B\to K$
vers celle de la flèche composée~$A\to B\to K$. 

\medskip
Il résulte immédiatement de l'une ou l'autre de ces définitions que si~$y$
est un point de~$\spec B$ d'image~$x$
sur~$\spec A$, le morphisme~$\phi$ induit un plongement
$\kappa(x)\hookrightarrow \kappa(y)$ tel que le diagramme

$$\xymatrix{
B\ar[rr]^{g\mapsto g(y)}&&{\kappa(y)}\\
A\ar[rr]_{f\mapsto f(x)}\ar[u]^\phi&&{\kappa(x)}
\ar@{^{(}->}[u]
}$$
commute. En particulier, on 
a pour tout~$f\in A$ l'équivalence
$$f(x)=0\iff\phi(f)(y)=0.$$

\medskip
En conséquence, on a~$\psi^{-1}(V(E))=V(\phi(E)))$
pour toute partie~$E\subset A$, et 
$\psi^{-1}(D(f))=D(\phi(f))$ pour tout~$f\in A$. 
Il s'ensuit que l'application~$\psi$
est
{\em continue}. 

\medskip
Ainsi, $A\mapsto \spec A$ apparaît comme un foncteur contravariant de la catégorie
des anneaux vers celle des espaces topologiques. 

\deux{ex-corps-morspec}
{\bf Un premier exemple.}
Soit~$A\to K$ un morphisme d'un anneau~$A$ vers un corps~$K$, soit~$x$ le point
correspondant de~$\spec A$ et soit~$\xi$ l'unique point de~$\spec K$. La flèche~$A\to K$
induit une application~$\spec K \to \spec A$. L'image de~$\xi$ est simplement 
par définition la classe du morphisme~$A\to K$, c'est-à-dire~$x$ ; et la flèche 
$\kappa(x)\hookrightarrow \kappa(\xi)=K$ est l'injection canonique de~$\kappa(\xi)$
dans~$K$. 

\deux{loc-quot}
Soit~$\phi \colon A\to B$ un morphisme d'anneaux. Les assertions suivantes
sont équivalentes :

\medskip
i) il existe une partie multiplicative~$S$ de~$A$ et un idéal~$I$ de~$A$ tels que la~$A$-algèbre~$B$
s'identifie à~$S^{-1}A/(I\cdot S^{-1}A)$. 

ii) tout élément de~$B$ est de la forme~$\phi(a)/\phi(s)$ où~$a$ et~$s$ appartiennent à~$A$ et où~$\phi(s)$
est inversible dans~$B$.

\medskip
En effet, i)$\Rightarrow$ii) est évidente ; pour~ii)$\Rightarrow$i), nous laissons au lecteur le soin de vérifier
qu'on peut prendre~$S=\phi^{-1}(B\ti)$ et~$I={\rm Ker}\;\phi$. 

\deux{spec-loc-quot}
Soit~$\phi : A\to B$ un morphisme d'anneaux vérifiant les conditions équivalentes
i) et~ii) du~\ref{loc-quot}
ci-dessus, et soient~$S$ et~$I$ comme dans~i). Soit~$\psi \colon \spec B \to \spec A$
la flèche induite par~$\phi$. Nous allons démontrer que~$\psi$ induit un homéomorphisme
$$\spec B\simeq \{x\in \spec A, x\in V(I)\;{\mathrm et}\;s(x)\neq 0\;\;\forall s\in S\},$$
et que~$\kappa(\psi(y))=\kappa(y)$ pour tout~$y\in \spec B$. 

\trois{prop-univ-loc-quot}
On vérifie immédiatement que~$B$ représente le foncteur covariant
de~$\ann$ dans~$\ens$ qui envoie~$C$ sur le sous-ensemble de~$\hom_{\ann}(A,C)$
constitué des morphisme~$f$ tels que $f(a)=0$ pour tout~$a\in I$
et tels que~$f(s)$ soit inversible pour tout~$s\in S$. 

\trois{inj-loc-quot}
{\em La flèche~$\psi : \spec B\to \spec A$ induit
un homéomorphisme de~$\spec B$ sur~$\psi(\spec B)$.}
 En effet, soit~$b\in B$ ; écrivons~$b=\phi(a)/\phi(s)$ avec~$a\in A, s\in A$ et~$\phi(s)\in B\ti$. 
Soit~$y\in \spec B$ et soit~$x$ son image sur~$\spec A$. On a
les équivalences

$$b(y)=0$$
$$\iff (\phi(a)/\phi(s))(y)= 0$$
$$\iff \phi(a)(y)= 0$$
$$\iff a(x)=0.$$

On en déduit que le noyau de~$b\mapsto b(y)$, qui caractérise entièrement le point~$y$, 
ne dépend que de~$x$ ; autrement dit, $y$ est le seul antécédent de~$x$ sur~$\spec B$, et~$\psi$
est injective. 

Par ailleurs, la chaîne d'équivalence ci-dessus entraîne que~$y\in D(b)$ si et seulement
si~$x\in D(a)$. En conséquence, $\psi(D(b))=D(a)\cap \psi(\spec B)$, et l'injection continue
$\spec B\to \psi(\spec B)$ induite par~$\psi$ est dès lors ouverte ; c'est donc un homéomorphisme, 
ce qu'il fallait démontrer. 

\trois{im-loco-quot}
{\em Description de~$\psi(\spec B)$.}
La propriété universelle énoncée au~\ref{prop-univ-loc-quot}
signifie qu'un morphisme d'anneaux~$f\colon A\to C$ se factorise
{\em via}
$\phi\colon A \to B$ si et seulement si~$f(a)=0$ pour tout~$a\in I$
et~$f(s)\in C\ti$ pour tout~$s\in S$ (et qu'une telle factorisation, si elle existe,
est unique). En appliquant cette assertion dans le cas où~$C$ est un corps, et en la retraduisant
dans le langage des spectres, on voit qu'un point~$x$ de~$\spec A$ appartient
à~$\psi(\spec B)$ si et seulement si~$x\in V(I)$ et~$s(x)\neq 0$ pour tout~$s\in S$. 

\trois{corpsresi-loc-quot}
{\em La flèche~$\psi$ préserve
le corps résiduel des points.} 
Commençons par une remarque d'ordre général. Soit~$C\to K$ un morphisme d'un anneau~$C$ vers
un corps~$K$, soit~$\got p$ son noyau et soit~$x$ le point correspondant de~$\spec C$. Comme
$\kappa(x)={\rm Frac}\;C/\got p$, c'est le plus petit sous-corps de~$K$ contenant l'image de~$C$. 

\medskip
Soit maintenant~$y\in \spec B$ et soit~$x$ son image
sur~$\spec A$. Pour montrer que~$\kappa(y)=\kappa(x)$, il suffit de s'assurer que 
le plus petit sous-corps~$F$ de~$\kappa(y)$ contenant l'image de~$A$ par la flèche
composée~$A\to B\to \kappa(y)$ est égal à~$\kappa(y)$ lui-même. 

\medskip
Par définition, $F$ contient tous les éléments de la forme~$\phi(a)(y)$ pour~$a$
parcourant~$A$. Si~$s$ est un élément de~$A$
 tel que~$\phi(s)\in B\ti$, l'élément~$\phi(s)(y)$ de~$F$
 est non nul, et son inverse~$\phi(s)(y)^{-1}$
appartient à~$F$.
On voit donc que~$F$ contient tous les éléments de la forme~$(\phi(a)/\phi(s))(y)$ où~$a\in A$, 
$s\in A$ et où~$\phi(s)\in B\ti$ ; autrement dit, il contient tous les éléments de la forme~$b(y)$
où~$b\in B$, et coïncide de ce fait avec~$\kappa(y)$. 

\deux{exem-loc-quot}
{\bf Exemples}. Nous allons décliner~\ref{spec-loc-quot}
dans un certain nombre
de cas particuliers importants. Soit~$A$ 
un anneau. 

\trois{spec-asuri}
Soit~$I$ un idéal de~$A$. En appliquant 
\ref{spec-loc-quot}
avec~$S=\{1\}$, on voit que le morphisme
quotient~$A\to A/I$ induit un homéomorphisme préservant
les corps résiduels 
de~$\spec A/I$ sur le fermé~$V(I)$ de~$\spec A$. 

\medskip
Notons qu'on peut avoir~$V(I)=\spec A$ sans que
l'idéal~$I$ soit nul : 
comme~$V(I)=\bigcap_{f\in I}V(f)$, le fermé~$V(I)$ est 
plus précisément égal
à~$\spec A$ si et seulement si~$V(f)=\spec A$ pour tout~$f\in I$, 
c'est-à-dire si et seulement si~$I$ est constitué d'éléments nilpotents
(\ref{null-speca}).

Dans ce cas, la flèche~$A\to A/I$ induit
en vertu de ce qui précède
un homéomorphisme
$\spec A/I\simeq \spec A$. 

\trois{spec-af}
Soit~$f$ un élément de~$A$. Si~$x\in \spec A$, 
on a~$f(x)\neq 0$ si et seulement si~$f^n(x)\neq 0$ pour tout~$n$. 
En appliquant \ref{spec-loc-quot}
avec~$S=\{f^n\}_n$
et~$I=\{0\}$, on voit que le morphisme
de localisation~$A\to A_f$ induit 
un homéomorphisme préservant
les corps résiduels 
de~$\spec A_f$ sur l'ouvert~$D(f)$ de~$\spec A$. 

\trois{spec-ap}
Soit~$\got p$ un idéal premier de~$A$ et soit~$x$ le point correspondant
de~$\spec A$. En appliquant \ref{spec-loc-quot}
avec~$S=A\setminus \got p$ et~$I=\{0\}$, on voit que le morphisme
de localisation~$A\to A_{\got p}$ induit un homéomorphisme préservant
les corps résiduels 
de~$\spec A_{\got p}$ sur
$$\{y\in \spec A, f(y)\neq 0\;\;\forall f\notin \got p\}.$$
Ce dernier ensemble peut se décrire d'une manière un peu plus géométrique, en remarquant
que l'implication $(f\notin \got p)\Rightarrow (f(y)\neq 0)$ équivaut à sa contraposée
$(f(y)=0)\Rightarrow f\in \got p.$ En se rappelant que~$f\in \got p\iff f(x)=0$, 
et en utilisant la description de~$\overline{\{y\}}$ donnée au~\ref{pas-separe},
on en déduit que l'image de~$\spec A_{\got p}$ est précisément l'ensemble
des {\em générisations}
de~$x$, c'est-à-dire des éléments~$y$ tels que~$x\in \overline{\{y\}}$. 

\trois{spec-kappax}
Soit~$\got p$ un idéal premier de~$A$ et soit~$x$ le point correspondant
de~$\spec A$. En appliquant \ref{spec-loc-quot}
avec~$S=A\setminus \got p$ et~$I=\got p$, 
on voit que le morphisme
canonique~$A\to {\rm Frac}\; A/\got p$ induit un homéomorphisme préservant
les corps résiduels 
de~$\spec {\rm Frac}\;A/\got p$ sur
$$\{y\in \spec A, f(y)\neq 0\;\;\forall f\notin \got p\;\;{\rm et}\;f(y)=0\;\;\forall f\in \got p\},$$
qui n'est autre que l'ensemble des points~$y$ tels que le noyau de l'évaluation en~$y$ soit
exactement~$\got p$ ; c'est donc le singleton~$\{x\}$. 

\medskip
On retrouve ainsi l'exemple~\ref{ex-corps-morspec}
dans le cas particulier du morphisme canonique~$A\to \kappa(x)$. 

\deux{fibre-mor-spec}
{\bf Fibres d'une application entre spectres.}
Soit~$\phi \colon A \to B$ un morphisme d'anneaux et soit~$\psi \colon \spec B\to \spec A$
l'application continue induite. Soit~$x\in \spec A$, et soit~$\got p$ l'idéal premier correspondant. 
On se propose de donner une description de la fibre~$\psi^{-1}(x)$.

\medskip
Posons~$S=A\setminus \got p$. Le point~$x$ peut se décrire comme
$$\{y\in \spec A, f(y)\neq 0\;\;\forall f\notin \got p\;\;{\rm et}\;f(y)=0\;\;\forall f\in \got p\},$$
({\em cf.}~\ref{spec-kappax}
ci-dessus). En conséquence, 
$$\psi^{-1}(x)=\{z\in \spec B, \phi(f)(z)\neq 0\;\;\forall f\notin \got p\;\;{\rm et}\;\phi(f)(z)=0\;\;\forall f\in \got p\},$$
que l'on peut récrire
$$\{z\in \spec B, g(z)\neq 0\;\;\forall g\in \phi(S)\;\;{\rm et}\;g(z)=0\;\;\forall g\in \phi(\got p)\}.$$

Il résulte alors de~\ref{spec-loc-quot} que~$B\to \phi(S)^{-1}B/(\phi(\got p)\cdot \phi(S)^{-1}B)$
induit un homéomorphisme préservant les corps résiduels 
entre~$\spec B/(\phi(\got p)\cdot \phi(S)^{-1}B)$ et~$\psi^{-1}(x)$. Par ailleurs, la commutation du produit
tensoriel à la localisation et au quotient garantit que
$$\phi(S)^{-1}B/(\phi(\got p)\cdot \phi(S)^{-1}B)\simeq B\otimes_A S^{-1}A/(\got p\cdot S^{-1}A)$$
$$=B\otimes_A{\rm Frac}\;A/\got p=B\otimes_A \kappa(x).$$

\medskip
Récapitulons : on a finalement
montré que~$B\to B\otimes_A \kappa(x)$ induit un homéomorphisme 
$$\spec B\otimes_A \kappa(x) \simeq \psi^{-1}(x)$$ qui préserve les corps résiduels. 

\deux{spectre-prod}
{\bf Spectre d'un produit.}
Soient~$A$ et~$B$ deux anneaux. On note
respectivement~$e$ et~$f$ les idempotents~$(1,0)$ et~$(0,1)$
de~$A\times B$. 

\trois{decomp-prod}
On a~$ef=0$ et~$e+f=1$ ; on en déduit aussitôt
que~$\spec (A\times B)$ est la réunion disjointe
des ouverts fermés~$D(e)=V(f)$ et~$D(f)=V(e)$. 

\trois{prod-ann-ferm}
Le noyau de la projection~$A\times B\to A$ est égal à~$(f)$ ; en conséquence, 
cette projection induit un homéomorphisme préservant les corps
résiduels~$\spec A \simeq V(f)$ ; de même, la seconde 
projection induit un homéomorphisme préservant 
les corps résiduels~$\spec B\simeq V(e)$. 

\trois{prod-ann-ouv}
On peut décrire ces homéomorphismes 
d'une manière en quelque sorte duale
de la précédente. 
Pour cela, on vérifie
(l'exercice est laissé au lecteur) que
la projection~$A\times B\to A$ identifie~$A$ au localisé
$(A\times B)_e$ ; en conséquence, cette projection induit un homéomorphisme
préservant les corps
résiduels~$\spec A \simeq D(e)$. De même  la seconde 
projection induit un homéomorphisme préservant 
les corps résiduels~$\spec B\simeq D(f)$. 

\medskip
Bien entendu, tout ceci est compatible avec~\ref{prod-ann-ferm}
puisque~$D(e)=V(f)$ et~$V(f)=D(e)$. 

\trois{resume-prodann}
En résumé, on a donc un homéomorphisme
canonique
$$\spec(A\times B)\simeq \spec A\coprod \spec B$$
modulo lequel~$\spec A=D(e)=V(f)$ et~$\spec B=V(e)=D(f)$. 

\deux{ex-fibres-spec}
{\bf Un exemple.}
Soit~$\psi : \spec \CC[T]\to \spec \RR[T]$ le morphisme
induit par le plongement naturel~$\RR[T]\hookrightarrow \CC[T]$. Nous allons
étudier ses fibres. D'après~\ref{fibre-mor-spec},
la fibre de~$\psi$ en~$x$ s'identifie pour tout~$x\in \spec \RR[T]$ 
au spectre de l'anneau 
$\CC[T]\otimes_{\RR[T]}\kappa(x)$.

\trois{fibre-eta-rtct}
{\em La fibre générique.}
Soit~$\eta$ le point générique de~$\spec \RR[T]$. Par ce qui précède, la fibre~$\psi^{-1}(\eta)$ 
s'identifie au spectre de
$$ \CC[T]\otimes_{\RR[T]}\RR(T)=\RR[T][U]/(U^2+1)\otimes_{\RR[T]}\RR(T)
=\RR(T)[U]/(U^2+1)=\CC(T).$$
La fibre~$\psi^{-1}(\eta)$ contient donc un unique
point~$\xi$ dont le corps résiduel est~$\CC(T)$, et le morphisme d'évaluation~$f\mapsto f(\xi)$ n'est autre
que l'inclusion~$\CC[T]\hookrightarrow \CC(T)$ ; en conséquence,~$\xi$ est
le point générique de~$\spec \CC[T]$. 

\trois{fibre-r-rtrc}
{\em La fibre en un point naïf.}
Soit~$a\in \RR$. Il lui correspond un point fermé~$x$
de~$\spec \RR[T]$ de corps
résiduel~$\RR$, qui est précisément
le lieu d'annulation 
de~$T-a$ ; l'évaluation~$f\mapsto f(x)$ s'identifie
à l'évaluation classique~$f\mapsto f(a)$. 
La fibre~$\psi^{-1}(x)$ est donc le lieu d'annulation
de~$T-a$ sur~$\spec \CC[T]$, qui est également réduit à un point fermé~$y$
de corps résiduel~$\CC$, le morphisme~$f\mapsto f(y)$ 
correspondant à l'évaluation classique en~$a$. 
En quelque sorte, on peut dire que~$x$ est~«$a$ vu comme point réel»,
et que~$y$ est~«$a$ vu comme point complexe». 

\medskip
Vérifons la compatibilité avec la description «tensorielle»
de~$\psi^{-1}(x)$. Celle-ci assure 
que~$\psi^{-1}(x)$ s'identifie naturellement
au spectre de~$\CC[T]\otimes_{\RR[T]}\kappa(x)$,
soit encore à celui de~$\CC[T]\otimes_{\RR[T]}\RR[T]/(T-a)$, et finalement à~$\spec \CC[T]/(T-a)$. 
Comme~$\CC[T]/(T-a)$ est isomorphe à~$\CC$ par l'évaluation en~$a$, on retrouve ce qu'on attendait :
$\psi^{-1}(x)$ contient un unique point~$y$ de corps résiduel~$\CC$, et~$f\mapsto f(y)$ coïncide avec
l'évaluation classique en~$a$. 

\trois{fibre-c-rtrc}
{\em La fibre en un point fermé non naïf.}
Soit~$P$ un polynôme irréductible
de degré~$2$ sur~$\RR$. Il lui correspond un point fermé~$x$
qui est précisément
le lieu d'annulation 
de~$P$, et dont le corps résiduel~$\RR[T]/P$ s'identifie
à~$\CC$ une fois choisie (arbitrairement)
une racine~$\alpha$ de~$P$ ; l'évaluation~$f\mapsto f(x)$ s'identifie
alors à l'évaluation classique~$f\mapsto f(\alpha)$.

La fibre~$\psi^{-1}(x)$ est donc le lieu d'annulation
de~$P=(T-\alpha)(T-\bar \alpha)$ sur~$\spec \CC[T]$, et est dès
lors
égale
à~$V(T-\alpha)\cup V(T-\bar \alpha)$. Elle comprend
deux points fermés~$y$ et~$z$ de corps résiduel~$\CC$,
le morphisme~$f\mapsto f(y)$ 
correspondant à l'évaluation classique en~$\alpha$, et le morphisme~$f\mapsto f(z)$ 
à l'évaluation classique en~$\bar \alpha$.  

Intuitivement, $\alpha$ et~$\bar \alpha$ sont indiscernables sur~$\RR$ et ne définissent 
donc
qu'un seul point fermé sur~$\spec \RR[T]$ 
(un peu «gros» : son
corps résiduel est de degré~$2$ sur le corps de base~$\RR$) ; une fois les scalaires
étendus à~$\CC$ elles deviennent discernables et définissent
deux points fermés distincts et «naïfs» sur~$\spec \CC[T]$ : chacun a un corps résiduel égal au
nouveau corps de base~$\CC$. 

\medskip
Vérifions la compatibilité avec la description «tensorielle»
de~$\psi^{-1}(x)$. Celle-ci assure 
que~$\psi^{-1}(x)$ s'identifie naturellement
au spectre de~$\CC[T]\otimes_{\RR[T]}\kappa(x)$,
soit encore à celui de~$\CC[T]\otimes_{\RR[T]}\RR[T]/P$, et finalement à~$\spec \CC[T]/P$. 
Comme~$\CC[T]/P$ est isomorphe par
le lemme chinois à~$\CC\times \CC$ {\em via}
les évaluations en~$\alpha$ et~$\bar \alpha$,  on retrouve ce qu'on attendait :
en vertu de ce qui précède et de~\ref{resume-prodann},
$\psi^{-1}(x)$ contient deux points fermés~$y$ et~$z$ de corps résiduel~$\CC$,
le morphisme~$f\mapsto f(y)$ 
correspondant à l'évaluation classique en~$\alpha$, et le morphisme~$f\mapsto f(z)$ 
à l'évaluation classique en~$\bar \alpha$. 

\section{Description de~$\spec \ZZ[T]$ et~$\spec k[S,T]$ lorsque~$k$ est algébriquement clos}
\markboth{Le spectre comme espace topologique}{$\spec \ZZ[T]$  et $\spec k[S,T]$}

\subsection*{Le spectre de~$\ZZ[T]$} 

\deux{intro-speczt}
Le but de ce qui suit est de décrire le spectre de~$\ZZ[T]$ ; c'est un exemple qu'il est fondamental de bien comprendre
et méditer, car il offre un excellent échantillon des bizarreries et curiosités schématiques. 

Pour le décrire, nous allons recourir
à une stratégie très fréquente en géométrie : nous allons utiliser
une application de source~$\spec \ZZ[T]$
dont nous comprenons bien le but et les fibres. 

Cette application sera simplement la flèche $\psi\colon \spec \ZZ[T]\to \spec \ZZ$ 
induite par l'unique morphisme de~$\ZZ$ dans~$\ZZ[T]$. Son but est $\spec \ZZ$, que nous avons déjà décrit ;
et 
pour tout~$x\in \spec \ZZ$, la fibre~$\psi^{-1}(x)$ s'identifie au spectre de~$\ZZ[T]\otimes_{\ZZ}\kappa(x)$,
c'est-à-dire à~$\spec \kappa(x)[T]$, que nous avons décrit également.  

\deux{speczt-fibregen}
{\bf La fibre générique.}
Soit~$\eta$ le point générique de~$\spec \ZZ$. La fibre générique~$\psi^{-1}(\eta)$ 
s'identifie à~$\spec \QQ[T]$ et possède donc deux types de points. 

\trois{gen-gen-zt}
{\em Le point générique.}
Désignons par~$\xi_\eta$ le point générique de~$\spec \QQ[T]\simeq \psi^{-1}(\eta)$. 

\medskip
Lorsqu'on voit~$\xi_\eta$ comme un 
point de~$\spec \QQ[T]$, son corps résiduel est~$\QQ(T)$,
et l'évaluation
$f\mapsto f(\xi_\eta)$ est simplement le plongement~$\QQ[T]\hookrightarrow \QQ(T)$. 

\medskip
En conséquence, lorsqu'on voit~$\xi_\eta$ comme appartenant 
à~$\spec \ZZ[T]$, son corps résiduel est~$\QQ(T)$, et l'évaluation 
$f\mapsto f(\xi_\eta)$ est la flèche composée
$$\ZZ[T]\hookrightarrow \QQ[T]\hookrightarrow \QQ(T).$$
L'idéal premier correspondant à~$\xi_\eta$ est le noyau de cette dernière, c'est-à-dire
l'idéal nul. L'adhérence~$\overline{\{\xi_\eta\}}$ est alors égale à~$V(0)$, 
qui n'est autre que~$\spec \ZZ[T]$ tout entier. 

\trois{ferm-gen-zt}
{\em Les points fermés.}
Soit~$P$ un polynôme irréductible unitaire à coefficients dans~$\QQ[T]$. Il définit
un point fermé~$y_{\eta,P}$ sur~$\spec \QQ[T]\simeq \psi^{-1}(\eta)$. 

\medskip
Lorsqu'on voit~$y_{\eta,P}$ comme point de~$\spec \QQ[T]$, c'est le lieu des zéros de~$P$, son corps
résiduel est~$\QQ[T]/P$, et l'évaluation~$f\mapsto f(y_{\eta,P})$ est la flèche 
quotient~$\QQ[T]\to \QQ[T]/P$. 

\medskip
En conséquence, lorsqu'on voit~$y_{P,\eta}$ comme appartenant 
à~$\spec \ZZ[T]$, son corps résiduel est~$\QQ[T]/P$, et l'évaluation 
$f\mapsto f(y_{P,\eta})$ est la flèche composée
$$\ZZ[T]\hookrightarrow \QQ[T]\to \QQ[T]/P.$$
L'idéal premier correspondant à~$y_{P,\eta}$ est le noyau de cette dernière
flèche. Un raisonnement
fondé sur la factorialité de~$\ZZ$ (et que nous laissons au lecteur) assure que ce noyau
est l'idéal~$(P_0)$, où~$P_0$ est le produit de~$P$ par le plus petit multiple commun 
des dénominateurs de ses coefficients (écrits sous forme irréductible) ; c'est un polynôme
appartenant à~$\ZZ[T]$ dont le contenu (le plus grand diviseur commun des coefficients)
vaut~$1$. 

\medskip
L'adhérence~$\overline{\{y_{P,\eta}\}}$ est donc égale à~$V(P_0)$ ; nous en dirons
quelques mots un peu plus loin. 

\deux{speczt-fibreferm}
{\em Les fibres fermées.}
Soit
$p$ un nombre premier et soit~$x_p$ le
point fermé correspondant de~$\spec \ZZ$. 
La fibre fermée~$\psi^{-1}(x_p)$ s'identifie à~$\spec \FF_p[T]$
et possède donc deux types 
de points. 

\trois{gen-ferm-zt}
{\em Le point générique.}
Désignons par~$\xi_p$ le point générique de~$\spec \FF[T]\simeq \psi^{-1}(x_p)$.

\medskip
Lorsqu'on voit~$\xi_p$ 
comme un point de~$\spec \FF_p[T]$, son corps résiduel est~$\FF_p(T)$, et le
morphisme d'évaluation
$f\mapsto f(\xi_p)$ est simplement le plongement~$\FF_p[T]\hookrightarrow \FF_p(T)$. 

\medskip
En conséquence, lorsqu'on voit~$\xi_p$ comme appartenant 
à~$\spec \ZZ[T]$, son corps résiduel est~$\FF_p(T)$, et l'évaluation 
$f\mapsto f(\xi_p)$ est la flèche composée
$$\ZZ[T]\to \FF_p[T]\hookrightarrow \FF_p(T).$$
L'idéal premier correspondant à~$\xi_p$ est le noyau de cette dernière, 
c'est-à-dire~$(p)$. L'adhérence~$\overline{\{\xi_p\}}$ est alors égale à~$V(p)$, 
qui n'est autre que~$\psi^{-1}(x_p)$ (puisque~$x_p$ s'identifie
lui-même à~$V(p)\subset \spec \ZZ$). 

{\em Remarque.} On aurait
pu voir directement sans calculer le noyau
de l'évaluation que~$\overline{\{\xi_p\}}=\psi^{-1}(x_p)$,
puisque~$\xi_p$ est dense dans~$\psi^{-1}(x_p)$ et puisque cette dernière est fermée
dans~$\spec \ZZ[T]$. 

\trois{ferm-ferm-zt}
{\em Les points fermés.}
Soit~$P$ un polynôme irréductible unitaire à coefficients dans~$\FF_p[T]$. Il définit
un point fermé~$y_{p,P}$ sur~$\spec \FF_p[T]\simeq \psi^{-1}(x_p)$.  

\medskip
Lorsqu'on voit~$y_{p,P}$ comme un point de~$\spec \FF_p[T]$, 
c'est le lieu ses zéros de~$P$, son corps
résiduel est
le corps fini $\FF_p[T]/P$ (de cardinal~$p^{\deg P}$)
et l'évaluation en~$y_{p,P}$ est la flèche quotient~$\FF_p[T]\to \FF_p[T]/P$.  

\medskip
En conséquence, lorsqu'on voit~$y_{p,P}$ comme appartenant 
à~$\spec \ZZ[T]$, son corps résiduel est~$\FF_p[T]/P$. L'évaluation 
$f\mapsto f(y_{p,P})$ est la flèche composée
$$\ZZ[T]\to \FF_p[T]\to \FF_p[T]/P,$$ qui est surjective puisque composée
de deux surjections, et le point~$y_{p,P}$ est fermé (ce qui était d'ailleurs
évident
{\em a priori}, puisqu'il est fermé dans une fibre fermée). 

L'idéal maximal correspondant à~$y_{p,P}$
est le noyau de~$f\mapsto f(y_{p,P})$. On vérifie
aussitôt que si~$P^\sharp$ désigne un relevé quelconque de~$P$ dans~$\ZZ[T]$, 
ledit noyau est engendré par~$p$ et~$P^\sharp$ ; le point~$y_{p,P}$ est en conséquence
égal à~$V(p,P^\sharp)$. 

\medskip
Par exemple, $V(7, T-3)=y_{7, T-3}\in \psi^{-1}(x_7)$ ;
modulo l'identification~$\psi^{-1}(x_7)\simeq \spec \FF_7[T]$, le point~$y_{7,T-3}$ 
est le point de~$\spec \FF_7[T]$ défini
par l'annulation de~$T-3$, 
c'est à-dire le
point naïf correspondant à l'élément~$3$ de~$\FF_7$. 
Le corps résiduel de~$y_{7,T-3}$ est donc~$\FF_7$, et l'évaluation en~$y_{7,T-3}$ envoie un polynôme~$P\in \ZZ[T]$ 
sur la classe modulo~$7$ de~$P(3)$. 

\deux{adh-vp0}
{\bf Retour à l'étude de~$\overline{\{y_{\eta, P}\}}$, où~$P$ est un polynôme
irréductible de~$\QQ[T]$.}
Nous reprenons les notations~$P$ et~$P_0$ du~\ref{ferm-gen-zt}, et allons
décrire un peu plus précisément l'adhérence~$V(P_0)$ de~$y_{\eta, P}$, en regardant sa trace sur
chacune des fibres. Comme~$y_{\eta, P}$ est fermé dans~$\psi^{-1}(\eta)$, on a~$V(P_0)\cap \psi^{-1}(\eta)=\{y_{\eta, P}\}$. 

Écrivons~$P_0=\sum a_i T^i$. Soit~$p$ un nombre premier ; nous noterons~$a\mapsto \bar a$
la réduction modulo~$p$.  L'intersection de~$V(P_0)$ avec~$\psi^{-1}(x_p)\simeq \spec \FF_p[T]$ 
s'identifie au fermé~$V(\sum \overline{a_i} T^i)$ de~$\spec \FF_p[T]$. 

\trois{adh-vp0-fibre}
Par définition de~$P_0$, les~$a_i$ sont globalement premiers entre eux ; en particulier ils ne peuvent être
tous nuls modulo~$p$, et~$\sum \overline{a_i} T^i$ est donc un élément
{\em non nul}
de~$\FF_p[T]$ ; en conséquence, $V(P_0)\cap \psi^{-1}(x_p)$ est un ensemble
fini de points fermés. 

Cet ensemble peut être vide : c'est le cas si et seulement si~$\sum \overline{a_i} T^i$
est inversible, c'est-à-dire (compte-tenu du fait qu'il est non nul) si et seulement
si~$\overline{a_i}=0$ pour tout~$i\geq 1$, ou encore si et seulement si~$p$ divise~$a_i$
pour tout~$i\geq 1$. 

\trois{adhvp0-conclu}
Il résulte de ce qui précède
que l'ensemble des nombres premiers~$p$ tels que~$\overline{\{y_{\eta, P}\}}\cap \psi^{-1}(x_p)=\varnothing$
est {\em fini} ; autrement dit, $\overline{\{y_{\eta,P}\}}$ rencontre presque toutes les fibres fermées de~$\psi$. 

En particulier,~$\overline{\{y_{\eta, P}\}}$ n'est pas réduit au singleton~$\{y_{\eta, P}\}$ et le point~$y_{\eta, P}$ n'est
pas fermé. Il découle alors de toute l'étude menée ci-dessus
que les points fermés de~$\spec \ZZ[T]$ sont
exactement les points fermés de ses fibres fermées au-dessus de~$\spec \ZZ$, ou encore ses points
à corps résiduel fini. 

\deux{xP-caspart}
{\bf Étude de~$\overline{\{y_{\eta, P}\}}$ : deux exemples explicites.}

\trois{zt2+1}
{\em Le cas où~$P=T^2+1$.}
Dans ce cas~$P_0=T^2+1$ aussi. Soit~$p$ un nombre premier. Nous allons
décrire l'intersection de~$\overline{\{y_{\eta, T^2+1}\}}=V(T^2+1)$ avec~$\psi^{-1}(x_p)\simeq \spec \FF_p[T]$. 
Elle s'identifie à~$V(T^2+1)\subset \spec \FF_p[T]$. On distingue trois cas.

\medskip
$\bullet$ {\em Supposons que~$p=1\;\text{mod}\;4$.}
Le polynôme~$T^2+1$ est alors irréductible dans~$\FF_p[T]$.
En conséquence, $\overline{\{y_{\eta, T^2+1}\}}\cap \psi^{-1}(x_p)$
consiste en un point fermé dont le corps résiduel est~$\FF_p[T]/(T^2+1)$ (qui compte~$p^2$ éléments) ; 
en tant que point de~$\spec \ZZ[T]$, il est égal à~$V(p,T^2+1)$. 

$\bullet$ {\em Supposons que~$p=-1\;\text{mod}\;4$.}
Le polynôme~$T^2+1$ de~$\FF_p[T]$ 
s'écrit alors~$(T-a)(T+a)$ pour un certain~$a\in \FF_p\ti$. En conséquence, $\overline{\{y_{\eta, T^2+1}\}}\cap \psi^{-1}(x_p)$
est le sous-ensemble~$V(T-a)\cup V(T+a)$ de~$\spec \FF_p[T]$. 
Il consiste en deux points fermés de corps résiduel~$\FF_p$, à savoir les points naïfs
correspondant aux racines~$a$ et~$(-a)$. Si~$\alpha$ désigne n'importe quel
entier de classe modulo~$p$ égale à~$a$, les deux points en question, vus comme
appartenant à~$\spec \ZZ[T]$, sont respectivement égaux à~$V(p,T-\alpha)$
et~$V(p,T+\alpha)$. 

$\bullet$ {\em Supposons que~$p=2$.}
Le polynôme~$T^2+1$ de~$\FF_2[T]$ 
est égal à~$(T-1)^2$. En conséquence, $\overline{\{y_{\eta, T^2+1}\}}\cap \psi^{-1}(x_2)$
est le sous-ensemble~$V(T-1)$
de~$\spec \FF_2[T]$. 
Il consiste en un seul point fermé de corps résiduel~$\FF_2$ : le point naïf qui correspond
à~$1$. En tant point de~$\spec \ZZ[T]$, il est égal à~$V(2,T-1)$.

\trois{zt-demi}
{\em Le cas où~$P=T-(1/2)$.}
Dans ce cas~$P_0=2T-1$. Soit~$p$ un nombre premier. Nous allons
décrire l'intersection de~$\overline{\{y_{\eta, T-(1/2)}\}}=V(2T-1)$ avec~$\psi^{-1}(x_p)\simeq \spec \FF_p[T]$. 
Elle s'identifie à~$V(2T-1)\subset \spec \FF_p[T]$. On distingue deux cas.

\medskip
$\bullet$ {\em Supposons que~$p$ est impair.}
Dans ce cas~$(1/2)$ existe dans~$\FF_p$.
En conséquence,~$\overline{\{y_{\eta, T-(1/2)}\}}\cap \psi^{-1}(x_p)$
est le sous-ensemble~$V(T-(1/2))$
de~$\spec \FF_p[T]$. Il consiste en un seul point fermé de corps résiduel~$\FF_p$ : le point naïf qui correspond
à~$(1/2)$.  En tant point de~$\spec \ZZ[T]$, il est égal à~$V(p, T-\alpha)$ où~$\alpha$ est n'importe
quel entier tel que~$2\alpha$ soit congru à~$1$ modulo~$p$.  

$\bullet$ {\em Supposons que~$p=2$.}
Le polynôme~$2T-1$ de~$\FF_2[T]$
est alors égal à~$1$. En conséquence, $\overline{\{y_{\eta, T-(1/2))}\}}\cap \psi^{-1}(x_2)$
est le sous-ensemble~$V(1)$
de~$\spec \FF_2[T]$, qui est
{\em vide.}

\subsection*{Le spectre de~$k[S,T]$} 

\deux{intro-speckst}
On fixe un corps {\em algébriquement clos}~$k$.
Nous allons décrire dans ce qui suit le spectre de~$k[S,T]$, 
par une démarche parallèle à celle suivie
{\em supra}
à propos de~$\spec \ZZ[T]$. On note~$\psi\colon \spec k[S,T]\to \spec k[S]$ 
induite par le plongement~$k[S]\hookrightarrow k[S,T]$. Pour
tout point~$x$
de~$\spec k[S,T]$, la
fibre~$\psi^{-1}(x)$ s'identifie au spectre de~$k[S][T]\otimes_{k[S]}\kappa(x)$,
c'est-à-dire à~$\spec \kappa(x)[T]$. 

\medskip
Mentionnons avant d'entamer l'étude détaillée de~$\spec k[S,T]$ que celui-ci
est la variante schématique du plan affine, et que~$\psi \colon \spec k[S,T]\to \spec k[S]$
est quant à lui l'avatar de la première projection. 

\deux{speckst-fibregen}
{\bf La fibre générique.}
Soit~$\eta$ le point générique de~$\spec k[S]$. La fibre générique~$\psi^{-1}(\eta)$ 
s'identifie à~$\spec k(S)[T]$ et possède donc deux types de points.

\trois{gen-gen-kst}
{\em Le point générique.}
Désignons par~$\xi_\eta$ le point générique de~$\spec k(S)[T]\simeq \psi^{-1}(\eta)$. 

\medskip
Lorsqu'on voit~$\xi_\eta$ comme point de~$\spec k(S)[T]$, 
son corps résiduel est~$k(S,T)$, et le
morphisme d'évaluation
$f\mapsto f(\xi_\eta)$ est simplement le plongement~$k(S)[T]\hookrightarrow k(S,T)$. 

\medskip
En conséquence, lorsqu'on voit~$\xi_\eta$ comme appartenant 
à~$\spec k[S,T]$, son corps résiduel est~$k(S,T)$, et l'évaluation 
$f\mapsto f(\xi_\eta)$ est la flèche composée
$$k[S,T]\hookrightarrow k(S)[T]\hookrightarrow k(S,T).$$
L'idéal premier correspondant à~$\xi_\eta$ est le noyau de cette dernière, c'est-à-dire
l'idéal nul. L'adhérence~$\overline{\{\xi_\eta\}}$ est alors égale à~$V(0)$, 
qui n'est autre que~$\spec k[S,T]$ tout entier. 

\trois{ferm-gen-kst}
{\em Les points fermés.}
Soit~$P\in k(S)[T]$ un polynôme irréductible unitaire. Il définit
un point fermé~$y_{\eta,P}$ sur~$\spec k(S)[T]\simeq \psi^{-1}(\eta)$. 

\medskip
Lorsqu'on voit~$y_{\eta,P}$ comme un point
de~$\spec k(S)[T]$ c'est le lieu des zéros de~$P$, son corps
résiduel est~$k(S)[T]/P$ et l'évaluation en~$y_{\eta,P}$ est la 
flèche quotient~$k(S)[T]\to k(S)[T]/P$. 

\medskip
En conséquence, lorsqu'on voit~$y_{\eta,P}$ comme appartenant 
à~$\spec k[S,T]$, son corps résiduel est~$k(S)[T]/P$, et l'évaluation 
$f\mapsto f(y_{\eta,P})$ est la flèche composée
$$k[S,T]\hookrightarrow k(S)[T]\to k(S)[T]/P.$$
L'idéal premier correspondant à~$y_{\eta,P}$ est le noyau de cette dernière
flèche. Un raisonnement
fondé sur la factorialité de~$k[S]$ assure que ce noyau
est l'idéal~$(P_0)$, où~$P_0$ est le produit de~$P$ par le plus petit multiple commun 
des dénominateurs de ses coefficients (écrits sous forme irréductible) ; c'est un polynôme
appartenant à~$k[S][T]$ dont le contenu (le plus grand diviseur commun des coefficients)
vaut~$1$. 

\medskip
L'adhérence~$\overline{\{y_{\eta,P}\}}$ est donc égale à~$V(P_0)$ ; nous en dirons
quelques mots un peu plus loin. 

\deux{speckst-fibreferm}
{\em Les fibres fermées.}
Soit
$\lambda\in k$ et soit~$x_{S-\lambda}$ le
point fermé correspondant de~$\spec k[S]$.
La fibre fermée~$\psi^{-1}(x_{S-\lambda})$ s'identifie
au spectre de
l'anneau~$k[S,T]\otimes_{k[S]}\kappa(x_{S-\lambda})\simeq k[T]$ (ce dernier isomorphisme
étant celui qui envoie~$S$ sur~$\lambda$). 
Elle possède donc deux types de points. 

\trois{gen-ferm-kst}
{\em Le point générique.}
Désignons par~$\xi_{S-\lambda}$ le point générique de~$\spec k[T]\simeq \psi^{-1}(x_{S-\lambda})$. 

\medskip
Lorsqu'on voit~$\xi_{S-\lambda}$ comme un point de~$\spec k[T]$, son corps résiduel est~$k(T)$, et 
le morphisme d'évaluation
$f\mapsto f(\xi_{S-\lambda})$ est simplement le plongement~$k[T]\hookrightarrow k(T)$. 

\medskip
En conséquence, lorsqu'on voit~$\xi_{S-\lambda}$ comme appartenant 
à~$\spec k[S,T]$, son corps résiduel est~$k(T)$, et l'évaluation 
$f\mapsto f(\xi_{S-\lambda})$ est la flèche composée
$$\xymatrix{
{k[S,T]}\ar[rr]^{S\mapsto \lambda}&&{k[T]}\ar@{^{(}->}[r]&{ k(T)}}.$$
L'idéal premier correspondant à~$\xi_{S-\lambda}$ est le noyau de cette dernière, 
c'est-à-dire~$(S-\lambda)$. L'adhérence~$\overline{\{\xi_{S-\lambda}\}}$ est alors égale à~$V(S-\lambda)$, 
qui n'est autre que~$\psi^{-1}(x_{S-\lambda})$ (puisque~$x_{S-\lambda}$ s'identifie
lui-même à~$V(S-\lambda)\subset \spec k[S]$). 

{\em Remarque.} On aurait
pu voir directement sans calculer le noyau
de l'évaluation que~$\overline{\{\xi_{S-\lambda}\}}=\psi^{-1}(x_{S-\lambda})$,
puisque~$\xi_{S-\lambda}$ est dense dans~$\psi^{-1}(x_{S-\lambda})$ et puisque cette dernière est fermée
dans~$\spec k[S,T]$. 

\trois{ferm-ferm-kst}
{\em Les points fermés.}
Soit~$\mu$ un élément de~$k$. Il définit 
un point fermé naïf~$y_{S-\lambda, T-\mu}$ sur~$\spec k[T]\simeq \psi^{-1}(x_{S-\lambda})$.

\medskip
Si l'on voit~$y_{S-\lambda, T-\mu}$ comme un point de~$\spec k[T]$, c'est le lieu des zéros de~$T-\mu$, son corps
résiduel est~$k$, et l'évaluation en~$y_{S-\lambda, T-\mu}$ est l'évaluation
classique en~$\mu$. 

\medskip
En conséquence, lorsqu'on voit~$y_{S-\lambda, T-\mu}$ comme appartenant 
à~$\spec k[S,T]$, son corps résiduel est~$k$. L'évaluation 
$f\mapsto f(y_{S-\lambda, T-\mu})$ est la flèche composée
$$\xymatrix{
{k[S,T]}\ar[rr]^{S\mapsto \lambda}&&{k[T]}\ar[rr]^{T\mapsto \mu}&& k},$$ qui coïncide avec l'évaluation
classique~$f\mapsto f(\lambda, \mu)$ ; le point~$y_{S-\lambda, T-\mu}$
est donc le point fermé naïf
correspondant au couple~$(\lambda, \mu)$ de~$k^2$, l'idéal maximal associé
est le noyau~$(S-\lambda, T-\mu)$
de cette évaluation, et~$y_{S-\lambda, T-\mu}=V(S-\lambda, T-\mu)$.

\deux{adh-vp0-k}
{\bf Retour à l'étude de~$\overline{\{y_{\eta,P}\}}$, où~$P$ est un polynôme
irréductible de~$k(S)[T]$.}
Nous reprenons les notations~$P$ et~$P_0$ du~\ref{ferm-gen-kst}, et allons
décrire un peu plus précisément l'adhérence~$V(P_0)$ de~$x_P$, en regardant sa trace sur
chacune des fibres. 

Comme~$y_{\eta,P}$ est fermé dans~$\psi^{-1}(\eta)$, on a~$V(P_0)\cap \psi^{-1}(\eta)=\{y_{\eta,P}\}$. 
Écrivons~$P_0=\sum a_i T^i$. Soit~$\lambda\in k$.
L'intersection de~$V(P_0)$ avec
la fibre~$\psi^{-1}(x_{S-\lambda})\simeq \spec k[T]$ 
s'identifie au fermé~$V(\sum a_i(\lambda) T^i)$ de~$\spec k[T]$. 

\trois{adh-vp0-fibre-k}
Par définition de~$P_0$, les~$a_i$ sont globalement premiers entre eux ; en particulier ils ne peuvent être
tous nuls modulo~$S-\lambda$, et~$\sum a_i(\lambda)T^i$ est donc un élément
{\em non nul}
de~$k[T]$ ; en conséquence, $V(P_0)\cap \psi^{-1}(x_{S-\lambda})$ est un ensemble
fini de points fermés. 

Cet ensemble peut être vide : c'est le cas si et seulement si~$\sum a_i(\lambda)T^i$
est inversible, c'est-à-dire (compte-tenu du fait qu'il est non nul) si et seulement
si~$a_i(\lambda)=0$ pour tout~$i\geq 1$, ou encore si et seulement si~$S-\lambda$ divise~$a_i$
pour tout~$i\geq 1$. 

\trois{adhvp0-conclu}
Il résulte de ce qui précède
que l'ensemble des éléments~$\lambda$
de~$k$ tels que~$\overline{\{y_{\eta, P}\}}\cap \psi^{-1}(x_{S-\lambda})=\varnothing$
est {\em fini} ; autrement dit, $\overline{\{y_{\eta,P}\}}$ rencontre presque toutes les fibres fermées de~$\psi$,
et donc une infinité de telles fibres
(le corps~$k$ est algébriquement clos, et partant infini). 

En particulier,~$\overline{\{y_{\eta,P}\}}$ n'est pas réduit au singleton~$\{y_{\eta, P}\}$ 
et le point~$y_{\eta, P}$ n'est
pas fermé. Il découle alors de toute l'étude menée ci-dessus
que les points fermés de~$\spec k[S,T]$ sont
exactement les points fermés de ses fibres fermées au-dessus de~$\spec k[S]$, ou encore ses points
«naïfs», c'est-à-dire
ceux dont le corps résiduel est~$k$, ce qu'on savait déjà par ailleurs (\ref{kbarre-x0}). 

\deux{xP-caspart}
{\bf Étude de~$\overline{\{y_{\eta, P}\}}$ : deux exemples explicites.}

\trois{zt2+1}
{\em Le cas où~$P=T^2-S$ et où~$k$ est de caractéristique différente de~$2$.}
Dans ce cas~$P_0=T^2-S$ aussi. Soit~$\lambda\in k$. 
Nous allons
décrire l'intersection de~$\overline{\{y_{\eta, T^2-S}\}}=V(T^2-S)$ avec
la fibre~$\psi^{-1}(x_{S-\lambda})\simeq \spec k[T]$. 
Elle s'identifie à~$V(T^2-\lambda)\subset \spec k[T]$. On distingue deux cas.

\medskip
$\bullet$ {\em Supposons que~$\lambda\neq 0$.}
Le polynôme~$T^2-\lambda$ de~$k[T]$ 
est alors scindé à racines simples (car la caractéristique
de~$k$ est différente de 2) ; il s'écrit $(T-\sqrt \lambda)(T+\sqrt \lambda)$ où l'on désigne par~$\sqrt \lambda$
l'une des deux racines carrées de~$\lambda$.  En conséquence, $\overline{\{y_{\eta, T^2-S}\}}\cap \psi^{-1}(x_{S-\lambda})$
est le sous-ensemble~$V(T-\sqrt \lambda)\cup V(T+\sqrt \lambda)$ de~$\spec k[T]$. 
Il consiste en deux points fermés de corps résiduel~$k$, à savoir les points naïfs
correspondant aux racines~$\sqrt \lambda$ et~$(-\sqrt \lambda)$.
Vus comme points de~$\spec k[S,T]$, 
ce sont les points naïfs correspondants aux couples~$(\lambda, \sqrt \lambda )$ et~$(\lambda, -\sqrt \lambda)$. 

$\bullet$ {\em Supposons que~$\lambda=0$.}
Dans ce cas, $\overline{\{y_{\eta, T^2-S}\}}\cap \psi^{-1}(x_S)$
est le sous-ensemble~$V(T)$
de~$\spec k[T]$. 
Il consiste en un seul point fermé de corps résiduel~$k$ : le point naïf qui correspond
à~$0$. Vu comme point de~$\spec k[S,T]$, 
c'est le point naïf correspondant à l'origine~$(0,0)$.

\trois{kst-hyperbole}
{\em Le cas où~$P=T-(1/S)$.}
Dans ce cas~$P_0=ST-1$. Soit~$\lambda$
appartenant à~$k$.  
Nous allons
décrire l'intersection de~$\overline{\{y_{\eta, T-(1/S)}\}}=V(ST-1)$ avec
la fibre~$\psi^{-1}(x_{S-\lambda})\simeq \spec k[T]$. 
Elle s'identifie à~$V(\lambda T-1)\subset \spec k[T]$. On distingue deux cas.

\medskip
$\bullet$ {\em Supposons que~$\lambda\neq 0$.}
Dans ce cas~$\overline{\{y_{\eta, T-(1/S)}\}}\cap \psi^{-1}(x_{S-\lambda})$
est le sous-ensemble~$V(T-(1/\lambda))$
de~$\spec k[T]$. Il consiste en un seul point fermé de corps résiduel~$k$ : le point naïf qui correspond
à~$(1/\lambda)$.  Vu comme point de~$\spec k[S,T]$, 
c'est le point naïf correspondant au couple~$(\lambda,1/\lambda)$. 

$\bullet$ {\em Supposons que~$\lambda =0$.}
Dans ce cas,~$\overline{\{y_{\eta, T-(1/S)}\}}\cap \psi^{-1}(x_{S-\lambda})$
est le sous-ensemble~$V(1)$
de~$\spec k[T]$, qui est
{\em vide.}

\deux{comment-kst}
{\bf Récapitulation.}
On déduit de ce qui précède
que~$\spec k[S,T]$
comprend trois types de points. 

\trois{naifs-kst}
Il y a tout d'abord les points naïfs, qui sont exactement les points fermés
de~$\spec k[S,T]$, ou encore ceux de corps résiduel~$k$ ; ils constituent un ensemble en bijection naturelle avec~$k^2$. 
Si~$(\lambda, \mu)\in k^2$ le point fermé associé peut se décrire comme~$V(T-\lambda, S-\mu)$
et il correspond à l'idéal maximal~$(T-\lambda, S-\mu)$ ; l'évaluation en ce point
est l'évaluation usuelle des polynômes en~$(\lambda, \mu)$. 

\trois{point-gen-kst}
Il y a ensuite le point générique que nous avons noté~$\xi_\eta$, dont l'adhérence est~$\spec k[S,T]$ tout entier, 
et dont le corps résiduel est~$k(S,T)$. L'idéal correspondant à~$\xi_\eta$ est l'idéal nul, et l'évaluation
en~$\xi_\eta$ est le plongement~$k[S,T]\hookrightarrow k(S)[T]$.

\trois{points-inter}
Il y a enfin une famille de points
«intermédiaires» : ceux que nous avons notés
$\xi_{S-\lambda}$ où~$\lambda \in k$ et~$y_{\eta,P}$, où~$P$
est un polynôme irréductible de~$k(S)[T]$. Leur point commun est le suivant : 
l'idéal premier correspondant à chacun d'eux est engendré par un polynôme irréductible
de~$k[S,T]$ (il s'agit de~$S-\lambda$ pour~$\xi_{S-\lambda}$, et 
de celui que nous avons noté~$P_0$
pour~$y_{\eta,P}$ -- le lecteur vérifiera qu'on obtient ainsi {\em tous}
les polynômes irréductibles de~$k[S,T]$). 

Si~$z$ est l'un de ces points et si~$Q$ désigne le polynôme irréductible de~$k[S,T]$
qui lui correspond, il résulte
de~\ref{gen-ferm-kst}, \ref{adh-vp0-k} et~\ref{adh-vp0-fibre-k}
que~$\overline{\{z\}}=V(Q)$ ne comprend, hormis
le point~$z$ lui-même, que des points naïfs ; en identifiant ces derniers à des éléments de~$k^2$, on
a donc
$$\overline{\{z\}}=\{z\}\cup\{(\lambda, \mu)\in k^2, \;Q(\lambda, \mu)=0\}.$$
L'ensemble~$E:=\{(\lambda, \mu)\in k^2, \;Q(\lambda, \mu)=0\}$, qui n'est autre 
que la courbe algébrique «naïve»
d'équation~$Q=0$, est toujours infini : 

\medskip
$\bullet$ si~$z=\xi_{S-\lambda}$ pour un certain~$\lambda$ alors~$Q=(S-\lambda)$ 
et~$E=\{(\lambda, \mu)\}_{\mu \in k}$ ; 

$\bullet$ si~$z=y_{\eta,P}$ pour un certain~$P$ alors~$Q=P_0$ et l'on a vu
plus haut (\ref{adh-vp0-k} {\em et sq.})
que l'ensemble $E_\lambda:=\{\mu\in k \;\;{\rm t.q.}\;(\lambda,\mu)\in E\}$ est
fini pour tout~$\lambda$, et non vide pour presque tout~$\lambda$.

\medskip
Quant au corps résiduel~$\kappa(z)={\rm Frac}\;k[S,T]/Q$, 
il se décrit plus précisément comme suit (\ref{ferm-gen-kst}, \ref{gen-ferm-kst}) : 

$\bullet$ si~$z=\xi_{S-\lambda}$ pour un certain~$\lambda$
alors~$\kappa(z)\simeq k(T)$ ; 

$\bullet$ si~$z=y_{\eta,P}$ pour un certain~$P$ alors
$\kappa(z)=k(S)[T]/P$. 

\medskip
Dans les deux cas, $\kappa(z)$ est une extension de~$k$ 
de degré de transcendance égal à~$1$. 

\deux{analog-kst-zt}
{\bf Quelques commentaires.} Ce qui précède met bien en évidence l'analogie
entre~$\spec \ZZ[T]$ et~$\spec k[S,T]$. Le lecteur pourra
d'ailleurs se convaincre que nous aurions pu nous contenter
d'un paragraphe coiffant ces deux exemples, en étudiant~$\spec A[T]$ 
{\em via}
le morphisme~$\psi \colon \spec A[T]\to \spec A$ induit par~$A\hookrightarrow A[T]$, où~$A$
est un anneau principal {\em ayant un ensemble infini d'éléments irréductibles}. 

\trois{infini-irred}
La condition sur l'existence d'une infinité d'éléments irréductibles sert simplement
à assurer que si~$y_{\eta,P}$ est un point fermé de la fibre générique~$\psi^{-1}(\eta)$, 
son adhérence rencontre au moins une fibre fermée et n'est en particulier pas réduite à~$\{y_{\eta,P}\}$. 
Nous allons esquisser ici un contre-exemple à ce fait lorsque~$A$ n'a qu'un nombre fini d'irréductibles ; 
les détails sont laissés au lecteur. 

\medskip
Soit~$p$ un nombre premier. On vérifie que le localisé~$\ZZ_{(p)}$ de~$\ZZ$ est principal, 
et a un unique élément irréductible, à savoir~$p$. L'adhérence~$\overline{\{y_{\eta, T-(1/p)}\}}$ 
ne rencontre alors pas {\em la}
fibre fermée de~$\psi \colon \spec \ZZ_{(p)}[T]\to \spec \ZZ_{(p)}$, et~$y_{\eta,T-(1/p)}$ 
apparaît ainsi comme un point fermé de~$\spec \ZZ_{(p)}[T]$
{\em qui est situé sur la fibre générique de~$\psi$}. 
L'évaluation correspondante n'est autre que
la flèche
naturelle de~$\ZZ_{(p)}[T]$ vers~$\QQ$ qui envoie~$T$ sur~$1/p$, 
et qui est 
effectivement surjective puisque~$\QQ=\ZZ_{(p)}[1/p]$. 

\trois{zt-geom}
Cette analogie entre~$\spec \ZZ[T]$ et~$\spec k[S,T]$ permet
de penser en termes géométriques à~$\ZZ[T]$, qui pouvait apparaître
comme de nature davantage algébrique ou arithmétique\footnote
{Insistons sur le fait que, contrairement à une impression qu'on peut avoir à première vue
(par exemple face à la profusion de points étranges sur un objet censé jouer le rôle
du plan affine), la théorie des schémas est faite pour 
favoriser l'intuition géométrique à propos d'objets sur
lesquels elle semblait inopérante
à première vue, et pas pour la chasser
lorsqu'elle est naturellement présente ! }. Nous attirons
par exemple l'attention du lecteur sur la similitude des exemples~\ref{zt-demi} 
et~\ref{kst-hyperbole} : le fermé~$V(2T-1)$ de~$\spec \ZZ[T]$
ressemble beaucoup à l'hyperbole~$V(ST-1)\subset \spec k[S,T]$. Cette dernière rencontre
toutes les droites verticales, c'est-à-dire les fibres fermées de~$\psi$, à l'exception de la fibre
en l'origine : comme~$S$ s'annule en l'origine, l'hyperbole ne coupe pas la fibre correspondante 
-- en fait, elle la rencontre plus précisément «à l'infini», en un sens
que nous préciserons plus loin lorsque nous aurons introduit la géométrie projective. 

Le même phénomène vaut pour~$V(2T-1)\subset \spec \ZZ[T]$ : puisque~$2$, vu comme
fonction sur~$\spec \ZZ$, s'annule en~$x_2$, le fermé~$V(2T-1)$ ne rencontre pas la fibre
$\psi^{-1}(x_2)$ -- et là encore, nous verrons plus bas qu'il la rencontre en fait «à l'infini»,
ce qui permettra au lecteur imaginatif de penser à cette
fibre comme à une asymptote de~$V(2T-1)$. 

\section{Compléments sur la topologie de~$\spec A$}
\markboth{Le spectre comme espace topologique}{Compléments sur la topologie de~$\spec A$}

\subsection*{Idéaux saturés et fermés de Zariski}

\deux{intro-idsat}
Soit~$A$ un anneau et soit~$I$ un idéal de~$A$. On note~$\sqrt I$ l'image réciproque
du nilradical de~$A/I$ par la flèche quotient~$A\to A/I$ ; c'est un idéal de~$A$, qu'on peut
également décrire comme
l'ensemble des~$a\in A$ pour  lesquels il existe~$n\in \NN$ tel que~$a^n\in I$. 
On dit que~$\sqrt I$ est le {\em radical}
de~$I$.

\trois{a-inclus-arad}
Les faits suivants découlent immédiatement
de la définition : si~$I$ est un idéal de~$A$ alors~$I\subset \sqrt I$ et~$\sqrt{\sqrt I}=\sqrt I$ ; si~$J$
est un idéal de~$A$ tel que~$I\subset J$ alors~$\sqrt I \subset \sqrt J$ ; l'idéal~$\sqrt{(0)}$ est le nilradical
de~$A$. 

\trois{coince-rad}
Il résulte de~\ref{a-inclus-arad}
que si~$J$ est un idéal de~$A$ tel que
$$I\subset J\subset \sqrt I$$ alors~$\sqrt J=\sqrt I.$

\deux{def-idsat}
Nous dirons qu'un idéal~$I$ de~$A$ est
{\em saturé}\footnote{Cette terminologie
n'est pas particulièrement
standard, mais il ne semble pas
exister d'adjectif universellement utilisé pour qualifier un tel idéal.} s'il est égal à son radical ; cela revient
à demander que~$A/I$ soit réduit. 

\trois{ex-id-sat}
Si~$I$ est un idéal de~$A$ alors~$\sqrt I$ est saturé (\ref{a-inclus-arad}). 

\trois{id-prem-sat}
Si~$\got p$ est un idéal premier de~$A$ il est saturé ($A/\got p$ est intègre, et {\em a fortiori}
réduit).  

\deux{ferme-speca}
La formule~$I\mapsto V(I)$ définit une surjection décroissante (pour l'inclusion)
de l'ensemble des idéaux de~$A$ vers l'ensemble des fermés de~$\spec A$. 

\medskip
Par ailleurs, la formule
$$F\mapsto \sch I(F):=\{f\in A\;\;{\rm t.q.}\;\;\forall x\in F\; \; f(x)=0\}$$
définit une application décroissante de l'ensemble des fermés de~$\spec A$ vers l'ensemble des idéaux de~$A$. 

\medskip
Il résulte immédiatement
des définitions que $I\subset \sch I(V(I))$ pour tout idéal~$I$ de~$A$,
et que~$F\subset V(\sch I(F))$ pour tout fermé~$F$
de~$\spec A$. 

\trois{if-sat}
Soit~$F$ un fermé de~$\spec A$, et soit~$f\in A$ tel que~$f^n\in \sch I(F)$ pour un certain entier~$n$. On a
alors~$f^n(x)=0$ pour tout~$x\in F$, et partant~$f(x)=0$ pour tout~$x\in \sch F$ ; ainsi, $f\in \sch I(F)$ et~$\sch I(F)$ est saturé. 

\trois{ivi-rad}
Soit~$I$ un idéal de~$A$ et soit~$f\in A$. La flèche~$A\to A/I$ induit un homéomorphisme
de~$\spec A/I$ sur~$V(I)$. En conséquence on a pour tout~$f\in A$ l'équivalence entre les assertions suivantes : 

\medskip
i) $f\in \sch I(V(I))$, {\em i.e.}
$f$ s'annule en tout point~$x$ de~$V(I)$ ;

ii) la classe~$\bar f$ de~$f$ modulo~$I$ s'annule en tout point de~$\spec A/I$. 

\medskip
Mais cette dernière condition revient à demander que~$\bar f$ soit nilpotente dans~$A/I$, 
donc que~$f\in \sqrt I$. On a ainsi démontré que
$$\sch I(V(I))=\sqrt I.$$

\deux{lemme-idsat-ferm}
{\bf Lemme.}
{\em Les flèches~$I\mapsto V(I)$ et~$F\mapsto \sch I(F)$ établissent une bijection
décroissante d'ensembles
ordonnés entre l'ensemble des idéaux
{\em saturés}
de~$A$ et l'ensemble des fermés de~$\spec A$.}

\medskip
{\em Démonstration.}
Soit~$I$ un idéal saturé de~$A$. On a alors~$\sch I(V(I))=\sqrt I=I$ (la première égalité est due à~\ref{ivi-rad}, 
la seconde à l'hypothèse de saturation). 

\medskip
Soit~$F$ un fermé de~$\spec A$ ; écrivons~$F=V(J)$ pour un certain idéal~$J$ de~$A$. 
On sait que~$\sch I(F)$ est saturé (\ref{if-sat}). On a~$F\subset V(\sch I(F))$. Par ailleurs 
comme~$F=V(J)$ il vient~$J\subset \sch I(F)$ et donc~$V(\sch I(F))\subset V(J)=F$. Ainsi,
$V(\sch I(F))=F$, ce qui achève la démonstration.~$\Box$ 

\deux{comment-vi-if}
{\bf Commentaire.}
Le lemme~\ref{lemme-idsat-ferm}
ci-dessus affirme en particulier
que la restriction de~$I\mapsto V(I)$ à l'ensemble des idéaux saturés est injective. 
Mais on peut également en déduire une condition nécessaire et suffisante pour
que deux idéaux~$I$ et~$J$ (non nécessairement saturés) vérifient l'égalité~$V(I)=V(J)$. 
En effet, comme~$F\mapsto \sch I(F)$ est injective en vertu de~{\em loc. cit.}, ce sera le cas
si et seulement si~$\sch I(V(I))=\sch I(V(J))$, c'est-à-dire si et seulement si~$\sqrt I=\sqrt J$
(\ref{ivi-rad}). 

C'est vrai en particulier lorsque~$J=0$. On a donc~$V(I)=V(0)=\spec A$ si et seulement si~$\sqrt I$
est égal au nilradical de~$A$ ; on voit immédiatement que c'est le cas si et seulement si~$I$ lui-même
est contenu dans le nilradical de~$A$, et l'on retrouve ainsi ce qui avait été mentionné en~\ref{spec-asuri}.

\subsection*{Le cas d'une algèbre de type fini sur un corps algébriquement clos}

\deux{constr-kalgclos}
Supposons maintenant que~$A$ est une algèbre de type fini sur un corps algébriquement clos~$k$, 
et choisissons un isomorphisme de~$k$-algèbres
~$A\simeq k[T_1,\ldots, T_n]/(P_1,\ldots, P_r)$.
On reprend
les notations de~\ref{cas-k-alg}
{\em et sq.} : on pose~$X=\spec A$, et l'on note~$X(k)$ l'ensemble~$\hom_k(A,k)$, qui coïncide avec l'ensemble~$X_0$
des points fermés de~$X$. L'ensemble~$X(k)$
s'identifie par ailleurs à celui des~$n$-uplets~$(x_1,\ldots, x_n)\in k^n$ en lesquels les~$P_j$ s'annulent.

\deux{top-zariski-xk}
On munit~$X(k)$ de la topologie induite par celle de~$X$, qu'on appelle encore topologie de Zariski ; ses fermés sont les
parties de la forme~$V(E)\cap X(k)$ où~$E$ est une partie de~$A$ (on peut d'ailleurs 
se limiter aux idéaux de~$A$) et une base d'ouverts de~$X(k)$ est formée des parties
de la forme~$D(f)\cap X(k)$ où~$f\in A$. 

\medskip
Si l'on voit~$X(k)$ comme
un ensemble de~$n$-uplets, alors pour tout~$E\subset A$ et tout~$f\in A$ on a
$$E\cap X(k)=\{(x_1,\ldots, x_n)\in X(k)\;{\rm t.q.}\;g(x_1,\ldots,x_n)=0\;\forall \;g\in E\}$$
et
$$D(f)\cap X(k)=\{(x_1,\ldots, x_n)\in X(k)\;{\rm t.q.}\;f(x_1,\ldots,x_n)\neq 0\
\}.$$

\deux{top-const}
Soit~$T$ un espace topologique. On notera~$\sch C(T)$ le plus petit
sous-ensemble de~$\sch P(T)$ contenant les fermés, les ouverts, et qui est stable
par unions finies, intersections finies et passage au complémentaire. On vérifie
aussitôt que~$\sch C(T)$ est l'ensemble des parties de~$T$
de la forme
$\bigcup_{i\in I}U_i\cap F_i$ où~$I$ est un ensemble fini, où les~$U_i$ sont des ouverts et où
les~$T_i$
sont des fermés.

\deux{const-x-xk}
{\bf Proposition.}
{\em L'application~$C\mapsto C(k):=C\cap X(k)$ induit une bijection de~$\sch C(X)$ 
sur~$\sch C(X(k))$.}

\medskip
{\em Démonstration.} Par définition de la topologie induite, $C(k)\in \sch C(X(k))$
pour tout~$C\in \sch C(X)$, et toute partie appartenant à~$\sch C(X(k))$ est de la forme~$C(k)$
pour une telle~$C$. Il reste donc à s'assurer que~$C\mapsto C(k)$ est injective. 

\trois{vide-const-k}
{\em Un cas particulier.}
Soit~$C\in \sch C(X)$ telle que~$C(k)=\varnothing$ ; nous allons montrer que~$C=\varnothing$. 
Comme~$C$ est une union finie de parties de la forme~$U\cap F$, où~$U$ est ouvert et~$F$ fermé, 
on peut supposer que~$C=U\cap F$. Le fermé~$F$ s'écrit~$V(I)$ pour un certain~$I$, et s'identifie donc
à~$\spec A/I$. Quitte à remplacer~$A$ par~$A/I$, on peut supposer~$F=\spec A$ et~$C=U$. Dans ce cas~$C$ est réunion 
d'ouverts de la forme~$D(f)$, avec~$f\in A$, et l'on est ainsi ramené au cas où~$C=D(f)$. L'ouvert~$C$ s'identifie
alors au spectre de la~$k$-algèbre de type fini~$A_f$. Comme~$C(k)=\emptyset$, le spectre de~$A_f$ n'a pas de
point fermé, ce qui veut dire que~$A_f$ est nulle et que~$C=\spec A_f=\varnothing$. 

\trois{vide-const-casgen}
{\em Le cas général.}
Soient~$C$ et~$D$ deux parties appartenant
à~$\sch C(X)$ telles que~$C(k)=D(k)$. Soit~$C'$ l'intersection
de~$C$ et du complémentaire de~$D$, et soit~$D'$ l'intersection de~$D$ et du complémentaire
de~$C$ ; les parties~$C'$ et~$D'$ appartiennent à~$\sch C(X)$. 
Par hypothèse~$C'(k)=D'(k)=\varnothing$ ; le cas particulier traité au~\ref{vide-const-k}
ci-dessus assure alors que~$C'=D'=\varnothing$, et donc que~$C=D$.~$\Box$

\deux{ideaux-fermes-k}
{\bf Quelques conséquences.}

\trois{inclu-x-xk}
Soient~$C$ et~$D$ deux parties appartenant à~$\sch C(X)$. On a les équivalences
$$C\subset D \iff C\cap (X\setminus D)=\varnothing\iff C(k)\cap (X(k)\setminus D(k))=\varnothing \iff C(k)\subset D(k)$$
(la deuxième équivalence découle de la proposition~\ref{const-x-xk}
ci-dessus, les autres sont tautologiques). 

\trois{ferme-x-xk}
Soit~$F$ un fermé de Zariski de~$X(k)$ ; on notera~$\sch I(F)$ l'idéal de~$A$ formé des 
fonctions~$f$ qui s'annulent en tout point de~$F$. 

\medskip
Le fermé~$F$ est de la forme~$G(k)$, où~$G$ est un fermé de Zariski de~$X$
qui est uniquement déterminé d'après la proposition~\ref{const-x-xk}
ci-dessus. Si~$f\in A$ on a les équivalences
$$f\in \sch I(F)\iff G(k)\subset V(f)(k) \iff G\subset V(f)\iff f\in \sch I(G),$$
la deuxième équivalence résultant de~\ref{inclu-x-xk}. Ainsi~$\sch I(F)=\sch I(G)$. 

\trois{bij-ferme-ideal-k}
En combinant~\ref{ferme-x-xk}, 
la proposition~\ref{const-x-xk} et le lemme~\ref{lemme-idsat-ferm}, 
on voit que
$I\mapsto V(I)(k)$ et~$F\mapsto \sch I(F)$ établissent une bijection entre l'ensemble des idéaux
saturés de~$A$ et l'ensemble des fermés de Zariski de~$X(k)$. 

\deux{comment-x-xk}
La proposition~\ref{const-x-xk}
et ses conséquences signalées en~\ref{ideaux-fermes-k}
{\em et sq.}
expliquent pourquoi l'on peut, pour un grand nombre de questions, se
passer du langage des schémas lorsqu'on fait de la géométrie
algébrique sur un corps algébriquement clos : dans ce contexte, 
une partie définissable de manière algébrique (c'est-à-dire par des conditions d'annulation
ou de non-annulation de polynômes) est connue sans ambiguïté dès qu'on connaît
ses points naïfs, et l'on ne perd donc
pas grand-chose à ne considérer que lesdits points. 

\subsection*{Espaces topologiques irréductibles, composantes irréductibles, dimension de Krull}

\deux{def-espirr}
{\bf Définition.}
Soit~$X$ un espace topologique. On dit que~$X$ est {\em irréductible}
si~$X$ est non vide et si pour tout couple~$(Y,Z)$ de fermés de~$X$ tels que
$X=Y\cup Z$ on a~$X=Y$ ou~$X=Z$. 

\trois{irr-bourbaki}
On peut donner une définition 
bourbakiste de l'irréductibilité -- analogue à celle de l'intégrité 
donnée
en~\ref{ann-integre-bourb} -- en disant que~$X$ est irréductible si et seulement si toute union
finie de fermés stricts de~$X$ est stricte ; cela force en particulier la réunion
{\em vide} de tels fermés (qui est l'ensemble vide) à être stricte, et donc~$X$ à être non vide. 

\trois{def-duale-irr}
Il est tautologique qu'un espace topologique~$X$ est irréductible
si et seulement si~$X\neq \varnothing$ et si tout ouvert 
non vide de~$X$ est dense dans~$X$ (passer au complémentaire
dans la définition initiale). 

On en déduit que si~$X$ est un espace topologique irréductible, 
tout ouvert non vide de~$X$ est encore irréductible. 

\trois{irr-adh}
Soit~$X$ un espace topologique. Si~$X$ possède une partie dense irréductible, 
il est irréductible : c'est immédiat. 

\medskip
Comme un singleton est trivialement irréductible, tout espace topologique possédant
un point générique (c'est-à-dire dense) est irréductible. 

\trois{irr-connexe}
Il résulte de la définition qu'un espace topologique irréductible est en particulier connexe. 

\deux{comment-irred}
La notion d'espace irréductible n'a guère d'intérêt lorsqu'on s'intéresse
aux  espaces topologiques usuels. On démontre par exemple
aisément (le lecteur est invité à le faire à titre d'exercice)
qu'un espace topologique séparé est irréductible
si et seulement si c'est un singleton. 

Elle est par contre extrêmement utile en géométrie algébrique, 
qui manipule des espaces à la topologie assez grossière et très combinatoire. 
La proposition suivante en est une bonne illustration.

\deux{prop-irred-sp}
{\bf Proposition.}
{\em Soit~$A$ un anneau. Les assertions suivantes sont équivalentes. 

\medskip
i) $\spec A$ a un unique point générique. 

ii) $\spec A$ a un point générique. 

iii) $\spec A$ est irréductible. 

iv) Le quotient~$A_{\rm red}$
de~$A$ par son nilradical est intègre. }

\medskip
{\em Démonstration.} 
On sait que la flèche quotient~$A\to A_{\rm red}$
induit un homéomorphisme~$\spec A_{\rm red}\simeq \spec A$ (\ref{spec-asuri}). 
On peut donc remplacer~$A$ par~$A_{\rm red}$, et ainsi supposer que~$A$ est réduit, 
c'est-à-dire que~$A=A_{\rm red}$. 

\medskip
Il est clair que~i)$\Rightarrow$ii)$\Rightarrow$iii). 

\medskip
Supposons que~iii) soit vraie, et montrons que~$A=A_{\rm red}$ est intègre. 
Comme~$\spec A$ est irréductible, il est non vide et~$A$ est donc non nul. 

Soient~$f$ et~$g$ deux éléments de~$A$ tels que~$fg=0$. On a alors
$\spec A=V(fg)=V(f)\cup V(g)$. 

Par irréductibilité, il vient
$\spec A =V(f)$ ou~$\spec A=V(g)$. En conséquence, $f$ est nilpotente
ou~$g$ est nilpotente ; comme~$A$ est réduit, on a~$f=0$ ou~$g=0$ et~$A$ est intègre.

\medskip
Supposons que~iv) soit vraie, et montrons~i). Soit~$x\in \spec A$ et soit~$\got p$
l'idéal premier correspondant. Le point~$x$ est générique si et seulement si tout idéal premier
de~$A$ contient~$\got p$ ; mais il est clair que cette propriété est vérifiée par~$(0)$ (qui est premier
car~$A$ est intègre) et par lui seul, d'où~i).~$\Box$ 

\deux{fermes-irred}
Soit~$A$ un anneau.

\trois{bij-points-fermirr}
Soit~$F$ un fermé de~$\spec A$. On a~$F=V(I)$ pour
un certain idéal~$I$ de~$A$, et~$F$ est homéomorphe à~$\spec A/I$. Par la proposition~\ref{prop-irred-sp}
ci-dessus, 
$F$ est irréductible si et seulement si il admet un point générique, lequel est alors unique. 

\medskip
On en déduit que~$x\mapsto \overline{\{x\}}$ établit une bijection entre~$\spec A$ et l'ensemble
de ses fermés irréductibles, la réciproque envoyant un fermé~$F$ sur son unique point générique. 

\trois{bij-idprem-fermirr}
Soit~$I$ un idéal saturé de~$A$. Par définition, $A/I$ est réduit ; il s'ensuit
d'après~{\em loc. cit.}
que
$V(I)\simeq \spec A/I$ est irréductible si et seulement si~$A/I$ est intègre, c'est-à-dire si et seulement si~$I$ est premier. 
En conséquence, $\got p\mapsto V(\got p)$ et~$F\mapsto \sch I(F)$ 
établissent une bijection décroissante d'ensembles ordonnés entre l'ensemble
des idéaux premiers de~$A$ et celui des fermés irréductibles de~$\spec A$. Bien entendu, cette bijection 
est simplement la traduction en termes d'idéaux premiers de la bijection
décrite au~\ref{bij-points-fermirr}
ci-dessus par la formule plus
géométrique~$x\mapsto \overline{\{x\}}$. 

\deux{def-sobre}
{\bf Intermède culturel.}
On dit qu'un espace topologique est
{\em sobre}
si chacun de ses fermés irréductibles admet un et un seul point générique. On vient
de voir que le spectre d'un anneau est sobre.

\trois{sobrification}
Tout espace topologique admet une {\em sobrification} ; nous laissons
au lecteur le soin de la définir, et de la construire -- suivant la philosophie habituelle
: {\em just do it !}

\trois{sobrification-topos}
On démontre que si~$X$ et~$Y$
sont deux espaces topologiques, les catégories des faisceaux (d'ensembles)
sur~$X$ et sur~$Y$ sont équivalentes si et seulement si~$X$ et~$Y$ ont même sobrification. 

\deux{points-schem-interp}
{\bf Interprétation des points schématiques en
termes classiques.}
Soit~$k$ un corps algébriquement clos et soit~$A$
une~$k$-algèbre de type fini ; posons~$X=\spec A$.

\trois{irred-irresk}
Il résulte
de la proposition~\ref{const-x-xk}
et de~\ref{inclu-x-xk}
que~$F\mapsto F(k)$ établit une bijection entre l'ensemble
des fermés irréductibles de~$X$ et l'ensemble
des fermés irréductibles de~$X(k)$ ; en conséquence, 
$x\mapsto \overline{\{x\}}(k)$ établit une bijection entre~$X$ et l'ensemble
des fermés irréductibles de~$X(k)$. 

\medskip
Insistons à ce propos sur le fait que tous les points
de l'espace topologique~$X(k)$ sont fermés ; en conséquence, 
un fermé irréductible non singleton de~$X(k)$ (c'est-à-dire,
par ce qui précède, un fermé de la forme~$\overline{\{x\}}(k)$ avec~$x\notin X(k)$)
n'a pas de point générique dans~$X(k)$. 

\trois{naif-schema}
On voit donc 
que topologiquement, 
on peut construire~$X$ 
à partir de~$X(k)$ en rajoutant un point
générique par fermé irréductible non singleton. 
L'espace~$X$ apparaît
ainsi comme la sobrification de~$X(k)$. 

\trois{interret-point-gen}
Soit~$x$ un point de~$X$ et soit~$G$
le fermé irréductible~$\overline{\{x\}}(k)$. 
Pour tenter d'appréhender intuitivement le point~$x$, 
on pourra se référer au slogan suivant : {\em une propriété
(raisonnable) est vraie en~$x$ si et seulement si elle est vraie
sur un ouvert de Zariski non vide de~$G$.}

C'est par exemple le cas
pour l'annulation ou la non-annulation des fonctions : 
si~$f\in A$ et si~$f(x)=0$ alors~$f$ est nulle en tout point de~$G$ ; 
et si~$f(x)\neq 0$ l'ouvert~$D(f)\cap G$ est non vide, puisque
$D(f)\cap \overline{\{x\}}\neq \varnothing$ (il contient~$x$). 
Mais ce le sera aussi pour un grand nombre d'autres propriétés
plus subtiles. 

\medskip
Le passage
de la variété naïve~$X(k)$ au schéma~$X$
s'apparente ainsi à une déclinaison d'un procédé très répandu en mathématiques : 
on estime souvent avoir intérêt, pour des raisons de confort psychologique, 
à inventer un objet conceptuellement un peu compliqué pour pouvoir 
remplacer des énoncés de la forme «{\em il existe un ensemble}
sur lequel telle
propriété est vraie»
par des énoncés plus agréables de la forme 
«telle propriété est vraie {\em en tel point}». 

Donnons un exemple de ce type de démarche, 
qui vous est certainement très familier : 
{\em la construction du corps des réels}. 
On remplace
l'objet assez simple~$\QQ$ par l'objet plus compliqué~$\RR$, 
mais on gagne en simplicité des assertions. Par exemple si~$P\in \QQ[T]$, 
on a équivalence entre~«$P(\sqrt 2)>0$»
et~«il existe un entier~$M>0$ et un entier~$N>0$ tel que pour tout nombre rationnel positif~$r$
satisfaisant les inégalités~$2-1/M\leq r^2\leq 2+1/M$ on ait~$P(r)>1/N$».

\subsection*{Espaces noethériens et composantes
irréductibles}

\deux{def-noeth}
{\bf Définition.}
Soit~$X$ un espace topologique. On vérifie sans peine que les assertions suivantes 
sont équivalentes : 

\medskip
i) tout ensemble non vide de fermés de~$X$ admet un élément minimal (pour l'inclusion) ; 

ii) toute suite décroissante de fermés de~$X$ est stationnaire. 

\medskip
Lorsqu'elles sont satisfaites, on dit que~$X$
est {\em noethérien}. 

\deux{prop-base-noeth}
Il découle immédiatement de la définition
que tout fermé d'un espace topologique noethérien est encore noethérien. 

Nous invitions par ailleurs le lecteur à démontrer qu'un espace topologique~$X$
est noethérien si et seulement si tous ses ouverts sont quasi-compacts. 

\deux{ex-esp-noeth}
{\bf Exemple.}
Soit~$A$ un anneau noethérien. Tout ensemble non vide d'idéaux de~$A$ admet un élément maximal ; 
c'est en particulier vrai lorsqu'on se restreint aux ensembles d'idéaux saturés, et il découle
alors du lemme~\ref{lemme-idsat-ferm}
que l'espace topologique~$\spec A$ est noethérien. 

\medskip
\deux{contre-ex-noeth}
En fait, le lemme~\ref{lemme-idsat-ferm}
garantit précisément que
si~$A$ est un anneau,~$\spec A$ est noethérien si et seulement
si tout ensemble non vide d'idéaux saturés de~$A$ a un élément maximal. C'est
une propriété
{\em a priori}
plus faible que la noethérianité, et l'on peut effectivement
construire un exemple d'anneau non noethérien à spectre noethérien. 
Par exemple, soit~$k$ un corps, soit~$I$ l'idéal de~$k[X_i]_{i\in \NN}$
engendré par tous les monômes de degré~$2$, et soit~$A$
le quotient~$k[X_i]_{i\in \NN}/I$. 
Nous laissons au lecteur le soin de vérifier les points suivants : 

$\bullet$ l'idéal~$(\overline{X_i})_{i\in \NN}$ est le seul idéal premier de~$A$, et $\spec A$
est donc un singleton, évidemment noethérien ; 

$\bullet$ l'idéal~$(\overline{X_i})_{i\in \NN}$ n'est pas de type fini, et~$A$ n'est donc pas noethérien. 

\deux{lemme-dec-irr}
{\bf Lemme.}
{\em Soit~$X$ un espace topologique
noethérien. Il existe un ensemble
fini~$E$ de fermés
irréductibles de~$X$
possédant les propriétés suivantes : 

\medskip

i) les fermés appartenant à~$E$ sont deux à deux non comparables pour l'inclusion ; 

ii) $X$ est la réunion des fermés appartenant à~$E$. 

\medskip
De plus si~$E$ est un tel ensemble, 
tout fermé irréductible de~$X$ est contenu dans (au moins)
un fermé appartenant à~$E$, et~$E$
apparaît ainsi comme l'ensemble
des fermés irréductibles maximaux de~$X$. Il est donc unique, et ses
éléments sont appelés les {\em composantes irréductibles}
de~$X$.
}

\medskip
{\em Démonstration.}
On procède en plusieurs étapes. 

\trois{prem-etap-compirr}
{\em Première étape : $X$ est une réunion finie de fermés irréductibles.}
On suppose que ce n'est pas le cas, et l'on note~$F$ l'ensemble des fermés de~$X$
qui ne sont pas réunion finie de fermés irréductibles. Par hypothèse, $X\in F$ et~$F\neq \emptyset$ ; comme~$X$
est noethérien, $F$ admet un élément minimal~$Y$. 

\medskip
Comme~$Y\in F$, il n'est pas réunion finie de fermés irréductibles ; en particulier, il est non vide (sinon, 
ce serait la réunion vide de tels fermés) et non irréductible. Il existe donc deux fermés stricts~$Z$ et~$T$
de~$Y$ tels que~$Y=Z\cup T$. Comme~$Z$ et~$T$ sont strictement contenus dans~$Y$, la minimalité de~$Y$
implique qu'ils n'appartiennent pas à~$F$. Chacun d'eux est donc une union finie de fermés irréductibles, et il en va
dès lors de même de~$Y$, ce qui est absurde. 

\trois{fin-preuve-noetherirr}
Il existe donc un ensemble fini~$E$ de fermés irréductibles de~$X$ satisfaisant~ii). Si~$Y$ et~$Z$ sont deux fermés
appartenant à~$E$ avec~$Y\subsetneq Z$, alors~$X=\bigcup_{T\in E,T\neq Y}T$, et on peut donc retirer~$Y$ de~$E$
sans altérer la validité de~ii). En recommençant l'opération autant de fois que nécessaire, on obtient bien un ensemble
fini de fermés irréductibles de~$X$ qui satisfait~i) et~ii). 

\trois{conclu-irr}
Soit maintenant~$E$ un tel ensemble et soit~$Y$ un fermé irréductible de~$X$. Comme~$X$ est réunion
des éléments de~$E$, le fermé~$Y$ est réunion de ses fermés~$Y\cap Z$ où~$Z$ parcourt~$E$. Comme~$Y$
est irréductible, l'un au moins de ces fermés n'est pas strict ; il existe donc~$Z\in E$ tel que~$Y\cap Z=Y$, 
c'est-à-dire tel que~$Y\subset Z$, ce qui achève la démonstration.~$\Box$

\deux{ex-compirr}
{\bf Exemple.}
Soit~$k$ un corps et soit~$f$ un élément non nul 
de~$k[T_1,\ldots, T_n]$. Le fermé
de Zariski~$V(f)$ de~$\spec k[T_1,\ldots, T_n]$ s'identifie
à~$\spec k[T_1,\ldots, T_n]/(f)$. L'anneau~$ k[T_1,\ldots, T_n]/(f)$
étant noethérien, l'espace topologique~$V(f)$ est lui aussi noethérien, 
et est en conséquence justiciable du lemme précédent. Nous allons 
décrire explicitement
ses composantes irréductibles. 

\medskip
Comme~$f$ est un élément non nul de l'anneau factoriel~$k[T_1,\ldots,T_n]$, 
il s'écrit
comme un produit fini~$\prod P_i^{n_i}$ où les~$P_i$ sont
des polynômes irréductibles deux à deux non associés et les~$n_i$ des entiers~$>0$. 
On a alors
$$V(f)=\bigcup V(P_i^{n_i})=\bigcup V(P_i).$$

\medskip
Fixons~$i$. Comme~$P_i$ est un polynôme irréductible de l'anneau factoriel
$k[T_1,\ldots, T_n]$, l'idéal~$(P_i)$ est premier, et~$V(P_i)$ est donc
irréductible d'après~\ref{bij-idprem-fermirr}.  

\medskip
Par ailleurs, si~$i$ et~$j$ sont deux indices distincts, $(P_i)$ et~$(P_j)$ sont
non comparables pour l'inclusion (puisque~$P_i$ et~$P_j$ ne le sont pas pour
la divisibilité), et~$V(P_i)$ et~$V(P_j)$ ne le sont donc
pas non plus d'après~{\em loc. cit.}

\medskip
En conséquence, les~$V(P_i)$ sont exactement les composantes irréductibles de~$V(f)$. 

\subsection*{Dimension de Krull}

\deux{def-dimkrulltop}
{\bf Définition.}
Soit~$X$ un espace topologique. La {\em dimension de Krull}
de~$X$ est la borne supérieure de l'ensemble des entiers~$n$ pour lesquels
il existe une chaîne strictement croissante

$$X_0 \subsetneq X_1\subsetneq\ldots \subsetneq X_n$$ où
les~$X_i$ sont des fermés irréductibles de~$X$. 

\trois{comment-krulltopo}
{\em Commentaires.}
L'ensemble d'entiers dont on prend la borne supérieure
dans la définition ci-dessus est vide
si et seulement si~$X$ n'a pas de fermés irréductibles, 
ce qui signifie que~$X=\varnothing$ : dans le cas contraire, 
il existe~$x\in X$, et~$\overline{\{x\}}$ est un fermé irréductible
de~$X$.  

Si~$X=\varnothing$, sa dimension de Krull est donc égale à~$-\infty$
(voir la note de bas de page au paragraphe~\ref{valeurs-krull}). Sinon, c'est un élément
de~$\NN\cup\{+\infty\}$, qui vaut~$+\infty$ si et seulement si~$X$ possède des chaînes strictement
croissantes arbitrairement longues de fermés irréductibles.

\trois{krull-geom}
La dimension de Krull n'est pas une notion pertinente 
pour comprendre les espaces topologiques usuels. Ainsi, 
comme les seuls fermés irréductibles de~$\RR^n$ sont les points, 
la dimension de Krull de~$\RR^n$ est nulle quel que soit~$n$. 

\medskip
Elle est en revanche tout à fait adaptée aux espaces topologiques
rudimentaires qui interviennent en géométrie algébrique, et correspond
alors parfaitement à l'idée intuitive qu'on se fait de la dimension. 

\medskip
Par exemple, 
dire qu'un espace~$X$~est de dimension de Krull égal à 2 signifie que les chaînes
strictement croissantes
de fermés irréductibles de~$X$ les plus longues
qu'on puisse trouver
comportent trois éléments (la numérotation 
commence à~$0$). Or c'est ce que l'on attend, indépendamment de la définition 
précise donnée à ce terme, d'une surface algébrique, sur laquelle une telle chaîne
doit être constituée d'un point, d'une courbe irréductible, et de la surface elle-même. 

\trois{dim-top}
Signalons qu'il existe une notion bien plus fine de dimension
en topologie générale qui est pertinente pour tout ce qui
est peu ou prou modelé sur~$\RR$ ; par exemple,
$\RR^n$ est de dimension~$n$. Mais nous n'en aurons pas besoin
dans le cadre de ce cours, et ne donnerons pas sa définition. 

\deux{dim-krull-ann}
Soit~$A$ un anneau. Il résulte des définitions 
et du dictionnaire entre fermés
irréductibles de~$\spec A$ et idéaux premiers de~$A$
que la
dimension de Krull de l'espace topologique~$\spec A$ est égale
à la dimension de Krull de~$A$. 

\medskip
Le théorème~\ref{prop-dimkrull-transc}
implique alors que si~$A$ est une algèbre intègre et de type fini sur
un corps~$k$, la dimension de Krull de~$\spec A$ est égale au
degré de transcendance de~${\rm Frac}\;A$ sur~$k$. 

\chapter{La notion de schéma}

\section{La catégorie des schémas}
\markboth{La notion de schéma}{La catégorie des schémas}

\subsection*{Le spectre comme espace localement annelé}

\deux{intro-spec-locann}
Soit~$A$ un anneau. Le but de ce qui suit est de
munir l'espace topologique~$\spec A$ d'une structure d'espace
localement annelé, c'est-à-dire de construire un faisceau d'anneaux
sur~$\spec A$ dont les fibres soient des anneaux locaux. 

\medskip
Nous aurons en fait également
besoin plus loin
de considérer sur l'espace
localement annelé~$\spec A$ des faisceaux de modules d'un certain type. 
Aussi avons-nous choisi de donner une construction
qui fournit directement ces faisceaux de modules, dont notre faisceau d'anneaux
apparaîtra simplement comme un cas particulier. 

\deux{def-mpref}
Soit~$M$ un~$A$-module et soit~$U$ un ouvert
de~$\spec A$. On note~$S(U)$
l'ensemble des éléments~$f$ de~$A$ qui ne s'annulent pas sur~$U$. C'est une partie 
multiplicative de~$A$, égale à~$A\ti$ si~$U=\spec A$. On désigne par~$M_{\rm pref}(U)$ le~$S(U)^{-1}A$-module~$S(U)^{-1}M$.
Si~$V$ est un ouvert contenu dans~$U$, on dispose d'une application~$A$-linéaire naturelle~$M_{\rm pref}(U)\to M_{\rm pref}(V)$, 
et~$M_{\rm pref}$
apparaît ainsi comme un préfaisceau de~$A$-modules sur~$\spec A$.

\trois{apref-alg}
On remarque que~$A_{\rm pref}$ hérite quant à lui d'une structure plus riche :
c'est un
préfaisceau de~$A$-algèbres, et
$M_{\rm pref}$ est de manière naturelle
un~$A_{\rm pref}$-module. 

\trois{mpref-fibres}
{\em Les fibres de~$M_{\rm pref}$.}
Fixons~$x\in \spec A$, et soit~$\Sigma$ l'ensemble des éléments
de~$A$ non nuls en~$x$ ; notons que~$\Sigma$
est précisément le complémentaire dans~$A$
de l'idéal premier~$\got p$
correspondant à~$x$. 

\medskip
L'ensemble des voisinages ouverts de~$x$, muni de l'ordre
{\em opposé}
à celui de l'inclusion, est filtrant ; si~$U$ et~$V$ sont deux
voisinages ouverts de~$x$ avec~$V\subset U$ alors~$S(U)\subset S(V)$, et la 
réunion des~$S(U)$ lorsque~$U$ parcourt l'ensemble
des voisinages ouverts de~$x$ est égale à~$\Sigma$ (si~$f\in \Sigma$, 
alors~$f\in S(D(f))$ ). 

\medskip
On déduit alors de~\ref{mp-limind-mf}
que~$M_{{\rm pref},x}$ s'identifie à~$\Sigma^{-1}M$, c'est-à-dire
à~$M_{\got p}$.  
En particulier,
$A_{{\rm pref},x}$ est l'anneau local~$A_{\got p}$. 

\trois{mpref-loc-mf}
Nous allons faire une remarque de nature technique
à propos de la description locale des sections de~$M_{\rm pref}$, qui nous servira
par la suite. 

Soit~$U$ un ouvert de~$\spec A$, soit~$\sigma$ une section de~$M_{\rm pref}$
sur~$U$ et soit~$x\in U$. On peut écrire~$\sigma=\mu/s$ avec~$\mu\in M$ et~$s\in S(U)$. 
Soit~$g$ tel que~$D(g)$ soit un voisinage ouvert de~$x$ dans~$U$. Comme~$s$ ne s'annule
pas sur~$U$, on a~$D(g)\subset D(s)$ et donc~$D(g)=D(sg)$. La restriction de~$\sigma$ à~$D(sg)$
s'écrit~$\mu/s=\mu g/sg$. Ainsi, on a montré l'existence d'un élément~$f$ de~$A$ tel que~$x\in D(f)\subset U$
et tel que~$\sigma|_{D(f)}$ soit de la forme~$m/f$ avec~$m\in M$. 

\deux{def-mtilde}
On note~$\red M$ le faisceau associe au préfaisceau~$M_{\rm pref}$.
Le faisceau~$\red A$ est un faisceau de~$A$-algèbres, et~$\red M$ 
est un~$\red A$-module. 

\trois{fibres-redm}
Soit~$x\in \spec A$ et soit~$\got p$ l'idéal premier correspondant. La fibre~$\red M_x$ 
est canoniquement isomorphe à~$M_{{\rm pref}, x}$, c'est-à-dire à~$M_{\got p}$. 

\trois{fibres-reda}
En particulier, la fibre de~$\red A$ en~$x$ s'identifie à
l'anneau local~$A_{\got p}$. 

\trois{mor-mf-mtilde}
Soit~$f\in A$. Comme~$f\in S(D(f))$, on dispose d'une application~$A$-linéaire
naturelle~$M_f\to S(D(f))^{-1}M=M_{\rm pref}(D(f))$ (qui envoie~$m/f$ vu comme 
appartenant à~$M_f$ sur~$m/f$ vu comme appartenant à~$S(D(f))^{-1}M$). 
Par composition avec la flèche canonique~$M_{\rm pref}(D(f))\to D(f)$, 
on obtient une application~$A$-linéaire~$M_f\to \red M(D(f))$. 

\deux{annonce-faisceau-struct}
Le théorème qui suit est absolument fondamental. 
Il donne une description 
explicite de~$\red M$ qui
permet de manipuler effectivement ce dernier
(et qui s'appliquera en particulier à~$\red A$). 
Il est à la 
base de toute la théorie des schémas. 
Comme vous allez le voir, sa preuve met en jeu deux types d'ingrédients : 

\medskip
$\bullet$ les propriétés générales des localisés, et notamment la condition de nullité d'une fraction
(ou d'égalité de deux fractions) ; 

$\bullet$ une petite astuce de calcul, pas difficile mais absolument cruciale, qui permet
d'exhiber par une formule explicite un élément de~$M$ répondant à certaines conditions. 

\deux{theo-faisc-struct}
{\bf Théorème.}
{\em Pour tout~$f\in A$ l'application~$A$-linéaire~$M_f\to \red M(D(f))$
est bijective. En particulier, $M\to \red M(\spec A)$ est bijective (prendre~$f$
égal à~$1$).}

\medskip
{\em Démonstration.}
Nous allons commencer par une remarque qui permet de simplifier un peu la démonstration. Soit~$f\in A$
et soit~$U$ un ouvert de~$D(f)$. 
Le morphisme~$A\mapsto A_f$ identifie de manière naturelle~$\spec A_f$ à~$D(f)$. Modulo cette identification,
il résulte aisément des définitions et du fait que~$f\in S(U)$
que~$M_{\rm pref}(U)$
s'identifie à~$M_{f,\rm pref}(U)$, où~$M_{f,\rm pref}$
est le préfaisceau sur~$\spec A_f$ associé au~$A_f$-module~$M_f$. En conséquence, 
$\red M|_{D(f)}$ s'identifie au faisceau~$\red{M_f}$ sur~$\spec A_f$.

\medskip
Si l'on montre l'injectivité (resp. la surjectivité)
de~$M\to \red M(\spec A)$, ceci entraînera donc l'injectivité
(resp. la surjectivité)
de~$M_f\to \red M(D(f))$ pour tout~$f$, en appliquant le résultat
à l'anneau~$A_f$ en lieu et place de~$A$. 

\trois{structural-inj}
{\em Preuve de l'injectivité de~$M\to \red M(\spec A)$, et partant
de
l'injectivité~$M_f\to \red M(D(f))$ pour tout~$f$}. 
Soit~$m$ un élément de~$M$ dont l'image dans~$\red M(\spec A)$ est nulle. 
Les germes de cette
image sont alors tous nuls, ce qui signifie
que l'image de~$m$ dans~$M_{\got p}$ est nulle pour tout idéal premier
$\got p$ de~$M$ (\ref{fibres-redm}). En conséquence, $m=0$ 
d'après le lemme~\ref{lemme-null-msi}. 

\trois{structural-surj}
{\em Preuve de la surjectivité de~$M\to \red M(\spec A)$, et partant
de
la surjectivité de~$M_f\to \red M(D(f))$ pour tout~$f$}. 
Soit~$\sigma$
appartenant à~$\red M(\spec A)$ ; nous allons construire un élément~$m\in M$
d'image égale à~$\sigma$. 

\medskip
En vertu de~\ref{mpref-loc-mf},
de la quasi-compacité de~$\spec A$, 
et du fait que toute section de~$\red M$ provient
localement d'une section de~$M_{\rm pref}$, 
il existe une famille finie~$(f_i)_{i\in I}$ d'éléments de~$A$ 
et une famille~$(m_i)_{i\in I}$ d'éléments de~$M$ tels
que les propriétés suivantes soient satisfaites : 

\medskip
a) $\spec A=\bigcup D(f_i)$ ; 

b) pour tout~$i$, la restriction~$\sigma|_{D(f_i)}$
provient de la section~$m_i/f_i$ de~$M_{\rm pref}(D(f_i))$. 

\medskip
Soient~$i$ et~$j$ deux indices. Par construction, la restriction~$\sigma|_{D(f_if_j)}$
provient de la section~$m_i/f_i=f_jm_i/(f_if_j)$ de~$M_{\rm pref}(D(f_if_j))$,
et également de la section~$f_im_j/(f_if_j)$ de ce dernier. 
La flèche~$M_{f_if_j}\to \red M(D(f_if_j))$
est injective d'après ce qu'on a déjà montré au~\ref{structural-inj} ; en conséquence,
on a~$f_im_j/(f_if_j)=f_jm_i/(f_if_j)$ dans~$M_{(f_if_j)}$. Cela signifie qu'il existe un entier~$N$
tel que~$(f_if_j)^N(f_jm_i-f_im_j)=0$ ; notons que comme l'ensemble d'indices est fini,
l'entier~$N$ peut être choisi indépendamment de~$i$ et~$j$. 

\medskip
Nous allons simplifier un peu ces égalités. On commence par remarquer
que~$$(f_if_j)^N(f_jm_i-f_im_j)=0$$
$$\iff f_i^Nf_j^{N+1}m_i=f_j^Nf_{i}^{N+1}m_j,$$
qui n'est autre que l'égalité classique des produits en croix
pour les fractions~$f_i^Nm_i/(f_i^{N+1})$ et~$f_j^Nm_j/(f_j^{N+1})$. 

On remplace alors pour tout~$i$
la fonction~$f_i$ par~$f_i^{N+1}$ (ce qui ne change
pas~$D(f_i)$) et l'élément $m_i$ par~$f_i^Nm_i$. Les conditions
a) et~b) ci-dessus restent vérifiées, et l'on a de plus
$f_im_j-f_jm_i=0$ pour tout~$(i,j)$. 

\medskip
Nous en venons au cœur de la preuve : la construction d'un antécédent
de~$\sigma$. Remarquons pour commencer que si~$m\in M$, le
fait que son image dans~$\red M(\spec A)$ 
soit égale à~$\sigma$ peut se tester sur chacun des ouverts~$D(f_i)$ 
(puisque~$\red M$ est un faisceau). 
On en déduit, en considérant pour tout~$i$
le diagramme commutatif
$$\xymatrix{
M\ar[r]\ar[d]&{\red M(\spec A)}\ar[d]
\\
{M_{f_i}}\ar[r]&{\red M(D(f_i))}
},$$
que pour que l'image de~$m$ dans~$\red M(\spec A)$ soit égale
à~$\sigma$, il suffit que l'image de~$m$ dans~$M_{f_i}$
soit égale à~$m_i/f_i$ pour tout~$i$. 

\medskip
Nous allons maintenant utiliser l'hypothèse
que les~$D(f_i)$ recouvrent~$\spec A$. Cela signifie
qu'il existe~$(a_i)\in A^I$ tel que~$\sum a_if_i=1$. 
Intervient alors l'astuce de calcul que nous avons
évoquée plus haut : on
{\em pose}
$m=\sum_{j\in I} a_j m_j$
et nous allons vérifier que l'image de~$m$ dans~$M_{f_i}$
est égale à~$m_i/f_i$ pour tout~$i$. 

\medskip
Fixons~$i$. On a
$$f_im=f_i\sum_j a_j m_j=\sum_j a_j(f_i m_j)$$
$$=\sum_j a_j(f_jm_i)\;\;\;\text{(~car}\;f_im_j=f_jm_i\;\;\text{pour tout}\;j\text{)}$$
$$=(\sum_j a_jf_j)m_i=m_i.$$
En conséquence, l'image de~$m$ dans~$M_{f_i}$ est égale à~$m_i/f_i$,
ce qui achève la démonstration.~$\Box$

\deux{speca-locann-maintenant}
Muni du faisceau d'anneaux~$\red A$, l'espace topologique~$\spec A$ devient un espace
{\em localement}
annelé, dont le faisceau structural~$\red A$
sera désormais noté~$\sch O_{\spec A}$. Les faits suivant sont
des reformulations de~\ref{fibres-reda}
et du théorème~\ref{theo-faisc-struct}
pour~$M=A$ (nous nous permettons
de noter les isomorphismes canoniques comme des égalités) : 

\medskip
$\bullet$ si~$x$ est un point de~$\spec A$ correspondant
 à l'idéal~$\got p$ alors~$\sch O_{\spec A,x}=A_{\got p}$ ; 
 
 $\bullet$ si~$f\in A$ alors~$\sch O_{\spec A}(D(f))=A_f$ ; en particulier, $\sch O_{\spec A}(\spec A)=A$. 
 
 \deux{compatible-notations}
 {\bf Compatibilité des notations.}
 Soit~$x\in \spec A$
 et soit~$\got p$ l'idéal premier qui lui
 correspond. 
 Nous faisons face
 {\em a priori}
 à un conflit de notations : on dispose
 d'un corps~$\kappa(x)$ et d'un morphisme~$f\mapsto f(x)$
 de~$\sch O_{\spec A}(\spec A)=A$ dans~$\kappa(x)$ fournis par la théorie
 générale des espaces
 localement annelés (\ref{eval-u}) ; et d'un
 corps~$\kappa(x)$ et d'un morphisme~$f\mapsto f(x)$
 de~$A$ dans~$\kappa(x)$ définis directement
 (\ref{spec-def}).  Rappelons
 en quoi ils consistent. 
  
 \trois{kappax-fx-locann}
 {\em Définitions dans le contexte des espaces localement annelés.}
 Le corps~$\kappa(x)$ est alors le corps
 résiduel de l'anneau local~$\sch O_{\spec A, x}$, c'est-à-dire
 $A_{\got p}/\got pA_{\got p}={\rm Frac}\;A/\got p$. 
 Et l'évaluation est la composée de la flèche naturelle
 $$\sch O_{\spec A}(\spec A)=A\to \sch O_{\spec A,x}=A_{\got p}$$ et de
la flèche quotient~$\sch O_{\spec A, x}\to \kappa(x)=A_{\got p}\to
A_{\got p}/\got pA_{\got p}={\rm Frac}\;A/\got p$. Elle s'identifie à la flèche canonique
de~$A$ dans~${\rm Frac}\;A/\got p$. 

\trois{kappax-fx-spec}
 {\em Définition directe.}
 On a~$\kappa(x)={\rm Frac}\;A/\got p$ et
 la flèche~$f\mapsto f(x)$ est la flèche canonique
 $A\to {\rm Frac}\;A/\got p$. 
 
 \trois{conclu-pasconflit}
 Tout va donc pour le mieux : il n'y a
 {\em a posteriori}
 plus de conflit. 
 
 \subsection*{Les schémas : définition et premières propriétés}

 \deux{def-schema}
 {\bf Définition.}
 Un
 {\em schéma}
 est un espace localement annelé~$(X,\sch O_X)$ qui est localement
 isomorphe au spectre d'un anneau. Cela signifie plus précisément que pour tout~$x\in X$
 il existe un voisinage ouvert~$U$ de~$x$ dans~$X$, un anneau~$A$ et un isomorphisme
 d'espaces localement annelés~$(U,\sch O_U)\simeq (\spec A,\sch O_{\spec A})$.

 \medskip
 Un {\em morphisme de schémas}
 de~$(Y,\sch O_Y)$ vers~$(X,\sch O_X)$ est simplement un
 morphisme d'espaces localement annelés de $(Y,\sch O_Y)$ vers~$(X,\sch O_X)$.
 En d'autres termes, on définit 
 la catégorie des schémas comme une sous-catégorie pleine
 de la catégorie des espaces annelés. 
 
 \deux{comment-defschem}
Considérer les espaces localement
 annelés localement isomorphes à des espaces
 d'un certain type préalablement fixé peut s'avérer utile
 dans bien d'autres contextes géométriques. 
 
 Par exemple, 
 il est loisible de définir la catégorie des variétés différentielles comme
 la sous-catégorie pleine de la catégorie des espaces localement
 annelés en~$\RR$-algèbres
 formée des espaces localement isomorphes
 à un ouvert~$U$ de~$\RR^n$ (pour un certain~$n$) 
 muni de son faisceau des fonctions~$\sch C^\infty$. Cela dit, 
 en pratique, les géomètres différentiels préfèrent utiliser des cartes et atlas ; 
 c'est essentiellement une affaire de goût, et nous laissons au lecteur le soin
 de vérifier l'équivalence des deux points de vue.

 \deux{pre-exem-schemas}
 {\bf Premiers exemples.}

\trois{speca-schem}
Si~$A$ est un anneau, $\spec A$ est un schéma par définition.  Un schéma
qui est isomorphe à~$\spec A$ pour un certain~$A$ sera qualifié
d'{\em affine}.  

\trois{df-schemaff}
Soit~$A$ un anneau et soit~$f\in A$. On a signalé
au début de la preuve du théorème~\ref{theo-faisc-struct}
que la restriction de~$\sch O_{\spec A}=\red A$
à l'ouvert~$D(f)$ s'identifie, modulo l'homéomorphisme
$\spec A_f\simeq D(f)$ induit par~$A\to A_f$, au faisceau
$\red{A_f}=\sch O_{\spec A_f}$. En conséquence, 
l'ouvert~$D(f)$ de l'espace localement annelé~$\spec A$ est un schéma, 
et même un schéma affine : il est canoniquement isomorphe
à~$\spec A_f$. 

\medskip
Il s'ensuit que~$\spec A$ possède une base d'ouverts
qui sont des schémas affines, que l'on qualifiera plus
simplement d'{\em ouverts affines} ; ce fait s'étend immédiatement 
à un schéma quelconque. Il en résulte que tout ouvert d'un schéma
est un schéma. 

\deux{imm-ouverte-schem}
{\bf Les immersions ouvertes.}

\trois{prop-univ-schemeouv}
{\em La propriété universelle d'un ouvert.}
Les faits suivants se déduisent de~\ref{imm-ouverte-esplocann}.  

\medskip
Si~$U$ est un ouvert d'un schéma~$X$, l'inclusion
$j :U\hookrightarrow X$ est sous-jacente à un morphisme
de schémas noté encore~$j$, tel que pour tout schéma~$Y$ et tout 
morphisme~$\psi : Y\to X$ vérifiant~$\psi(Y)\subset U$, il existe
un unique morphisme de schémas~$\chi : Y\to U$ tel que~$j\circ \chi=\psi$. 
Le couple~$(U,j)$ représente donc le foncteur qui envoie un schéma~$Y$
sur l'ensemble des morphismes de~$Y\to X$ dont l'image est contenue
dans~$U$. 

\medskip
En d'autres termes, toute factorisation
{\em ensembliste}
d'un morphisme de schémas par un ouvert est automatiquement
{\em morphique.}

\medskip
Si nous avons insisté sur cette propriété, qui semble tellement évidente qu'on omet 
souvent de l'expliciter, c'est parce qu'elle
ne vaut pas sauf exceptions 
pour les
{\em fermés}
d'un schéma. 
Plus précisément, 
nous verrons plus loin
qu'un fermé~$Z$ d'un schéma~$X$ admet toujours
{\em au moins une} 
structure naturelle de schéma, mais qu'on ne peut pas en général
en trouver une
par laquelle se factorise tout morphisme~$\psi : Y\to X$
vérifiant~$\psi(Y)\subset Z$. 

\trois{def-immouverte-schem}
On dira qu'un morphisme de schémas~$Y\to X$ est une
{\em immersion ouverte}
s'il induit un isomorphisme entre~$Y$ et un ouvert de~$X$. 

\deux{mor-k-schem}
Soit~$X$ un schéma, soit~$K$ un corps et soit~$\xi$
l'unique point de~$\spec K$. Soit~$x\in X$ et soit~$\lambda$ un 
plongement de~$\kappa(x)$ dans~$K$. Nous allons
montrer qu'il existe un unique
morphisme~$\psi: \spec K\to X$
tel que~$\psi(\xi)=x$
et tel que la flèche induite~$\kappa(x)\hookrightarrow \kappa(\xi)=K$
soit égale à~$\lambda$. 

\trois{unique-psi-Kx}
{\em Unicité de~$\psi$.}
Soit~$\psi$ comme ci-dessus. Comme~$\psi(\xi)=x$,
le morphisme~$\psi$ est uniquement déterminé
ensemblistement. Soit~$U$ un ouvert de~$X$.
Il reste à s'assurer que
la flèche
$\psi^*: \sch O_X(U)\to \sch O_{\spec K}(\psi^{-1}(U))$
induite par~$\psi$ est elle aussi uniquement déterminée
par le couple~$(x,\lambda)$. 

\medskip
Si~$x\notin U$
on a~$\psi^{-1}(U)=\varnothing$  
et~$\psi^*(f)=0$ pour toute~$f\in \sch O_X(U)$. 

Si~$x\in U$ alors~$\psi^{-1}(U)=\{\xi\}$ et
on déduit du diagramme commutatif

$$\xymatrix{
{\sch O_{\spec K}(\psi^{-1}(U))=K}\ar[rrr]^{\rm Id=(g\mapsto g(\xi))}&&&
{K=\kappa(\xi)}
\\
{\sch O_X(U)}\ar[u]^{\psi^*}\ar[rrr]^{f\mapsto f(x)}&&&
\kappa(x)\ar@{^{(}->}[u]_{\lambda}
}$$
que~$\psi^*(f)= \lambda(f(x))$
pour toute~$f\in \sch O_X(U)$,
d'où notre assertion. 

\trois{const-psi}
{\em Existence de~$\psi$.} 
On s'inspire des égalités dont on vient
de voir qu'elles sont {\em nécessairement}
vérifiées. 

On décrit
tout d'abord~$\psi$
ensemblistement en posant~$\psi(\xi)=x$. 
Soit maintenant~$U$ un ouvert de~$X$. Nous définissons
un morphisme~$\psi^*$
de~$\sch O_X(U)$
vers~$\sch O_{\spec K}(\psi^{-1}(U))$ comme suit : 

\medskip
$\bullet$ si~$x\notin U$ on a~$\psi^{-1}(U)=\varnothing$ et l'on pose~$\psi^*(f)=0$ 
pour tout~$f\in \sch O_X(U)$  ; 

$\bullet$ si~$x\in U$ on a~$\psi^{-1}(U)=\spec K$, et l'on pose
$$\psi^*(f)=\lambda(f(x))\in K=\sch O_{\spec K}(\spec K).$$
Ces formules étant compatibles
aux restrictions, on obtient ainsi un morphisme de~$\spec K$
vers~$X$, dont il est immédiat qu'il satisfait les conditions requises. 

\deux{mor-Kx-conclu}
En vertu de ce qui précède, 
se donner un
morphisme de~$\spec K$ vers~$X$ revient à se donner
un point~$x$ de~$X$ et un plongement~$\kappa(x)\hookrightarrow K$. 

\medskip
En particulier, à tout point~$x$ de~$X$ est associé
un morphisme canonique~$\spec \kappa(x)\to X$ : celui
qui correspond au couple~$(x,{\rm Id}_{\kappa(x)})$. 

\deux{top-schemas-intro}
{\bf Quelques propriétés topologiques
des schémas.}
Nous allons nous contenter d'énoncés très généraux -- il est difficile d'être
plus précis
sans hypothèses spécifiques sur
les schémas en jeu. Le lemme~\ref{lemme-x-irred}
et le~\ref{irred-points-schemgen}
étendent des résultats précédemment démontrés pour les
spectres d'anneaux (prop.~\ref{prop-irred-sp}
et~\ref{bij-points-fermirr}). 

\trois{q-comp-schem}
Comme tout schéma affine est quasi-compact 
(en tant qu'espace topologique), un schéma
est quasi-compact si et seulement
si il est réunion {\em finie}
d'ouverts affines. 

\trois{ex-schem-pasqcomp}
Insistons sur le fait
qu'un schéma n'a aucune raison d'être quasi-compact en général. Par exemple, 
soit~$(X_i)_{i\in I}$ une famille quelconque de schémas. La somme disjointe~$\coprod X_i$
des~$X_i$ dans la catégorie des espaces localement annelés est un schéma : cette propriété
est en effet par définition locale, et~$\coprod X_i$ est recouvert par ses ouverts~$X_i$
qui sont des schémas. Si les~$X_i$ sont tous non vides et si l'ensemble~$I$ est infini, $\coprod X_i$
n'est pas quasi-compact. 

\medskip
Notons que comme les schémas forment une sous-catégorie pleine de la catégorie des espaces
localement annelé, $\coprod X_i$ est la somme disjointe des~$X_i$ {\em dans la catégorie des schémas}. 

\medskip
{\em Remarque :}
nous invitons le lecteur à vérifier
que le raisonnement suivi 
aux paragraphes~\ref{spectre-prod}
{\em et sq.}
montre en réalité l'existence d'un {\em isomorphisme
de schémas}
$$\spec (A\times B)\simeq \spec A \coprod \spec B.$$

\trois{lemme-x-irred}
{\bf Lemme.}
{\em Soit~$X$
un schéma et soit~$Y$ un fermé de~$X$. 
Les assertions suivantes sont équivalentes : 

\medskip
i) $Y$ possède un unique point générique ; 

ii) $Y$ possède un point générique ; 

iii) $Y$ est irréductible.}

\medskip
{\em Démonstration.}
Il est clair que~i)$\Rightarrow$ii)$\Rightarrow$iii). Supposons
maintenant que~iii)
soit vraie. Comme~$Y$ est irréductible,
il est non vide, et il existe donc un
ouvert affine $U$
de~$X$ tel que~$U\cap Y$ soit non vide. L'espace~$Y$
étant irréductible, il en va de même de son ouvert non vide~$U\cap Y$.
Comme celui-ci est par ailleurs un fermé du schéma affine~$U$, 
il est lui-même homéomorphe à un schéma affine 
et possède donc un unique point générique~$\eta$ (prop.~\ref{prop-irred-sp}, {\em cf.}
aussi~\ref{bij-points-fermirr}). Le point~$\eta$ est dense dans~$U\cap Y$,
lequel est dense dans~$Y$ par irréductibilité de ce dernier. En conséquence, $\eta$ est dense
dans~$Y$ : c'en est un point générique. 

\medskip
Il reste à s'assurer de l'unicité de~$\eta$. Soit~$\xi$ un point
générique de~$Y$. La densité de~$\xi$ dans~$Y$ signifie que~$\xi$ appartient à tout ouvert non vide de~$Y$,
et en particulier à~$Y\cap U$. Le point~$\xi$ qui est dense dans~$Y$
l'est {\em a fortiori}
dans~$Y\cap U$ ; par unicité du point générique de~$U$, il vient alors~$\xi=\eta$.~$\Box$ 

\trois{irred-points-schemgen}
On en déduit que~$x\mapsto \overline{\{x\}}$ induit une bijection entre~$X$
et l'ensemble de ses fermés irréductibles. 

\trois{rem-ferme-schema}
{\em Remarque.}
Nous avons énoncé et démontré le lemme~\ref{lemme-x-irred}
ci-dessus pour tout fermé de~$X$, et pas seulement 
pour le schéma~$X$ lui-même. Mais cette généralité est en réalité illusoire : en effet,
comme nous l'avons déjà mentionné au~\ref{prop-univ-schemeouv}, nous verrons
plus bas que tout fermé d'un schéma possède (au moins) une structure de schéma. 

\subsection*{Morphismes vers le spectre d'un anneau}

\deux{espann-speca}
{\bf Lemme.}
{\em Soit~$X$ un espace annelé et soit~$A$ un anneau. Pour tout morphisme
$\phi : A\to \sch O_X(X)$ et toute application continue~$\chi \colon X\to \spec A$
il existe {\em au plus un}
morphisme d'espaces 
annelés~$\psi \colon X\to \spec A$ dont l'application continue
sous-jacente coïncide avec~$\chi$ et tel que le morphisme
$$A=\sch O_{\spec A}(\spec A)\to \sch O_X(X)$$ induit par~$\psi$ soit égal à~$\phi$.}

\medskip
{\em Démonstration}.
Soit~$\psi \colon X\to \spec A$
un morphisme
d'espaces annelés comme dans l'énoncé. 
Comme les ouverts de la forme~$D(f)$
forment une base de la topologie 
de~$\spec A$, il suffit de s'assurer que pour tout~$f$
appartenant à~$A$,
la flèche
$$A_f=\sch O_{\spec A}(D(f))\to \sch O_X(\chi^{-1}(D(f)))$$
induite par~$\psi$
est uniquement
déterminé par~$\phi$. 

Soit~$f\in A$. Posons~$U=\chi^{-1}(D(f))$. 
Soit~$\rho\colon A\to \sch O_X(U)$
la composée de~$\phi \colon A\to \sch O_X(X)$ et de la
flèche de restriction
$\sch O_X(X)\to \sch O_X(U)$. 
Le diagramme
commutatif 

$$\xymatrix{
{\sch O_X(X)}\ar[rr]&&{\sch O_X(U)}
\\
A\ar[u]^\phi\ar[rr]\ar[urr]^\rho&&A_f\ar[u]
}$$
montre que la flèche
composée~$A\to A_f\to \sch O_X(U)$ est 
égale à~$\rho$. Mais la propriété universelle du localisé~$A_f$ assure qu'il y a {\em au plus}
un morphisme de~$A_f$ vers~$\sch O_X(U)$ dont la composée avec~$A\to \sch O_X(U)$ est égale
à~$\rho$ (elle assure aussi qu'un tel morphisme existe si et seulement si~$\rho(f)$ est inversible) ;
ainsi, le morphisme $A_f\to \sch O_X(U)$ est uniquement déterminé
par~$\phi$, ce qu'on souhaitait établir.~$\Box$

\deux{locann-speca}
{\bf Théorème.}
{\em Soit~$X$ un espace localement annelé et soit~$A$ un anneau. Pour tout morphisme
$\phi : A\to \sch O_X(X)$ il existe un unique morphisme d'espaces localement
annelés~$\psi  : X\to \spec A$ tel que le morphisme
$$A=\sch O_{\spec A}(\spec A)\to \sch O_X(X)$$ induit par~$\psi$ soit égal à~$\phi$.}

\medskip
{\em Démonstration.} Nous allons, comme souvent, commencer par
établir l'unicité ; puis nous nous inspirerons des conditions nécessaires qui auront
été dégagées à cette occasion pour exhiber un morphisme satisfaisant les conditions requises. 

\trois{locann-speca-unique}
{\em Preuve de l'unicité 
de~$\psi$.}
Soit~$\psi$ comme dans l'énoncé. Par hypothèse, 
il induit le morphisme~$\phi$ entre les anneaux de sections globales. 
Pour montrer qu'il est uniquement déterminé il suffit alors, en vertu du lemme~\ref{espann-speca}, 
de montrer qu'il est uniquement déterminé {\em ensemblistement}. 

\medskip
Soit~$x\in X$ et soit~$f\in A$. L'image de~$f$ dans~$\sch O_X(X)$
est par hypothèse égale
à
$\phi(f)$. Comme~$\psi$ est un morphisme d'espaces localement
annelés, on a
$$\phi(f)(x)=0\iff f(\psi(x))=0.$$
Ainsi, $\psi(x)$ est {\em nécessairement}
le point correspondant au noyau de~$f\mapsto \phi(f)(x)$
(ou encore au morphisme~$f\mapsto \phi(f)(x)$ de~$A$ dans~$\kappa(x)$).

\trois{locann-speca-exist}
{\em Existence de~$\psi$}. Commençons par le définir ensemblistement : pour tout~$x\in X$, 
on note~$\psi(x)$ le point de~$\spec A$ correspondant au noyau de
la flèche~$f\mapsto \phi(f)(x)$
((ou encore au morphisme~$f\mapsto \phi(f)(x)$ de~$A$ dans~$\kappa(x)$).
Cela signifie
que l'on a pour tout~$f\in A$ l'équivalence $$\phi(f)(x)=0\iff f(\psi(x))=0.$$ Il s'ensuit
que~$\psi^{-1}(D(f))=D(\phi(f))$, et~$\psi$ est en conséquence continue. 

\medskip
Pour faire de~$\psi$ un morphisme d'espaces annelés, il faut maintenant 
se donner un morphisme~$\sch O_{\spec A}\to \psi_*\sch O_X$ ; nous reprenons les notations
de~\ref{def-mpref}
{\em et sq.}

Soit~$U$ un ouvert de~$\spec A$. Si~$s\in S(U)$, il résulte de la définition de~$\psi$
que~$\phi(s)$ est inversible sur~$\psi^{-1}(U)$. La flèche composée

$$\xymatrix
{A\ar[rr]^\phi&&{\sch O_X(X)}\ar[r]&{\sch O_X(\psi^{-1}(U))}}$$
envoie donc~$S$ dans l'ensemble des éléments inversibles
de~$\sch O_X(\psi^{-1}(U))$. Elle induit dès
lors un morphisme
de~$S(U)^{-1}A$ vers~$\sch O_X(\psi^{-1}(U))$. Cette construction
étant compatible
aux restrictions, elle définit un morphisme
de préfaisceaux d'anneaux de~$A_{\rm pref}$
vers~$\psi_*\sch O_X$. Comme~$\psi_*\sch O_X$
est un faisceau,
ce morphisme induit 
un morphisme de faisceau d'anneaux 
de~$\red A=\sch O_{\spec A}$ vers~$\psi_*\sch O_X$
dont l'effet sur les sections globales coïncide par construction
avec~$\phi$. 

\medskip
Il reste à s'assurer que le morphisme d'espaces 
annelés~$\psi$ est bien un morphisme d'espaces {\em localement}
annelés. Soit~$x\in X$ et soit~$g$ une section de~$\sch O_{\spec A}$
définie sur un voisinage de~$\psi(x)$. On peut toujours supposer que
ce voisinage est de la forme~$D(f)$ avec~$f(\psi(x))\neq 0$, auquel cas
la fonction~$g$ est de la forme
$a/f$ ; on a~$\psi^{-1}(D(f))=D(\phi(f))$, et le diagramme
commutatif
$$\xymatrix{
{\sch O_X(X)}\ar[rr]&&{\sch O_X(D(\phi(f)))}\\
A\ar[u]^\phi\ar[rr]&&A_f\ar[u]}$$
assure que l'image de~$g$ dans~$\sch O_X(\psi^{-1}(D(f)))$
est égale à~$\phi(a)/\phi(f)$ On a les équivalences
$$g(\psi(x))=0\iff a(\psi(x))=0\iff \phi(a)(x)=0\iff [\phi(a)/\phi(f)](x)=0,$$
ce qu'il fallait démontrer.~$\Box$

\deux{imm-ouv-specann}
{\bf Exemple.}
Soit~$A$ un anneau, soit~$U$ un ouvert
de~$\spec A$, et soit~$\rho$ la flèche 
de restriction de~$A$ vers~$\sch O_U(U)$. 
Le morphisme induit
par l'immersion ouverte~$U\hookrightarrow A$ 
au niveau des anneaux de sections globales est égal à~$\rho$. 
En conséquence, la
flèche~$U\to \spec A$ associée par le théorème~\ref{locann-speca}
au morphisme~$\rho$ coïncide
avec l'immersion ouverte~$U\hookrightarrow \spec A$. 

\deux{apll-spectres}
Soit~$\phi  \colon A \to B$
un morphisme d'anneaux. 
Le théorème~\ref{locann-speca}
permet d'associer à~$\phi$
un morphisme~$\psi \colon \spec B\to \spec A$,
à savoir l'unique morphisme
qui induit~$\phi$ entre leurs anneaux de sections globales. 
Les faits suivants découlent de la description explicite
de~$\psi$, donnée en~\ref{locann-speca-exist}
lors 
de la preuve de
{\em loc. cit.}

\trois{deja-defini}
Ensemblistement, $\psi$
coïncide avec
l'application continue construite au~\ref{spec-fonct} :
si~$\got q$ un idéal premier de~$B$, et si~$y$ est
le point correspondant de~$\spec B$, alors~$\psi(y)$ correspond à l'idéal premier~$\phi^{-1}(\got q)$ de~$A$.

\trois{psi-def-attendue}
Soit~$f\in F$. On a~$\psi^{-1}(D(f))=D(\phi(f))$
et le morphisme
$$[A_f=\sch O_{\spec A}(D(f))]\to [B_{\phi(f)}=\sch O_{\spec B}(D(\phi(f))]$$
induit par~$\psi$ est le morphisme canonique $A_f\to B_{\phi(f)}$ déduit de~$\phi$. 

\trois{psi-def-fibres}
Soit~$\got q$ un idéal premier de~$B$, et soit~$y$ le point correspondant de~$\spec B$. 
Le morphisme
$$[A_{\phi^{-1}(\got q)}=\sch O_{\spec A,\psi(y)}] \to[ B_{\got q}=\sch O_{\spec B,y}]$$
induit par~$\psi$ est le morphisme canonique $A_{\phi^{-1}(\got q)}\to B_{\got q}$ déduit de~$\phi$. 

\deux{cas-identi}
Soit~$A$ un anneau. Le morphisme~${\rm Id}_{\spec A}$ induit
l'identité sur l'anneau des sections globales de~$\spec A$ ; c'est
donc le morphisme
associé à~${\rm Id}_A$.

\medskip
Soient~$\phi : A\to B$ et~$\phi' : B\to C$  deux morphismes
d'anneaux. La composée des morphismes~$\spec C\to \spec B$
et~$\spec B\to \spec A$ respectivement associés à~$\phi'$ et~$\phi$
est un morphisme de~$\spec C$ vers~$\spec A$ qui induit
le morphisme~$\phi'\circ \phi$ entre leurs anneaux de
sections globales : c'est donc
le morphisme
de~$\spec C$
vers~$\spec A$ associé à~$\phi'\circ \phi$. 

\medskip
Notons que ces fait pourraient aussi se déduire des descriptions
explicites des morphismes évoqués, 
que nous avons fournies en~\ref{deja-defini}--\ref{psi-def-fibres}. 

\deux{fonct-spec-locann}
Ainsi, $A\mapsto \spec A$ apparaît comme un foncteur
contravariant de la catégorie des anneaux vers celle des espaces localement annelés,
qui induit pour tout couple~$(A,B)$ d'anneaux une 
bijection
$$\hom_{\ann}(A,B)\simeq \hom_{\mathsf{Sch}}(\spec B,\spec A).$$
Il s'ensuit (exercice~\ref{equiv-pf-es})
que ce foncteur établit une anti-équivalence\footnote
{Le préfixe «anti»
fait référence au fait que c'est un fonceur contravariant.}
entre la catégorie des anneaux et celle des schémas
{\em affines}. On vérifie immédiatement que~$X\mapsto O_X(X)$
en est un quasi-inverse. 

\deux{interp-fonct-locann}
Nous allons maintenant pouvoir donner diverses interprétations
fonctorielles du théorème~\ref{locann-speca}. 

\trois{reformul-locann-speca}
Soit~$X$ un espace localement annelé. Le théorème~\ref{locann-speca}
assure en particulier l'existence d'un unique morphisme~$\chi : X\to \spec \sch O_X(X)$
induisant~${\rm Id}_{\sch O_X(X)}$
sur les anneaux de sections globales, dont nous dirons que c'est le morphisme
{\em canonique}
de~$X$ vers~$\spec \sch O_X(X)$ (notons que si~$X=\spec A$ on a~$\sch O_X(X)=A$
et~$\chi$ est alors nécessairement égal à l'identité, par unicité). 

\medskip
Soit~$\phi$ un morphisme
de~$A$ vers~$\sch O_X(X)$. Il induit un morphisme
$\theta$
de~$\spec \sch O_X(X)$ vers~$\spec A$ ;
la composée~$\theta \circ \chi$
est une flèche
de~$X$
vers~$\spec A$ qui induit le morphisme~$\phi$
sur les anneaux de sections globales : c'est donc le morphisme~$\psi$
dont le théorème~\ref{locann-speca}
assure l'existence et l'unicité. 

\medskip
Ainsi, $\psi$ se factorise
{\em via}
$\chi$. Cette factorisation est unique : si~$\theta'$ est un morphisme
de~$\spec \sch O_X(X)$ vers~$\spec A$ tel que~$\theta'\circ \chi=\psi$ on
vérifie aussitôt que~$\theta'$ induit le morphisme~$\phi$
sur les anneaux de sections globales, et il est dès lors égal à~$\theta$. 

\medskip
En conséquence, le foncteur covariant 
de la catégorie des schémas affines vers les ensembles
qui envoie~$Y$ sur~$\hom_{\mathsf{Sch}}(X,Y)$ est représentable
par le couple~$(\spec \sch O_X(X), \chi \colon X\to \spec \sch O_X(X))$. 

\trois{reformul-locann-adj}
Le théorème~\ref{locann-speca}
affirme que pour tout couple~$(X,A)$ formé d'un espace
localement annelé et d'un anneau, l'application naturelle naturelle
$\hom_{\mathsf{Esp-loc-ann}}(X,\spec A)\to \hom_{\ann}(A,\sch O_X(X))$
est bijective. On vérifie immédiatement que cette bijection est fonctorielle en~$X$ et~$A$. 

\medskip
Autrement dit, {\em $A\mapsto \spec A$ est adjoint
au foncteur~$X\mapsto \sch O_X(X)$ qui va des espaces localement annelés vers les anneaux}. 

\medskip
Notons que les foncteurs en jeu sont
contravariants ; c'est pour cela que nous n'avons
pas précisé «adjoint à droite»
ou «adjoint à gauche», aucune des deux terminologies ne semblant plus naturelle que l'autre. 

\deux{conclu-spec-esplocann}
{\bf La place des spectres au sein des espaces localement annelés.} 
Il résulte du~\ref{reformul-locann-adj}
ci-dessus que la notion de spectre est naturelle
dès lors qu'on s'intéresse aux espaces localement annelés généraux,
puisque~$A\mapsto \spec A$ est adjoint à~$X\mapsto \sch O_X(X)$. 

Ce fait permet en un sens de penser à~$\spec A$ comme à {\em l'espace
localement annelé le plus général d'anneau des sections globales égal à~$A$}. 

\deux{comment-adj-local}
{\bf Remarque culturelle.}
La notion de spectre, et sa caractérisation fonctorielle
donné au~\ref{conclu-spec-esplocann}
ci-dessus, se généralisent comme suit. Le
foncteur d'inclusion de la catégorie des espaces localement annelés dans celle
des espaces annelés admet un adjoint à droite~$\mathfrak {Sp}$, qui est une sorte
de «spectre étalé» ; et si~$A$ est un anneau, $\spec A$ s'identifie
à~$\mathfrak {Sp}(\{*\},\underline A)$. 
Nous ne sous servirons pas de~$\mathfrak{Sp}$, mais nous allons
indiquer sa construction pour le lecteur intéressé. 
Soit~$(X,\sch O_X)$ un espace annelé. 

\trois{ens-mathfraksp}
{\em Description ensembliste}. 
L'ensemble sous-jacent à $\mathfrak{Sp}(X)$
est l'union disjointe~$\coprod_{x\in X}\spec \sch O_{X,x}$, munie
d'une application naturelle~$p$ vers~$X$ (le spectre de~$\sch O_{X,x}$ est 
placé au-dessus de~$x$). 

\trois{top-mathfraksp}
{\em La topologie de~$\mathfrak{Sp}(X)$.}
Pour tout ouvert~$U$ de~$X$ et tout
élément~$f$ de~$\sch O_X(U)$,
on note~$\got D(U,f)$ le sous-ensemble de~$\mathfrak{Sp}(X)$ 
égal à la réunion, pour~$x$ parcourant~$U$, des ouverts~$D(f_x)\subset \spec \sch O_{X,x}=p^{-1}(x)$. 
On munit~$\mathfrak{Sp}(X)$ de la topologie engendrée par les~$\got D(U,f)$, pour
laquelle~$p$ est continue. 

\trois{locann-mathfraksp}
{\em Le faisceau structural de~$\mathfrak{Sp}(X)$.}
Soit~$V$ un ouvert de~$\mathfrak{Sp}(X)$. Pour tout
ouvert~$U$ de~$X$ contenant~$p(V)$, notons~$S_V(U)$ l'ensemble
des sections~$f$ de~$\sch O_X$ sur~$U$ telles que~$V\subset \got D(U,f)$. C'est une partie multiplicative 
de~$\sch O_X(U)$  ; soit~$\Lambda(V)$ la limite inductive des~$S_V(U)^{-1}\sch O_X(U)$,
où~$U$ parcourt l'ensemble des ouverts de~$X$ contenant~$p(V)$. 

\medskip
On munit alors~$\mathfrak{Sp}(X)$
du faisceau associé au préfaisceau~$V\mapsto \Lambda(V)$. On vérifie que~$\mathfrak{Sp}(X)$
est localement annelé, que~$p\colon \mathfrak {Sp}(X)\to X$
est un morphisme d'espaces annelés, et que pour tout espace localement annelé~$Y$, l'application
$\phi\mapsto p\circ \phi$ établit une bijection
entre~$\hom_{\mathsf{Esp-loc-ann}}(Y,\mathfrak{Sp}(X))$ et~$\hom_{\mathsf{Esp-ann}}(Y,X)$.

\subsection*{Un critère d'affinité, et un premier
contre-exemple}

\deux{intro-etre-affine}
{\em A priori}, le fait d'être affine ne semble
pas être une propriété facile à vérifier pour un schéma donné~$X$ : cela
signifie en effet qu'{\em il existe}
un anneau~$A$ et un isomorphisme~$X\simeq \spec A$.

Mais le lemme ci-dessous assure qu'en réalité,
cela peut s'exprimer de façon directe, 
sans quantificateur existentiel toujours un peu désagréable. 

\deux{critere-affine}
{\bf Lemme.}
{\em Soit~$X$ un schéma. Il est affine si et seulement 
si le morphisme canonique~$X\to \spec \sch O_X(X)$ est un isomorphisme.}

\medskip
{\em Démonstration.}
Si~$X\to \spec \sch O_X(X)$ est un isomorphisme, $X$ est affine par définition. Supposons
réciproquement que~$X$ soit affine, donc qu'il existe un isomorphisme~$X\simeq \spec A$. 
On dispose alors d'un diagramme
commutatif
$$\xymatrix{
X\ar[rrr]\ar[d]_\simeq &&&{\spec \sch O_X(X)}\ar[d]^\simeq\\
{\spec A}\ar[rrr]^(.35){\rm Id}&&&{\spec A=\spec \sch O_{\spec A}(\spec A)}
}$$
qui montre que~$X\to \spec \sch O_X(X)$ est un isomorphisme.~$\Box$

\deux{ex-non-affine}
{\bf Un premier exemple de schéma non affine.}
Soit~$k$ un corps. Nous allons introduire une notation
que nous utiliserons dans toute la suite du cours : on désigne
par~$\Aff^n_k$ 
le schéma~$\spec k[T_1,\ldots, T_n]$. Il est muni d'un morphisme naturel
vers~$\spec k$, induit par le plongement~$k\hookrightarrow k[T_1,\ldots, T_n]$ ; on dit
que~$\Aff^n_k$ est {\em l'espace affine de dimension~$n$ sur le corps~$k$.}

\medskip
Pour le moment, nous allons travailler avec
le plan affine~$\Aff^2_k$. Soit~$U$ l'ouvert complémentaire
dans~$\Aff^2_k$
de l'origine~$(0,0)$, vue comme point fermé de~$\Aff^2_k$ . 

\trois{oxu-polyn}
{\em Détermination de~$\sch O_{\Aff^2_k}(U)$.}
Comme~$(0,0)=V(T_1,T_2)$, l'ouvert~$U$
 est la réunion de~$D(T_1)$ et~$D(T_2)$. On a

$$\sch O_{\Aff^2_k}(D(T_1))=k[T_1,T_2]_{T_1}=\left\{\frac P{T_1^n}, P\in k[T_1,T_2], n\in \NN\right\}\subset k(T_1,T_2),$$
$$\sch O_{\Aff^2_k}(D(T_2))=k[T_1,T_2]_{T_2}=\left\{\frac P{T_2^n}, P\in k[T_1,T_2], n\in \NN\right\}\subset k(T_1,T_2)$$
et
$$\sch O_X(D(T_1)\cap D(T_2))=\sch O_X(D(T_1T_2))=k[T_1,T_2]_{T_1T_2}$$
$$=\left\{\frac P{T_1^nT_2^n}, P\in k[T_1,T_2], n\in \NN\right\}\subset k(T_1,T_2).$$

Se donner une section de~$\sch O_{\Aff^2_k}$ sur~$U$ revient à se donner un couple~$(t_1,t_2)$ formé d'une section~$t_1$ de~$\sch O_{\Aff^2_k}$
sur~$D(T_1)$ et d'une section~$t_2$ de~$\sch O_{\Aff^2_k}$ sur~$D(T_2)$, telles que~$t_1|_{D(T_1T_2)}=t_2|_{D(T_1T_2)}$. Par ce qui précède, cela revient à se donner
une fraction~$R\in K(T_1,T_2)$ pouvant à la fois s'écrire sous la forme~$P/T_1^n$ et
sous la forme~$Q/T_2^m$ où~$P$ et~$Q$ appartiennent
à~$k[T_1,T_2]$. Un argument élémentaire d'arithmétique des anneaux factoriels assure qu'un tel~$R$ appartient nécessairement à~$k[T_1,T_2]=\sch O_{\Aff^2_k}(\Aff^2_k)$.
Autrement dit,
la flèche de restriction~$\sch O_{\Aff^2_k}(\Aff^2_k)\to \sch O_{\Aff^2_k}(U)$ est un isomorphisme. 

\trois{oxu-nonaff}
La flèche canonique
de~$U$ vers~$\spec \sch O_{\Aff^2_k}(U)$ s'identifie
par ce qui précède
(et en vertu de l'exemple~\ref{imm-ouv-specann})
à l'immersion ouverte
de~$U$ dans~$\Aff^2_k$,  qui n'est pas
un isomorphisme (elle n'est déjà pas ensemblistement surjective, puisque son image ne contient pas l'origine). 
Il en résulte que~$U$ n'est pas affine.

\section{Recollement de schémas, construction des
produits fibrés}
\markboth{La notion de schéma}{Recollements, produits fibrés}

\subsection*{Recollements de schémas}

\deux{x-schema}
Soit~$X$ un schéma. Un
{\em schéma sur~$X$}, ou encore un~{\em $X$-schéma}, 
est un couple~$(Y,\phi)$ où~$Y$ est un schéma
et~$\phi$ un morphisme de~$Y$ vers~$X$. Un morphisme
entre de~$X$-schémas de~$(Y,\phi)$ vers~$(Z,\psi)$ 
est un morphisme de~$Y$ vers~$Z$ tel que le diagramme
$$\xymatrix{
Y\ar[rr]\ar[dr]_\phi&&Z\ar[dl]^\psi\\
&X&}$$ commute. Bien entendu, le plus souvent on omettra
de mentionner explicitement~$\phi$ et~$\psi$, et l'on parlera
de {\em $X$-morphisme}
de~$Y$ vers~$Z$.

\trois{a-schema}
Si~$A$ est un anneau on dira le plus souvent~«$A$-schéma»
et~«$A$-morphisme»
au lieu de~«$\spec A$-schéma»
et~«$\spec A$-morphisme».

\deux{limind-schem-intro}
Soit~$\sch D=((X_i), (E_{ij})$ un diagramme
dans la catégorie des schémas, et soit~$X$ sa limite inductive
dans la catégorie des espaces localement annelés. 
{\em Si~$X$ 
a le bon goût d'être un schéma}, c'est {\em a fortiori}
la limite inductive de~$\sch D$ dans la catégorie 
des schémas (puisque celle-ci est une sous-catégorie pleine
de la catégorie des espaces localement annelés). 

\deux{limind-schem}
Nous allons maintenant donner quelques exemples
de diagrammes
dans la catégorie
des schémas dont la limite inductive dans~$\mathsf{Esp}\text{-}\mathsf{loc}\text{-}\mathsf{ann}$
est effectivement un schéma, et qui sont donc justiciables de la remarque précédente. 

\medskip
Dans ce paragraphe, «diagramme» signifiera «diagramme {\em dans la catégorie
des schémas}», et les limites inductives des diagrammes seront calculées dans la catégorie
des espaces localement annelés. 

\trois{somme-disj-schem-bis}
Si~$\sch D$ est un diagramme sans flèches~$(X_i)$, alors~$\limind \sch D=\coprod X_i$
est un schéma : on l'a déjà signalé en~\ref{ex-schem-pasqcomp}. 

\trois{diag-limind-autom-schem}
Soit~$I$ un ensemble ordonné
et 
soit~$\sch D$ un diagramme commutatif
induit par un foncteur
de~$I$ dans la catégorie des schémas (\ref{desc-idiagramme}). 
Pour tout
indice~$i$ appartenant à~$I$ on note~$X_i$ l'objet correspondant de~$\sch D$ ; pour tout
couple~$(i,j)\in I$ 
avec~$i\leq j$, on note~$f_{ij}$ la flèche~$X_i \to X_j$ de~$\sch D$.

\medskip
Nous faisons les hypothèses suivantes : 

i) les~$f_{ij}$ sont toutes des immersions ouvertes ; 

ii) pour tout~$(i,j)\in I$ et tout~$\ell$ majorant~$i$ et~$j$, 
l'ouvert~$f_{i\ell}(X_i)\cap f_{j\ell}(X_j)$ de~$X_{\ell}$ est la réunion 
des~$f_{\alpha\ell}(X_\alpha)$ pour~$\alpha$ minorant~$i$ et~$j$. 

\medskip
Nous laissons le lecteur vérifier les assertions ci-dessous (qui valent
d'ailleurs pour tout diagramme de ce type dans la catégorie des {\em espaces annelés}). 

1) Les flèches naturelles~$\lambda_i \colon X_i\to \limind \sch D$ sont des immersions ouvertes. 

2) Pour tout~$(i,j)\in I^2$, l'ouvert~$\lambda_i(X_i)\cap \lambda_j(X_j)$ de~$\limind \sch D$ est la réunion 
des~$\lambda_\alpha(X_\alpha)$ pour~$\alpha$ minorant~$i$ et~$j$.

\medskip
La famille~$(\lambda_i(X_i))$ est un recouvrement ouvert de~$\sch D$, 
et~$\lambda_i(X_i)$ est un schéma pour tout~$i$ (puisque c'est le cas de~$X_i$). 
L'espace localement annelé~$\limind \sch D$ est donc un schéma. On dira
qu'il est obtenu par
{\em recollement des~$X_i$ le long des~$f_{ij}$}. 

\trois{limind-recoll-schema-cocy}
Soit~$I$ un ensemble (quelconque) d'indices.  Pour tout~$i$, soit~$X_i$ un 
schéma. Pour tout couple~$(i,j)$ avec~$i\neq j$ on se donne un ouvert~$(X_{ij})$ de~$X_i$, et
un isomorphisme~$\iota_{ij}\colon X_{ij}\simeq X_{ji}$. On suppose que les~$\iota_{ij}$ satisfont 
les conditions suivantes : 

\medskip
i) $\iota_{ij}=\iota_{ji}^{-1}$ pour tout~$(i,j)$ avec~$i\neq j$ ; 

ii) si~$(i,j,k)$ sont trois indices deux à deux distincts alors
$$\iota_{ij}(X_{ik}\cap X_{ij})=X_{ji}\cap X_{jk}$$
et
$$\iota_{jk}\circ \iota_{ij}=\iota_{ik},$$
les deux membres étant vus comme des isomorphismes
de~$X_{ij}\cap X_{ik}$ sur~$X_{ki}\cap X_{kj}$. 

\medskip
Soit~$\sch D$ le diagramme dont
les objets sont les~$X_i$ et les~$X_{ij}$, et dont les flèches
sont les isomorphismes~$\iota_{ij}$ et les immersions ouvertes~$X_{ij}\hookrightarrow X_i$. 

\medskip
Nous laissons le lecteur vérifier les assertions ci-dessous (qui valent
d'ailleurs pour tout diagramme de ce type dans la catégorie des {\em espaces annelés}). 

1) Les flèches naturelles~$\lambda_i \colon X_i\to \limind \sch D$ sont des immersions ouvertes. 

2) Pour tout~$(i,j)\in I^2$ avec~$i\neq j$,
on a 
$$\lambda_i(X_{ij})=\lambda_j(X_{ji})=\lambda_i(X_i)\cap \lambda_j(X_j),$$
et~$\iota_{ij}=\lambda_j^{-1}\circ \lambda_i$, où~$\lambda_i$ et~$\lambda_j$ sont respectivement
vus dans cette égalité
comme des isomorphismes de~$X_{ij}$ et~$X_{ji}$ sur~$\lambda_i(X_i)\cap \lambda_j(X_j)$.  

\medskip
\medskip
La famille~$(\lambda_i(X_i))$ est un recouvrement ouvert de~$\sch D$, 
et~$\lambda_i(X_i)$ est un schéma pour tout~$i$ (puisque c'est le cas de~$X_i$). 
L'espace localement annelé~$\limind \sch D$ est donc un schéma. On dira
qu'il est obtenu par
{\em recollement des~$X_i$ le long des~$\iota_{ij}$}.

\trois{rapp-liminid-yschema}
{\em Remarque.}
Nous aurons plusieurs fois l'occasion d'appliquer ce qui précède, en le combinant
implicitement avec la remarque suivante. 

\medskip
Soit~$\sch D=((X_i), (E_{ij}))$ un diagramme dans la catégorie des schémas
et soit~$Y$ un schéma. Supposons que chaque~$X_i$ soit muni d'une structure de~$Y$-schéma,
et que les éléments de~$E_{ij}$ soient pour tout~$(i,j)$ des~$Y$-morphismes. Si $\limind \sch D$ 
existe dans la catégorie des schémas, elle hérite d'une structure naturelle de~$Y$-schéma est s'identifie
à la limite inductive de~$\sch D$ dans la catégorie
des~$Y$-schémas : c'est un fait complètement général, {\em cf.}~\ref{sur-s-ind}.  

\subsection*{La droite projective et la droite affine avec origine dédoublée}

\deux{intro-copies-a1}
Soit~$k$ un corps. Nous
allons travailler dans ce qui suit avec deux copies~$X$ et~$Y$ 
de la droite affine~$\Aff^1_k$, vue comme
un~$k$-schéma de façon évidente.
Pour éviter de les confondre, nous
écrirons~$X=\spec k[T]$ et~$Y=\spec k[S]$. On note~$U$ l'ouvert~$D(T)$ de~$X$,
et~$V$ l'ouvert~$D(S)$ de~$Y$. On a~$U\simeq \spec k[T,T^{-1}]$ et~$V\simeq \spec k[S,S^{-1}]$. 
On note~$i$ l'immersion ouverte de~$U$ dans~$X$, et~$j$ celle de~$V$ dans~$Y$.

\deux{droite-proj-recoll}
{\bf La droite projective.}
L'isomorphisme de~$k$-algèbres
$$\begin{array}{ccc}
k[S,S^{-1}]&\to& k[T,T^{-1}]\\
S&\mapsto &T^{-1}\\
S^{-1}&\mapsto &T\end{array}$$
induit un isomorphisme
de~$k$-schémas~$\psi : U\to V$. 

\medskip
On note~$\PP^1_k$ le~$k$-schéma
obtenu par recollement de~$X$ et~$Y$ le long
de~$\psi$ et~$\psi^{-1}$, défini en~\ref{limind-recoll-schema-cocy} -- notez que
la condition~ii) de
{\em loc cit.} est ici vide puisqu'on ne recolle que deux ouverts. Le~$k$-schéma
$\PP^1_k$ 
est également appelé
la~{\em droite projective sur~$k$.}
Nous allons maintenant la décrire plus avant, en déclinant dans ce
cas particulier les énoncés~1) et~2) de~{\em loc. cit.}

\trois{topo-droite-proj}
La droite projective~$\PP^1_k$ est réunion de deux ouverts affines~$X'$ et~$Y'$
(les images de~$X$ et~$Y$). Chacun d'eux est une copie de la droite affine : on 
a~$X'\simeq \spec k[T]$ et~$Y'\simeq k[S]$. Leur intersection~$X'\cap Y'$
est égale à~$D(T)$ en tant qu'ouvert de~$X'$, et donc à~$\spec k[T,T^{-1}]$ ; elle
est égale à~$D(S)$ en tant qu'ouvert de~$Y'$, et donc à~$\spec k[S,S^{-1}]$. L'isomorphisme
entre~$D(T)\subset X'$ et~$D(S)\subset Y'$ induit par leurs identifications avec~$X'\cap Y'$ 
est celui fourni par le morphisme de~$k$-algèbres qui envoie~$S$ sur~$T^{-1}$. 

\medskip
Le complémentaire de~$X'\cap Y'$ dans~$\PP^1_k$ est constitué de deux points fermés : l'origine~$V(T)$
de~$X'\simeq \Aff^1_k$, et l'origine~$V(S)$ de~$Y'\simeq \Aff^1_k$. Comme~$S=T^{-1}$ sur l'ouvert~$X'\cap Y'$, 
il est raisonnable, si l'on décide (par exemple) de privilégier la variable~$T$, 
de noter ces deux points en question
d'une part~$0$ (pour le point de~$X'$ d'équation~$T=0$, comme il se doit), et d'autre part~$\infty$, 
pour le point de~$Y'$ d'équation~$S=0$ à laquelle on pense comme «~$T^{-1}=0$». 
On peut donc voir~$\PP^1_k$ comme la droite affine à laquelle on a rajouté
un point fermé de corps résiduel~$k$ (donc un point naïf, si l'on veut) «à l'infini». 

\medskip
Si l'on effectue la construction analogue en topologie en remplaçant~$\Aff^1_k$ par 
$\RR$, on obtient un cercle, la droite réelle se recollant par ses deux bouts sur
le point à l'infini ; si l'on remplace~$\Aff^1_k$  par~$\CC$, on obtient une sphère.
Dans les deux cas, l'espace
construit apparaît comme une compactification de l'espace de départ (celle d'Alexandrov, en l'occurrence).

\medskip
Cela reste vrai
{\em mutatis mutandis}
dans le cadre des schémas : nous avons déjà évoqué 
l'existence d'un avatar schématique de la compacité
que nous rencontrerons plus loin, la {\em propreté} ; et nous verrons
à cette occasion que~$\PP^1_k$ est un~$k$-schéma propre. 

\trois{fonc-droite-pro}
Déterminons maintenant la~$k$-algèbre
des sections globales de~$\sch O_{\PP^1_k}$. Se donner un élément de~$\sch O_{\PP^1_k}(\PP^1_k)$,
c'est se donner
une fonction sur~$X'$ et une fonction sur~$Y'$ dont les restrictions à~$X'\cap Y'$ coïncident. 
Autrement dit, cela revient à se donner un polynôme~$P\in k[T]$ et un polynôme~$Q\in k[S]$ dont
les images dans~$\sch O_{\PP^1_k}(X'\cap Y')=k[T,T]^{-1}$ coïncident. La restriction de~$S$ à~$X'\cap Y'$
étant égale à~$T^{-1}$, cette condition de coïncidence signifie simplement que
l'on a~$P=Q(T^{-1})$ dans l'anneau~$k[T,T^{-1}]$. 
La seule possibilité pour qu'un polynôme en~$T$ soit égal à un polynôme en~$T^{-1}$ est évidemment que les deux
soient constants, et il vient~$\sch O_{\PP^1_k}(\PP^1_k)=k$. 

\trois{rem-princ-maximum}
{\em Remarque}. On a signalé au~\ref{topo-droite-proj}
que l'on pouvait penser à~$\PP^1_k$ comme à un objet compact, et l'on a par
ailleurs vu au~\ref{fonc-droite-pro}
que les seules
fonctions globales sur~$\PP^1_k$ sont les constantes : c'est un cas particulier
de l'avatar schématique du {\em principe du maximum}
de la géométrie analytique complexe. 

\trois{p1k-pasaffine}
Comme~$\sch O_{\PP^1_k}(\PP^1_k)=k$, son spectre est réduit à un point. 
Le morphisme canonique~$\PP^1_k\to \spec \sch O_{\PP^1_k}(\PP^1_k)$ n'est
donc évidemment pas un isomorphisme, et~$\PP^1_k$ n'est dès lors pas affine. 

\medskip
En fait, {\em ce n'est même pas un ouvert d'un schéma affine} : en effet, on sait 
d'après le~\ref{reformul-locann-speca}
que tout morphisme de~$\PP^1_k$ vers un schéma affine se factorise par
~$\spec  \sch O_{\PP^1_k}(\PP^1_k)$, et a donc pour image ensembliste un point. 

\deux{droite-dedoubl}
{\bf La droite affine avec origine dédoublée.}
L'isomorphisme de~$k$-algèbres
$$\begin{array}{ccc}
k[S,S^{-1}]&\to& k[T,T^{-1}]\\
S&\mapsto &T\\
S^{-1}&\mapsto &T^{-1}\end{array}$$
induit un isomorphisme
de~$k$-schémas~$\chi : U\to V$. 

\medskip
De manière analogue à ce qui a été fait au~\ref{droite-proj-recoll} 
on définit le recollement de~$X$
et~$Y$ le long de~$\chi$ et~$\chi^{-1}$. 
Cette limite sera notée~$\DD_k$ (ce n'est pas une notation standard). 
Nous allons maintenant la décrire plus avant, en déclinant dans ce
cas particulier les faits mentionnés au~\ref{limind-recoll-schema-cocy}. 

\trois{topo-droite-dedoubl}
Le~$k$-schéma~$\DD_k$  est réunion de deux ouverts affines~$X''$ et~$Y''$
(les images de~$X$ et~$Y$). Chacun d'eux est une copie de la droite affine : on 
a~$X''\simeq \spec k[T]$ et~$Y''\simeq k[S]$. Leur intersection~$X''\cap Y''$
est égale à~$D(T)$ en tant qu'ouvert de~$X''$, et donc à~$\spec k[T,T^{-1}]$ ; elle
est égale à~$D(S)$ en tant qu'ouvert de~$Y''$, et donc à~$\spec k[S,S^{-1}]$. L'isomorphisme
entre~$D(T)\subset X''$ et~$D(S)\subset Y''$ induit par leurs identifications avec~$X''\cap Y''$ 
est celui fourni par le morphisme de~$k$-algèbres qui envoie~$S$ sur~$T$. 

\medskip
Le complémentaire de~$X''\cap Y''$ dans~$\DD_k$ est constitué de deux points fermés : l'origine~$V(T)$
de~$X''\simeq \Aff^1_k$, et l'origine~$V(S)$ de~$Y''\simeq \Aff^1_k$. Mais
comme~$S=T$ sur l'ouvert~$X''\cap Y''$, 
l'origine~$V(S)$ n'est pas cette fois-ci «rejetée à l'infini» : tout se passe 
plutôt comme si l'on avait {\em dédoublé}
l'origine classique en une origine dans~$X''$ et une autre dans~$Y''$. 

\medskip
Si l'on effectuait la construction analogue en topologie en remplaçant~$\Aff^1_k$ par 
$\RR$ (resp.~$\CC$), on obtiendrait «une droite réelle (resp. complexe) avec origine dédoublée»
qui n'est pas un espace séparé : tout voisinage de la première origine rencontre
tout voisinage de la seconde. 

\medskip
Nous verrons plus loin qu'il existe une notion de~$k$-schéma séparé, et que~$\DD_k$ 
n'est justement pas séparé. Mais cette notion n'est pas purement topologique (la topologie
de Zariski n'est de toute façon presque jamais séparée). 

\trois{fonc-droite-dedoubl}
Déterminons maintenant la~$k$-algèbre
des sections globales de~$\sch O_{\DD_k}$. Se donner un élément de~$\sch O_{\DD_k}(\DD_k)$ c'est se donner
une fonction sur~$X''$ et une fonction sur~$Y''$ dont les restrictions à~$X''\cap Y''$ coïncident. 
Autrement dit, cela revient à se donner un polynôme~$P\in k[T]$ et un polynôme~$Q\in k[S]$ dont
les images dans~$\sch O_{\DD_k}(X''\cap Y'')=k[T,T]^{-1}$ coïncident. La restriction de~$S$ à~$X''\cap Y''$
étant égale à~$T$, cette condition de coïncidence signifie simplement que
le poyôme~$Q$
est égal à~$P(S)$. 
La restriction à~$X''$ induit donc un isomorphisme~$\sch O_{\DD_k}(\DD_k)=k[T]$ ; de même, 
la restriction  à~$Y''$ induit un isomorphisme~$\sch O_{\DD_k}(\DD_k)=k[S]$, l'isomorphisme
entre~$k[T]$
et~$k[S]$
induit par ces deux identifications étant celui qui envoie~$T$ sur~$S$.

\trois{deddouble-pasaffine}
Soit~$P$ un polynôme appartenant à~$k[T]$, vu comme fonction
globale sur~$\DD_k$. Par définition, son évaluation en l'origine de~$X''$
(donnée par~$T=0$)
est égale à~$P(0)$. Quant à son évaluation en l'origine de~$Y''$ (donnée par~$S=0$),
elle s'obtient en substituant~$S$ à~$T$, puis en faisant~$S=0$ ; c'est donc encore~$P(0)$. 

\medskip
Il s'ensuit que le morphisme canonique~$p : \DD_k\to \spec \sch O_{\DD_k}(\DD_k)$ envoie
les deux origines sur le même point de~$\spec \sch O_{\DD_k}(\DD_k)\simeq \spec k[T]=\Aff^1_k$, à savoir
l'origine de~$\Aff^1_k$. En particulier, $p$ n'est pas un isomorphisme et~$\DD_k$ n'est
dès lors pas affine.

\medskip
En fait, {\em ce n'est même pas un ouvert d'un schéma affine} : en effet, on sait 
d'après le~\ref{reformul-locann-speca}
que tout morphisme de~$\DD_k$ vers un schéma affine se factorise par~$p$, 
et envoie donc les deux origines sur le même point. 

\subsection*{Produits fibrés de schémas}

\deux{intro-prodfib-schem}
Le but de ce qui suit est de montrer que la catégorie des schémas
admet des produits fibrés, et d'en donner une description (raisonnablement)
explicite. Nous allons commencer par une remarque qui jouera un rôle crucial
pour recoller nos constructions locales.

\trois{prod-fib-ouv-ouv}
Soient~$f : Y\to X$ et~$g: Z\to X$ deux morphismes de schémas. Supposons que l'on sache
que~$Y\times_XZ$ existe ; notons~$p$ sa projection sur~$Y$, et~$q$
sa projection sur~$Z$. 
Soient~$U, V$ et~$W$ des ouverts de~$X,Y$ et~$Z$ respectivement, tels que
$f(V)\subset U$ et~$g(W)\subset U$. 

\medskip
On vérifie alors immédiatement, en combinant
propriétés universelles des produits fibrés et propriétés 
universelles des immersions ouvertes (\ref{prop-univ-schemeouv})
que le produit fibré~$V\times_UW$ existe, et plus
qu'il s'identifie
canoniquement à l'ouvert~$p^{-1}(V)\cap q^{-1}(W)$ de~$Y\times_XZ$
(on note qu'il ne dépend pas de~$U$). 

\trois{prod-fib-ouv}
{\em Un cas particulier.}
Soit~$f: Y\to X$ un morphisme de schémas. Le produit fibré~$Y\times_XX$ existe et
s'identifie tautologiquement à~$Y$. Il résulte alors de~\ref{prod-fib-ouv-ouv}
que pour tout ouvert~$U$ de~$X$ le produit fibré~$Y\times_XU$ existe et s'identifie
à l'ouvert~$f^{-1}(U)$ de~$Y$ (on pourrait aussi le démontrer directement, là encore
à l'aide de la propriété universelle des immersions ouvertes). 

\deux{const-prod-fib}
{\bf Construction des produits fibrés.}
Nous allons procéder en plusieurs étapes. 

\trois{prod-tens-prod-fib}
{\em Produits fibrés dans la catégorie des schémas affines.}
Soit~$A$ un anneau et soient~$B$ et~$C$
deux~$A$-algèbres. Dans la catégorie des anneaux,
la somme amalgamée de~$B$ et~$C$
le long de~$A$ existe : ce n'est autre que
le produit tensoriel~$B\otimes_A C$. Comme
le foncteur $D\mapsto \spec D$ établit une anti-équivalence entre la catégorie des anneaux et celle
des schémas affines, le produit fibré~$\spec B\times_{\spec A}\spec C$
existe dans la catégorie des schémas {\em affines}
et s'identifie à~$\spec (B\otimes_AC)$. 

\trois{prod-aff-prod-fib}
{\em Produits fibrés de schémas affines dans la catégorie des schémas.}
Soient~$Y\to X$ et~$Z\to X$ des morphismes de schémas affines. 
Notons $Y\times_XZ$ leur produit fibré dans la catégorie
des schémas {\em affines}, 
qui existe vertu du~\ref{prod-tens-prod-fib}
ci-dessus. 

\medskip
Soit~$T$ un schéma quelconque. 
On dispose de bijections canoniques fonctorielles en~$T$ 

$$\hom (T, Y\times_XZ)\simeq \hom (\spec \sch O_T(T), Y\times_XZ)$$
$$\simeq \hom(\spec \sch O_T(T), Y)\times_{\hom(\spec \sch O_T(T), X)}\hom(\spec \sch O_T(T), Z)$$
$$\simeq \hom(T,Y)\times_{\hom(T,X)}\hom(T,Z)$$
(la seconde provient de la définition de~$Y\times_XZ$ comme produit fibré dans la catégorie des schémas affines, 
et la première et la troisième du fait que tout
morphisme de~$T$ vers un schéma affine se factorise {\em canoniquement}
par le schéma affine~$\spec \sch O_T(T)$ d'après~\ref{reformul-locann-speca}). 
Il en résulte que~$Y\times_XZ$ est également le produit fibré de~$Y$ et~$Z$ au-dessus de~$X$ dans la catégorie
de {\em tous les schémas}. 

\trois{prod-fib-casgen}
{\em Produits fibrés : le cas général.}
Soient~$f:Y\to X$ et~$g:Z\to X$ des morphismes de schémas. 
Notons~$I$ l'ensemble des triplets~$(U,V,W)$ où~$U$ (resp.~$V$, resp.~$W$)
est un ouvert
affine de~$X$ (resp.~$Y$, resp.~$Z$) et où~$f(V)\subset U$ et~$g(W)\subset U$. 
On munit
l'ensemble~$I$ de l'ordre partiel pour lequel
on a~$(U,V,W)\leq (U',V',W')$ si~$U\subset U', V\subset V'$ et~$W\subset W'$. 

\medskip
Pour tout~$(U,V,W)\in I$, le produit fibré~$V\times_UW$ existe dans la catégorie
des schémas en vertu de~\ref{prod-aff-prod-fib}. Il résulte par ailleurs de~\ref{prod-fib-ouv-ouv}
que si~$(U,V,W)$ et~$(U',V',W')$ dont deux
éléments de~$I$ avec~$(U,V,W)\leq (U',V',W')$, il existe une immersion ouverte canonique
de~$V\times_UW$ dans~$V'\times_{U'}W'$. 

\medskip
La collection des schémas~$V\times_UW$ pour~$(U,V,W)$ parcourant~$I$, et des immersions
ouvertes que nous venons d'évoquer, constitue un diagramme du type mentionné au~\ref{diag-limind-autom-schem}. 
On peut donc procéder au recollement des~$V\times_UW$
le long desdites immersions. Nous laissons le lecteur vérifier que le schéma ainsi obtenu est le produit
fibré~$Y\times_XZ$, la preuve étant essentiellement formelle modulo les deux faits suivants : 

\medskip
$\bullet$ tout schéma est réunion d'ouverts affines ; 

$\bullet$ si~$T$ et~$T'$ sont deux schémas, $\Omega \mapsto \hom(\Omega, T')$ 
définit un {\em faisceau}
sur~$T$ (cela découle immédiatement de la définition d'un morphisme de schémas). 

\medskip
Il résulte de cette construction 
que~$Y\times_XZ$ est réunion de ses ouverts
affines~$V\times_UW$ pour~$(U,V,W)$ parcourant~$I$. 

\deux{comment-prodfib}
{\bf Quelques commentaires.}

\trois{rem-objfinal}
{\em Produits fibrés, objet final, produits cartésiens.}
Soit~$X$ un schéma. En vertu du théorème~\ref{locann-speca},
l'ensemble~$\hom(X,\spec \ZZ)$ est en bijection avec~$\hom_{\ann}(\ZZ,\sch O_X(X))$, lequel
est un singleton. En conséquence, $\spec \ZZ$ est l'objet final de la catégorie des schémas. 

\medskip
Il s'ensuit au vu de~\ref{const-prod-fib}
{\em et sq.} que le produit cartésien de deux schémas~$Y$ et~$Z$
existe toujours : c'est leur produit fibré au-dessus de~$\spec \ZZ$. 

\trois{prod-fib-esp-sj}
L'une des difficultés techniques et psychologiques de la théorie des schémas est
que l'espace topologique sous-jacent à un produit fibré {\em n'est pas}, en général, 
le produit fibré des espaces topologiques sous-jacents. Donnons un contre-exemple très simple. 

\medskip
Le produit fibré $\spec \CC\times_{\spec \RR}\spec \CC$ est égal
à~$\spec \CC\otimes_\RR \CC$. On a
$$\CC\otimes_{\RR}\CC=\CC\otimes_{\RR}\RR[T]/T^2+1\simeq \CC[T]/T^2+1$$
$$=\CC[T]/(T-i)(T+i)\simeq \CC[T]/(T-i)\times \CC[T]/(T+i)\simeq \CC\times \CC.$$

Le produit fibré  $\spec \CC\times_{\spec \RR}\spec \CC$ s'identifie donc
au spectre de~$\CC\times \CC$, c'est à dire à~$\spec \CC\coprod \spec \CC$,
qui comprend deux points. Mais comme~$\spec \CC$ et~$\spec \RR$ sont deux singletons, le produit fibré
{\em topologique}
de~$\spec \CC$ par lui-même au-dessus de~$\spec \RR$ est un singleton. 

\trois{intuition-produit-fibre}
Si~$Y\to X$ et~$Z\to X$ sont deux morphismes de schémas, il y a deux manières
d'envisager le produit fibré~$Y\times_X Z$ (aucune n'est meilleure que l'autre, tout dépend du contexte).
On peut bien sûr y penser comme à un objet symétrique en~$Y$ et~$Z$. Mais on peut également privilégier
l'un des deux facteurs, disons ~$Y$, et voir~$Y\times_XZ$
comme {\em un schéma qui est à~$Z$ ce que~$Y$ est
à~$X$} ; on traduit cette idée en disant que~$Y\times_XZ$ est {\em le~$Z$-schéma déduit du~$X$-schéma~$Y$
par changement de base de~$X$ à~$Z$}. 
Illustrons ces deux visions du produit fibré par des exemples.

\medskip
Soit~$k$ un corps. On a
$$\Aff^1_k\times_{\spec k}\Aff^1_k
=\spec k[S]\otimes_k k[T]=\spec k[S,T]=\Aff^2_k.$$
On se trouve ainsi face à un bon exemple de conception
«symétrique» du produit fibré : le produit 
de la droite affine par elle-même (sur un corps de base fixé) est égal au plan affine. 

\medskip
Donnons-nous maintenant une extension~$L$ de~$k$. On a alors
$$\Aff^1_k\times_{\spec k}\spec L=\spec k[T]\otimes_k L=\spec L[T]=\Aff^1_L.$$
Ici, c'est plutôt la seconde conception qui s'impose : $\Aff^1_L$ est à~$L$ ce que~$\Aff^1_k$ est à~$k$ ; ou encore,
si l'on préfère, le~$L$-schéma~$\Aff^1_L$ se déduit du~$k$-schéma~$\Aff^1_k$ par changement de base de~$k$ à~$L$. 

\trois{aff-base-qcque}
{\em Généralisation des exemples ci-dessus.}
Pour tout~$n$, on note~$\Aff^n_{\ZZ}$ le schéma~$\spec \ZZ[T_1,\ldots, T_n]$. Pour tout couple
$(n,m)$ d'entiers on a
$$\Aff^n_{\ZZ}\times_{\spec \ZZ}\Aff^m_{\ZZ}=\spec \ZZ[T_1,\ldots, T_n]\otimes_{\ZZ}\ZZ[S_1,\ldots, S_m]$$
$$=\spec \ZZ[T_1,\ldots, T_n,S_1,\ldots, S_m]=\Aff^{n+m}_{\ZZ}.$$
Si~$X$ est un schéma quelconque, on pose~$\Aff^n_X=\Aff^n_{\ZZ}\times_{\spec \ZZ}X$ ; on dit
que c'est {\em l'espace affine de dimension~$n$ {\em relatif}
sur~$X$}. 
Si~$A$ est un anneau, on écrira
le plus souvent~$\Aff^n_A$ au lieu de~$\Aff^n_{\spec A}$ ; on a~$\Aff^n_A=\spec A[T_1,\ldots, T_n]$. 
Lorsque~$A$ est un corps, cette notation est compatible avec celle précédemment introduite. 

\medskip
Soit~$X$ un schéma, soit~$Y$ un~$X$-schéma, et soient~$n$ et~$m$ deux entiers. 
On a~$$\Aff^n_X\times_X Y=(\Aff^n_{\ZZ}\times_{\spec \ZZ}X)\times_X Y=\Aff^n_{\ZZ}\times_{\spec \ZZ}Y=\Aff^n_Y$$
et
$$\Aff^n_X\times_X\Aff^m_X=\Aff^n_{\ZZ}\times_{\spec \ZZ} X\times_X (X\times_{\spec \ZZ}\Aff^m_{\ZZ})$$
$$=\Aff^n_{\ZZ}\times_{\spec \ZZ}\Aff^m_{\ZZ}\times_{\spec \ZZ}X=\Aff^{n+m}_{\ZZ}\times_{\spec \ZZ}X=\Aff^{n+m}_X.$$

\deux{fibre-et-produit}
{\bf Structure de schéma sur une fibre.}
Soit~$\psi : Y\to X$ un morphisme de schémas et soit~$x$
un point de~$X$ ; on dispose
d'un morphisme canonique~$\spec \kappa(x)\to X$ (\ref{mor-Kx-conclu}). 
Soit~$p: Y\times_X\spec \kappa(x)\to Y$ la première projection. Le diagramme
commutatif
$$\xymatrix{
{Y\times_X \spec \kappa(x)}\ar[r]
\ar[d]_p&{\spec \kappa(x)}\ar[d]\\
Y\ar[r]^\psi&X}$$
assure que l'image de~$p$ est contenue dans~$\psi^{-1}(x)$. 

\trois{p-homeo-fibre}
{\em Le morphisme~$p$ induit un homéomorphisme
$$Y\times_X\spec \kappa(x)\simeq \psi^{-1}(x).$$}
En effet, grâce à~\ref{prod-fib-ouv-ouv}, on peut raisonner localement
sur~$Y$ et restreindre~$X$ au voisinage de~$x$, ce qui autorise à supposer~$Y$
et~$X$ affines, auquel cas l'assertion voulue découle\footnote{Modulo
le fait suivant, dont nous laissons la vérification au lecteur : si~$A$
est un anneau et~$x$ un point de~$\spec A$, le morphisme
$\spec \kappa(x)\to \spec A$ induit par l'évaluation~$A\to \kappa(x)$
coïncide avec celui défini au~\ref{mor-Kx-conclu}.}
de~\ref{fibre-mor-spec}. 

\trois{struct-schema-fibre}
Cet identification topologique
~$\psi^{-1}(x)\simeq Y\times_X \spec \kappa(x)$
permet de munir la fibre~$\psi^{-1}(x)$
d'une structure de~$\kappa(x)$-schéma. 

\medskip
Que les fibres des morphismes soient elles-mêmes 
des objets de la théorie a de multiples avantages. Nous allons
en mentionner un, particulièrement important : la présence
éventuelle de nilpotents dans le faisceau structural 
du schéma~$\psi^{-1}(x)$ permet de détecter de manière naturelle
et élégante les phénomènes de multiplicité. 

\medskip
Donnons un exemple. Soit~$k$ un corps de caractéristique
différente de 2, et soit~$X$
le~$k$-schéma $\spec k[U,V,T]/(U^2+TV^2-T)$. La flèche naturelle
de~$k[T]$
dans~$k[U,V,T]/(U^2+TV^2-T)$
induit un morphisme~$p$
de~$X$
vers~$\Aff^1_k$. 

Soit~$x\in \Aff^1_k$. Le~$\kappa(x)$-schéma~$\psi^{-1}(x)$ est
égal à
$$\spec \kappa(x)[U,V]/(U^2+T(x)V^2-T(x)).$$
Si~$T(x)\neq 0$, c'est-à-dire si
le point~$x$ n'est pas l'origine de
la droite affine~$\Aff^1_k$, 
le polynôme~$U^2+T(x)V^2-T(x)$ de~$\kappa(x)[U,V]$ est irréductible
(car~$\kappa(x)$ est de caractéristique différente de~$2$), et l'anneau des fonctions
globales sur le schéma~$\psi^{-1}(x)$ est intègre. 

Si~$x$ est l'origine $T(x)=0$, et~$U^2+T(x)V^2-T(x)=U^2$ ; l'anneau des fonctions
globales sur~$\psi^{-1}(x)$ est alors égal à~$k[U,V]/U^2$, {\em et n'est pas réduit.}

\medskip
Cette apparition de nilpotents est la manifestation rigoureuse de ce qu'on décrirait
informellement de la manière suivante : {\em la famille de coniques
affines d'équations~$U^2+TV^2-T=0$, 
dépendant du paramètre~$T$, dégénère en une droite {\em double}
lorsque~$T=0$.}

\trois{fibres-du-produit}
Soient~$\psi : Y\to X$ et~$Z\to X$ deux morphismes de schémas, soit~$z$ un point
de~$Z$ et soit~$x$ son image sur~$X$.  La fibre de~$Y\times_X Z$ en~$z$
s'identifie à
$$(Y\times_X Z)\times_Z \spec \kappa(z)=Y\times_X \spec \kappa(z)$$
$$=(Y\times_X \spec \kappa(x))\times_{\spec \kappa(x)}\spec \kappa(z)
=\psi^{-1}(x)\times_{\spec \kappa(x)}\spec
\kappa(z).$$

\deux{esp-sj-et-prodfb}
Nous allons maintenant essayer de décrire précisément la différence
entre l'espace sous-jacent au produit fibré et le produit fibré des espaces
sous-jacents. Pour ce faire, il va être commode de noter~$|X|$ l'espace topologique
sous-jacent à un schéma~$X$. 

\trois{lemme-espsj}
{\bf Lemme.}
{\em Soient~$\psi : Y\to X$ et~$\chi : Z\to X$ deux morphismes de schémas. Il existe une
application continue surjective naturelle
$$\pi : |Y\times_XZ|\to |Y|\times_{|X|}|Z|.$$ Si~$(y,z)\in  |Y|\times_{|X|}|Z|$
et si~$x$ désigne l'image commune de~$y$ et~$z$ sur~$X$ alors
$$\pi^{-1}(y,z)\simeq \spec \kappa(y)\times_{\spec \kappa(x)}\spec \kappa(z)=\spec (\kappa(y)\otimes_{\kappa(x)}\kappa(z)).$$}

\medskip
{\em Démonstration.}
Le diagramme commutatif
$$\xymatrix{
{Y\times_XZ}\ar[r]\ar[d]&Z\ar[d]\\
Y\ar[r]&X}$$ en induit un
$$\xymatrix{
{|Y\times_XZ|}\ar[r]\ar[d]&|Z|\ar[d]\\
|Y|\ar[r]&|X|}$$
dans la catégorie des espaces topologiques qui, en vertu
de la propriété universelle du produit fibré, induit lui-même
une application continue
$$\pi :  |Y\times_XZ|\to |Y|\times_{|X|}|Z|.$$

\medskip
Donnons-nous maintenant~$x, y$ et~$z$ comme dans l'énoncé, et notons~$p$ et~$q$
les projections respectives de~$Y\times_XZ$ sur~$Y$ et~$X$. Par définition, l'espace topologique
$\pi^{-1}(y,z)$ est l'intersection~$p^{-1}(y)\cap q^{-1}(z)$. 

Il résulte de~\ref{fibres-du-produit}
que~$q^{-1}(z)$ s'identifie à~$\psi^{-1}(x)\times_{\spec \kappa(x)}\spec \kappa(z)$. En
réappliquant~{\em loc. cit.}
à la projection
de~$\psi^{-1}(x)\times_{\spec \kappa(x)}\spec \kappa(y)$
sur~$\psi^{-1}(x)$, on voit que

$$\pi^{-1}(y,z)=p^{-1}(y) \cap q^{-1}(z)\simeq \spec \kappa(y) \times_{\spec \kappa(x)}\spec \kappa(z).$$

Pour conclure, il reste à établir la surjectivité de~$\pi$, c'est-à-dire à s'assurer que
la fibre~$\pi^{-1}(y,z)\simeq \spec (\kappa(y)\otimes_{\kappa(x)}\kappa(z))$ 
est non vide, c'est-à-dire encore que l'anneau $\kappa(y)\otimes_{\kappa(x)}\kappa(z)$
est non nul. Mais c'est immédiat, puisqu'il s'agit du produit tensoriel de deux espaces vectoriels
non nuls sur le corps~$\kappa(x)$. 

\trois{prod-non-vide}
Soient~$\psi : Y\to X$ et~$\chi : Z\to X$ deux morphismes de schémas.
On déduit
du lemme ci-dessus que~$Y\times_XZ$ est vide si et seulement
si~$|Y|\times_{|X|}|Z|$ est vide, c'est-à-dire si et seulement si~$\psi(Y)\cap \chi(Z)= \varnothing$.

\medskip
Comme
la vacuité d'un spectre équivaut à la nullité
de l'anneau correspondant, la
remarque ci-dessus permet de donner une interprétation géométrique d'un phénomène
{\em a priori}
purement algébrique : si~$A$ est un anneau et~$B$ et~$C$ sont deux~$A$-algèbres alors~$B\otimes_AC$
est nul si et seulement si les images de~$\spec B$ et~$\spec C$ sur~$\spec A$ sont disjointes. 

\section{Faisceaux quasi-cohérents}
\markboth{La notion de schéma}{Faisceaux quasi-cohérents}

\subsection*{Faisceaux quasi-cohérents sur un schéma affine}

\deux{intro-fascqc}
Soit~$A$ un anneau, 
soit~$X$ le spectre de~$A$
et soit~$M$ un~$A$-module. 
Aux paragraphes~\ref{intro-spec-locann}
{\em et sq.} nous avons
construit un faisceau~$\red M$
sur~$X$. Lorsque~$M$
est égal à~$A$, c'est ce 
faisceau que nous avons utilisé pour faire de~$X$ un espace localement annelé ; 
en général, $\red M$ est un~$\sch O_X$-module. 

\trois{rapp-prop-mtilde}
Ce faisceau~$\red M$ possède les propriétés suivantes (\ref{fibres-redm},
th.~\ref{theo-faisc-struct}) :

\medskip
$\bullet$ pour tout~$f\in A$ on a $\red M(D(f))=M_f$ ; en particulier
$\red M(X)=M$ (prendre~$f
$ égal à~$1$) ; 

$\bullet$ si~$x$ est un point de~$X$ correspondant à un idéal premier~$\got p$ alors
$$\red M_x=M_{\got p}=A_{\got p}\otimes_A M=\sch O_{X,x}\otimes_AM.$$

\trois{fonct-mtilde}
Cette construction est fonctorielle. Plus précisément, 
soit~$u$ une
application~$A$-linéaire
de~$M$ vers
un~$A$-module~$N$. 
Pour tout ouvert~$U$ de~$X$, elle induit
une application~$S(U)^{-1}A$-linéaire de
$S(U)^{-1}M$ vers~$S(U)^{-1}N$, c'est-à-dire
de~$M_{\rm pref}(U)$ vers~$N_{\rm pref}(U)$
(nous utilisons les notations de~\ref{intro-spec-locann}
{\em et sq.}). Ces flèches étant compatibles aux restrictions
lorsque~$U$ varie, elles définissent 
un morphisme de
préfaisceau de~$M_{\rm pref}$
vers~$N_{\rm pref}$, puis un morphisme
de~$\red u : \red M$ vers~$\red N$ en passant
aux faisceaux associés.
Si~$f\in A$ 
la flèche~$\red u(D(f))\colon \red M(D(f))\to \red N(D(f))$
est l'application~$M_f\to N_f$ déduite 
de~$u$ par localisation, et l'on a un résultat analogue concernant les fibres. 
En particulier, on retrouve~$u$ à partir de~$\red u$ en considérant
le morphisme induit sur les modules de sections globales. 

\medskip
Donnons-nous maintenant
un morphisme~$v$ 
de~$\red M$ vers~$\red N$. Le morphisme~$v$ induit
par passage aux modules des sections
globales une application~$A$-linéaire $w$
de~$M$ vers~$N$. Pour tout~$f\in A$, la commutativité
du diagramme
$$\xymatrix{
M\ar[rr]^{w=v(X)}\ar[d]&&N\ar[d]\\
{M_f}\ar[rr]_{v(D(f))}&&N_f}$$
montre que~$v(D(f))$ coïncide avec l'application
déduite de~$w$ par localisation, c'est-à-dire encore avec
l'application~$M_f\to N_f$ déduite de~$\red w : \red M\to \red N$.

\medskip
Comme les~$D(f)$ forment une base d'ouverts
de~$X$, il vient~$v=\red w$.

\trois{autre-const-redm}
Nous allons proposer une construction de~$\red M$ qui diffère un peu
de celle donnée initialement, et est plus naturelle du point de vue des espaces 
localement annelés (si nous ne l'avons pas utilisée lorsque nous
 lorsque nous avons défini~$\red M$, c'est parce que
 nous n'avions alors pas encore muni~$X$ d'une structure d'espace localement annelé). 
 Nous reprenons les notations de~\ref{intro-spec-locann}
{\em et sq.} 
 
 \medskip
Soit~$M^\sharp$
le préfaisceau sur~$X$ donné par la formule~$U\mapsto  \sch O_X(U)\otimes_A M$ Si~$U$ est un ouvert de~$X$,
l'application de restriction~$A\to \sch O_X(U)$ envoie les éléments de~$S(U)$ sur des éléments
inversibles 
de~$\sch O_X(U)$ ; elle induit donc un morphisme de~$S(U)^{-1}A$ vers~$\sch O_X(U)$, et partant 
un morphisme de~$M_{\rm pref}(U)$ vers~$M^\sharp(U)$. Cette construction est
compatible aux restrictions, 
et l'on obtient ainsi un morphisme de préfaisceaux~$M_{\rm pref}\to M^\sharp$. 

\medskip
Par commutation du produit
tensoriel aux limites inductives, la fibre de~$M^\sharp$ en un point~$x$
de~$X$ s'identifie canoniquement à~$\sch O_{X,x}\otimes_A M$
(c'est une conséquence formelle de l'exactitude à droite du produit tensoriel),
c'est-à-dire
à~$M_{\got p}$ si~$\got p$ est l'idéal premier correspondant à~$x$. Le morphisme
de préfaisceaux~$M_{\rm pref}\to M^\sharp$ induit donc un isomorphisme au niveau des fibres, et
partant un isomorphisme entre les faisceaux associés. Ainsi, $\red M$ est le faisceau
associé à~$M^\sharp$. 

\trois{conse-msharp}
Soit~$U$ un ouvert affine de~$X$. La restriction de~$M^\sharp$
à~$U$ est le préfaisceau~$V\mapsto M\otimes_A V=(M\otimes_A \sch O_X(U))\otimes_{\sch O_X(U)}\sch O_X(V)$ ;
c'est donc le préfaisceau~$(M\otimes_X \sch O_X(U))^\sharp.$ En conséquence, 
$\red M|_U=\red{(M\otimes_X \sch O_X(U))}$ ; on dispose en particulier d'un isomorphisme naturel
$\red M(U)\simeq M\otimes_A\sch O_X(U)$, ce qui n'était jusqu'alors connu que lorsque~$U$ est de la forme~$D(f)$ 
(\ref{rapp-prop-mtilde}).

\deux{def-qc}
{\bf Définition.}
On dit qu'un~$\sch O_X$-module~$\sch F$ sur~$X$ est 
{\em quasi-cohérent}
s'il existe un~$A$-module~$M$ et un isomorphisme~$\red M\simeq \sch F$. 

\trois{comment-quasi}
{\em Commentaire sur la terminologie.}
Il existe également une notion de faisceau {\em cohérent} : 
c'est un faisceau quasi-cohérent satisfaisant certaines conditions de finitude, aisées
à énoncer si~$A$ est noethérien mais plus délicates en général. Nous n'en aurons pas besoin dans ce
cours. 

\trois{prem-exemple}
{\em Premiers exemples.}
Le faisceau nul (qui est égal à~$\red{\{0\}}$) et le
faisceau~$\sch O_X=\red A$ sont quasi-cohérents. 

\trois{equiv-cat}
On déduit de~\ref{rapp-prop-mtilde}
et~\ref{fonct-mtilde}
que~$M\mapsto \red M$ établit une équivalence entre la catégorie
des~$A$-modules et celle des~$\sch O_X$-modules quasi-cohérents, et que~$\sch F\mapsto \sch F(X)$
en est un quasi-inverse. 

\deux{fqque-fqc}
Soit~$\sch F$ un~$\sch O_X$-module quelconque sur~$X$. Pour tout ouvert~$U$
de~$X$, la restriction~$\sch F(X)\to \sch F(U)$ induit un morphisme~$\sch F(X)\otimes_A \sch O_X(U)\to \sch F(U)$,
c'est-à-dire un morphisme~$\sch F(X)^\sharp(U) \to \sch F(U)$, avec les notations de~\ref{autre-const-redm}. 
On obtient ainsi un morphisme de préfaisceaux~$\sch F(X)^\sharp \to \sch F$, qui induit un morphisme de
faisceaux~$\red{\sch F(X)}\to \sch F$. Il est clair que~$\sch F$ est quasi-cohérent si et seulement si ce morphisme
est un isomorphisme. Cela peut se tester sur les fibres, et revient donc à demander
que~$\sch F(X)\otimes_A\sch O_{X,x}\to \sch F_x$ soit un isomorphisme pour tout~$x\in X$. 

\trois{intuit-qc}
Quelle est la signification intuitive de la quasi-cohérence ? On peut y penser comme
à une propriété de {\em stabilisation}. En effet, 
on déduit de~\ref{fqque-fqc}
qu'un~$\sch O_X$-module
$\sch F$ est quasi-cohérent
si et seulement si pour tout ouvert affine~$U$ de~$X$,
le morphisme
naturel~$\sch O_X(U)\otimes_{\sch O_X(X)}\sch F(X)\to \sch F(U)$ est un isomorphisme. 
Autrement dit, lorsqu'on passe de~$X$ à~$U$, le module des sections
de~$\sch F$ change aussi peu que possible : il subit simplement une extension des scalaires,
ce qui veut dire
en quelque sorte qu'aucune nouvelle section ou aucune nouvelle relation
entre sections déjà existantes
ne surgissent {\em ex nihilo}.

\trois{suite-exacte-qc}
Soit~$0\to M'\to M\to M''\to 0$ un diagramme dans la catégorie
des~$A$-modules. D'après le lemme~\ref{lemme-ex-msi}, 
ce diagramme est une suite exacte si et seulement si
$0\to M_{\got p}'\to M_{\got p}\to M''_{\got p}\to 0$
est une suite exacte pour tout idéal premier~$\got p$ de~$A$.

\medskip
Compte-tenu de~\ref{rapp-prop-mtilde}
et du fait que les propriétés d'exactitude faisceautique se détectent
fibre à fibre, on
déduit de ce qui précède que
le diagramme~$0\to M'\to M\to M''\to 0$ 
est une suite exacte si et seulement si~$0\to \red M'\to \red M\to \red M''\to 0$ est une suite
exacte de faisceaux. 

\medskip
Insistons sur une conséquence frappante de cet énoncé : {\em lorsqu'on le restreint à la catégorie
des faisceaux quasi-cohérents}, le foncteur des sections globales~$\sch F \mapsto \sch F(X)$ est exact
(rappelons qu'en général, il est seulement exact {\em  à gauche} : la surjectivité peut poser des problèmes). 

\trois{qc-stab-quot}
Soit~$\sch F\to \sch G$ un morphisme entre faisceaux quasi-cohérents sur~$X$. Son
noyau, son conoyau et son image (faisceautiques !) sont alors quasi-cohérents. 

\medskip
En effet, soient~$M$ et~$N$ les modules de sections globales
de~$\sch F$ et~$\sch G$ ; on a~$\sch F\simeq \red M, \sch G\simeq \red N$, et le morphisme~$\sch F \to \sch G$ est induit par
une application~$A$-linéaire~$M\to N$. Si~$P$ (resp.~$Q$) désigne le noyau (resp. conoyau)
de celle-ci, 
il résulte de~\ref{suite-exacte-qc}
que~$\red P$ (resp.~$\red Q$) est le noyau (resp. conoyau)
de~$\sch F\to \sch G$, d'où notre assertion en ce qui concerne le noyau et le conoyau ;
on traite le cas de l'image en remarquant simplement que c'est le noyau du conoyau. 

\trois{qc-somme-dir}
Soit~$\sch D=((\sch F_i), (E_{ij})$ un diagramme
dans la catégorie des faisceaux quasi-cohérents sur~$X$, et soit~$\sch F$
sa limite inductive dans la catégorie des~$\sch O_X$-modules. 
Pour tout~$i$, soit~$M_i$ le module des sections globales de~$\sch F_i$. 
Le diagramme~$\sch D$ induit un diagramme~$((M_i), (E'_{ij}))$
dans la catégorie des~$A$-modules ; soit~$M$ sa limite inductive. 
On déduit de la commutation du produit tensoriel aux limites
inductives et de~\ref{autre-const-redm}
que~$\red M$ s'identifie à la limite inductive 
des~$\red {M_i}$. En conséquence, $\sch F=\red M$
et~$\sch F$ est quasi-cohérent. 

\medskip
Notons un cas particulier important : une somme directe quelconque de faisceaux quasi-cohérents est quasi-cohérents. 

\trois{qc-limproj}
Par contre, on prendra garde qu'en général une limite projective
de faisceaux quasi-cohérents (calculée dans la catégorie des~$\sch O_X$-modules) n'est pas quasi-cohérent. 
Donnons un exemple simple. Soit~$\sch F$ le~$\sch O_{\spec \ZZ}$-module~$\sch O_{\spec \ZZ}^{\NN}$. 
Nous allons montrer qu'il n'est pas quasi-cohérent. 

\medskip
On a~$\sch F(\ZZ)=\ZZ^{\NN}$
et~$\sch F((D(2)))=(\ZZ[1/2])^{\NN}.$
Mais la flèche canonique $$\ZZ[1/2]\otimes_{\ZZ}\ZZ^{\NN}\to (\ZZ[1/2])^{\NN}$$ n'est pas un isomorphisme : en effet, 
nous invitons
le lecteur à vérifier qu'elle identifie $\ZZ[1/2]\otimes_{\ZZ}\ZZ^{\NN}$ au sous-module
strict de~$ (\ZZ[1/2])^{\NN}$
constitué des suites {\em à dénominateurs bornés.}
En conséquence, $\sch F$ n'est pas quasi-cohérent. 

\medskip
Si l'on reprend le langage un peu imagé utilisé au~\ref{intuit-qc},
on peut dire que lorsqu'on passe de~$\spec \ZZ$ à son ouvert
affine~$D(2)$, une flopée de nouvelles sections de
~$\sch F$ surgissent
{\em ex nihilo} : toutes les suites à dénominateurs non bornés.

\deux{stabilite-qc}
{\bf Fonctorialité.}
Soit~$\phi : A\to B$ un morphisme d'anneaux.
Notons~$Y$ le spectre de~$B$, et soit~$\psi : Y\to X$
le morphisme
de schémas  correspondant
à~$\phi$. Soit~$M$ un~$A$-module et soit~$N$ un~$B$-module. 

\trois{qcoh-imdir}
{\em Le~$\sch O_X$-module~$\psi_*\red N$ s'identifie à~$\red{_AN}$, où~$_AN$ désigne~$N$ vu comme~$A$-module.
En particulier, le faisceau $\psi_*\red N$ est quasi-cohérent.} 
En effet, on a~$\psi_*\red N(X)=\red N(Y)=N$. Soit maintenant~$U$ un ouvert affine
de~$X$, d'anneau des fonctions~$C$. On a
$$\psi_*\red N(U)=\red N(\psi^{-1}(U))=\red N(Y\times_X U) =N\otimes_B(B\otimes_AC)=N\otimes_A C=\red{_AN}(U)$$
Les ouverts affines formant une base de la topologie de~$X$, le morphisme
canonique de~$\red{_AN}=\red{\psi_*\red N(X)}$ vers~$\psi_*\red N$ est un isomorphisme, ce qu'on souhaitait établir. 

\trois{qcoh-imr}
{\em Le~$\sch O_Y$-module~$\psi^*\red M$ s'identifie à~$\red{B\otimes_AM}$
et est en particulier quasi-cohérent.}
En effet, on dispose par construction d'une
application~$B$-linéaire
~$B\otimes_A M\to \psi^*\red M(Y)$ ; il suffit alors
de montrer que la flèche composée
$$u : \red{B\otimes_AM}\to \red{\psi^*\red M(Y)}\to \psi^*\red M$$
est un isomorphisme. Soit~$y\in Y$ et soit~$x$ son image sur~$X$. 
La fibre de~$\red{B\otimes_AM}$ en~$y$ s'identifie
à~$\sch O_{Y,y}\otimes_B(B\otimes_AM)=\sch O_{Y,y}\otimes_A M$. 

\medskip
Quant à la fibre de~$\psi^*\red M$ en~$y$, elle s'identifie 
à
$$\sch O_{Y,y}\otimes_{\sch O_{X,x}}\red M_x=\sch O_{Y,y}\otimes_{\sch O_{X,x}}(\sch O_{X,x}\otimes_A M)=\sch O_{Y,y}\otimes_A M.$$
En conséquence, la flèche~$u$ induit un isomorphisme au niveau des fibres, et est 
de ce fait elle-même un isomorphisme.

\trois{rem-imm-ouv-qc}
{\em Remarque.}
Ce qu'on a vu plus haut au~\ref{conse-msharp}
est un cas particulier de ce qui précède : celui où~$Y$ est un ouvert affine
de~$X$. 

\subsection*{Caractère local de la quasi-cohérence, faisceaux
quasi-cohérents sur un schéma quelconque}

\deux{intro-qcoh-loc}
Nous allons maintenant établir un résultat fondamental, qui n'a rien d'évident au
vu des définitions : le fait que la quasi-cohérence est une propriété {\em locale.}
Comme vous allez le voir, la preuve n'est pas triviale, même si elle ne repose
{\em in fine}
que sur la sempiternelle condition d'égalité de fractions dans un module localisé.

\deux{theo-qc-loc}
{\bf Théorème.}
{\em Soit~$A$ un anneau et soit~$X$ son spectre. Soit~$(U_i)_{i\in I}$ un recouvrement de~$X$
par des ouverts affines, et soit~$\sch F$ un~$\sch O_X$-module. Les assertions suivantes sont équivalentes : 

\medskip
1) $\sch F$ est quasi-cohérent ; 

2) $\sch F|_{U_i}$ est quasi-cohérent pour tout~$i$.}

\medskip
{\em Démonstration.}
L'implication~1)$\Rightarrow$2) découle directement de~\ref{conse-msharp}. Supposons
maintenant que~2) soit vraie, et montrons~1). Toujours grâce à~\ref{conse-msharp}, on peut 
raffiner le recouvrement~$(U_i)$ sans altérer~2). Cela autorise à supposer que~$U_i$ est de
la forme~$D(f_i)$ pour tout~$i$ ; comme~$X$ est quasi-compact, on peut également
faire l'hypothèse que l'ensemble d'indices~$I$ est fini. 

\medskip
Pour montrer que~$\sch F$ est quasi-cohérent, nous allons vérifier
que le morphisme naturel~$\red{\sch F(X)}\to \sch F$ est un isomorphisme. Cela peut se tester localement ; 
il suffit donc de s'assurer que~$\red {\sch F(X)}|_{U_i}\to \sch F|_{U_i}$ est un isomorphisme pour tout~$i$.

\medskip
Fixons~$i$, et montrons que~$\red {\sch F(X)}|_{U_i}\to \sch F|_{U_i}$ est un isomorphisme. 
Les~$\sch O_{U_i}$-modules $\red {\sch F(X)}|_{U_i}$
et~$\sch F|_{U_i}$ étant quasi-cohérents, cette dernière condition 
se teste sur les sections globales : il suffit donc
de s'assurer que
$$\red{\sch F(X)}(U_i)=\sch F(X)_{f_i}\to \sch F(U_i)$$ est un isomorphisme.

\trois{qc-loc-inj}
{\em Injectivité de~$\sch F(X)_{f_i}\to \sch F(U_i)$.}
Soit~$s/f_i^r$ un élément de~$\sch F(X)_{f_i}$ dont l'image
dans~$\sch F(U_i)$ est nulle. Comme~$f_i$ est inversible sur~$U_i=D(f_i)$, 
la restriction de~$s$ à~$U_i$ est nulle. 

\medskip
Soit~$j\in I$ 
(ce qui suit est trivial si~$j=i$, mais
on n'exclut pas ce cas). La restriction de~$s$ à~$U_i\cap U_j$ est
{\em a fortiori}
nulle,
ce que l'on peut récrire $(s|_{U_j})|_{U_i\cap U_j}=0$.
Le faisceau~$\sch F|_{U_j}$
est quasi-cohérent, et~$U_i\cap U_j$ est l'ouvert~$D(f_i)$ du schéma affine~$U_j$ ; en conséquence, 
$\sch F(U_i\cap U_j)$ s'identifie à~$\sch F(U_j)_{f_i}$. L'annulation de~$s|_{U_j}$ sur~$U_i\cap U_j$ signifie
alors qu'il existe~$N$ tel que~$f_i^Ns|_{U_j}=0$ ; comme l'ensemble~$I$ est fini, on peut choisir~$N$ de sorte que cette égalité 
vale pour tout~$j$. Puisque~$\sch F$ est un faisceau, il vient
$f_i ^Ns=0$, et la fraction~$s/f_i^r\in \sch F(X)_{f_i}$ est nulle, ce qui achève de montrer l'injectivité requise. 

\trois{qc-loc-surj}
{\em Surjectivité de~$\sch F(X)_{f_i}\to \sch F(U_i)$.}
Soit~$\sigma\in \sch F(U_i)$. Il s'agit de montrer qu'il existe un entier~$N$ tel que~$f_i^N\sigma$
se prolonge
en une section globale de~$\sch F$. 

\medskip
Soit~$j\in I$
(ce qui suit est trivial si~$j=i$, mais
on n'exclut pas ce cas). Le faisceau~$\sch F|_{U_j}$
est quasi-cohérent, et~$U_i\cap U_j$ est l'ouvert~$D(f_i)$ du schéma affine~$U_j$ ; en conséquence, 
$\sch F(U_i\cap U_j)$ s'identifie à~$\sch F(U_j)_{f_i}$. Il s'ensuit qu'il existe~$n$ tel que la restriction
$f_i^n\sigma|_{U_i\cap U_j}$ se prolonge en une section~$\sigma_j$ de~$\sch F$ sur~$U_j$ ; 
comme l'ensemble~$I$ est fini, on peut choisir
un entier~$n$ convenant pour tout~$j$. 

\medskip
Soient~$j$ et~$j'$ deux éléments de~$I$. Les restrictions de~$\sigma_j|_{U_j\cap U_{j'}}$ et~$\sigma_{j'}|_{U_j\cap U_{j'}}$ 
à~$U_j\cap U_{j'}\cap U_i$ coïncident (elles sont toutes deux égales à la restriction de~$f_i^n\sigma$). 
Comme~$U_j\cap U_{j'}$ est l'ouvert affine~$D(f_{j'})$ de~$U_j$, le faisceau~$\sch F|_{U_j\cap U_{j'}}$ est quasi-cohérent ; puisque
$U_i\cap U_j\cap U_{j'}$ est l'ouvert~$D(f_i)$ du schéma affine~$U_j\cap U_{j'}$, 
$\sch F(U_i\cap U_j\cap U_{j'})$ s'identifie à~$\sch F(U_j\cap U_{j'})_{f_i}$. Il existe donc un entier~$\ell$
tel que~$f_i^\ell \sigma_j|_{U_j\cap U_{j'}}=f_i^\ell\sigma_{j'}|_{U_j\cap U_{j'}}$ ; comme l'ensemble~$I$ est fini, on peut choisir
un entier~$\ell$
convenant pour tout~$(j,j')$. 

\medskip
Posons~$N=n+\ell$. Par construction, la famille des~$f_i^\ell\sigma_j$ se recolle (pour~$j$ variable)
en une section globale~$s$ de~$\sch F$, et~$s|_{U_i}=f_i^N\sigma$.~$\Box$ 

\deux{def-qc}
{\bf Proposition-définition.}
{\em Soit~$X$ un schéma et soit~$\sch F$ un~$\sch O_X$-module. Les assertions suivantes sont
équivalentes :

\medskip
1) pour tout ouvert affine~$U$ de~$X$, le~$\sch O_U$-module~$\sch F|_U$ est quasi-cohérent ;

2) il existe un recouvrement~$(U_i)$ de~$X$ par des ouverts affines tels
que~$\sch F|_{U_i}$ soit un~$\sch O_{U_i}$-module quasi-cohérent pour tout~$i$.

\medskip
Lorsqu'elles sont satisfaite, on dit que~$\sch F$ est {\em quasi-cohérent.}}

\medskip
{\em Démonstration.} 
Il est clair que~1)$\Rightarrow$2). Réciproquement, supposons que~2) soit satisfaite, et soit~$U$ un 
ouvert affine de~$X$. Il existe un recouvrement ouvert~$(V_j)$
de~$U$ par des ouverts affines tel que~$V_j$
soit contenu pour tout~$j$ dans~$U_{i(j)}\cap U$ pour un certain indice~$i(j)$. 

Pour tout~$j$, le faisceau~$\sch F|_{V_j}$ est la restriction du faisceau quasi-cohérent~$\sch F|_{U_{i(j)}}$, et est donc
quasi-cohérent. Il résulte alors du théorème~\ref{theo-qc-loc}
que~$\sch F|_U$ est quasi-cohérent.~$\Box$ 

\deux{fais-quas-sorites}
{\bf Premières propriétés.}
Soit~$X$ un schéma. 

\trois{ox-qc}
Le faisceau nul ainsi que le structural~$\sch O_X$ sont quasi-cohérents en vertu de~\ref{prem-exemple}.  

\trois{qc-proploc}
Soit~$\sch F$ un~$\sch O_X$-module. Il résulte immédiatement de la définition que si~$\sch F$ est quasi-cohérent,
sa restriction à tout ouvert de~$X$ l'est encore ; et que s'il existe un recouvrement ouvert~$(U_i)$ de~$X$ tel que~$\sch F|_{U_i}$ 
soit quasi-cohérent pour tout~$i$, alors~$\sch F$ est quasi-cohérent.

%
\trois{qc-stable}
On déduit de~\ref{qc-stab-quot}
et~\ref{qc-somme-dir}
que le noyau, le conoyau 
et l'images
d'un morphisme de~$\sch O_X$-modules quasi-cohérents
sont quasi-cohérents, et qu'une 
limite inductive de~$\sch O_X$-modules quasi-cohérents est quasi-cohérente ; en particulier, 
une somme directe de~$\sch O_X$-modules quasi-cohérents est quasi-cohérente.

\trois{imrec-qc}
Soit~$\psi\colon Y\to X$ un morphisme de schémas et soit~$\sch F$ un~$\sch O_X$-module quasi-cohérent. 
L'image réciproque~$\psi^*\sch F$ de~$\sch F$ sur~$Y$ est alors un~$\sch O_Y$-module quasi-cohérent. 
En effet, comme la propriété est locale, on peut raisonner localement sur~$Y$ et~$X$, et se ramener ainsi au cas
où tous les deux sont affines,
pour lequel l'assertion requise a été démontrée en~\ref{qcoh-imr}.

\deux{loc-libre-qch}
Soit~$X$ un schéma. Si~$\sch F$ est un~$\sch O_X$-module localement
libre de rang fini, il est quasi-cohérent : c'est une conséquence immédiate du caractère local
de la quasi-cohérence, et de la quasi-cohérence de~$\sch O_X^m$ pour tout~$m$. 

\medskip
Supposons maintenant que~$X=\spec A$ pour un certain anneau~$A$, et soit~$M$ un~$A$-module. 
On déduit de la quasi-compacité de~$\spec A$ que le faisceau quasi-cohérent~$\red M$ est localement
libre de rang fini si et seulement si il existe une famille {\em finie} $(f_i)$
d'éléments de~$A$ et une famille finie~$(n_i)$ d'entiers tels que~$\red M|_{D(f_i)}$
soit isomorphe à~$\sch O_{D(f_i)}^{n_i}$ 
pour tout~$i$, c'est-à-dire encore tels que~$M_{f_i}\simeq A_{f_i}^{n_i}$ pour tout~$i$. 
On déduit alors de~\ref{theo-proj-fond}
que~$\red M$ est localement libre de rang fini si et seulement si~$M$ est projectif
et de type fini.

\deux{qc-imdir}
{\bf Image directe d'un faisceau
quasi-cohérent.}
La situation est moins simple que pour l'image réciproque : comme nous allons le voir,
l'image directe d'un faisceau quasi-cohérent n'est pas quasi-cohérente en général. Nous allons
commencer par énoncer une condition de finitude suffisante -- un peu rébarbative -- pour qu'elle le soit, 
puis nous donnerons un contre-exemple dans un cas assez simple où cette condition n'est pas
remplie. 

\trois{imdir-qc}
Soit~$A$ un anneau, soit~$X$ son spectre et soit~$\psi : Y\to X$ un morphisme de schémas. Soit~$\sch F$
un faisceau quasi-cohérent sur~$Y$. On suppose
que~$Y$ satisfait la propriété suivante : 

\medskip
($*$)~{\em il existe un recouvrement fini~$(U_i)$ de~$Y$ par des ouverts affines tels que~$U_i\cap U_j$ soit pour tout~$(i,j)$
une réunion finie d'ouverts affines}. 

\medskip
Nous allons démontrer que sous cette hypothèse, $\psi_*\sch F$ est un~$\sch O_X$-module quasi-cohérent. Pour cela, 
on écrit chacun des~$U_i\cap U_j$ comme réunion finie d'ouverts affine~$V_{ij\ell}$. On 
écrit~$U_i=\spec A_i$ et~$V_{ij\ell}=\spec B_{ij\ell}$ pour tout~$i,j,\ell$.

Comme~$\sch F$ est un faisceau, on a une suite exacte 

$$\xymatrix{0\ar[r]&{ \sch F(Y)}\ar[r]&{\prod_i \sch F(U_i)}\ar[rrrr]^(0.44){(s_i)_i\mapsto \left(s_i|_{V_{ij\ell}}-s_j|_{V_{ij\ell}}\right)_{ij\ell}}&&&&{\prod_{i,j,\ell}\sch F(V_{ij\ell})}}.$$

Soit~$f\in A$. On note~$Y', U'_i$ et~$V'_{ij\ell}$ les produits
fibrés 
$$Y\times_X D(f), U_i\times_X D(f) \;\text{et}\;V_{ij\ell}\times_X D(f),$$
et le~$f$ en indice fera référence à la localisation {\em en tant que~$A$-module}. 
La~$A$-algèbre~$A_f$ est plate. Par ailleurs, les produits intervenants dans la suite exacte ci-dessus comprennent
par hypothèse un nombre fini de facteurs ; ce sont donc également des sommes directes, et de ce fait ils commutent au produit tensoriel. 
Il en résulte qu'en appliquant~$A_f\otimes_A \bullet$ à la suite précédente, on obtient
une suite exacte 
$$\xymatrix{0\ar[r]&{\sch F(Y)_f}\ar[r]&{\prod_i \sch F(U_i)_f}\ar[rrrr] ^(0.44){(s_i)_i\mapsto \left(s_i|_{V_{ij\ell}}-sj|_{V_{ij\ell}}\right)_{ij\ell}}
&&&&{\prod_{i,j,\ell}\sch F(V_{ij\ell})_f}}.$$
Fixons $i,j$ et~$\ell$. On a
$$\sch F(U_i)_f=A_f\otimes_A \sch F(U_i) =(A_f\otimes_A A_i)\otimes_{A_i} \sch F(U_i)=\sch F(U'_i)$$
car~$U'_i=\spec (A_f\otimes_A A_i)$
et car~$\sch F|_{|U_i}$ est quasi-cohérent. 
On a de même
$$\sch F(V_{ij\ell})_f=A_f\otimes_A \sch F(V_{ij\ell}) =(A_f\otimes_A B_{ij\ell})\otimes_{B_{ij\ell}} \sch F(V_{ij\ell})=\sch F(V'_{ij\ell}).$$

La suite exacte ci-dessus se récrit donc
$$\xymatrix{0\ar[r]&{\sch F(Y)_f}\ar[r]&{\prod_i \sch F(U'_i)}\ar[rrrr]^(0.44){(s_i)_i\mapsto \left(s_i|_{V'_{ij\ell}}-s_j|_{V'_{ij\ell}}\right)_{ij\ell}}&&&&{\prod_{i,j,\ell}\sch F(V'_{ij\ell})}}.$$

Les~$U'_i$ forment
un recouvrement ouvert de~$Y'$, et les~$V'_{ij\ell}$ forment pour tout~$i$
un recouvrement ouvert de~$U_i$. En utilisant une fois encore le fait que~$\sch F$ faisceau,  on
déduit de ce qui précède que la flèche naturelle de~$\sch F(Y)_f$ vers~$\sch F(Y')=\sch F(Y\times_XD(f))=\psi_*\sch F(D(f))$ est un isomorphisme. 
Comme les~$D(f)$
forment une base d'ouverts de~$X$, le morphisme
naturel de~$\sch O_X$-modules~$\red {_A\sch F(Y)}=\red{\psi_*\sch F(X)}\to \psi_*\sch F$ est un isomorphisme, 
et~$\psi_*\sch F$ est quasi-cohérent. 

\trois{contre-ex-imdir-qc}
{\bf Un contre-exemple.}
Dans le raisonnement suivi ci-dessus, l'hypothèse~$(*)$ joue un rôle crucial : elle
permet de décrire le module des sections globales du faisceau~$\sch F$ 
par une suite exacte mettant en jeu des produits {\em finis}
de modules, produits qui sont donc des sommes directes, et commutent dès lors au produit tensoriel. 
Nous allons maintenant
donner un exemple simple de situation où~$(*)$ est prise en défaut,
et où il existe un faisceau quasi-cohérent dont l'image directe n'est pas quasi-cohérente. 

\medskip
Soit~$T$ la
somme disjointe
de copies de~$\spec \ZZ$ paramètrée par
$\NN$, et soit
$\psi$ l'unique
morphisme de~$T$ vers~$\spec \ZZ$.
Le schéma~$T$ n'est manifestement pas quasi-compact, 
et ne peut donc être réunion finie d'ouverts affines ; en conséquence, l'hypothèse~$(*)$ 
n'est pas vérifiée. Nous allons montrer que~$\psi_*\sch O_T$
{\em n'est pas}
quasi-cohérent. 

\medskip
Pour tout ouvert~$U$ de~$\spec \ZZ$, l'ouvert~$\psi^{-1}(U)$
est somme disjointe de copies de~$U$ paramètrées par~$\NN$, et l'on a donc
On a~
$$\psi_*\sch O_T(U)=\sch O_T(\psi^{-1}(U))=(\sch O_{\spec \ZZ}(U))^{\NN}.$$
Ainsi, $\psi_*\sch O_T=\sch O_{\spec \ZZ}^\NN$, 
dont on a vu au~\ref{qc-limproj}
qu'il n'est pas
quasi-cohérent.

\trois{imdir-qc-casgen}
Les résultats de~\ref{imdir-qc}
(qui concernait le cas d'un schéma de
base affine) peuvent se globaliser comme suit. 
Soit~$\psi : Y\to X$ un morphisme de schémas et soit~$\sch F$ un~$\sch O_Y$-module
quasi-cohérent. Supposons que tout point de~$X$
admette un voisinage ouvert affine~$U$ tel que~$\psi^{-1}(U)$ satisfasse~$(*)$ ;
en vertu du caractère local de la quasi-cohérence,
le~$\sch O_X$-module~$\psi_*\sch F$ est
alors quasi-cohérent. 

\subsection*{Faisceaux quasi-cohérents d'idéaux et fermés}

\deux{def-vi-qc}
Soit~$X$ un schéma et soit~$\sch I\subset \sch O_X$ un faisceau quasi-cohérent d'idéaux. 
On note~$V(\sch I)$ le sous-ensemble de~$X$ constitué des points~$x$ possédant la propriété suivante : {\em 
pour tout voisinage ouvert~$U$ de~$x$ et toute fonction~$f$
appartenant à~$\sch I(U)\subset \sch O_X(U)$ on a~$f(x)=0$.}

\trois{vi-ferme-casaff}
Supposons que~$X$ est le spectre d'un anneau~$A$ ; le faisceau~$\sch I$ est alors
égal à~$\red I$ pour un certain idéal~$I$ de~$A$. Nous allons montrer
que~$V(\sch I)$ est égal au fermé~$V(I)$. 

\medskip
Comme~$I=\sch I(X)$ on a par définition~$V(\sch I) \subset V(I)$. 
Réciproquement, soit~$x\in V(I)$ et soit~$U$ et~$f$ comme dans~\ref{def-vi-qc} ; il s'agit
de montrer que~$f(x)=0$. Quitte à restreindre~$U$ on peut le supposer de la forme~$D(g)$ avec~$g(x)\neq 0$. 
Dans ce cas~$f\in \red I(D(g))=I_g$, ce qui veut dire que~$f$ s'écrit~$a/g^n$ avec~$a\in I$ et~$n\geq 0$. 
Comme~$a\in I$ et~$x\in V(I)$ on a~$a(x)=0$ et~$f(x)=0$, ce qu'on souhaitait. 

\trois{vi-ferme-casgen}
On ne suppose plus~$X$ affine. Il découle du~\ref{vi-ferme-casaff}
ci-dessus que~$V(\sch I)\cap U$ est égal au fermé~$V(\sch I(U))$
de~$U$ pour tout ouvert affine~$U$
de~$X$. Comme être fermé est une propriété locale, $V(\sch I)$ est un fermé
de~$X$. 

\trois{vi-support}
On peut également caractériser $V(\sch I)$ comme le {\em support}
du~$\sch O_X$-module quasi-cohérent~$\sch O_X/\sch I$, c'est-à-dire
comme l'ensemble des points~$x$ de~$X$ tel que la fibre
$(\sch O_X/\sch I)_x$ soit {\em non nulle}. 

\medskip
En effet, soit~$x\in X$. La suite exacte
$$0\to \sch I\to \sch O_X\to \sch O_X/\sch I\to 0$$
induit une suite exacte de~$\sch O_{X,x}$-modules
$$0\to \sch I_x\to \sch O_{X,x}\to (\sch O_X/\sch I)_x\to 0.$$ 
La fibre~$(\sch O_X/\sch I)_x$ s'identifie donc à~$\sch O_{X,x}/\sch I_x$. Elle est nulle
si et seulement si il existe un élément~$f\in \sch I_x$ qui n'appartient pas à l'idéal
maximal de~$\sch O_{X,x}$, c'est-à-dire une section~$f$ de~$\sch I$ définie
au voisinage de~$x$ et telle que~$f(x)\neq 0$. Autrement dit, 
$$(\sch O_X/\sch I)_x=\{0\}\iff x\notin V(\sch I),$$
comme annoncé. 

\trois{rem-support-oxi}
{\em Remarque.}
Notons~$j$ l'inclusion~$V(\sch I)\hookrightarrow X$. Par ce qui précède,
$(\sch O_X/\sch I)|_{X\setminus V(\sch I)}$ est nul. Le lecteur est invité
à démontrer que cela équivaut à dire que l'homomorphisme
canonique
$$\sch O_X/\sch I\to j_*j^{-1}\sch O_X/\sch I$$ est un isomorphisme. 

\deux{digr-schem-red}
Nous nous proposons maintenant de démonter que réciproquement,
tout fermé d'un schéma~$X$ est de la forme~$V(\sch I)$ pour un certain faisceau d'idéaux
quasi-cohérent~$\sch I$ sur~$X$. Pour cela, il est nécessaire de faire une petite digression
et d'introduire la notion de schéma {\em réduit}. 

\deux{lemme-sch-red}
{\bf Lemme-définition.}
{\em Soit~$X$ un schéma. Les assertions suivantes sont équivalentes : 

\medskip
i) pour tout~$x\in X$, l'anneau local~$\sch O_{X,x}$ est réduit ; 

ii) pour tout ouvert~$U$ de~$X$, l'anneau~$\sch O_X(U)$ est réduit ; 

iii) il existe un recouvrement~$(U_i)$ de~$X$ par des ouverts affines tels que~$\sch O_X(U_i)$ soit
réduit pour tout~$i$.

\medskip
Lorsqu'elles sont satisfaites, on dit que~$X$ est {\em réduit}.}

\medskip
{\em Démonstration.}
Supposons que~i) est vraie, et soit~$U$ un ouvert de~$X$. Soit~$f$
un élément nilpotent de~$\sch O_X(U)$.
Pour tout~$x\in U$, le germe de~$f$ en~$x$ est nilpotent et donc nul d'après l'hypothèse~i). 
Ainsi, $f$ est nulle, et~ii) est vraie. 

\medskip
Il est clair que~ii) entraîne~iii). Supposons maintenant que~iii) est vraie, et soit~$x\in X$. Il existe~$i$ tel
que~$x\in U_i$. Par hypothèse, $U_i$ s'écrit~$\spec A$ pour un certain anneau réduit~$A$. L'anneau local
$\sch O_{X,x}=\sch O_{U_i,x}$ est de la forme~$A_{\got p}$, où~$\got p$ est un idéal premier de~$A$ ; il est en
conséquence réduit
d'après
le lemme~\ref{coroll-reduit}.~$\Box$  

\deux{retour-ferme-viqc}
Revenons maintenant au problème qui nous intéresse. On se donne un schéma~$X$ et un fermé~$F$ de~$X$. Soit~$\sch I(F)$ 
le sous-faisceau d'idéaux de~$\sch O_X$ défini par la formule
$$U\mapsto \{f\in \sch O_X(U), f(x)=0\;\;\forall x\in F\cap U\}.$$
Nous allons démontrer que~$\sch I(F)$ est quasi-cohérent et que~$F=V(\sch I(F))$. Ces deux propriétés étant locales, on peut supposer
que~$X$ est affine ; c'est donc le spectre d'un anneau~$A$, et~$F=V(I)$ pour un certain idéal~$I$ de~$A$ que l'on peut choisir saturé. 

\trois{schi-f-i}
Nous allons montrer que~$\sch I(F)=\red I$, ce qui assurera la quasi-cohérence de~$\sch I(F)$. En vertu de~\ref{vi-ferme-casaff},
on a~$V(\red I)=V(I)=F$, ce qui montre que~$\red I\subset \sch I(F)$ ; nous allons établir l'inclusion
réciproque. 

\medskip
Soit donc un ouvert~$U$ de~$X$ et~$f$ un élément de~$\sch I(F)(U)$ ; il s'agit de vérifier que~$f\in \red I(U)$. On peut s'en assurer
localement, et donc supposer que~$U$ est de la forme~$D(g)$ pour un certain~$g\in A$. 
La fonction~$f$ s'annule alors par hypothèse en tout point du
fermé~$F\cap D(g)$ de~$D(g)=\spec A_g$, qui n'est autre que~$V(I\cdot A_g)=V(I_g)$. 

\medskip
Puisque~$V(I_g)$ s'identifie à~$\spec A_g/I_g$, l'image~$\bar f$ de~$f$
dans~$A_g/I_g$ s'annule en tout point de~$\spec A_g/I_g$, ce 
qui veut dire qu'elle est nilpotente. Mais comme~$I$ est saturé, $A/I$ est réduit, et~$\spec A/I$ est donc un schéma réduit ; en conséquence,
$A_g/I_g=\sch O_{\spec A/I}(D(\bar g))$ est réduit, et son élément nilpotent~$\bar f$ est dès lors nul. Il s'ensuit~$f\in I_g=\red I(D(g))$,
ce qu'on souhaitait. 

\trois{schi-f-i-conclu}
On a donc~$\sch I(F)=\red I$, et partant~$V(\sch I(F))=V(\red I)$. Mais d'après~\ref{vi-ferme-casaff}, ce dernier est égal à~$V(I)=F$, ce qui achève
la démonstration.  

\trois{rem-schi-f-i-conclu}
{\em Remarque.}
Il résulte immédiatement des définitions que si~$\sch J$ est un faisceau d'idéaux quasi-cohérent tel que~$V(\sch J)=F$ alors~$\sch J$
est contenu dans~$\sch I(F)$ ; 
autrement dit, $\sch I(F)$ est {\em le plus grand}
faisceau quasi-cohérent d'idéaux définissant~$F$.

\section{Morphismes affines}
\markboth{La notion de schéma}{Morphismes affines}

\subsection*{Spectre d'une algèbre quasi-cohérente et morphismes
affines}

\deux{def-alg-sc}
Soit~$X$ un schéma. Une~$\sch O_X$-algèbre~$\sch A$ est dite quasi-cohérente 
si le~$\sch O_X$-module~$\sch A$ est quasi-cohérent. Donnons quelques exemples. 

\trois{ox-i-algqc}
Si~$\sch I$ est un faisceau quasi-cohérent sur~$X$, le quotient~$\sch O_X/\sch I$ est une~$\sch O_X$-algèbre
quasi-cohérente. 

\trois{algqc-casaff}
Supposons que~$X$ est le spectre d'un anneau~$A$. Si~$B$ est une~$A$-algèbre,
le faisceau quasi-cohérent~$\red B$ hérite d'une structure
naturelle de~$\sch O_X$-algèbre. Il découle de~\ref{equiv-cat}
que~$B\mapsto \red B$ induit une équivalence entre la catégorie des~$A$-algèbres et celle des~$\sch O_X$-
algèbres quasi-cohérentes, dont~$\sch B\mapsto \sch B(X)$ est un quasi-inverse. 

\deux{construction-spec-algqc}
{\bf Le spectre d'une algèbre quasi-cohérente}. Le but de ce qui suit est de donner
une variante globale (ou relative, ou faisceautique, comme on voudra) du foncteur
$A\mapsto \spec A$. On fixe un schéma~$X$, et une~$\sch O_X$-algèbre quasi-cohérente~$\sch A$. 

\trois{contr-loc-spec-algqc}
Pour tout ouvert affine~$U$ de~$X$, on pose~$Y_U=\spec \sch A(U)$ ; c'est un schéma affine. Comme~$\sch A (U)$ 
est une~$\sch O_X(U)$-algèbre, le schéma~$Y_U$ est fourni avec un morphisme naturel
$Y_U\to \spec \sch O_X(U)=U\hookrightarrow X$. 

\trois{pre-recollement-spec-au}
Soient maintenant~$U$ et~$V$ deux ouverts affines de~$X$ tels que~$V\subset U$. Comme~$\sch A$ est quasi-cohérente, 
on a~$\sch A(V)=\sch O_X(V)\otimes_{\sch O_X(U)}\sch A(U)$ ; il vient~$Y_V=Y_U\times_U V$. Il existe donc une immersion
ouverte naturelle~$\iota_{VU}$
de~$Y_V$ dans~$Y_U$, laquelle est un~$X$-morphisme.

\trois{recollement-spec-au}
Le diagramme constitué des~$Y_U$ et des immersions~$\iota_{VU}$ est du type décrit au~\ref{diag-limind-autom-schem} ; on peut
donc recoller les~$Y_U$ le long des~$\iota_{VU}$ ; on obtient un~$X$-schéma que l'on appelle le {\em spectre de la~$\sch O_X$-algèbre quasi-cohérente}
$\sch A$ et que l'on note~$\spec \sch A$. Soit~$\pi$ le morphisme~$\spec \sch A\to X$. 

\medskip
Soit~$U$ un ouvert affine de~$X$. On a
par construction une identification naturelle
$$\pi^{-1}(U)=\spec \sch A\times_X U\simeq \spec \sch A(U),$$
modulo laquelle le morphisme~$\pi^{-1}(U)\to U=\spec \sch O_X(U)$ est induit par la flèche
structurale~$\sch O_X(U)\to \sch A(U)$. Il en résulte un isomorphisme
$$\pi_*(\sch O_{\pi^{-1}(U)})\simeq \red{(\sch A(U))}\simeq \sch A|_U.$$ 
La formation de ces isomorphismes commute aux restrictions, et elle
induit donc (en vertu du fait que~$\underline{\text{Isom}}(\pi_*\sch O_{\spec \sch A}, \sch A)$ est un faisceau)
un isomorphisme naturel~$\pi_*\sch O_{\spec \sch A}\simeq \sch A$.

\deux{algqc-ex-prop}
{\bf Exemples et premières propriétés.}

\trois{spec-relatif-casaff}
{\em Un exemple trivial.}
Soit~$A$ un anneau et soit~$B$ une~$A$-algèbre. Par construction, le~$A$-schéma $\spec \red B$ s'identifie à~$\spec B$. 

\trois{esp-aff-rel-moraff}
Soit~$X$ un schéma. 
Le faisceau~$\sch A :=U\mapsto \sch O_X(U)[T_1,\ldots, T_n]$ est de manière naturelle une~$\sch O_X$-algèbre
quasi-cohérente. Il découle immédiatement des définitions que le~$X$-schéma~$\spec \sch A$ s'identifie
à~$\Aff^n_X$.

\trois{algqc-contrav}
Soit~$X$ un schéma. Si~$\sch A$ et~$\sch B$ sont deux~$\sch O_X$-algèbres
quasi-cohérentes, tout morphisme~$\sch A \to \sch B$ induit un~$X$-morphisme
$\spec \sch B\to \spec \sch A$ (le définir au-dessus des ouverts affines de~$X$ 
par fonctorialité contravariante du spectre classique, et recoller). 

\medskip
La flèche~$\sch A\mapsto \spec \sch A$ apparaît ainsi
de manière naturelle comme 
un foncteur contravariant 
de la catégorie des~$\sch O_X$-algèbres quasi-cohérentes vers celle des~$X$-schémas.

\trois{algqc-pleine-fid}
Soient~$\sch A$ et~$\sch B$ deux~$\sch O_X$-algèbres quasi-cohérentes
et soient~$p$ et~$q$
les morphismes respectifs de~$\spec \sch A$ et~$\spec \sch B$ vers~$X$. 
Soit~$\lambda$
l'application naturelle $$\lambda\colon \hom_{\sch O_X\text{-}\mathsf{Alg}}(\sch A,\sch B)\to \hom_X(\spec \sch B,\spec \sch A)$$
définie au~\ref{algqc-contrav}
ci-dessus. 
Soit~$\phi \colon \spec \sch B\to \spec \sch A$
un~$X$-morphisme. Il induit un morphisme~$\sch O_{\sch \spec \sch A}\to \phi_*\sch O_{\spec \sch B}$ 
puis, par application de~$p_*$, un morphisme

$$\sch A\simeq p_*\sch O_{\spec \sch A}\to p_*\circ \phi_*\sch O_{\spec \sch B}=q_*\sch O_{\spec \sch B}\simeq \sch B.$$

On obtient ainsi une application~$\mu$ de~$\hom_X(\spec \sch B,\spec \sch A)$
vers~$\hom_{\sch O_X\text{-}\mathsf{Alg}}(\sch A,\sch B)$, et l'on vérifie aisément que~$\lambda$
et~$\mu$ sont 
des bijections
réciproques l'une de l'autre (c'est une propriété locale sur~$X$, ce
qui permet de se ramener au cas où tout le monde est affine,
dans lequel c'est une reformulation de~\ref{fonct-spec-locann}). 

\deux{y-specpioy}
Soit~$\pi \colon Y\to X$ un morphisme 
de schémas. 

\trois{constr-y-specpioy}
Soit~$U$ un ouvert affine de~$X$ ; on dispose
d'un morphisme naturel~$\pi^{-1}(U)\to \spec \sch O_Y(\pi^{-1}(U))$ et
pour tout ouvert affine~$V$ de~$U$, d'un diagramme commutatif

$$
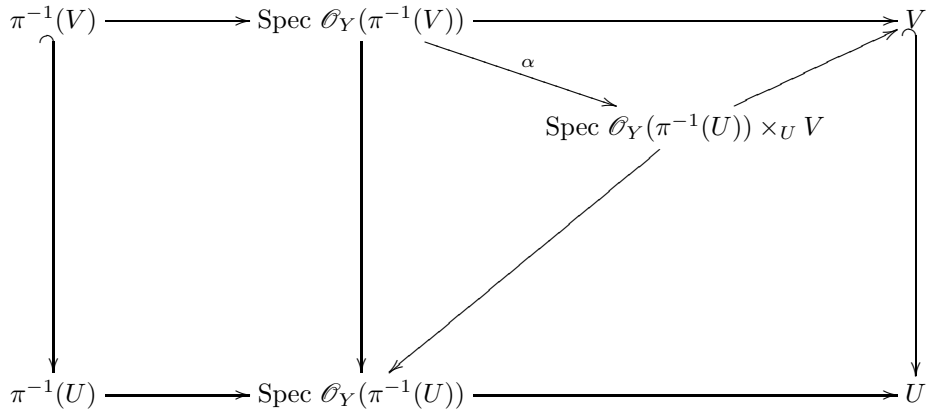
\xymatrix{
{\pi^{-1}(V)}\ar@{_{(}->}[dddd]\ar[rr]&&
{\spec \sch O_Y(\pi^{-1}(V))}\ar[dddd]\ar[rr]\ar[rd]^\alpha&&V\ar@{_{(}->}[dddd]
\\&&&{\spec \sch O_Y(\pi^{-1}(U))\times_U V}\ar[ru]\ar[lddd]&
\\
\\
\\
{\pi^{-1}(U)}\ar[rr]&&{\spec \sch O_Y(\pi^{-1}(U))}\ar[rr]&&U
}$$
dans lequel la flèche~$\alpha$ est un isomorphisme dès que
la~$\sch O_X$-algèbre~$\pi^*\sch O_Y$
est quasi-cohérente. {\em Sous cette dernière hypothèse}, 
on dispose donc d'un diagramme commutatif 
$$
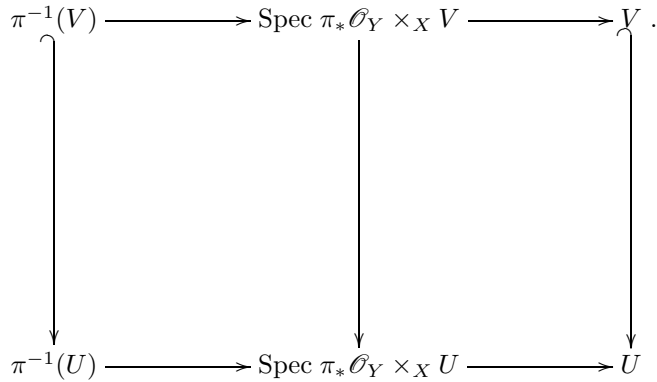
\xymatrix{
{\pi^{-1}(V)}\ar@{_{(}->}[dddd]\ar[rr]&&
{\spec \pi_*\sch O_Y\times_X V}\ar[dddd]\ar[rr]&&V\ar@{_{(}->}[dddd]
\\
\\
\\
\\
{\pi^{-1}(U)}\ar[rr]&&{\spec \pi_*\sch O_Y\times_X U}\ar[rr]&&U
}.$$

Par recollement de ces diagrammes pour~$(U,V)$ variables, on obtient
un morphisme de~$X$-schémas~$Y\to \spec \pi_*\sch O_Y$.

\trois{cas-y-specpioy}
Supposons que~$Y=\spec \sch A$ pour une certaine~$\sch O_X$-algèbre quasi-cohérente~$\sch A$. Dans
ce cas~$\pi_*\sch O_Y$ est quasi-cohérente et s'identifie plus précisément à~$\sch A$ (\ref{recollement-spec-au}) ; la construction
du~\ref{constr-y-specpioy}
ci-dessus fournit alors un
morphisme
de~$X$-schémas~$\spec \sch A\to \spec \sch A$, dont on vérifie aussitôt que c'est l'identité. 

\deux{def-morph-aff}
{\bf Proposition-définition.}
{\em Soit~$\pi \colon Y\to X$ un morphisme de schémas. Les assertions suivantes sont équivalentes. 

\medskip
i) La~$\sch O_X$-algèbre~$\pi_*\sch O_Y$ est quasi-cohérente, et~$Y\to \spec \pi_*\sch O_Y$ est un isomorphisme. 

ii) Il existe une~$\sch O_X$-algèbre quasi-cohérente~$\sch A$ et un~$X$-isomorphisme
$$Y\simeq \spec \sch A.$$

iii) Pour tout ouvert affine~$U$ de~$X$, le schéma~$\pi^{-1}(U)$ est affine. 

iv) Il existe un recouvrement~$(U_i)$ de~$X$
par des ouverts affines tels que~$\pi^{-1}(U_i)$
soit un schéma affine pour tout~$i$.

\medskip
Lorsqu'elles sont satisfaites, on dit que~$\pi$ est {\em affine},
ou que~$Y$ est {\em relativement affine}
sur~$X$.}

\medskip
{\em Démonstration.}
Il est clair que~i)$\Rightarrow$ii), et~ii)$\Rightarrow$i)
d'après~\ref{cas-y-specpioy}. Il découle
de la construction
même de~$\spec \sch A$ que~ii)$\Rightarrow$iii) (on l'a déjà signalé en~\ref{recollement-spec-au}),
et~iii)$\Rightarrow$iv)
est évident. 

\medskip
Supposons maintenant que~iv) est vraie, et montrons que les conditions équivalentes~i) et~ii)
sont vérifiées. Elles sont de nature locale sur~$X$ (c'est l'énoncé~i) qui le montre) ; il suffit 
par conséquent
de démontrer
qu'elles sont vraies sur chaque~$U_i$. Fixons donc~$i$. Par hypothèse, $U_i$ et~$\pi^{-1}(U_i)$ sont
affines ; il s'ensuit, 
en vertu de~\ref{spec-relatif-casaff},
que~$\pi^{-1}(U_i)\to U_i$ satisfait~ii) (et partant~i)),
ce qui achève la démonstration.~$\Box$

\deux{comment-morph-aff}
{\bf Commentaires.}
L'aspect le plus spectaculaire de la proposition précédente
est l'équivalence entre~iii)
et~iv) : s'il existe un recouvrement 
de~$X$ par des ouverts affines dont l'image réciproque par~$\pi$ 
est affine, 
alors~$\pi^{-1}(U)$ est affine
pour {\em tout}
ouvert affine~$U$ 
de~$X$. 

Le lecteur sera peut-être étonné qu'un résultat aussi fort ait une
preuve aussi courte et d'apparence très formelle. Mais 
lorsqu'on invoque
le caractère local de~i), on invoque en particulier le caractère local de la
quasi-cohérence ; et ce dernier 
est 
lui-même fondé sur
le théorème~\ref{theo-qc-loc} 
dont la démonstration met en jeu des
arguments non triviaux,  à base de calcul de fractions dans les modules localisés. 
C'est donc là qu'il se «passe vraiment quelque chose». 

\deux{reformu-moraff}
On peut reformuler~\ref{recollement-spec-au},
\ref{algqc-contrav}
et~\ref{algqc-pleine-fid}
en disant que pout tout schéma~$X$, le foncteur~$\sch A\mapsto \spec \sch A$
établit une anti-équivalence entre la catégorie des~$\sch O_X$-algèbres quasi-cohérentes et celle
des~$X$-schémas relativement affines, dont~$(\phi \colon Y\to X) \mapsto \phi_*\sch O_Y$ 
est un quasi-inverse.

\deux{prod-fib-aff}
{\bf Stabilité du caractère affine par composition et changement de base.}

\trois{aff-compos}
Soient~$\phi \colon Z\to Y$ et~$\psi \colon Y\to X$ deux morphismes affines. 
La composée~$\psi \circ \phi$
est alors affine : c'est immédiat en utilisant la condition
équivalente~iii) de la proposition~\ref{def-morph-aff}
ci-dessus. 

\trois{aff-changebase}
Soit~$X$ un schéma. Soit~$\phi \colon Y\to X$ un morphisme
{\em affine}, et soit~$\psi \colon Z\to X$ un morphisme. La projection
$\pi \colon Y\times_X Z\to Z$ est alors affine. Plus
précisément si~$Y=\spec \sch A$ pour une certaine
$\sch O_X$-algèbre quasi-cohérente~$\sch B$,
le~$Z$-schéma~$Y\times_X Z$ 
s'identifie alors à~$\spec \psi^*\sch B$. 

\medskip
En effet, la question est locale sur~$Z$, et {\em a fortiori}
sur~$X$. On peut donc supposer tout d'abord~$X$ affine, puis~$Z$ affine ; comme~$\phi$
est affine, $Y$ est affine. Soient~$A, B$ et~$C$
les anneaux correspondant respectivement aux schémas affines~$X,Y$ et~$Z$. 
On a alors~$\sch B= \red B$, 
et~
$$Y\times_X Z=\spec (B\otimes_A C)=\spec \red{B\otimes_A C}=\spec \psi^* \red B=\spec \psi^*\sch B.$$

\subsection*{Les immersions fermées}
Nous allons maintenant présenter une première classe
absolument fondamentale de morphismes affines : les {\em immersions fermées}.

\deux{def-imm-ferm}
Soit~$X$ un schéma et soit~$\phi : Y\to X$ un morphisme. Les conditions suivantes sont équivalentes : 

\medskip
i) le morphisme~$\phi$
est affine et le morphisme structural~$\sch O_X\to \phi_*\sch O_Y$ est surjectif ; 

ii) il existe un faisceau quasi-cohérent d'idéaux~$\sch I$ sur~$X$ tel que
le~$X$-schéma~$Y$
soit isomorphe à~$\spec \sch O/\sch I$. 

\medskip
En effet si~i) est vraie, le noyau~$\sch I$ 
de~$\sch O_X\to \phi_*\sch O_Y$
est un faisceau quasi-cohérent d'idéaux et~$\phi_*\sch O_Y\simeq \sch O_X/\sch I$,
d'où~ii) puisque~$Y=\spec \phi_*\sch O_Y$ ; et si~ii) est vraie on
a~$\phi_*\sch O_Y=\sch O/\sch I$, d'où~i). 

\medskip
Lorsque ces conditions sont satisfaites, on dit que~$\phi$ est une 
{\em immersion fermée.}
Notons que l'idéal~$\sch I$ de~ii) est alors uniquement déterminé : il est
{\em nécessairement} égal au noyau de~$\sch O_X\to \phi_*\sch O_Y=\sch O/\sch I$. 
Nous dirons que~$\sch I$ est {\em le faisceau d'idéaux associé à~$\phi$},
et inversement que~$\phi$ est {\em l'immersion fermée associée à~$\sch I$}. 

\deux{prem-prop-imm-f}
{\bf Exemples et premières propriétés.} 

\trois{imm-f-local}
Soit~$\phi \colon Y\to X$
un morphisme de schémas et soit~$(U_i)$ un recouvrement ouvert de~$X$. 
On déduit de la caractérisation d'une immersion fermée
par la propriété~i) ci-dessus 
que~$\phi$ est une immersion fermée si et seulement si~$\phi^{-1}(U_i)\to U_i$ est une 
immersion fermée pour tout~$i$.

\trois{imm-f-affine}
Soit~$A$ un anneau. La caractérisation des
immersions fermées par la propriété~ii)
ci-dessus assure qu'un morphisme~$\phi \colon Y\to \spec A$ est une immersion
fermée si et seulement si le~$A$-schéma~$Y$ est de la forme~$\spec A/I$
pour un certain idéal~$I$ de~$A$ ; le faisceau~$\sch I$ associé à~$\phi$ 
est alors égal à~$\red I$. 
On voit que dans ce cas, $Y\to \spec A$ induit un homéomorphisme entre~$Y$ et
le fermé~$V(I)=V(\sch I)$ de~$\spec A$. 

\trois{imm-f-ferme}
Soit~$\phi \colon Y\to X$ une immersion fermée et soit~$\sch I\subset \sch O_X$ le faisceau
d'idéaux associés. Le morphisme~$\phi$ induit un homéomorphisme~$Y\simeq V(\sch I)$ : cette assertion
est en effet locale, ce qui permet de se ramener au cas affine traité au~\ref{imm-f-affine}
ci-dessus.

\trois{imm-f-ferme-moinsun}
Soit~$\phi \colon Y\to X$ une immersion fermée et soit~$\sch I\subset \sch O_X$ le faisceau
d'idéaux associés. Le morphisme naturel~$\phi^{-1}(\sch O_X/\sch I)\to \sch O_Y$ est alors un
isomorphisme. Pour le voir, on raisonne localement sur~$X$, ce qui permet de se ramener au cas affine. 
Le caractère bijectif de la flèche étudiée se vérifie alors fibres à fibres, et revient à l'assertion suivante
d'algèbre commutative : soit~$A$ un anneau, soit~$I$ un idéal de~$A$ et soit~$\got p$ un idéal premier
de~$A$ contenant~$I$ ; la flèche
canonique~$A_{\got p}/IA_{\got p}\to (A/I)_{\got p/I}$
est un isomorphisme.

\trois{imm-f-local-conseq}
Soit~$\phi \colon Y \to X$ un morphisme de schémas
et soit~$(U_i)$ un recouvrement de~$X$ par des ouverts affines. 
Il découle de~\ref{imm-f-local}
et~\ref{imm-f-affine}
que pour que~$\phi$ soit une immersion fermée, il 
faut et il suffit que 
pour tout~$i$, le~$U_i$-schéma~$\phi^{-1}(U_i)$ soit de la forme~$\spec \sch O_X(U_i)/I_i$
pour un certain idéal~$I_i$ de~$\sch O_X(U_i)$ ; et
que si c'est le cas alors pour
{\em tout}
ouvert affine~$U$
de~$X$, le~$U$-schéma $\phi^{-1}(U)$ est de la forme~$\spec \sch O_X(U)/I$
pour un certain idéal~$I$ de~$\sch O_X(U)$. Le lecteur s'assurera que ce
passage des seuls ouverts d'un
recouvrement affine 
fixé à {\em tous}
les ouverts affines
repose {\em in
fine}
là encore sur le caractère local de la quasi-cohérence (théorème~\ref{theo-qc-loc}). 

\deux{prop-uni-immf}
{\bf Propriété universelle d'une
immersion fermée.}
Soit~$\phi \colon Y \to X$ une immersion fermée, et soit~$\sch I\subset \sch O_X$ 
le faisceau d'idéaux correspondants. Le faisceau~$\sch I$ est inclus dans
le (et même égal au) noyau de~$\sch O_X\to \phi_*\sch O_Y$. 

\medskip
\trois{prop-uni-immf-alamain}
Soit~$Z$ un schéma et soit~$\psi \colon Z \to X$ un morphisme tel que~$\sch I$ soit 
contenu dans le noyau de~$\sch O_X\to \psi_*\sch O_Z$. Il existe alors un unique
morphisme~$\chi \colon Z\to Y$ tel que le diagramme 
$$\xymatrix{
Z\ar[rrd]^\psi\ar[rd]_\chi&&
\\&Y\ar@{^{(}->}[r]_\phi&X}$$

commute. 
En effet, comme~$U\mapsto \hom(U,Y)$ est un faisceau sur~$Z$, on peut raisonner localement, 
et donc supposer~$X$ et~$Z$ affines ; soient~$A$ et~$B$ les anneaux correspondants. Comme~$\phi$ est l'immersion
fermée associée à~$\sch I$, le~$A$-schéma
$Y$ est égal à~$\spec A/I$, où~$I$ est l'idéal~$\sch I(X)$ de~$A$. Notre hypothèse sur~$\psi$ signifie 
simplement
que~$I\subset {\rm Ker}\;(A\to B)$, et le résultat voulu est une simple reformulation du fait que~$A\to B$
se factorise alors de manière unique par~$A/I$.

\trois{prop-uni-immf-enonce}
On peut résumer conceptuellement ce
qui précède en disant que~$(Y,\phi)$
représente le foncteur
$$Z\mapsto \{\psi \in \hom(Z,X), \;\;{\sch I}\subset {\rm Ker}\;(\sch O_X\to \psi_*\sch O_Z)\}.$$

\deux{prod-fib-immf}
{\bf Bon comportement des immersions fermées par composition et changement de base.}

\trois{composition-immf}
Soient~$\psi \colon Z\to Y$ et~$\phi \colon Y\to X$
deux immersions fermées. La composée
~$\phi \circ \psi$ est alors une immersion fermée  : c'est
par exemple une
conséquence immédiate de~\ref{imm-f-local-conseq}
et du fait que la composée de deux morphismes d'anneaux surjectifs est encore
une surjection. 

\trois{changebase-immf}
Soit~$\phi \colon Y\to X$ une immersion
fermée induite par un faisceau quasi-cohérent
d'idéaux~$\sch I$ sur~$X$, et soit~$\psi \colon Z\to X$
un morphisme de schémas. 

\medskip
La flèche~$\sch I\hookrightarrow \sch O_X$ induit
une flèche~$\psi^*\sch I\to \sch O_Z$ (qui n'a pas
de raison d'être injective) ; son image est un
faisceau quasi-cohérent d'idéaux~$\sch J$ de~$\sch O_Z$,
et il résulte des définitions que~$V(\sch J)=\psi^{-1}(V(\sch I))$.

\medskip
Il découle de
la description du foncteur représenté par~$(Y,\phi)$ ({\em cf.}~\ref{prop-uni-immf-enonce})
que le produit fibré~$Y\times_XZ$ représente le foncteur

$$T\mapsto \{\chi \in \hom(T,Z), \sch I\subset {\rm Ker}\; (\sch O_X\to \psi_*\chi_*\sch O_T)\}.$$
La condition~$\sch I\subset {\rm Ker}\; (\sch O_X\to \psi_*\chi_*\sch O_T)$ équivaut à dire
que la flèche composée~$\sch I\to \sch O_X\to \psi_*\chi_*\sch O_T$ est nulle, c'est-à-dire encore
par adjonction que la flèche composée~$\psi^*\sch I\to \sch O_Z\to \chi_*\sch O_T$ est nulle ; 
mais c'est le cas si et seulement si~$\sch J\subset {\rm Ker}\;(\sch O_Z\to \chi_*\sch O_T)$. En conséquence, 
$Y\times_X Z\to Z$ est l'immersion fermée associée à~$\sch J$. 

\medskip
On peut donner une deuxième preuve moins yonedesque de cette assertion. Il suffit de se ramener en raisonnant
localement au cas où tous les schémas en jeu sont affines, et de remarquer qu'elle est alors une simple reformulation 
du fait que si~$A$ est un anneau, $B$ une~$A$-algèbre et~$I$ un idéal de~$A$
alors~$B\otimes_A A/I\simeq B/IB$. 

\medskip
Mentionnons pour conclure ce paragraphe un cas particulier intéressant : supposons
que~$\psi \colon Z\to X$ se factorise par~$Y\to X$ ; d'après la propriété universelle d'une immersion
fermée assure, 
cela signifie que~$\sch I$
est contenu dans~${\rm Ker}\;(\sch O_X\to \psi_*\sch O_Z)$, et cette
factorisation est alors unique. 

Comme~$\sch I$ est contenu dans~${\rm Ker}\;(\sch O_X\to \psi_*\sch O_Z)$,
la flèche~$\psi^*\sch I\to \sch O_Z$ est nulle par adjonction, et $\sch J$ est dès lors nul. 
La projection~$Y\times_X Z\to Z$ s'identifiant à l'immersion fermée
définie par~$\sch J$, il vient
$Y\times_X Z=Z$.

\deux{prod-fib-deuximmf}
{\bf Produit fibré de deux immersions fermées.}
Soient~$\phi \colon Y\to X$ et~$\psi \colon Z\to X$ deux immersions 
fermées, respectivement associées à des faisceaux d'idéaux~$\sch I$ et~$\sch J$.

\trois{desc-prodfib-vi-vj}
On déduit de la description des foncteurs représentés par~$(Y,\phi)$ et~$(Z,\psi)$ que
$Y\times_XZ$ représente le foncteur
$$T\mapsto \{\chi \in \hom(T,X), \sch I\subset {\rm Ker}\;(\sch O_X\to \chi_*\sch O_T)\;\;\text{et}\;\,\sch J\subset {\rm Ker}\;(\sch O_X\to \chi_*\sch O_T)\}$$
$$=\{\chi \in \hom(T,X), \sch I+\sch J\subset {\rm Ker}\;(\sch O_X\to \chi_*\sch O_T)\},$$
où la somme~$\sch I+\sch J$ est par définition le
faisceau cohérent d'idéaux égal à l'image
de la flèche canonique~$\sch I\oplus \sch J\to \sch O_X$. En conséquence, $Y\times_XZ\to X$ est l'immersion
fermée associée à~$\sch I+\sch J$. 

\trois{fact-imm-ferm}
L'immersion~$\psi$ se factorise par~$\phi$
si et seulement si~
$$\sch J\subset {\rm Ker}\;(\sch O_X\to \phi_*\sch O_Y)={\rm Ker}\;(\sch O_X\to \sch O_X/\sch I)=\sch I.$$
Si c'est le cas, on déduit de la remarque faite à la fin de~\ref{changebase-immf}
que~$Z\times_X Y\simeq Z$ ; il découle alors de {\em loc. cit.}
que~$Z\to Y$ est l'immersion fermée
associée au faisceau quasi-cohérent d'idéaux sur~$Y$ 
engendré par~$\phi^*\sch J$. 

\deux{sous-schem-ferm}
{\bf La notion de sous-schéma fermé.}
Soit~$X$ un schéma, soit~$F$ un fermé de~$X$
et soit~$\sch I\subset \sch O_X$ un faisceau quasi-cohérent d'idéaux
tel que~$V(\sch I)$
soit égal à~$F$ (on sait qu'il en existe 
au moins un, 
{\em cf}.~\ref{retour-ferme-viqc}). 

\medskip
L'immersion fermée~$\spec \sch O_X/\sch I\to X$ induit un homéomorphisme
entre~$\spec \sch O_X/\sch I$ et~$F$, et permet donc de munir par transport 
de structure l'espace topologique~$F$ d'une structure de schéma ; un tel schéma
est appelé le {\em sous-schéma fermé 
de~$X$ défini par~$\sch I$}, et l'on dit que~$F$ est son
{\em support}. 

\medskip
On a signalé plus haut 
(rem.~\ref{rem-support-oxi})
que la restriction de~$\sch O_X/\sch I$ 
à~$X\setminus F$ est nulle, et que cela se traduit
en disant que~$\sch O_X/\sch I$ s'identifie à~$i_*i^{-1}\sch O_X/\sch I$, où~$i$
est l'inclusion de~$F$ dans~$X$. Le faisceau~$\sch O_X/\sch I$ provient donc
du faisceau d'anneaux~$i^{-1}\sch O_X/\sch I$ sur~$F$, et il résulte de~\ref{imm-f-ferme-moinsun}
que la structure de schéma dont on a muni~$F$ est précisément
induite par~$i^{-1}\sch O_X/\sch I$. 

\deux{intro-struct-reduite}
Soit~$X$ un schéma. L'ensemble des (classes d'isomorphie de)
sous-schémas fermés de~$X$ est en bijection avec l'ensemble des faisceaux-quasi-cohérents
d'idéaux sur~$X$. On munit l'ensemble des sous-schémas fermés de~$X$ de la relation d'ordre
pour laquelle~$F\leq G$ si et seulement si~$F\hookrightarrow X$ se factorise
par~$G\hookrightarrow X$. Si~$\sch I$ et~$\sch J$ désignent les faisceaux
d'idéaux respectivement associés à~$F$ et~$G$, cela revient à demander
que~$\sch J\subset \sch I$ (\ref{fact-imm-ferm}).

\trois{struct-reduite}
{\em La structure réduite.}
Si~$F$ est un fermé de~$X$, il existe un plus petit sous-schéma fermé~$F_{\rm red}$
de~$X$ de support~$F$ : celui qui est induit par~$\sch I(F)$ ({\em cf.}-~\ref{retour-ferme-viqc}
{\em et sq.}). Nous allons montrer que~$F_{\rm red}$ est réduit, et que c'est le seul sous-schéma fermé
réduit de support~$F$. 

\medskip
Fixons un recouvrement~$(U_i)$ de~$X$ par des ouverts affines et soit~$\sch J$ un faisceau quasi-cohérent d'idéaux de support~$F$.
Le sous-schéma
fermé défini par~$\sch J$ est réduit si et seulement si~$\sch O_X(U_i)/\sch J(U_i)$ est réduit pour tout~$i$, ce qui revient
à demander que l'idéal $\sch J(U_i)$ de~$\sch O_X(U_i)$ soit saturé pour tout~$i$ ; mais cela signifie précisément
que~$\sch J=\sch I(F)$ (\ref{schi-f-i}). 

\trois{prop-univ-fred}
Soit~$T$ un schéma et soit~$\phi \colon T \to X$
un morphisme. Si~$\phi$
admet une factorisation (nécessairement unique)
par~$F_{\rm red}$ alors~$\phi(T)\subset F$.

\medskip
Faisons maintenant l'hypothèse que~$\phi(T)\subset F$,
et supposons de surcroît 
que~$T$ est {\em réduit}. 
Nous allons montrer que~$\phi$ se factorise par~$F_{\rm red}$. 

Soit~$U$ un ouvert de~$X$ et soit~$f\in \sch I(F)(U)$.
Comme~$\phi(T)\subset F$, 
la fonction~$\phi^*f\in \sch O_T(\phi^{-1}(U))$ 
s'annule en tout point de~$U$. Cela implique
que sa
restriction à tout ouvert affine de~$\phi^{-1}(U)$ est nilpotente, donc nulle {\em puisque~$T$ est réduit}. 
Autrement dit, $\phi^*f=0$, ce qui revient à dire 
que l'image de~$f$ dans~$\phi_*\sch O_T$ est nulle ;
ainsi, $\sch I(F)\subset  {\rm Ker}\;(\sch O_X\to \phi_*\sch O_T)$,
ce qui garantit que~$\phi$ se factorise par~$F_{\rm red}$. 

\medskip
En d'autres termes,~$(F_{\rm red}, F_{\rm red}\hookrightarrow X)$
représente le foncteur covariant de la catégorie des schémas
{\em réduits}
vers celle des ensembles qui envoie~$T$ sur
$$\{\phi\in \hom(T,X), \phi(T)\subset F\}.$$

\trois{x-red}
On peut appliquer ce qui précède lorsque~$F=X$.
Soit~$\sch J$
l'idéal des fonctions s'annulant en tout point de~$X$ ; ce sont exactement les fonctions
dont la restriction à tout ouvert affine est nilpotente -- on vérifie aussitôt
que cela vaut encore pour leur restriction à tout ouvert quasi-compact, et 
nous vous laissons construire un exemple d'une telle fonction qui ne serait pas elle-même nilpotente, sur un schéma
non quasi-compact. 

Le schéma~$X_{\rm red}$
est par définition le sous-schéma fermé défini par~$\sch J$ ; son espace
topologique sous-jacent est~$X$ tout entier ; on dit que c'est le {\em schéma réduit associé à~$X$}. 
Si~$T$ est un schéma réduit, tout morphisme de~$T$ vers~$X$ se factorise
canoniquement
par~$X_{\rm red}$. 

\deux{comment-sous-schem}
Soit~$X$ un schéma et soit~$F$ un fermé
de~$X$. On prendra garde qu'en général, il n'existe pas de structure 
de sous-schéma fermé sur~$F$ telle que tout morphisme
$T\to X$ se factorisant {\em ensemblistement}
par~$F$ se factorise par la structure en question  -- la
propriété universelle de~$F_{\rm red}$
décrite au~\ref{prop-univ-fred}
ne concerne que les schémas {\em réduits}.

\medskip
En effet, supposons qu'il existe une telle structure, et soit~$Y$
le sous-schéma fermé correspondant. Si~$Z$ est un autre sous-schéma fermé de~$X$
de support~$F$, l'image de~$Z\hookrightarrow X$ est contenue dans~$F$ et se factorise
donc par~$Y$ ; autrement dit, $Y$ est nécessairement {\em le plus grand sous-schéma fermé}
de support~$F$. Or l'ensemble des sous-schémas fermés de support~$F$ n'a pas, en général, de plus
grand élément ; ou, si l'on préfère, l'ensemble des faisceaux quasi-cohérents d'idéaux de lieu des zéros~$F$
n'a pas forcément de plus petit élément.

\trois{pas-de-plugrand-scf}
Ainsi, soit~$p$ un nombre premier. Un faisceau quasi-cohérent d'idéaux sur~$\spec \ZZ$ a pour lieu des zéros
le singleton~$\{x_p\}$ si et seulement
si il est de la forme~$\sch I_n:=\red{(p^n)}$ pour un certain~$n>0$ ; la famille~$(\sch I_n)$ étant strictement décroissante, 
elle n'a pas de plus petit élément.

\trois{un-plusgrand-scf}
Il peut toutefois arriver qu'il existe un plus petit idéal quasi-cohérent de lieu des zéros~$F$. Par exemple, supposons que le fermé
$F$ soit également ouvert, et soit~$G$ l'ouvert fermé complémentaire. Soit~$f\in \sch O_X(X)$ la fonction telle que~$f|_G=1$ et~$f|_F=0$. 
Soit~$\sch I$ l'idéal quasi-cohérent égal à l'image de la flèche~$\sch O_X\to \sch O_X, a\mapsto af$. Nous laissons le lecteur vérifier que~$\sch I$
est le plus petit idéal quasi-cohérent de lieu des zéros égal à~$F$, et que la structure de sous-schéma fermé qu'il définit sur~$F$ est sa structure 
d'ouvert, par laquelle se factorise tout morphisme dont l'image ensembliste est contenue dans~$F$. 

\deux{sous-schf-intuition}
Soit~$X$ un schéma et soit~$F$ un fermé de~$X$. Intuitivement les différentes structures de sous-schéma fermé
de support~$F$ codent
les différentes manières d'envisager~$F$ «avec multiplicités» (la structure~$F_{\rm red}$ correspondant au cas sans multiplicités) ; 
ou, si l'on préfère,
les différentes manières d'épaissir~$F$ infinitésimalement
à l'intérieur de~$X$, la relation~$Y\leq Z$ entre
deux sous-schémas fermés de support~$F$
pouvant alors 
s'interpréter comme «$Y$ est moins épais que~$Z$». 

\deux{ex-schf}
{\bf Exemples}. 
Soit~$k$ un corps. Posons~$X=\spec k[S,T]$, et soit~$F$ le fermé~$V(S)$ de~$X$. Nous allons
décrire différentes structures de sous-schéma fermé sur~$F$. Comme~$X$ est affine, 
un faisceau cohérent d'idéaux 
de~$X$ de
lieu des zéros~$F$ est simplement un idéal~$I$ de~$k[S,T]$ tel que~$V(I)=F$. 

\trois{ex-schf-red}
{\em La structure réduite}.
On a~$V(S)=F$. 
L'anneau~$k[S,T]/S\simeq k[T]$ est intègre, et {\em a fortiori}
réduit. La structure induite sur~$F$ par l'idéal~$(S)$
est donc la structure réduite
$F_{\rm red}$. 

\medskip
{\em Remarque.}
Comme~$k[T]$ est principal les ouverts de~$F_{\rm red}$
sont exactement les~$D(P)$ pour~$P\in k[T]$ ; et
cela vaut également
pour n'importe quel sous-schéma fermé
de~$X$ de support~$F$, puisque c'est une assertion purement topologique.

\trois{ex-sch-pasred1}
Soit~$I$ l'idéal~$(S^2)$. On a~$V(I)=F$, et~$I$ induit donc
une structure de schéma~$F_1$
sur~$F$. 
L'anneau quotient~$k[S,T]/S^2$ possède un élément
nilpotent d'ordre 2, à savoir~$\overline S$. La structure correspondante n'est en conséquence pas réduite. 
Moralement, la droite~$F$ a été 
un peu épaissie pour devenir une «droite double», 
et ce de façon relativement uniforme : si~$P$
est un élément de~$k[T]$ 
on a en effet
$$\sch O_{F_1}(D(P))=k[S,T]_{(P)}/S^2=k[T]_{(P)}[S]/S^2=\sch O_{F_{\rm red}}(D(P))[S]/S^2.$$
On voit notamment que la restriction de~$\overline S$ à tout ouvert non vide de~$F_1$ est encore nilpotente
d'ordre 2.

\trois{compos-imm}
{\em Composantes immergées}.
Il peut exister des façons
plus subtiles d'épaissir~$F$. Par exemple, 
soit~$J$ l'idéal~$(S^2, ST)$. On a~$V(J)=F$, et~$J$
induit donc une structure de schéma~$F_2$ sur~$F$. 
Le quotient
$k[S,T]/(S^2, ST)$ n'est pas réduit : la fonction~$\overline S$ est nilpotente d'ordre~$2$. 
Le schéma~$F_2$ est en conséquence un épaississement de~$F$, mais moins «uniforme» que~$F_1$. 
On a en effet
$$\sch O_{F_2}(D(T))=k[S,T]_{(T)}/(S^2, ST)$$
$$=\underbrace{k[T]_{(T)}[S]/(S^2,ST)=k[T]_{(T)}[S]/S}_{
\text{car}\;T\;\text{est inversible dans}\;k[T]_{(T)}}=k[T]_{(T)}=\sch O_{F_{\rm red}}(D(T)).$$
L'ouvert dense
$D(T)$ de~$F_2$ est ainsi {\em réduit} :  l'épaississement disparaît dès qu'on retire l'origine, c'est donc elle qui en un sens
porte toute la multiplicité de la situation. On dit que l'origine est une {\em composante immergée}
du schéma~$F_2$. 

\medskip
{\em Remarque.}
Notez bien que la restriction de la fonction nilpotente non nulle~$\overline S$ à l'ouvert réduit~$D(T)$ de~$F_2$ est nulle : contrairement
à ce qu'on pourrait croire naïvement, une fonction sur un schéma peut s'annuler en restriction à un ouvert dense sans être nulle. 

\trois{compos-imm-subt}
Soit~$J'$ l'idéal ~$(S^3, S^2T)$. On a~$V(J')=F$, et~$J'$
induit donc une structure de schéma~$F_3$ sur~$F$. 
La fonction~$\overline S$ est alors nilpotente d'ordre~3
sur~$F_3$. Par ailleurs
$$\sch O_{F_3}(D(T))=k[S,T]_{(T)}/(S^3, S^2T)$$
$$=\underbrace{k[T]_{(T)}[S]/(S^3,S^2T)=k[T]_{(T)}[S]/S^2}_{
\text{car}\;T\;\text{est inversible dans}\;k[T]_{(T)}}=\sch O_{F_1}(D(T)).$$
Ainsi, la restriction de~$\overline S$ à l'ouvert~$D(T)$ de~$F_3$ est désormais
nilpotente {\em d'ordre 2 seulement} (et il en ira de même de sa restriction à n'importe quel
ouvert non vide de~$D(T)$, d'après~\ref{ex-sch-pasred1}). Cet exemple combine donc 
les phénomènes décrits en~\ref{ex-sch-pasred1}
et~\ref{compos-imm} : 
on peut y penser comme à un épaississement global du fermé~$F$, de multiplicité «générique» égale à~$2$, 
possédant un surcroît de multiplicité porté par l'origine. 

\subsection*{Morphismes finis}

\deux{def-qc-loctf}
{\bf Lemme.}
{\em Soit~$X$ un schéma et soit~$\sch F$ un~$\sch O_X$-module quasi-cohérent. Les assertions suivantes
sont équivalentes : 

\medskip
i) pour tout ouvert affine~$U$ de~$X$, le~$\sch O_X(U)$-module~$\sch F(U)$ est de type fini ; 

ii) il existe un recouvrement~$(U_i)$ de~$X$ par des ouverts affines tels que le~$\sch O_X(U_i)$-module
$\sch F(U_i)$ soit de type fini pour tout~$i$.}

\medskip
{\em Démonstration.} 
Il est évident que~i)$\Rightarrow$ii). Supposons que~ii) soit vraie, 
et soit~$U$ un ouvert affine de~$X$ ; posons~$A=\sch O_X(U)$. Soit~$x\in U$ ; il existe~$i$
tel que~$x\in U_i$, et il existe donc~$f\in A$ tel que~$f(x)\neq 0$ et tel que l'ouvert~$D(f)$ de~$U$
soit contenu dans~$U_i$. On a alors~$\sch F(D(f))=\sch O_X(D(f))\otimes_{\sch O_X(U_i)}\sch F(U_i)$ 
(puisque~$\sch F$ est quasi-cohérent), et~$\sch F(D(f))$ est donc un~$A_f$-module de type fini, 
qui est égal à~$\sch F(U)_f$, là encore par quasi-cohérence de~$\sch F$.

\medskip
Par quasi-compacité de~$U$ il existe~$f_1,\ldots, f_r\in A$ telles que~$U=\bigcup D(f_i)$ ou, 
si l'on préfère, telles que~$(f_1,\ldots, f_n)=A$, et telles que~$\sch F(U)_{f_i}$ soit un~$A_{f_i}$-module
de type fini pour tout~$i$. Le lemme~\ref{type-fini-si}
assure alors que~$\sch F(U)$ est un~$A$-module de type fini.~$\Box$ 

\deux{def-morph-fin}
{\bf Définition.}
Un morphisme de schéma~$\phi \colon Y\to X$ est dit {\em fini}
s'il est affine et si le~$\sch O_X$-module quasi-cohérent~$\phi_*\sch O_Y$ satisfait les conditions équivalentes
du lemme~\ref{def-qc-loctf}. On dira également si c'est le cas que~$Y$ est un {\em $X$-schéma fini}. 

\deux{reform-morph-fin}
Le lemme~\ref{def-qc-loctf}
(combiné à la proposition~\ref{def-morph-aff})
assure que pour tout morphisme de schémas~$\phi \colon Y\to X$, les conditions suivantes sont équivalentes :

\medskip
i) $\phi$ est fini ; 

ii) pour tout ouvert affine~$U$ de~$X$, le schéma~$\phi^{-1}(U)$ est affine et~$\sch O_Y(\phi^{-1}(U))$ est une~$\sch O_X(U)$-algèbre finie ; 

iii) il existe un recouvrement~$(U_i)$ de~$X$ par des ouverts affines tels que pour tout~$i$, le schéma~$\phi^{-1}(U_i)$
soit affine et~$\sch O_Y(\phi^{-1}(U_i))$ soit une~$\sch O_X(U_i)$-algèbre finie. 

\deux{ex-morph-fin}
{\bf Exemples et premières propriétés.}

\trois{imm-ferm-fini}
Si~$A$ est  un anneau et~$I$ un idéal de~$A$, la~$A$-algèbre~$A/I$ est finie ; il s'ensuit que
toute immersion fermée est un morphisme fini.

\trois{anneaux-entiers-cdn}
Soit~$K$ un corps de nombres et soit~$\got O_K$ l'anneau des entiers de~$K$ ; comme~$\got O_K$ est une~$\ZZ$-algèbre
finie (en tant que~$\ZZ$-module, il est libre de rang
égal à~$[K:\QQ]$), le morphisme~$\spec \got O_K\to \spec \ZZ$ est fini.

\trois{oxt-modp}
Soit~$X$ un schéma et soit~$P=T^n+\sum_{i\leq n-1}a_i T^i$ un polynôme
unitaire à coefficients dans~$\sch O_X(X)$. Le faisceau
$$\sch O_X[T]/P:=U\mapsto \sch O_X(U)[T]/\left(T^n+\sum_{i\leq n-1}a_i|_U T^i\right)$$
est une~$\sch O_X$-algèbre quasi-cohérente, qui comme~$\sch O_X$-module
satisfait visiblement les conditions équivalentes du lemme~\ref{def-qc-loctf}. En conséquence, 
$\spec \sch O_X[T]/P$ est un~$X$-schéma fini. 

\trois{comp-morph-fin}
La composée de deux morphismes finis est un morphisme fini : c'est une conséquence immédiate
de~\ref{reform-morph-fin}
et du fait que si~$A$ est un anneau, si~$B$ est une~$A$-algèbre finie et si~$C$ est une~$B$-algèbre finie
alors~$C$ est une~$A$-algèbre finie. 

\trois{fini-changebase}
Soit~$X$ un schéma, soit~$Y$ un~$X$-schéma fini et soit~$Z$ un~$X$-schéma. Le produit
fibré~$Y\times_X Z$ est alors un~$Z$-schéma fini. En effet, en raisonnant localement et en utilisant
une fois encore~\ref{reform-morph-fin}, cette assertion se ramène au fait connu suivant : si~$A$ est un anneau, si~$B$
est une~$A$-algèbre finie et si~$C$ est une~$A$-algèbre alors~$C\otimes_A B$ est une~$C$-algèbre finie. 

\deux{morph-fin-ferm}
{\bf Proposition.}
{\em Soit~$\phi \colon Y \to X$ un morphisme fini. Il est alors {\em fermé},
c'est-à-dire que $\phi(Z)$ est un fermé de~$X$ pour tout fermé~$Z$ de~$Y$.}

\medskip
{\em Démonstration.}
Soit~$Z$ un fermé de~$Y$. Munissons-le d'une structure quelconque de sous-schéma fermé. 
La composée $Z\hookrightarrow Y\to X$ est alors encore un morphisme fini, et notre proposition
revient à montrer que son image est fermée. Quitte à remplacer~$Y\to X$ par ce morphisme, 
on peut donc supposer que~$Z=Y$. 

\medskip
Être un fermé étant une propriété locale, on peut supposer que~$X$ est affine. On a donc~$X=\spec A$
pour un certain anneau~$A$, et~$Y$ est alors le spectre d'une~$A$-algèbre finie~$B$. Soit~$I$ le noyau de~$A\to B$
et soit~$T$ le spectre de~$A/I$. Le morphisme~$\phi \colon Y \to X$ admet une factorisation canonique
$Y\to T \hookrightarrow X$, correspondant à la factorisation~$A\to A/I\hookrightarrow B$ au niveau des anneaux. 

\medskip
Comme~$A/I\hookrightarrow B$ est une injection qui fait de~$B$ une~$A/I$-algèbre finie, et
{\em a fortiori}
entière, le lemme de {\em going-up}
(lemme~\ref{going-up}, {\em cf.}
notamment
la remarque qui précède sa démonstration)
assure que~$\spec B\to \spec A/I$ est surjective, c'est-à-dire que~$Y\to T$ est surjective. Puisque~$T\hookrightarrow X$
a pour image le fermé~$V(I)$, il vient~$\phi(Y)=V(I)$.~$\Box$ 

\section{Morphismes de type fini}
\markboth{La notion de schéma}{Morphismes de type fini}

\subsection*{Définition, exemples, premières propriétés}

\deux{prop-sch-tf}
{\bf Proposition.}
{\em Soit~$A$ un anneau et soit~$X$ un~$A$-schéma possédant la propriété suivante : il existe un recouvrement 
ouvert~$(X_i)$ de~$X$ tel que pour tout~$i$ le~$A$-schéma~$X_i$ soit égal à~$\spec A_i$ pour une certaine~$A$-algèbre
{\em de type} fini~$A_i$. Sous ces hypothèses, pour tout ouvert affine~$U$ de~$X$ la~$A$-algèbre~$\sch O_X(U)$ est de type fini.}

\medskip
{\em Démonstration.}
Soit~$U$ un ouvert affine de~$X$ et soit~$B$ la~$A$-algèbre~$\sch O_X(U)$. 

\trois{recouvre-u-schtf}
Soit~$x\in U$. Il appartient à~$X_i$ pour un certain~$i$. Il existe alors~$a\in A_i$ tel que
l'ouvert~$V:=D(a)$ de~$X_i$ soit contenu dans~$U\cap X_i$ et contienne~$x$. 
On a~$\sch O_X(V)=(A_i)_a=A_i[T]/(aT-1)$ ; en conséquence, $\sch O_X(V)$ est une~$A$-algèbre de type fini. 

\medskip
Il existe~$f\in B$ tel que l'ouvert~$D(f)$ de~$U$ soit contenu dans~$V$ et contienne~$x$. Cet ouvert
est {\em a fortiori}
égal à l'ouvert~$D(f)$ de~$V$, et son algèbre des fonctions est donc égale à~$\sch O_X(V)_f=\sch O_X(V)[T]/(fT-1)$ ; 
c'est donc encore une~$A$-algèbre de type fini. 

\medskip
Par quasi-compacité de~$U$ on en déduit qu'il existe une famille finie $(f_1,\ldots, f_n)$ d'éléments de~$B$ tels que
$U=\bigcup D(f_i)$ et tel que~$\sch O_X(D(f_i))=B_{f_i}$ soit pour tout~$i$ une~$A$-algèbre de type fini. 

\trois{conclu-tpfinilocal}
Puisque~$U=\bigcup D(f_i)$ il existe une famille~$(b_i)$ d'éléments
de~$B$ tels que~$\sum b_i f_i=1$. Par ailleurs, il existe par hypothèse pour tout~$i$
une famille finie d'éléments de~$B_{f_i}$ engendrant celle-ci comme~$A$-algèbre ; 
on les écrit~$\beta_{i1}/f_i^{n_{i1}},\ldots, \beta_{ir_i}/f_i^{n_{ir_i}}$ où les~$\beta_{ij}$ appartiennent
à~$B$. Soit~$C$ la sous-$A$-algèbre de~$B$ engendrée par les~$b_i$, les~$f_i$ et les~$\beta_{ij}$ ; nous allons
montrer que~$B=C$, ce qui achèvera la démonstration. 

\medskip
Soit~$b\in B$. Fixons~$i$. Par choix des~$\beta_{ij}$, on peut écrire l'élément~$b/1$ de~$B_{f_i}$ 
comme un polynôme à coefficients dans~$A$ en les~$\beta_{ij}/f_i^{n_j}$ ; la condition d'égalité entre fractions
entraîne alors qu'il existe~$N$ tel que~$f_i^Nb$ soit un polynôme à coefficients dans~$A$
en~$f_i$ et les~$\beta_{ij}$ ; en particulier, $f_i^Nb$ appartient à~$C$, et cela reste vrai si l'on augmente l'exposant~$N$. 

\medskip
Il existe en conséquence~$N$ tel que~$f_i^mb\in C$ pour tout~$i$ dès que~$m\geq N$. 
On a~$b =(\sum_{1\leq i\leq n} b_i f_i)^{nN}b$. Lorsqu'on développe cette expression, on trouve une somme de termes 
de la forme~$b_1^{e_1}\ldots b_n^{e_n}f_1^{e_1}\ldots f_n^{e_n}b$ avec~$\sum e_i =nN$. Dans un tel terme, il existe nécessairement~$i_0$ tel que
$e_{i_0}\geq N$. On a alors
$$b_1^{e_1}\ldots b_n^{e_n}f_1^{e_1}\ldots f_n^{e_n}b=\left(\prod_{i\neq i_0}b_i^{e_i}f_i^{e_i}\right)\cdot b_{i_0}^{e_{i_0}}\cdot (f_{i_0}^{e_{i_0}}b)\in C,$$
et~$b$ appartient donc à~$C$.~$\Box$ 

\deux{rem-loc-tf}
{\bf Remarque}. La proposition ci-dessus affirme en particulier que si~$A$ est un anneau et si~$U$
est un ouvert affine de~$\spec A$ alors~$\sch O_{\spec A}(U)$ est une~$A$-algèbre de type fini. Cela n'avait rien d'évident
{\em a priori}, sauf quand~$U$ est de la forme~$D(f)$ car alors~$\sch O_{\spec A}(U)=A_f=A[T]/(fT-1)$. 

\deux{pro-mor-tf}
{\bf Proposition.}
{\em Soit~$\phi \colon Y\to X$ un morphisme de schémas. 

\begin{itemize}
\item[A)] Les assertions suivantes sont équivalentes :

\medskip
\begin{itemize} 
\item[i)] pour tout ouvert affine~$U$ de~$X$, le schéma~$\phi^{-1}(U)$ admet un recouvrement ouvert
fini 
par des spectres de~$\sch O_X(U)$-algèbres de type fini ;

\item[ii)] il existe un recouvrement~$(U_i)$ de~$X$ par des ouverts affines tels que
$\phi^{-1}(U_i)$ admette pour tout~$i$
un recouvrement ouvert
fini 
par des spectres de~$\sch O_X(U_i)$-algèbres de type fini.

\end{itemize}
\medskip
\item[B)] Si les assertions ci-dessus sont satisfaites alors pour tout ouvert affine~$U$
de~$X$ et tout ouvert affine~$V$ de~$\phi^{-1}(U)$, la~$\sch O_X(U)$-algèbre
$\sch O_Y(V)$ est de type fini. 
\end{itemize}}

\medskip
{\em Démonstration.}
Commençons par montrer~A).
Il est clair que~i)$\Rightarrow$ii). Supposons que~ii)
est vraie, et montrons~i). Soit~$U$ un ouvert affine de~$X$. Par quasi-compacité de~$U$, 
il existe~$(f_1,\ldots, f_r)\in \sch O_X(U)$ tels que les~$D(f_j)$ recouvrent~$U$ et tels
que pour tout~$j$ il existe~$i(j)$ vérifiant $D(f_j)\subset U_{i(j)}\cap U$. Il suffit pour conclure
de démontrer que pour tout~$j$ le schéma $\phi^{-1}(D(f_j))$ admet 
un recouvrement ouvert
fini 
par des spectres de~$\sch O_X(U)$-algèbres de type fini.

\medskip
Fixons~$j$, et écrivons~$f$ au lieu de~$f_j$ et~$i$ au lieu de~$i(j)$. Le schéma~$\phi^{-1}(U_i)$ 
admet un recouvrement ouvert fini~$(V_\alpha)$ où chaque~$V_\alpha$ est le spectre
d'une~$\sch O_X(U_i)$-algèbre de type fini~$A_\alpha$. 

\medskip
Le schéma~$\phi^{-1}(D(f))=\phi^{-1}(U_i)\times_{U_i} D(f)$ est réunion de
ses ouverts
$$V_\alpha\times_{U_i} D(f)=\spec A_\alpha\otimes_{\sch O_X(U_i)}\sch O_X(D(f)).$$
Pour tout~$\alpha$, la~$\sch O_X(D(f))$-algèbre~$A_\alpha\otimes_{\sch O_X(U_i)}\sch O_X(D(f))$
est de type fini ; comme~$\sch O_X(D(f))=\sch O_X(U)[T]/(fT-1)$, elle
est également de type fini
sur~$\sch O_X(U)$, ce qui achève la démonstration
de~A).

\medskip
L'assertion~B) découle
quant à elle
immédiatement
de la proposition~\ref{prop-sch-tf},
appliquée au schéma~$\phi^{-1}(U)$.~$\Box$ 

\deux{def-mor-tpf}
{\bf Définition.}
On dit qu'un morphisme de schémas~$\phi \colon Y\to X$ est {\em de type fini}
s'il satisfait les conditions
équivalentes~i) et~ii)
de la proposition~\ref{pro-mor-tf}
ci-dessus. On dit parfois aussi que $Y$ est {\em de type fini sur~$X$},
ou bien est un {\em $X$-schéma de type fini}.  

\deux{pro-mor-tpf}
{\bf Exemples et premières propriétés}. 

\trois{fini-tpf}
Il résulte immédiatement des définitions qu'un morphisme fini est de type fini ; c'est en particulier le cas des immersions fermées. 

\trois{aff-tf}
Si~$A$ est un anneau et~$B$ une~$A$-algèbre alors~$\spec B\to \spec A$ est de type fini si et seulement si~$B$ est de type fini
comme~$A$-algèbre : la condition est en effet suffisante par définition, et nécessaire en vertu de l'assertion~B)
de la proposition~\ref{pro-mor-tf}. 

\trois{anfini}
Pour tout entier~$n$ et tout schéma~$X$, le schéma~$\Aff^n_X$ est de type fini sur~$X$. En effet, la propriété est par définition locale sur~$X$, ce qui permet
de se ramener au cas où celui-ci est affine, auquel cas c'est immédiat car~$\Aff^n_A$ est égal à~$\spec A[T_1,\ldots, T_n]$ pour tout anneau~$A$. 

\trois{compose-tf}
La composée de deux morphismes de type fini est de type fini : on le déduit immédiatement
de leur caractérisation {\em via}
la propriété~i) de la proposition~\ref{prop-sch-tf}. 

\trois{tpf-changebase}
Soit~$X$ un schéma, soit~$Y$ un~$X$-schéma de type fini et soit~$Z$ un~$X$-schéma.
Le~$Z$-schéma $Y\times_X Z$ est alors de type fini. En effet, on peut raisonner localement sur~$Z$, 
et {\em a fortiori}
sur~$X$ ; cela autorise à supposer~$X$ et~$Z$
affines. Dans ce cas, $Y$ possède un recouvrement ouvert
affine fini~$(V_i)$ tel que~$\sch O_Y(V_i)$ soit une~$\sch O_X(X)$-algèbre de type fini pour tout~$i$. Le schéma~$Y\times_XZ$
est alors réunion de ses ouverts affines~$V_i\times_X Z$ ; pour tout~$i$, la~$\sch O_Z(Z)$-algèbre
$\sch O_{Y\times_X Z}(V_i\times_X Z)$ est égale à
$\sch O_Y(V_i)\otimes_{\sch O_X(X)}\sch O_Z(Z)$ et est donc
de type fini, d'où notre assertion. 

\subsection*{Schémas de type fini sur un corps}

\deux{intro-tf-k}
Soit~$k$ un corps, et soit~$X$ un~$k$-schéma de type fini. 

\trois{x-tf-k}
Par définition, $X$ est recouvert par un nombre fini d'ouverts affines de la forme~$\spec A$ où~$A$
est une~$k$-algèbre de type fini. Une telle~$A$ étant noethérienne, son spectre est noethérien 
(\ref{def-noeth}
{\em et sq.}). Il s'ensuit aisément que
l'espace topologique
$X$ est lui-même noethérien. 

\trois{ouv-aff-xtf}
On déduit de l'assertion~B) de la proposition~\ref{pro-mor-tf}, ou directement de la proposition~\ref{prop-sch-tf},
que si~$U$ est un ouvert affine de~$X$ alors~$\sch O_X(U)$ est une~$k$-algèbre de type fini.
Insistons sur l'importance 
que~$U$ soit {\em affine} : il existe des contre-exemples lorsqu'il ne l'est pas, même sur~$\CC$. 

\trois{ouv-gen-xtf}
Comme~$X$ est noethérien, ses ouverts sont tous quasi-compacts et donc encore de type fini sur~$k$ d'après~\ref{ouv-aff-xtf}. 

\trois{points-fermes-xtf}
Soit~$x\in X$. Son corps résiduel~$\kappa(x)$ est une extension finie de~$k$. Le point~$x$ est fermé si et seulement si il l'est dans
tout ouvert affine le contenant ; on déduit alors de~\ref{pts-fermes-a}
que~$x$ est fermé si et seulement si~$\kappa(x)$ est une extension finie de~$k$.

\trois{points-fermes-attention}
Notons une conséquence importante de ce qui précède : si~$x\in X$ et si~$U$ est un ouvert de~$X$ contenant~$x$
alors~$x$ est fermé dans~$X$ si et seulement si il l'est dans~$U$ : les deux propriétés équivalent en effet à la finitude de~$\kappa(x)$
sur~$k$, puisque~$U$ est lui aussi un~$k$-schéma de type fini. 

\medskip
Cette équivalence vous paraît peut-être anodine, mais il n'en est rien. Par exemple, soit~$S$ le spectre de~$k[[t]]$. L'anneau~$k[[t]]$ a deux
idéaux premiers : l'idéal nul et l'idéal maximal~$(t)$. En conséquence~$S$ comprend deux points : le point générique~$\eta$ et un unique
point fermé~$s$.
Le point~$\eta$ est ouvert et dense. Il est évidemment fermé dans l'ouvert~$\{\eta\}$, mais n'est pas fermé
dans~$S$. 

\trois{exist-pointferme}
Nous proposons au lecteur de montrer en exercice que tout schéma quasi-compact non vide possède un point fermé. 
Mais ici, 
on peut voir directement que si~$X\neq \varnothing$ il possède un point fermé  : en effet, il existe dans ce cas un ouvert affine
non vide~$U$ de~$X$, lequel possède un point fermé~$x$, qui est également fermé dans~$X$ d'après la remarque
du~\ref{points-fermes-attention}
ci-dessus. 

\deux{point-x-extk}
Soit~$L$
une extension de~$k$. L'ensemble~$\hom_k(\spec L,X)$ s'identifie
à l'ensemble des couples~$(x,\iota \colon \kappa(x)\hookrightarrow L)$ où~$x\in X$
et où~$\iota$ est un~$k$-plongement ({\em cf.}~\ref{mor-k-schem} {\em et sq.}). Cet ensemble
se note aussi~$X(L)$ ; ses éléments sont aussi appelés les~{\em $L$-points}
de~$X$. 

\trois{point-x-krat}
En particulier, $X(k)$ est l'ensemble des couples
$(x, \iota \colon \kappa(x)\hookrightarrow k)$. Mais si~$x\in X$
l'ensemble des~$k$-plongements de~$\kappa(x)$ dans~$k$ est facile
à décrire : c'est~$\{{\rm Id}_k\}$ si~$\kappa(x)=k$ et~$\varnothing$
sinon. En conséquence, l'ensemble des~$k$-points
de~$X$ s'identifie à l'ensemble
des points schématiques de~$X$ de corps résiduel~$k$. 

\trois{extension-k-attention}
On prendra garde que 
si~$L$ est une extension stricte de~$k$, 
le lien entre~$L$-point et point schématique est 
en général plus subtil : 
un point fermé donné $x$ peut
supporter différents~$L$-points de~$X$. Ainsi, 
supposons que~$k=\RR$, que~$X=\Aff^1_{\RR}$ et que~$x$ est 
le point~$V(T^2+1)$ (où~$T$ est la fonction coordonnée).
Le corps~$\kappa(x)$
est alors égal à~$\RR[T]/T^2+1$, et il existe
deux~$\RR$-morphismes de~$\kappa(x)$ dans~$\CC$, 
à savoir~$T\mapsto i$ et~$T\mapsto (-i)$. Le point~$x$ est
donc le support de deux~$\CC$-points distincts de~$X$. 

\trois{point-xaff-extk}
Si le schéma~$X$ est affine, on peut l'écrire sous la forme
$\spec k[T_1,\ldots, T_n]/(P_1,\ldots, P_r)$ pour une certaine
famille~$(P_j)$ de polynômes en~$n$ variables, 
et l'on a alors
$$X(L)=\hom_k(\spec L,X)
=\hom_k(k[T_1,\ldots, T_n]/(P_1,\ldots, P_r), L)$$
$$\simeq \{(x_1,\ldots, x_n)\in L^n, \;P_j(x_1,\ldots, x_n)=0\;\forall j\}$$ 
(ainsi, notre notation~$X(L)$ est 
compatible dans ce cas avec celle introduite en~\ref{def-xl}). On voit que~$X(L)$
est l'ensemble des~$L$-points de la variété algébrique naïve
qui correspond à~$X$, ce qui explique le choix de l'expression~«$L$-point»
en dépit des ambiguïtés dues à l'autre sens du mot «point»
({\em cf.}~\ref{extension-k-attention}). 

\trois{xtf-kbarre-galois}
On ne suppose plus~$X$ affine. Soit~$\bar k$ une clôture algébrique de~$k$ et soit
$G$ le groupe~${\rm Gal}(\bar k/k)$. Le groupe~$G$ agit sur~$\bar k$, donc sur~$\spec \bar k$, 
et donc par composition sur~$\hom_k(\spec \bar k, X)=X(\bar k)$. Si~$U$ est un ouvert affine de~$X$, 
l'ensemble~$U(\bar k)$ est simplement le sous-ensemble de~$X(\bar k)$ formé des~$\bar k$-points dont le point schématique
sous-jacent est situé sur~$U$ ; il est dès lors stable sous l'action de~$G$. 

\medskip
On déduit alors de~\ref{action-g}
que l'ensemble~$X_0$ des points fermés de~$X$ s'identifie naturellement à~$X(\bar k)/G$ ; on retrouve ainsi l'injection
de~$X(k)$ dans~$X_0$ mentionnée au~\ref{point-x-krat}. 

\deux{constr-xtf}
Supposons pour ce paragraphe que~$k$
est algébriquement clos ; dans ce cas~$X(k)$ s'identifie au sous-ensemble~$X_0$ de~$X$,
et on le munit de la topologie induite. Rappelons (\ref{top-const})
que si~$T$ est un espace topologique, on note~$\sch C(T)$ l'ensemble
des parties de~$T$ de la forme~$\bigcup_{1\leq i\leq n}U_i\cap F_i$ où les~$U_i$ sont ouverts et les~$F_i$ fermés. Les faits
suivants se déduisent des résultats correspondants déjà établis dans le cas affine (\ref{constr-kalgclos}
{\em et sq.}, \ref{points-schem-interp}
{\em et sq.}). 

\trois{const-bij-xtf}
La flèche~$C\mapsto C(k):=C\cap X(k)$ établit une bijection entre~$\sch C(X) $ et~$\sch C(X(k))$ ; si~$C$ et~$D$ sont deux
éléments de~$\sch C(X)$ alors~$C\subset D$ si et seulement si~$C(k)\subset D(k)$ ; si~$C\in \sch X$ alors~$C$ est un fermé irréductible
si et seulement si~$C(k)$ est un fermé irréductible de~$X(k)$. 

\trois{fermes-irred-xtf}
L'application~$x\mapsto \overline{\{x\}}(k)$ établit une bijection entre~$X$ et l'ensemble des fermés irréductibles de~$X(k)$ ; en conséquence, 
$X$ s'obtient en rajoutant à~$X(k)$ (dont tous les points sont fermés)
un point générique par fermé irréductible non singleton. En d'autres termes, $X$ est la sobrification de~$X(k)$. 

\deux{dimkrull-xtf}
{\bf Dimension de Krull du schéma~$X$.}
On ne suppose plus que~$k$ est algébriquement clos. Nous allons
montrer que si~$X\neq \varnothing$ sa dimension de Krull
est finie.

\trois{intro-krull-xtf}
Comme~$X$ est noethérien,
il possède une décomposition~$X=\bigcup_{1\leq i\leq n}X_i$
en composantes irréductibles (lemme~\ref{lemme-dec-irr}). Nous laissons
le lecteur vérifier que l'on a alors 
$$\dim_{\rm Krull}X=\sup_i \dim_{\rm Krull}X_i.$$ Il suffit donc de démontrer que la dimension de Krull
de chacune des~$X_i$ est finie ; on s'est ainsi ramené au cas où~$X$ est irréductible. 

\trois{krull-cas-irred}
Comme~$X_{\rm red}\hookrightarrow X$ est une immersion fermée elle est de type fini, et~$X_{\rm red}$
est donc un~$k$-schéma de type fini ; de plus, $X_{\rm red}\hookrightarrow X$ induit un homéomorphisme
entre les espaces topologiques sous-jacents. On peut donc remplacer~$X$ par~$X_{\rm red}$ et le supposer
réduit. 

\medskip
Le schéma irréductible~$X$ possède un unique point générique~$\xi$. Soit~$U$ un ouvert affine non vide de~$X$. 
Il est irréductible ; comme~$X$ est réduit, $\sch O_X(U)$ est réduit, et donc intègre. C'est une~$k$-algèbre de type fini, 
et la dimension de Krull de~$U$ est égale au degré de transcendance de~${\rm Frac}\;\sch O_X(U)$ sur~$k$
(th.~\ref{prop-dimkrull-transc}).

\medskip
Le point~$\xi$ est l'unique
point générique de~$U$, et l'on a
$${\rm Frac}(\sch O_X(U))=\sch O_{U,\xi}=\sch O_{X,\xi}.$$
Le corps
${\rm Frac}(\sch O_X(U))$ s'identifie donc à~$\sch O_{X,\xi}$ et ne dépend
en particulier pas de~$U$. On l'appelle le {\em corps des fonctions}
de~$X$. Soit~$d$ son degré de transcendance sur~$k$ ; nous allons montrer que la dimension 
de Krull de~$X$ est égale à~$d$. Par ce qui précède, nous avons déjà que~$d$ est la 
dimension de Krull de tout ouvert {\em affine}
non vide de~$X$. 

\medskip
{\em La dimension de Krull de~$X$ est majorée par~$d$}.
En effet, soit
$$F_0\subsetneq F_1\subsetneq\ldots\subsetneq F_n$$ une suite de fermés irréductibles de~$X$. Comme~$F_0$ est
irréductible, il est non vide et rencontre donc un ouvert affine~$U$ de~$X$ ; celui-ci rencontre
{\em a fortiori}
chacun des~$F_i$. 

Pour tout~$i$, l'intersection~$F_i\cap U$ est un ouvert non vide de l'espace irréductible~$F_i$, il 
est donc irréductible et dense dans~$F_i$ ; il vient~$F_i=\overline{U\cap F_i}$. On en déduit
que les ensembles~$U\cap F_i$ sont deux à deux distincts, puisque les~$F_i$ le sont. 
La suite 
$$U\cap F_0\subsetneq U\cap F_1\subsetneq\ldots\subsetneq U\cap F_n$$
est ainsi une chaîne strictement croissante de fermés irréductibles de~$U$ ; comme
celui-ci est de dimension de Krull égale à~$d$, on a~$n\leq d$, et~$\dim_{\rm Krull}X \leq d$. 

\medskip
{\em La dimension de Krull de~$X$ est minorée par~$d$.}
Soit~$U$ un ouvert affine non vide de~$X$. Il est de dimension de Krull~$d$ ; en conséquence, 
il existe une chaîne strictement croissante
$$G_0\subsetneq G_1\subsetneq\ldots\subsetneq G_d$$ 
de fermés irréductibles de~$U$. Pour tout~$i$, le fermé~$\overline{G_i}$ de~$X$
est irréductible, et son intersection avec~$U$ est égale à~$G_i$. On en déduit
que les fermés~$\overline{G_i}$ sont deux à deux distincts, puisque les~$G_i$ le sont. 
La suite 
$$\overline{G_0}\subsetneq \overline{G_1}\subsetneq\ldots\subsetneq \overline{G_d}$$
est ainsi une chaîne strictement croissante de fermés irréductibles de~$X$ ; 
on a donc~$\dim_{\rm Krull}X \geq d$, ce qui termine la preuve. 

\trois{comment-dimkrull-tf}
{\em Remarque.}
Si~$U$ est un ouvert non vide quelconque de~$X$ on a
l'égalité~$\sch O_{U,\xi}=\sch O_{X,\xi}$ ; par
ce qui précède, il s'ensuit que~$\dim_{\rm Krull}U=d$ : la dimension de Krull de tout 
ouvert non vide de~$X$
coïncide avec celle de~$X$. 

\medskip
Une fois encore, cette remarque n'est pas anodine. Pour le voir,
considérons le spectre~$S$ de~$k[[t]]$ que nous avons
décrit au~\ref{points-fermes-attention} ;
soit~$s$ son point fermé et soit~$\eta$ son point générique. Les seuls fermés
irréductibles de~$S$ sont~$\{s\}$ et~$S$, et l'on a évidemment~$\{s\}\subsetneq S$ ; 
en conséquence, la dimension de Krull de~$S$ est~$1$, mais celle 
de son ouvert dense~$\{\eta\}$ est égale à~$0$. 

\section{Le foncteur des points d'un schémas, ou la revanche 
du point de vue ensembliste}\label{FONCTPOINT}
\markboth{La notion de schéma}{Le foncteur des points
d'un schéma}

\deux{intro-yoneda-schemas}
Nous avons jusqu'à maintenant considéré les schémas comme des espaces localement annelés. Cette approche offre
l'avantage de décalquer, dans une certaine mesure, l'intuition géométrique classique : elle permet de parler d'ouverts
et de fermés, d'évaluer les fonctions... Toutefois, elle présente en regard des inconvénients assez lourds : on doit accepter
que le corps résiduel varie avec le point considéré, qu'une fonction puisse s'annuler ponctuellement partout sans être pour autant 
nulle, que l'espace sous-jacent au produit fibré ne soit pas le produit fibré des espaces sous-jacents, etc. 

\medskip
Mais on peut penser à un schéma autrement : comme n'importe quel objet de n'importe quelle 
catégorie, il est en vertu du lemme de Yoneda entièrement déterminé par le foncteur qu'il représente. 
Et si tautologique que soit ce constat, nous allons voir sur plusieurs exemples qu'il peut parfois présenter
un côté rafraîchissant : il permet en effet
dans un certain nombre de cas
de revenir en un sens à la définition première 
de la géométrie algébrique, à savoir
l'étude des
{\em ensembles}
de solutions d'équations polynomiales, et des
{\em applications}
polynomiales entre iceux. 

\subsection*{Premiers exemples}

\deux{suite-intro-yoneda-schema}
Soit~$S$ un schéma ; dans ce qui suit, nous allons travailler dans la catégorie
$S\text{-}\mathsf{Sch}$
des~$S$-schémas 
(qui est celle des schémas tout courts lorsque~$S$ est égal à~$\spec \ZZ$). Soit~$X$ un~$S$-schéma.  Pour
tout~$S$-schéma~$T$, nous noterons~$X(T)$ l'ensemble~$\hom_S(T,X)$, et nous écrirons souvent
par abus~$X(A)$ plutôt que~$X(\spec A)$ ; notons que ces conventions
sont compatibles avec la notation~$X(L)$ introduite plus haut (\ref{point-x-extk}) En d'autres termes, $T\mapsto X(T)$ est le foncteur
contravariant de~$S\text{-}\mathsf{Sch}$
dans~$\ens$
représenté par~$X$. 
On dit parfois que~$X(T)$ est l'ensemble des~{\em $T$-points}
(ou des~{\em $A$-points}
si~$X=\spec A$)
du~$T$-schéma~$X$.

\trois{x-egal-sonfoncteur}
Le lemme de Yoneda assure que~$X$ est entièrement déterminé par~$X\mapsto X(T)$, et que
se donner un
morphisme~$Y\to X$ dans $S\text{-}\mathsf{Sch}$
revient à se donner un morphisme
entre les foncteurs~$T\mapsto Y(T)$ et~$T\mapsto X(T)$. Il n'y a donc
aucun inconvénient à identifier, si on le juge utile, un
schéma~$X$ au foncteur~$T\mapsto X(T)$, 
ce qui justifie {\em a posteriori}
la notation~$X(T)$. 

\trois{test-seulement-affine}
Soit~$X$ un~$S$-schéma. Si~$T$ est un~$S$-schéma il est recouvert par des ouverts affines, et~$U\mapsto \hom_S(U,X)$
est un faisceau sur~$X$. Il s'ensuit que le {\em foncteur}~$X$ est entièrement déterminé par sa restriction
à 
la catégorie~$S\text{-}\mathsf{Aff}$
des~$S$-schémas qui sont affines (dans l'absolu, pas relativement à~$S$). 

Pour la même raison, si~$Y$ est un~$X$-schéma,
tout morphisme
entre les foncteurs~$Y|_{S\text{-}\mathsf{Aff}}$ et~$X|_{S\text{-}\mathsf{Aff}}$
s'étend d'une unique manière en un morphisme de foncteurs de~$Y$ vers~$X$. 

\medskip
Il n'y a donc aucun inconvénient, si l'on préfère travailler avec des schémas affines, à 
se contenter de voir 
un~$S$-schéma~$X$ comme un foncteur contravariant de~$S\text{-}\mathsf{Aff}$
dans~$\ens$. 

\trois{prod-fib-xt}
Soient~$X$ et~$S'$ deux~$S$-schémas, et soit~$T$ un~$S'$-schéma, que l'on peut voir comme
un~$S$-schéma par composition avec~$S'\to S$. Nous vous laissons vérifier qu'il existe une bijection naturelle
$$X(T)\simeq (X\times_S S')(T)$$
(à gauche, $X$ est vu comme foncteur de~$S\text{-}\mathsf{Sch}$
vers~$\ens$ ; à droite, $X\times_S S'$ est vu comme foncteur de~$S'\text{-}\mathsf{Sch}$
vers~$\ens$). Remarquez qu'il s'agit simplement d'une déclinaison 
de l'énoncé~\ref{rest-cs}, dont nous vous avions
déjà proposé la preuve en exercice.

\deux{exemple-foncteur-schema}
{\bf Un exemple.}
Soit~$A$ un anneau et soit~$(P_j)_{j\in J}$ une famille
de polynômes appartenant à~$A[U_1,\ldots, U_n]$. Posons
$$X=\spec A[U_1,\ldots, U_n]/(P_j)_j.$$
Soit~$T$ un~$A$-schéma. 
L'ensemble~$X(T)$ est égal à~$\hom_A(T,X)$, 
c'est-à-dire à~$\hom_A(A[U_1,\ldots, U_n]/(P_j)_j, \sch O_T(T))$. 
Or cet ensemble est lui-même en bijection naturelle, {\em via}
la flèche~$\phi \mapsto (\phi(\overline{U_1},\ldots, \phi(\overline {U_n}))$,
avec l'ensemble des~$n$-uplets~$(t_1,\ldots, t_n)$ de~$\sch O_T(T)^n$
tels que~$P_j(t_1,\ldots, t_n)=0$ pour tout~$j$. On dispose donc d'une bijection
fonctorielle en~$T$

$$X(T)\simeq \{(t_1,\ldots, t_n)\in \sch O_T(T)^n\;\;\text{t.q.}\;P_j(t_1,\ldots,t_n)=0\;\;\forall j\},$$
qui envoie un morphisme~$\psi$ sur~$(\psi^*\overline{U_1},\ldots, \psi^*\overline{U_n})$. 
Le {\em foncteur}~$X$ est donc tout simplement le foncteur
«ensemble des~$n$-uplets solutions du systèmes d'équations polynomiales~$(P_j)_j$» : c'est 
le retour annoncé du point vue naïf ou ensembliste sur la géométrie algébrique. 

\trois{cas-part-ox-mora1}
{\em Un cas particulier.}
Supposons que~$n=1$ et que la famille des~$P_j$ est vide (une variable,
pas d'équations). On a alors~$X=\Aff^1_A$, et l'on dispose par ce qui précède
pour tout~$A$-schéma~$T$ d'une bijection~$\hom_A(T,\Aff^1_A)\simeq \sch O_T(T)$,
fonctorielle en~$T$, qui envoie un morphisme~$\psi$ sur~$\psi^*\overline U$ (où~$U$ est ici
la fonction coordonnée sur~$\Aff^1_A$). On a donc une bonne raison supplémentaire de penser
à~$\sch O_T(T)$ comme à l'anneau des fonctions sur~$T$ : les éléments de~$\sch O_T(T)$
«sont»
exactement les~$A$-morphismes de~$T$ vers la droite affine. 

\trois{foncteur-ouvert-df}
Revenons au cadre général décrit au~\ref{exemple-foncteur-schema}, 
et soit~$f$
appartenant à~$A[U_1,\ldots, U_n]$ ; nous allons décrire fonctoriellement l'ouvert~$X':=D(\bar f)$
de~$X$. Comme~$(A[U_1,\ldots, U_n]/(P_j)_j)_{\bar f}=A[U_1,\ldots, U_n,V]/((P_j)_j, Vf-1)$, on déduit de
ce qui précède que l'on a pour tout~$A$-schéma~$T$ une bijection naturelle
entre~$X'(T)$ et
$$\{(t_1,\ldots, t_n,s)\in \sch O_T(T)^n\;\;\text{t.q.}\;P_j(t_1,\ldots,t_n)=0\;\;\forall j\;\text{et}\;f(t_1,\ldots, t_n)s=1\}.$$
La condition~$f(t_1,\ldots, t_n)s=1$ peut se récrire
«$f(t_1,\ldots, t_n)$ est inversible et~$s$ est son inverse», d'où une bijection
naturelle entre~$X'(T)$ et
$$ \{(t_1,\ldots, t_n)\in \sch O_T(T)^n\;\;\text{t.q.}\;P_j(t_1,\ldots,t_n)=0\;\;\forall j\;\text{et}\;f(t_1,\ldots, t_n)\in \sch O_T(T)^\times\}.$$

Attention donc : du point de vue fonctoriel, la condition~«$f\neq 0$»
se traduit par~«$f$ inversible». 
Observez d'ailleurs à ce propos
que la condition d'être non nulle n'est de toutes
façons pas fonctorielle : un morphisme d'anneaux peut
avoir un noyau non trivial !

\trois{mor-schem-fonct}
Donnons-nous maintenant une famille~$(Q_\ell)_{\ell \in \Lambda}$ de polynômes
appartenant à~$A[V_1,\ldots, V_m]$, et posons
$$Y=\spec A[V_1,\ldots, V_m]/(Q_\ell)_\ell.$$ 
Le {\em foncteur}
$Y$ est égal en vertu de ce qui
précède à
$$T\mapsto \{(s_1,\ldots, s_m)\in \sch O_T(T)^m\;\;\text{t.q.}\;Q_\ell(s_1,\ldots,s_m)=0\;\;\forall \ell\}.$$ 

\medskip
Soit~$\psi \colon Y \to X$ un morphisme de~$A$-schémas ; nous allons décrire
le morphisme correspondants entre les foncteurs~$Y$ et~$X$. Le morphisme~$\psi$
est un élément de~$X(Y)$, et correspond 
dès lors par ce qui précède à
un~$n$-uplet~$(g_1,\ldots, g_n)$ d'éléments de~ $A[V_1,\ldots, V_m]/(Q_\ell)_\ell$
tels que~$P_j(g_1,\ldots, g_n)=0$. Pour tout~$i$, choisissons un polynôme~$G_i\in A[V_1,\ldots, V_m]$
relevant~$g_i$. 

\medskip
En reprenant l'ensemble des constructions, on voit 
aisément que le morphisme
induit par~$\psi$ entre les foncteurs~$Y$ et~$X$
est donné par la formule 
$$(s_1,\ldots, s_m)\mapsto (G_1(s_1,\ldots, s_m), \ldots, G_n(s_1,\ldots, s_m)).$$
Notez que tout est consistant : cette application ne dépend
bien
que des
$g_i$ et pas du choix des~$G_i$
(parce que le~$m$-uplet~$(s_1,\ldots, s_m)$ appartient à~$Y(T)$), et le~$n$-uplet
de droite appartient bien à~$X(T)$ (parce que chaque~$P_j$ s'annule en~$(g_1,\ldots, g_n)$). 

\medskip
Ainsi, lorsqu'on considère un morphisme entre les~$A$-schémas~$Y$ et~$X$ comme un morphisme
entre les foncteurs correspondants, on obtient une {\em application polynomiale} : là encore, 
on retombe sur la géométrie algébrique naïve ou ensembliste. 

\subsection*{Traduction schématique d'énoncés naïfs}

\deux{intro-naif-trad}
Bien que le but de la théorie des schémas soit
de simplifier la vie
des géomètres algébristes, on peut avoir à première vue l'impression
qu'elle la complique singulièrement. Pour se convaincre qu'il n'en est rien, 
il est important de bien comprendre que les énoncés «naïfs»
se transposent aisément, et le plus souvent 
de façon automatique, dans ce nouveau contexte. 
Nous allons illustrer ce propos par un exemple. 

\deux{intro-parametrage-cercle}
Soit~$k$ un corps algébriquement clos de caractéristique différente de~$2$.
Nous allons partir d'un énoncé de géométrie algébrique classique sur le corps~$k$, qui traduit
le fait qu'on peut paramétrer le «cercle»
sur~$k$ en faisant tourner une droite non verticale de pente~$t$ autour de~$(-1,0)$ et en considérant 
son deuxième point d'intersection avec le cercle : 
{\em les applications
polynomiales
$$t\mapsto \left(\frac{1-t^2}{1+t^2}, \frac{2t}{1+t^2}\right)\;\;\text{et}\;\;(x,y)\mapsto \frac {y}{x+1}$$
établissent un isomorphisme algébrique entre
$$\{t\in k, 1+t^2\neq 0\}\;\;\text{et}\;\{(x,y)\in k^2, x^2+y^2=1\;\;\text{et}\;x+1\neq 0\}.$$}
Nous allons donner deux traductions
de ce fait dans le langage des schémas, la première en termes d'espaces annelés et la seconde en termes 
de foncteurs des points, restreints aux schémas affines 
pour simplifier (\ref{test-seulement-affine}). 

\trois{traduc-espann}
{\em La première traduction.}
Posons
$$A=(k[x,y]/(x^2+y^2-1))_{x+1}\;\;\text{et}\;B=\spec k[t]_{t^2+1}.$$
On dispose de deux morphismes de~$k$-algèbres 
$$A\to B, \;\;\;x\mapsto \frac {1-t^2}{1+t^2}, y\mapsto \frac{2t}{1+t^2}$$ 

et~$$B\to A, \;\;t\mapsto \frac y {x+1}$$ qui sont des bijections réciproques l'une de l'autre, 
et qui induisent donc deux isomorphismes réciproques l'un de l'autre entre~$\spec B$ et~$\spec A$.

\trois{traduc-foncteur}
{\em La seconde traduction.}
Les foncteurs
$$ Y\colon R\mapsto \{t\in R, 1+t^2\in R^\times\}$$
$$\text{et}\;\;X\colon R \mapsto \{(x,y)\in R^2, x^2+y^2=1, (x+1)\in R^\times\}$$
de~$k\text{-}\mathsf{Alg}$
dans~$\ens$ sont représentables par des~$k$-schémas (affines), 
et
les formules $$t\mapsto \left(\frac{1-t^2}{1+t^2}, \frac{2t}{1+t^2}\right)\;\;\text{et}\;\;(x,y)\mapsto \frac {y}{x+1}$$
définissent
deux isomorphismes de foncteurs réciproques
l'un de l'autre entre~$Y$ et~$X$.

\medskip
On remarque que cette traduction fonctorielle est plus élémentaire
que la précédente, dans la mesure où elle est un décalque 
presque direct des énoncés naïfs (pensez toutefois à remplacer
partout «$\neq 0$»
par «inversible»), sans la contorsion psychologique consistant
à passer par les morphismes d'algèbres en changeant le sens des flèches.

\subsection*{Schémas en groupes} 

\deux{intro-groupe}
{\bf La notion d'objet en groupes dans une catégorie}. 
Un groupe est un ensemble~$G$ muni d'une application de~$G\times G$ vers~$G$ qui est 
associative, possède un élément neutre (nécessairement unique), et pour laquelle
tout élément admet un symétrique. 

\trois{reformu-groupe}
{\em Un petit jeu un peu loufoque}. Nous allons
traduire ce qui précède en termes purement catégoriques, 
{\em en nous interdisant de faire référence aux éléments}. 
Un groupe est donc un ensemble~$G$ muni des données supplémentaires suivantes. 
\medskip

$\bullet$ Un morphisme (d'ensembles !)~$\mu \colon G\times G\to G$
tel que les flèches composées
$$\xymatrix{
{G\times G\times G}\ar[rrr]^{(\mu\circ {\rm pr}_{12}, {\rm pr}_3)}&&&{G\times G}\ar[r]^(.53)\mu&G}$$
et
$$\xymatrix{
{G\times G\times G}\ar[rrr]^{({\rm pr}_1, \mu\circ {\rm pr}_{23})}&&&{G\times G}\ar[r]^(.53)\mu&G}$$
coïncident. 

$\bullet$ Un morphisme~$e \colon \got f\to G$, où~$\got f$ est l'objet final de~$\ens$ 
(c'est-à-dire «le» singleton)
tel que les applications composées
$$\xymatrix{
G\ar[rrr]^{(e\circ \pi, {\rm Id}_G)}&&&{G\times G}\ar[r]^(.53)\mu&G}$$
et~$$\xymatrix{
G\ar[rrr]^{({\rm Id}_G, e\circ \pi)}&&&{G\times G}\ar[r]^(.53)\mu&G},$$
où~$\pi$ est l'unique morphisme de~$G$ vers~$\got f$, soient toutes deux
égales à~${\rm Id}_G$. 

$\bullet$ Un morphisme~$i \colon G\to G$ tel que
les applications composées
$$\xymatrix{
G\ar[rrr]^{(i, {\rm Id}_G)}&&&{G\times G}\ar[r]^(.53)\mu&G}$$
et~$$\xymatrix{
G\ar[rrr]^{({\rm Id}_G,i)}&&&{G\times G}\ar[r]^(.53)\mu&G}$$
soient toutes deux égales à~$e\circ \pi$. 

\trois{obj-groupes}
{\em Première définition d'un objet en groupes.}
Soit maintenant~$\mathsf C$ une catégorie. On suppose que~$\mathsf C$ a
un objet final et que le produit cartésien de deux objets existe toujours dans~$\mathsf C$ 
(on peut de manière équivalente requérir que les produits cartésiens de familles finies
d'objets de~$\mathsf C$ existent dans~$\mathsf C$, l'objet final étant alors le produit vide). 
Soit~$\got f$ l'objet final de~$\mathsf C$. Un
{\em objet en groupes}
dans la catégorie~$\mathsf C$ est un objet~$G$ de~$\mathsf C$ muni d'un morphisme
$\mu \colon G\times G\to G$, d'un morphisme~$e\colon \got f\to G$, et d'un morphisme
$i\colon G\to G$ tels que les axiomes catégoriques
du~\ref{reformu-groupe}
soient satisfaits {\em verbatim}. 

\trois{obj-groupes-yoneda}
{\em Seconde définition d'un objet en groupes.}
Il résulte immédiatement
du lemme de Yoneda qu'un objet en
groupes de~$\mathsf C$ est un objet~$G$ 
de~$\mathsf C$ muni, pour tout objet~$T$ de~$\mathsf C$, 
d'une structure de groupe fonctorielle en~$T$ sur~$\hom_{\mathsf C}(T,G)$ ; ou, 
si l'on préfère, d'une factorisation de~$T\mapsto \hom_{\mathsf C}(T,G)$ 
{\em via}
le foncteur d'oubli de~$\gp$ dans~$\ens$. 

\trois{groupe-en-gp}
{\em Exercice}. Montrez qu'un groupe en groupes est un groupe abélien.

\deux{def-schem-en-gpes}
Soit~$S$ un schéma. Il y a d'après ce qui précède
deux manières de se donner une structure de~$S$-schémas en groupes
sur un~$S$-schéma~$G$ : on peut ou bien se donner pour tout~$S$-schéma~$T$
une structure de groupe fonctorielle en~$T$ sur~$G(T)$, ou bien se donner trois morphismes
de~$S$-schémas
$$\mu \colon G\times G \to G, \;e\colon S\to G\;\;\text{et}\;i\colon G\to G$$
satisfaisant
les axiomes requis. En général, la première méthode est nettement plus simple ; nous allons
l'illustrer dans un instant par un exemple. 

\medskip
Avant cela, faisons une remarque. Supposons que~$S$ et~$G$ soient affines, disons~$S=\spec A$
et~$G=\spec B$. Se donner un triplet~$(\mu, e, i)$ comme ci-dessus revient alors 
à se donner un triplet de morphismes de~$A$-algèbres
$$\lambda \colon B\to B\otimes_AB, \;\;\epsilon \colon A\to B\;\;\text{et}\;j \colon B\to B$$
satisfaisant les axiomes «duaux»
de ceux imposés à~$(\mu, e,i)$, que nous vous laissons expliciter ; on dit qu'un tel triplet~$(\lambda, \epsilon, j)$ fait de~$B$
une~{\em $A$-algèbre de Hopf}. 

\deux{exemple-gmz}
{\bf Exemple de schéma en groupes : le groupe multiplicatif sur~$\ZZ$}.
Soit~${\mathbb G}_m$ le~$\ZZ$-schéma~$\spec \ZZ[U,U^{-1}]$. Nous allons le munir d'une structure de~$\ZZ$-schéma en groupes, 
que nous allons définir de deux façons
différentes. 

\trois{gmz-yoneda}
{\em Définition fonctorielle}. 
Soit~$T$ un schéma (ou un~$\ZZ$-schéma, c'est la même chose). On déduit de
\ref{foncteur-ouvert-df}
que~$\psi \mapsto \psi^*U$ établit une bijection fonctorielle en~$T$
entre~${\mathbb G}_m(T)$ et~$\sch O_T(T)^\times$. On voit donc
immédiatement que~${\mathbb G}_m(T)$ hérite d'une structure de groupe
fonctorielle en~$T$, qui fait de~${\mathbb G}_m$ un~$\ZZ$-schéma en groupes
appelé pour des raisons évidentes le {\em groupe multiplicatif}
(sur~$\spec \ZZ$). 

\trois{gmz-hopf}
{\em Définition par une structure d'algèbre de Hopf}. 
Soit~$\lambda$ le morphisme d'anneaux de~$\ZZ[U,U]^{-1}$
dans~$\ZZ[V,V^{-1}]\otimes_{\ZZ}\ZZ[W,W^{-1}]=\ZZ[V,W, V^{-1}, W^{-1}]$ qui
envoie~$U$ sur~$VW$ ; soit~$\epsilon$ le morphisme d'anneaux
de~$\ZZ[U,U^{-1}]$ vers~$\ZZ$ qui envoie~$U$ sur~$1$, et soit~$j$ le morphisme
d'anneaux de~$\ZZ[U,U^{-1}]$ dans lui-même qui échange~$U$ et~$U^{-1}$. 

Nous vous laissons vérifier que~$(\lambda, \epsilon, j)$ satisfait les axiomes
des algèbres de Hopf, et que la structure de schéma en groupes que
ce triplet induit dès lors
sur~${\mathbb G}_m$ est la même que celle définie {\em supra}. 

\trois{conclu-gmz}
Vous observez donc sur cet exemple le phénomène que nous avions annoncé : la définition
d'une structure de schémas en groupes {\em via}
le foncteur des points est en général la plus naturelle. 

\chapter{Schémas projectifs}

\section{Le schéma $\proj B$}
\markboth{Schémas projectifs}{Le schéma $\proj B$}

\subsection*{Un peu d'algèbre graduée}
\deux{def-gradue}
{\bf Définition.}
Un {\em anneau gradué}
est un anneau~$B$ munit d'une décomposition 
en somme directe de groupes abéliens~$B=\bigoplus_{n\in \ZZ}B_n$
telle que~$1\in B_0$ et telle que~$B_n\cdot B_m\subset B_{n+m}$ pour tout~$(m,n)$. On dit que~$B$
est gradué
{\em en degrés positifs}
si~$B_n=\{0\}$ pour tout~$n<0$ (ce sera le plus souvent le cas ici).

\trois{gradue-pasgradue}
Tout anneau (usuel)~$A$ peut être vu comme un anneau gradué en posant~$A_0=A$ et~$A_n=0$ si~$n\neq 0$. 

\trois{bzeto-sousanneau}
Soit~$B=\bigoplus B_n$ un anneau gradué. Nous dirons que le sommande~$B_n$ est 
l'ensemble des {\em éléments homogènes de degré~$n$}
de~$B$ (notez que le degré d'un élément homogène non nul est uniquement déterminé,
et que~$0$ est homogène de tout degré). 

Il résulte  immédiatement des définitions
que~$B_0$ est un sous-anneau de~$B$, et que chaque~$B_n$
est un sous-$B_0$-module de~$B$.

\trois{algebre-gradue}
Si~$A$ est un anneau (usuel), une {\em $A$-algèbre graduée}
est un anneau gradué~$B$ muni d'un morphisme 
de~$A$ dans~$B_0$. 

\trois{exemple-gradue}
{\em Exemple.}
Si~$A$ est un anneau et~$n\in\NN$~alors~$A[T_0,\ldots, T_n]$ a une structure 
naturelle de~$A$-algèbre
graduée en degrés positifs : pour tout~$d\in \NN$, l'ensemble de ses éléments homogènes 
de degré~$d$ est le~$A$-module engendré par les~$\prod T_i^{e_i}$ avec~$\sum e_i=d$. 

\trois{ideal-homogene}
Soit~$B=\bigoplus B_n$ un anneau gradué. Si~$I$ est un idéal de~$B$, on vérifie immédiatement que les conditions
suivantes sont équivalentes : 

\medskip
i) $I=\bigoplus (I\cap B_n)$ ; 

ii) l'idéal~$I$ possède une famille génératrice constituée d'éléments homogènes ; 

iii) pour tout~$b\in B$, on a~$b\in I$ si et seulement si c'est le cas
de chacune de
ses composantes homogènes.

\medskip
Lorsqu'elles sont satisfaites, on dit que~$I$ est {\em homogène}. La décomposition
$$B/I=\bigoplus B_n/(I\cap B_n)$$ fait alors de~$B/I$ un anneau gradué ; si~$B$
est gradué en degrés positifs, il en va de même de~$I$. 

\trois{homogene-premier}
Soit~$B$ un anneau gradué et soit~$I$ un idéal homogène de~$B$. Nous vous invitons
à vérifier que~$I$ est premier si et seulement si~$I\neq \{0\}$ et si
$$(ab\in I)\Rightarrow a\in I\;\text{ou}\;b\in I$$
pour tout couple~$(a,b)$ d'éléments
{\em homogènes}
de~$B$. 

\trois{s-mult-homogene}
Soit~$B=\bigoplus B_n$ un anneau gradué et soit~$S$ une partie multiplicative
de~$B$ {\em constituée d'éléments homogènes}. L'anneau~$S^{-1}B$ hérite alors
d'une graduation naturelle, pour laquelle~$(S^{-1}B)_n$ est l'ensemble des éléments
pouvant s'écrire sous la forme~$\frac a s$ avec~$a\in B_{m+n}$ et~$s\in S\cap B_m$
pour un certain~$m$. Notez que même si~$B$ est graduée en degrés positifs ce
n'est pas forcément le cas
de~$S^{-1}B$ ; ainsi, si~$0\notin S$ et si~$s$ est un élément de~$S$ de degré~$m>0$
alors~$\frac 1 s$ est un élément non nul et homogène de degré~$(-m)$ de~$S^{-1}B$. 

\trois{prod-tens-gradue}
Soit~$A$ un anneau, soit~$B$ une~$A$-algèbre graduée et soit~$C$
une~$A$-algèbre. La décomposition~$B=\bigoplus B_n$ 
induit une décomposition
$$C\otimes_A B=\bigoplus (C\otimes_A B_n)$$
qui fait de~$C\otimes_A B$ une~$C$-algèbre graduée. 

\deux{mor-ann-gradue}
Soient~$B$ et~$C$ deux anneaux gradués et soit~$d$ un entier.
Un morphisme $\phi \colon B\to C$ est dit 
{\em homogène de degré~$d$}
si~$\phi(B_n)\subset C_{nd}$
pour tout~$n$. Si~$B$ a une structure d'algèbre graduée
sur un certain anneau~$A$, un tel morphisme fait de~$C$ une~$A$-algèbre graduée. 
Donnons maintenant quelques exemples.

\trois{ex-mor-gradue}
Soit~$A$ un anneau, et soient~$d,n$ et~$m$ trois entiers. Soit~$(P_0,\ldots, P_n)$ une famille
de polynômes homogènes de degré~$d$ appartenant à~$A[T_0,\ldots, T_m]$. L'unique morphisme
de~$A$-algèbres de~$A[S_0,\ldots, S_n]$ dans~$A[T_0,\ldots, T_m]$ qui envoie~$S_i$ sur~$P_i$ pour tout~$i$
est homogène
de degré~$d$. 

\trois{id-homogene-mor}
Soit~$B$ un anneau gradué et soit~$I$ un idéal homogène de~$B$. Le
morphisme quotient~$B\to B/I$
est homogène de degré~$1$. 

\trois{mult-homogene-mor}
Soit~$B$ un anneau gradué et soit~$S$ une partie multiplicative de~$B$ constituée d'éléments homogènes. 
La flèche naturelle~$B\to S^{-1}B$ est homogène de degré~$1$.

\trois{prod-tens-gradue-mor}
Soit~$A$ un anneau, soit~$B$ une~$A$-algèbre graduée et soit~$C$ une~$A$-algèbre. 
Le morphisme naturel de~$A$-algèbres $B\to C\otimes_A B$ est alors homogène de degré~$1$. 

\subsection*{Construction de~$\proj B$}

\deux{intro-proj}
Soit~$B$ un anneau gradué en degrés positifs. On désigne
par~$B\pos$
l'idéal homogène~$\bigoplus_{n>0} B_n$
de~$B$, par~$B\homog$ l'ensemble des éléments
homogènes de~$B$, et par~$B\posh$ l'ensemble des 
éléments homogènes de~$B\pos$ ; en d'autres termes,
$$B\homog=\bigcup_{n\geq 0}B_n
\;\;\text{et}\;\;B\posh=\bigcup_{n>0}B_n.$$

\deux{def-projb}
Le but de ce qui suit est d'associer à l'anneau gradué~$B$
un schéma, par une construction analogue à celle du spectre
mais un peu plus compliquée. 

\trois{def-projb-ensemble}
{\em Définition ensembliste}. On note~$\proj B$ l'ensemble des idéaux premiers homogènes de~$B$
{\em ne contenant pas~$B\pos$}. 

\trois{def-projb-topologie}
{\em La topologie de~$\proj B$.}
Si~$I$ est un idéal homogène de~$B$, on note~$V(I)$ l'ensemble des~$\got p\in \proj B$ tel que~$I\subset \got p$
(ou, ce qui revient au même, tels que~$\got p$ contienne tous les éléments
homogènes de~$I$).
On vérifie immédiatement que les sous-ensembles de~$\proj B$
de la forme~$V(I)$ sont les fermés d'une topologie, dont on munit~$\proj B$. Il 
résulte aussitôt des définitions que si~$\got p\in \proj B$, 
son adhérence dans~$\proj B$ est égale à~$V(\got p)$ ; et qu'une
partie de~$\proj B$ est ouverte si et seulement si elle est 
réunions de parties de la 
forme~
$$D(f):=\{\got p \in \proj B, f\notin \got p\}$$ où~$f\in B\homog.$

\medskip
Soit~$(f_i)$ une famille d'éléments
de~$B\posh$ engendrant~$B\pos$. Si~$\got p\in \proj B$
il ne contient pas~$B\pos$, et ne saurait donc contenir tous les~$f_i$. 
Par conséquent, $\proj B=\bigcup_i D(f_i)$.

\trois{rem-topologie-projb}
{\em Remarque}. Par définition, $\proj B$ est un sous-ensemble de~$\spec B$. 
Si~$f\in B\homog$, l'ouvert~$D(f)$ de~$\proj B$
est
simplement
l'intersection
de l'ouvert~$D(f)$ de~$\spec B$ avec~$\proj B$. 

Soit maintenant~$g\in B$ ; écrivons~$g=\sum_{n\leq N} g_n$ avec~$g_n\in B_n$ pour tout~$n$. 
Soit~$\got p\in \proj B$.
En vertu de la caractérisation d'un idéal homogène par la propriété~iii)
de~\ref{ideal-homogene}, on a
l'équivalence 
$$g\notin \got p \iff \exists i\; \; g_i\notin \got p.$$ 
En conséquence, $D(g)\cap \proj B$ est l'ouvert~$\bigcup_{i\leq n} D(g_i)$ de~$\proj B$. 

\medskip
Il découle de ce qui précède 
que la topologie de~$\proj B$ coïncide avec la topologie induite par celle de~$\spec B$. 

\trois{def-projb-fasci}
{\em Le faisceau d'anneaux sur~$\proj B$}.
Si~$U$ est un ouvert de~$\proj B$, on note $S\homog(U)$ 
l'ensemble des éléments~$f$ de~$B\homog$
tels que~$U\subset D(f)$. C'est une partie multiplicative
de~$B$ constituée d'éléments homogènes ; l'anneau~$S\homog(U)^{-1}B$
hérite donc d'après~\ref{s-mult-homogene}
d'une graduation naturelle 
et en particulier d'un sous-anneau~$(S\homog(U)^{-1}B)_0$. 

\medskip
La flèche~$U\mapsto (S\homog(B)^{-1}U)_0$ est un préfaisceau d'anneaux sur~$\proj B$ ; 
on note~$\sch O_{\proj B}$ le 
faisceau associé. 

\deux{theo-projb-schema}
{\bf Théorème.}
{\em On conserve les notations
de~\ref{intro-proj}
{\em et sq.}

\medskip
i) 
Pour tout~$f\in B\posh$ l'espace
annelé~$(D(f),\sch O_{\proj B}|_{D(f)})$ est canoniquement
isomorphe à~$\spec (B_f)_0$.

ii)
L'espace annelé~$(\proj B, \sch O_{\proj B})$ est un schéma.}

\medskip
{\em Démonstration}. 
On sait que~$\proj B=\bigcup_{f\in B\posh} D(f)$ (\ref{def-projb-topologie}) ; l'assertion~ii)
est dès lors
une conséquence immédiate de~i), et du fait qu'être un schéma est, pour un espace annelé
donné, une propriété locale. 
Il suffit donc de démontrer~i). 

\medskip
Par hypothèse, $f\in B_r$ pour un certain~$r$ {\em strictement positif}. 
Pour tout~$\got p$
appartenant à~$\proj B$ et tout~$n$, on pose~$\got p_n=\got p\cap B_n$. 

\trois{df-bf0-appli}
{\em Construction d'une application~$D(f)\to \spec (B_f)_0$.}
L'inclusion continue
$\proj B\hookrightarrow \spec B$ envoie l'ouvert~$D(f)$ de~$\proj B$ 
dans l'ouvert~$D(f)$ de~$\spec B$, qui lui-même s'identifie à~$\spec B_f$. L'inclusion
de~$(B_f)_0$ dans~$B_f$ induit par ailleurs une application continue~$\spec B_f\to \spec (B_f)_0$. 
Par composition, on obtient une application continue
de~$D(f)$ vers~$\spec (B_f)_0$. 
Donnons-en
une description explicite. Soit~$\got p\in \proj B$
et soit~$\got q$ son image dans~$\spec (B_f)_0$. 
Il résulte alors de~\ref{ideaux-smoinsuna}
que pour tout entier~$m$ et tout~$a\in B_{rm}$, 
l'élément~$\frac a {f^m}$ de~$(B_f)_0$
appartient à~$\got q$ si et seulement si~$a\in \got p_{rm}$. 

\trois{rem-pre-dfb0}
{\em Une remarque}. 
Soit~$n$ un entier. On note~$\lambda(n)$ l'ensemble
des couples~$(d,\delta)$ avec~$d>0$ tels que~$dn=r\delta$. Le fait que~$r>0$
par hypothèse a deux conséquences : $\lambda(n)$ contient~$(r,n)$ (et est en particulier non vide) ; 
et si~$(d,\delta)\in \lambda(n)$ alors~$\delta$ est égal à~$nd/r$, et est donc entièrement déterminé par~$d$. 

\medskip
Soient maintenant~$n\in \NN, a\in B_n$ et~$(d_0,\delta_0)\in \lambda(n)$ ; soit~$\got q$
un idéal premier de~$(B_f)_0$. Les assertions
suivantes sont équivalentes : 

\medskip
1) $\frac {a^d}{f^\delta}\in \got q$ pour tout~$(d,\delta)\in \lambda(n)$ ; 

2)  $\frac {a^{d_0}}{f^{\delta_0}}\in \got q$. 

\medskip
Il est en effet clair que~1)$\Rightarrow$2). Supposons maintenant que~2)
soit vraie, et soit~$(d,\delta)\in \lambda(n)$ ; nous allons montrer que
$\frac {a^d}{f^\delta}\in \got q$. 
Les couples~$(dd_0, \delta d_0)$ 
et~$(dd_0,\delta_0d)$ appartiennent tous deux à~$\lambda(n)$, ce qui implique que~$\delta d_0=\delta_0d$. 
On a alors
$$\left(\frac{a^{d_0}}{f^{\delta_0}}\right)^d=\frac{a^{dd_0}}{f^{d\delta_0}}=\frac{a^{dd_0}}{f^{\delta d_0}}
=\left(\frac{a^d}{f^\delta}\right)^{d_0}.$$ 
Comme~$d>0$ on a~$\left(\frac{a^d}{f^\delta}\right)^{d_0}
=\left(\frac{a^{d_0}}{f^{\delta_0}}\right)^d\in \got q$
et~$\frac{a^d}{f^\delta}\in \got q$ puisque~$\got q$ est premier.

\trois{df-bf0-inj}
{\em L'application~$D(f)\to \spec (B_f)_0$ est injective}. 
Soit~$\got q$
un idéal premier de~$(B_f)_0$ ; il s'agit de montrer qu'il est l'image d'au plus un élément~$\got p$
de~$\proj B$.

\medskip
Soit donc~$\got p$ un antécédent de~$\got q$ dans~$D(f)$. Pour montrer que~$\got p$ est uniquement déterminé, 
il suffit de s'assurer que~$\got p_n$ est uniquement déterminé pour tout~$n$, puisque~$\got p$ est homogène. 
Fixons donc~$n$ et soit~$a\in B_n$. Si~$a\in \got p_n$ alors~$a^d\in \got p_{dn}$ pour tout~$d>0$, 
et~$\frac{a^d}{f^\delta}$
appartient donc à~$\got q$
pour tout~$(d,\delta)\in \lambda(n)$. 

Réciproquement, supposons que~$\frac{a^d}{f^\delta}\in \got q$
pour tout~$(d,\delta)\in \lambda(n)$, et choisissons un tel~$(d,\delta)$. 
On a alors~$a^d\in \got p$, d'où~$a\in \got p_n$ puisque~$\got p$
est premier. En conséquence, 
$\got p_n$ est {\em nécessairement}
l'ensemble des éléments~$a\in B_n$ tels que~$\frac{a^d}{f^\delta}$
appartienne à~$\got q$
pour tout~$(d,\delta)\in \lambda(n)$, d'où l'unicité de~$\got p$. 

\trois{df-bf0-surj}
{\em L'application~$D(f)\to \spec (B_f)_0$ est surjective}.
On s'inspire comme il se doit de la formule qu'on a exhibée lors de la preuve 
de l'injectivité. 

\medskip
Soit
donc~$\got q$ un idéal premier de~$(B_f)_0$. 
Pour tout entier~$n$, 
on définit~$\got p_n$ comme l'ensemble des~$a\in B_n$ tels que~$\frac {a^d}{f^\delta}\in \got q$
pour tout~$(d,\delta)\in \lambda(n)$ ; notons qu'en vertu de la remarque~\ref{rem-pre-dfb0},
il suffit de s'assurer que c'est le cas pour {\em un}
tel~$(d,\delta)$. On pose~$\got p=\sum \got p_n$. Nous allons tout d'abord montrer que~$\got p$
est un élément de l'ouvert~$D(f)$ de~$\proj B$. 

\medskip
Vérifions
pour commencer
que la somme de deux éléments de~$\got p$ appartient à~$\got p$.
On peut raisonner composante homogène par composante homogène. 
Soient donc~$n\in \NN$ et~$a$ et~$b$
deux éléments
de~$\got p_n$ ; nous allons
prouver que~$(a+b)\in \got p_n$. 

Soit~$(d,\delta)\in \lambda(n)$. On a alors~$(2d,2\delta)\in \lambda(n)$, et il suffit de montrer
que~$\frac{(a+b)^{2d}}{f^{2\delta}}\in \got q$. 
Or lorsqu'on développe~$\frac{(a+b)^{2d}}{f^{2\delta}}$, on obtient une somme de termes 
qui sont de la forme 
$\frac{a^ib^j}{f^{2\delta}}$ avec~$i+j=2d$. Dans un tel terme, l'un des deux entiers~$i$ et~$j$ est au moins égal à~$d$, 
et le terme en question est donc multiple de~$\frac {a^d}{f^\delta}$ ou de~$\frac {b^d}{f^\delta}$,
et il appartient en conséquence à~$\got q$, puisque~$a$ et~$b$ appartiennent
à~$\got p_n$. Il en résulte que~$\frac{(a+b)^{2d}}{f^{2\delta}}\in \got q$, ce qu'on souhait
établir. 

\medskip
Il est immédiat que~$0\in \got p$ et que~$\got p$ est stable par multiplication externe par les éléments de~$B$ ; en conséquence,
$\got p$ est un idéal de~$B$, qui est homogène par sa forme même. Puisque~$\got q$ est premier, $1=\frac f f \notin \got q$, et~$f\notin
\got p$ ; en particulier~$\got p\neq B$. 

\medskip
Montrons que l'idéal homogène strict~$\got p$ de~$B$ est premier. Soient~$n$ et~$m$ deux entiers et soient~$a\in B_n$ et~$b\in B_m$ tels
que~$ab\in \got p$. 
On a alors~$\frac{(ab)^r}{f^{n+m}}=\frac{a^r}{f^n}\cdot \frac {b^r}{f^m}\in \got q$, et donc~$\frac{a^r}{f^n}\in \got q$
ou~$\frac{b^r}{f^m}\in \got q$ puisque~$\got q$ est premier ; ainsi, $a\in \got p_n$ ou~$b\in \got p_m$, ce qu'il fallait démontrer. 

\medskip
On a vu ci-dessus que~$f\notin \got p$, et~$\got p$ appartient
dès lors à l'ouvert~$D(f)$ de~$\proj B$ (notons que le fait que~$f\notin \got p$ garantit
que~$\got p$ ne contient pas~$B\pos$). 

\medskip
Il reste 
à
prouver
que l'image~$\got r$ de~$\got p$  sur~$\spec (A_f)_0$ est égale à~$\got q$.
Soit~$n$ un entier et soit~$a\in B_{rn}$. Il résulte de
la définition de~$\got p$ et de la description explicite de~$\got r$ en fonction
de ce dernier qu'on a les équivalences

$$\frac a {f^n}\in \got r\iff a\in \got p_{rn}\iff \frac a {f^n}\in \got q,$$
ce qui assure
que~$\got r=\got q$ et achève 
de démontrer
que~$D(f)\to \spec (B_f)_0$
est surjective. 

\trois{df-bf0-homeo}
{\em La bijection continue~$D(f)\to \spec (B_f)_0$
est un homéomorphisme.}
Soit~$n\in \NN$ et soit~$a\in B_n$. Choisissons
un couple~$(d,\delta)$ dans~$\lambda(n)$. Soit~$\got p\in D(f)$ 
et soit~$\got q$ son image sur~$\spec (B_f)_0$. Il résulte des descriptions
de~$\got q$ en fonction de~$\got p$ (\ref{df-bf0-appli})
et de~$\got p$ en fonction de~$\got q$ (\ref{df-bf0-surj})
que~$a\in \got p$ si et seulement si~$\frac {a^d}{f^\delta}\in \got q$. 
En conséquence, la bijection continue~$D(f)\simeq \spec (B_f)_0$
identifie le fermé~$V(a)\cap D(f)$ de~$D(f)$ au fermé~$V(\frac{a^d}{f^\delta})$
de~$\spec (B_f)_0$ ; par passage au complémentaire, elle identifie également
l'ouvert~$D(a)\cap D(f)$ de~$D(f)$ à l'ouvert~$D(\frac{a^d}{f^{\delta}})$ de~$\spec (B_f)_0$. 

\medskip
Puisque les ouverts de~$D(f)$ de
la forme~$D(a)\cap D(f)$ avec~$a\in b\homog$
constituent une base de la topologie de~$D(f)$, 
la flèche~$D(f)\to \spec (B_f)_0$ est un homéomorphisme.

\trois{df-bf0-espann}
{\em L'homéomorphisme~$D(f)\simeq \spec (B_f)_0$ est sous-jacent
à un isomorphisme d'espaces localement annelés}. Soit~$U$ un ouvert de~$D(f)$,
et soit~$V$ son image sur~$\spec (B_f)_0$. Notons~$S(V)$ l'ensemble des éléments
de~$(B_f)_0$ qui ne s'annulent en aucun point de~$V$.

\medskip
Commençons par une remarque que nous allons utiliser implicitement
plusieurs fois dans la suite. Soient~$n\in\NN$ et~$a\in B_{rn}$. 
Si~$\got p\in U$ et si~$\got q$ 
désigne son image sur~$V$, alors~$\frac{a}{f^n}$
appartient à~$\got q$ si et seulement si~$a$ appartient 
à~$\got p$. 
Il en résulte que~$\frac a{f^n}\in S(V)$ si et seulement si~$a\in S\homog(U)$. 

\medskip
Par hypothèse, $f\in S\homog(U)$ ; en conséquence, on dispose
d'un morphisme d'anneaux de~$B_f$ vers~$S\homog(U)^{-1}B$, qui est 
par construction 
homogène de degré~$1$ et envoie en particulier~$(B_f)_0$ dans~$(S\homog(U)^{-1}B)_0$. 
Il découle de la remarque précédente que ce morphisme envoie~$S(V)$ dans l'ensemble
des éléments inversibles de~$(S\homog(U)^{-1}B)_0$ ; il se factorise
dès lors par une flèche
de~$S(V)^{-1}(B_f)_0$ dans~$(S\homog(U)^{-1}B)_0$, dont nous allons montrer qu'elle est bijective. 

\medskip
{\em Preuve de l'injectivité}. Soient~$n$ et~$m$ deux entiers, et soient~$a$ et~$b$
appartenant respectivement à~$B_{rn}$ et~$b\in B_{rm}$
tels que~$\frac{a}{f^n}\in S(V)$,
et tels que l'image de~$\left(\frac a{f^n}\right)^{-1}\cdot \frac b{f^m}$ dans~$(S\homog(U)^{-1}B)_0$ soit nulle. Cette image est égale
à~$\frac{bf^n}{af^m}$ ; dire qu'elle est nulle signifie qu'il existe
un élément~$s\in S\homog(U)$ tel que~$sbf^n=0$. Soit~$\ell$ tel que~$s\in B_\ell$
et soit~$(d,\delta)\in \lambda(\ell)$. Comme~$s$
appartient à~$S\homog(U)$, on a~$\frac {s^d} {f^\delta}\in S(V)$. 
L'égalité~$sbf^n=0$ implique que l'élément~$\frac {s^d} {f^\delta} \cdot \frac b {f^m}$ de~$(B_f)_0$ est nul, et donc que l'élément~$\frac b{f^m}$
de~$S(V)^{-1}(B_f)_0$ est nul ; l'élément~$\left(\frac a{f^n}\right)^{-1}\cdot \frac b{f^m}$ 
de~$S(V)^{-1}(B_f)_0$
est {\em a fortiori}
nul, ce qu'il fallait démontrer. 

\medskip
{\em Preuve de la surjectivité}. 
Soit~$n$ un entier, soit~$a\in B_n$ et soit~$s$
appartenant à~$S\homog(U)\cap B_n$. Soit~$(d,\delta)\in \lambda(n)$ ; 
notons que~$\frac {s^d}{f^\delta}\in S(V)$ puisque~$s\in S\homog(U)$. On a
dans l'anneau~$(S\homog(U)^{-1}B)_0$ les égalités

$$\frac a s =\frac {as^{d-1}}{s^d}=\frac {f^\delta}{s^d}\cdot \frac{as^{d-1}}{f^\delta},$$
et~$\frac a s$ est donc l'image de l'élément~$\left(\frac {s^d}{f^\delta}\right)^{-1}\frac{as^{d-1}}{f^\delta}$
de~$S(V)^{-1}(B_f)_0$.

\medskip
Par ce qui précède, la restriction du préfaisceau
$U\mapsto (S\homog(U)^{-1}B)_0$ à
l'ouvert~$D(f)$ s'identifie, {\em via}
l'homéomorphisme~$D(f)\simeq \spec (B_f)_0$, au préfaisceau~$V\mapsto S(V)^{-1}(B_f)_0$. 
En conséquence,~$\sch O_{\proj B}|_{D(f)}$ s'identifie au faisceau associé à~$V\mapsto S(V)^{-1}(B_f)_0$,
qui n'est autre que~$\sch O_{\spec (B_f)_0}$ ; ceci achève la démonstration du théorème.~$\Box$ 

\deux{comment-projb}
{\bf Premières propriétés de~$\proj B$}. 

\trois{projb-nonvide}
Le schéma~$\proj B$ étant
la réunion des~$D(f)$ pour~$f$ parcourant~$B\posh$, il est vide
si et seulement si~$D(f)=\varnothing$
pour tout tel~$f$. En vertu du théorème~\ref{theo-projb-schema},
cela revient à demander que~$(B_f)_0$ soit nul pour tout~$f\in B\posh$.
Or comme~$B(f)_0$ est un sous-anneau de~$B_f$, on a~$1=0$ dans~$(B_f)_0$ 
si et seulement si c'est le cas dans~$B_f$ ; en d'autres termes, $(B_f)$ est nul
si et seulement si~$(B_f)_0$ est nul. 

\medskip
Ainsi, $\proj B$ est vide si et seulement si~$B_f$ est nul pour tout~$f\in B\posh$, 
c'est-à-dire encore si et seulement si tout élément de~$B\posh$ est nilpotent. 

\trois{anneau-loc-projb}
Soit~$\got p\in \proj B$,
et soit~$\Sigma$ l'ensemble des éléments de~$B\homog$
n'appartenant pas à~$\got p$. C'est une partie multiplicative
de~$B$, constituée par définition d'éléments homogènes. Nous 
laissons le lecteur vérifier que l'anneau
local~$\sch O_{\proj B,\got p}$ s'identifie
à~$(\Sigma^{-1}B)_0$.

\trois{cs-projbreduit}
Supposons que~$B$ soit réduit. Dans ce cas, $B_f$ est réduit pour tout~$f$
appartenant à~$B$,
et en particulier pour tout~$f\in B\posh$. Pour un tel~$f$, le sous-anneau~$(B_f)_0$ de~$B_f$
est alors lui aussi réduit, et le schéma~$D(f)=\spec (B_f)_0$ est donc réduit. 
Comme~$\proj B$ est la réunion des~$D(f)$ pour~$f$ parcourant~$B\posh$, il est réduit. 

\trois{cs-projbintegre}
Supposons que~$B$ soit intègre et que~$B\posh$ ne soit pas réduit
au singleton~$\{0\}$. L'idéal homogène~$\{0\}$ de~$B$ est alors 
premier et ne contient pas~$B\posh$ ; c'est donc un point
de~$\proj B$, dont l'adhérence dans~$\proj B$ est égale à~$V({0})$, c'est-à-dire
à~$\proj B$ tout entier. En conséquence, $\proj B$ est irréductible.

\trois{b-algr-proj-asch}
Supposons que~$B$ soit une algèbre
graduée sur un certain anneau~$A$. 
Le préfaisceau
$$U\mapsto (S\homog(U)^{-1}B)_0$$
(avec les notations de~\ref{def-projb-fasci} )
est alors de manière naturelle un préfaisceau de~$A$-algèbres, et il en va de même
de son faisceau associé. Le schéma~$\proj B$ hérite
par ce biais d'une structure naturelle de~$A$-schéma. 

\subsection*{Fonctorialité partielle de la construction}

\deux{projb-fonct}
Soient~$B$ et~$C$
deux anneaux gradués, soit~$d$ un entier
strictement positif
et soit~$\phi \colon B\to C$
un morphisme d'anneaux homogène de degré~$d$. Les choses ne se passent 
pas aussi bien que pour les spectres d'anneaux puisque~$\phi$ n'induit pas
en général
un morphisme de~$\proj C$ {\em tout entier} vers~$\proj B$ : comme
nous allons le voir, 
il peut être nécessaire de se restreindre à un ouvert. 

\trois{fonct-projb-defensemble}
Soit~$\got q\in \proj C$. On vérifie immédiatement que l'idéal
premier~$\phi^{-1}(\got q)$ de~$B$ est homogène. Il appartient à~$\proj B$
si et seulement si il ne contient pas~$B\posh$, ce qui signifie précisément
qu'il existe~$f\in B\posh$ tel que~$\phi(f)\notin \got q$. 

\medskip
Soit~$\Omega$ l'ouvert de~$\proj C$ égal à la réunion des~$D(\phi(f))$ pour~$f$ parcourant
$B\posh$ ; par ce qui précède, on dispose d'une application naturelle~$\psi$
de~$\Omega$ vers~$\proj B$. 

\medskip
Soit~$f\in B\posh$. Il résulte de notre construction que~$\psi^{-1}(D(f))=D(\phi(f))$
(ce qui montre la continuité de~$\psi$)
et que le diagramme

$$\xymatrix{{D(\phi(f))}\ar@{^{(}->}[rr]\ar[d]_\psi&&{\spec C_{\phi(f)}}\ar[d]\ar[r]&{\spec( C_{\phi(f)})_0}\ar[d]\\
{D(f)}\ar@{^{(}->}[rr]&&{\spec B_f}\ar[r]&{\spec(B_f)_0}}$$commute.
Il s'ensuit que modulo les homéomorphismes
canoniques
$$D(\phi(f))\simeq \spec (C_{\phi(f)})_0\;\;\text{et}\;\;D(f)\simeq \spec (B_f)_0,$$
la restriction de~$\psi$ à~$D(\phi(f))$ est simplement l'application continue
naturelle~$\spec (C_{\phi(f)})_0)\to \spec (B_f)_0$. 

\trois{fonct-projb-defespann}
Soit~$U$ un ouvert de~$\proj B$. Si~$f\in S\homog(U)$ alors~$\phi(f)$
appartient à~$S\homog(\psi^{-1}(U))$. 
Le morphisme~$\phi$ induit donc un morphisme de~$S\homog(U)^{-1}B$ vers~$S\homog(\psi^{-1}(U))^{-1}C$, 
dont on voit aussitôt qu'il est homogène de degré~$d$. En particulier, il envoie~$(S\homog(U)^{-1}B)_0$ vers~$(S\homog(\psi^{-1}(U))^{-1}C)_0$. 
En faisant varier~$U$, on obtient ainsi un morphisme de préfaisceaux
$$[U\mapsto (S\homog(U)^{-1}B)_0]\to \psi_*[V\mapsto (S\homog(V)^{-1}C)_0]$$ puis, par passage
aux faisceaux associés, un morphisme~$\psi^*$
de~$\sch O_{\proj B}$ vers~$\psi_*\sch O_{\proj C}|_\Omega$ ; la donnée de~$\psi^*$
fait de~$\psi$ un morphisme d'espaces annelés de~$\Omega$ vers~$\proj B$.

\trois{descr-expl-mor}
Soit~$f\in B\posh$. Par construction, la flèche
$$\sch O_{\proj B}(D(f))=(B_f)_0\to \sch O_{\proj C}(D(\phi(f))=(C_{\phi(f)})_0$$ 
induite par~$\psi^*$ est simplement le morphisme naturel
de~$(B_f)_0$ vers~$(C_{\phi(f)})_0$ induit par~$\phi$. 

\medskip
On en déduit, grâce à~\ref{fonct-projb-defensemble}
et
au lemme~\ref{espann-speca},
que modulo les
isomorphismes canoniques de
schémas
$$D(\phi(f))\simeq \spec (C_{\phi(f)})_0\;\;\text{et}\;\;D(f)\simeq \spec (B_f)_0,$$
le morphisme d'espaces annelés~$\psi|_{D(\phi(f))}\colon D(\phi(f))\to D(f)$ coïncide avec le morphisme 
naturel~$\spec (C_{\phi(f)})_0\to \spec (B_f)_0$. 

\medskip
Ceci valant pour tout~$f\in B\posh$, le morphisme d'espaces annelés~$\psi$ est un morphisme d'espaces localement annelés 
(cette propriété est en effet de nature locale), c'est-à-dire un morphisme de schémas.

\trois{rem-projb-moraff}
{\em Remarque}.
Comme~$\psi^{-1}(D(f))=D(\phi(f))$ pour tout~$f\in B\posh$,
le morphisme~$\psi$ est affine.

\deux{exemple-projb-omegavide}
L'ouvert~$\Omega$ de définition du morphisme~$\psi$ ci-dessus est la réunion 
des~$D(\phi(f))$ pour~$f\in B\posh$ ; on vérifie immédiatement
qu'on peut se contenter de faire parcourir à~$f$ une partie de~$B\posh$ engendrant l'idéal~$B\pos$. 

\medskip
On prendra garde que~$\Omega$ n'a en général
aucune raison d'être égal à~$\proj C$ tout entier (il l'est toutefois dans deux situations
particulières importantes que nous étudierons ci-dessous en~\ref{ihomog-ferme}
{\em et sq.}, ainsi qu'en~\ref{proj-changebase}).
Il peut même, comme
le montre l'exemple détaillé
au~\ref{proj-kx}
ci-dessous, être {\em vide}
sans que~$\proj C$ le soit.

\deux{proj-kx}
Soit~$A$ un anneau. 
Puisque~$T$ engendre~$A[T]\posh$, le~$A$-schéma
$\proj A[T]$ s'identifie à son ouvert~$D(T)$, 
c'est-à-dire d'après le théorème~\ref{theo-projb-schema}
au spectre de~$(A[T]_{(T)})_0=A[T,T^{-1}]_0=A$. Autrement dit, 
$\proj A[T]\to \spec A$ est un isomorphisme (en particulier, $\proj A[T]$
est non vide dès que~$A$ est non nul). 

\medskip
Soit~$\phi$ le morphisme de~$A$-algèbres de~$A[T]$ dans lui-même qui envoie~$T$ sur~$0$. Il est 
homogène de degré~$1$ (et d'ailleurs de n'importe quel autre degré). Il induit donc un morphisme~$\psi$ 
d'un ouvert~$\Omega$ de~$\proj A[T]$ vers~$\proj A[T]$. L'idéal~$A[T]\pos$ étant engendré par~$T$, on a 
$$\Omega =D(\phi(T))=D(0)=\varnothing.$$

\deux{ihomog-ferme}
Soit~$B$ un anneau gradué et soit~$I$ un idéal homogène de~$B$. Le morphisme
quotient~$B\to B/I$ étant homogène de degré~$1$, il induit un morphisme~$\psi \colon \Omega \to \proj B$,
où~$\Omega$ est la réunion 
des~$D(\bar f)$ pour~$f$ parcourant~$B\posh$. Or par définition de la graduation de~$B/I$, 
l'ensemble des~$\bar f$ pour~$f$ parcourant~$B\posh$ est égal à~$(B/I)\posh$ ; il s'ensuit
que~$\Omega$ est la réunion des~$D(g)$ pour~$g$ parcourant~$(B/I)\posh$, et donc que~$\Omega$
est égal à~$\proj B/I$ tout entier.

\trois{ihomog-immferm}
Soit~$n>0$ et soit~$f\in B_n$. L'image réciproque de~$D(f)= \spec B_f$ sur~$\proj B/I$ est l'ouvert~$D(\bar f)= \spec ((B/I)_{\bar f})_0$. 
Il est immédiat que~$IB_f$ est un idéal homogène de~$B_f$, et que
l'isomorphisme canonique~$(B/I)_{\bar f}\simeq (B_f)/(IB_f)$ est homogène de degré~$1$. Par conséquent,
$$D(\bar f)= \spec ((B_f)/(IB_f))_0=\spec (B_f)_0/(IB_f\cap (B_f)_0).$$

La flèche~$D(\bar f)\to D(f)$ induite par~$\psi$ est donc
l'immersion fermée associée à
l'idéal
$IB_f\cap (B_f)_0$ de~$(B_f)_0$, idéal qui est simplement l'ensemble des éléments
de la forme~$\frac a {f^r}$ avec~$r\in \NN$ et~$a\in I\cap B_{rn}$. 

%

\trois{conclu-projbi-immersion}
Comme être une immersion fermée
est une propriété locale sur le but, on déduit de ce qui précède
que  le morphisme~$\psi \colon \proj B/I\to \proj B$
est une immersion fermée. Nous allons déterminer son image. 
 
 \medskip
Soit~$\got q\in \proj B/I$. Par définition, $\psi(\got q)$ est l'image
réciproque de~$\got q$ dans~$B$, qui 
est un idéal premier homogène de~$B$ {\em ne contenant pas~$B\posh$}
(c'est précisément ce que signifie l'égalité~$\Omega=\proj B/I$). 

On en déduit que pour tout~$\got q\in \proj B/I$, 
l'idéal~$\psi(\got q)$ de~$B$ contient~$I$. 
Inversement, si~$\got p$ est un élément de~$\proj B$ qui contient~$I$, on vérifie
aussitôt que~$\got q:=\got p/I$ est un idéal premier homogène de~$B/I$ ne contenant pas~$(B/I)\pos$, c'est-à-dire un élément
de~$\proj B/I$, et que~$\psi(\got q)=\got p$. 

\medskip
L'immersion fermée~$\psi \colon \proj B/I\hookrightarrow \proj B$ a donc pour image~$V(I)$.

\deux{projb-aschema-mor}
Soit~$A$ un anneau et soient~$B$
et~$C$ deux~$A$-algèbres graduées. Soit~$d$ un entier et soit~$\phi \colon B \to C$
un morphisme de~$A$-algèbres homogène de degré~$d$. Il induit d'après~\ref{projb-fonct}
{\em et sq.} un morphisme~$\psi$ d'un ouvert~$\Omega$
de~$\proj C$ vers~$\proj B$. On vérifie 
sans peine
que les morphismes d'anneaux intervenant dans la construction
de~$\psi$ sont des morphismes de~$A$-algèbres ; en conséquence, $\psi$ est un morphisme de~$A$-schémas. 

\deux{proj-changebase}
Soit~$A$ un anneau, soit~$B$ une~$A$-algèbre graduée, et soit~$C$ une~$A$-algèbre. 
Le schéma~$\proj (C\otimes_A B)$ est un~$C$-schéma d'après~\ref{b-algr-proj-asch}. 
On dispose
par ailleurs
d'un~$A$-morphisme naturel 
$$\phi \colon B\to C\otimes_A B$$ qui est homogène de degré~$1$, et induit donc un~$A$-morphisme~$\psi$ 
d'un ouvert~$\Omega$ de~$\proj C\otimes_A B$ vers~$\proj B$. Cet ouvert~$\Omega$ est la réunion
des~$D(\phi(f))$ pour~$f$
parcourant~$B\posh$. Mais comme~${\phi(f)}_{f\in B\posh}$ engendre~$(C\otimes_A B)\pos=C\otimes_A(B\pos)$, 
l'ouvert~$\Omega$ est en fait égal à~$\proj (C\otimes_A B)$ tout entier.  

\medskip
Le
diagramme commutatif
$$\xymatrix{
{\proj(C\otimes_AB)}\ar[d]_\psi\ar[r]&{\spec C}\ar[d]\\
{\proj B}\ar[r]&{\spec A}}$$
définit un morphisme~$\proj(C\otimes_AB)\to \proj B\times_{\spec A}\spec C$. 
Nous allons montrer que c'est un isomorphisme ; on peut 
pour ce faire raisonner localement sur~$\proj B$. Soit~$f\in B\posh$. 
On a 
$$D(f)\simeq \spec (B_f)_0\;\;\text{et}\;\;\psi^{-1}(D(f))=D(\phi(f))\simeq \spec ((C\otimes_AB)_{\phi(f)})_0.$$
On vérifie immédiatement que l'isomorphisme
$C\otimes_AB_f\simeq (C\otimes_AB)_{\phi(f)}$ est homogène de degré~$1$ ; il vient
$$C\otimes_A ((B_f)_0)=(C\otimes_AB_f)_0\simeq ((C\otimes_AB)_{\phi(f)})_0.$$
Par conséquent, $\psi^{-1}(D(f))\simeq D(f)\times_{\spec A}\spec C$, d'où notre assertion. 

\section{Le schéma~$\PP^n_A$}
\markboth{Schémas projectifs}{Le schéma~$\PP^n_A$}

\deux{def-pna}
{\bf Définition}.
Soit~$A$ un anneau et soit~$n\in \NN$. On note~$\PP^n_A$
le~$A$-schéma~$\proj A[T_0,\ldots, T_n]$. 

\trois{changebase-pna}
Si~$B$ est une~$A$-algèbre, on dispose d'après~\ref{proj-changebase}
d'un isomorphisme canonique~$\PP^n_B\simeq \PP^n_A\times_{\spec A}\spec B$. En particulier, 
on a~$\PP^n_A=\PP^n_{\ZZ}\times_{\spec \ZZ}\spec A$. 

\medskip
Si~$x\in \spec A$, la fibre de~$\PP^n_A$ en~$x$ s'identifie à~$\PP^n_A\times_{\spec A}\spec \kappa(x)$, 
c'est-à-dire à~$\PP^n_{\kappa(x)}$. 

\trois{ex-p0a}
{\em Premiers exemples}.
Si~$A=\{0\}$ alors~$\PP^n_A$ est vide pour tout~$n$. Si~$A\neq \{0\}$
l'idéal~$A[T_0,\ldots, T_n]\pos$ de~$A[T_0,\ldots, T_n]$ n'est pas constitué d'éléments
nilpotents (par exemple, les~$T_i$ ne sont pas nilpotents), et~$\PP^n_A$ est donc non vide
(\ref{projb-nonvide}).  
 
\medskip
Le~$A$-schéma~$\PP^0_A$ est égal
par définition à~$\proj A[T]$. Il résulte alors de~\ref{proj-kx}
que~$\PP^0_A\to \spec A$ est un isomorphisme. 

\trois{pnx-def-rem}
On définit plus généralement, pour tout schéma~$X$, le~$X$-schéma~$\PP^n_X$
comme étant égal à~$\PP^n_{\ZZ}\times_{\spec \ZZ}X$ ; en vertu de~\ref{changebase-pna},
cette définition est compatible avec la précédente lorsque~$X$ est affine. 

\trois{compat-pna-p1k}
Soit~$k$ un corps. La notation~$\PP^1_k$ semble 
{\em a priori}
désigner deux~$k$-schémas différents : 

\medskip
$\bullet$ celui construit par recollement au~\ref{droite-proj-recoll} ; 

$\bullet$ celui défini au~\ref{def-pna}
ci-dessus, à savoir~$\proj k[T_0,T_1]$. 

\medskip
Mais nous verrons un peu plus loin 
au~\ref{conclu-p1-compatible}
que ce conflit de notations n'est qu'apparent et que ces deux~$k$-schémas coïncident. 

\deux{rem-spectre-qcgrad}
{\bf Remarque}.
Le foncteur $A\mapsto \spec A$ admet
une variante globale ou relative, à savoir la formation du spectre d'une
$\sch O_X$-algèbre quasi-cohérente sur un schéma~$X$. On peut de même définir
(mais nous n'en aurons pas besoin) une variante
globale ou relative de~$B\mapsto \proj B$, en introduisant la notion de~$\sch O_X$-algèbre quasi-cohérente
{\em graduée}, et en associant un~$X$-schéma~$\proj \sch B$ à une telle $\sch O_X$-algèbre~$\sch B$ ; nous laissons
le lecteur intéressé deviner puis écrire en détail la construction de~$\proj \sch B$. 

\medskip
On peut grâce à cette notion définir directement~$\PP^n_X$, sans produit fibré : c'est~$\proj \sch O_X[T_0,\ldots, T_n]$. 

\deux{desc-pna-casgeneral}
{\bf Les cartes affines standard de~$\PP^n_A$}. Soit~$A$
un anneau et soit~$n$ un entier. 
L'idéal~$A[T_0, \ldots,T_n]\pos$ de~$A[T_0,\ldots, T_n]$ 
est engendré par les~$T_i$. En conséquence, 
$\PP^n_A$ est la réunion de ses ouverts~$U_i:=D(T_i)$. 

\trois{ann-localis-pna}
Si~$(e_i)_{0\leq i\leq n}$
est une famille d'entiers, nous noterons~$A\left[T_0,\ldots, T_n, \frac{1}{\prod T_i^{e_i}}\right]$
la~$A$-algèbre graduée~$A[T_0,\ldots, T_n]_{\prod T_i^{e_i}}$. On vérifie immédiatement
que le morphisme homogène de degré~$1$ naturel de
$A\left[T_0,\ldots, T_n, \frac{1}{\prod T_i^{e_i}}\right]$ vers~$A\left[T_0,\ldots, T_n, \frac{1}{\prod T_i}\right]$ est injectif ; 
nous nous permettrons donc de considérer implicitement~$A\left[T_0,\ldots, T_n, \frac{1}{\prod T_i^{e_i}}\right]$
comme une sous-algèbre graduée de~$A\left[T_0,\ldots, T_n, \frac{1}{\prod T_i}\right]$. 

\trois{desc-dti-pna}
Soit~$i\in \{0,\ldots, n\}$. En vertu de l'assertion~i)
du théorème~\ref{theo-projb-schema}, on dispose d'un~$A$-isomorphisme canonique
$$U_i\simeq \spec A\left[T_0,\ldots, T_n, \frac{1}{T_i}\right]_0.$$ On montre aisément
que~$A\left[T_0,\ldots, T_n, \frac{1}{T_i}\right]_0$ est simplement 
l'algèbre de polynômes
en~$n$-variables~$A\left[\frac{T_\ell}{T_i}\right]_{\ell \neq i}$. On a donc
$$U_i\simeq \spec A\left[\frac{T_\ell}{T_i}\right]_{\ell \neq i}\simeq \Aff^n_A.$$ 

 \trois{bord-carte-aff-pn}
Fixons~$i$. Le fermé complémentaire de l'ouvert~$D(T_i)$ est par définition 
$V(T_i)$. Le quotient de~$A[T_0,\ldots, T_n]$ par son idéal homogène~$(T_i)$ 
est simplement l'anneau gradué~$A[T_j]_{j\neq i}$. La flèche quotient~$A[T_0,\ldots, T_n]\to A[T_j]_{j\neq i}$
induit en vertu de~\ref{ihomog-ferme}
{\em et sq.}
une immersion fermée
$$\proj A[T_j]_{j\neq i} \hookrightarrow A[T_0,\ldots, T_n]$$ d'image~$V(T_i)$,
qui permet de munir celui-ci d'une structure de sous-schéma fermé. 
Comme~$A[T_j]_{j\neq i}$ est isomorphe à~$A[T_0,\ldots, T_{n-1}]$
si~$n\geq 1$
(par renumérotation), le fermé~$V(T_i)$ muni de la structure
en question est isomorphe à~$\PP^{n-1}_A$ dès que~$n\geq 1$ (il est en 
particulier non vide dès que~$A$ est non nul). 
Si~$n=0$ alors~$T_0$ engendre~$A[T_0]\pos$, et~$V(T_0)=\varnothing$.

\trois{pna-copies-ana}
Le~$A$-schéma~$\PP^n_A$ est ainsi réunion de~$n+1$ copies de l'espace
affine relatif~$\Aff^n_A$ ; c'est en particulier
un~$A$-schéma {\em de type fini}. 
Soient~$i$ et~$j$ deux éléments de~$\{0,\ldots, n\}$. 
L'intersection~$U_i\cap U_j$ est l'ouvert~$D(T_iT_j)$ de~$\PP^n_A$ ; en utilisant 
une fois encore le théorème~\ref{theo-projb-schema}, on voit que
$$U_i\cap U_j\simeq \spec (A[T_0,\ldots, T_n]_{T_iT_j})_0.$$
En tant qu'ouvert de~$U_i=\spec  A\left[\frac{T_\ell}{T_i}\right]_{\ell \neq i}$, 
l'intersection~$U_i\cap U_j=U_i\cap D(T_j)$ est égale en vertu de~\ref{df-bf0-homeo}
à
$$D\left(\frac{T_j}{T_i}\right)=\spec A\left[\frac{T_\ell}{T_i}, \frac{T_i}{T_j}\right]_{\ell \neq i}.$$
De même, en tant qu'ouvert de~$U_j=\spec  A\left[\frac{T_\ell}{T_j}\right]_{\ell \neq j}$, 
l'intersection~$U_i\cap U_j$ est égale à
$$D\left(\frac{T_i}{T_j}\right)=\spec A\left[\frac{T_\ell}{T_j}, \frac{T_j}{T_i}\right]_{\ell \neq j}.$$
Ces deux dernières descriptions de~$U_i\cap U_j$ pourraient tout aussi bien se déduire de
la première et des égalités
$$(A[T_0,\ldots, T_n]_{T_iT_j})_0=A\left[\frac{T_\ell}{T_i}, \frac{T_i}{T_j}\right]_{\ell \neq i}=
A\left[\frac{T_\ell}{T_j}, \frac{T_j}{T_i}\right]_{\ell \neq j}$$
entre sous-anneaux de~$A\left[T_0,\ldots, T_n, \frac{1}{\prod T_i}\right]$. 

\deux{fonctions-pna}
{\bf Les fonctions globales sur~$\PP^n_A$}. 
Il résulte de~\ref{desc-pna-casgeneral}
{\em et sq.}
que~$\sch O_{\PP^n_A}(\PP^n_A)$ 
est le sous-anneau de~$A\left[T_0,\ldots, T_n,\frac 1 {\prod T_i}\right]$ constitué
des éléments qui peuvent s'écrire pour tout~$i$ comme un polynôme en les variables
$\frac{T_j}{T_i}$ pour~$j\neq i$. 

\medskip
En tant que~$A$-module, $A\left[T_0,\ldots, T_n,\frac 1 {\prod T_i}\right]$ est libre et admet
pour base la famille des monômes de la forme~$\prod T_i^{e_i}$ où les~$e_i$ appartiennent à~$\ZZ$. 
Fixons~$i$. Un élément de~$A\left[T_0,\ldots, T_n,\frac 1 {\prod T_i}\right]$ peut s'écrire
comme un polynôme en les variables
$\frac{T_j}{T_i}$ pour~$j\neq i$ si et seulement si son écriture dans la base évoquée ne fait intervenir que des monômes
de degré total nul où seul~$T_i$ est autorisé à avoir un exposant négatif.

\medskip
Il est immédiat que seules les constantes peuvent satisfaire cette condition pour tout~$i$ ; en conséquence, la~$A$-algèbre
$\sch O_{\PP^n_A}(\PP^n_A)$ 
est égale à~$A$. 

\medskip
On en déduit que~$\PP^n_A$ n'est pas affine dès que~$A\neq\{0\}$ et dès que~$n\geq 1$. En effet, le morphisme
structural~$\PP^n_A\to \spec \sch O_{\PP^n_A}(\PP^n_A)$ est égal à~$\PP^n_A\to \spec A$ ; or
si~$A\neq\{0\}$ et~$n\geq 1$, cette flèche n'est pas un isomorphisme. Pour le voir, on choisit~$x\in \spec A$
(ce qui est possible puisque~$A\neq\{0\}$ et l'on remarque que la fibre de~$\PP^n_A\to \spec A$ en~$x$, qui s'identifie
à~$\PP^n_{\kappa(x)}$, n'est pas réduite à un singleton (elle contient par exemple une copie de~$\Aff^n_{\kappa(x)}$, et
donc de l'ensemble~$\kappa(x)^n$ de ses points naïfs). 

\medskip
Notons par contre que si~$A=\{0\}$ alors~$\PP^n_A$ est vide (et en particulier affine), et que~$\PP^0_A\simeq \spec A$
est toujours affine.

\deux{desc-alternative-pna}
Grâce à~\ref{desc-pna-casgeneral}, 
on peut donner une deuxième description de~$\PP^n_A$, par recollement 
de cartes affines et sans faire appel à la construction «$\proj$», que nous allons maintenant
esquisser. 

\trois{prepare-recollement-pna}
Pour tout~$i$ compris entre~$0$ et~$n$, on se donne 
une famille~$(\tau_{\ell i})_{0\leq \ell \leq n , \ell \neq i}$ d'indéterminées
et l'on pose~$X_i=\spec A[\tau_{\ell i}]_{\ell\neq i}$ (pour faire 
le lien avec~\ref{pna-copies-ana}, il faut penser que~$\tau_{\ell i}=\frac{T_\ell}{T_i}$). 

\medskip
Pour tout couple~$(i,j)$ d'indices avec~$i\neq j$, le
morphisme~$\phi_{ij}$ de~$A$-algèbres

$$A\left[\tau_{\ell j}, \frac 1{\tau_{ij}}\right]\to A\left[\tau_{\ell i}, \frac 1{\tau_{ji}}\right],$$

$$\tau_{\ell j}\mapsto \tau_{\ell i} \cdot \frac 1 {\tau_{ji}}\;\;(\ell \neq i),\;\;\tau_{ij}\mapsto \frac 1 {\tau_{ji}}$$
est un isomorphisme de réciproque~$\phi_{ji}$. 

Pour tout couple d'indices~$(i,j)$ avec~$i\neq j$, 
on note~$\iota_{ij}$ l'isomorphisme
de~$X_{ij}:=D(\tau_{ji})\subset U_i$ vers~$X_{ji}$ induit
par~$\phi_{ij}$. 

\trois{conclu-recollement-pna}
Il n'est pas difficile de voir que la famille
des~$X_i$, des~$X_{ij}$ et des~$\iota_{ij}$
satisfait 
les conditions~i) et~ii) de~\ref{limind-recoll-schema-cocy}. Cela autorise
à procéder comme expliqué
dans~{\em loc. cit.}
au recollement des~$X_i$ le long
des isomorphismes~$\iota_{ij}$ ; nous vous laissons
vérifier que le~$A$-schéma ainsi obtenu est isomorphe à~$\PP^n_A$
(l'ouvert~$X_i$ s'envoyant sur la carte affine~$D(T_i)$). 

\trois{conclu-p1-compatible}
{\em Le cas particulier où~$n=1$}. 
Par ce qui précède, le~$A$-schéma~$\PP^1_A$
peut se décrire
comme le recollement des
copies~$X_0=\spec A[\tau_{10}]$ et~$X_1=\spec A[\tau_{01}]$ de~$\Aff^1_A$, 
le long des isomorphismes réciproques l'un de l'autre
$$\iota_{10}\colon D(\tau_{10})=\spec A[\tau_{10}, \tau_{10}^{-1}]\to \spec A[\tau_{01},\tau_{01}^{-1}]=D(\tau_{01})$$
$$\text{et}\; \iota_{01}\colon D(\tau_{01})=\spec A[\tau_{01}, \tau_{01}^{-1}]\to \spec A[\tau_{10},\tau_{10}^{-1}]=D(\tau_{10})$$
donnés par les formules
$$\iota_{10}(\tau_{10})=\tau_{01}^{-1}\;\;\text{et}\;\iota_{01}(\tau_{01})=\tau_{10}^{-1}.$$

\medskip
Lorsque~$A$ est un corps, on retrouve très précisément la construction du~\ref{droite-proj-recoll}. 
Il s'ensuit que pour tout corps~$k$, les deux définitions concurrentes du~$k$-schéma
~$\PP^1_k$, à savoir celle de~\ref{droite-proj-recoll}
et celle de~\ref{def-pna},
coïncident. 

\deux{homogene-deshomogene}
{\bf Homogénisation et déshomogénisation.}
Fixons un indice~$i$, et posons~$\tau_j=\frac{T_j}{T_i}$ pour~$j\neq i$. L'ouvert affine~$U_i=D(T_i)$
de~$\PP^n_A$ s'identifie à~$\spec A[\tau_j]_{j\neq i}$ (\ref{desc-dti-pna}). 

\trois{deshomo-ferme}
Soit~$I$ un idéal homogène de~$A[T_0,\ldots, T_n]$, et soit~$B$ la~$A$-algèbre
graduée quotient~$A[T_0,\ldots, T_n]/I$. 
Soit~$(f_\ell)_{\ell}$ une famille génératrice de~$I$ constituée d'éléments homogènes et non nuls ; pour tout~$\ell$, on note~$d_\ell$
le degré de~$f_\ell$. 

\medskip
Pour tout~$\ell$, posons~$g_\ell=f_\ell/T_i^\ell$. 
On peut voir~$g_\ell$ comme le «déshomogénéisé»
de~$f_\ell$ (relativement à la variable~$T_i$) ; il s'obtient à partir de~$f_\ell$ en remplaçant~$T_i$ par~$1$
et~$T_j$ par~$\tau_j$ pour tout~$j\neq i$. 

\medskip
Il résulte de~\ref{ihomog-ferme}
{\em et sq.}
que la flèche quotient~$A[T_0,\ldots, T_n]\to B$ induit une immersion fermée
$$\proj B\hookrightarrow \PP^n_A$$ d'image~$V(I)$, 
et que
l'immersion fermée induite~$\proj B\times_{\PP^n_A}U_i\hookrightarrow U_i$ est définie
par l'idéal~$J:=(g_\ell)_{\ell}$ de~$A[\tau_j]$. En particulier, $V(I)\cap U_i=V(J)$
(ce qu'on aurait pu déduire directement de~\ref{df-bf0-homeo}).

\trois{homo-ferme}
Inversement, soit~$J$ un idéal de~$A[\tau_j]_{j\neq i}$.
Choisissons une famille~$(g_\ell)_{\ell}$
de générateurs de~$J$ constituée d'éléments non nuls. Pour tout~$\ell$, notons~$d_\ell$ le degré de~$g_\ell$, et~$f_\ell$
le polynôme de~$A[T_0,\ldots, T_n]$ déduit de~$g_\ell$
par «homogénéisation». Plus précisément, si l'on écrit~$g_\ell=\sum_{(e_j)_j}a_{(e_j)}\prod \tau_j^{e_j}$
alors~$f_\ell=\sum_{(e_j)}a_{(e_j)}T_i^{d_\ell-\sum e_j}\prod T_j^{e_j}$. Par construction, $f_\ell$ est non nul
et homogène
de degré~$\ell$, et~$g_\ell$ est son déshomogénéisé.
Soit~$I$ l'idéal (homogène)
de~$A[T_0,\ldots, T_n]$ engendré par les~$f_\ell$, et soit~$B$ la~$A$-algèbre
graduée quotient~$A[T_0,\ldots, T_n]/I$.  
On est maintenant exactement dans la situation considérée au~\ref{deshomo-ferme}
ci-dessus : $A[T_0,\ldots, T_n]\to B$ induit une immersion
fermée~$\proj B\hookrightarrow \PP^n_A$ d'image~$V(I)$, 
et l'immersion fermée~$\proj B\times_{\PP^n_A}U_i\hookrightarrow U_i$
est celle définie par l'idéal $J$ (en particulier, on a l'égalité~$V(I)\cap U_i=V(J)$). 

\trois{exemple-homo-deshomo}
{\em Exemples}. Supposons que~$A=\CC$, que~$n=2$ et que~$i=0$. 

\medskip
{\em Déshomogénéisation}. L'idéal homogène
$(T_0T_1-T_2^2+2  T_3^2, T_1^3-iT_0T_1T_2)$ de~$\CC[T_0,T_1,T_2]$ 
induit une immersion fermée~$X\hookrightarrow \PP^2_{\CC}$, 
et~$X\times_{\PP^2_{\CC}}U_0\hookrightarrow U_0$ est l'immersion fermée 
définie par l'idéal~$(\tau_1-\tau_2^2+2\tau_3^2, \tau_1^3-i\tau_1 \tau_2)$ de~$\CC[\tau_1,\tau_2]$
(on a déshomogénéisé les équations de~$X$). 

{\em Homogénéisation}. L'idéal~$(\tau_1^2-3\tau_1+\tau_2^4, \tau_1^3-\tau_2+7)$ de~$\CC[\tau_1,\tau_2]$
induit une immersion fermée~$Y\hookrightarrow U_0$. Si~$Z\hookrightarrow \PP^2_{\CC}$ 
désigne l'immersion fermée induite par l'idéal homogène~$(T_0^2T_1^2-3T_0^3T_1+T_2^4, T_1^3-T_0^2T_2+7T_0^3)$
de~$\CC[T_0,T_1,T_2]$ on a un isomorphisme de~$U_0$-schémas
$Z\times_{\PP^2_{\CC}}U_0 \simeq Y$ (on a homogénéisé les équations de~$Y$).

\deux{imm-ferm-pna}
{\bf Proposition.}
{\em Toute immersion fermée de but~$\PP^n_A$ 
est de la forme~$\proj A[T_0,\ldots, T_n]/I\hookrightarrow \PP^n_A$ pour un 
certain idéal homogène~$I$ de~$A[T_0,\ldots, T_n]$.}

\medskip
{\em Démonstration}. Soit~$\theta$ une immersion fermée de but~$\PP^n_A$. 
Pour tout~$i$, le morphisme~$\theta^{-1} (U_i)\hookrightarrow U_i$ est une immersion fermée, 
associée à un idéal~$I_i$ de~$A\left[\frac{T_j}{T_i}\right]_{j\neq i}$. Pour tout couple~$(i,j)$, les
idéaux~de~$(A[T_0, \ldots, T_n]_{(T_iT_j)})_0$ engendrés par~$I_i$ et~$I_j$
coïncident (ils sont tous deux égaux à l'idéal définissant l'immersion~$\theta^{-1}(U_i\cap U_j)\hookrightarrow U_i\cap U_j$). 
Dans ce qui suit, toutes les fractions que nous allons manipuler seront 
considérées comme vivant dans~$A\left[T_0,\ldots, T_n, \frac 1 {\prod T_i}\right]$ 
({\em cf.}~\ref{ann-localis-pna}). 

\medskip
Il s'agit de montrer l'existence d'un idéal homogène~$I$ de~$A[T_0,\ldots, T_n]$ tel que pour tout~$i$, 
l'idéal~$I_i$ soit l'ensemble des éléments de la forme~$\frac {f}{T_i^d}$ où~$d\in \NN$ et où~$f$ est un élément
de~$I$ homogène de degré~$d$. 
Soit~$I$ l'idéal homogène de~$A[T_0,\ldots, T_n]$ tel que pour tout entier~$d$
on ait
$$I\cap A[T_0,\ldots,  T_n]_d=\left\{f \in  A[T_0,\ldots,  T_n]_d, \;\;\frac f {T_i^d}\in T_i\;\forall i\right\}\;;$$
nous allons
montrer qu'il répond au problème posé. 

\medskip
Compte-tenu de la définition de même de~$I$, il suffit de prouver l'assertion qui suit : {\em pour
tout
$i$ et tout
élément~$\alpha$ de~$I_i$, il existe un entier~$d$ et un élément~$f$ de~$I$, homogène de degré~$d$, tel 
que~$\alpha =\frac f {T_i^d}$.}

\medskip
Soit donc~$i\in \{0,\ldots, n\}$ et~$\alpha \in I_i$. On peut écrire~$\alpha$ sous la forme~$\frac g {T_i^\delta}$ pour un certain entier~$\delta$
et un certain polynôme~$g$ homogène de degré~$\delta$. 

Soit~$j\in  \{0,\ldots, n\}$ (ce qui suit est trivial si~$j=i$, mais nous n'avons pas
de raison d'exclure ce cas). Vu comme appartenant à~$(A[T_0,\ldots, T_n]_{T_iT_j})_0$, l'élément~$\alpha$ appartient
à l'idéal engendré par~$I_j$ ; cela signifie qu'il est de la forme~$\beta_j \frac {T_j^{r_j}}{T_i^{r_j}}$ pour un certain entier~$r_j$ et un certain~$\beta_j \in I_j$. 
On peut toujours si besoin est agrandir~$r_j$ : en effet, on a pour tout~$s\geq 0$ l'égalité
$$\beta_j \frac {T_j^{r_j}}{T_i^{r_j}}=\left(\beta_j \frac {T_i^s}{T_j^s}\right)\cdot \frac{T_j^{r_j+s}}{T_i^{r_j+s}},$$
et~$\beta_j  \frac {T_i^s}{T_j^s}$ appartient évidemment encore à l'idéal~$I_j$. 
Il s'ensuit (comme il n'y a qu'un nombre fini d'indices)
qu'il existe~$r$ et, pour tout~$j$, un élément~$\beta_j\in I_j$ tels que~$\alpha=\beta_j \frac {T_j^r}{T_i^r}$
pour tout~$j$. On peut écrire chacun des~$\beta_j$ sous la forme~$\frac {h_j}{T_j^s}$ où~$s$ est un entier
et~$h_j$ un polynôme homogène de degré~$s$ (il est clair que~$s$ peut être choisi indépendamment de~$j$ : au besoin,
il n'y a qu'à multiplier le numérateur et le dénominateur par une même puissance de~$T_j$). 

\medskip
On a finalement pour tout~$j$ l'égalité~$$\frac g {T_i^\delta}=\frac {h_j} {T_j^s}\cdot \frac {T_j^r}{T_i^r},$$
et donc~$$\frac {gT_i^r}{T_j^{\delta+r}}=\frac{h_j T_i^\delta}{T_j^{s+\delta}}\in I_j,$$
ce qui montre que~$gT_i^r \in I$. On conclut en remarquant que
$$\alpha =\frac g {T_i^\delta}=\frac {gT_i^r}{T_i^{\delta+r}}.\;\;\Box$$

\deux{esp-proj-corps}
{\bf Quelques propriétés de l'espace projectif sur un corps}. 
Soit~$k$ un corps et soit~$n$ un entier. 

\trois{pnk-integre-reduit}
L'anneau~$k[T_0,\ldots, T_n]$ est intègre, et~$k[T_0,\ldots, T_n]_+$ est non nul (il contient
au moins~$T_0$). Il découle alors de~\ref{cs-projbreduit} 
et~\ref{cs-projbintegre}
que~$\PP^n_k$ est irréductible et réduit. Ses ouverts non vides sont donc denses -- c'est 
en particulier le cas des~$D(T_i)$.

\trois{dim-pnk-n}
Comme~$\PP^n_k$ est irréductible, sa dimension de Krull
est égale en vertu de~\ref{krull-cas-irred}
à celle de n'importe quel de ses ouverts non vide. Puisque~$D(T_0)\simeq \Aff^n_k$, 
la dimension de Krull de~$\PP^n_k$ est égale à~$n$. 

\trois{corps-fonct-pnk}
Le corps des fonctions de~$\PP^n_k$ (\ref{krull-cas-irred})
est quant à lui égal au corps
des fractions de l'anneau des fonctions de n'importe lequel de ses ouverts affines non vides. 
Il est en particulier égal à
$${\rm Frac}\;\sch O_{\PP^n_k}(D(T_0))={\rm Frac}\;k\left[\frac{T_j}{T_0}\right]_{j\neq 0}=k\left(\frac{T_j}{T_0}\right)_{j\neq 0}.$$
On vérifie aussitôt que
ce corps peut se décrire indépendamment du choix d'une carte affine :  c'est l'ensemble
des éléments de~$k(T_0,\ldots, T_n)$ qui admettent une écriture
comme quotient de deux polynômes homogènes
de même degré. 

\trois{points-fermes-pnk}
Soit~$U$ un ouvert de~$\PP^n_k$ et soit~$x\in U$. On déduit de~\ref{points-fermes-xtf}
et~\ref{points-fermes-attention}
que les assertions suivantes sont équivalentes : 

\medskip
i) $x$ est fermé dans~$\PP^n_k$ ; 

ii) $x$ est fermé dans~$U$ ; 

iii) $\kappa(x)$ est une extension finie de~$k$. 

\trois{exoirred-deshomo}
{\em Exercice}. Fixons~$i$ entre~$0$ et~$n$, et posons~$\tau_j=\frac{T_j}{T_i}$ pour tout
indice~$j\neq i$. 
Montrez que les opérations d'homogénéisation
et de déshomogénéisation relatives à la variable~$T_i$ mettent en bijection l'ensemble
des (classes d'équivalence de)
polynômes irréductibles de~$k[\tau_j]_{j\neq i}$ et l'ensemble des (classes d'équivalence de)
polynômes irréductibles de~$k[T_0,\ldots, T_n]$ différents de~$T_i$. 

\trois{topologie-p1k}
{\em On suppose que~$n=1$}. La droite projective~$\PP^1_k$ est réunion de ses deux ouverts
denses~$D(T_0)$ et~$D(T_1)$. Chacun 
d'eux est isomorphe à~$\Aff^1_k$, ce qui implique qu'il 
n'est constitué que de points fermés et d'un point générique -- qui est nécessairement
par densité le point générique de~$\PP^1_k$ -- 
et que ses fermés stricts sont précisément les ensembles finis de points fermés. 

\medskip
On en déduit aisément que~$\PP^1_k$ est elle-même constituée d'un point générique
et de points fermés, et que ses fermés stricts sont précisément les ensembles finis de points fermés.

\section{Le foncteur des points de~$\PP^n_A$}
\markboth{Schémas projectifs}{Le foncteur des points de~$\PP^n_A$}

\deux{notations-points-pn}
Soit~$A$ un anneau. Si~$X$ et~$S$ sont deux~$A$-schémas, nous noterons
comme d'habitude~$X(S)$ l'ensemble~$\hom_{A\text{-}\mathsf{Sch}}(S,X)$.

\deux{annonce-fonct-pna}
On fixe un entier~$n$. Le but de cette section est de donner une description explicite et relativement
maniable du foncteur~$S\mapsto \PP^n_A(S)$ de~$A\text{-}\mathsf{Sch}$ vers~$\ens$. 

Pour tout~$i$, on note~$U_i$ l'ouvert~$D(T_i)$ de~$\PP^n_A$, et~$V_i$ l'ouvert affine~$D(T_i)$
de~$\Aff^{n+1}_A=\spec A[T_0,\ldots, T_n]$. La réunion~$V$ des~$V_i$ est un ouvert de~$\Aff^{n+1}_A$. 

\subsection*{Description partielle du foncteur des points de~$\PP^n_A$ : points
donnés par une famille de fonctions}

\deux{fleche-v-pn}
{\bf Construction d'un morphisme~$V\to \PP^n_A$.}

\trois{fleche-vi-pn}
Soit~$i\in \{0,\ldots, n\}$. 
L'inclusion~$A\left[\frac{T_j}{T_i}\right]_{j\neq i}\hookrightarrow A\left[T_0,\ldots, T_n,\frac 1 {T_i}\right]$ induit un 
morphisme~$\Psi_i \colon V_i\to U_i$. 

\medskip
Soit~$j\in \{0,\ldots, n\}$. L'intersection~$U_i\cap U_j$ est l'ouvert~$D(\frac{T_j}{T_i})$ de~$U_i$ ; son image réciproque
$\Psi_i^{-1}(U_i\cap U_j)$ est l'ouvert~$D(\frac{T_j}{T_i})$ de~$V_i$, qui s'identifie à son ouvert~$D(T_j)$, puisque~$T_i$ est inversible
sur~$U_i$. Autrement dit, 
on a l'égalité~$\Psi_i^{-1}(U_i\cap U_j)=V_i\cap V_j$.
Par construction, le morphisme
de schémas affines $\Psi_i|_{V_i\cap V_j}\colon V_i\cap V_j\to U_i \cap U_j$ 
est induit par l'inclusion $$(A\left[T_0,\ldots, T_n\right]_{T_iT_j})_0\hookrightarrow A\left[T_0,\ldots, T_n,\frac 1{T_iT_j}\right].$$

\medskip
En échangeant~$i$ et~$j$, on voit que~$\Psi_j^{-1}(U_i\cap U_j)=V_i\cap V_j$, et que le morphisme
de schémas affines $\Psi_j|_{V_i\cap V_j}\colon V_i\cap V_j\to U_i \cap U_j$
est induit par l'inclusion $$(A\left[T_0,\ldots, T_n\right]_{T_jT_i})_0\hookrightarrow A\left[T_0,\ldots, T_n,\frac 1{T_jT_i}\right].$$ 
{\em Il coïncide donc avec~$\Psi_i|_{V_i\cap V_j}$}. 

\trois{conclu-fleche-v-pna}
Il découle de ce qui précède
que les morphismes~$\Psi_i$ se recollent en un morphisme~$\Psi \colon V\to \PP^n_A$, qui
possède la propriété suivante : {\em pour tout~$i$, on a~$\Psi^{-1}(U_i)=V_i$, et la flèche~$\Psi|_{V_i}\colon V_i\to U_i$ est induite
 par l'inclusion~$A\left[\frac{T_j}{T_i}\right]_{j\neq i}\hookrightarrow A\left[T_0,\ldots, T_n,\frac 1 {T_i}\right]$}. 
 
 \deux{desc-fonct-v-pna}
 {\bf Description fonctorielle du morphisme~$\Psi$}. 
 
 \trois{fonct-vi-et-v}
 Soit~$S$ un~$A$-schéma.
 L'ensemble~$\Aff^{n+1}_A(S)$ s'identifie canoniquement, par la flèche~$\chi \mapsto (\chi^*T_i)$, 
 à ~$\sch O_S(S)^{n+1}$. Soit~$(f_0,\ldots, f_n)$ appartenant à~$\sch O_S(S)^{n+1}$ et soit~$\chi \colon S \to \Aff^{n+1}_A$
 le morphisme correspondant. Comme~$\chi^*T_i$
 est égal à~$f_i$, on voit que pour tout~$s\in S$, on a équivalence entre~$f_i(s)=0$ et~$T_i(\chi(s))=0$. 
 
 \medskip
 Il s'ensuit que~$\chi$ se factorise par~$V_i$ pour un certain~$i$ (resp. par~$V$) si et seulement si~$f_i$ est inversible
 (resp. si et seulement si pour tout~$s\in S$ l'une au moins des~$f_j$ est inversible en~$s$). 
 
 \medskip
 En d'autres termes, $V_i(S)$ s'identifie au sous-ensemble de~$\sch O_S(S)^{n+1}$ formé des familles~$(f_j)$
 avec~$f_i$ inversible, et~$V(S)$ à celui formé des familles~$(f_j)$ telles
 que~$S=\bigcup D(f_j)$. 
 
 \trois{notation-crochet-pna}
 Le morphisme~$\Psi$ induit un morphisme de foncteurs
 $\chi \mapsto \Psi \circ \chi$ -- qui le caractérise entièrement -- de~$S\mapsto V(S)$ vers~$S\mapsto \PP^n_A(S)$. 
 Si~$S$ est un~$A$-schéma et si~$(f_0,\ldots, f_n)$ est un élément de~$V(S)$, on notera~$[f_0:f_1:\ldots:f_n]$ son image dans~$\PP^n_A(S)$ par ce morphisme. 
 
 \trois{test-appartient-ui}
 Fixons un indice~$i$. Soit~$S$ un~$A$-schéma, soit~$(f_0,\ldots,f_n)\in V(S)$ et soit~$\chi \colon S \to V$ le morphisme
 correspondant. L'image réciproque de~$U_i$ 
 par la flèche $\Psi \circ \chi$ (flèche qui n'est autre que l'élément~$[f_0:\ldots :f_n]$ de~$\PP^n_A(S)$)
 est égale à l'image réciproque 
 de~$\Psi^{-1}(U_i)=V_i$ par~$\chi$, c'est-à-dire à~$D(f_i)$ 
 en vertu de~\ref{fonct-vi-et-v}.

 \medskip
 En particulier, on voit que~$[f_0: \ldots:f_n]$ appartient au sous-ensemble~$U_i(S)$ de~$\PP^n_A(S)$ (constitué des
 morphismes qui se factorisent ensemblistement par~$U_i$) si et seulement si~$S=D(f_i)$, c'est-à-dire 
 si et seulement si~$f_i$ est inversible, ou encore
 si et seulement si~$(f_0,\ldots, f_n)\in V_i(S)$. 
 
 \trois{desc-crochet-ui}
 Fixons~$i$, soit~$S$ un~$A$-schéma et soit~$(f_0,\ldots,f_n)\in V_i(S)$. D'après~\ref{test-appartient-ui}, l'élément~$[f_0:\ldots: f_n]$
 de~$\PP^n_A(S)$ appartient à~$U_i(S)$. 
   
 \medskip
 On a~$U_i=\spec A\left[\frac{T_j}{T_i}\right]_{j\neq i}$. L'application~$\chi \mapsto \left(\chi^*\left(\frac{T_j}{T_i}\right)\right)_j$ permet donc
 d'identifier~$U_i(S)$ à
 l'ensemble des familles~$(g_j)_{0\leq j\leq n, j\neq i}$ d'éléments de~$\sch O_S(S)$. Comme~$\Psi|_{V_i}\colon V_i\to U_i$ est 
 induit par l'inclusion~$A\left[\frac{T_j}{T_i}\right]_{j\neq i}\hookrightarrow A\left[T_0,\ldots, T_n,\frac 1 {T_i}\right]$, 
 l'élément~$[f_0:f_1:\ldots :f_n]$ de~$U_i(S)$ correspond au~$n$-uplet~$(f_j/f_i)_j$ par l'identification ci-dessus. 
On en déduit deux faits importants. 
 
 \medskip
 $\bullet$ L'application~$(f_0,\ldots, f_n)\mapsto [f_0:f_1:\ldots:f_n]$ de~$V_i(S)$ dans~$U_i(S)$ est surjective : 
 en effet, si~$(g_j)_{j\neq i}$ est un élément de~$U_i(S)$, il est par ce qui précède égal à $[g_0:\ldots:g_{i-1}:1:g_{i+1}:\ldots:g_n]$. 
 
 $\bullet$ Si~$(f_0,\ldots, f_n)$ et~$(g_0,\ldots, g_n)$ sont deux éléments de~$V_i(S)$ alors
 $$[f_0:f_1:\ldots :f_n]=[g_0:g_1:\ldots :g_n]$$ si et seulement si $f_j/f_i=g_j/g_i$ pour tout~$j\neq i$. On vérifie aussitôt que cela
 revient à demander qu'il existe~$\lambda \in \sch O_S(S)\ti$ tel que~$f_j=\lambda g_j$ pour tout~$j$.

 \trois{relation-generale-crochet}
 On désigne toujours par~$S$
 un~$A$-schéma, et l'on se donne deux éléments
 $(f_0,\ldots, f_n)$ et~$(g_0,\ldots, g_n)$ de~$\sch O_S(S)^{n+1}$. Le but de ce qui suit est de
 montrer que~$$[f_0:\ldots :f_n]=[g_0:\ldots:g_n]$$ 
 si et seulement si il existe~$\lambda\in \sch O_S(S)\ti$ telle que~$f_i=\lambda g_i$ pour tout~$i$, et qu'une telle~$\lambda$
 est nécessairement unique dans ce cas.  
 
 \medskip
 Supposons qu'il existe une telle~$\lambda$. On a alors immédiatement les égalités~$D(f_i)=D(g_i)$ pour tout~$i$, et sur~$D(f_i)$ la fonction~$\lambda$ est nécessairement 
 égale à~$f_i/g_i$, ce qui montre déjà son unicité puisque les~$D(f_i)$ recouvrent~$S$ par définition de~$V(S)$. Pour cette même raison il suffit, pour montrer
 que les éléments~$[f_0:\ldots f_n]$ et~$[g_0:\ldots:g_n]$ de~$\PP^n_A(S)$ coïncident, de
 prouver que c'est le cas des éléments
 $[f_0|_{S_i}:\ldots:f_n|_{S_i}]$ et~$[g_0|_{S_i}:\ldots:g_n|_{S_i}]$ de~$\PP^n_A(S_i)$ pour tout~$i$, où l'on a posé~$S_i=D(f_i)=D(g_i)$.

Fixons donc~$i$. Les fonctions~$f_i$ et~$g_i$ sont inversibles sur~$S_i$, et les~$S_i$-points
$(f_0|_{S_i},\ldots, f_n|_{S_i})$ et~$(g_0|_{S_i}, \ldots, g_n|_{S_i})$
appartiennent en conséquence à~$V_i(S_i)$. 
En vertu de~\ref{desc-crochet-ui}, 
l'existence de la fonction~$\lambda$ entraîne  
l'égalité
$$[f_0|_{S_i}:\ldots:f_n|_{S_i}]=[g_0|_{S_i}:\ldots:g_n|_{S_i}],$$
qui est ce qu'on voulait. 

\medskip
Réciproquement, supposons que~$[f_0:f_1:\ldots :f_n]=[g_0:g_1:\ldots :g_n]$. Fixons~$i$. 
D'après~\ref{test-appartient-ui},
les ouverts~$D(f_i)$ et~$D(g_i)$ de~$T$ sont tous deux égaux à l'image réciproque de~$U_i$
par~$[f_0:f_1:\ldots :f_n]=[g_0:g_1:\ldots :g_n]$ ; en conséquence, ils coïncident ; posons~$S_i=D(f_i)=D(g_i)$. 
Comme~$f_i$ et~$g_i$ sont inversibles sur~$S_i$, l'égalité
$$[f_0|_{S_i}:\ldots :f_n|_{S_i}]=[g_0|_{S_i}:\ldots :g_n|_{S_i}]$$
implique en vertu de~\ref{desc-crochet-ui}
qu'il existe une fonction inversible~$\lambda_i$ sur~$S_i$ telle que~$f_j|_{S_i}=\lambda_i g_j|_{S_i}$ pour tout~$j$. 
L'assertion d'unicité déjà établie entraîne que~$\lambda_i|_{S_i\cap S_j}=\lambda_j|_{S_i\cap S_j}$
pour tout~$(i,j)$, et les~$\lambda_i$
se recollent ainsi en une fonction inversible~$\lambda$ qui possède la propriété requise. 

\deux{comment-crochet}
Récapitulons : le morphisme~$\Psi \colon V\to \PP^n_A$
induit pour tout~$A$-schéma~$S$
une application~$(f_0,\ldots, f_n)\mapsto [f_0:\ldots:f_n]$ de~$V(S)$
vers~$\PP^n_A(S)$ dont on a décrit au~\ref{relation-generale-crochet}
le~«noyau», c'est-à-dire les conditions 
sous lesquelles deux éléments ont même image : il faut et il suffit qu'ils satisfassent
la relation de~«colinéarité inversible». 

\medskip
Il est par contre difficile en général de décrire son image
que nous noterons~$\PP_A^{n,\sharp}(S)$. 
Indiquons
tout de même quelques faits à son sujet. 

\trois{unionui-image}
Soit~$S$ un~$A$-schéma. Pour tout~$i$, le sous-ensemble
~$U_i(S)$ de~$\PP^n_A(S)$ est contenu dans~$\PP_A^{n,\sharp}(S)$ : 
c'est une simple conséquence de la surjectivité de
l'application~$V_i(S)\to U_i(S)$
(\ref{desc-crochet-ui}). 

\medskip
{\em Supposons que~$\PP^n_A(S)$ soit la réunion des~$U_i(S)$} (nous verrons
un peu plus bas que
cette condition 
est effectivement vérifiée lorsque~$S$ est le spectre d'un corps, et plus généralement d'un anneau local). On a alors en vertu de ce qui précède 
l'égalité~$\PP^{n,\sharp}_A(S)=\PP_A^n(S)$. 
Notez par ailleurs que comme~$V_i(S)$ est l'image réciproque
de~$U_i(S)$ pour tout~$i$ (\ref{test-appartient-ui}), l'égalité~$\PP^n_A(S)=\bigcup U_i(S)$ implique
que~$V(S)=\bigcup V_i(S)$. On dispose donc d'une bijection canonique
entre~$\PP^n_A(S)$ et le quotient
de
$$V(S)=\bigcup V_i(S)=\{(f_0,\ldots, f_n)\in \sch O_S(S)^{n+1},\;\exists i,\;f_i\in \sch O_S(S)\ti\}$$
par la relation de colinéarité inversible. 

\trois{pdiese-pastout}
En général, $\PP_A^{n,\sharp}(S)\subsetneq \PP^n_A(S)$. Donnons un exemple. On suppose que~$A\neq\{0\}$ et que~$n\geq 1$, et l'on pose
$S=\PP^n_A$. Comme~$\sch O_{\PP^n_A}(\PP^n_A)=A$ d'après~\ref{fonctions-pna}, le sous-ensemble
$\PP^{n,\sharp}_A(\PP^n_A)$ est simplement constitué d'éléments de la forme~$[a_0:\ldots:a_n]$ où les~$a_i$ 
appartiennent à~$A$, et où les~$D(a_i)$ recouvrent~$\spec A$. 

\medskip
Or {\em ${\rm Id}_{\PP^n_A}$ ne peut pas être de cette forme}.   En effet, soit~$(a_0,\ldots, a_n)$ comme ci-dessus, 
et soit~$x$ un point de~$\spec A$ (comme~$A$ est non nul, son spectre est non vide). Il existe~$i$ tel que~$a_i(x)\neq 0$, 
et~$a_i$ est donc inversible sur toute la fibre~$\PP ^n_{\kappa(x)}$ de~$\PP^n_A$ en~$x$. Il s'ensuit que le morphisme
$[a_0:\ldots:a_n]$ envoie toute la fibre~$\PP ^n_{\kappa(x)}$ sur l'ouvert~$D(T_i)$ de~$\PP ^n_{\kappa(x)}$, qui est strict
car~$n\geq 1$ (\ref{bord-carte-aff-pn}). En conséquence, $[a_0:\ldots:a_n]$ ne peut être égal à~${\rm Id}_{\PP^n_A}$. 

\deux{pndiese-corps-local}
Nous allons toutefois donner deux exemples fondamentaux dans lesquels~$\PP_A^{n,\sharp}(S)$ est égal à~$\PP^n_A(S)$. Le premier
d'entre eux s'avérera être un cas particulier du second, mais nous avons choisi de le traiter séparément
au vu de son importance. 

\trois{pndiese-corps}
Soit~$k$ une~$A$-algèbre qui est un corps. Comme~$\spec k$ ne comprend
qu'un point, tout morphisme de~$\spec k$ vers~$\PP^n_A$ a nécessairement une image contenue dans~$U_i$ pour un certain~$i$. 
Il s'ensuit que~$\PP^n_A(k)=\bigcup U_i(k)$. On déduit alors de~\ref{unionui-image}
que~$\PP^{n,\sharp}_A(k)=\PP^n_A(k)$, et plus précisément que~$\PP^n_A(k)$ s'identifie 
naturellement au quotient
de
$k^{n+1}\setminus\{(0,\ldots, 0)\}$ par la relation de colinéarité inversible ; on retrouve ainsi la description classique
ou naïve de l'espace projectif. 

\trois{pndiese-anneaulocal-lemme}
{\bf Lemme.}
{\em Soit~$B$ un anneau local. Si~$W$ est un ouvert de~$\spec B$ qui contient son unique
point fermé alors~$W=\spec B$}. 

\medskip
{\em Démonstration}. 
Soit~$x$ le point fermé de~$\spec B$. Comme~$W$ contient~$x$, il existe
$f\in B$ tel que~$D(f)\subset W$ et tel que~$f(x)\neq 0$. Mais cette dernière condition signifie que~$f$
n'appartient pas à l'idéal maximal de~$B$, et donc que~$f$ est inversible. En conséquence~$D(f)$ est égal à~$\spec B$ tout entier, et
il en va {\em a fortiori}
de même de~$W$.~$\Box$

\trois{pndiese-anneaulocal}
Soit maintenant~$B$ une~$A$-algèbre locale et soit~$x$ le point fermé de~$\spec B$. Soit~$\chi$
un morphisme de~$\spec B$ vers~$\PP^n_A$. Il existe
un indice~$i$ tel que~$\chi(x)\in U_i$. L'image réciproque~$\chi^{-1}(U_i)$
est donc un ouvert de~$\spec B$ qui contient~$x$ ; d'après le lemme~\ref{pndiese-anneaulocal-lemme}
ci-dessus, c'est~$\spec B$ tout entier, ce qui veut dire que~$\chi(\spec B)\subset U_i$. Il s'ensuit que~$\PP^n_A(B)=\bigcup U_i(B)$.
On déduit alors de~\ref{unionui-image}
que~$\PP^{n,\sharp}_A(B)=\PP^n_A(B)$, et plus précisément que~$\PP^n_A(B)$ s'identifie 
naturellement au quotient
de
$$\{(b_0,\ldots, b_n)\in B^{n+1},\;\exists i,\; b_i\in B\ti\}$$ 
par la relation de colinéarité inversible. 

\subsection*{Quelques exemples}

\deux{points-pn-corps}
On suppose pour ce paragraphe que l'anneau~$A$ est un corps, que nous préférons noter~$k$. Soit~$x\in \PP^n_k(k)$
(on peut voir~$x$ aussi bien comme un morphisme de~$\spec k$ vers~$\PP^n_k$ que comme
un point schématique de~$\PP^n_k$ de corps résiduel~$k$, {\em cf.}~\ref{point-x-krat} ; dans ce qui suit, nous utiliserons implicitement
ces deux interprétations).
D'après~\ref{pndiese-corps}, 
il existe un~$(n+1)$-uplet~$(a_0,\ldots, a_n)$ d'éléments {\em non tous nuls}
de~$k$ tel que~$x=[a_0:\ldots:a_n]$. Pour tout~$i$, on a~$x\in U_i(k)$ si et seulement si~$a_i\neq 0$
(\ref{test-appartient-ui}).

\trois{coordonnees-ui}
Soit~$i$ tel que~$a_i\neq 0$, c'est-à-dire encore tel que~$x\in U_i(k)$. 
On a un isomorphisme~$U_i\simeq \spec k\left[\frac{T_j}{T_i}\right]_{j\neq i}$ qui permet d'identifier
$U_i(k)$ à l'ensemble des $n$-uplets~$(b_j)_{0\leq j\leq n, j\neq i}$ d'éléments de~$k$. D'après~\ref{desc-crochet-ui},
le~$n$-uplet qui correspond à~$x$ 
est~$(a_j/a_i)_{j\neq i}$. 

\trois{coordonnees-ui-p1}
Déclinons ces faits
dans le cas particulier où~$n=1$. La droite projective~$\PP^1_k$ est alors
réunion de deux cartes affines~$U_0$ et~$U_1$, respectivement munies des fonctions coordonnées~$\tau_1:=T_1/T_0$ et~$\tau_0:=T_0/T_1$. 
Le point~$x$ s'écrit~$[a_0:a_1]$, où~$(a_0,a_1)\in k^2\setminus\{(0,0)\}$. 

\medskip
Si~$a_0\neq 0$ alors~$x\in U_0(k)$ ; c'est  
le point naïf d'équation~$\tau_1=a_1/a_0$
de la carte~$U_0$. 

Si~$a_1\neq 0$ alors~$x\in U_1(k)$ ; c'est le point naïf d'équation~$\tau_0=a_0/a_1$ de
la carte affine~$U_1$. 

Si~$a_0=0$ alors~$x\notin U_0(k)$. C'est le point naïf d'équation~$\tau_0=0$ de la carte affine~$U_1$ ou, si l'on préfère, 
le «point à l'infini» relativement à la coordonnée~$\tau_1$. 

Si~$a_1=0$ alors~$x\notin U_1(k)$. C'est le point naïf d'équation~$\tau_1=0$ de la carte affine~$U_0$ ou, si l'on préfère, 
le «point à l'infini» relativement à la coordonnée~$\tau_0$. 

\medskip
Notons que le point fermé~$x$ de~$\PP^1_k$
est égal à~$V(a_1T_0-a_0T_1)$ : il suffit en effet de vérifier
que~$\{x\}\cap U_0=V(a_1-a_0\tau_1)$ et~$\{x\}\cap U_1=V(a_1\tau_0-a_0)$, et cela découle
aussitôt de ce qui précède. 

\deux{etude-section-general}
On ne suppose
plus que~$A$ est un corps. Le~$A$-schéma~$\PP^1_A$ est réunion de deux
cartes affines~$U_0$ et~$U_1$, respectivement munies des fonctions coordonnées~$\tau_1:=T_1/T_0$ et~$\tau_0:=T_0/T_1$. 

\trois{section-p1a}
Soit~$s\in \PP^{1,\sharp}_A(A)$. Par définition, $s$ est une section du morphisme structural~$\PP^1_A\to \spec A$, 
section qui est de la forme
$[a_0:a_1]$ où~$a_0$ et~$a_1$ sont deux éléments de~$A$ tels que~$\spec A=D(a_0)\cup D(a_1)$, c'est-à-dire encore tels que l'idéal
$(a_0,a_1)$ de~$A$ soit égal à~$A$. 

\medskip
Soit~$x\in \spec A$. L'image~$s(x)$ est alors
(par fonctorialité de toutes les constructions) le~$\kappa(x)$-point $[a_0(x), a_1(x)]$ de la fibre~$\PP^1_{\kappa(x)}$ de~$\PP^1_A$ en~$x$
(comme~$\spec A=D(a_0)\cup D(a_1)$ on a bien~$a_0(x)\neq 0$ ou~$a_1(x)\neq 0$). 
On a donc
$$s(\spec A)\cap \PP^1_{\kappa(x)}=s(x)=\underbrace{V(a_1(x)T_0-a_0(x)T_1)}_{\text{fermé~de~}\;\PP^1_{\kappa(x)}}=\underbrace{V(a_1T_0-a_0T_1)}
_{\text{fermé~de}\;\PP^1_A}\cap \PP^1_{\kappa(x)}.$$ 

\medskip
Ceci valant pour tout~$x$, il vient~$s(\spec A)=V(a_1T_0-a_0T_1)$.

\trois{section-21}
Supposons que~$A=\ZZ$, que~$a_0=2$ et que~$a_1=1$. Soit~$y\in \spec \ZZ$. Son image~$s(y)$ 
est le~$\kappa(y)$-point~$[1:2]$ de~$\PP^1_{\kappa(x)}$. On distingue maintenant deux cas. 

\medskip
$\bullet$ Si~$2(y)\neq 0$, c'est-à-dire
si~$y\neq x_2$,  alors~$s(y)$ est le~$\kappa(y)$-point d'équation~$\tau_1=(1/2)$ de la carte
affine~$U_0\cap \PP^1_{\kappa(y)}$. 

$\bullet$ Si~$2(y)=0$, c'est-à-dire
si~$y=x_2$, alors~$s(y)$ est l'origine de la carte affine~$U_1\cap \PP^1_{\FF_2}$, c'est-à-dire
encore le point à l'infini de~$\PP^1_{\FF_2}$ relatif à la coordonnée~$\tau_1$. 

\medskip
L'image~$s(\spec \ZZ)$ est le fermé~$V(2T_1-T_0)$ de~$\PP^1_{\ZZ}$, qui est irréductible. En effet, $\spec \ZZ$
est irréductible et il en va donc de même de son image par n'importe quelle application continue -- on peut aussi si l'on préfère
remarquer qu'il existe un
homéomorphisme naturel
$$V(2T_1-T_0)\simeq \proj \ZZ[T_0,T_1]/(2T_1-T_0)=\proj \ZZ[T_1]\simeq \spec \ZZ,$$ et utiliser là encore l'irréductibilité de~$\spec \ZZ$. 

\medskip
{\em Remarque}. La description explicite de~$s(y)$ donnée ci-dessus pour tout~$y$
appartenant à~$\spec \ZZ$
montre que~$V(2T_1-T_0)$ est la réunion disjointe de son ouvert non vide~$U_0\cap V(2T_1-T_0)$
(qui est le fermé~$V(2\tau_1-1)$ de~$U_0\simeq \Aff^1_{\ZZ}$) et du point à l'infini (relativement
à la coordonnée~$\tau_1$) de la fibre~$\PP^1_{\FF_2}$. Comme
le fermé~$V(2T_1-T_0)$ est irréductible,
son ouvert non vide~$V(2\tau_1-1)\subset U_0$ en est une partie dense, ce qui veut dire que
$V(2T_1-T_0)=\overline{V(2\tau_1-1)}$. 

\medskip
Or le fermé~$V(2\tau_1-1)$ de~$U_0\simeq \Aff^1_{\ZZ}$ a déjà été étudié 
en~\ref{zt-demi}. Nous avions signalé qu'il ne rencontrait pas~$\Aff^1_{\FF_2}$, 
et mentionné en~\ref{zt-geom}
qu'il l'intersectait en fait moralement «à l'infini». Cette assertion
un peu vague a désormais sa traduction rigoureuse :  nous venons
en effet de voir que dans~$\PP^1_{\ZZ}$, l'adhérence 
de~$V(2\tau_1-1)$ est précisément la réunion de ce dernier et du point à l'infini de~$\PP^1_{\FF_2}$.  

\trois{section-23}
Supposons que~$A=\ZZ$, que~$a_0=2$ et que~$a_1=3$ (notons que~$(2,3)$ engendre bien~$\ZZ$
puisque~$2$ et~$3$ sont premiers entre eux). Soit~$y\in \spec \ZZ$. Son image~$s(y)$ 
est le~$\kappa(y)$-point~$[2:3]$ de~$\PP^1_{\kappa(y)}$. On distingue maintenant trois cas.  

\medskip
$\bullet$ Si~$2(y)\neq 0$ et~$3(y)\neq 0$  c'est-à-dire
si~$y\notin\{x_2,x_3\}$ est le~$\kappa(y)$-point d'équation~$\tau_1=(3/2)$ de la carte
affine~$U_0\cap \PP^1_{\kappa(y)}$, et le~$\kappa(y)$-point d'équation~$\tau_0=(2/3)$
de la carte affine~$U_1\cap \PP^1_{\kappa(y)}$.

$\bullet$ Si~$2(y)=0$, c'est-à-dire
si~$y=x_2$ alors~$s(y)$ est l'origine de la carte affine~$U_1\cap \PP^1_{\FF_2}$, c'est-à-dire
encore le point à l'infini de~$\PP^1_{\FF_2}$ relatif à la coordonnée~$\tau_1$. 

$\bullet$ Si~$3(y)=0$, c'est-à-dire
si~$y=x_2$ alors~$s(y)$ est l'origine de la carte affine~$U_0\cap \PP^1_{\FF_3}$, c'est-à-dire
encore le point à l'infini de~$\PP^1_{\FF_3}$ relatif à la coordonnée~$\tau_0$. 

\medskip
L'image~$s(\spec \ZZ)$ est le fermé~$V(2T_1-3T_0)$ de~$\PP^1_{\ZZ}$, qui est irréductible
puisque~$\spec \ZZ$ est irréductible. 

\subsection*{Les faisceaux~$\sch O(d)$}

\deux{2intro-foncteur-points}
Si~$k$ est une~$A$-algèbre qui est un corps, on a vu que~$\PP^n_A(k)$ possède une description agréable : 
c'est le quotient de~$k^{n+1}\setminus \{(0,\ldots, 0)\}$ par la relation de colinéarité inversible. 

\trois{desc-pnk-pasbon}
Malheureusement, cette description ne se généralise pas telle quelle aux~$A$-schémas 
quelconques. On a certes construit, pour tout~$A$-schéma~$S$, 
un sous-ensemble naturel~$\PP^{n,\sharp}_A(S)$ de~$\PP^n_A(S)$ qui s'identifie au quotient
de
$$V(S)=\{(f_0,\ldots, f_n)\in \sch O_S(S)^{n+1}, \bigcup D(f_i)=S\}$$ par la relation de colinéarité 
inversible. Mais on a signalé qu'en général, $\PP^{n,\sharp}_A(S)$ est un sous-ensemble {\em strict}
de~$\PP^n_A(S)$ (\ref{pdiese-pastout}). 

\trois{desc-pnk-bon}
Nous allons maintenant donner une description légèrement différente de~$\PP^n_A(k)$
qui aura l'avantage de bien se généraliser à un~$A$-schéma quelconque. 

\medskip
Soit~$\mathsf C$ la catégorie définie comme suit. Ses objets sont les familles~$(L, s_0,\ldots, s_n)$ 
où~$L$ est un~$k$-espace vectoriel de dimension~$1$, et où les~$s_i$ sont des éléments
non tous nuls de~$L$. Un morphisme de~$(L, (s_i))$ vers~$(L', (s'_i))$ est une {\em bijection}
linéaire~$\phi \colon L\to L'$ qui envoie~$s_i$ sur~$s'_i$ pour tout~$i$ (les morphismes de~$\mathsf C$
sont donc tous des isomorphismes). 

Nous allons montrer qu'il existe une bijection naturelle entre~$\PP^n_A(k)$ et l'ensemble~$\mathsf I$
des classes d'isomorphie
d'objets de~$\mathsf C$. 

\medskip
Soit~$x\in \PP^n_A(k)$. Écrivons~$x=[a_0:\ldots:a_n]$ où les~$a_i$ sont des scalaires non tous nuls. 
La classe d'isomorphie de l'objet~$(k, (a_i))$ de~$\mathsf C$ ne dépend alors que de~$x$, et pas du choix des~$a_i$. En effet, 
si~$\lambda\in k\ti$, l'homothétie de rapport~$\lambda$ est un isomorphisme de~$(k, (a_i))$ sur~$(k, (\lambda a_i))$. On a ainsi
défini une application de~$\PP^n_A(k)$ vers~$\mathsf I$. 

\medskip
Réciproquement, soit~$(L, (s_i))$ un objet de~$\mathsf C$. Choisissons une base de~$L$, c'est-à-dire un isomorphisme~$L\simeq k$ ; pour tout~$i$,
notons~$a_i$ l'image de~$s_i$ sous cette bijection. Comme deux isomorphismes entre~$L$ et~$k$ «diffèrent»
simplement d'une homothétie de rapport inversible, l'élément~$[a_0:\ldots:a_n]$ de~$\PP^n_A(k)$ ne dépend que de~$(L, (s_i))$, et même 
que de la classe d'isomorphie de ce dernier. On a ainsi construit une application de~$\mathsf I$ vers~$\PP^n_A(k)$. 

\medskip
On vérifie aussitôt que les deux flèches~$\PP^n_A(k)\to \mathsf I$ et~$\mathsf I\to \PP^n_A(k)$ ainsi construites sont des bijections
réciproques l'une de l'autre. 

\deux{intro-def-od}
Pour pouvoir généraliser la description de~$\PP^n_A(k)$ donnée ci-dessus à un~$A$-schéma 
quelconque, il va être nécessaire d'introduire un~$\sch O_{\PP^n_A}$-module localement libre de rang~$1$ particulier, qui sera noté~$\sch O(1)$. 
Nous allons en fait plus généralement définir pour tout~$d\in \ZZ$ un~$\sch O_{\PP^n_A}$-module~$\sch O(d)$. 

\trois{def-od}
Soit~$d\in \ZZ$. Rappelons que si~$U$ est un ouvert de~$\PP^n_A$, on note~$S\homog(U)$ l'ensemble
des polynômes homogènes~$f\in A[T_0,\ldots, T_n]$ tels que~$U\subset D(f)$. 
La flèche
$$U\mapsto (S\homog(U)^{-1}A[T_0,\ldots, T_n])_d$$ est un préfaisceau sur~$\PP^n_A$, 
qui est de manière naturelle un module sur le préfaisceau d'anneaux
$U\mapsto (S\homog(U)^{-1}A[T_0,\ldots, T_n])_0.$

\medskip
On note~$\sch O(d)$ le faisceau associé à~$U\mapsto (S\homog(U)^{-1}A[T_0,\ldots, T_n])_d$. C'est de manière
naturelle un module sur le faisceau d'anneaux~$\sch O(0)$ associé à~$(S\homog(U)^{-1}A[T_0,\ldots, T_n])_0$,
qui n'est autre que~$\sch O_{\PP^n_A}$ par définition de ce dernier. 

\trois{od-ui}
Soit~$i\in \{0,\ldots, n\}$. Pour tout ouvert~$U\subset U_i=D(T_i)$, l'application
$f\mapsto T_i^df$ induit une bijection de~$S\homog(U)^{-1}A[T_0,\ldots, T_n])_0$ vers~$S\homog(U)^{-1}A[T_0,\ldots, T_n])_d$,
de réciproque
$g\mapsto T_i^{-d}g$. 

Il s'ensuit que~$f\mapsto T_i^df$ induit un isomorphisme~$\sch O_{U_i}\simeq \sch O(d)|_{U_i}$, 
de réciproque
$g\mapsto T_i^{-d}g$. 

\medskip
Comme les~$U_i$ recouvrent~$\PP^n_A$, le~$\sch O_{\PP^n_A}$-module~$\sch O(d)$ est localement libre de rang~$1$.

\trois{sec-globales-od}
Soit~$i\in \{0,\ldots, n\}$. Il découle de~\ref{od-ui}
et du fait que
$$\sch O_{U_i}(U_i)=A\left[T_0,\ldots, T_n, \frac 1 {T_i}\right]_0$$
que~$\sch O(d)(U_i)$ s'identifie naturellement à~$A\left[T_0,\ldots, T_n, \frac 1 {T_i}\right]_d$. 

\medskip
On en déduit, par des raisonnements et calculs analogues à ceux suivis en~\ref{pna-copies-ana}
et~\ref{fonctions-pna}, que la flèche canonique~$A[T_0,\ldots, T_n]_d\to \sch O(d)(\PP^n_A)$
est un isomorphisme. 
Énonçons tout de suite quelques conséquences de ce fait. 

\medskip
$\bullet$ {\em Supposons~$d<0$}.  Le~$A$-module
$\sch O(d)(\PP^n_A)$ s'identifie à l'ensemble des polynômes en~$T_0,\ldots, T_n$ qui sont
homogènes de degré~$d$ ; il est donc nul. Comme~$\sch O_{\PP^n_A}(\PP^n_A)=A$, 
on voit que~$\sch O(d)$ est non trivial (comme~$\sch O_{\PP^n_A}$-module localement libre de rang~$1$)
dès que~$A\neq\{0\}$. 

$\bullet$ {\em Supposons~$d\geq 0$}.  Le~$A$-module~$\sch O(d)(\PP^n_A)$ s'identifie à l'ensemble des polynômes en~$T_0,\ldots, T_n$ qui sont
homogènes de degré~$d$ ; il est donc libre de rang~$r(n,d):=\left(\begin{array}cn+d\\d\end{array}\right)$ (faites l'exercice !). Si~$n\geq 1$ et~$d>0$
on vérifie aussitôt que~$r(d,n)>r(0,n)=1$ ; il s'ensuit que~$\sch O(d)$ 
est non trivial 
dès que~$n\geq 1, d>0$ et~$A\neq\{0\}$. 

\trois{vsec-vhomog-compat}
Supposons que~$d\geq 0$ et soit~$f\in A[T_0,\ldots, T_n]_d$. On peut en vertu de~\ref{sec-globales-od}
la voir comme un élément de~$\sch O(d)(\PP^n_A)$. Soit~$F$ le lieu des zéros
de~$f$
{\em vue comme section globale de~$\sch O(d)$}. Soit~$i\in \{0,\ldots n\}$. On a vu au~\ref{od-ui}
que~$g\mapsto T_i^{-d}g$ établit un isomorphisme entre~$\sch O(d)|_{U_i}$ et~$\sch O_{U_i}$. 
Cet isomorphisme envoie~$f$ sur la fonction~$\frac f {T_i^d}$. En conséquence, $F\cap U_i$ est le lieu des
zéros de la fonction~$\frac f {T_i^d}$, lequel s'identifie à~$V(f)\cap U_i$ (d'après~\ref{df-bf0-homeo}). 

\medskip
Ceci valant pour tout~$i$, il vient~$F=V(f)$.  Le lieu des zéros~$V(f)$ de~$f$ au sens {\em ad hoc}
de la géométrie projective est donc en fait son lieu des zéros comme section de~$\sch O(d)$. 

\medskip
Soit~$g\in A[T_0,\ldots, T_n]_d$. Le quotient~$\frac f g$ est une fonction bien définie
sur l'ouvert affine~$D(g)$. Comme on a évidemment $f =\frac f g \cdot g$ dans
le~$\sch O_{\PP^n_A}(D(g))$-module~$\sch O(d)(D(g))$, on voit que~$\frac f g$ est 
bien le quotient des sections~$f$ et~$g$ de~$\sch O(d)$, au sens de~\ref{recap-sectinv}.  

\medskip
{\em Nous utiliserons
ces faits implicitement dans toute la suite du texte}. 

\trois{od1-od2}
Soient~$d_1$ et~$d_2$ deux entiers relatifs. Pour tout ouvert~$U$ de~$\PP^n_A$, le produit
définit une application bilinéaire
de
$$S\homog(U)^{-1}A[T_0,\ldots, T_n]_{d_1}\times S\homog(U)^{-1}A[T_0,\ldots, T_n]_{d_2}$$
vers~$S\homog(U)^{-1}A[T_0,\ldots, T_n]_{d_1+d_2},$
d'où par passage au produit tensoriel et faisceautisation un morphisme

$$\sch O (d_1)\otimes_{\sch O_{\PP^n_A}}\sch O(d_2)\to \sch O(d_1+d_2).$$

Nous allons montrer qu'il s'agit d'un {\em isomorphisme}. Pour cela, il suffit de raisonner localement ; on peut donc fixer~$i\in \{0,\ldots, n\}$
et établir l'assertion requise sur la carte~$U_i$. On déduit de~\ref{od-ui}
que~$T_i^{d_1}$ (resp.~$T_i^{d_2}$) est une section inversible de~$\sch O(d_1)|_{U_i}$ (resp. de~$\sch O(d_2)|_{U_i}$). En conséquence, 
$T_i^{d_1}\otimes T_i^{d_2}$ est une section inversible de~$(\sch O(d_1)\otimes_{\sch O_{\PP^n_A}}\sch O(d_2))|_{U_i}$. 

Le morphisme ci-dessus envoie~$T_i^{d_1}\otimes T_i^{d_2}$
sur~$T_i^{d_1+d_2}$ qui est elle-même d'après~{\em loc. cit.}
une section inversible de~$\sch O(d_1+d_2)|_{U_i}$. L'assertion requise s'ensuit aussitôt. 

\subsection*{Description complète du foncteur des points $\PP^n_A$}

\deux{def-lt}
Soit~$S$ un schéma. On note~$\mathsf L_S$ la catégorie définie comme suit. Ses objets sont les familles
$(\sch L, s_0,\ldots, s_n)$ où~$\sch L$ est un~$\sch O_S$-module localement libre de rang~$1$ et où les~$s_i$
sont des sections globales de~$\sch L$ telles que~$S=\bigcup D(s_i)$. Si~$(\sch L, (s_i))$ et~$(\sch L', (s'_i))$ sont
deux objets de~$\mathsf L_S$, un morphisme de~$(\sch L, (s_i))$ vers~$(\sch L', (s'_i))$ est un isomorphisme de~$\sch L$
sur~$\sch L'$
qui envoie~$s_i$ sur~$s'_i$ pour tout~$i$ (ainsi, tout morphisme de~$\mathsf L_S$ est un isomorphisme). 

\trois{cas-part-lk}
Si~$k$ est un corps, la catégorie~$\mathsf L_{\spec k}$ s'identifie à la catégorie~$\mathsf C$ 
définie au~\ref{desc-pnk-bon}. 

\trois{exemple-o1}
Comme~$\PP^n_A=\bigcup D(T_i)$, la famille~$(\sch O(1), (T_i))$ est un objet de~$\mathsf L_{\PP^n_A}$. 

\trois{fonctorialite}
Soit~$\Psi \colon S'\to S$ un morphisme de schémas et soit~$(\sch L, (s_i))$ un objet de~$\mathsf L_S$. Il est immédiat
que~$(\Psi^*\sch L, (\Psi^*s_i))$ est un objet de~$\mathsf L_{S'}$. 

\deux{fibre-morphisme-pn}
Soit~$S$ un~$A$-schéma et soit~$(\sch L, (s_i))$ un objet de~$\mathsf L_S$. Nous allons lui associer un~$A$-morphisme
de~$S$ vers~$\PP^n_A$. 

\medskip
Soit~$U$ un ouvert de~$S$. Supposons qu'il existe un indice~$i$ tel que~$s_i|_{U}$ soit inversible. 
Comme la fonction~$s_i/s_i$ est égale à~$1$ (et est en particulier inversible),
la famille~$(s_j/s_i)_{0\leq j\leq n}$ de fonctions sur~$U$ donne lieu à un~$A$-morphisme
$$[s_0/s_i:\ldots: s_n/s_i]$$
de~$U$ vers~$\PP^n_A$, et même vers~$U_i$.  

Si~$j$ est un autre indice tel que
$s_j$ soit inversible sur~$U$, on a pour tout~$\ell$ l'égalité~$s_\ell/s_i=(s_j/s_i)s_\ell/s_j$, et les deux morphismes
$$[s_0/s_i:\ldots: s_n/s_i]\;\;\text{et}\;\;[s_0/s_j:\ldots: s_n/s_j]$$ de~$U$ vers~$\PP^n_A$
coïncident donc. On a ainsi construit un~$A$-morphisme $\chi_U\colon U\to \PP^n_A$ qui ne dépend d'aucun choix. Il est 
immédiat que si~$V$ est un ouvert de~$U$ alors~$\chi_V=\chi_U|_V$. 

\medskip
Comme~$S=\bigcup D(s_i)$, les ouverts~$U$ de~$S$ sur lequel l'une au moins des~$s_i$ est inversible recouvrent~$S$. 
Lorsque~$U$ parcourt l'ensemble desdits ouverts, les morphismes~$\chi_U$ se recollent en un~$A$-morphisme
$\chi \colon S \to \PP^n_A$. Il est immédiat que ce morphisme
ne dépend que de la classe d'isomorphie de~$(\sch L, (s_i))$ (un isomorphisme entre deux objets de~$\mathsf L_S$
préservant les quotients des sections concernées). 

\medskip
Il résulte de la définition de~$\chi$ ainsi que de~\ref{test-appartient-ui}
que pour tout~$(i,j)$, l'ouvert~$\chi|_{D(s_j)}^{-1}(U_i)$ de~$D(s_j)$ est égal à~$D(s_i/s_j)$, 
c'est-à-dire à~$D(s_i)\cap D(s_j)$. En fixant~$i$ et faisant varier~$j$, il vient
$\chi^{-1}(U_i)=D(s_i)$.

Comme~$\chi|_{D(s_i)}=[s_0/s_i:\ldots:s_n/s_i]$ on déduit de~\ref{desc-crochet-ui}
que~$(\chi|_{D(s_i)})^*(T_j/T_i)=(s_j/s_i)/(s_i/s_i)=s_j/s_i$
pour tout~$j\neq i$ (et c'est d'ailleurs vrai trivialement aussi pour~$j=i$).  

\medskip
Le morphisme~$\chi$ sera noté~$[s_0:\ldots:s_n]$.

\deux{rem-fibre-morphisme-pn}
{\bf Remarque.}
Soit~$S$ un~$A$-schéma et soient~$(f_0,\ldots, f_n)$ 
des fonctions sur~$S$ telles que~$S=\bigcup D(f_i)$. 
La famille~$(\sch O_S, (f_i))$ est alors un objet de~$\mathsf L_S$, et il est immédiat
que le~$A$-morphisme~$[f_0:\ldots, f_n]$ défini ci-dessus
coïncide avec celui que nous notions~$[f_0:\ldots:f_n]$ jusqu'à présent. Il n'y a donc
pas de conflits de notations. 

\deux{annonce-fonct-pn}
Nous sommes maintenant en mesure de donner une description précise du foncteur~$S\mapsto \PP^n_A(S)$
qui généralisera ce qui a été fait au~\ref{desc-pnk-bon}. 

\deux{theo-fonct-pna}
{\bf Théorème.}
{\em Soit~$S$ un~$A$-schéma. Les flèches
$$\chi \mapsto (\chi^*\sch O(1), (\chi^*T_i))$$
et
$$(\sch L, (s_i))\mapsto [s_0:\ldots:s_n]$$ établissent
une bijection fonctorielle en~$S$ entre~$\PP^n_A(S)$ et l'ensemble
des classes d'isomorphie d'objets de~$\mathsf L_S$.}

\medskip
{\em Démonstration}. 
Les flèches de l'énoncé constituent clairement deux applications fonctorielles en~$S$. Il reste à s'assurer
qu'elles sont réciproques l'une de l'autre. 

\trois{mor-vers-fib}
Soit~$\chi$ un~$A$-morphisme de~$S$ vers~$\PP^n_A$. Nous allons montrer
que le morphisme~$[\chi^*T_0,\ldots, \chi^*T_n]$ coïncide avec~$\chi$. 
C'est une propriété qu'il suffit de vérifier localement ; nous allons donc nous assurer
qu'elle est vraie sur
l'ouvert~$S_i=\chi^{-1}(U_i)=D(\chi^*T_i)$ pour tout~$i$, ce qui permettra de conclure. 

\medskip
Soit~$i\in \{0,\ldots, n\}$. Par définition, la restriction de~$[\chi^*T_0,\ldots, \chi^*T_n]$
à~$S_i$
est l'élément
$$[\chi^*T_0/\chi^*T_i:\ldots: \chi^*T_n/\chi^*T_i]=[\chi^*(T_0/T_i):\ldots:\chi^*(T_n/T_i)],$$
de~$U_i(S_i)$. 

\medskip
On sait par ailleurs d'après~\ref{desc-crochet-ui} 
que~$\theta \mapsto (\theta^*(T_j/T_i))_{j\neq i}$ établit une bijection entre~$U_i(S_i)$ 
et l'ensemble des~$n$-uplets~$(g_j)_{j\neq i}$ de fonctions sur~$S_i$ ; et que modulo cette bijection,~$[\chi^*(T_0/T_i):\ldots: \chi^*(T_n/T_i)]$
correspond au~$n$-uplet $(\chi^*(T_j/T_i)/\chi^*(T_i/T_i))_j=(\chi^*(T_j/T_i))_j$. Il en résulte immédiatement
que~$[\chi^*(T_0/T_i):\ldots: \chi^*(T_n/T_i)]=\chi|_{S_i}$, comme annoncé.

\trois{fib-vers-mor}
Soit~$(\sch L, (s_i))$ un objet de~$\mathsf L_S$.
Posons
$$\chi=[s_0:\ldots:s_n].$$
Nous allons
montrer que~$(\chi^*\sch O(1), (\chi^*T_i))$ est isomorphe à~$(\sch L, (s_i))$. 

\medskip
Fixons~$i$. On sait que~$\chi^{-1}(U_i)=D(s_i)$ (\ref{fibre-morphisme-pn}). Sur cet ouvert, $s_i$ est une section
inversible de~$\sch L$, et~$\chi^*T_i$ est une section inversible de~$\sch O(1)$. Il existe donc
un unique isomorphisme~$\ell_i \colon \sch L|_{U_i}\simeq \chi^*\sch O(1)|_{U_i}$ qui envoie~$s_i$
sur~$\chi^*T_i$. 
Soit~$j\in \{0,\ldots, n\}$. On a les égalités
$$\ell_i(s_j)=\ell_i((s_j/s_i)s_i)=(s_j/s_i)\ell_i(s_i)=(s_j/s_i)\chi^*T_i=\chi^*(T_j/T_i)\chi^*T_i=\chi^*T_j$$
(pour l'avant dernière égalité, {\em cf}. \ref{fibre-morphisme-pn}). 

\medskip
On voit en particulier que~$\ell_i|_{D(s_i)\cap D(s_j)}$ est l'unique isomorphisme de~$\sch L|_{D(s_i)\cap D(s_j)}$
sur~$\chi^*\sch O(1)|_{D(s_i)\cap D(s_j)}$ qui envoie~$s_j$ sur~$\chi^*T_j$ ; il coïncide donc nécessairement
avec~$\ell_j|_{D(s_i)\cap D(s_j)}$. 

\medskip
On en déduit que les isomorphismes~$\ell_i$ se recollent en un isomorphisme~$\ell \colon \sch L \to \chi^*\sch O(1)$. On a vu 
au cours de la preuve que l'égalité~$\ell(s_j)=\chi^*T_j$ valait pour tout~$j$ sur chacun des~$U_i$. Deux sections d'un faisceau qui sont
localement égales le sont globalement, et l'on a donc~$\ell(s_j)=\chi^*T_j$ pour tout~$j$, ce qui achève la démonstration.~$\Box$ 

\deux{comment-fonct-pn}
{\bf Commentaires.}
Soit~$S$ un~$A$-schéma. 

\trois{fibre-meme-mor}
Soit~$\sch L$ un~$\sch O_S$-module localement libre de rang~$1$, 
et soient~$(s_i)_{0\leq i\leq n}$ et~$(t_i)_{0\leq i\leq n}$ deux familles de sections globales de
$\sch L$ telles que
$$S=\bigcup D(s_i)=\bigcup D(t_i).$$ Ces familles définissent
deux morphismes~$[s_0:\ldots: s_n]$ et~$[t_0:\ldots:t_n]$ de~$S$
vers~$\PP^n_A$. 

\medskip
En vertu du théorème~\ref{theo-fonct-pna}, 
ces deux morphismes coïncident si et seulement si il existe un automorphisme~$\ell$ de~$\sch L$
envoyant~$s_i$ sur~$t_i$ pour tout~$i$. Mais les automorphismes de~$\sch L$ sont précisément les homothéties
de rapport inversible ; en conséquence, 
$$[s_0:\ldots :s_n]=[t_0:\ldots:t_n]$$ si et seulement si il existe~$\lambda \in \sch O_S(S)\ti$ 
tel que~$t_i=\lambda s_i$ pour tout~$i$ : on généralise ainsi~\ref{relation-generale-crochet}
au cas des morphismes définis par une famille de sections de {\em n'importe
quel}
$\sch O_S$-module localement libre de rang~$1$ ({\em cf}. remarque~\ref{rem-fibre-morphisme-pn}) 

\trois{cas-identite-pastriv}
Les faits suivants se déduisent du théorème~\ref{theo-fonct-pna}
et de la remarque~\ref{rem-fibre-morphisme-pn} : 

$\bullet$ si~$(\sch L, (s_i))$ est un objet de~$\mathsf L_S$, le morphisme~$[s_0:\ldots:s_n]$ appartient
à~$\PP_A^{n,\sharp}(S)$ si et seulement si~$\sch L$ est trivial ; 

$\bullet$ si~$\chi \colon S\to \PP^n_A$ est un~$A$-morphisme, il appartient à~$\PP_A^{n,\sharp}(S)$ si et seulement
si~$\chi^*\sch O(1)$ est trivial. 

\medskip
On peut ainsi donner un autre éclairage sur le contre-exemple~\ref{pdiese-pastout} : 
on a évidemment~${\rm Id}_{\PP^n_A}^*\sch O(1)=\sch O(1)$, et l'on a signalé au~\ref{sec-globales-od}
que~$\sch O(1)$ n'est pas trivial dès que~$A\neq\{0\}$ et~$n\geq 1$. 

\trois{cas-pndiese-pn}
Il résulte de~\ref{cas-identite-pastriv}
que pour que~$\PP^n_A(S)=\PP^{n,\sharp}_A(S)$, il suffit que~tout~$\sch O_S$-module
localement libre de rang~$1$ soit trivial. Citons trois cas dans lesquels cette dernière propriété est satisfaite
(en ce qui concerne les deux premiers, nous avions déjà démontré directement l'égalité~$\PP^n_A(S)=\PP^{n,\sharp}_A(S)$
en~\ref{pndiese-corps}
et~\ref{pndiese-anneaulocal}).  

\medskip
$\bullet$ Le cas où~$S$ est le spectre d'un corps (c'est évident). 

$\bullet$ Le cas où~$S$ est le spectre d'un anneau local (c'est dû au lemme~\ref{pndiese-anneaulocal-lemme}). 

$\bullet$ Le cas où~$S$ est le spectre d'un anneau principal (c'est une conséquence de~\ref{loc-libre-qch}
et
du corollaire~\ref{coro-proj-prin}).

\medskip
On déduit notamment de ce dernier exemple que~$\PP^n_{\ZZ}(\ZZ)=\PP_{\ZZ}^{n,\sharp}(\ZZ)$. 

\deux{fin-chapitre-fonctpn}
{\bf Fonctorialité des différentes constructions}. 
Dans tout ce chapitre, nous avions fixé un anneau de base~$A$, et la plupart du temps nous avons omis de le mentionner
explicitement dans les notations (excepté pour~$\PP^n_A$). En toute rigueur,
nous aurions dû parler de l'ouvert~$V_A$ de~$\Aff^{n+1}_A$, du morphisme
$\Psi_A \colon V_A \to \PP^n_A$, et des faisceaux~$\sch O(d)_A$. Adoptons pour un instant ces conventions, plus précises (mais également un peu plus lourdes). 
Soit~$B$ une~$A$-algèbre. On vérifie alors sans peine (nous vous laissons le faire en exercice)
que toutes nos constructions se comportent bien par extension des scalaires de~$A$ à~$B$. Plus précisément :

\medskip
$\bullet$ $\Psi_B \colon V_B\to \PP^n_B$ se déduit de~$V_A\to \PP^n_A$ par produit fibré avec~$\spec B$ au-dessus de~$\spec A$ ; 

$\bullet$ le~$\sch O_{\PP^n_B}$-module~$\sch O(1)_B$
s'identifie à~$\pi^*\sch O(1)_A$ où~$\pi$ est le morphisme canonique~$\PP^n_B\to \PP^n_A$ ; 

$\bullet$ si~$S$ est un~$A$-schéma, si~$(\sch L, (s_i))$ est un objet de~$\mathsf L_S$, 
et si~$p$ désigne le morphisme canonique de~$S_B:=S\times_{\spec A}\spec B$
vers~$S$ alors
le morphisme
$$[p^*s_0:\ldots :p^*s_n]\colon S_B\to \PP^n_B$$ 
se déduit de~$[s_0:\ldots:s_n]$ par produit fibré avec~$\spec B$ au-dessus de~$\spec A$. 

\section{Quelques exemples de morphismes en géométrie projective}
\markboth{Schémas projectifs}{Quelques morphismes}

\deux{intro-mor-proj}
Jusqu'à maintenant, nous avons vu une seule méthode de construction
de morphismes en géométrie projective, consistant à exploiter la fonctorialité
(partielle) du schéma~$\proj B$ en l'anneau gradué~$B$ (\ref{projb-fonct}
{\em et sq.}). Le moins qu'on puisse dire
est qu'elle n'est pas particulièrement engageante, et le but de cette section 
est d'en proposer d'autres, fondées sur la description explicite du foncteur des points
de l'espace projectif (th.~\ref{theo-fonct-pna}). 

Nous observerons à cette occasion une nouvelle
manifestation de la philosophie dégagée à la section~\ref{FONCTPOINT}. En effet, comme nous le verrons, 
cette approche fonctorielle
permettra peu ou prou de retrouver le point de vue naïf ou ensembliste sur la géométrie projective, selon lequel 
les sous-variétés sont les lieux des zéros de systèmes d'équations polynomiales homogènes, et les morphismes des
applications définies par des formules polynomiales homogènes ; à une petite subtilité près toutefois : il faudra considérer
les polynômes homogènes en un sens «tensoriel» et non «multiplicatif» -- pour la signification précise
de cette 
remarque, {\em cf.}~\ref{polynom-interp-tens}
{\em infra}. 

\deux{notations-mor-pro}
{\bf Quelques notations}. 
On fixe pour toute la suite de la section
un anneau~$A$. 

\trois{notation-lns} Si~$S$ est un schéma et~$n$ un entier, nous noterons~$\ml n S$ la catégorie
que nous avions simplement notée~$\mathsf L_S$ au~\ref{def-lt} (l'entier~$n$ était alors fixé une fois pour toutes ; 
ce ne sera pas le cas dans cette section, et il est donc préférable de le faire figurer explicitement dans les notations). 
Nous désignerons par~$\isol n S$ l'ensemble des classes d'isomorphie d'objets de~$\ml n S$. 

Pour tout~$n\in \NN$ et tout~$A$-schéma~$S$, le théorème~\ref{theo-fonct-pna}
fournit une bijection~$$\PP^n_A(S)\simeq \isol n S$$ fonctorielle en~$S$. 

\trois{polynom-interp-tens}
{\em Interprétation tensorielle d'un polynôme homogène}. 
Soient~$n$ et~$d$ deux entiers et soit~$P$ un 
polynôme homogène de degré~$d$ appartenant à~$A[T_0,\ldots, T_n]$. 
Écrivons~$P=\sum_{(e_i)}a_{(e_i)}\prod_i T_i^{e_i}$, où~$(e_i)$ parcourt
la famille des~$(n+1)$-uplets d'entiers de somme égale à~$d$, et où les~$a_{(e_i)}$ 
sont des scalaires (évidemment presque tous nuls). 

\medskip
Soit~$S$ un~$A$-schéma, soit~$\sch L$ un~$\sch O_S$-module localement
libre de rang~$1$ et soit~$(s_i)_{0\leq i\leq n}$ une famille de sections globales de~$\sch L$. 
On pose alors
$$P(s_0,\ldots,s_n)=\sum_{(e_i)}a_{(e_i)}\bigotimes_i s_i^{\otimes e_i}\in \sch L^{\otimes d}(S).$$

Remarquons que si~$\sch L=\sch O_X$ alors~$\sch L^{\otimes d}$ s'identifie à~$\sch O_X$
{\em via}
la multiplication des fonctions, et que modulo cette identification $P(s_0,\ldots, s_n)$ a son sens habituel. 

\deux{reinterp-famille-homog}
{\bf Morphisme donné par une famille de polynômes homogènes}. 
Soient~$n, m$ et~$d$ trois entiers, avec~$d>0$. Soit~$(P_0,\ldots, P_n)$ une famille de polynômes
homogènes de degré~$d$ en les variables~$(S_0,\ldots, S_m)$. 

\trois{defpsi-mor-pnm}
Soit~$\phi$ l'unique
morphisme de~$A$-algèbres de~$A[T_0,\ldots, T_n]$ dans~$A[S_0,\ldots, S_m]$ 
qui envoie~$T_i$ sur~$P_i$ pour tout~$d$. Il est homogène de degré~$d$, et induit donc
en vertu de~\ref{projb-fonct}
{\em et sq.}
un morphisme de~$A$-schémas $\psi \colon \Omega \to \PP^n _A$ où~$\Omega$ est 
l'ouvert~$\bigcup D(P_i)$ de~$\PP^m_A$. 

\trois{defchi-mor-pnm}
Par ailleurs, $(\sch O(d)|_{\Omega}, (P_i|_\Omega))$ est un objet de~$\ml n \Omega$, et définit
donc lui-même un morphisme de~$A$-schémas~$\chi \colon \Omega \to \PP^n_A$. 

\trois{lemme-egalite-chi-psi}
{\bf Lemme}. 
{\em Les morphismes~$\psi$ et~$\chi$ de~$\Omega	$ vers~$\PP^n_A$ sont égaux}. 

\medskip
{\em Démonstration}. 
Pour tout~$i$, on pose~$U_i=D(T_i)\subset \PP^n_A$. 
Il résulte de~\ref{fonct-projb-defensemble}
et~\ref{test-appartient-ui}
que l'on a pour tout~$i$
l'égalité~$\psi^{-1}(U_i)=\chi^{-1}(U_i)=D(P_i)$. Il suffit pour conclure 
de montrer que pour tout~$i$, les morphismes 
de~$D(P_i)$ vers~$U_i$ induits par~$\psi$ et~$\chi$ coïncident.  

\medskip
Soit donc~$i\in \{0,\ldots, n\}$. Pour montrer que~$\psi|_{D(P_i)}\colon D(P_i)\to U_i$ et
$\chi|_{D(P_i)}\colon D(P_i)\to U_i$ coïncident, il suffit de s'assurer que~$\chi^*(T_j/T_i)=\psi^*(T_j/T_i)$
pour tout~$j\neq i$ (\ref{desc-crochet-ui}). Or
pour tout~$j\neq i$ on a~$\chi^*(T_j/T_i)=P_j/P_i$ d'après~\ref{desc-crochet-ui}, et
$\psi^*(T_j/T_i)=P_j/P_i$ d'après~\ref{descr-expl-mor}.~$\Box$  

\subsection*{Immersions ouvertes et fermées}

\deux{imm-ouv-pn-fonct}
{\bf Immersions ouvertes}. Soient~$n$ et~$d$ deux entiers, et soit~$f$ un élément homogène
de degré~$d$
de~$A[T_0,\ldots, T_n]$. 

\trois{df-sousens-isolns}
Soit~$S$ un~$A$-schéma. L'ensemble~$D(f)(S)$ est de manière
naturelle un sous-ensemble de~$\PP^n_A(S)$ (c'est l'ensemble des morphismes de~$A$-schémas
de~$S$ vers~$\PP^n_A$ dont l'image ensembliste est contenue dans~$D(f)$). 

\medskip
Modulo l'identification canonique de~$\PP^n_A(S)$ à~$\isol n S$, l'ensemble~$D(f)(S)$ apparaît
dès lors
comme un sous-ensemble de~$\isol n S$. Le but du lemme qui suit est d'en donner une description,
aussi proche que possible de l'intuition ensembliste qui fait de~$D(f)$ le lieu de non-annulation, ou
plus exactement d'inversibilité, de~$f$.

\trois{df-pn-intuitif}
{\bf Lemme.}
{\em Le sous-ensemble~$D(f)(S)$ de~$\isol n S$ est constitué des classes d'objets
$(\sch L, (s_i))$ tels que la section~$f(s_0,\ldots, s_n)$ de~$\sch L^{\otimes d}$ soit inversible}. 

\medskip
{\em Démonstration}. Soit~$\psi \colon S \to \PP^n_A$ un morphisme. L'élément de~$\isol n S$
auquel il correspond est la classe de~$(\psi^*\sch O(1), (\psi^*T_i))$. Pour tout entier~$i$ entre~$0$
et~$n$, posons~$s_i=\psi^*T_i, U_i=D(T_i)$ 
et~$S_i=\psi^{-1}(U_i)=D(s_i)$. On a~$\psi(S)\subset D(f)$ si et seulement si~$\psi(S_i)$
est contenu dans~$D(f)\cap U_i$ pour tout~$i$. 

\medskip
Fixons~$i$. L'ouvert~$D(f)\cap U_i$ de~$U_i=\spec A\left[\frac{T_j}{T_i}\right]$ est égal à~$D\left(\frac f{T_i^d}\right)$. En conséquence,
$\psi(S_i)$ s'envoie dans~$D(f)\cap U_i$ si et seulement si~$\psi^*\frac f{T_i^d}$
appartient à~$\sch O_{S_i}(S_i)\ti$. Or
l'élément$\psi^*\frac f{T_i^d}$ de~$\sch O_{S_i}(S_i)$ est égal à~$f(s_0,\ldots, s_n)/(s_i^{\otimes d})$, et est donc inversible
si et seulement si~$f(s_0, \ldots, s_n)|_{S_i}$ est  inversible.

\medskip
En conséquence, $\psi(S)\subset D(f)$ si et seulement si~$f(s_0, \ldots, s_n)|_{S_i}$ est inversible pour tout~$i$, c'est-à-dire
si et seulement si~$f(s_0,\ldots,s_n)$ est inversible.~$\Box$ 

\deux{imm-ferm-pn-fonct}
{\bf Immersions fermées}. Soit~$n$ un entier et soit~$I$ un idéal homogène
de~$A[T_0,\ldots, T_n]$. Donnons-nous une famille
génératrice~$(g_\ell)$ de~$I$, où chaque~$g_\ell$
est homogène d'un certain degré~$d_\ell$. Soit~$Z$ le sous-schéma fermé~$\proj A[T_0,\ldots, T_n]/I$
de~$\PP^n_A$. 

\trois{vi-sousens-isolns}
Soit~$S$ un~$A$-schéma. L'ensemble~$Z(S)$ est de manière
naturelle un sous-ensemble de~$\PP^n_A(S)$ (c'est l'ensemble des morphismes~$\psi$ de~$A$-schémas
qui se factorisent par~$Z$, c'est-à-dire encore qui sont tels
que~$\psi^*a=0$ pour toute section~$a$ du faisceau quasi-cohérent
d'idéaux définissant~$Z$). 

\medskip
Modulo l'identification canonique de~$\PP^n_A(S)$ à~$\isol n S$, l'ensemble~$Z(S)$ apparaît
dès lors
comme un sous-ensemble de~$\isol n S$. Le but du lemme qui suit est d'en donner une description,
aussi proche que possible de l'intuition ensembliste qui fait de~$Z$ le lieu des zéros des~$g_\ell$. 

\trois{vi-pn-intuitif}
{\bf Lemme.}
{\em Le sous-ensemble~$Z(S)$ de~$\isol n S$ est constitué des classes d'objets
$(\sch L, (s_i))$ tels que~$g_\ell(s_0,\ldots, s_n)=0$ pour tout~$\ell$.}

\medskip
{\em Démonstration}. Soit~$\psi \colon S \to \PP^n_A$ un morphisme. L'élément de~$\isol n S$
auquel il correspond est la classe de~$(\psi^*\sch O(1), (\psi^*T_i))$. Pour tout entier~$i$ entre~$0$
et~$n$, posons~$s_i=\psi^*T_i, U_i=D(T_i)$ 
et~$S_i=\psi^{-1}(U_i)=D(s_i)$. Le morphisme~$\psi$ se factorise par~$Z$ si et seulement et seulement si~$\psi|_{S_i}\colon S_i \to U_i$
se factorise par~$Z_i:=Z\times_{\PP^n_A}U_i$ pour tout~$i$. 

\medskip
Fixons~$i$. Le sous-schéma fermé~$Z_i$
de~$U_i$ est défini par l'idéal~$I_i$ de~$A\left[\frac{T_j}{T_i}\right]_{j\neq i}$ engendré par les~$\frac{g_\ell}{T_i^{d_\ell}}$. 
En conséquence,
$\psi|_{S_i}$ se factorise par~$Z_i$ si et seulement si~$\psi^*\frac {g_\ell}{T_i^{d_\ell}}=0$ pour tout~$\ell$. 
Or 
l'élément~$\psi^*\frac {g_\ell}{T_i^{d_\ell}}$ de~$\sch O_{S_i}(S_i)$ est égal pour tout~$\ell$
à~$g_\ell(s_0,\ldots, s_n)/(s_i^{\otimes d_\ell})$, et il est donc nul 
si et seulement si~$g_\ell(s_0, \ldots, s_n)|_{S_i}=0$. 

\medskip
En conséquence, $\psi$ se factorise par~$Z$ si et seulement si on a pour tout~$\ell$ et tout~$i$ l'égalité~$g_\ell(s_0, \ldots, s_n)|_{S_i}=0$, 
ce qui revient à demander que~$g_\ell(s_0, \ldots, s_n)=0$ pour tout~$\ell$.~$\Box$ 

\subsection*{Un plongement de~$\PP^1_A$ dans~$\PP^2_A$}

\deux{plonge-conique}
Le but de ce qui suit est de construire
dans le cadre schématique
l'immersion fermée de~$\PP^1_A$ dans~$\PP^2_A$ qui est naïvement
donnée
par la formule~$[s_0:s_1]\mapsto [s_0^2:s_0s_1:s_1^2]$. Comme nous allons le voir, 
le point de vue «foncteur des points» permet de la {\em définir rigoureusement}
par cette même
formule
-- ou plus exactement par sa déclinaison tensorielle. 

\deux{def-mor-conique}
Soit~$S$ un~$A$-schéma, et soit~$(\sch L, (s_0,s_1))$ un objet de~$\ml 1 S$. 
Il est immédiat que~$\Phi(\sch L , (s_0, s_1)):=(\sch L^{\otimes 2}, (s_0^{\otimes 2}, s_0\otimes s_1, s_1^{\otimes 2}))$ 
est un objet de~$\ml 2 S$, dont la classe d'isomorphie ne dépend que de celle de~$(\sch L, (s_0, s_1))$. On a ainsi défini
une application encore notée~$\Phi$, fonctorielle en~$S$, de~$\isol 1 S$ vers~$\isol 2 S$, et partant un morphisme~$\psi \colon \PP^1_A\to \PP^2_A$. 

\deux{prop-p1-conique}
{\bf Proposition}. {\em Le morphisme~$\psi$
induit un isomorphisme
$$\PP^1_A\simeq \proj A[T_0, T_1, T_2]/(T_0T_2-T_1^2)\hookrightarrow \PP^2_A.$$}

\medskip
{\em Démonstration}. Nous allons faire un usage intensif des constructions de faisceaux
localement libres de rang~$1$ à partir de cocycles, détaillées en~\ref{cocy-cobo}
{\em et sq.} Soit~$S$ un~$A$-schéma, et soit~$\mathsf D$ le sous-ensemble
de~$\isol 2 S$ formé des classes d'objets~$(\sch M, (t_0, t_1,t_2))$ tels que~$t_0\otimes t_2=t_1^{\otimes 2}$. En vertu
du lemme~\ref{vi-pn-intuitif}, il suffit de démontrer que 
l'application~$\Phi \colon \isol 1 S \to \isol 2 S$ induit une bijection~$\isol 1 S \simeq \mathsf D$. 
Il est immédiat que~$\Phi(\isol 1 S)\subset \mathsf D$. 

\trois{conique-def-theta}
Soit~$(\sch M, (t_0, t_1,t_2))$ un objet de~$\ml 2 S$
tel que~$t_0\otimes t_2=t_1^{\otimes 2}$. Cette égalité assure que si~$t_0$ et~$t_2$ s'annulent en un point de~$S$, 
il en va de même de~$t_1$, ce qui est absurde. En conséquence, la réunion de~$S_0:=D(t_0)$ et~$S_2:=D(t_2)$ est égale
à~$S$. Soit~$(f_{ij})$ le cocycle subordonné au recouvrement de~$S$ par~$S_0$ et~$S_2$, défini
par la formule~$f_{20}=t_1/t_0$ (remarquons : que~$t_1/t_0$ appartient bien à~$\sch O_S(S_0\cap S_2)\ti$ car~$t_1^{\otimes 2}=t_0\otimes t_2$ ; et que
les autres~$f_{ij}$ s'obtiennent grâce aux relations de cocycle)

\medskip
Soit~$\sch L$
le~$\sch O_S$-module localement libre de rang~$1$ obtenu en tordant~$\sch M$ par le cocyle~$(f_{ij})$
(\ref{cocycle-vers-fibre}). Rappelons brièvement ce que cela signifie. Les restrictions~$\sch L|_{S_0}$ et~$\sch L|_{S_2}$ 
s'identifient respectivement à~$\sch M|_{S_0}$ et~$\sch M|_{S_2}$, mais la conditions de coïncidence de deux sections
sur~$S_0\cap S_2$ est tordue : si l'on se donne un ouvert~$U_0$ de~$S_0$ et une section~$a_0\in \sch M(U_0)=\sch L(U_0)$, 
ainsi qu'un ouvert~$U_2$ de~$S_2$ et une section~$a_2$ de~$\sch M(U_2)=\sch L(U_2)$, les sections~$a_0$ et~$a_2$ {\em du
faisceau~$\sch L$}
coïncident sur~$U_0\cap U_2$ si et seulement si~$a_2=f_{20}a_0=(t_1/t_0)a_0$ dans~$\sch M(U_0\cap U_2)$. 

On prendra garde que lorsqu'on travaille sur un ouvert~$V\subset S_0\cap S_2$, il y a deux manières différentes d'identifier~$\sch L|_{V}$
à~$\sch M|_{V}$, selon qu'on voit~$V$ comme contenu dans~$S_0$ ou dans~$S_2$ ; il importe, lorsqu'on doit effectuer les calculs, de bien préciser
laquelle de ces deux identifications on utilise, et surtout de ne pas les mélanger indûment.

\medskip
Par définition de~$\sch L$,  les sections~$t_0|_{S_0}$ et~$t_1|_{S_2}$
de~$\sch M$
se recollent en une section globale~$s_0$ de~$\sch L$.
Comme~$t_0\otimes t_2=t_1^{\otimes 2}$, on a dans l'anneau~$\sch O_S(S_0\cap S_2)$ l'égalité
$$\frac {t_2}{t_1}=\frac {t_1}{t_0},$$
et les sections~$t_1|_{S_0}$ et~$t_2|_{S_2}$
de~$\sch M$
se recollent donc en une section globale~$s_0$ de~$\sch L$. Comme la section~$t_0$ de~$\sch M$ ne s'annule pas sur~$S_0$, la section~$s_0$ de~$\sch L$ ne s'annule pas sur~$S_0$ ; comme la section~$t_2$ de~$\sch M$ ne s'annule pas sur~$S_2$, la section~$s_1$ de~$\sch L$ ne s'annule pas sur~$S_2$. 
En conséquence, $(\sch L, (s_0, s_1))$ est un objet de~$\ml 1 S$ que nous noterons~$\Theta(\sch M, (t_0, t_1,t_2))$. Sa classe d'isomorphie ne
dépend visiblement que de la classe
d'isomorphie de~$(\sch M, (t_0, t_1,t_2))$, et l'on note encore~$\Theta$ l'application de~$\mathsf D$ vers~$\isol 1 S$ induite par ce procédé. Nous allons montrer
que~$\Theta \circ \Phi={\rm Id}_{\isol 1 S}$ et que~$\Phi\circ \Theta={\rm Id}_{\mathsf D}$, ce qui permettra de conclure. 

\trois{conique-fonct-sens1}
{\em Montrons que~$\Theta \circ \Phi={\rm Id}_{\isol 1 S}$}.
Soit~$(\sch L, (s_0,s_1))$ un objet de~$\ml 1 S$, et posons~$(\sch N, (\sigma_0, \sigma_1))=\Theta (\sch L^{\otimes 2}, (s_0^{\otimes 2}, s_0\otimes s_1, s_1^{\otimes 2}))$. 
Nous allons prouver que~$(\sch L, (s_0,s_1))\simeq (\sch N, (\sigma_0, \sigma_1))$. D'après la construction détaillée de~$\Theta$ au~\ref{conique-def-theta},
les restrictions de~$\sigma_0$ à~$D(s_0^{\otimes 2})=D(s_0)$ et~$D(s_1^{\otimes 2})=D(s_1)$ sont respectivement «égales»
à~$s_0^{\otimes 2}$ et~$s_0\otimes s_1$ ; et les restrictions  de~$\sigma_1$ à~$D(s_0)$ et~$D(s_1)$ sont
respectivement égales à~$s_0\otimes s_1$ et~$s_1^{\otimes 2}$.

\medskip
La section~$\sigma_0$ est inversible sur~$D(s_0)$ et la section~$\sigma_1$ est inversible
sur~$D(s_1)$. En conséquence, il existe deux isomorphismes
$$\ell_0 \colon \sch L|_{D(s_0)}\simeq \sch N|_{D(s_0)}\;\;\text{et}\;\;\ell_1 \colon \sch L|_{D(s_1)}\simeq \sch N|_{D(s_1)}$$ tels que~$\ell_0(s_0)=\sigma_0$ et~$\ell_1(s_1)=\sigma_1$
(chacun d'eux est {\em caractérisé}
par l'égalité correspondante). 
Par définition de~$\sigma_0$ et~$\sigma_1$ on a sur~$D(s_0)\cap D(s_1)\subset D(s_0)$ l'égalité
$$\frac{\sigma_0}{\sigma_1}=\frac{s_0^{\otimes 2}}{s_0\otimes s_1}=\frac {s_0}{s_1},$$ ce qui montre que les restrictions de~$\ell_0$ et~$\ell_1$ à
$D(s_0)\cap D(s_1)$ coïncident, et permet de les recoller en un isomorphisme~$\ell\colon \sch L\simeq \sch N$. Il reste à s'assurer que~$\ell(s_0)=\sigma_0$ 
et~$\ell(s_1)=\sigma_1$. L'égalité~$\ell(s_0)=\sigma_0$ est vraie par définition 
de~$\ell$ sur~$D(s_0)$, et l'on a sur~$D(s_1)$ les égalités
$$\ell(s_0)=\frac {s_0}{s_1}\ell(s_1)=\frac {s_0}{s_1}\sigma_1=\frac{s_0}{s_1}s_1^{\otimes 2}=s_0\otimes s_1=\sigma_0.$$
En conséquence, $\ell(s_0)=\sigma_0$ sur~$S$ tout entier. 

\medskip
L'égalité~$\ell(s_1)=\sigma_1$ est vraie par définition 
de~$\ell$ sur~$D(s_1)$, et l'on a sur~$D(s_0)$ les égalités
$$\ell(s_1)=\frac {s_1}{s_0}\ell(s_0)=\frac {s_1}{s_0}\sigma_0=\frac{s_1}{s_0}s_0^{\otimes 2}=s_0\otimes s_1=\sigma_1.$$
En conséquence, $\ell(s_1)=\sigma_1$ sur~$S$ tout entier. 

\trois{conique-fonct-sens2}
{\em Montrons que~$\Phi\circ \Theta={\rm Id}_{\mathsf D}$}.
Soit~$(\sch M, (t_0, t_1, t_2))$ un objet de~$\ml 2 S$ dont la classe
appartient à~$\mathsf D$, et soit~$(\sch L, (s_0, s_1))$ son image par~$\Theta$, dont nous utiliserons
la définition et les propriétés élémentaires établies au~\ref{conique-def-theta}, en reprenant
les notations~$S_0=D(t_0)$ et~$S_2=D(t_2)$
de~{\em loc. cit.}.
Il s'agit
de prouver que~$(\sch L^{\otimes 2}, (s_0^{\otimes 2}, s_0\otimes s_1, s_1^{\otimes 2}))$ est isomorphe
à~$(\sch M, (t_0, t_1, t_2))$.  

\medskip
La section~$s_0$ de~$\sch L$
étant inversible sur~$S_0$, il existe un
unique isomorphisme~$\ell_0 \colon \sch M|_{S_0}\simeq \sch L^{\otimes 2}|_{S_0}$ qui envoie
$t_0$ sur~$s_0^{\otimes 2}$ ; de même, il existe un
unique isomorphisme~$\ell_2  \colon \sch M|_{S_2}\simeq \sch L^{\otimes 2}|_{S_2}$ qui envoie~$t_2$
sur~$s_1^{\otimes 2}$. 

\medskip
Par définition des sections~$s_i$ et en vertu de l'égalité~$t_0\otimes t_2=t_1^{\otimes 2}$
on a sur~$S_0\cap S_2\subset S_0$ l'égalité 
$$\frac{s_0^{\otimes 2}}{s_1^{\otimes 2}}=\frac{t_0^{\otimes 2}}{t_1^{\otimes 2}}=\frac {t_0}{t_2},$$
et les restrictions de~$\ell_0$ et~$\ell_2$ à~$S_0\cap S_2$ coïncident donc, ce qui permet de les recoller en un
isomorphisme~$\ell\colon \sch M\simeq \sch L^{\otimes 2}$. Il reste à calculer~$\ell(t_0), \ell(t_1)$ et~$\ell(t_2)$. 

\medskip
$\bullet$ On sait que~$\ell(t_0)=s_0^{\otimes 2}$ sur~$S_0$. Sur~$S_2$, on a les égalités
$$\ell(t_0)=\frac{t_0}{t_2}\ell(t_2)=\frac{t_0}{t_2}s_1^{\otimes 2}=\frac{t_0}{t_2}{t_2^{\otimes 2}}=t_0\otimes t_2=t_1^{\otimes 2}=s_0^{\otimes 2},$$
et donc~$\ell(t_0)=s_0^{\otimes 2}$ sur~$S$ tout entier. 

\medskip
$\bullet$ On sait que~$\ell(t_2)=s_1^{\otimes 2}$ sur~$S_2$. Sur~$S_0$, on a les égalités
$$\ell(t_2)=\frac{t_2}{t_0}\ell(t_0)=\frac{t_2}{t_0}s_0^{\otimes 2}=\frac{t_2}{t_0}{t_0^{\otimes 2}}=t_0\otimes t_2=t_1^{\otimes 2}=s_1^{\otimes 2},$$
et donc~$\ell(t_2)=s_1^{\otimes 2}$ sur~$S$ tout entier. 

\medskip
$\bullet$ Sur~$S_0$, on a les égalités
$$\ell(t_1)=\frac {t_1}{t_0}\ell(t_0)=\frac {t_1}{t_0} s_0^{\otimes 2} =\frac {t_1}{t_0}t_0^{\otimes 2}=t_0\otimes t_1=s_0\otimes s_1\;;$$
sur~$S_2$, on a les égalités
$$\ell(t_1)=\frac {t_1}{t_2}\ell(t_2)=\frac {t_1}{t_2} s_1^{\otimes 2} =\frac {t_1}{t_2}t_2^{\otimes 2}=t_1\otimes t_2=s_0\otimes s_1.$$
En conséquence, $\ell(t_1)=s_0\otimes s_1$ sur~$S$ tout entier, ce qui achève la démonstration.~$\Box$ 

\subsection*{Les plongements de Segre et de Veronese}

\deux{mor-segre-motivation}
{\bf Le plongement de Segre}.
Soient~$n$ et~$m$ deux entiers. Le but de ce qui suit est de construire, dans le contexte
schématique, une immersion fermée~$\PP^n_A\times_A\PP^m_A\hookrightarrow \PP^{nm+n+m}_A$, 
appelée le {\em plongement de Segre}, et qui est donnée naïvement par la formule
$$([s_0:\ldots:s_n], [t_0:\ldots :t_m])\mapsto [s_it_j]_{0\leq i\leq n, 0\leq j\leq m}$$
(notez qu'il y a à droite~$(n+1)(m+1)=nm+n+m+1$ coordonnées, et que le but
est donc bien~$\PP^{nm+n+m}$).  Là encore, 
le point de vue  «foncteur des points» permet de la {\em définir rigoureusement}
par cette même
formule
-- ou plus exactement par sa déclinaison tensorielle. 

\deux{def-mor-segre}
Soit~$\sch S$ un~$A$-schéma, soit~$(\sch L, (s_i))$ un objet de~$\ml n S$ et soit~$(\sch M, (t_j))$ un
objet de~$\ml m S$. Il est immédiat que
$$\Phi(\;(\sch L, (s_i)), (\sch M, (t_j))\;):=(\sch L\otimes \sch M, (s_i\otimes t_j)_{i,j})$$ est un objet
de~$\ml {nm+n+m} S$, dont la classe d'isomorphie ne dépend que de celles de~$(\sch L, (s_i))$ et~$(\sch M, (t_j))$. 
On obtient ainsi une application fonctorielle en~$S$, notée encore~$\Phi$, 
de~$\isol n S\times \isol m S$ vers~$\isol {nm+m+n} S$, et partant un morphisme
de~$A$-schémas 
$$\psi \colon \PP^n_A\times_A\PP^m_A\to \PP^{nm+n+m}_A,$$
appelé {\em morphisme de Segre}. 

\deux{prop-mor-segre}
{\bf Proposition.}
{\em Le morphisme~$\psi$ induit un isomorphisme
$$\PP^n_A\times_A\PP^m_A \simeq \proj A[\Sigma_{ij}]_{i,j}/(\Sigma_{ij}\Sigma_{i'j'}-\Sigma_{ij'}\Sigma_{i'j})_{i\neq i', j\neq j'}\hookrightarrow \PP^{nm+n+m}_A.$$}

\medskip
{\em Démonstration}. Posons~${\bf n}=\{0,\ldots, n\}$ et~${\bf m}=\{0,\ldots, m\}$. 
Soit~$S$ un~$A$-schéma. Notons~$\mathsf D$ le sous-ensemble de~$\isol {nm+n+m}S$ formé des classes d'objets
$(\sch N, (\sigma_{ij}))$ tels que l'on ait
$$\sigma_{ij}\otimes \sigma_{i'j'}=\sigma_{ij'}\otimes\sigma_{i'j}$$ pour tout~$(i,i',j,j')$ avec~$i\neq i'$ et~$j\neq j'$ 
(remarquez que lorsque~$i=i'$ ou~$j=j'$ l'égalité est automatiquement vérifiée). En vertu
du lemme~\ref{vi-pn-intuitif}, il suffit de démontrer que 
l'application~$\Phi \colon \isol n S\times \isol m S\to \isol {nm+n+m}S$ induit une bijection~$\isol n S \times \isol m S\simeq \mathsf D$. 
Il est immédiat que~$\Phi(\isol n S\times \isol m S)\subset \mathsf D$. 

\trois{def-theta-segre}
Soit~$(\sch N, (\sigma_{ij}))$ un objet de~$\ml {nm+n+m}S$ tel que~$\sigma_{ij}\otimes \sigma_{i'j'}=\sigma_{ij'}\otimes\sigma_{i'j}$ pour tout~$(i,i',j,j')$. 
Pour tout~$(i,j)$ on pose~$S_{ij}=D(\sigma_{ij})$. Remarquez qu'en vertu des équations satisfaites par les~$\sigma_{ij}$, on a~$S_{ij}\cap S_{i'j'}=S_{ij'}\cap S_{i'j}$
pour tout~$(i, i', j, j')$ ; nous utiliserons implicitement ce fait dans tout ce qui suit. 

\medskip
Soient~$i$ et~$i'$ deux éléments de~$\bf n$ et soient~$j$ et~$j'$ deux éléments de~$\bf m$.
On pose~$$f_{iji'j'}=\sigma_{ij}/\sigma_{ij'}=\sigma_{i'j}/\sigma_{i'j'}$$
et~$$g_{iji'j'}=\sigma_{ij}/\sigma_{i'j}=\sigma_{ij'}/\sigma_{i'j'}\;\;$$
ce sont des fonctions inversibles sur~$S_{ij}\cap S_{i'j'}$. Un calcul immédiat (qui utilise les équations
satisfaites par les~$\sigma_{ij}$) montre que~$(f_{iji'j'})$ et~$(g_{iji'j'})$ sont deux cocycles
subordonnés au recouvrement~$(S_{ij})$. On note~$\sch L$ (resp.~$\sch M$)
le~$\sch O_S$-module localement libre de rang~$1$ obtenu en tordant~$\sch N$ avec le cocycle~$(f_{iji'j'})$
(resp.~$(g_{iji'j'})$). 

\medskip
Soit~$i\in \bf n$. Pour tout~$(i', j)\in {\bf n}\times{\bf m}$, posons~$\lambda^i_{i'j}=\sigma_{ij}|_{S_{i'j}}$. 
On vérifie immédiatement  que 
$$\sigma_{ij}=\frac{\sigma_{i'j}}{\sigma_{i'j'}}\sigma_{ij'}$$ sur~$S_{i'j}\cap S_{i''j'}$
pour tout~$(i', i'', j, j')$, ce qui entraîne que les~$\lambda^i_{i'j}$ se recollent, pour~$(i',j)$ variables, 
en une section globale~$s_i$ de~$\sch L$. 

\medskip
Soit~$j\in \bf m$. Pour tout~$(i, j')\in {\bf n}\times{\bf m}$, posons~$\mu^j_{ij'}=\sigma_{ij}|_{S_{ij'}}$. 
On vérifie immédiatement que
$$\sigma_{ij}=\frac{\sigma_{ij'}}{\sigma_{i'j'}}\sigma_{i'j}$$ sur~$S_{ij'}\cap S_{i'j''}$
pour tout~$(i, i', j', j'')$, ce qui entraîne que les~$\mu^j_{ij'}$ se recollent, pour~$(i,j')$ variable, 
en une section globale~$t_j$ de~$\sch M$

\medskip
Par construction, $S_{ij}\subset D(s_i)$ et~$S_{ij}\subset D(t_j)$ pour tout~$(i,j)$. Il s'ensuit
que~$(\sch L, (s_i))$ est un objet de~$\ml n S$, et que~$(\sch M, (t_j))$ est un objet de~$\ml m S$ ; la
classe d'isomorphie 
du couple
$$\Theta(\sch N, (\sigma_{ij})):=(\;(\sch L, (s_i)), (\sch M, (t_j))\;)$$ ne dépend manifestement que de celle de~$(\sch N, (\sigma_{ij}))$. 
Il s'ensuit que~$\Theta$ induit une application, notée encore~$\Theta$,   
de~$\mathsf D$
vers~$\isol n S \times \isol m S.$ Nous 
allons prouver que~$\Theta \circ \Phi={\rm Id}_{\isol n S\times \isol m S}$
et~$\Phi\circ \Theta ={\rm Id}_{\mathsf D}.$

\trois{iso-segre-sens1}
{\em Montrons que~$\Theta \circ \Phi={\rm Id}_{\isol n S\times \isol m S}$}. 
Soit~$(\sch L, (s_i))$ un objet de~$\ml n S$ et soit~$(\sch M, (t_j))$ un objet de~$\ml m S$. 
Posons
$$(\;(\sch E, (\xi_i)), (\sch F, (\eta_j))\;)=\Theta(\sch L\otimes \sch M, (s_i\otimes t_j)).$$ Il s'agit
de prouver que~$(\sch E, (\xi_i))\simeq (\sch L, (s_i))$ et que~$(\sch F, (\eta_j))\simeq (\sch M, (t_j))$. 
Pour tout~$(i,j)$ on pose~$S_{ij}=D(s_i\otimes t_j)=D(s_i)\cap D(t_j)$. Il résulte de la définition de~$\Theta$ que pour tout~$(i, i',j)$, 
la restriction de~$\xi_i$ à~$S_{i'j}$ est «égale»
à~$s_i\otimes t_j$, et que pour tout~$(i, j, j')$ la restriction de~$\eta_j$ à~$S_{ij'}$ est «égale»
à~$s_i\otimes t_j$. 

\medskip
Fixons~$(i,j)$. La section~$\xi_i$ de~$\sch E$ est inversible sur~$S_{ij}$, et il en va de même de la section~$s_i$
de~$\sch L$. Il existe donc un unique isomorphisme~$\ell_{ij}$ de~$\sch L|_{S_{ij}}$ sur~$\sch E|_{S_{ij}}$ envoyant~$s_i$
sur~$\xi_i$. 

\medskip
Pour tout~$(i,i',j,j')$ on a les égalités
$$\frac {\xi_i}{\xi_{i'}}=\frac{s_i\otimes t_j}{s_{i'}\otimes t_j}=\frac {s_i}{s_{i'}}$$
dans l'anneau des fonctions de~$S_{ij}\cap S_{i'j'}\subset S_{ij}$, ce qui implique que~$\ell_{ij}$ et~$\ell_{i'j'}$ coïncident
sur~$S_{ij}\cap S_{i'j'}$. Les~$\ell_{ij}$ se recollent donc en un isomorphisme~$\ell \colon \sch L\to \sch E$. Il reste à calculer
$\ell(s_i)$ pour tout~$i$. 

\medskip
Fixons~$i$. Soit~$(i', j)\in {\bf n}\times {\bf m}$. Sur l'ouvert~$S_{i'j}$, on a les égalités
$$\ell(s_i)=\frac {s_i}{s_{i'}}\ell(s_{i'})=\frac {s_i}{s_{i'}}\xi_{i'}=\frac {s_i}{s_{i'}}s_{i'}\otimes t_j=s_i\otimes t_j=\xi_i.$$
Il s'ensuit que~$\ell(s_i)=\xi_i$ sur tout~$S$, ce qui achève de montrer que~$(\sch E, (\xi_i))$ et~$(\sch L, (s_i))$
sont isomorphes. 

\medskip
La démonstration que~$(\sch F, (\eta_j))\simeq (\sch M, (t_j))$
est la même {\em mutatis mutandis}.

\trois{iso-segre-sens2}
{\em Montrons que~$\Phi\circ \Theta ={\rm Id}_{\mathsf D}.$}
Soit~$(\sch N, (\sigma_{ij}))$ un objet de~$\ml {nm+n+m} S$
tel que~$\sigma_{ij}\otimes \sigma_{i'j'}=\sigma_{ij'}\otimes\sigma_{i'j}$ pour tout~$(i,i',j,j')$. 
Posons
$$(\;(\sch L, (s_i)), (\sch M, (t_j))\;)=\Theta(\sch N, (\sigma_{ij})).$$ Il s'agit
de montrer que~$(\sch N, (\sigma_{ij}))\simeq (\sch L\otimes \sch M, s_i \otimes t_j).$ Pour tout~$(i,j)$, on note
$S_{ij}$ l'ouvert~$D(\sigma_{ij})$. Il résulte de la définition de~$\Theta$ que pour tout~$(i,i',j)$
la restriction de~$s_i$ à~$S_{i'j}$ est «égale»
à~$\sigma_{ij}$, et que pour tout~$(i,j,j')$ la restriction de~$t_j$ à~$S_{ij'}$ est «égale»
à~$\sigma_{ij}$. 

\medskip
Fixons~$(i,j)$. Le produit tensoriel~$s_i\otimes t_j$ est inversible sur~$S_{ij}$. Il existe donc un unique
isomorphisme~$\ell_{ij} \colon \sch N|_{S_{ij}}\simeq (\sch L\otimes \sch M)|_{S_{ij}}$ tel que~$\ell_{ij}(\sigma_{ij})=s_i\otimes t_j$. 

\medskip
Pour tout~$(i,j,i',j')$ on a les égalités
$$\frac{s_i\otimes t_j}{s_{i'}\otimes t_{j'}}=\frac{\sigma_{ij}\otimes \sigma_{ij}}{\sigma_{i'j}\otimes \sigma_{ij'}}
=\frac{\sigma_{ij}\otimes \sigma_{ij}}{\sigma_{ij}\otimes \sigma_{i'j'}}=\frac{\sigma_{ij}}{\sigma_{i'j'}}$$
dans l'anneau des fonctions de~$S_{ij}\cap S_{i'j'}\subset S_{ij}$, ce qui montre que les restrictions de~$\ell_{ij}$ et~$\ell_{i'j'}$
à~$S_{ij}\cap S_{i'j'}$ coïncident. Les~$\ell_{ij}$ se recollent donc en un isomorphisme
$\ell\colon \sch N\simeq \sch L\otimes \sch M$ ; il reste à calculer~$\ell(\sigma_{ij})$ pour tout~$(i,j)$. 

\medskip
Fixons~$(i,j)$. Soit~$(i',j')\in {\bf n}\times {\bf m}$. On a sur l'ouvert~$S_{i'j'}$ les égalités

$$\ell(\sigma_{ij})=\frac{\sigma_{ij}}{\sigma_{i'j'}}\ell(\sigma_{i'j'})=\frac{\sigma_{ij}}{\sigma_{i'j'}}s_{i'}\otimes t_{j'}
=\frac{\sigma_{ij}}{\sigma_{i'j'}}\sigma_{i'j'}\otimes \sigma_{i'j'}$$
$$=\sigma_{ij}\otimes \sigma_{i'j'}=\sigma_{ij'}\otimes \sigma_{i'j}=s_i\otimes t_j.$$
Il s'ensuit que~$\ell(\sigma_{ij})=s_i\otimes t_j$ sur~$S$ tout entier, ce qui termine la preuve
que~$(\sch N, (\sigma_{ij}))\simeq (\sch L\otimes \sch M, s_i \otimes t_j)$ et achève la démonstration
de la proposition.~$\Box$ 

\deux{def-veronese}
{\bf Le plongement de Veronese}. Soit~$n\in \NN$ et soit
$$\psi \colon \PP^n_A\times_A \PP^n_A\hookrightarrow \PP^{n^2+2n}_A$$ le plongement 
de Segre (définie en~\ref{def-mor-segre}, voir aussi la proposition~\ref{prop-mor-segre}). Le {\em morphisme de Veronese}
$\chi  \colon \PP^n_A\to \PP^{n^2+2n}_A$
est la flèche composée
$$\xymatrix{
{\PP^n_A}\ar[rr]^(0.4){({\rm Id}, {\rm Id})}&&{\PP^n_A\times_A\PP^n_A}\ar[r]^\psi
&\PP^{n^2+2n}_A}$$
(il est donné, en termes naïfs, par la formule~$[s_0:\ldots :s_n]\mapsto [s_is_j]_{ij}$). 

\deux{prop-veronese}
{\bf Proposition.}
{\em Le morphisme de Veronese
$\chi \colon \PP^n_A\to \PP^{n^2+2n}_A$
s'identifie à l'immersion fermée
$$\proj A[\Sigma_{ij}]/[(\Sigma_{ij}\Sigma_{i'j'}-\Sigma_{ij'}\Sigma_{i'j})_{i\neq i', j\neq j'}, (\Sigma_{ij}-\Sigma_{ji})_{i\neq j}]\hookrightarrow \PP^{n^2+2n}.$$}

\medskip
{\em Démonstration.}
Soit~$S$ un~$A$-schéma. L'application~$\isol n S\to \isol n S\times \isol n S$ induite par~$({\rm Id}, {\rm Id})\colon \PP^n_A\to \PP^n_A\times_A\PP^n_A$ est simplement la diagonale~$\lambda\mapsto (\lambda, \lambda)$, qui identifie
$\isol n S$ au sous-ensemble~$\Delta$ de~$\isol n S\times \isol n S$ constitué des couples dont les deux composantes sont égales. 

\medskip
{\em Nous reprenons les notations~$\Phi, \Theta, \mathsf D$ de~\ref{def-mor-segre} et de la démonstration de la
proposition.~\ref{prop-mor-segre} -- notez simplement que maintenant~$m=n$}. Soit~$\mathsf E$ le sous-ensemble
de~$\isol {n^2+2n}S$ constitué des classes d'objets
$(\sch N, (\sigma_{ij}))$ tels que l'on ait
$$\sigma_{ij}\otimes \sigma_{i'j'}=\sigma_{ij'}\otimes\sigma_{i'j}$$ pour tout~$(i,i',j,j')$ 
et~$\sigma_{ij}=\sigma_{ji}$ pour tout~$(i,j)$ avec~$i\neq j$ (notez que si~$i=j$, l'égalité
est automatiquement vérifiée).  En vertu
du lemme~\ref{vi-pn-intuitif}, il suffit de démontrer que 
l'application~$\Phi \colon \isol n S\times \isol n S\to \isol {n^2+2n}S$ induit une bijection~$\Delta \simeq \mathsf E$. 
Il est immédiat que~$\Phi(\Delta)\subset \mathsf E$, et l'on sait d'après la preuve de la proposition~\ref{prop-mor-segre}
que~$\Phi$ induit une bijection~$\ml n S\simeq \mathsf D$ de réciproque~$\Theta$. Il suffit donc pour conclure de 
vérifier que~$\Theta(\mathsf E)\subset \Delta$. 

\medskip
Soit donc~$(\sch N, (\sigma_{ij}))$ un objet de~$\ml {n^2+2n} S$
tel que 
$$\sigma_{ij}\otimes \sigma_{i'j'}=\sigma_{ij'}\otimes\sigma_{i'j}$$ pour tout~$(i,i',j,j')$ 
et~$\sigma_{ij}=\sigma_{ji}$ pour tout~$(i,j)$. Pour tout~$(i,j)$, on note~$S_{ij}$
l'ouvert~$D(\sigma_{ij})$ ; il résulte de nos hypothèses que~$S_{ij}=S_{ji}$ pour tout~$(i,j)$. 
Posons
$$(\;(\sch L, (s_i)), (\sch M, (t_j))\;)=\Theta(\sch N, (\sigma_{ij})).$$ Il s'agit de prouver
que~$(\sch L, (s_i))\simeq (\sch M, (t_j))$. On rappelle que pour tout~$(i, i',j)$, 
la restriction de~$s_i$ à~$S_{i'j}$ est «égale»
à~$\sigma_{ij}$, et que pour tout~$(i, j, j')$ la restriction de~$t_j$ à~$S_{ij'}$ est «égale»
à~$\sigma_{ij}$. 

\medskip
Soit~$(i,j)$ un couple d'indices. La section~$s_i$ de~$\sch L$ est inversible sur~$S_{ij}$. Comme~$S_{ij}=S_{ji}$, 
il en va de même de la section~$t_i$ de~$\sch M$. En conséquence, il existe un
unique isomorphisme~$\ell_{ij}\colon \sch L|_{S_{ij}}\simeq \sch M|_{S_{ij}}$
tel que~$\ell_{ij}(s_i)=t_i$. 

\medskip
Pour tout~$(i,i', j, j')$ on a dans l'anneau~$\sch O_S(S_{ij}\cap S_{i'j'})$
les égalités
$$\underbrace{\frac {s_i}{s_{i'}}=\frac {\sigma_{ij}}{\sigma_{i'j}}}_{\text{calcul~effectué~dans}\;S_{ij}}
=\;\;\;\;\frac{\sigma_{ij'}}{\sigma_{i'j'}}
\;\;\;\;=\underbrace{\frac{\sigma_{j'i}}{\sigma_{j'i'}}=\frac {t_i}{t_{i'}}}_{\text{calcul~effectué~dans}\;S_{i'j'}=S_{j'i'}},$$
ce qui montre que les restrictions de~$\ell_{ij}$ et~$\ell_{i'j'}$ à l'ouvert~$S_{ij}\cap S_{i'j'}$ coïncident. 
Il en résulte que la famille des~$(\ell_{ij})$ se recolle en un isomorphisme~$\ell \colon \sch L \simeq \sch M.$

\medskip
Soit~$i\in\{0, \ldots, n\}$ et soient~$(i', j)$ deux entiers compris entre~$0$ et~$n$. On a sur l'ouvert~$S_{i'j}$
les égalités
$$\ell(s_i)=\frac {s_i}{s_{i'}}\ell(s_{i'})=
\underbrace{\frac {s_i}{s_{i'}}t_{i'}=\frac{\sigma_{ij}}{\sigma_{i'j}}t_{i'}}_
{\text{calcul~effectué~dans}\;S_{ij}}=
\;\;\underbrace{\frac{\sigma_{ji}}{\sigma_{ji'}}t_{i'}=\frac{\sigma_{ji}}{\sigma_{ji'}}\sigma_{ji'}=\sigma_{ji}=t_i}_
{\text{calculs~effectués~dans}\;S_{ij}=S_{ji}}.$$

Il s'ensuit que~$\ell(s_i)=t_i$ sur~$S$ tout entier, ce qui prouve
que~$(\sch L, (s_i))$ est isomorphe à~$(\sch M, (t_j))$
et achève la démonstration.~$\Box$ 

\trois{coro-imm-ferme-pn}
{\bf Corollaire}. 
{\em Le morphisme~$({\rm Id}, {\rm Id})\colon \PP^n_A\to \PP^n_A\times_A\PP^n_A$
est une immersion fermée.}

\medskip
{\em Démonstration.}
Notons~$\delta$ le morphisme en question. La proposition~\ref{prop-veronese}
assure que~$\psi \circ \delta$ est une immersion fermée, et~$\psi$ est elle-même une immersion fermée
d'après la proposition~\ref{prop-mor-segre}.  Il résulte alors de~\ref{fact-imm-ferm}
que~$\delta$ est une immersion fermée.~$\Box$ 

\trois{exemple-p1-p1}
{\em Exemple.} Nous allons décliner les propositions~\ref{prop-mor-segre}
et~\ref{prop-veronese} lorsque~$n=m=1$ ; notez que dans ce cas~$n^2+2n=3$. 
Identifions~$\PP^3_A$ à~$\proj A[\Sigma_{00}, \Sigma_{01}, \Sigma_{10}, \Sigma_{1,1}]$.

\medskip
$\bullet$ Le plongement de Segre identifie~$\PP^1_A\times_A \PP^1_A$ à la «quadrique»
de~$\PP^3_A$ définie par l'idéal homogène~$(\Sigma_{00}\Sigma_{11}-\Sigma_{01}\Sigma_{10})$.

$\bullet$
Le plongement de Veronese identifie~$\PP^1_A$ au sous-schéma fermé de~$\PP^3_A$
défini par l'idéal homogène~$(\Sigma_{00}\Sigma_{11}-\Sigma_{01}\Sigma_{10}, \Sigma_{01}-\Sigma_{10})$. 

\section{Séparation et propreté}
\markboth{Schémas projectifs}{Séparation et propreté}

\subsection*{Morphismes séparés}

\deux{intro-separation}
La notion naïve (purement topologique) de séparation n'a guère d'intérêt
en théorie des schémas, faute d'être suffisamment discriminante : en effet, comme
on a eu l'occasion de le voir, les schémas ne sont presque jamais topologiquement séparés.

\trois{annonce-exist-separation}
Il existe toutefois, comme on va le voir, une «bonne»
notion de séparation en théorie des schémas, qui est conforme à l'intuition
-- par exemple, si~$k$ est un corps, les~$k$-schémas~$\Aff^n_k$ et~$\PP^n_k$
sont séparés, mais la droite affine avec origine dédoublée~$\DD_k$ vue aux~\ref{droite-dedoubl}
{\em et sq.}
ne l'est pas. 

\trois{separation-topologie} 
On peut donner en topologie deux définitions d'un espace séparé (leur équivalence est immédiate, nous vous
laissons la vérifier). 

\medskip
i) Un espace topologique~$X$ est séparé si pour tout couple~$(x,y)$ de points de~$X$ avec~$x\neq y$,
il existe un voisinage ouvert~$U$ de~$x$ dans~$X$ et un voisinage ouvert~$V$ de~$y$ dans~$X$ tels que~$U\cap V=\varnothing$. 

ii) Un espace topologique~$X$ est séparé si
la diagonale
$\{(x,x)\}_{x\in X}$ est un sous-ensemble fermé de~$X\times X$. 

\medskip
Pour ce qui nous intéresse ici, la définition~ii) est meilleure : comme nous le verrons, elle se décalque
très naturellement en géométrie
algébrique et fournit
la bonne notion de séparation dans ce contexte -- alors qu'à notre connaissance, 
il n'existe pas de façon pertinente de «schématiser»
la définition~i). 
 
\deux{def-immersion}
{\bf Définition.}
Soit~$\phi \colon Y\to X$ un morphisme de schémas. On dit que
$\phi$ est une
{\em immersion}
s'il existe un ouvert~$\Omega$ de~$X$ tel que~$\phi$ induise une immersion
fermée~$Y\hookrightarrow \Omega$.

\trois{exemple-immersion}
{\em Exemples.}
Il résulte immédiatement de la définition que les immersions ouvertes et les immersions fermées
sont des cas particuliers d'immersions (la terminologie choisie est donc cohérente). Il n'est pas difficile
de construire des immersions qui ne soient ni ouvertes ni fermées. Donnons-nous par exemple un corps~$k$,
soit~$\Omega$ l'ouvert~$D(S)$ de~$\Aff^2_k=\spec k[S,T]$ et soit~$Y$ le sous-schéma fermé
de~$\Omega$ défini par l'idéal~$(T)$. Par définition, la flèche composée~$Y\hookrightarrow \Omega \hookrightarrow \Aff^2_k$
est une immersion, mais elle n'est ni ouverte ni fermée, car son image est égale à~$D(S)\cap V(T)$ et n'est ni ouverte ni 
fermée. 

\trois{immersion-monomor}
Soit~$\phi \colon Y\to X$ un morphisme de schémas. Si~$\phi$ est une immersion
alors~$Y\times_XY\simeq Y$. En effet, il existe par hypothèse un ouvert~$\Omega$
de~$Y$ tel que~$\phi$ se factorise par une immersion fermée~$Y\hookrightarrow \Omega$. 
On sait que~$Y\times_XY$ est alors égal à~$Y\times_\Omega Y$, et ce dernier s'identifie à~$Y$
en vertu de~\ref{changebase-immf}.  

\trois{prop-gen-immersion}
On vérifie sans peine que la composée de deux immersions est une immersion, 
et que le fait d'être une immersion est stable par changement de base. 

\trois{lemme-immer-fermetop}
{\bf Lemme}. 
{\em Soit~$\phi \colon Y\to X$ un morphisme de schémas. On suppose
que~$\phi$ est une immersion. Pour que~$\phi$ soit une immersion fermée, 
il faut et il suffit que~$\phi(Y)$ soit un fermé de~$X$.}

\medskip
{\em Démonstration}. Si~$\phi$ est une immersion fermée, $\phi(Y)$ est un fermé de~$X$. 
Réciproquement, supposons que~$\phi(Y)$ soit un fermé de~$X$ et soit~$U$ son ouvert complémentaire. Comme~$\phi$ est une immersion, 
il existe un ouvert~$\Omega$ de~$X$ tel que~$\phi$ induise une immersion fermée~$Y\hookrightarrow \Omega$. 
En conséquence, $$Y\times_X\Omega \to \Omega =\phi^{-1}(\Omega)\to \Omega=Y\to \Omega$$ est une immersion fermée. Par ailleurs,
$$Y\times_XU\to U=\phi^{-1}(U)\to U=\varnothing \to U$$ est aussi une immersion fermée. Comme~$\phi(Y)\subset \Omega$, les ouverts~$U$ et~$\Omega$
recouvrent~$Y$ ; puisque le fait d'être une immersion fermée est une propriété locale sur le but, $\phi$ est une immersion fermée.~$\Box$ 

\trois{remm-ouverttop-immer}
{\em Remarque.}
L'assertion analogue pour les immersions ouvertes est {\em fausse} : par exemple, si~$X$ est un schéma non réduit, $X_{\rm red}\hookrightarrow X$
est une immersion (fermée) dont l'image est~$X$, mais ce n'est pas une immersion ouverte, car sinon ce serait un isomorphisme et~$X$ serait réduit. 

\deux{immer-diagonale}
{\bf Définition.} Soit~$X$ un schéma et soit~$Y$ un~$X$-schéma. La {\em diagonale}
(du morphisme~$Y\to X$, ou du~$X$-schéma~$Y$, ou de~$Y$ au-dessus de~$X$)
est la flèche~
$$\delta \colon \xymatrix{Y\ar[rr]^(0.4){({\rm Id}, {\rm Id})}&&{Y\times_XY}}.$$

\trois{diagonale-cas-affine}
Supposons que~$X$ est le spectre d'un anneau~$A$, et~$Y$ celui d'une~$A$-algèbre~$B$. 
La diagonale~$\delta$ est alors
induite par le morphisme d'anneaux~$B\otimes_AB\to B$
correspondant au couple~$({\rm Id}_B, {\rm Id}_B)$, qui n'est autre que
la «multiplication»~$b\otimes \beta \mapsto b\beta$. Celle-ci est manifestement surjective ; en conséquence,
$\delta$ est une immersion fermée. 

\trois{diagonale-cas-gen}
On ne suppose plus que~$X$ et~$Y$ sont affines. Soit~$\Omega$ la réunion des ouverts
de~$Y\times_XY$ qui sont de la forme~$V\times_UV$ où~$V$ est un ouvert affine de~$Y$ et~$U$
un ouvert affine de~$X$ contenant l'image de~$V$. 

\medskip
Donnons-nous un tel couple~$(U,V)$. Comme~$V\times_UV$ est l'intersection des images
réciproques de~$V$ par les deux projections de~$Y\times_XY$ sur~$Y$, le produit fibré
$$Y\times_{Y\times_XY}(V\times_UV)=\delta^{-1}(V\times_UV)\to V\times_UV$$
est simplement 
la diagonale~$V\to V\times_UV$ du morphisme~$V\to U$. En vertu de~\ref{diagonale-cas-affine}, 
c'est une immersion fermée. 

\medskip
Puisqu'être une immersion fermée est une propriété locale sur le but d'un morphisme, $\delta^{-1}(\Omega)\to \Omega$ 
est une immersion fermée. Par ailleurs, soit~$z\in Y\times_XY$. Par définition, les images de~$\delta(z)$ 
par les deux projections sur~$Y$ sont égales à un même point~$y$ ; soit~$x$ l'image de~$y$ sur~$X$. Choisissons
un voisinage 
affine~$U$ de~$x$ dans~$X$ et un voisinage affine~$V$ de~$y$ dans $Y\times_XU$ ; 
par construction, le point~$\delta(z)$ appartient
à~$V\times_UV \subset \Omega$. Ainsi, $\delta^{-1}(\Omega)=Y$. La diagonale~$\delta$
se factorise donc par une immersion fermée~$Y\hookrightarrow \Omega$ ; par conséquent, 
{\em $\delta$ est une immersion.}

\deux{def-mor-separe}
{\bf Définition.}
Soit~$Y\to X$ un morphisme de schémas. On dit que $Y\to X$ est {\em séparé}
si la diagonale~$Y\hookrightarrow Y\times_XY$ (qui est une immersion d'après le~\ref{diagonale-cas-gen}
ci-dessus) est une immersion fermée. On dit parfois aussi
que~$Y$ est un {\em $X$-schéma séparé}, ou que~$Y$ est {\em séparé sur~$X$}. 

\medskip
Un schéma est dit
{\em séparé}
s'il est séparé sur~$\spec \ZZ$. 

\deux{propriete-separation}
Soit~$Y\to X$ un morphisme de schémas et soit~$\delta \colon Y\hookrightarrow Y\times_XY$ 
l'immersion diagonale. 

\trois{affine-separe}
Si~$Y$ et~$X$ sont affines, il  découle de~\ref{diagonale-cas-affine}
que~$Y\to X$ est séparé. 

\trois{separe-test-topologique}
En général, comme~$\delta$ est une immersion, on déduit du lemme~\ref{lemme-immer-fermetop}
que~$Y\to X$ est séparé si et seulement si~$\delta(Y)$ est fermé
dans~$Y\times_XY$. 

\trois{separe-prop-locale}
Soit~$U$ un ouvert de~$X$ et soit~$V$ son image réciproque sur~$Y$. La flèche~$Y\times_{Y\times_XY} (V\times_UV)$ 
s'identifie à l'immersion diagonale~$V\hookrightarrow V\times_UV$ (\ref{diagonale-cas-gen}). Si~$Y\to X$ est séparé,
cette flèche est donc une immersion fermée, ce qui veut dire que~$V$ est séparé sur~$U$. 

Si~$X$ possède un recouvrement ouvert~$(U_i)$ tel que~$Y\times_X U_i\to U_i$ soit séparé pour tout~$i$, il résulte
de ce qui précède et du caractère local (au but) de la propriété d'être une immersion fermée que~$\delta$ est une immersion fermée, 
et donc que~$Y\to X$ est séparé. {\em La séparation est donc une propriété locale sur le but}. 

\trois{mor-affine-separe}
Il résulte de~\ref{affine-separe}
et~\ref{separe-prop-locale}
que si le morphisme~$Y\to X$ est affine, il est séparé. 

\trois{separe-changebase}
Soit~$X'\to X$ un morphisme ; posons~$Y'=Y\times_X X'$. On
vérifie sans peine que la diagonale~$Y'\hookrightarrow Y'\times_{X'}Y'$
s'identifie à la flèche canonique
$$Y\times_{Y\times_XY}(Y'\times_{X'}Y')\to Y'\times_{X'}Y'.$$
C'est donc une immersion fermée dès que~$\delta$ est une immersion fermée. Autrement dit,
si~$Y\to X$ est séparé alors~$Y'\to X'$ est séparé. 

\trois{separe-composition}
Soit~$Z\to Y$ un morphisme de schémas. Supposons que~$Z\to Y$ et~$Y\to X$
soient séparés. Nous allons montrer qu'il en va alors de même de la flèche
composée~$Z\to X$. 

\medskip
La flèche diagonale~$Z\to Z\times_XZ$ est composée de~$Z\to Z\times_YZ$, qui est une immersion fermée
par hypothèse, et de~$Z\times_YZ\to Z\times_XZ$. Il suffit dès lors
pour conclure de montrer que
$$Z\times_YZ\to Z\times_XZ$$ est une immersion fermée. Nous allons montrer que cette dernière flèche
s'identifie naturellement à
$$Y\times_{Y\times_X Y}(Z\times_XZ)\to Z\times_X Z,$$
ce qui permettra de conclure puisque~$Y\to Y\times_XY$ est par hypothèse
une immersion fermée.  

\medskip
On cherche donc à établir que le carré
commutatif
$$\xymatrix{
{Z\times_YZ}\ar[d]_p\ar[r]^q&{Z\times_XZ}\ar[d]^\pi
\\
Y\ar[r]_(0.4)\delta&{Y\times_XY}}$$
est cartésien (c'est-à-dire qu'il identifie
le terme en haut à gauche au produit
fibré des trois autres). Face à ce genre de problème, 
{\em le seul 
réflexe sain est d'invoquer le lemme de Yoneda pour 
se ramener à l'assertion ensembliste correspondante}
(pour un exemple de raisonnement détaillé de ce type, {\em cf.}~\ref{red-cas-ens}). 

\medskip
On suppose donc pour un instant que le carré ci-dessus vit dans la catégorie des ensembles, et nous allons
montrer qu'il est cartésien. Appelons~$g$ la flèche~$Z\to Y$, et~$f$ la flèche~$Y\to X$. 
Par définition, le produit~$Y\times_XY$ est l'ensemble des couples~$(y,y')\in Y^2$ tels que~$f(y)=f(y')$, le produit~$Z\times_XZ$ est l'ensemble
des couples~$(z,z')\in Z^2$ tels que~$f(g(z))=f(g(z'))$, et le produit~$Z\times_YZ$ est l'ensemble
des couples~$(z,z')\in Z^2$ tels que~$g(z))=g(z')$. Les flèches du diagramme sont données par les formules suivantes : 

\medskip
$\bullet$ $q(z,z')=(z,z')$ ; 

$\bullet$ $\pi(z,z')=(g(z), g(z'))$ ; 

$\bullet$ $p(z,z')=g(z)=g(z')$ ; 

$\bullet$ $\delta(y)=(y,y)$. 
 

\medskip
Il s'agit maintenant
de s'assurer que pour tout triplet~$(y, z, z')\in Y\times (Z\times_XZ)$ tel que~$\delta(y)=\pi(z,z')$ il existe un 
unique élément de~$Z\times_Y Z$ dont l'image par~$p$ est égale à~$y$ et l'image par~$q$ à~$(z,z')$. L'unicité
est claire : étant donnée la formule qui définit~$q$, si un tel élément existe, ce ne peut être que~$(z,z')$. 
Il reste à s'assurer que celui-ci convient. Mais l'égalité~$\delta(y)=\pi(z,z')$ signifie que~$g(z)=y$ et~$g(z')=y$, 
ce qui signifie précisément que~$(z,z')\in Z\times_YZ$ et que~$p(z,z')=y$ ; on a de plus~$q(z,z')=(z,z')$, 
ce qui termine la démonstration. 

\deux{separe-exemple-contrex}
{\bf Exemples et contre-exemples}. 

\trois{rappel-affine-separe}
Nous avons déjà vu que les morphismes affines sont séparés
(\ref{mor-affine-separe}). En particulier, pour tout schéma~$X$ et tout entier~$n$, 
le schéma~$\Aff^n_X$ est séparé sur~$X$. 

\trois{imm-separe}
Toute immersion est séparée, puisque sa diagonale est un {\em isomorphisme}
d'après~\ref{immersion-monomor}.  

\trois{projectif-separe}
Pour tout schéma~$X$ et tout entier~$n$, le schéma~$\PP^n_X=\PP^n_{\ZZ}\times_{\ZZ}X$ est séparé sur~$X$. 
En effet, en vertu de~\ref{separe-changebase}, il suffit de traiter le cas où~$X$
est égal à~$\spec \ZZ$, auquel cas
c'est une conséquence directe du corollaire~\ref{coro-imm-ferme-pn}
(qui établit d'ailleurs en fait directement la séparation de~$\PP^n_A$ sur~$\spec A$ pour tout anneau~$A$). 

\trois{mor-projectif}
Un morphisme de schémas~$Y\to X$ est dit
{\em quasi-projectif}
(resp. {\em projectif})
si pour tout~$x\in X$ il existe un voisinage ouvert~$U$ de~$x$ dans~$X$
et un entier~$n$ tel que le morphisme~$Y\times_XU\to U$ se factorise par une
immersion (resp. une immersion fermée)~$Y\times_XU\hookrightarrow \PP^n_U$.
Il résulte de~\ref{imm-separe}, \ref{projectif-separe},
~\ref{separe-composition}
et~\ref{separe-prop-locale}
que tout morphisme quasi-projectif est séparé ; c'est {\em a fortiori}
le cas de tout morphisme projectif. 

\trois{dedouble-nonsepare}
Soit~$k$ un corps et soit~$\DD_k$ la droite affine avec origine dédoublée construite aux~\ref{droite-dedoubl}
{\em et sq.} Elle n'est pas séparée sur~$k$ ; nous allons esquisser une démonstration de ce fait. 

\medskip
Redonnons brièvement ici la description de~$\DD_k$. Elle est réunion de deux ouverts affines~$X=\spec k[S]$ et~$Y=\spec k[T]$. Leur intersection
est égale à~$D(S)=\spec k[S,S^{-1}]$ en tant qu'ouvert de~$X$, et à~$D(T)=\spec k[T, T^{-1}]$ en tant qu'ouvert de~$Y$. L'isomorphisme entre ces deux identifications
est induit par l'isomorphisme de~$k$-algèbres
$$k[S,S^{-1}]\simeq k[T, T^{-1}], \;S \mapsto T.$$

\medskip
On dispose d'un morphisme naturel~$\pi \colon X\coprod Y\to \DD_k$. 
Complétons le diagramme 
$$\xymatrix{
&{ (X\coprod Y)\times_k (X\coprod Y)}\ar[d]\\
\DD_k\ar[r]&{\DD_k\times_k\DD_k}}$$ en un carré cartésien

$$\xymatrix{
\Delta\ar[r]\ar[d]&{ (X\coprod Y)\times_k (X\coprod Y)}\ar[d]\\
\DD_k\ar[r]&{\DD_k\times_k\DD_k}}$$ 
dont la flèche
horizontale supérieure est une immersion, et une immersion fermée si~$\DD_k$ est séparée.

\medskip
Moralement, on peut penser à~$\Delta$ comme au {\em graphe de la relation
d'équivalence qui a permis de définir~$\DD_k$ à partir de~$X\coprod Y$}, mais nous allons
maintenant en donner
une description rigoureuse. Commençons par observer que nous sommes dans la même situation formelle qu'au~\ref{separe-composition} ; il s'ensuit que l'immersion
$$\Delta\hookrightarrow \left(X\coprod Y\right)\times_k \left(X\coprod Y\right)$$
s'identifie à
$$\left(X\coprod Y\right)\times_{\DD_k}\left(X\coprod Y\right)\to \left(X\coprod Y\right)\times_k \left(X\coprod Y\right).$$
La source et le but de cette flèche admettent chacun une décomposition en quatre ouverts disjoints, et cette flèche préserve
ces décompositions. 
On se retrouve donc avec quatre immersions différentes à considérer. 

\medskip
a) {\em L'immersion~$X\times_{\DD_k}X\to X\times_k X$.} Comme~$X$ est un ouvert de~$\DD_k$, 
le produit fibré~$X\times_{\DD_k}X$ s'identifie à~$X$, et l'immersion étudiée est donc l'immersion
diagonale~$X\hookrightarrow X\times_kX$, qui est fermée puisque~$X=\spec k[T]$ est affine (on peut évidemment
la calculer directement et voir qu'elle s'identifie à~$\spec k[T_1,T_2]/(T_1-T_2)\hookrightarrow \spec k[T_1,T_2]$). 

\medskip
b) {\em L'immersion~$Y\times_{\DD_k}Y\to Y\times_k Y$.} Pour la même raison, c'est l'immersion
diagonale~$Y\hookrightarrow Y\times_k Y$, et elle est fermée. 

\medskip
c) {\em L'immersion~$X\times_{\DD_k}Y\to X\times_kY$}. Comme~$X$ et~$Y$ sont des ouverts de~$\DD_k$, le terme de
gauche est l'intersection de~$X$ et~$Y$ dans~$\DD_k$, laquelle est s'identifie à~$\spec k[S, S^{-1}]\simeq \spec k[T,T^{-1}]$ 
(l'isomorphisme envoyant~$T$ sur~$S$). Quant au produit fibré~$X\times_k Y$, c'est de façon naturelle
le spectre de~$k[S,T]$. La flèche entre les deux est donnée par le morphisme d'algèbres~$k[S,T]\to k[S,S^{-1}]$ qui envoie
$S$ et~$T$ sur~$S$. {\em Il n'est pas surjectif et cette flèche n'est donc pas une immersion fermée}. 

Donnons quelques
précisions. Le morphisme~$k[S,T]\to k[S,S^{-1}]$
est la composée de~$k[S,T]\hookrightarrow k[S, S^{-1},T]$ et de la surjection~$k[S,S^{-1},T]\to k[S,S^{-1}]$ qui envoie~$T$ sur~$S$
et a pour noyau~$(S-T)$.
L'immersion $X\times_{\DD_k}Y\to X\times_kY$ est donc égale à la composée de l'immersion fermée~$V(S-T)\cap D(S)\hookrightarrow D(S)$ 
(où le fermé~$V(S-T)\cap D(S)$ de~$D(S)$ est muni de sa structure réduite, d'anneau associé~$k[S,S^{-1}]$), et de l'immersion ouverte~$D(S)\hookrightarrow \spec k[S,T]$. 
On voit bien que son image n'est pas fermée : c'est~$V(S-T)\cap D(S)=V(S-T)\setminus\{(0,0)\}$, c'est-à-dire la diagonale {\em épointée}
(cela traduit le fait qu'on a identifié chaque point de~$X$  {\em à l'exception de l'origine}
au point correspondant de~$Y$). 

\medskip
d) {\em L'immersion~$Y\times_{\DD_k}X\to Y\times_kX$}. Elle se décrit exactement comme l'immersion considérée au~c) ; elle a également pour image
la diagonale épointée, et n'est donc pas fermée.

\medskip
On voit donc qu'en raison de~c) et~d), l'immersion
$$\Delta\hookrightarrow \left(X\coprod Y\right)\times_k \left(X\coprod Y\right)$$
n'est pas fermée ; il s'ensuit que~$\DD_k$ n'est pas séparée. 

\trois{comment-p1-separe}
{\em Remarque}. On sait que le~$k$-schéma~$\PP^1_k$ est séparé (\ref{projectif-separe}).
On peut par ailleurs en donner une construction par recollement, analogue à celle utilisée pour
définir~$\DD_k$ (\ref{droite-proj-recoll}
{\em et sq.}) : si l'on reprend les notations ci-dessus, la seule
différence avec le cas de la droite à origine dédoublée
réside dans le fait que l'isomorphisme entre les deux identifications
$$X\cap Y\simeq \spec k[S,S^{-1}]\;\text{et}\;
X\cap Y\simeq \spec k[T, T^{-1}]$$ est induit par le morphisme d'algèbres qui envoie~$T$ sur~$S^{-1}$ (et non pas~$S$).  
Supposons que l'on cherche, dans ce nouveau contexte, à décrire l'immersion 
$$\Delta\hookrightarrow \left(X\coprod Y\right)\times_k \left(X\coprod Y\right).$$
Tout se passe comme ci-dessus jusqu'au point~b) inclus, mais une différence fondamentale apparaît au point~c) : 
l'immersion $X\times_{\PP^1_k}Y\to X\times_kY$ est alors donnée par le morphisme d'algèbres de~$k[S,T]$ vers~$k[S,S^{-1}]$
qui envoie~$T$ sur~$S^{-1}$ et qui est {\em surjectif}, de noyau~$(TS-1)$ ; c'est donc une immersion fermée, d'image
l'hyperbole~$V(ST-1)$ (et qui induit la structure réduite sur celle-ci, d'anneau associé
$k[S,S^{-1}]$). 

\deux{annonce-lemme-interaff}
Nous allons terminer ces considérations sur la séparation par un lemme facile qui a son intérêt, même si nous ne nous
en servirons pas dans la suite. 

\deux{lemme-interaff}
{\bf Lemme.}
{\em Soit~$A$ un anneau et soit~$X$ un~$A$-schéma séparé ; soient~$U$ et~$V$ deux ouverts affines de~$X$. 
L'intersection~$U\cap V$ est affine.}

\medskip
{\em Démonstration.}
Comme~$X$ est séparé sur~$A$, e morphisme diagonal~$\delta \colon X\hookrightarrow X\times_AX$ est une immersion fermée. 
Soient~$p$ et~$q$ les deux projections de~$X\times_AX$ vers~$X$. L'ouvert~$U\times_AV$ de~$X\times_AX$ est affine
puisque~$U, V$ et~$\spec A$ le sont, et il est égal à~$p^{-1}(U)\cap q^{-1}(V)$. Il s'ensuit
que~$\delta^{-1}(U\times_AV)=U\cap V$, et~$\delta$ induit donc une
immersion fermée~$U\cap V\hookrightarrow U\times_AV$. 
Comme~$U\times_AV$ est affine, on en déduit que~$U\cap V$ est affine.~$\Box$ 

\deux{contrex-interaffine}
{\bf Remarque}. Donnons un contre-exemple
au lemme~\ref{lemme-interaff}
ci-dessus lorsque l'hypothèse de séparation
n'est pas satisfaite. 
Soit~$k$ un corps.
On peut définir, par un procédé analogue à celui utilisé pour définir~$\DD_k$
que nous avons rappelé au~\ref{dedouble-nonsepare}
ci-dessus, le {\em plan affine avec origine dédoublée}.
C'est un~$k$-schéma qui n'est pas séparé. Il est réunion de deux ouverts ouverts affines~$X$ et~$Y$. 
Chacun d'eux est isomorphe à~$\Aff^2_k$, et leur intersection s'identifie (comme ouvert de~$X$ aussi bien que de~$Y$
à~$\Aff^2_k\setminus \{(0,0)\}$ ; elle n'est donc pas affine (\ref{ex-non-affine} {\em et sq.}).

\subsection*{Morphismes propres}

\deux{motiv-propre}
Il en va de la compacité comme de la séparation : il ne semble pas raisonnable, au vu
de la la grossièreté de la topologie
de Zariski en général, d'espérer une notion satisfaisante de compacité en géométrie algébrique qui 
soit définissable en termes purement topologique. Songez par exemple que sur un corps~$k$, 
les espaces topologiques~$\Aff^1_k$ et~$\PP^1_k$ sont {\em homéomorphes}
(c'est une conséquence immédiate de~\ref{topologie-p1k}, et du fait que l'ensemble des points
fermés de~$\Aff^1_k$ et celui de~$\PP^1_k$ sont de même cardinal infini) ; or quelque soit
le sens que l'on donne à l'adjectif «compact», la décence exige que~$\Aff^1_k$ ne le soit pas et
que~$\PP^1_k$ le soit. 

\deux{topologie-compact}
Soit~$X$ un espace topologique séparé et localement compact. On démontre que~$X$ est compact si et seulement
si pour tout espace localement compact~$Y$, la projection~$X\times Y\to Y$ est {\em fermée}
(cela signifie que l'image d'un fermé est fermé). L'expérience a montré que c'est cette caractérisation de la compacité
qui se prête le mieux à une transposition dans le monde des schémas -- sous le nom de {\em propreté}. 
Pour pouvoir définir celle-ci, nous allons avoir besoin d'une première notion, celle de morphisme
{\em universellement fermé}. 

\deux{def-univ-ferm}
{\bf Définition}. Un morphisme de schémas~$Y\to X$
est dit {\em universellement fermé}
si pour tout~$X$-schéma $X'$, l'application 
continue~$Y\times_XX'\to X'$ est fermée. 

\deux{exemple-uniferm}
{\bf Exemples et contre-exemples}. 

\trois{fini-uniferme}
Soit~$Y\to X$ un morphisme fini de schémas ; il est universellement fermé. En effet, 
soit~$X'$ un~$X$-schéma. Le morphisme~$Y\times_XX'\to X'$ est fini, 
et est en conséquence fermé (prop.~\ref{morph-fin-ferm}). 

\trois{a1-pas-uniferme}
Soit~$k$ un corps.
Le morphisme~$\Aff^1_k\to \spec k$ n'est pas universellement fermé (notez qu'il 
est par contre fermé, et qu'il n'a pas grand mérite puisque~$\spec k$ est un point). 

En effet, 
la projection~$p\colon \Aff^2_k\to \Aff^1_k$ par rapport à la seconde
variable n'est pas fermée, puisque
l'image par~$p$ de l'hyperbole
$V(T_1T_2-1)$ est l'ouvert~$D(T_2)$, qui n'est pas fermé. 

\deux{prop-base-univferme}
Soit~$Y\to X$ un morphisme de schémas.

\trois{change-base-uniferme}
Si~$Y\to X$ est universellement fermé, alors pour tout~$X$-schéma~$X'$, 
le morphisme~$Y\times_XX'\to X'$ est universellement fermé : c'est une conséquence immédiate
de la définition, qui {\em impose}
la stabilité par changement de base. 

\trois{composition-uniferme}
Soit~$Z\to Y$ un morphisme. Supposons que~$Z\to Y$ et~$Y\to X$ soient universellement fermés ; 
la flèche composée~$Z\to Y\to X$ est alors universellement fermée. 

\medskip
En effet, soit~$X'$ un
$X$-schéma. Posons
$$Y'=Y\times_XX'\;\;\text{et}\;\;Z'=Z\times_XX'=Z\times_YY'.$$
Comme~$Y\to X$ est universellement fermé, $Y'\to X'$ est fermé. Comme~$Z\to Y$ est universellement fermé, 
$Z'\to Y'$ est fermé. Il est immédiat que la composée de deux applications fermées est fermée ; en conséquence, $Z'\to X'$ est fermé,
et~$Z\to X$ est universellement fermée. 

\trois{uniferme-local}
Soit~$Y\to X$ un morphisme. Supposons qu'il existe un recouvrement ouvert~$(U_i)$ de~$X$ tel que~$Y\times_X U_i\to U_i$ 
soit universellement fermé pour tout~$i$. Dans ce cas, $Y\to X$ est universellement fermé. 

En effet, soit~$X'$ un~$X$-schéma. Pour tout~$i$, posons~$U'_i=X'\times_XU_i$ ; la famille~$(U'_i)$ est un recouvrement
ouvert de~$X'$. Pour tout~$i$, le morphisme
$$(Y\times_X U_i)\times_{U_i} U'_i=(Y\times_X X')\times_{X'}U'_i\to U'_i$$ est fermé,
puisque~$Y\times_X U_i\to U_i$ est universellement fermé. Il s'ensuit immédiatement, compte-tenu du fait qu'être fermé est, pour un sous-ensemble
de~$X'$, une propriété locale, que~$Y\times_XX'\to X'$ est fermé, d'où notre assertion. 

\deux{topologie-separe-univferme}
Il est bien connu en topologie générale que si~$\phi$ est une
application continue d'un espace topologique compact~$Y$ vers un espace
topologique {\em séparé}~$Z$ alors~$\phi(Y)$ est une partie fermée de~$Z$. Nous allons
énoncer un avatar de ce résultat dans le monde des schémas. 

\deux{uniferme-image}
{\bf Lemme.}
{\em Soit~$Y\to X$ un morphisme de schémas universellement fermé, 
soit~$Z$ un~$X$-schéma séparé et soit~$\phi \colon Y\to Z$ un~$X$-morphisme. L'image~$\phi(Y)$
est un fermé de~$Z$.}

\medskip
{\em Démonstration}.
Le morphisme~$\phi$ peut s'écrire comme la flèche composée

$$\xymatrix{Y\ar[rrr]^{({\rm Id}, \phi)}&&&{Y\times_X Z}\ar[r]&Z}.$$
Comme~$Y\to X$ est universellement fermé, $Y\times_XZ\to Z$ est fermé. Il suffit
donc pour conclure de s'assurer que~$({\rm Id}, \phi)\colon Y\to Y\times_X Z$ a une image fermée. Nous allons
pour ce faire montrer que le carré commutatif 

$$\xymatrix{Y\ar[d]_\phi\ar[rrr]^{({\rm Id}, \phi)}&&&{Y\times_X Z}\ar[d]^{(\phi\circ p, q)}\\
Z\ar[rrr]^(0.42){\delta}&&&{Z\times_XZ}}$$
(où~$\delta$ est la diagonale de~$Z\to X$, et où~$p$ et~$q$ sont les projections respectives de~$Y\times_XZ$ vers~$Y$
et~$Z$) est {\em cartésien} : sa flèche du bas étant une immersion fermée en vertu de l'hypothèse de séparation faite sur le~$X$-schéma~$Z$, 
il en résultera que~$({\rm Id}, \phi)\colon Y \to Y\times_XZ$ est une immersion fermée, et en particulier a une image fermée,
ce qui achèvera la preuve. 

\medskip
 Là encore, nous nous ramenons grâce au lemme de Yoneda au cas d'un diagramme analogue dans la catégorie des ensembles.  
 Soit~$(z,y,z')\in Z\times Y \times_X Z$ tel que~$\delta(z)=(\phi\circ p, q)(y,z')$, c'est-à-dire tel que~$(z,z)=(\phi(y), z')$, ou encore tel que~$z'=z=\phi(y)$.
 Il s'agit de montrer 
 qu'il existe un unique élément~$\eta\in Y$ tel que~$({\rm Id}, \phi)(\eta)=(y,z')$ (c'est-à-dire tel que~$\eta=y$ et~$\phi(\eta)=z'$)
 et tel que~$\phi(\eta)=z$. Or il est immédiat 
 que~$\eta=y$ est 
 solution du problème, et est la seule.~$\Box$ 
 
\deux{def-mor-propre}
{\bf Définition}. Un morphisme de schémas~$Y\to X$
est dit
{\em propre}
s'il est séparé, de type fini et universellement fermé. On dira également
que $Y$ est un {\em $X$-schéma propre}
ou que~$Y$ est {\em propre sur ~$X$}. 

\trois{prem-pro-morpropres}
Les propriétés, pour un morphisme, d'être séparé, d'être de type fini, et d'être universellement fermé
sont stables par composition, par changement de base, et sont locales sur le but ; il en va donc de même pour la propreté. 

\trois{fini-propre}
{\em Un premier exemple}. 
Tout morphisme fini est de type fini, est affine donc séparé (\ref{mor-affine-separe}), 
et est universellement fermé (\ref{fini-uniferme}). Autrement dit, tout morphisme fini est propre.
En particulier, une immersion fermée est propre (notez que le caractère universellement
fermé  peut être établi directement dans ce cas, alors que pour les morphismes
finis généraux il fait {\em in fine}
appel au lemme de {\em going-up}). 

\deux{theo-prop-propre}
{\bf Théorème.}
{\em Soit~$X$ un schéma et soit~$n$ un entier. Le morphisme~$\PP^n_X\to X$ est propre.}

\medskip
{\em Démonstration}. Comme~$\PP^n_X=\PP^n_{\ZZ}\times_{\ZZ}X$, et comme la propreté est stable
par changement de base (\ref{prem-pro-morpropres}), il suffit de traiter le cas où~$X=\spec \ZZ$. Le morphisme
$\PP^n_{\ZZ}\to \ZZ$ est de type fini, et est séparé ({\em cf}.~\ref{projectif-separe}, ou
directement le corollaire
\ref{coro-imm-ferme-pn}). 

\medskip
Il reste à s'assurer que~$\PP^n_{\ZZ}\to \spec \ZZ$ est universellement fermé, c'est-à-dire
que~$\PP^n_Y\to Y$ est fermé pour tout~$Y$. Le fait, pour une partie d'un schéma~$Y$, d'être
fermée dans~$Y$ est une propriété locale ; il en résulte qu'on peut supposer que~$Y$ est le spectre d'un anneau~$A$.

\trois{preliminaires-vi-ferme}
Soit~$I$ un idéal homogène de~$A[T_0,\ldots, T_n]$.
Nous allons démontrer que l'image de~$V(I)\subset \PP^n_A$ sur~$\spec A$ est
fermée, ce qui permettra de conclure. L'immersion fermée~$\proj A[T_0,\ldots, T_n]/I\hookrightarrow \PP^n_A$
induit un homéomorphisme~$\proj A[T_0,\ldots, T_n]/I\simeq V(I)$. Il suffit donc de vérifier que l'image
de~$\proj A[T_0,\ldots, T_n]/I$
sur~$\spec A$ est fermée ; nous allons plus précisément prouver que son complémentaire~$U$
est ouvert. 

\trois{condition-xu-viferme}
Soit~$x\in \spec A$. Notons~$J(x)$ l'idéal de~$\kappa(x)[T_1,\ldots, T_n]$
engendré par l'image de~$I$. La fibre de~$\proj A[T_0,\ldots, T_n]/I$
en~$x$  s'identifie à
$$\proj (A[T_0,\ldots, T_n]/I)\otimes_A\kappa(x))=\proj \kappa(x)[T_0,\ldots, T_n]/J(x).$$
Dire que~$x\in U$
signifie que la fibre en question est {\em vide}, c'est-à-dire, 
en vertu de~\ref{projb-nonvide}, que tout élément homogène de degré strictement positif de~$\kappa(x)[T_0,\ldots, T_n]/J(x)$
est nilpotent. Comme l'idéal~$(\kappa(x)[T_0,\ldots, T_n]/J(x))\pos$
est engendré par~$(T_0,\ldots, T_n)$, cela revient à demander que les~$T_i$ soient nilpotents dans~$\kappa(x)[T_0,\ldots, T_n]/J(x)$, 
ou encore que~$(\kappa(x)[T_0,\ldots, T_n]/J(x))_d$
soit nul pour un certain~$d$.

\trois{vi-ferme-suite}
Soit~$d$ un entier On pose~$I_d=I\cap A[T_0,\ldots, T_n]_d$, et on définit le~$A$-module~$Q_d$ par la suite exacte
$$0\to I_d\to A[T_0,\ldots,T_n]_d\to Q_d\to 0.$$ Pour tout~$x\in \spec A$,
l'exactitude à droite du produit tensoriel garantit l'exactitude de la suite
$$I_d\otimes_A\kappa(x)\to \kappa(x)(T_0,\ldots, T_n]_d\to Q_d\otimes_A \kappa(x)\to 0,$$ ce qui montre 
que~$Q_d\otimes_A \kappa(x)\simeq (\kappa(x)[T_0,\ldots, T_n]/J(x))_d$. Par conséquent, il découle
de~\ref{preliminaires-vi-ferme}
que
le point~$x$ appartient à~$U$ si et seulement si il existe~$d\in \NN$ tel que~$Q_d\otimes_A \kappa(x)=\{0\}$. 

\trois{vi-ferme-fin}
Soit~$d$ un entier.
Si~$x\in \spec A$, l'espace vectoriel~$Q_d\otimes_A \kappa(x)$ s'identifie à~$\red{Q_d}\otimes \kappa(x)$,
au sens de~\ref{f-kappax} (où~$\red{Q_d}$ désigne le faisceau quasi-cohérent sur~$\spec A$ associé à~$Q_d$). 
Comme~$Q_d$ est de type fini d'après
sa définition, il s'écrit comme un quotient de~$A^m$ pour un certain~$m$, et~$\red{Q_d}$ 
s'écrit dès lors comme un quotient de~$\sch O_{\spec A}^m$. Il résulte alors du corollaire~\ref{coro-naka-geom}
(voir aussi les commentaire qui le suivent en~\ref{comment-zero-ouvert})
que l'ensemble~$V_d$ des points~$x$ de~$\spec A$ tels que~$Q_d\otimes_A\kappa(x)=\red{Q_d}\otimes \kappa(x)=\{0\}$ est
ouvert. 

\medskip
L'ensemble~$U$ étant égal en vertu de~\ref{vi-ferme-suite}
à la réunion des~$V_d$ pour~$d\in \NN$, il est ouvert, ce qui achève la démonstration.~$\Box$ 

\deux{coro-morproj-propre}
{\bf Corollaire}.
{\em Tout morphisme projectif (\ref{mor-projectif})
est propre}.

\medskip
{\em Démonstration}. C'est une conséquence immédiate du théorème~\ref{theo-prop-propre}
ci-dessus, du fait que les immersions fermées sont propres (\ref{fini-propre}), et du bon comportement
de la propreté  à divers égards (\ref{prem-pro-morpropres}).~$\Box$

\subsection*{Un «principe du maximum»
en géométrie algébrique}
\deux{annonce-principe-max}
Nous nous proposons pour terminer ce cours d'établir une variante algébrique du principe du
maximum de la géométrie complexe. Nous aurons besoin du lemme suivant. 

\deux{lemme-schema-integre}
{\bf Lemme.}
{\em Soit~$X$ un schéma irréductible et réduit. L'anneau~$\sch O_X(X)$ est intègre.}

\medskip
{\em Démonstration}. 
Comme~$X\neq \varnothing$ (puisqu'il est irréductible), l'anneau~$\sch O_X(X)$ est non nul
(\ref{locann-non-vide}). Soient~$f$ et~$g$ deux éléments de~$\sch O_X(X)$ tels que~$fg=0$. 
On a alors~$X=V(fg)=V(f)\cup V(g)$, et comme~$X$ est irréductible il vient
$X=V(f)$ ou~$X=V(g)$. Supposons par exemple que~$X=V(f)$. La restriction de~$f$ à tout ouvert affine de~$X$
est alors nilpotente, donc nulle puisque~$X$ est réduit ; il s'ensuit que~$f=0$. On a de même~$g=0$ si~$X=V(g)$, ce qui
achève la preuve.~$\Box$

\deux{theo-princ-max}
{\bf Théorème}.
{\em Soit~$k$ un corps et soit~$X$ un~$k$-schéma propre, irréductible et réduit. L'anneau~$\sch O_X(X)$ est une extension finie
de~$k$ ; en particulier, $\sch O_X(X)=k$ si~$k$ est algébriquement clos}. 

\medskip
{\em Démonstration}. 
Comme~$X$ est irréductible et réduit, l'anneau~$\sch O_X(X)$ est intègre d'après le lemme~\ref{lemme-schema-integre}.
Soit~$f\in \sch O_X(X)$. Elle induit un~$k$-morphisme~$\psi \colon X\to \Aff^1_k$, caractérisé par le fait que~$\psi^*T=f$
({\em cf}.~\ref{cas-part-ox-mora1}). La composée de~$\psi$ et de l'immersion ouverte~$\Aff^1_k\hookrightarrow \PP^1_k$ 
(obtenue en identifiant~$\Aff^1_k$ à l'une des deux cartes affines standard de~$\PP^1_k$) est un~$k$-morphisme
$X\to \PP^1_k$. Comme le~$k$-schéma~$X$ est propre, et en particulier universellement fermé, et comme~$\PP^1_k$
est séparé sur~$k$,  l'image~$\psi(X)$ est fermée dans~$\PP^1_k$. Étant par ailleurs contenue dans~$\Aff^1_k$, cette image
est nécessairement un ensemble fini de points fermés, et consiste finalement en un unique point fermé~$x$
car~$X$ est irréductible. 

\medskip
Le point fermé~$x$ de~$\Aff^1_k$ correspond à un polynôme irréductible~$P$
de~$k$ ; comme~$X$ est réduit, $\psi$ se factorise par~$\{x\}_{\rm red}=\spec k[T]/P$, 
ce qui veut dire que~$P(f)=P(\psi^*T)=\psi^*P(T)=0$. 

\medskip
L'anneau~$\sch O_X$ est ainsi une~$k$-algèbre intègre dont tous les éléments sont entiers sur~$k$ ; c'est donc un corps
(lemme~\ref{alg-ent-corps}). Il reste à s'assurer qu'elle est de type fini sur~$k$. 

Le~$k$-schéma~$X$ est propre, et donc de type fini. Il est non vide car irréductible, et possède donc
un point fermé~$y$. L'évaluation en~$y$ est un~$k$-morphisme de~$\sch O_X(X)$ dans~$\kappa(y)$, 
injectif puisque~$\sch O_X(X)$ est un corps. Comme~$[\kappa(y):k]<+\infty$, le corps~$\sch O_X(X)$ est une
extension finie de~$k$.

 

\end{document}